\documentclass[12pt,a4paper]{article}
\usepackage[utf8]{inputenc}
\usepackage[T1]{fontenc}
\usepackage{amsmath}
\usepackage{amsthm}
\usepackage{amsfonts}
\usepackage{amssymb}
\usepackage{mathtools}
\usepackage{mathrsfs}
\usepackage{url}
\usepackage{mathdots}
\usepackage[english]{babel}
\usepackage{geometry}
\usepackage{verbatim}
\usepackage{hyperref}
\usepackage{float}
\usepackage{authblk}
\usepackage{enumitem}
\usepackage{tikz-cd}
\usepackage{extarrows}
\usepackage{comment}
\usepackage{bbm}
\usepackage[style=alphabetic, backend=bibtex, backref=true]{biblatex}
\usepackage{csquotes}
\bibliography{IHAg}

\hypersetup{
    bookmarks=true,         % show bookmarks bar?
    unicode=false,          % non-Latin characters in Acrobat<E2><80><99>s bookmarks
    pdftoolbar=true,        % show Acrobat<E2><80><99>s toolbar?
    pdfmenubar=true,        % show Acrobat<E2><80><99>s menu?
    pdffitwindow=false,     % window fit to page when opened
    pdfstartview={FitH},    % fits the width of the page to the window
  pdftitle={The Euler characteristic of \(\ell\)-adic local systems on \(\mathcal{A}_n\)},    % title
    pdfauthor={Olivier Taïbi},     % author
    colorlinks=true,       % false: boxed links; true: colored links
    linkcolor=blue,          % color of internal links
    citecolor=green,        % color of links to bibliography
    filecolor=green,      % color of file links
    urlcolor=cyan}           % color of external links

\newcommand{\op}{\operatorname}
\newcommand{\Q}{\mathbb{Q}}
\newcommand{\A}{\mathbb{A}}
\newcommand{\R}{\mathbb{R}}
\providecommand{\C}{\mathbb{C}}
\renewcommand{\C}{\mathbb{C}}
\newcommand{\N}{\mathbb{N}}
\newcommand{\Qp}{\mathbb{Q}_p}
\newcommand{\Qpbar}{\overline{\mathbb{Q}}_p}
\newcommand{\etabar}{\overline{\eta}}
\newcommand{\Fpbar}{\overline{\mathbb{F}}_p}
\newcommand{\sbar}{\overline{s}}
\newcommand{\Qbar}{\overline{\mathbb{Q}}}
\newcommand{\Qellbar}{\overline{\mathbb{Q}_{\ell}}}
\newcommand{\Qell}{\mathbb{Q}_{\ell}}

\newcommand{\Fp}{\mathbb{F}_p}
\newcommand{\F}{\mathbb{F}}
\newcommand{\Fq}{\mathbb{F}_q}
\newcommand{\Zp}{\mathbb{Z}_p}
\newcommand{\Zq}{\mathbb{Z}_q}
\newcommand{\Zell}{\mathbb{Z}_{\ell}}
\newcommand{\Z}{\mathbb{Z}}

\newcommand{\Sets}{\mathrm{Sets}}
\newcommand{\Perf}{\mathrm{Perf}}
\newcommand{\Sh}{\mathrm{Sh}}
\newcommand{\id}{\mathrm{id}}
\newcommand{\ctf}{\mathrm{ctf}}
\newcommand{\diag}{\mathrm{diag}}
\newcommand{\Sym}{\operatorname{Sym}}
\newcommand{\univ}{\mathrm{univ}}
\newcommand{\ad}{\mathrm{ad}}
\newcommand{\even}{\mathrm{even}}
\newcommand{\Tot}{\operatorname{Tot}}
\newcommand{\Hecke}{\operatorname{Hecke}}
\newcommand{\PreH}{\operatorname{PreH}}
\newcommand{\Ob}{\operatorname{Ob}}
\newcommand{\Rep}{\operatorname{Rep}}

\newcommand{\colim}{\operatorname{colim}}
\newcommand{\sm}{\mathrm{sm}}
\newcommand{\Frob}{\operatorname{Frob}}
\newcommand{\Spec}{\operatorname{Spec}}
\newcommand{\BC}{\operatorname{BC}}

\mathchardef\mhyphen="2D
\newcommand{\corr}{\mathrm{corr}\mhyphen}
\newcommand{\sign}{\mathrm{sign}}

\newcommand{\Frac}{\mathrm{Frac}}

\newcommand{\AJ}{\mathrm{AJ}}
\newcommand{\sesi}{\mathrm{ss}}
\newcommand{\extr}{\mathrm{extr}}
\newcommand{\GSpin}{\mathrm{GSpin}}
\newcommand{\std}{\mathrm{std}}
\newcommand{\Lie}{\operatorname{Lie}}

\newcommand{\IH}{\mathrm{IH}}
\newcommand{\gen}{\mathrm{gen}}
\newcommand{\Spin}{\mathrm{Spin}}
\newcommand{\Sp}{\mathrm{Sp}}
\newcommand{\GSp}{\mathrm{GSp}}
\newcommand{\spin}{\mathrm{spin}}
\newcommand{\Sfrak}{\mathfrak{S}}
\newcommand{\afrak}{\mathfrak{a}}
\newcommand{\efrak}{\mathfrak{e}}

\newcommand{\gfrak}{\mathfrak{g}}
\newcommand{\hfrak}{\mathfrak{h}}
\newcommand{\kfrak}{\mathfrak{k}}
\newcommand{\lfrak}{\mathfrak{l}}
\newcommand{\mfrak}{\mathfrak{m}}
\newcommand{\nfrak}{\mathfrak{n}}
\newcommand{\pfrak}{\mathfrak{p}}
\newcommand{\tfrak}{\mathfrak{t}}
\newcommand{\ufrak}{\mathfrak{u}}

\newcommand{\wfrak}{\mathfrak{w}}
\newcommand{\zfrak}{\mathfrak{z}}
\newcommand{\gfrakhat}{\widehat{\mathfrak{g}}}
\newcommand{\lfrakhat}{\widehat{\mathfrak{l}}}
\newcommand{\Tr}{\operatorname{tr}}
\newcommand{\Cent}{\mathrm{Cent}}
\newcommand{\Ad}{\operatorname{Ad}}
\newcommand{\unr}{\mathrm{unr}}
\newcommand{\nonendo}{\mathrm{ne}}
\newcommand{\Sat}{\mathrm{Sat}}
\newcommand{\disc}{\mathrm{disc}}
\newcommand{\ab}{\mathrm{ab}}
\newcommand{\lin}{\mathrm{lin}}
\newcommand{\her}{\mathrm{her}}

\newcommand{\adm}{\mathrm{adm}}
\newcommand{\adj}{\mathrm{adj}}
\newcommand{\der}{\mathrm{der}}
\newcommand{\cont}{\mathrm{cont}}

\newcommand{\sico}{\mathrm{sc}}
\newcommand{\Out}{\mathrm{Out}}
\newcommand{\GL}{\mathrm{GL}}
\newcommand{\GLbf}{\mathbf{GL}}
\newcommand{\SLbf}{\mathbf{SL}}
\newcommand{\PGLbf}{\mathbf{PGL}}
\newcommand{\SL}{\mathrm{SL}}
\newcommand{\PGL}{\mathrm{PGL}}
\newcommand{\SO}{\mathrm{SO}}
\newcommand{\SU}{\mathrm{SU}}
\newcommand{\End}{\operatorname{End}}
\newcommand{\Aut}{\operatorname{Aut}}
\newcommand{\Hom}{\operatorname{Hom}}

\newcommand{\RuHom}{\mathrm{R\underline{Hom}}}
\newcommand{\Gal}{\mathrm{Gal}}
\newcommand{\GalQ}{\mathrm{Gal}_{\mathbb{Q}}}
\newcommand{\GalQell}{\mathrm{Gal}_{\mathbb{Q}_{\ell}}}
\newcommand{\GalQp}{\mathrm{Gal}_{\mathbb{Q}_p}}
\newcommand{\Res}{\operatorname{Res}}
\newcommand{\res}{\operatorname{res}}
\newcommand{\cores}{\operatorname{cores}}
\newcommand{\Ind}{\mathrm{Ind}}
\newcommand{\ind}{\mathrm{ind}}

\newcommand{\Cscr}{\mathscr{C}}
\newcommand{\Wscr}{\mathscr{W}}
\newcommand{\Acal}{\mathcal{A}}
\newcommand{\Bcal}{\mathcal{B}}
\newcommand{\Ccal}{\mathcal{C}}

\newcommand{\Ecal}{\mathcal{E}}
\newcommand{\Gcal}{\mathcal{G}}
\newcommand{\Hcal}{\mathcal{H}}
\newcommand{\Lcal}{\mathcal{L}}
\newcommand{\Mcal}{\mathcal{M}}
\newcommand{\Ncal}{\mathcal{N}}
\newcommand{\Ocal}{\mathcal{O}}
\newcommand{\Rcal}{\mathcal{R}}
\newcommand{\Scal}{\mathcal{S}}
\newcommand{\Tcal}{\mathcal{T}}
\newcommand{\Ucal}{\mathcal{U}}
\newcommand{\Xcal}{\mathcal{X}}
\newcommand{\Dbb}{\mathbb{D}}
\newcommand{\PsCo}{\mathcal{PC}}
\newcommand{\Pcal}{\mathcal{P}}
\newcommand{\Qcal}{\mathcal{Q}}
\newcommand{\ICcal}{\mathcal{IC}}
\newcommand{\ICcalt}{\tilde{\mathcal{IC}}}
\newcommand{\IC}{\mathrm{IC}}
\newcommand{\Fcal}{\mathcal{F}}
\newcommand{\AFcal}{\mathcal{AF}}
\newcommand{\Std}{\mathrm{Std}}
\newcommand{\Hbf}{\mathbf{H}}
\newcommand{\Hhat}{\widehat{\mathbf{H}}}
\newcommand{\Gbf}{\mathbf{G}}
\newcommand{\Cbf}{\mathbf{C}}
\newcommand{\Dbf}{\mathbf{D}}
\newcommand{\Kbf}{\mathbf{K}}
\newcommand{\Rbf}{\mathbf{R}}
\newcommand{\Obf}{\mathbf{O}}
\newcommand{\SObf}{\mathbf{SO}}
\newcommand{\Spbf}{\mathbf{Sp}}
\newcommand{\GSpbf}{\mathbf{GSp}}
\newcommand{\PGSpbf}{\mathbf{PGSp}}
\newcommand{\GSObf}{\mathbf{GSO}}
\newcommand{\PGSObf}{\mathbf{PGSO}}
\newcommand{\Gbfad}{\mathbf{G}_{\mathrm{ad}}}
\newcommand{\Ghat}{\widehat{\mathbf{G}}}
\newcommand{\Lhat}{\widehat{\mathbf{L}}}
\newcommand{\Mhat}{\widehat{\mathbf{M}}}
\newcommand{\That}{\widehat{\mathbf{T}}}
\newcommand{\nuhat}{\widehat{\nu}}
\newcommand{\thetahat}{\widehat{\theta}}
\newcommand{\Pbf}{\mathbf{P}}
\newcommand{\Abf}{\mathbf{A}}
\newcommand{\Lbf}{\mathbf{L}}
\newcommand{\Mbf}{\mathbf{M}}
\newcommand{\Qbf}{\mathbf{Q}}
\newcommand{\Nbf}{\mathbf{N}}
\newcommand{\Bbf}{\mathbf{B}}
\newcommand{\Tbf}{\mathbf{T}}
\newcommand{\Ubf}{\mathbf{U}}
\newcommand{\Zbf}{\mathbf{Z}}
\newcommand{\Zhat}{\widehat{\mathbb{Z}}}
\newcommand{\cpsc}{c_{p, \mathrm{sc}}}
\newcommand{\Mpsi}{\mathcal{M}_{\psi}}
\newcommand{\Mpsisc}{\mathcal{M}_{\psi, \mathrm{sc}}}
\newcommand{\GMpsisc}{\mathcal{GM}_{\psi, \mathrm{sc}}}
\newcommand{\dpsitau}{\dot{\psi}_{\tau}}
\newcommand{\dpsitausc}{\dot{\psi}_{\tau, \mathrm{sc}}}
\newcommand{\Psit}{\widetilde{\Psi}}
\newcommand{\taut}{\widetilde{\tau}}
\newcommand{\pr}{\mathrm{pr}}
\newcommand{\prDS}{\mathrm{pr}_{\mathrm{DS}}}

\providecommand{\Re}{\operatorname{Re}}

\newcommand{\ol}[1]{\overline{#1}}
\newcommand{\ul}[1]{\underline{#1}}
\newcommand{\wt}[1]{\widetilde{#1}}

\newtheorem{theo}{Theorem}[subsection]
\newtheorem{lemm}[theo]{Lemma}
\newtheorem{coro}[theo]{Corollary}
\newtheorem{defi}[theo]{Definition}
\newtheorem{prop}[theo]{Proposition}
\newtheorem{prodef}[theo]{Proposition-Definition}
\newtheorem{rema}[theo]{Remark}
\newtheorem{exam}[theo]{Example}
\newtheorem{assu}[theo]{Assumption}

\newtheorem{theointro}{Theorem}
\newtheorem{conjintro}{Conjecture}

\numberwithin{equation}{subsection}

\begin{document}

\baselineskip=16pt

\author{Olivier Taïbi}

\title{The Euler characteristic of \(\ell\)-adic local systems on \(\mathcal{A}_n\)}

\maketitle

\begin{abstract}
  We study the Euler characteristic of \(\ell\)-adic local systems on the moduli stack \(\Acal_n\) of principally polarized abelian varieties of dimension \(n\) associated to algebraic representations of \(\GSpbf_{2n}\), as virtual representations of the absolute Galois group of \(\Q\) and the unramified Hecke algebra of \(\GSpbf_{2n}\).
  To this end we take the last steps of the Ihara-Langlands-Kottwitz method to compute the intersection cohomology of minimal compactifications of Siegel modular varieties in level one, following work of Kottwitz and Morel, proving an unconditional reformulation of Kottwitz' conjecture in this case.
  This entails proving the existence of \(\GSpin\)-valued Galois representations associated to certain level one automorphic representations for \(\PGSpbf_{2n}\) and \(\SObf_{4n}\).
  As a consequence we prove the existence of \(\GSpin\)-valued Galois representations associated to level one Siegel eigenforms, a higher genus analogue of theorems of Deligne (genus one) and Weissauer (genus two).
  Using Morel's work and Franke's spectral sequence we derive explicit formulas expressing the Euler characteristic of compactly supported cohomology of automorphic \(\ell\)-adic local systems on Siegel modular varieties in terms of intersection cohomology.
  Specializing to genus three and level one, we prove an explicit conjectural formula of Bergström, Faber and van der Geer for the compactly supported Euler characteristic in terms of spin Galois representations associated to level one Siegel cusp forms.
  Specializing to trivial local systems we give explicit formulas for the number of points of \(\Acal_n\) over finite fields for all \(n \leq 7\).
\end{abstract}

\newpage
\setcounter{tocdepth}{2}
\tableofcontents
\newpage

\section{Introduction}

For an integer \(n \geq 1\) let \(\Acal_n\) be the moduli stack of \(n\)-dimensional principally polarized abelian varieties, a smooth Deligne-Mumford stack over \(\Z\) of relative dimension \(n(n+1)/2\).
We recall the precise moduli problem in Section \ref{sec:def_An}.
For a prime number \(\ell\) there is a natural functor \(\Fcal\) from the category of finite-dimensional algebraic representations of \(\GSpbf_{2n,\Qell}\) (``conformal symplectic group'' over \(\Qell\)) to the category of \(\ell\)-adic sheaves on \(\Acal_{n,\Z[\ell^{-1}]} := \Acal_n \times_{\Z} \Z[\ell^{-1}]\) (see Section \ref{sec:local_systems}).
Denote by \(\GalQ := \Gal(\Qbar/\Q)\) the absolute Galois group of \(\Q\).
Denote by \(\Hcal^{\unr}(\GSpbf_{2n})\) the (unramified) Hecke algebra of \(\GSpbf_{2n}(\A_f)\) in level \(\GSpbf_{2n}(\Zhat)\) (with rational coefficients).
The main goal of this paper is to prove a formula expressing, in the Grothendieck group of finite-dimensional continuous \(\Qell\)-representations of the absolute Galois group \(\GalQ\) with commuting action of \(\Hcal^{\unr}(\GSpbf_{2n})\), the Euler characteristics
\begin{equation} \label{eq:intro_e_c_An_FV}
  e_c(\Acal_{n,\Qbar}, \Fcal(V)) := \sum_{i=0}^{n(n+1)} (-1)^i \left[ H^i_c(\Acal_{n,\Qbar}, \Fcal(V)) \right]
\end{equation}
in terms of \(\ell\)-adic Galois representations associated to certain automorphic representations.
This problem reduces to the case where the representation \(V\) of \(\GSpbf_{2n,\Qell}\) is irreducible, which we assume for the rest of this introduction.
For any prime \(p \neq \ell\) the virtual Galois representation \eqref{eq:intro_e_c_An_FV} is unramified at \(p\) and its restriction to the decomposition group at \(p\) equals \(e_c(\Acal_{n,\Fpbar}, \Fcal(V))\) (see Proposition \ref{pro:spe_gal_hecke_coho}).
When \(V\) is trivial the Grothendieck-Lefschetz trace formula tells us that knowing the Euler characteristic \(e_c(\Acal_{n,\Fpbar}, \Qell) \in K_0(\Rep_{\Qell}^{\cont}(\GalQ))\) (forgetting the action of \(\Hcal^{\unr}(\GSpbf_{2n})\)) is equivalent to knowing the (weighted) counts
\[ |\Acal_n(\mathbb{F}_{p^m})| := \sum_{(A,\lambda) \in \Acal_n(\mathbb{F}_{p^m}) / \sim} |\Aut(A,\lambda)|^{-1} \]
for all integers \(m \geq 1\).
As a first application of our main results and \cite[Theorem 9.3.3]{CheLan} we obtain the following explicit formulas.

\begin{theointro} \label{thmintro:card_An_Fq_small_n}
  For \(1 \leq n \leq 6\) and any prime number \(\ell\) the virtual representation \(e_c(\Acal_{n,\Qbar}, \Qell)\) of \(\GalQ\) (forgetting the Hecke action) is Tate, equivalently there exists a polynomial \(P_n \in \Z[X]\) such that for any prime power \(q\) we have \(|\Acal_n(\Fq)| = P_n(q)\).
  More precisely we have
  \begin{align*}
    |\Acal_1(\Fq)| &= q, \\
    |\Acal_2(\Fq)| &= q^3 + q^2, \\
    |\Acal_3(\Fq)| &= q^6 + q^5 + q^4 + q^3 + 1, \\
    |\Acal_4(\Fq)| &= q^{10} + q^9 + q^8 + 2q^7 + q^6 + q^5 + q^4 + q,  \\
    |\Acal_5(\Fq)| &= q^{15} + q^{14} + q^{13} + 2q^{12} + 2q^{11} + 2q^{10} + 2q^9 + 2q^8 + q^7 + q^6 + q^5 + q^3 + q^2, \\
    |\Acal_6(\Fq)| &= q^{21} + q^{20} + q^{19} + 2q^{18} + 2q^{17} + 4q^{16} + 4q^{15} + 4q^{14} + 5q^{13} + 4q^{12} \\
                   &\ \ \ + 4q^{11} + 3q^{10} + 2q^9 + q^8 + 2q^7 + 2q^6 + q^5 + q^4 + q^3 + 1.
  \end{align*}

  For \(n=7\) this Euler characteristic is not Tate, more precisely for any prime power \(q=p^m\) we have
  \begin{align*}
    |\Acal_7(\Fq)| =& \ q^{28} + q^{27} + q^{26} + 2q^{25} + 2q^{24} + 3q^{23} + 4q^{22} + 4q^{21} + 4q^{20} + 6q^{19} + 7q^{18} \\
                    & \ + 8q^{17} + 7q^{16} + 6q^{15} + 5q^{14} + 4q^{13} + 4q^{12} + 2q^{11} + 3q^{10} + 4q^9 + 3q^8 \\
                    & \ + 3q^7 + q^6 + q^5 + 2q^4 + q \\
                    & \ + (q^6 + q^5 + q^4 + q^3 + 1)  \times a(p^m)
  \end{align*}
  where the family of integers \((a(p^m))_{m \geq 0}\) is defined by the equality in \(\Z[[T]]\)
  \[ \sum_{m \geq 0} a(p^m) T^m = \frac{3 - 2 \tau(p^2) T + p^{11} \tau(p^2) T^2}{1 - \tau(p^2) T + p^{11} \tau(p^2) T^2 - p^{33} T^3} \]
  where \(\tau(p^2)\) is a coefficient of the \(q\)-expansion of the Ramanujan \(\Delta\) function
  \[ \sum_{m \geq 0} \tau(m) q^m = q \prod_{m \geq 1} (1-q^m)^{24}. \]
\end{theointro}

These formulas seem to be new for \(n>3\), see \cite[Theorem 8.1]{BFG} for \(n=3\) (see also \cite{Hain_A3}).
In principle one can give similar explicit formulas for all \(n \leq 12\) (but with more complicated ingredients than just \(\Delta_{11}\)), see Remark \ref{rem:card_An_Fq_upto12}.

For a non-trivial representation \(V\), via the geometric construction of \(\Fcal(V)\) using the universal abelian variety there is a similar interpretation of the Euler characteristic \(e_c(\Acal_{n,\Fpbar}, \Fcal(V))\) using point counting: see \cite[\S 8]{BFG}.
Irreducible representations of \(\GSpbf_{2n,\Qell}\) come by extension of scalars from irreducible representations of \(\GSpbf_{2n,\Q}\), and this interpretation implies in particular that for an irreducible representation \(V\) of \(\GSpbf_{2n,\Q}\) the traces (for \(n \in \Z\))
\[ \Tr(\Frob_p^n \,|\, e_c(\Acal_{n,\Qbar}, \Fcal(\Qell \otimes_\Q V))) \]
are all rational and do not depend on the choice of \(\ell \neq p\).

Our results are motivated by \cite[Conjecture 7.1]{BFG}, which conjectures an explicit expression for \eqref{eq:intro_e_c_An_FV} (forgetting the Hecke action) in the case \(n=3\) in terms of \(2^g\)-dimensional ``spin'' \(\ell\)-adic Galois representations conjecturally associated to Siegel eigenforms for \(\Spbf_{2g}(\Z)\), \(g \leq 3\).
(In fact their conjecture is at the level of motives over \(\Q\), and we are only considering \(\ell\)-adic realizations.)
Bergström, Faber and van der Geer arrived at this conjectural formula using explicit point counts.
This conjecture follows similar results in genus \(\leq 2\) recalled below.
We first explain the application of our main results to the existence of spin Galois representations for level one Siegel modular forms and to Conjecture 7.1 of \cite{BFG}, before explaining our main results concerning the cohomology of intersection complexes on the minimal compactification of \(\Acal_n\) and explicit formulas relating it to compactly supported cohomology on \(\Acal_n\).

\subsection{Spin Galois representations for level one Siegel modular forms}
\label{sec:intro_SMF}

As a corollary of one of our main results we prove the existence of \(\GSpin\)-valued \(\ell\)-adic Galois representations in higher genus, as we now explain.
Applying Schur functors to the Hodge bundle of \(\Acal_{n,\C}\) (an \(n\)-dimensional vector bundle) and taking global sections ``vanishing at infinity'' yields, for an irreducible algebraic representation of \(\GLbf_n\) parametrized by its highest weight \(\ul{k} = (k_1 \geq \dots \geq k_n)\), the finite-dimensional vector space \(S_{\ul{k}}(\Spbf_{2n}(\Z))\) of Siegel cusp forms of weight \(\ul{k}\) and level \(\Spbf_{2n}(\Z)\) (the precise definition is recalled in \cite{vdG_123}).
It is endowed with an action\footnote{More precisely we incorporate the extra factor \(\eta(\gamma)^{\sum \lambda_i - g(g+1)/2}\) in \cite[Definition 8]{vdG_123}.} of the (commutative) Hecke algebra \(\Hcal^{\unr}(\GSpbf_{2n})\).
Characters of \(\Hcal^{\unr}(\GSpbf_{2n})_\C := \C \otimes_\Q \Hcal^{\unr}(\GSpbf_{2n})\) correspond via the Satake isomorphism to families, indexed by the set of all prime numbers, of semi-simple conjugacy classes in \(\GSpin_{2n+1}(\C)\).
In particular to an eigenform \(f\) is associated a family \((c_p(f))_p\) of such conjugacy classes.
For better rationality properties (see \cite[\S 8]{Gross_Satake} for details) it is convenient to consider the family \((c_p^\mathrm{arith}(f))_p\) defined by \(c_p^\mathrm{arith}(f) = p^{n(n+1)/4} c_p(f)\) instead.
There is a morphism (of split connected reductive groups over \(\Q\)) \(\GSpin_{2n+1} \to \SO_{2n+1}\) with kernel identified to \(\GL_1\), as well as a natural morphism \(\beta: \GSpin_{2n+1} \to \GL_1\) which is \(t \mapsto t^2\) on \(\GL_1 \simeq Z(\GSpin_{2n+1})\).
For an eigenform \(f \in S_{\ul{k}}(\Spbf_{2n}(\Z))\) and a prime number \(p\) we have \(\beta(c_p^\mathrm{arith}(f)) = p^{\sum_{i=1}^n k_i - n(n+1)/2}\).

We will be particularly interested in two irreducible representations of the algebraic group \(\GSpin_{2n+1}\):
\begin{itemize}
\item the \(2n+1\)-dimensional standard representation \(\Std\), which factors through \(\SO_{2n+1}\),
\item the \(2^n\)-dimensional spin representation \(\spin\), which maps \(z \in \GL_1 \simeq Z(\GSpin_{2n+1})\) to \(z \, \id\).
\end{itemize}
Putting together known results due to many mathematicians (Theorem \ref{thm:existence_rho_SO}, applied using Arthur's endoscopic classification for \(\Spbf_{2n}\): see \cite[\S 9]{ChRe} or \cite[\S 5]{Taibi_dimtrace}) we know that for any field isomorphism\footnote{As usual only the restriction of this isomorphism to the algebraic closure of \(\Q\) plays a role.} \(\iota: \C \simeq \Qellbar\) and any eigenform \(f \in S_{\ul{k}}(\Spbf_{2n}(\Z))\) of weight \(\ul{k}\) satisfying \(k_n \geq n+1\), there exists a unique continuous semi-simple morphism \(\rho_{f,\iota}^{\SO}: \GalQ \to \SO_{2n+1}(\Qellbar)\) which is unramified away from a finite set of prime numbers and such that for almost all primes \(p\) the semi-simplification of \(\rho_{f,\iota}^{\SO}(\Frob_p)\) (here \(\Frob_p\) denotes the geometric Frobenius element) is equal to the projection of \(\iota(c_p^\mathrm{arith}(f))\) along \(\GSpin_{2n+1} \to \SO_{2n+1}\).
We even know that \(\rho_{f,\iota}^{\SO}\) is unramified away from \(\ell\) and that this relation holds at all primes \(p \neq \ell\), and that \(\rho_{f,\iota}^{\SO}\) is crystalline at \(\ell\).
As a corollary of one of our main results we obtain the existence of spin Galois representations associated to level one Siegel eigenforms, in the form of the following theorem.

\begin{theointro}[Corollary \ref{cor:GSpin_Gal_rep_Siegel}] \label{thmintro:Gal_rep_SMF}
  Let \(n \geq 1\) and \(k_1 \geq \dots \geq k_n \geq n+1\) be integers, and denote \(\ul{k} = (k_1, \dots, k_n)\).
  Let \(f \in S_{\ul{k}}(\Spbf_{2n}(\Z))\) be an eigenform.
  Let \(\ell\) be a prime number, and choose \(\iota: \C \simeq \Qellbar\).
  \begin{enumerate}
  \item There exists a unique continuous lift \(\rho_{f,\iota}^{\GSpin}: \GalQ \to \GSpin_{2n+1}(\Qellbar)\) of \(\rho_{f,\iota}^{\SO}\) which is unramified away from \(\ell\), crystalline at \(\ell\) and which satisfies \(\beta \circ \rho_{f,\iota}^{\GSpin} = \chi_\ell^{n(n+1)/2-\sum_{i=1}^n k_i}\) where \(\chi_\ell: \GalQ \to \Qell^\times\) is the \(\ell\)-adic cyclotomic character.
  \item For any prime number \(p \neq \ell\) the semi-simplification of \(\rho_{f,\iota}^{\GSpin}(\Frob_p)\) belongs to \(\iota(c_p^\mathrm{arith}(f))\).
  \item Any continuous semi-simple \(\rho: \GalQ \to \GSpin_{2n+1}(\Qellbar)\) unramified away from a finite set of primes and such that the semi-simplification of \(\rho_{f,\iota}^{\GSpin}(\Frob_p)\) belongs to \(\iota(c_p^\mathrm{arith}(f))\) for almost all \(p\) is conjugated to \(\rho_{f,\iota}^{\GSpin}\).
  \end{enumerate}
\end{theointro}

The first and third points are an easy consequence of results of Patrikis and Conrad (see Proposition \ref{pro:ex_lift_cond_one}) and the Kronecker-Weber theorem, so our main contribution to this theorem is the second point, which we will prove somewhat indirectly.
For \(n=1\) (resp.\ \(n=2\)) we have \(\GSpin_3 \simeq \GL_2\) (resp.\ \(\GSpin_5 \simeq \GSp_4\)) and Theorem \ref{thmintro:Gal_rep_SMF} is the level one case of \cite{Deligne_GalGL2} (resp.\ \cite{Weissauer_GalGSp4}).
For \(n>2\) the existence of these Galois representations is new.

\subsection{A conjecture of Bergström, Faber and van der Geer}
\label{sec:intro_BFG}

Following \cite{BFG} for a weight \(\ul{k} = (k_1, \dots, k_n)\) satisfying \(k_n \geq n+1\) we define\footnote{In fact for \(n=1\) following \cite{BFG} we will use a slightly different definition in the weight \(2\) case, see Section \ref{sec:BFG}.} the \(\ell\)-adic representation of \(\GalQ\)
\[ S[\ul{k}]_\ell = \bigoplus_f \spin \circ \rho_{f,\iota}^{\GSpin} \]
where the sum ranges over eigenforms \(f \in S_{\ul{k}}(\Spbf_{2n}(\Z))\).
One could show, from our method to construct the morphisms \(\rho_{f,\iota}^{\GSpin}\), that \(S[\ul{k}]_\ell\) does not depend on the choice of \(\iota\), as the notation suggests.
Bergström, Faber and van der Geer conjecture \cite[\S 5]{BFG} the existence of a motive \(S[\ul{k}]\) over \(\Q\) whose \(\ell\)-adic realization is isomorphic to \(S[\ul{k}]_\ell\).
For \(n \leq 3\) the fact that \(S[\ul{k}]_\ell\) may be defined over \(\Qell\) and does not depend on the choice of \(\iota\) follows from our results in Section \ref{sec:BFG}.
Let \(\lambda = (\lambda_1 \geq \dots \geq \lambda_n)\) be a dominant weight for \(\Spbf_{2n}\), corresponding to an irreducible algebraic representation \(V_\lambda\) of \(\Spbf_{2n,\Qell}\).
In the moduli interpretation there is a natural choice of extension \(V_{\lambda,0}\) of \(V_\lambda\) to \(\GSpbf_{2n,\Qell}\): letting \(z\) in the center \(\Zbf(\GSpbf_{2n,\Qell}) \simeq \GLbf_{1,\Qell}\) act by \(z^{-\sum_i \lambda_i} \, \id\).
As recalled in Section \ref{sec:local_systems} the local system \(\Fcal(V_{\lambda,0})\) is pure of weight \(\sum_i \lambda_i\) and is ``effective''.
The authors of \cite{BFG} define in \S 5 loc.\ cit.\ the extraneous contribution \(e_{n,\extr}(\lambda)_\ell \in K_0(\Rep_{\Qellbar}^{\cont}(\GalQ))\) by the equation
\[ e_c(\Acal_{n,\Qbar}, \Fcal(V_{\lambda,0})) = (-1)^{n(n+1)/2} S[\ul{k}]_\ell + e_{n,\extr}(\lambda)_\ell, \]
the idea being that \(e_c(\Acal_{n,\Qbar}, \Fcal(V_{\lambda,0}))\) should be equal to \((-1)^{n(n+1)/2} S[\ul{k}]_\ell\) up to ``smaller'' error terms (endoscopic or related to the boundary).
For \(n=1\) and \(\lambda_1>0\) we simply have \(e_{1,\extr}(\lambda)_\ell=-1\) (see \cite[Theorem 2.3]{BFG}).
For \(n=2\) Faber and van der Geer conjectured an explicit formula for \(e_{2,\extr}(\lambda)_\ell\) (in fact, of the conjectural virtual motive over \(\Q\) whose \(\ell\)-adic realization should be \(e_{2,\extr}(\lambda)_\ell\)) in terms of \(S[-]_\ell\) (in genus one), recalled in Section \ref{sec:BFG_genus_two} (see also \cite[\S 6.3]{BFG}).
This conjecture was later proved by Weissauer and van der Geer in the regular case (see the discussion after Conjecture 6.1 in \cite{BFG}), and by Petersen \cite{Petersen} in general.
For \(n=3\), using explicit point counts over finite fields the authors of \cite{BFG} conjectured an explicit formula for \(e_{3,\extr}(\lambda)_\ell\) (again, their conjectural formula is motivic), see Conjecture 7.1 loc.\ cit.

\begin{theointro}[Theorem \ref{thm:BFG}] \label{theointro:BFG}
  Conjecture 7.1 of \cite{BFG} holds true at the level of \(\ell\)-adic Galois representations.
\end{theointro}

We will see that no similar formula can be expected in genus \(>3\): the extraneous term is not expressed just in terms of lower-dimensional Galois representations \(S[-]_\ell\) alone.

\subsection{Intersection cohomology of local systems on the minimal compactification of \(\Acal_n\) and \(\GSpin\)-valued Galois representations}
\label{sec:intro_IH_GSpin_Gal_rep}

We now explain our main results, which hold for an arbitrary genus \(n \geq 1\).
Our first goal is to prove a special case of a conjecture of Kottwitz \cite{Kottwitz_AA} describing, for a Shimura datum \((\Gbf,\Xcal)\) with associated Shimura tower \((\Sh(\Gbf,\Xcal,K))_K\) (quasi-projective varieties over the reflex field \(E\)) and minimal compactifications \(\Sh(\Gbf,\Xcal,K) \hookrightarrow \Sh(\Gbf,\Xcal,K)^*\), and an algebraic representation \(V\) of \(\Gbf_{\Qell}\) with associated intersection complex \(\IC(V)\) on \(\Sh(\Gbf,\Xcal,K)^*\) (the intermediate extension of \(\Fcal(V)\)), the representations
\begin{equation} \label{eq:intro_IH}
  \IH^i(\Gbf, \Xcal, V) := \varinjlim_K H^i(\Sh(\Gbf,\Xcal,K)^*_{\ol{E}}, \IC^K(V))
\end{equation}
of \(\Gbf(\A_f) \times \Gal_E\), in terms of the conjectural global Langlands correspondence (more precisely, Arthur's conjectures and their \(\ell\)-adic realizations), which we now recall without going into full details.
The special case relevant to our situation is the one where \((\Sh(\Gbf,\Xcal,K))_K\) is the tower \((\Acal_{n,K})_K\) of Siegel modular varieties, \(K\) ranges over (neat) compact open subgroups of \(\Gbf=\GSpbf_{2n}(\A_f)\), and we take \(\GSpbf_{2n}(\Zhat)\)-invariants in \eqref{eq:intro_IH}.
In this case \(\Gbf\) is split and the reflex field \(E\) is simply \(\Q\).
These properties simplify the general discussion in \cite{Kottwitz_AA} a little, so we assume that they hold for the rest of this section.

\subsubsection{Kottwitz' conjecture in general}

We temporarily assume the existence of the\footnote{As for absolute Galois groups or Weil groups, the Langlands group should be associated to a choice of algebraic closure.} Langlands group \(L_\Q\) of \(\Q\), a topological group together with a continuous surjective morphism onto the Weil group \(W_\Q\) such that the kernel is compact and connected (i.e.\ a projective limit of compact connected Lie groups).
Among other extra data it should come with an embedding \(W_\R \to L_\Q\) (well-defined up to conjugacy), where \(W_\R\) is the Weil group of \(\R\) (an extension of \(\Gal(\C/\R)\) by \(\C^\times\)), such that the composition \(W_\R \to L_\Q \to W_\Q\) is the usual embedding.
Similarly at the non-Archimedean places we should have embeddings \(W_{\Qp} \times \SU(2) \hookrightarrow L_\Q\) for all primes \(p\).
A conjecture of Langlands predicts a bijection between isomorphism classes of irreducible continuous representations \(\varphi: L_\Q \to \GL_N(\C)\) and cuspidal automorphic representations for \(\GLbf_{N,\Q}\).

Kottwitz' conjecture involves Arthur-Langlands parameters \(\psi: L_\Q \times \SL_2(\C) \to \Ghat(\C)\), i.e.\ continuous semi-simple morphisms which are holomorphic on the factor \(\SL_2(\C)\) and whose centralizer \(C_\psi := \Cent(\psi, \Ghat)\) is finite modulo \(Z(\Ghat)\), whose restriction to \(W_\R\) (along the diagonal embedding, using \(W_\R \hookrightarrow \SL_2(\C), w \mapsto \diag(|w|^{1/2}, |w|^{-1/2})\)) has infinitesimal character (see Definition \ref{def:inf_char}) opposite to that of \(V\).
We denote by \(\Psi(\Gbf,V)\) the set of \(\Ghat(\C)\)-conjugacy classes of such parameters.
For \(\psi \in \Psi(\Gbf,V)\) the centralizer \(C_\psi\) is abelian and the finite group \(C_\psi/Z(\Ghat)\) is \(2\)-torsion.
To the Shimura datum \((\Gbf, \Xcal)\) is associated a conjugacy class of cocharacters \(\mu: \GLbf_{1,\C} \to \Gbf_\C\) and thus a representation \(r_{-\mu}: \Ghat \to \GL(Y)\) of the Langlands dual group \(\Ghat\) having extremal weight \(-\mu\).
We thus obtain a representation \(L_\Q \times \SL_2(\C) \times C_\psi \to \GL(Y)\), and we have a decomposition \(Y = \bigoplus_\nu Y_\nu\) where the sum ranges over the set \(N\) of characters \(\nu\) of \(C_\psi\) whose restriction to \(Z(\Ghat)\) is \(-\mu|_{Z(\Ghat)}\).
To \(\psi \in \Psi(\Gbf,V)\) is conjecturally associated a packet (multiset) \(\Pi_f(\psi)\) of irreducible representations of \(\Gbf(\A_f)\) together with a map
\begin{align*}
  \Pi_f(\psi) \times N & \longrightarrow \Z_{\geq 0} \\
  (\pi_f, \nu) & \longmapsto m(\psi, \pi_f, \nu),
\end{align*}
using Arthur's conjectures (see \cite[p.\ 200]{Kottwitz_AA}).
Denoting by \(d\) the dimension of the Shimura varieties \(\Sh(\Gbf, \Xcal, K)\) we consider the representations
\begin{equation} \label{eq:def_Y_psi_pi_f}
  Y(\psi, \pi_f) := \bigoplus_\nu |\cdot|^{-d/2} Y_\nu^{\oplus m(\psi, \pi_f, \nu)}
\end{equation}
of \(L_\Q \times \SL_2(\C)\) as representations of \(L_\Q\) using the diagonal embedding \(L_\Q \hookrightarrow L_\Q \times \SL_2(\C), g \mapsto (g, \diag(|g|^{1/2},|g|^{-1/2}))\).
These representations of \(L_\Q\) are algebraic, i.e.\ their restriction to \(W_\R\) decompose as direct sums of characters \(z \mapsto z^a \ol{z}^b\) with \(a,b \in \Z\).
This conjecturally implies that they are motivic, in particular for \(\iota: \C \simeq \Qellbar\) they should have \(\ell\)-adic realizations \(Y(\psi, \pi_f)_\iota\), vector spaces over \(\Qellbar\) with a continuous action of \(\GalQ\).
These realizations are characterized by the following compatibility: for almost all primes \(p \neq \ell\) the restriction of \(Y(\psi, \pi_f)\) to \(W_{\Qp} \times \SU(2)\) is trivial on \(I_{\Qp} \times \SU(2)\), where \(I_{\Qp}\) denotes the inertia subgroup of \(W_{\Qp}\), and for all such \(p\) the representation \(Y(\psi, \pi_f)_\iota\) should be unramified at \(p\) and we should have an equality of characteristic polynomials in\footnote{The coefficients of these polynomials are expected to be algebraic over \(\Q\).} \(\Qellbar[T]\)
\[ \det(T \id - \Frob_p \,|\, Y(\psi, \pi_f)_\iota) = \iota \left( \det(T \id - \Frob_p \,|\, Y(\psi, \pi_f)) \right). \]

\begin{conjintro}[{Kottwitz' conjecture \cite[\S 10]{Kottwitz_AA}, in our simplified setting}] \label{conj:Kottwitz}
  We should have an isomorphism of representations of \(\GalQ \times \Gbf(\A_f)\)
  \[ \bigoplus_i \Qellbar \otimes_{\Qell} \IH^i(\Gbf, \Xcal, V) \simeq \bigoplus_\psi \bigoplus_{\pi_f \in \Pi_f(\psi)} Y(\psi, \pi_f)_\iota \otimes_{\Qellbar} \iota(\pi_f) \]
  where \(\iota(\pi_f) = \Qellbar \otimes_{\iota,\C} \pi_f\).
\end{conjintro}

By purity this conjecture also characterizes the individual intersection cohomology groups \(\IH^i(\Gbf, \Xcal, V)\).
For example tempered parameters in \(\Psi(\Gbf, V)\), i.e.\ those parameters which are trivial on the factor \(\SL_2(\C)\), only contribute to middle degree (\(d\)) intersection cohomology.
Purity also implies that the Euler characteristic of \(\IH^\bullet(\Gbf, \Xcal, V)\) determines each \(\IH^i(\Gbf, \Xcal, V)\), a property which does not hold for compactly supported (or ordinary) cohomology.

\subsubsection{The Ihara-Langlands-Kottwitz method}

The Ihara-Langlands-Kottwitz method is a strategy to prove Conjecture \ref{conj:Kottwitz}, assuming the spectral expansion (``stable multiplicity formula'') of the stabilization of the trace formula for certain elliptic endoscopic groups of \(\Gbf\).
Very roughly, this strategy consists of three steps.
\begin{enumerate}
\item Using a generalization of the Grothendieck-Lefschetz trace formula (Deligne's conjecture, proved by Pink \cite{Pink_localterms} and Fujiwara \cite{Fujiwara_Deligne_conj}) and a group-theoretic description of points of Shimura varieties over finite fields, obtain an expression for the trace of a Frobenius element composed with a Hecke operator (satisfying certain assumptions) on
  \begin{equation} \label{eq:intro_e_c_G_X_V}
    e_c(\Gbf, \Xcal, V) := \sum_i (-1)^i \left[ \varinjlim_K H^i_c(\Sh(\Gbf, \Xcal, K), \Fcal(V)) \right]
  \end{equation}
  resembling the elliptic part of the trace formula for \(\Gbf\) (with a twist).
\item Express intersection cohomology in terms of compactly supported cohomology (for certain \(\ell\)-adic local systems on strata of minimal compactifications).
\item Stabilize (in the sense of the stabilization of the trace formula) the expression obtained by putting together the first two steps, in order to express the trace of a Frobenius element composed with a Hecke operator (again, satisfying certain assumptions) on
  \begin{equation} \label{eq:intro_e_IH_G_X_V}
    e_\IH(\Gbf, \Xcal, V) := \sum_i (-1)^i \left[ \IH^i(\Gbf, \Xcal, V) \right]
  \end{equation}
  in terms of the stabilization of the trace formula for certain elliptic endoscopic groups of \(\Gbf\).
\end{enumerate}

\subsubsection{The case of Siegel modular varieties}

For the case of Siegel modular varieties the first step was achieved by Kottwitz himself \cite{KottwitzPoints} (which concerns more generally PEL type Shimura varieties, see also \cite[\S 12]{Kottwitz_AA}), and the second and third step were achieved by Morel \cite{MorelSiegel1}, \cite{MorelSiegel2} (for the third step Kottwitz had already stabilized the ``elliptic part'' corresponding to the contribution of compactly supported cohomology of the open Shimura variety, see \cite[Theorem 7.2]{Kottwitz_AA}).
So it would seem that in the case of Siegel modular varieties we are already close to knowing Kottwitz' Conjecture \ref{conj:Kottwitz}.

To conclude however, we would like to know the spectral expansion of the stabilization of the trace formula for (certain) elliptic endoscopic groups of \(\GSpbf_{2n}\).
By \cite[Proposition 2.1.1]{MorelSiegel2} the relevant endoscopic groups are isomorphic to \(\Gbf(\Spbf_{2a} \times \SObf_{4b})\), where \(a+2b=n\) and for any \(\Q\)-algebra \(R\)
\[ \Gbf(\Spbf_{2a} \times \SObf_{4b})(R) = \{ (g_1, g_2) \in \GSpbf_{2a}(R) \times \GSObf_{4b}(R) \,|\, c(g_1) = c(g_2) \}, \]
denoting by \(c\) the similitude characters.
Of course the spectral expansion is not currently known in a strict sense, because the existence of the Langlands group \(L_\Q\) itself is still conjectural.
Endoscopic groups \(\Hbf\) as above are isogenous to products of split classical groups \(\Spbf_{2a}\) and \(\SObf_{4b}\).
For a classical group (such as \(\Spbf_{2a}\) and \(\SObf_{4b}\)) the Langlands dual group is also classical.
Conjugacy classes of parameters taking values in a classical group (in this case \(\SO_{2a+1}(\C)\) or \(\SO_{4a}(\C)\)) can be elementarily classified in terms of: the set of isomorphism classes of irreducible representations of \(L_\Q\), the duality map on this set, and for self-dual irreducible representations their type (symplectic or orthogonal) and determinant.
As recalled above isomorphism classes of irreducible representations of \(L_\Q\) are conjectured to be in bijection with cuspidal automorphic representations for general linear groups over \(\Q\), and the notions above (duality, determinant and symplectic/orthogonal type) admit translations on the automorphic side.
Using this observation Arthur formulated and proved a spectral expansion \cite[Corollary 3.4.2 and Theorem 4.1.2]{Arthur} which circumvents the hypothetical group \(L_\Q\), in terms of self-dual cuspidal automorphic representations for general linear groups.
Arthur's proof relies on the stabilization of the (twisted) trace formula.
Unfortunately it does not seem to be possible to proceed similarly for groups which are merely isogenous to classical groups.
Xu Bin \cite{Xu} obtained remarkable results towards a spectral expansion for groups such as \(\GSpbf_{2n}\), but a complete spectral expansion seems to be out of reach at the moment.

To understand the issue more concretely, but only conjecturally, consider a parameter \(\psi \in \Psi(\GSpbf_{2n}, V)\).
The condition at the real place implies that \(\Std \circ \psi\) decomposes as \(\bigoplus_{i=0}^r \psi_i\) for some \(r \geq 0\), where each \(\psi_i\) is irreducible, say of dimension \(m_i\), the \(\psi_i\)'s are non-isomorphic and exactly one of them has odd dimension, say \(\psi_0\).
For simplicity, and because this will be implied by the ``level one'' condition that we will eventually impose, assume that each character \(\det \psi_i\) is trivial.
Under this assumption we may find lifts \(\psi_i^{\GSpin}: L_\Q \times \SL_2(\C) \to \GSpin_{m_i}(\C)\), unique up to characters \(L_\Q \to Z(\GSpin_{m_i}(\C))\), and up to a character of \(L_\Q\) the representations \(Y_\nu\) occurring in the decomposition \eqref{eq:def_Y_psi_pi_f} decompose as
\begin{equation} \label{eq:intro_tens_prod_spin}
  (\spin \circ \psi_0^{\GSpin}) \otimes \bigotimes_{i=1}^r \spin^{\epsilon(\nu,i)} \circ \psi_i^{\GSpin}
\end{equation}
where for \(i>0\) the sign \(\epsilon(\nu,i)\) distinguishes one of the two half-spin representations \(\spin^{\pm}\) of \(\GSpin_{m_i}(\C)\).
Moreover (still conjecturally) for a given \(\pi_f\) the multiplicity \(m(\psi,\pi_f,\nu)\) vanishes except for exactly one value for \(\nu\), for which this multiplicity is one.
Each \(\psi_i\) should be the tensor product of a self-dual irreducible representation of \(L_\Q\) of dimension \(n_i\) and the irreducible algebraic representation of \(\SL_2(\C)\) of dimension \(d_i\).
The former should correspond to a self-dual cuspidal automorphic representation \(\pi_i\) for \(\GLbf_{n_i,\Q}\).
We use the notation \(\pi_i[d_i]\) for the pair \((\pi_i, d_i)\) to suggest this tensor product.
In the absence of the Langlands group \(L_\Q\) Arthur replaced the parameter \(\Std \circ \psi\) by the multiset \(\{ (\pi_i, d_i) | 0 \leq i \leq r \}\), that we simply denote by \(\pi_0[d_0] \oplus \dots \oplus \pi_r[d_r]\).
Among the difficulties in extending Arthur's stable multiplicity formula to groups such as \(\GSpbf_{2n}\), one has to find a way to (unconditionally) distinguish the various lifts \(\psi\) such that \(\Std \circ \psi\) corresponds to \(\bigoplus_i \pi_i[d_i]\), and for each such lift define a multiset \(\Pi_f(\psi)\) of irreducible representations of \(\GSpbf_{2n}(\A_f)\).
This has not been achieved in general.

\subsubsection{Adding the level one assumption}

Let us now restrict to the level one case: we only consider parameters \(\psi\) such that there exists \(\pi_f \in \Pi_f(\psi)\) which is everywhere unramified, i.e.\ the subspace of \(\GSpbf_{2n}(\Zhat)\) in \(\pi_f\) is non-zero.
This assumption should be equivalent to \(\psi\) factoring through the (conjectural) largest quotient \(L_\Z\) of \(L_\Q\) in which the groups \(I_{\Qp} \times \SU(2)\) become trivial (for all primes \(p\)).
Chenevier and Renard observed remarkable properties of this (conjectural) group \cite[Appendix B]{ChRe}, in particular it should decompose canonically as \(\R_{>0} \times L_\Z^1\) where \(L_\Z^1\) is a product of simply connected quasisimple compact Lie groups.
This gives us canonical choices for the lifts \(\psi_i^{\GSpin}\), by imposing that they factor through \(L_\Z\) and take values in \(\Spin_{m_i}(\C)\).
To \(\psi_i = \pi_i[d_i]\) is associated a family \((\tilde{c}_p(\psi_i))_p\) of semi-simple conjugacy classes in \(\GL_{m_i}(\C)\) (ranging over all primes \(p\)), defined by
\[ \tilde{c}_p(\psi_i) = c(\pi_{i,p}) \otimes \diag(p^{(d_i-1)/2}, p^{(d_i-3)/2}, \dots, p^{(1-d_i)/2}) \]
where \(c(\pi_{i,p})\) is the Satake parameter of \(\pi_{i,p}\).
Thus for our purpose, that is a particular case of the spectral expansion (``stable multiplicity formula'') for certain endoscopic groups of \(\GSpbf_{2n}\) in level one and for pseudo-coefficients of discrete series at the real place, we want to pin down the semisimple conjugacy classes
\[ \psi_i^{\GSpin}(\Frob_p, \diag(p^{1/2}, p^{-1/2})) \in \Spin_{m_i}(\C). \]
Since \(\psi_i^{\GSpin}\) is uniquely determined we will denote these by \(\cpsc(\psi_i)\).
These should satisfy \(\Std(\cpsc(\psi_i)) = \tilde{c}_p(\psi_i)\).
For \(i=0\) there is a unique semi-simple conjugacy class \(c_p(\psi_0)\) in \(\SO_{m_0}(\C)\) mapping to \(\tilde{c}_p(\psi_0)\), so this relation determines \(\cpsc(\psi_0)\) up to multiplication by \(Z(\Spin_{m_0}(\C)) \simeq \{ \pm 1 \}\).
For \(i>0\) the situation is more complicated, since in general there are two semisimple conjugacy classes in \(\SO_{m_i}(\C)\) mapping to \(\tilde{c}_p(\psi_i)\), exchanged by any element of \(\mathrm{O}_{m_i}(\C)\) having determinant \(-1\).
So in this case we first need to pin down the ``correct'' one \(c_p(\psi_i)\) and prove a spectral expansion for split groups \(\SObf_{4m}\) in level one and for pseudo-coefficients of discrete series at the real place (Proposition \ref{pro:refineSOeven}), slightly refining the specialization of Arthur's stable multiplicity formula in this setting (Arthur's results for even special orthogonal groups are all ``up to outer automorphism'').
We then pin down semi-simple conjugacy classes \(\cpsc(\psi_i)\) in spin groups and prove a spectral expansion for groups which are quotients of products of split groups \(\Spbf_{2a}\) and \(\SObf_{4b}\) under the same assumptions (Proposition \ref{pro:lift_stab_mult}).
These results rely on the stabilization of the trace formula and the (unconditional) automorphic counterpart to the property ``\(L^1_\Z\) is simply connected'' \cite[Proposition 4.4]{ChRe}.

Plugging this spectral expansion into Morel's formula \cite[Corollaire 5.3.3]{MorelSiegel2} and following Kottwitz' argument from \cite{Kottwitz_AA} in this unconditional setting, we obtain the following theorem.
If intersection cohomology does not vanish in level one then \(-1 \in \Zbf(\GSpbf_{2n})(\Qell)\) acts trivially on \(V\), so we may reduce to the case where \(V\) is a representation of \(\PGSpbf_{2n}\) by twisting (see Remark \ref{rema:twisting_sim}).

\begin{theointro}[weaker, vague version of Theorem \ref{theo:IH_explicit_crude}] \label{theointro:IH_explicit_crude}
  Let \(V\) be an irreducible algebraic representation of \(\PGSpbf_{2n}\).
  Then up to semi-simplification the representation
  \[ \Qellbar \otimes_{\Qell} \IH^\bullet(\GSpbf_{2n}, \Xcal, V)^{\GSpbf_{2n}(\Zhat)} \]
  of \(\GalQ \times \Hcal^{\unr}(\PGSpbf_{2n})\) is isomorphic to the sum of tensor products
  \[ \sigma_{\psi,\iota}^{\IH} \otimes_{\Qellbar} \iota(\chi_{f,\psi}) \]
  where
  \begin{itemize}
  \item \(\psi = \psi_0 \oplus \dots \oplus \psi_r\) ranges over Arthur's substitutes for global parameters for \(\Spbf_{2n}\) \cite[\S 1.4]{Arthur} which are unramified and of infinitesimal character determined by \(V\), in particular \(\psi_i = \pi_i[d_i]\) with \(\pi_i\) an everywhere unramified self-dual cuspidal automorphic representation for a general linear group,
  \item \(\chi_{f,\psi}\) is the character of \(\Hcal^{\unr}(\PGSpbf_{2n})\) corresponding to the images of \((\cpsc(\psi_i))_{0 \leq i \leq r}\) in \(\Spin_{2n+1}(\C)\) under the natural morphism
\[ \prod_{i=0}^r \Spin_{m_i}(\C) \longrightarrow \Spin_{2n+1}(\C), \]
  \item \(\sigma_{\psi,\iota}^{\IH}\) is a \(2^{n-r}\)-dimensional continuous semisimple representation of \(\GalQ\) over \(\Qellbar\) which is unramified away from \(\ell\) and such that for any prime \(p \neq \ell\) the semi-simplification of \(\sigma_{\psi,\iota}^{\IH}(\Frob_p)\) is conjugated to
    \begin{equation} \label{eq:intro_explicitIH_cc_Frob}
      \iota \left( p^{n(n+1)/4} \spin_{\psi_0}(\cpsc(\psi_0)) \otimes \spin_{\psi_1}^{u_1(\psi)}(\cpsc(\psi_1)) \otimes \dots \otimes \spin_{\psi_r}^{u_r(\psi)}(\cpsc(\psi_r)) \right)
    \end{equation}
    where \(u_i(\psi) \in \{ \pm 1 \}\) are explicit signs.
  \end{itemize}
\end{theointro}

Anticipating on the explicit relation between intersection and compactly supported cohomology (see \S \ref{sec:intro_Hc_IH}), we also deduce from a result of Faltings-Chai that each \(\sigma_{\psi,\iota}^{\IH}\) is crystalline at \(\ell\).

\subsubsection{\(\GSpin\)-valued Galois representations}

We expect from \ref{eq:intro_tens_prod_spin} that each \(\sigma_{\psi,\iota}^{\IH}\) decomposes as a tensor product of Galois representations obtained by composing \(\GSpin\)-valued Galois representations with (half-)spin representations, as suggested by \eqref{eq:intro_explicitIH_cc_Frob}, but this does not immediately follow from \eqref{eq:intro_explicitIH_cc_Frob}.
Proving that this holds (almost always) is the subject of Sections \ref{sec:odd_spin} and \ref{sec:even_spin}.
Let us mention the odd-dimensional case first (i.e.\ those parameters which can occur as \(\psi_0\) above).

\begin{theointro}[Theorem \ref{theo:odd_spin_rep}] \label{theointro:odd_spin_rep}
  Let \(\psi = \pi[d]\) be a self-dual Arthur-Langlands parameter of \emph{odd} orthogonal type, everywhere unramified and with algebraic regular infinitesimal character (i.e.\ its eigenvalues are distinct integers, see Definitions \ref{def:ICcal} and \ref{def:Psi}).
  There exists a continuous semisimple morphism \(\rho_{\psi,\iota}^{\GSpin}: \GalQ \to \GSpin_{2n+1}(\Qellbar)\) unramified away from \(\ell\) and crystalline at \(\ell\) and such that for any prime \(p \neq \ell\) the semi-simplification of \(\rho_{\psi,\iota}^{\GSpin}(\Frob_p)\) belongs to \(\iota(p^{n(n+1)/4} \cpsc(\psi))\).
  Moreover the conjugacy class of \(\rho_{\psi,\iota}^{\GSpin}\) admits the following characterizations.
  \begin{enumerate}
  \item If \(\rho: \GalQ \to \GSpin_{2n+1}(\Qellbar)\) is any continuous semisimple morphism such that for almost all primes \(p\), the semisimplifications of \(\rho|_{\GalQp}\) and \(\rho_{\psi,\iota}^{\GSpin}|_{\GalQp}\) are conjugate, then \(\rho\) is conjugate to \(\rho_{\psi,\iota}^{\GSpin}\).
  \item If \(\rho: \GalQ \to \GSpin_{2n+1}(\Qellbar)\) is a continuous morphism lifting the morphism \(\rho_{\psi,\iota}^{\SO}: \GalQ \to \SO_{2n+1}(\Qellbar)\) (Theorem \ref{thm:existence_rho_SO}), unramified away from \(\ell\) and crystalline at \(\ell\), then \(\rho\) is conjugate to \(\chi_\ell^N \rho_{\psi,\iota}^{\GSpin}\) for some integer \(N\).
  \end{enumerate}
\end{theointro}

We also prove in Proposition \ref{pro:rho_odd_GSpin_pi_2dp1} that in the non-tempered case, i.e.\ when \(\psi = \pi[d]\) with \(d>1\) (automatically odd), the morphism \(\rho_{\psi,\iota}^{\GSpin}\) can also be constructed from \(\rho_{\pi[1],\iota}^{\GSpin}\).
The basic idea to prove Theorem \ref{theointro:odd_spin_rep} is to combine the existence of \(\SO\)-valued Galois representations (recalled in Theorem \ref{thm:existence_rho_SO}), the existence of an essentially unique conductor one lift to \(\GSpin\), and to compare the composition of this lift with the spin representation with the representation \(\sigma_{\psi,\iota}^{\IH}\) constructed above.
This strategy is similar to the one used by Kret and Shin in \cite{KretShin_GSp}.

The even-dimensional case (in the preceding discussion, parameters which can occur as \(\psi_i\) for \(i>0\)) is more complicated because we do not know a priori the existence of a morphism \(\rho_{\psi,\iota}^{\SO}\) satisfying local-global compatibility at all primes \(p \neq \ell\) for the semisimple conjugacy classes \(c_p(\psi)\) in \(\SO_m(\C)\): a priori we only know local-global compatibility up to conjugation by \(\mathrm{O}_m(\Qellbar)\) (recalled in Theorem \ref{thm:existence_rho_O}).
We prove this refinement in almost all cases, and deduce the existence and uniqueness of \(\GSpin\)-valued Galois representations.

\begin{theointro}[Theorems \ref{thm:loc_glob_param_SO_even} and \ref{theo:even_GSpin_Gal_rep}] \label{theointro:even_spin_rep}
  Let \(\psi = \pi[d]\) be a self-dual Arthur-Langlands parameter of orthogonal type and \emph{even} dimension \(4n\), everywhere unramified and with algebraic regular infinitesimal character (i.e.\ its eigenvalues are distinct integers, see Definitions \ref{def:ICcal} and \ref{def:Psi}).
  Assume either \(n=1\), \(n\) even, \(d\) even, or that the infinitesimal character satisfies a regularity condition (see Definition \ref{def:bad_tau}).
  \begin{enumerate}
  \item There exists a continuous semisimple morphism \(\rho_{\psi,\iota}^{\SO}: \GalQ \to \SO_{4n}(\Qellbar)\) which is unramified away from \(\ell\) and crystalline at \(\ell\), such that for all primes \(p \neq \ell\) semi-simplification of \(\rho_{\psi,\iota}^{\SO}(\Frob_p)\) belongs to \(\iota(c_p(\psi))\).
    Any continuous semisimple morphism \(\rho: \GalQ \to \SO_{4n}(\Qellbar)\) satisfying this condition at almost all primes \(p\) is \(\SO_{4n}(\Qellbar)\)-conjugated to \(\rho_{\psi,\iota}^{\SO}\).
  \item There exists a continuous semisimple morphism \(\rho_{\psi,\iota}^{\GSpin}: \GalQ \to \GSpin_{4n}(\Qellbar)\) unramified away from \(\ell\) and crystalline at \(\ell\) and such that for any prime \(p \neq \ell\) the semi-simplification of \(\rho_{\psi,\iota}^{\GSpin}(\Frob_p)\) belongs to \(\iota(p^{n/2} \cpsc(\psi))\).
    It admits the same characterizations as in Theorem \ref{theointro:odd_spin_rep}.
  \end{enumerate}
\end{theointro}

For \(d>1\) even (resp.\ odd) we also prove in Proposition \ref{pro:Gal_GSpin_pi_2d} (resp.\ Proposition \ref{pro:rho_even_GSpin_pi_2dp1}) that \(\rho_{\psi,\iota}^{\GSpin}\) can also be constructed (in the \(d\) even case, up to outer automorphism) from the \(\GSp\)-valued (resp.\ \(\GSpin\)-valued) Galois representation associated to \(\pi\).
The proof of Theorem \ref{theointro:even_spin_rep} is rather indirect and relies heavily on \(\ell\)-adic families of automorphic representations for inner forms \(\Hbf\) of \(\PGSObf_{8m}\) which are split at all primes and definite (i.e.\ the Lie group \(\Hbf(\R)\) is compact).
Considering the contribution of the parameter \(1 \oplus \psi\) in Theorem \ref{theointro:IH_explicit_crude} we obtain the existence of a Galois representation \(\sigma_{\psi,\iota}^{\spin,\epsilon}\) which ought to be \(\spin^\epsilon \circ \rho_{\psi,\iota}^{\GSpin}\) (but we do not know the existence of \(\rho_{\psi,\iota}^{\GSpin}\) yet), for a sign \(\epsilon\) which we do not control (see Corollary \ref{coro:seed_ex_sigma} for details).
Assume first that \(n\) is even in Theorem \ref{theointro:even_spin_rep}, so that a definite inner form \(\Hbf\) of \(\PGSObf_{4n}\) split at all primes exists.
One can associate a level one automorphic representation \(\Pi\) for \(\Hbf\) to the parameter \(\psi\) (Example \ref{exam:mult_formula_definite}), and by \(\ell\)-adic interpolation, which is possible thanks to the fact that \(\Hbf\) is definite, we construct the ``other half-spin'' Galois representation \(\sigma_{\psi,\iota}^{\spin,-\epsilon}\).
By the main technical result of \cite{TaiCC} we can even \(\ell\)-adically interpolate \(\Pi\) by level one automorphic representations \(\Pi'\) (with associated parameter \(\psi' = \pi'[1]\)) such that the associated Galois representation has infinitesimally big image (the Lie algebra over \(\Qellbar\) generated by \(\Lie \Std(\rho_{\psi,\iota}^{\SO}(\GalQ))\) is maximal, i.e.\ equal to \(\mathfrak{so}_{4n}\)), in particular both representations \(\sigma_{\psi',\iota}^{\spin,\pm}\) are irreducible.
This allows us to deduce the first part of Theorem \ref{theointro:even_spin_rep} for \(\psi'\) and then for \(\psi\) (Corollary \ref{cor:loc_glob_param_SO8n}), and the second part follows as in the odd-dimensional case.
Part of this argument is similar to the one used by Kret and Shin in \cite{KretShin_GSO}, but they did not use \(\ell\)-adic interpolation and the ``infinitesimally big image'' condition is replaced by a weaker condition derived from their ``Steinberg at one place'' hypothesis.

The case where \(n\) is odd in Theorem \ref{theointro:even_spin_rep} is trickier, except for \(n=1\) because of the exceptional isomorphism \(\PGSObf_4 \simeq \PGLbf_2^2\).
For \(n>1\) odd we would like to apply the previous strategy to a level one automorphic representation \(\Pi\), for the inner form of \(\PGSObf_{4n+4}\) split at all places, corresponding to the parameter \(\psi \oplus \psi'\) where \(\psi'\) is \(4\)-dimensional.
Unfortunately we cannot find \(\psi'\) such that \(\Pi\) exists in all cases, whence the regularity condition on the infinitesimal character in Theorem \ref{theointro:even_spin_rep}.
See Section \ref{sec:using_endo} for details.

The proof of Theorem \ref{thmintro:Gal_rep_SMF} (Corollary \ref{cor:GSpin_Gal_rep_Siegel}) from Theorems \ref{theointro:odd_spin_rep} and \ref{theointro:even_spin_rep} is simply the observation that for the parameter \(\psi = \psi_0 \oplus \dots \oplus \psi_r\) corresponding to a Siegel cusp form, each \(\psi_i\) for \(i>0\) satisfies the assumption of Theorem \ref{theointro:even_spin_rep}.

\subsubsection{Tensor product decompositions in intersection cohomology}
\label{sec:intro_tensor_prod}

Going back to Theorem \ref{theointro:IH_explicit_crude} consider a parameter \(\psi = \psi_0 \oplus \dots \oplus \psi_r\).
We now have a \(\GSpin\)-valued Galois representation \(\rho_{\psi_0,\iota}^{\GSpin}\) by Theorem \ref{theointro:odd_spin_rep}, yielding a Galois representation \(\sigma_{\psi,\iota}^{\spin} := \spin \circ \rho_{\psi,\iota}^{\GSpin}\) (which is actually equal to the representation \(\sigma_{\psi,\iota}^{\IH}\) found in intersection cohomology).
For the other constituents of \(\psi\), namely \(\psi_i\) for \(i>0\), we only have \(\GSpin\)-valued Galois representations \(\rho_{\psi_i,\iota}^{\GSpin}\), and thus Galois representations \(\sigma_{\psi_i,\iota}^{\spin,\epsilon} := \spin^{\epsilon} \circ \rho_{\psi_i,\iota}^{\GSpin}\) for both values of the sign \(\epsilon\), under the assumption in Theorem \ref{theointro:even_spin_rep}.
Otherwise we only have \(\sigma_{\psi_i,\iota}^{\spin,-} := \sigma_{1 \oplus \psi, \iota}^{\IH}\).
It turns out that in the latter case the sign \(u_i(\psi)\) appearing in Theorem \ref{theointro:IH_explicit_crude} is always \(-1\), so this complication does not prevent us from proving in Theorem \ref{thm:sigma_is_tensor} that up to a Tate twist the representation \(\sigma_{\psi,\iota}^{\IH}\) is isomorphic to
\[ \sigma_{\psi_0,\iota}^{\spin} \otimes \sigma_{\psi_1,\iota}^{\spin,u_1(\psi)} \otimes \dots \otimes \sigma_{\psi_r,\iota}^{\spin,u_r(\psi)}. \]
This concludes the proof of our unconditional version of Kottwitz' conjecture (Conjecture \ref{conj:Kottwitz}) for level one Siegel modular varieties.

\subsection{Compactly supported cohomology}
\label{sec:intro_Hc_IH}

In Section \ref{sec:formula_e_c} we finally come back to the Euler characteristics \(e_c(\Acal_{n,\Qbar}, \Fcal(V))\) \eqref{eq:intro_e_c_An_FV}.
Now that we have a precise description of intersection cohomology, it is natural to try to express compactly supported cohomology in terms of intersection cohomology.
For this purpose it turns out that there is no significant simplification to be gained from restricting to level one, so we work in arbitrary level.
Denote by \((\GSpbf_{2n}, \Xcal_n)\) the Shimura datum corresponding to the Shimura tower \((\Acal_{n,K})_K\).
We prove in Theorem \ref{theo:Hc_from_IH} an explicit formula expressing, in a suitable Grothendieck group of representations of \(\GSpbf_{2n}(\A_f) \times \GalQ\), the Euler characteristic \(e_c(\GSpbf_{2n}, \Xcal_n, V)\) \eqref{eq:intro_e_c_G_X_V} in terms of (via parabolic induction)
\begin{enumerate}
\item the Euler characteristics \(e_\IH(\GSpbf_{2n'}, \Xcal_{n'}, V')\) for \(n' \leq n\) and certain explicit irreducible representations \(V'\) of \(\GSpbf_{2n',\Qell}\),
\item certain virtual representations \(e(\GLbf_1, a)\) (for \(a \in \Z\)) of \(\A_f^\times\) and \(e_{(2)}(\GLbf_2, a, b)\) (for \(a,b \in \C\) satisfying \(a-b \in \Z\)) of \(\GLbf_2(\A_f)\), respectively related to algebraic Hecke characters \(\Q^\times \backslash \A^\times \to \C^\times\) and (elliptic) modular forms (i.e.\ the representations \(\varinjlim_\Gamma S_k(\SL_2(\Z))\) of \(\GLbf_2(\A_f)\)), see Example \ref{exam:Euler_L2_coh_GL}.
\end{enumerate}

We now explain the steps that we take to arrive at Theorem \ref{theo:Hc_from_IH}.
For a connected reductive group \(\Gbf\) over \(\Q\) with maximal split central torus \(\Abf_\Gbf\) and a choice of maximal compact subgroup \(K_\infty\) of \(\Gbf(\R)\), to a finite-dimensional algebraic representation \(V\) of \(\Gbf\) are associated local systems (in \(\Q\)-vector spaces) \(\Fcal(V)\) on the locally symmetric spaces \(\Gbf(\Q) \backslash (\Gbf(\R)/K_\infty \Abf_\Gbf(\R)^0 \times \Gbf(\A_f)/K)\), where \(K\) ranges over neat compact open subgroups of \(\Gbf(\A_f)\).
We denote
\[ H^i(\Gbf, V) := \varinjlim_K H^i(\Gbf(\Q) \backslash (\Gbf(\R) / K_\infty \Abf_\Gbf(\R)^0 \times \Gbf(\A_f)/K), \Fcal(V)), \]
an admissible representation of \(\Gbf(\A_f)\) over \(\Q\), and
\[ e(\Gbf, V) := \sum_i (-1)^i [H^i(\Gbf, V)], \]
(in the Grothendieck group of admissible representations of \(\Gbf(\A_f)\) over \(\Q\)), and similarly for \(H^i_c(\Gbf, V)\) and \(e_c(\Gbf, V)\).
The cohomology groups \(H^i(\Gbf, V)\) may be identified with direct sums of group cohomology groups for certain arithmetic subgroups of \(\Gbf(\Q)\), and in many cases the compactly supported cohomology groups \(H^i_c(\Gbf, V)\) can be expressed by duality using ordinary cohomology groups \(H^i(\Gbf, V')\) (see Section \ref{sec:arith_gp_coh}).
Our starting point is a formula in the other direction derived from Morel's work \cite{MorelSiegel1} \cite{MorelBook} (see Corollary \ref{coro:IH_vs_Hc_Q}), which essentially expresses \(e_\IH(\GSpbf_{2n}, \Xcal_n, V)\) in terms (again, using parabolic induction) of
\begin{enumerate}
\item \(e_c(\GSpbf_{2n'}, \Xcal_{n'}, V')\) for \(n' \leq n\) and certain explicit irreducible representations \(V'\) of \(\GSpbf_{2n',\Qell}\),
\item \(e_c(\GLbf_{n'}, V')\) for \(n' \leq n\) and certain explicit irreducible representations \(V'\) of \(\GLbf_{n',\Q}\).
\end{enumerate}
The appearance of parabolic induction is implicit in Morel's work and we take this opportunity to make it explicit in Section \ref{sec:IH_vs_Hc_Morel}.
For this purpose we found convenient to introduce a new formulation for boundary cohomology, using a slight generalization of the notion of Shimura datum and Shimura varieties (Appendix \ref{sec:gen_Shim}).
This formulation incorporates arithmetic group cohomology and \(\ell\)-adic étale cohomology.

Our second step (Corollary \ref{coro:Euler_c_GL}) is to express the Euler characteristics \(e_c(\GLbf_{n'}, V')\) in terms of the (simpler) virtual representations \(e(\GLbf_1, a)\) and \(e_{(2)}(\GLbf_2, a, b)\).
In Section \ref{sec:Franke_formula} we recall Franke's filtration of the space of automorphic forms of a connected reductive group \(\Gbf\) over \(\Q\) \cite{Franke} and deduce (partially following Franke) in Corollary \ref{cor:Franke_Euler2} a formula expressing, in a suitable Grothendieck group of admissible representations of \(\Gbf(\A_f)\), the Euler characteristic \(e(\Gbf, V)\) in terms of Euler characteristics
\begin{equation} \label{eq:intro_e2}
  e_{(2)}(\Lbf, W) := \sum_i (-1)^i [\varinjlim_K H^i((\lfrak / \afrak_\Lbf, K_{\infty,L}), \Acal^2(\Lbf, \xi) \otimes W)]
\end{equation}
for \(\R\)-cuspidal Levi subgroups \(\Lbf\) of \(\Gbf\), where
\begin{itemize}
\item \(\lfrak = \C \otimes_\R \Lie \Lbf(\R)\) and \(\afrak_\Lbf = \C \otimes_\R \Lie \Abf_\Lbf(\R)\),
\item \(K_{\infty,\Lbf}\) is a maximal compact subgroup of \(\Lbf(\R)\),
\item \(W\) is a finite-dimensional representation of \(\Lbf(\R)\) with central character \(\xi^{-1}\) on \(\Abf_\Lbf(\R)^0\),
\item \(\Acal^2(\Lbf, \xi)\) is the space of automorphic forms \(f\) for \(\Lbf\) satisfying \(f(z \cdot) = \xi(z) f(\cdot)\) for \(z \in \Abf_\Lbf(\R)^0\) and whose restriction to \(\Lbf(\Q) \backslash \Lbf(\A)^1\) is square-integrable, where \(\Lbf(\A) = \Lbf(\A)^1 \times \Abf_\Lbf(\R)^0\) is the usual decomposition.
\end{itemize}
For \(\Gbf = \GLbf_{n'}\) any such Levi subgroup \(\Lbf\) is isomorphic to a product of \(\GLbf_1\)'s and \(\GLbf_2\)'s, and we deduce the following formula.

\begin{theointro}[Corollary \ref{coro:Euler_ord_GL}] \label{theointro:Euler_ord_GL}
  Let \(n \geq 1\).
  For \(a,b \in \Z_{\geq 0}\) such that \(a+2b=n\) denote by \(\Lbf_{a,b} \simeq \GLbf_1^a \times \GLbf_2^b\) the corresponding standard Levi subgroup of \(\GLbf_n\), and let \(\mathfrak{S}(a,b)\) be the subset of \(\mathfrak{S}_n\) consisting of \(\sigma\) such that
  \begin{enumerate}
  \item \(\sigma^{-1}(1) < \dots < \sigma^{-1}(a)\),
  \item \(\sigma^{-1}(a+1) < \sigma^{-1}(a+2)\), \dots, \(\sigma^{-1}(a+2b-1) < \sigma^{-1}(a+2b)\),
  \item \(\sigma^{-1}(a+1) < \sigma^{-1}(a+3) < \dots < \sigma^{-1}(a+2b-1)\).
  \end{enumerate}
  (In other words \(\Sfrak(a,b)\) parametrizes partitions of \(\{1,\dots,n\}\) into \(a\) singletons and \(b\) unordered pairs.)
  Consider a dominant weight \(\lambda = (\lambda_1 \geq \dots \geq \lambda_n)\) for \(\GLbf_n\) and let
  \[ \tau = (\tau_1 > \dots > \tau_n) := \lambda + \rho \]
  so that \(\tau_i = \lambda_i + \frac{n+1}{2} - i\).
  For \(\sigma \in \mathfrak{S}_n\) denote \(\sigma(\tau)_i = \tau_{\sigma^{-1}(i)}\) and \((\sigma \cdot \lambda)_i = \lambda_{\sigma^{-1}(i)} - \sigma^{-1}(i) + i\), i.e.\ \(\sigma \cdot \lambda = \sigma(\tau) - \rho\).
  Using notation introduced in Example \ref{exam:Euler_L2_coh_GL} (for \(e(\GLbf_1, -)\) and \(e_{(2)}(\GLbf_2, -, -)\)) we have, in the weak Grothendieck group \(K_0^{\Tr}(\Rep_{\C}^{\adm}(\GLbf_n(\A_f)))\) of admissible complex representations of \(\GLbf_n(\A_f)\) (see Definition \ref{def:Groth_adm}),
  \begin{multline*}
    e(\GLbf_n, V_{\lambda}) = \sum_{a+2b=n} (-1)^{a(a-1)/2} \sum_{\sigma \in \mathfrak{S}(a,b)} \epsilon(\sigma) \\
    \Ind_{\Lbf_{a,b}(\A_f)}^{\GLbf_n(\A_f)} \Big( \bigotimes_{i=1}^a e(\GLbf_1, (\sigma \cdot \lambda)_i + a + 1) |\cdot|_f^{(\sigma \cdot \lambda)_i + a + 1 - \sigma(\tau)_i} \\
    \otimes \bigotimes_{i=1}^b e_{(2)}(\GLbf_2, \sigma(\tau)_{a+2i-1}-1/2, \sigma(\tau)_{a+2i}+1/2) \Big)
  \end{multline*}
  where \(\Ind\) denotes normalized\footnote{This normalization is convenient because the image in \(K_0^{\Tr}(\Rep_{\C}^{\adm}(\GLbf_n(\A_f)))\) of such parabolically induced representations do not depend on the choice of a parabolic subgroup (with given Levi subgroup), but it introduces square roots of integers and is the reason for working with complex (rather than rational) coefficients here.
  It is easy to check a posteriori that these parabolically induced representations are naturally defined over \(\Q\).} parabolic induction.
\end{theointro}

Using Poincaré duality we easily deduce a similar formula for \(e_c(\GLbf_{n'}, V')\) in Corollary \ref{coro:Euler_c_GL}.
Plugging this formula into the previous formula (Corollary \ref{coro:IH_vs_Hc_Q}) and simplifying the resulting expression yields a formula for \(e_\IH(\GSpbf_{2n}, \Xcal_n, V)\) in terms of \(e_c(\GSpbf_{2n'}, \Xcal_{n'}, V')\), \(e(\GLbf_1, a)\) and \(e_{(2)}(\GLbf_2, a, b)\) (Theorem \ref{theo:IH_vs_Hc_simple}).

The third and final step consists of inverting this relation in Section \ref{sec:Hc_from_IH}, to obtain a (again, simplified) formula in the other direction.

\begin{theo}[Theorem \ref{theo:Hc_from_IH}] \label{theointro:Hc_from_IH}
  For an integer \(n \geq 1\), a dominant weight \(\lambda\) for \(\GSpbf_{2n}\) and a prime number \(\ell\) the Euler characteristic \(e_c(\GSpbf_{2n}, \Xcal_n, V_\lambda)\), in the weak Grothendieck group \(K_0^{\Tr}(\Rep_{\Qell}^{\adm,\cont}(\GSpbf_{2n}(\A_f) \times \GalQ))\) of admissible continuous representations of \(\GSpbf_{2n}(\A_f) \times \GalQ\) (see Definition \ref{def:Groth_adm}), is equal to
  \begin{multline*}
    \sum_{\substack{a,b \geq 0 \\ a+2b \leq n}} \sum_{w \in W(a,b,n)} (-1)^{a+b} \epsilon(w)  \\
    \ind_{\Pbf_{a,b,n}(\A_f)}^{\GSpbf_{2n}(\A_f)} \left( e_{(2)} \left( \GLbf_1^a \times \GLbf_2^b, (w \cdot \lambda)_{\lin} \right) \otimes e_\IH(\GSpbf_{2(n-a-2b)}, V_{(w \cdot \lambda)_{\her}}) \right).
  \end{multline*}
  where
  \begin{itemize}
  \item the subset \(W(a,b,n)\) of the Weyl group of \(\GSpbf_{2n}\) is defined at the beginning of Section \ref{sec:Hc_from_IH},
  \item \(\Pbf_{a,b,n}\) is the standard parabolic subgroup of \(\GSpbf_{2n}\) with Levi subgroup \(\GLbf_1^a \times \GLbf_2^b \times \GSpbf_{2(n-a-2b)}\), and the subscript \(\lin\) (resp.\ \(\her\)) corresponds to projecting to the \(\GLbf_1^a \times \GLbf_2^b\) (resp.\ \(\GSpbf_{2(n-a-2b)}\)) factor,
  \item \(\ind\) denotes (unnormalized) smooth parabolic induction.
  \end{itemize}
\end{theo}

Together with the special case of Kottwitz' conjecture proved earlier (Theorem \ref{theointro:IH_explicit_crude} and the tensor product decompositions explained in \S \ref{sec:intro_tensor_prod}) this gives (in principle) an explicit formula for the Euler characteristics \(e_c(\Acal_{n,\Qbar}, \Fcal(V))\).
Forgetting the Hecke action and translating using the exceptional isomorphism \(\PGSObf_4 \simeq \PGLbf_2^2\), we specialize this formula to the case of genus \(n \leq 3\) in Section \ref{sec:BFG} and finally prove Theorem \ref{theointro:BFG}.
In genus \(n>3\) some parameters \(\psi = \psi_0 \oplus \dots \oplus \psi_r\) in Theorem \ref{theointro:IH_explicit_crude} include even-dimensional components \(\psi_i\) of dimension \(\geq 8\), and their contributions to intersection cohomology cannot be expressed as a contribution to the representation \(S[\ul{k}]_\ell\) of \S \ref{sec:intro_BFG}.
In fact as \(n\) grows, and thus as \(r\) potentially grows, the difference between these two virtual representations ``increases'': not all parameters \(\psi\) (with suitable infinitesimal character) contribute to \(S[\ul{k}]_\ell\) (conditions recalled in Theorem \ref{thm:Siegel_formula}), and those that do contribute contribute, up to a Tate twist
\[ \bigoplus_{(\epsilon_i)_i \in \{ \pm 1 \}^r} \sigma_{\psi_0,\iota}^{\spin} \otimes \sigma_{\psi_1,\iota}^{\spin,\epsilon_1} \otimes \dots \otimes \sigma_{\psi_r,\iota}^{\spin,\epsilon_r}, \]
whereas all parameters \(\psi\) contribute to intersection cohomology, but each contributes only one tensor product (for \(\epsilon_i = u_i(\psi)\)).

\subsection{Acknowledgments}

I thank Sophie Morel for discussions related to her work, and Gaëtan Chenevier for his comments.
During this work I was supported by the ANR-19-CE40-0015-02 project COLOSS.

\subsection{Disclaimer}

The results of this paper rely on Arthur's endoscopic classification \cite{Arthur}, which is not yet completely unconditional.
Thanks to \cite{AGIKMS} the only remaining result to be proved is the (standard and non-standard) weighted fundamental lemma for Lie algebras over positive characteristic local fields (generalizing \cite{ChaLauII}).
On a positive note, we understand that Connor Halleck-Dub\'e is making good progress on this generalization, so one can be hopeful that the results of \cite{Arthur} will soon become unconditional.

\section{Notations and conventions}

\subsection{Class field theory}

For a prime number \(p\) we denote by \(\Frob_p\) the geometric Frobenius element of \(\Gal(\Fpbar/\Fp)\).
We normalize the reciprocity law in local class field theory by letting the geometric Frobenius element correspond to a uniformizer (e.g.\ \(\Frob_p\) corresponds to \(p\)).
For a prime number \(\ell\) we denote by \(\chi_\ell\) the \(\ell\)-adic cyclotomic character of \(\GalQ\), so that for any prime \(p \neq \ell\) we have \(\chi_\ell(\Frob_p) = p^{-1}\).

\subsection{Reductive groups and root data}
\label{sec:not_red_gps}

We use bold letters to denote reductive groups over global and local fields, and normal letters for their Langlands dual groups, that we will consider defined over \(\Qbar\).

Let \(n \geq 0\) be an integer and \(\Lambda\) a free \(\Z\)-module of rank \(2n\) endowed with a non-degenerate (over \(\Z\)) alternate bilinear form \(\langle \cdot, \cdot \rangle\).
The first reductive group of interest in this paper is the associated general symplectic group \(\Gbf\) (this notation is only temporary), that is the reductive group over \(\Z\) (in the sense of \cite[Exposé XIX, Définition 2.7]{SGA3-III}) defined by \(\Gbf(\Spec A ) = \{ (g, s) \in \GL(A \otimes_{\Z} V) \times A^{\times} | \langle g \cdot, g \cdot \rangle = s \langle \cdot, \cdot \rangle \}\).
Denote \(\nu : \Gbf \rightarrow \GLbf_1\), \((g,s) \mapsto s\) the similitude character.
Let \(J_n\) be the \(n\) by \(n\) ``antidiagonal'' matrix defined by \((J_n)_{i,j} = \delta_{i,n+1-j}\) (Kronecker delta).
There exists a basis of \(\Lambda\) in which the matrix of \(\langle \cdot, \cdot \rangle\) is \(\begin{pmatrix} 0 & J_n \\ -J_n & 0 \end{pmatrix}\).
This identifies \(\Gbf\) with the subgroup \(\GSpbf_{2n}\) of the matrix group \(\GLbf_{2n} \times \GLbf_1\), the restriction of the projection on the second factor being \(\nu\).
Let \(\Spbf_{2n}\) be the kernel of \(\nu: \GSpbf_{2n} \to \GLbf_1\).
Let \(\Tbf_{\GSpbf_{2n}}\) be the diagonal split maximal torus in \(\GSpbf_{2n}\), i.e.\ the subgroup consisting of \(t = (\diag(s t_1, \dots, s t_n, t_n^{-1}, \dots, t_1^{-1}), s)\) where \(s\) and the \(t_i\)'s belong to \(\GLbf_1\).
We will sometimes simply denote such a \(t\) by \((t_1, \dots, t_n, s)\), and similarly denote characters and cocharacters for \(\Tbf_{\GSpbf_{2n}}\) by tuples of integers.
The set of simple roots associated to the upper triangular Borel subgroup of \(\GSpbf_{2n}\) is \(\{\alpha_1, \dots, \alpha_n\}\) where \(\alpha_i(t) = t_i/t_{i+1}\) for \(i<n\) and \(\alpha_n(t) = s t_n^2\).
In particular \(X^*(\Tbf_{\GSpbf_{2n}})\) has a basis consisting of the simple roots \(\alpha_1, \dots, \alpha_n\) and \((\alpha_n + \nu)/2\) if \(n>0\).
If \(n=0\) then \(\nu\) is a basis of \(X^*(\Tbf_{\GSpbf_{2n}})\).
Let \(\Tbf_{\Spbf_{2n}} = \Spbf_{2n} \cap \Tbf_{\GSpbf_{2n}}\), a split maximal torus of \(\Spbf_{2n}\) whose group of characters admits as basis \(\alpha_1, \dots, \alpha_{n-1}, \alpha_n/2\) if \(n>0\).
We will sometimes need to consider symplectic and general symplectic group on the dual side, so we denote \(\Sp_{2n} = \Spbf_{2n,\Qbar}\), \(\Tcal_{\Sp_{2n}} = \Tbf_{\Spbf_{2n}, \Qbar}\), etc.

The Langlands dual group \(\widehat{\GSpbf_{2n}}\) is known to be isomorphic to \(\GSpin_{2n+1}\).
We will need an explicit identification.
Let \(\SO_{2n+1}\) be the special orthogonal group over \(\Qbar\) defined by
\(\SO_{2n+1}(R) = \{ M \in \SL_{2n+1}(R) \,|\, {}^t M J_{2n+1} M = J_{2n+1} \}\).
Just like for the symplectic group, this matrix realization is chosen so that
the diagonal and upper triangular subgroups of \(\SO_{2n+1}\) form a Borel pair,
which we denote \((\Tcal_{\SO_{2n+1}}, \Bcal_{\SO_{2n+1}})\).
We write \((z_1, \dots, z_n)\) for the element \(\diag(z_1, \dots, z_n, 1,
z_n^{-1}, \dots, z_1^{-1})\) of \(\Tcal_{\SO_{2n+1}}\).
We realise \(\SO_{2n+1}\) as the Langlands dual of \(\Spbf_{2n}\) by identifying
the simple root \(\alpha_i\) with a simple coroot \(\widehat{\alpha_i}\) for
\((\SO_{2n+1}, \Tcal_{\SO_{2n+1}}, \Bcal_{\SO_{2n+1}})\) as follows.
If \(i<n\) then \(\widehat{\alpha_i}(x) = (\dots, 1, x, x^{-1}, 1, \dots)\) where
\(x\) is the \(i\)-th coefficient and all but two coefficients are \(1\), and
\(\widehat{\alpha_n}(x) = (1, \dots, 1, x^2)\).
Let \(\Spin_{2n+1}\) be the simply connected cover of \(\SO_{2n+1}\), and let
\(\Tcal_{\Spin_{2n+1}}\) be the preimage of \(\Tcal_{\SO_{2n+1}}\).
The cocharacter group \(X_*(\Tcal_{\Spin_{2n+1}})\) equals \(\oplus_{i=1}^n \Z
\widehat{\alpha_i}\) and so an element of \(\Tcal_{\Spin_{2n+1}}\) can be
parametrized by an element \((z_1, \dots, z_n)\) of \(\Tcal_{\SO_{2n+1}}\) along
with \(s \in \GL_1\) satisfying \(s^2 = z_1 \dots z_n\).
We will simply denote such an element of \(\Tcal_{\Spin_{2n+1}}\) by \((z_1, \dots,
z_n, s)\).
For \(n>0\) define \(\GSpin_{2n+1}\) as the quotient of \(\Spin_{2n+1} \times \GL_1\)
by the diagonal subgroup \(\mu_2\), and define \(\GSpin_1 = \GL_1\).
Let \(\Tcal_{\GSpin_{2n+1}}\) be the image of \(\Tcal_{\Spin_{2n+1}} \times \GL_1\),
so that for \(n>0\) its group of points over any algebraically closed extension of
\(\Qbar\) is parametrized as
\begin{equation} \label{eq:param_TGSpin_odd}
  \left\{ (z_1, \dots, z_n, s, \lambda) \in \GL_1^{n+2} \,|\, s^2=z_1 \dots z_n
  \right\} / \langle (1, \dots, 1, -1, -1 ) \rangle.
\end{equation}

The identification of \(\GSpin_{2n+1}\) with the Langlands dual of \(\GSpbf_{2n}\)
is given by the above identification \(\SO_{2n+1} \simeq \widehat{\Spbf_{2n}}\)
along with the identification of \(\nu\) with the cocharacter \(\nuhat\) of
\(\Tcal_{\GSpin_{2n+1}}\) defined by \(\nuhat(x) = (1, \dots, 1, 1, x)\) (in
\(\Tcal_{\Spin_{2n+1}} \times \GL_1\)).

For \(n \geq 0\) let \(\Sp_{2n} = \Spbf_{2n} \times_{\Z} \Qbar\) and
\(\Tcal_{\Sp_{2n}} = \Tbf_{\Spbf_{2n}} \times_{\Z} \Qbar\), allowing us to use the
parametrization above.

Finally for \(n \geq 1\) let \(\Obf_{2n}\) be the (schematic) stabilizer in
\(\GLbf_{2n}\) of the quadratic form
\begin{align*}
  q: \Z^{2n} & \longrightarrow \Z \\
  (x_1, \dots, x_{2n}) & \longmapsto \sum_{i=1}^n x_i x_{2n+1-i}
\end{align*}
whose associated bilinear form \((x,y) \mapsto q(x+y)-q(x)-q(y)\) has Gram matrix
\(J_{2n}\).
Let \(\SObf_{2n}\) be the kernel of the Dickson morphism from \(\Obf_{2n}\) to the
constant group scheme \(\Z/2\Z\) over \(\Z\).
Let \(\Tbf_{\SObf_{2n}}\) be the diagonal split maximal torus in \(\SObf_{2n}\),
identified to \(\GLbf_1^n\) via
\begin{equation} \label{eq:param_T_SO_even}
  (t_1, \dots, t_n) \mapsto \diag(t_1, \dots, t_n, t_n^{-1}, \dots, t_1^{-1}).
\end{equation}
The set of simple roots associated to the upper triangular Borel subgroup of
\(\SObf_{2n}\) is \(\{ \alpha_1, \dots, \alpha_n \}\) where \(\alpha_i(t) =
t_i/t_{i+1}\) for \(i<n\) and \(\alpha_n(t) = t_{n-1} t_n\).
The Langlands dual group \(\widehat{\SObf_{2n}}\) is isomorphic to \(\SO_{2n} :=
\SObf_{2n} \times_{\Z} \Qbar\), and again we shall need an explicit
identification.
Letting \((\Tcal_{\SO_{2n}}, \Bcal_{\SO_{2n}})\) be the diagonal and upper
triangular Borel pair in \(\SO_{2n}\), it is easy to check that there is an
isomorphism between \(\widehat{\SObf_{2n}}\) and \(\SO_{2n}\) under which \(\alpha_i
\in X^*(\Tbf_{\SObf_{2n}})\) corresponds to the coroot \(\alpha_i^\vee \in
X_*(\Tcal_{\SO_{2n}})\).
Let \(\PGSObf_{2n}\) be the adjoint quotient of \(\SObf_{2n}\).
The Langlands dual \(\widehat{\PGSObf_{2n}}\) is thus identified to the simply
connected cover \(\Spin_{2n}\) of \(\SO_{2n}\).
The preimage \(\Tcal_{\Spin_{2n}}\) of \(\Tcal_{\SO_{2n}}\) in \(\Spin_{2n}\) is
parametrized as
\begin{equation} \label{eq:param_T_Spin_even}
  \left\{ (z_1, \dots, z_n, s) \in \GL_1^{n+1} \, \middle| \, s^2 = z_1 \dots z_n \right\}.
\end{equation}
As in the odd-dimensional case we define \(\GSpin_{2n}\) as the quotient of \(\Spin_{2n} \times \GL_1\) by the diagonally embedded \(\mu_2\) (on the first factor, the kernel of the projection \(\Spin_{2n} \to \SO_{2n}\)).
We have a morphism \(\GSpin_{2n} \to \SO_{2n}\) which is trivial on the second factor of \(\Spin_{2n} \times \GL_1\), and the points over any algebraically closed extension of \(\Qbar\) of the preimage \(\Tcal_{\GSpin_{2n}}\) of \(\Tcal_{\SO_{2n}}\) may be parametrized as in the odd case \eqref{eq:param_TGSpin_odd} as
\begin{equation} \label{eq:param_TGSpin_even}
  \left\{ (z_1, \dots, z_n, s, \lambda) \in \GL_1^{n+2} \,|\, s^2=z_1 \dots z_n  \right\} / \langle (1, \dots, 1, -1, -1 ) \rangle.
\end{equation}

For a group \(\Gbf\) as above we denote by \(\Std_{\Ghat}: \Ghat \to \GL_{N(\Ghat)}\) the standard representation of the dual group \(\Ghat\).

\section{The spectral side}

\subsection{Formal Arthur-Langlands parameters}

\begin{defi} \label{def:Calg}
  Let \(H\) be a complex connected reductive group, and let \(\hfrak\) be its Lie algebra.
  Let \(T\) be a maximal torus of \(H\), \(\tfrak\) its Lie algebra.
  We call a semi-simple \(H(\C)\)-orbit in \(\hfrak\) C-algebraic if it is represented by an element of \(\tfrak\) belonging to \(\rho + X_*(T)\) where \(2 \rho \in X_*(T)\) is the sum of the positive coroots (for the order defined by some choice of Borel subgroup of \(H\) containing \(T\)).
\end{defi}

\begin{defi} \label{def:ICcal}
  For \(n \geq 1\) and \(\Gbf = \mathbf{Sp}_{2n}\), we denote by \(\ICcal(\Gbf)\) the set of C-algebraic regular semisimple \(\Ghat(\C)\)-orbits in the Lie algebra \(\widehat{\gfrak} \simeq \mathfrak{so}_{2n+1}(\C)\) of \(\Ghat_{\C}\).
  Using the parametrization of \(\Tcal_{\SO_{2n+1}}\) given in the previous section, these are exactly the orbits of \((w_1, \dots, w_n) \in \Lie \Tcal_{\SO_{2n+1}}\) where \(w_1 > \dots > w_n > 0\) are integers.

  For \(n \geq 1\) and \(\Gbf=\SObf_{2n}\), write \(\theta\) for an outer automorphism of \(\Gbf\) induced by an element of \(\mathbf{O}_{2n}(\Z)\) of determinant \(-1\).
  Let \(\ICcal(\Gbf)\) be the set of C-algebraic regular semisimple \(\Ghat(\C)\)-orbits in \(\widehat{\gfrak}\) which are not invariant under \(\thetahat\), and let \(\ICcalt(\Gbf)\) be the set of \(\{ 1, \thetahat \}\)-orbits in \(\ICcal(\Gbf)\).
  More concretely these are the orbits of \((w_1, \dots, w_n) \in \Lie \Tcal_{\SO_{2n}}\) where \(w_1 > \dots > w_n > 0\) are integers.
  In order to be able to treat both cases simultaneously we also denote \(\ICcalt(\Spbf_{2n}) = \ICcal(\Spbf_{2n})\).
\end{defi}

For an integer \(m \geq 1\) a cuspidal automorphic representation \(\pi = \otimes'_v \pi_v\) for \(\PGLbf_{m,\Q}\), the infinitesimal character of \(\pi_\infty\) may be seen as a semisimple conjugacy class in \(\mathfrak{sl}_m(\C)\), or equivalently (considering its eigenvalues) as a multiset of cardinality \(m\) (orbit in \(\C^m\) under the symmetric group) \(\{ w_i(\pi_\infty) \,|\, 1 \leq i \leq m \}\).
We say that \(\pi\) has level one, or is everywhere unramified, if for all primes \(p\) we have \(\pi_p^{\PGLbf_m(\Zp)} \neq 0\).
We recall notation from \cite[p.275]{Taibi_dimtrace}
For \(n \geq 1\) and a family \((w_i)_{1 \leq i \leq n}\) satisfying \(w_i \in \Z\) and \(w_1 > \dots > w_n > 0\) we let \(O_e(w_1, \dots, w_n)\) (resp.\ \(O_o(w_1, \dots, w_n)\), resp.\ \(S(w_1, \dots, w_n)\)) be the set of self-dual level one cuspidal automorphic representations \(\pi\) for \(\PGLbf_{m,\Q}\) such that \(\pi_\infty\) has infinitesimal character \(\{ \pm w_i \,|\, 1 \leq i \leq n\}\) (resp.\ \(\{ \pm w_i \,|\, 1 \leq i \leq n \} \sqcup \{0\}\), resp.\ \(\{ \pm w_i \,|\, 1 \leq i \leq n\}\)).
We consider formal (unordered) finite sums \(\psi = \bigoplus_{i \in I} \pi_i[d_i]\) of pairs \((\pi_i, d_i)\) where \(\pi_i\) is in one of these three sets (i.e.\ \(\pi_i\) is a self-dual level one cuspidal automorphic representations for \(\PGLbf_{n_i,\Q}\) such that the eigenvalues of the infinitesimal character of \(\pi_{i,\infty}\) are distinct and either all integers or all in \(\frac{1}{2} + \Z\)) and \(d_i \geq 1\) are integers.
The notation \(\pi_i[d_i]\) (instead of \((\pi_i, d_i)\)) is meant to suggest the tensor product of the conjectural Langlands parameter of \(\pi_i\) with the irreducible \(d_i\)-dimensional representation of \(\SL_2\).
To such a formal sum \(\psi\) we associate an ``infinitesimal character'' which is the multiset of half-integers
\begin{equation} \label{eq:inf_char_sum_pii_di}
  \left\{ w_j(\pi_{i,\infty}) + \frac{d_i-1}{2} - k \,\middle|\, i \in I,\, 1 \leq j \leq n_i,\, 0 \leq k \leq d_i-1  \right\}
\end{equation}

\begin{defi} \label{def:Psi}
  For \(n \geq 1\) and \(\Gbf = \Spbf_{2n}\) or \(\SObf_{2n}\) and \(\taut \in \ICcalt(\Gbf)\) we let \(\Psit_\disc^{\unr, \taut}(\Gbf)\) denote the set of formal global Arthur-Langlands parameters \(\psi = \oplus_{i \in I} \pi_i[d_i]\) for \(\Gbf\) which are unramified at all finite places and have infinitesimal character equal to the image of \(\taut\) under the standard representation of \(\widehat{\gfrak}\).
  In particular this infinitesimal character is a genuine set, i.e.\ the terms in \eqref{eq:inf_char_sum_pii_di} are distinct.
  Let \(\Psit^{\unr, \taut}_{\disc, \nonendo}(\Gbf)\) be the subset of \(\Psit^{\unr, \taut}_{\disc}(\Gbf)\) consisting of all \(\psi\) as above for which \(|I|=1\), i.e.\ \(\psi = \pi[d]\).
\end{defi}

These sets \(\Psit_\disc^{\unr, \taut}(\Gbf)\) may thus be described combinatorially in terms of the sets \(O_e(\dots)\), \(O_o(\dots)\) and \(S(\dots)\) introduced above (see \cite[p.310]{Taibi_dimtrace} for more details).
Note that general substitutes for discrete parameters were defined in \cite[\S 1.4]{Arthur} (see also the review in \cite[\S 2]{Taibi_mult}.
Note that the definition of the present article differs from that of \cite{Taibi_dimtrace}: we find it more convenient to use infinitesimal characters \(\tau\) in this article, whereas we used dominant weights \(\lambda\) for \(\Gbf\) in loc.\ cit.
The two simply differ by half the sum of the positive roots for \(\Gbf\).

\begin{rema} \label{rem:SO2mod4_no_param}
  For \(n\) odd we always have \(\Psit_\disc^{\unr, \taut}(\SObf_{2n}) = \emptyset\): see \cite[Proposition 3.6]{ChRe} or \cite[3.\ p.309]{Taibi_dimtrace}.
  In particular the sets \(O_e(w_1, \dots, w_n)\) are empty for all odd \(n\) and all integers \(w_1 > \dots > w_n\).
\end{rema}

\begin{rema}[{\cite[Remark 4.1.6]{Taibi_dimtrace}}] \label{rem:no_param_unless_triv_res_center}
  For \(\Gbf = \Spbf_{2n}\) (resp.\ \(\SObf_{4n}\)) and \(\tau = (w_1, \dots, w_n) \in \ICcal(\Gbf)\) (resp.\ \(\taut = (w_1, \dots, w_{2n}) \in \ICcalt(\Gbf)\)) we have \(\Psit_\disc^{\unr, \taut}(\Gbf) = \emptyset\) unless \(\sum_{i=1}^n w_i \equiv n(n+1)/2 \mod 2\) (resp.\ \(\sum_{i=1}^{2n} w_i \equiv n \mod 2\)).
\end{rema}

\begin{defi} \label{def:inf_char}
  Let \(H\) be a connected complex reductive group, \(T\) a maximal torus of \(H\).
  \begin{enumerate}
  \item Let \(\varphi: \C^\times \to H(\C)\) be a continuous semisimple morphism.
    Up to conjugacy by \(H(\C)\) we may assume that \(\varphi\) takes values in \(T(\C)\), and then we have
    \[ \varphi(z) = z^{\tau} \ol{z}^{\tau'} := (z/|z|)^{\tau-\tau'} |z|^{\tau+\tau'} \]
    for uniquely determined \(\tau, \tau' \in \C \otimes_{\C} X_*(T) \simeq \Lie T\) satisfying \(\tau - \tau' \in X_*(T)\).
    We call the \(H(\C)\)-conjugacy class of \(\tau\) the infinitesimal character of \(\varphi\).
  \item For a continuous semisimple morphism \(\psi: \C^\times \times \SL_2(\C) \to H(\C)\) we define the infinitesimal character of \(\psi\) as that of
    \begin{align*}
      \varphi_{\psi}: \C^\times & \longrightarrow H(\C) \\
      z & \longmapsto \psi(z, \diag(|z|, |z|^{-1})).
    \end{align*}
  \end{enumerate}
\end{defi}

These definitions are meant to be applied to the case where \(H = \Ghat\) for some connected real reductive group \(\Gbf\) and \(\varphi\) (resp.\ \(\psi\)) is the restriction of a continuous semisimple morphism \(W_{\R} \to {}^L \Gbf\) (resp.\ \(W_{\R} \times \SL_2(\C) \to {}^L \Gbf\)).

\begin{defi} \label{def:Mpsi}
  \begin{enumerate}
  \item Let \(n \geq 1\), \(\Gbf = \Spbf_{2n}\) or \(\SObf_{4n}\) and \(\taut \in \ICcalt(\Gbf)\).
      For \(\psi \in \Psit^{\unr, \taut}_{\disc, \nonendo}(\Gbf)\), let \(\Mpsi\) be \emph{a} group isomorphic to \(\Ghat\), endowed with \(\tau_{\psi}\) a semisimple conjugacy class in \(\C \otimes_{\Qbar} \Lie \Mpsi\) such that \(\tau_{\psi}\) maps to \(\tilde{\tau}\) (this condition does not depend on the choice of an isomorphism \(\Mpsi \simeq \Ghat\)).

      Write \(\psi = \pi[d]\) where \(\pi = \pi_{\infty} \otimes \pi_f\) is a self-dual automorphic cuspidal representation of a general linear group and \(d \geq 1\).
      Let \(\psi_{\infty} : W_{\R} \times \SL_2(\C) \to \Mpsi(\C)\) be a continuous morphism, bounded on \(W_{\R}\) and algebraic on \(\SL_2(\C)\), such that:
      \begin{itemize}
        \item \(\Std_{\Gbf} \circ \psi_{\infty} : W_{\R} \times \SL_2(\C) \to \GL_{N(\Ghat)}(\C)\) is the localization of \(\psi\) at the real place, that is the tensor product of the Langlands parameter of \(\pi_{\infty}\) with the irreducible algebraic representation of \(\SL_2(\C)\) of dimension \(d\),
        \item the infinitesimal character of \(\psi_{\infty}\) is \(\tau_{\psi}\).
      \end{itemize}
      Such a \(\psi_{\infty}\) exists and is unique up to conjugation by \(\Mpsi(\C)\).

    \item Let \(n \geq 1\), \(\Gbf = \Spbf_{2n}\) or \(\SObf_{4n}\) and \(\taut \in \ICcalt(\Gbf)\).
      For \(\psi = \oplus_i \psi_i \in \tilde{\Psi}^{\unr, \tilde{\tau}}_{\disc}(\Gbf)\), let \(\Mpsi = \prod_i \mathcal{M}_{\psi_i}\), endowed with \(\tau_{\psi} = (\tau_{\psi_i})_i\).
      Let \(\psi_{\infty} : W_{\R} \times \SL_2(\C) \to \Mpsi(\C)\) be the morphism \(((\psi_i)_{\infty})_i\), well-defined up to conjugation by \(\Mpsi(\C)\).
      Let \(\varphi_{\psi_\infty} : W_{\R} \to \Mpsi(\C)\) be the composition of \(\psi_\infty\) with the morphism \(W_{\R} \to W_{\R} \times \SL_2(\C)\), \(w \mapsto (w, \diag(|w|^{1/2}, |w|^{-1/2}))\).
      By definition this parameter has infinitesimal character \(\tau_\psi\).

    \item
      Let \(n \geq 1\), \(\Gbf = \Spbf_{2n}\) or \(\SObf_{4n}\) and \(\tau \in \ICcal(\Gbf)\).
      For \(\psi = \oplus_i \psi_i \in \tilde{\Psi}^{\unr, \tilde{\tau}}_{\disc}(\Gbf)\), let \(\dpsitau : \Mpsi \rightarrow \Ghat\) be such that \(\Std_{\Gbf} \circ \dpsitau \simeq \oplus_i \Std_{\psi_i}\) and \(\mathrm{d} \dpsitau ( \tau_{\psi} ) = \tau\) (as semisimple conjugacy classes in \(\C \otimes_{\Qbar} \hat{\gfrak}\)).
      Then \(\dpsitau\) is uniquely determined up to composing with \(\Ad(g)\) for some \(g \in \Ghat(\Qbar)\), and \(\dot{\psi}_{\theta(\tau)} = \hat{\theta} \circ \dpsitau\).

      Let \(C_{\dpsitau} = \Cent( \dpsitau, \Ghat )\), a subgroup
      \footnote{More precisely, it is a finite group scheme over \(\Qbar\), but we
        will abusively identify it with its group of points over (any extension
      of) \(\Qbar\).}
      of the finite abelian \(2\)-torsion group \(C_{\dpsitau \circ
      \psi_{\infty}}\).
      Then \(\dpsitau\) induces an isomorphism \(Z(\Mpsi) \to C_{\dpsitau}\).
      In particular any other choice of \(\dpsitau\) yields a canonically
      isomorphic \(C_{\dpsitau}\), and for this reason we will often simply denote
      this group \(C_{\psi, \tau}\).
      This also shows that in the case where \(\Gbf = \SObf_{4n}\) we have a
      canonical isomorphism between \(C_{\psi, \tau}\) and \(C_{\psi,
      \thetahat(\tau)}\).
      For this reason we will often simply write \(C_{\psi}\) for \(C_{\psi,
      \tau}\).
    Let \(\Scal_\psi = C_\psi / Z(\Ghat)\).
    
  \end{enumerate}
\end{defi}

Note that the pair \((\Mpsi, \tau_{\psi})\) is well-defined up to an isomorphism which is unique up to \(\Mcal_{\psi,\ad}(\Qbar)\).
Of course this is already the case for \(\Mpsi\) (without the need for \(\tau_\psi\)) in the case where \(\Gbf = \mathbf{Sp}_{2n}\).

In the setting of the third part of Definition \ref{def:Mpsi} Arthur defined and element \(s_\psi \in C_\psi\) (see \cite[\S 1.4]{Arthur}) and a character \(\epsilon_\psi\) of \(\Scal_\psi\) (see (1.5.6) loc.\ cit.).

\begin{lemm} \label{lem:epsilon_psi_s_psi}
  Let \(n \geq 1\), \(\Gbf = \Spbf_{2n}\) (resp.\ \(\SObf_{4n}\)), \(\taut \in \ICcalt(\Gbf)\) and \(\psi = \oplus_{i \in I} \pi_i[d_i] \in \Psit_{\disc}^{\unr, \taut}(\Gbf)\).
  Let \(I_\even\) be the set of indices \(i \in I\) for which \(\pi_i[d_i]\) is even-dimensional (i.e.\ \(d_i\) is even or \(\pi_i\) is an automorphic representation for \(\GLbf_{n_i}\) where \(n_i\) is even).
  The set \(I \smallsetminus I_\even\) has one element (resp.\ is empty).
  For \(i \in I_\even\) let \(s_i \in Z(\Mpsi) = \prod_j Z(\Mcal_{\pi_j[d_j]})\) be the element which is trivial at all indices \(j \neq i\) and non-trivial at \(i\).
  \begin{enumerate}
  \item The group \(C_\psi\), seen as a vector space over \(\F_2\), admits \((s_i)_{i \in I_\even}\) as basis.
  \item The element \(s_\psi\) is the product of \(s_i\) over all indices \(i \in I\) such that \(d_i\) is even.
  \item We have
    \[ \epsilon_\psi(s_i) = \prod_{j \in I \smallsetminus \{i\}} \epsilon(1/2, \pi_i \times \pi_j)^{\min(d_i,d_j)} = \prod_{\substack{j \in I \smallsetminus \{i\} \\ d_i \not\equiv d_j \mod 2}} \epsilon(1/2, \pi_i \times \pi_j)^{\min(d_i,d_j)}. \]
  \end{enumerate}
\end{lemm}
\begin{proof}
  The first two assertions are clear, for the last one see \cite[(3.10)]{ChRe} and the following reference to \cite[Theorem 1.5.3 (b)]{Arthur} therein.
\end{proof}

\subsection{Stabilization of the trace formula}
\label{sec:STF}

We will use repeatedly the stabilization of the trace formula.
We state special cases that will be enough for our purposes.

For \(\Gbf\) a (connected) reductive group over \(\Qp\) (resp.\ \(\Q\)) and \(K\) a
compact open subgroup of \(\Gbf(\Qp)\) (resp.\ \(\Gbf(\A_f)\)) we denote
\(\Hcal(\Gbf(\Qp) // K)\) (resp.\ \(\Hcal(\Gbf(\A_f) // K)\)) the Hecke algebra in
level \(K\) with coefficients in \(\Q\).
Also denote \(\Hcal(\Gbf(\Qp))\) (resp.\ \(\Hcal_f(\Gbf)\)) the direct limit over
all compact open subgroups \(K\) (fixing a Haar measure on \(\Gbf(\Qp)\) resp.\
\(\Gbf(\A_f)\) identifies elements of this Hecke algebra with smooth compactly
supported functions).
For \(\Gbf\) connected reductive over \(\Z\) denote \(\Hcal_f^{\unr}(\Gbf) =
\Hcal(\Gbf(\A_f) // \Gbf(\Zhat))\).
For an extension \(F\) of \(\Q\) we add a subscript \(F\) to denote the analogous
Hecke algebras with coefficients in \(F\).
For \(\Gbf\) connected reductive over \(\R\) and \(K\) a maximal compact subgroup of
\(\Gbf(\R)\) denote \(\Hcal(\Gbf(\R), K)\) the Hecke algebra of smooth compactly
supported bi-\(K\)-finite distributions on \(\Gbf(\R)\).
When \(K\) plays no particular role we will often omit it from the notation.
For \(\Gbf\) connected reductive over \(\Q\) denote \(\Hcal(\Gbf) = \Hcal(\Gbf(\R))
\otimes \Hcal_f(\Gbf)_\C\), and let \(I(\Gbf)\) be the quotient of \(\Hcal(\Gbf)\) by
the subspace of distributions all of whose orbital integrals at semi-simple
regular elements of \(\Gbf(\Q_v)\) vanish, for any place \(v\) of \(\Q\).
Similarly define \(SI(\Gbf)\) by considering stable orbital integrals instead.
These quotients have obvious local analogues \(I(\Gbf(\Q_v))\) and
\(SI(\Gbf(\Q_v))\) at any place \(v\).
Finally for \(\Gbf\) connected and reductive over \(\Z\) we also denote
\(\Hcal^{\unr}(\Gbf) = \Hcal(\Gbf(\R)) \otimes \Hcal_f^{\unr}(\Gbf)_\C\).

For a split reductive group \(\Gbf\) over \(\Zp\), recall that the Satake morphism
\cite{Satake} \(\Sat_{\Gbf}\) is an isomorphism between \(\Hcal(\Gbf(\Qp) //
\Gbf(\Zp))_\C\) and the \(\C\)-algebra of invariant algebraic functions on the dual
group \(\Ghat\) (a basis over \(\C\) of this space of functions being given by
traces in irreducible representations of \(\Ghat\)).
Characters of this algebra are known to correspond to semi-simple conjugacy
classes in \(\Ghat(\C)\).
For example for \(n \geq 1\), \(\tilde{\tau} \in \tilde{\ICcal}(\mathbf{SO}_{4n})\)
and \(\psi \in \tilde{\Psi}^{\unr, \tilde{\tau}}_{\disc,
\nonendo}(\mathbf{SO}_{4n})\), we have an associated
\(\mathrm{O}_{4n}(\C)\)-conjugacy class of unramified Arthur parameters \(\psi_p:
\mathrm{W}_{\Qp}/\mathrm{I}_{\Qp} \times \SL_2(\C) \to \SO_{4n}(\C)\), and the
image of \((\Frob_p, \diag(p^{1/2}, p^{-1/2}))\) by \(\psi\) defines a \(\{1,
\thetahat \}\)-orbit \(\tilde{c}_p(\psi)\) of semisimple conjugacy classes in
\(\SO_{4n}(\C)\).
Similarly, for \(n \geq 1\), \(\tau \in \ICcal(\Spbf_{2n})\) and \(\psi \in
\Psit^{\unr, \tau}_{\disc, \nonendo}(\Spbf_{2n})\), we have an associated
semisimple conjugacy class \(c_p(\psi)\) in \(\mathrm{SO}_{2n+1}(\C)\).
So for \(\Gbf = \Spbf_{2n}\) or \(\SObf_{4n}\), \(\tau \in \ICcal(\Gbf)\) and \(\psi
\in \Psit^{\unr, \tau}_{\disc}(\Gbf)\) we have an associated orbit
\(\tilde{c}_p(\psi)\) of semisimple conjugacy classes in \(\Mpsi\) under the product
of the group of outer automorphisms of all even orthogonal factors of \(\Mpsi\).

\begin{rema} \label{rem:mult_one}
  As is used repeatedly in \cite{Arthur}, note that the results of \cite{JacSha_Euler2} and \cite{MoeWal_resGL} imply that \((\tilde{c}_p(\psi))_p\) determines \(\psi\) (even if we discard finitely many primes).
\end{rema}

We now state the special cases of the stabilization of the trace formula that we will use.
Let \(\Gbf\) be a split reductive group over \(\Z\) which is a quotient of a product \(\Gbf_{\mathrm{cl}}\) of groups isomorphic to \(\SObf_{2n}\) and \(\Spbf_{2n}\).
For \(\tau \in IC(\Gbf)\) let \(\Acal^2(\Gbf(\Q) \backslash \Gbf(\A) / \Gbf(\Zhat))_\tau\) be the space of square-integrable automorphic forms\footnote{The definition involves the maximal compact subgroup \(K\) of \(\Gbf(\R)\), but we omit it from the notation.} of level one which are eigenvectors for the infinitesimal character corresponding to \(\tau\).
Let \(I_{\disc, \tau}^\Gbf\) be the linear form on \(\Hcal^\unr(\Gbf)\) given by the
trace on \(\Acal^2(\Gbf(\Q) \backslash \Gbf(\A) / \Gbf(\Zhat))_\tau\) (by a
famous theorem of Harish-Chandra, if \(K_\infty\) is our chosen maximal compact
subgroup of \(\Gbf(\R)\) then any isotypical component of \(\Acal^2(\Gbf(\Q)
\backslash \Gbf(\A) / \Gbf(\Zhat))_\tau\) for some irreducible representation of
\(K_\infty\) is in fact finite-dimensional, so the trace we are considering is
simply that of an endomorphism of a finite-dimensional vector space).

Consider endoscopic data \(\efrak = (\Hcal, \Hbf, s, \xi)\) for \(\Gbf\).
For convenience we modify slightly the definition in \cite[\S 2.1]{KS}: instead
of letting \(\xi\) be an embedding of \(\Hcal\) in \({}^L \Gbf\), we let \(\Hcal\) be a
subgroup of \({}^L \Gbf\) and take \(\xi\) to be an isomorphism \(\Hhat \to \Hcal\).
We will mostly only need to consider elliptic endoscopic data which are
unramified at all finite places of \(\Q\).
Using the explicit description of all endoscopic data for special orthogonal and
symplectic groups and Minkowski's theorem on unramified extensions of \(\Q\) it is
easy to see that such endoscopic data are obtained as follows: \(s\) is a
semisimple element of \(\Ghat\) whose image in \(\widehat{\Gbf_{\mathrm{cl}}}\) has
order \(1\) or \(2\), \(\xi(\Hhat) = \Cent(s, \Ghat)^0\) and \(\Hcal = \xi(\Hhat)
\times W_F\).
(Not quite all of these give rise to \emph{elliptic} endosopic data: the condition for ellipticity is that the connected centralizer of \(s\) in \(\widehat{\Gbf_{\mathrm{cl}}}\) should not have any factor \(\SO_2\) except in factors of \(\widehat{\Gbf_{\mathrm{cl}}}\) which are themselves isomorphic to \(\SO_2\).)
Note in particular that any everywhere unramified elliptic endoscopic datum for
\(\Gbf\) is induced by one for \(\Gbf_{\ad}\), that \(\Hbf\) is also split and is also
a quotient of a product of symplectic and split even special orthogonal groups.
Thus there is an obvious way to extend \(\xi\) to an isomorphism \({}^L \Hbf \simeq
\Hcal\), and we will always use this L-embedding of \({}^L \Hbf\) into \({}^L \Gbf\).

We fix a global Whittaker datum \(\wfrak\) for \(\Gbf\), i.e.\ for a Borel pair
\((\Bbf, \Tbf)\) of \(\Gbf\) (defined over \(\Q\), and we may even take it to be
defined over \(\Z\)), denoting by \(\Ubf\) the unipotent radical of \(\Bbf\) and
\(\Ubf^\ab\) its largest abelian quotient (a vector group of dimension equal to
the rank of \(\Gbf\)), \(\wfrak\) is a generic unitary character of
\(\Ubf(\A)/\Ubf(\Q)\).
Denoting by \(\Tbf_\ad\) the image of \(\Tbf\) in the adjoint group of \(\Gbf\), the
group \(\Tbf(\Q)\) acts transitively on the set of such characters.
Since \(\Q^\times = \{ \pm 1 \} \times \bigoplus_{p \in \Pcal} p^\Z\), we see that
we may choose \(\wfrak\) such that for any prime number \(p\), the localization
\(\wfrak_p: \Ubf(\Qp) \to \C^\times\) is compatible with the \(\Zp\)-structure on
\(\Gbf_{\Qp}\) in the sense of \cite{CasSha} (see also \cite[\S
7]{Hales_unr_factors}).
In fact we even see that we may also impose that \(\wfrak_\infty\) lie in any
given \(\Gbf(\R)\)-orbit of Whittaker data for \(\Gbf_\R\).
The choice of \(\wfrak\) gives us, for any endoscopic datum \(\efrak = (\Hcal,
\Hbf, s, \xi)\) for \(\Gbf\) and any choice of L-embedding \({}^L \xi : {}^L \Hbf
\simeq \Hcal\) extending \(\xi\), a decomposition of the (canonical) global
transfer factor into a product over all places of \(\Q\) of normalized transfer
factors \cite[\S 5.3]{KS}.
In particular we have a notion of endoscopic transfer (\S 5.5 loc.\ cit.) from
\(I(\Gbf)\) to \(SI(\Hbf)\), which is the tensor product over all places of local
transfers.
Transfer is known to exist in general (in the real case see
\cite{Shelstad_temp1} for the case of Schwartz functions,
\cite{Shelstad_tw_geom_transf} for compactly supported functions, and
\cite[\S IV.3.4]{SFTT1}), but we will be mainly interested in particular cases
where it is somewhat explicit.

At a finite place \(p\) of \(\Q\) we will only use \(f_p \in \Hcal^{\unr}(\Gbf_{\Zp})\).
In this case endoscopic transfer is very explicit: it vanishes unless \(\efrak\)
is unramified.
If this holds, then transfer can be made explicit in terms of Satake transforms
(this statement is the version of the unramified fundamental lemma deduced in
\cite{Hales} from the case of the unit element for sufficiently large residual
characteristic proved by \cite{Ngo_LF}, \cite{Waldspurger_chgmt_car},
\cite{Waldspurger_pas_si_tordue}).
For global reasons we will only need the case where \(\Hbf\) is also split, and so
we recall the statement of the fundamental lemma only in this case.
Thus we may take the obvious embedding \({}^L \Hbf \to {}^L \Gbf\).
Choose a hyperspecial compact open subgroup (equivalently, a reductive model of
\(\Hbf\) over \(\Zp\)) \(\Hbf(\Zp)\) of \(\Hbf(\Qp)\).
We have a commutative diagram
\[ \begin{tikzcd}[column sep=5em]
  \Ocal(\Ghat)^{\Ghat} \arrow[r, "{\xi^*}"] & \Ocal(\Hhat)^{\Hhat} \\
  \Hcal^{\unr}(\Gbf_{\Zp}) \arrow[d] \arrow[u, "{\Sat_{\Gbf_{\Zp}}}" left,
    "{\sim}" right] & \Hcal^{\unr}(\Hbf_{\Zp}) \arrow[d] \arrow[u,
    "{\Sat_{\Hbf_{\Zp}}}" left, "{\sim}" right] \\
  I(\Gbf(\Qp)) \arrow[r] & SI(\Hbf(\Qp))
\end{tikzcd} \]
where the bottom horizontal map is endoscopic transfer.

At the Archimedean place of \(\Q\) we will mainly (but not exclusively) use the pseudo-coefficients \(f_\sigma(g) dg\) for discrete series representations \(\sigma\) of \(\Gbf(\R)\) given by \cite[Proposition 4, Corollaire]{CloDel}.
It follows from the endoscopic character relations \cite{Shelstad_temp2} and \cite[Lemma 5.3]{Shelstad_characters} that a suitable linear combination of pseudo-coefficients of discrete series representations of \(\Hbf(\R)\) is a transfer of \(f_\sigma(g) dg\).
In particular the transfer always vanishes (in \(SI(\Hbf_\R)\)) if \(\Hbf(\R)\) does not admit discrete series.
Note that Shelstad's lemma also implies that pseudo-coefficients of discrete series in the same L-packet have the same stable orbital integrals.
We denote by \(\PsCo(\Gbf, \tau)\) the (finite-dimensional) subspace of \(I(\Gbf_{\R})\) spanned by the pseudo-coefficients of discrete series representations having infinitesimal character \(\tau\).

Recall from \cite{Kottwitz_STFcusptemp} that associated to an elliptic endoscopic datum \(\efrak = (\Hbf, \Hcal, s, \xi)\) for \(\Gbf\) is a positive rational number \(\iota(\efrak) = \tau(\Gbf) \tau(\Hbf)^{-1} |\Out(\efrak)|\) (see also \cite[\S 2.1]{KS} for a definition of \(\Out(\efrak)\) with the formulation that we adopted here).
Since we will only need to consider groups \(\Gbf\) which are either split or inner forms of split groups, we simply have \(\tau(\Gbf) = |\pi_0(Z(\Ghat))|\), and similarly for \(\Hbf\).

\begin{theo}[Specialization of the stabilization of the trace formula: split case]
  \label{thm:STF_qs}
  For any group \(\Gbf\) as above, and any regular infinitesimal character \(\tau\), the linear form \(S_{\disc,\tau}^\Gbf\) recursively defined on \(\Hcal^\unr(\Gbf)\) by
  \[ S_{\disc, \tau}^\Gbf(f) = I_{\disc, \tau}^\Gbf(f) - \sum_{\efrak} \iota(\efrak) \sum_{\tau' \mapsto \tau} S_{\disc, \tau'}^\Hbf(f') \]
  where the first sum ranges over all non-trivial elliptic endoscopic data \(\efrak = (\Hbf, \Hcal, s, \xi)\) for \(\Gbf\) which are unramified at every prime and the second sum ranges over semisimple conjugacy classes in \(\widehat{\hfrak}\) mapping to \(\tau\) by the differential of \(\xi\), is stable.
\end{theo}
\begin{proof}
  This is a special case of the stabilization of the trace formula (see
  Global Theorem 2 and Lemma 7.3 (b) in \cite{ArthurSTF1}), refined by
  infinitesimal characters as in \cite[\S X.5.2]{SFTT2}.
  Note that in general the discrete part of the trace formula contains terms
  other than the trace in the discrete spectrum, but since \(\tau\) is regular the
  \(\tau\)-part of these other terms vanishes (see the argument on p.\ 268 of
  \cite{ArthurL2}).
  Note that any infinitesimal character \(\tau'\) which maps to a regular \(\tau\)
  is also regular.
\end{proof}

We will also need a special case of the stabilization of the trace formula for
certain non-quasisplit groups.
For the rest of this section \(\Gbf\) denotes an inner form of \(\PGSObf_{4n}\) which is split at all finite places of \(\Q\).
(In this paper we will only need the case where \(\Gbf(\R)\) is compact, which occurs if and only if \(n\) is even.)
The group \(\Gbf\) admits a reductive model over \(\Z\) (it may admit several non-isomorphic models if \(\Gbf(\R)\) is compact), and we fix such a model.
We also fix a realization of \(\Gbf\) as an inner form of \(\PGSObf_{4n}\), i.e.\ an isomorphism \(\Xi: \PGSObf_{4n,\Qbar} \simeq \Gbf_{\Qbar}\) such that for any \(\sigma \in \Gal(\Qbar/\Q)\) the automorphism \(c(\sigma) := \Xi^{-1} \sigma(\Xi)\) is inner.
In particular this defines a \(1\)-cocycle \(c \in Z^1(\Q, \PGSObf_{4n}(\Qbar))\).
Since the split group \(\PGSObf_{4n}\) is of adjoint type, the group
\(\PGSObf_{4n}(\Q)\) acts transitively on the set of global Whittaker data for
\(\PGSObf_{4n}\).
This (essentially unique) global Whittaker datum for \(\PGSObf_{4n}\) gives, for any endoscopic datum for \(\Gbf\), a decomposition of the adelic transfer factors as a product of local transfer factors: see \cite[\S 4.4]{Kalgri}.
For any prime number \(p\), the image of \(c\) in \(H^1(\Gal(\Qpbar/\Qp),
\PGSObf_{4n}(\Qpbar))\) is trivial, and this gives an isomorphism \(\Gbf_{\Qp}
\simeq \PGSObf_{4n}\), well-defined up to composition with an inner automorphism.
The local transfer factors at \(p\) for \(\Gbf\) are simply pulled back from those
for \(\PGSObf_{4n}\), in particular the transfer of unramified elements of the
Hecke algebra is given by the fundamental lemma as in the previous case.

\begin{theo}[Specialization of the stabilization of the trace formula: inner
  forms case]
  \label{thm:STF_inner}
  For such a group \(\Gbf\), for any \(f \in \Hcal^{\unr}(\Gbf)\) we have
  \[ I_{\disc, \tau}^\Gbf(f) = \sum_{\efrak} \iota(\efrak) \sum_{\tau' \mapsto
  \tau} S_{\disc, \tau'}^\Hbf(f') \]
  where the first sum ranges over all elliptic endoscopic data \(\efrak = (\Hbf,
  \Hcal, s, \xi)\) for \(\Gbf\) which are unramified at every prime.
\end{theo}
\begin{proof}
  This is deduced from the same references of the previous theorem, by the same
  argument.
\end{proof}

\subsection{Even orthogonal groups and outer automorphisms}

\begin{lemm} \label{lem:SO4n+2}
  For any \(n \geq 1\), letting \(\Gbf = \mathbf{SO}_{4n+2}\), for any algebraic
  regular infinitesimal character \(\tau\) for \(\Gbf\) which is not invariant under
  the outer automorphism of \(\Gbf\), the linear form \(S^{\Gbf}_{\disc,\tau}\)
  vanishes identically on \(\Hcal^{\unr}(\Gbf)\).
\end{lemm}
\begin{proof}
  Consider elliptic endoscopic data \(\mathfrak{e} = (\Hbf, \mathcal{H}, s, \xi)\)
  for \(\Gbf\) which are unramified at every finite place.
  We have \(\Hbf \simeq \mathbf{SO}_{4a} \times \mathbf{SO}_{4b+2}\) with \(a+b=n\)
  and \(b \geq 1\).
  By induction on \(n\) and using the definition (Theorem \ref{thm:STF_qs}), it is
  enough to show that \(I^{\Gbf}_{\disc, \tau}\) vanishes identically on
  \(\Hcal^{\unr}(\Gbf)\), i.e.\ that there is no discrete automorphic
  representation for \(\Gbf\) in level one and infinitesimal character \(\tau\).
  Now this follows from Arthur's multiplicity formula (\cite[Theorem
  1.5.2]{Arthur}, see \cite[Theorem 4.1.2]{Taibi_dimtrace} for the
  specialization to the everywhere unramified case) and the fact that there
  are no parameters with the appropriate infinitesimal character
  (Remark \ref{rem:SO2mod4_no_param}).
  Note that the multiplicity formula only describes the representations
  occurring in the discrete spectrum up to outer automorphism (at all places
  independently, i.e.\ \(\prod_v \Z/2\Z\)-orbits), but this is enough to imply
  vanishing of \(I^{\Gbf}_{\disc, \tau}\) since it is the trace in a genuine
  representation.
\end{proof}

Recall from \cite[Lemmas 3.7 and 3.15]{ChRe} or
\cite[pp.309--310]{Taibi_dimtrace} that for \(\Gbf = \Spbf_{2n}\) or \(\SObf_{4n}\),
for any \(\tau \in \ICcal(\Gbf)\) and any \(\psi \in \Psit_\disc^{\unr,
\taut}(\Gbf)\), the continuous semisimple morphism \(\dpsitau \circ \psi_\infty:
W_\R \times \SL_2 \to \Ghat\) (well-defined up to conjugation) is bounded on
\(W_\R\) (essentially by Clozel's purity lemma) and has C-algebraic infinitesimal
character, thus it is an Adams-Johnson parameter (see \cite[p.194]{Kottwitz_AA},
\cite[\S 5]{ArthurUnip}, \cite[\S 4.2.2]{Taibi_dimtrace}).
Adams and Johnson \cite{AdJo} attached a finite set \(\Pi^\AJ(\Gbf_\R, \dpsitau
\circ \psi_\infty)\) of irreducible unitary representations of \(\Gbf(\R)\) to such
a parameter, and since we have fixed a Whittaker datum for \(\Gbf\), we also get a
map from \(\Pi^\AJ(\Gbf_\R, \dpsitau \circ \psi_\infty)\) to the group of
characters of the finite abelian \(2\)-torsion group \(C_{\dpsitau \circ
\psi_\infty}/Z(\Ghat)\) (see the above-cited references for this formulation,
which does not appear in \cite{AdJo}).
This map will be denoted \(\pi \mapsto \langle \cdot, \pi \rangle\).
Moreover the linear form
\[ f(g) dg \mapsto \sum_{\pi \in \Pi^\AJ(\Gbf_\R, \dpsitau \circ \psi_\infty)} \langle s_\psi, \pi \rangle \Tr(\pi(f(g) dg)) \]
on \(I(\Gbf_\R)\), which we denote by \(\Lambda_{\dpsitau \circ \psi_\infty}^\AJ\), is stable, and the Adams-Johnson packet satisfies endoscopic character relations.
Adams-Johnson packets are closely related to discrete series L-packets, and we recall part of that connection.
Let \(\varphi_\tau: W_\R \to \Ghat\) be a discrete Langlands parameter having infinitesimal character \(\tau\) (such a parameter exists and is unique up to conjugation by \(\Ghat\)).
Recall from \cite[\S 4.2]{Taibi_dimtrace} that there is a certain quasi-split twisted Levi subgroup \(\Lbf^*_{\psi,\tau}\) of \(\Gbf_{\R}\) attached to \(\dpsitau \circ \psi_\infty\).

\begin{prop} \label{pro:tr_AJ_ps}
  For any \(\sigma\) in the L-packet \(\Pi(\Gbf_\R, \varphi_\tau)\), denoting by
  \(f_{\sigma}(g) dg\) a pseudo-coefficient for \(\sigma\), we have the following
  property.
  There exists exactly one \(\pi \in \Pi^\AJ(\Gbf_\R, \dpsitau \circ
  \psi_\infty)\) such that \(\Tr(\pi(f_{\sigma}(g) dg)) \neq 0\), and for this
  \(\pi\) we have
  \[ \Tr(\pi(f_{\sigma}(g) dg)) = \langle s_\psi, \pi \rangle (-1)^{q(\Lbf^*_{\psi,\tau})}.
  \]
  Assume that \(\dpsitau \circ \psi_\infty\) and \(\varphi_\tau\) are aligned, i.e.\ that the restrictions of  \(\dpsitau \circ \varphi_{\psi_\infty}\) and \(\varphi_\tau\) to \(\C^\times\) take values in the same maximal torus of \(\Ghat_{\C}\) and have the same holomorphic part (see Definition \ref{def:Mpsi}).
  This is always possible after conjugating by \(\Ghat(\C)\), and gives a canonical inclusion \(\Scal_{\dpsitau \circ \psi_\infty} \subset \Scal_{\varphi_\tau}\).
  We then have \(\langle \cdot, \pi \rangle = \langle \cdot, \sigma \rangle|_{\Scal_{\dpsitau \circ \psi_\infty}}\).

  In particular we have \(\Lambda_{\dpsitau \circ \psi_\infty}^\AJ(f_{\sigma}(g) dg) = (-1)^{q(\Lbf^*_{\psi,\tau})}\).
\end{prop}
\begin{proof}
  This follows from Johnson's resolution, see \cite[(4.2.1)]{Taibi_dimtrace},
  and the proof of Proposition 3.2.5 in \cite{Taibi_mult}.
  Note that this goes back to \cite{Kottwitz_AA} and \cite{ArthurUnip}, although
  the normalization by Whittaker datum was not available at the time.
\end{proof}

Recall that Arthur \cite{Arthur} also defined (in a more general setting) a
packet \(\Pi(\Gbf_\R, \dpsitau \circ \psi_\infty)\), which is a multi-set of
\(\theta\)-orbits of irreducible unitary representations of \(\Gbf\), and also
associated a character of \(C_{\dpsitau \circ \psi_\infty}/Z(\Ghat)\) to each
elements of this packet (all of this for more general parameters, and
characterized via twisted endoscopy for general linear groups).
Again the linear form \(\Lambda_{\dpsitau \circ \psi_\infty}\) defined as above
on \(I(\Gbf)^\theta\) is stable, and endoscopic character relations are satisfied.
The main result of \cite{AMR} is equivalent to the assertion that
\(\Lambda_{\dpsitau \circ \psi_\infty}\) is the restriction of \(\Lambda_{\dpsitau
\circ \psi_\infty}^\AJ\) to \(I(\Gbf)^\theta\), which implies that the packets
coincide and the associated characters are equal.
Note that in the orthogonal case where \(\theta\) is not trivial, we are assuming
that \(\tau\) is not fixed by \(\thetahat\) and so no element of \(\Pi^\AJ(\Gbf_\R,
\dpsitau \circ \psi_\infty)\) is fixed by \(\theta\), i.e.\ each orbit in
\(\Pi^\AJ(\Gbf_\R, \dpsitau \circ \psi_\infty)\) consists of two elements which
can be distinguished by their infinitesimal character.

\begin{lemm} \label{lem:pre_refineSOeven}
  Denote \(\Gbf = \SObf_{4n}\).
  Let \(\tau \in \ICcal(\Gbf)\) and \(\psi \in \Psit_\disc^{\unr, \taut}(\Gbf)\).
  Then for any \(\pi_\infty \in \Pi^\AJ(\Gbf_\R, \dpsitau \circ \psi_\infty)\)
  such that \(\langle \cdot, \pi_\infty \rangle|_{\Scal_{\psi}} = \epsilon_\psi\),
  there is a unique family \((c_p'(\Gbf, \psi, \pi_\infty))_{p \in \Pcal}\) of
  semisimple conjugacy classes in \(\Ghat\) such that for any prime \(p\) the
  conjugacy class \((c_p'(\Gbf, \psi, \pi_\infty))_{p \in \Pcal}\) belongs to the
  \(\thetahat\)-orbit \(\dpsitau(\tilde{c}_p(\psi))\) and the representation
  \(\pi_\infty \otimes \bigotimes'_p \pi_p\) of \(\Gbf(\A)\), where each \(\pi_p\) is
  unramified with Satake parameter \(c_p'(\Gbf, \psi, \pi_\infty)\), occurs in
  \(\Acal^2(\Gbf(\Q) \backslash \Gbf(\A))\).
  Moreover it occurs with multiplicity one.
\end{lemm}
\begin{proof}
  Of course this uses Arthur's endoscopic classification \cite[Theorem
  1.5.2]{Arthur}, specialized to level one (see \cite[Lemma
  4.1.1]{Taibi_dimtrace}).
  The parameter \(\psi\) does not have any odd-dimensional factor since \(\tau\) is
  not invariant under \(\thetahat\), and so in the notation of \cite{Arthur} we
  have \(m_\psi = 2\).
  Arthur's multiplicity formula (along with Remark \ref{rem:mult_one}) implies
  that there are exactly two irreducible representations \(\pi' = \pi'_\infty
  \otimes \pi'_f\) occurring in \(\Acal^2(\Gbf(\Q) \backslash \Gbf(\A))\) such that
  \(\pi'_\infty \in \{\pi_\infty, \pi_\infty^\theta\}\) and \(\pi'_f\) is everywhere
  unramified and at (almost) all prime numbers \(p\) its Satake parameter lies in
  \(\dpsitau(\tilde{c}_p(\psi))\).
  To conclude we only have to consider the global action of \(\theta\) on
  \(\Acal^2(\Gbf(\Q) \backslash \Gbf(\A))\) by \(f \mapsto (g \mapsto
  f(\theta(g)))\) which maps the isotypical spaces for an irreducible \(\pi'\) to
  that for \(\pi'^\theta\).
\end{proof}

\begin{prop} \label{pro:refineSOeven}
  There exists a unique family \((c_p(\psi))_{\psi, p}\), where \(\psi \in \Psit^{\unr, \taut}_{\disc, \nonendo}(\mathbf{SO}_{4n})\) for some \(n \geq 1\) and \(\tilde{\tau} \in \tilde{\ICcal}(\mathbf{SO}_{4n})\) and \(p \in \Pcal\), such that:
  \begin{itemize}
    \item
      \(c_p(\psi)\) is a semisimple conjugacy class in \(\Mpsi\) which belongs to
      the \(\{ 1, \hat{\theta} \}\)-orbit \(\tilde{c}_p(\psi)\),
    \item
      for any \(n \geq 1\), any \(\tau \in \ICcal(\mathbf{SO}_{4n})\), and any \(f(g)
      dg = f_{\infty} (g_\infty) d g_\infty \prod_p f_p(g_p) dg_p \in
      \Hcal^{\unr}(\mathbf{SO}_{4n})\),
      \begin{align*}
        S_{\mathrm{disc}, \tau}^{\mathbf{SO}_{4n}}(f(g) dg) = \sum_{\psi \in
          \tilde{\Psi}^{\unr, \tau}_{\disc}(\mathbf{SO}_{4n})} &
          \frac{\epsilon_{\psi}(s_{\psi})}{|\mathcal{S}_{\psi}|}
          \Lambda_{\dpsitau \circ \psi_{\infty}}^\AJ(f_{\infty}(g_\infty) d
          g_\infty) \\
        & \times \prod_p \mathrm{Sat}_{\mathbf{G}_{\Zp}}(f_p(g_p)
          dg_p)(\dpsitau( c_p(\psi))).
      \end{align*}
      Implicitly, we have written \(c_p(\psi)\) for \((c_p(\psi_i))_i\) for
      endoscopic \(\psi = \oplus_i \psi_i\).
  \end{itemize}
\end{prop}
\begin{proof}
  We proceed by induction on \(n\).
  Denote \(\Gbf = \SObf_{4n}\).
  To simplify notation in the proof we implicitly fix Haar measures on all
  groups that appear, thus identifying smooth compactly supported distributions
  (e.g.\ \(f(g) dg\)) with functions (e.g.\ \(f\)).
  We use the stabilization of the trace formula (Theorem \ref{thm:STF_qs}):
  \[ I^{\mathbf{G}}_{\mathrm{disc}, \tau}(f) = \sum_{\mathfrak{e}}
  \iota(\mathfrak{e}) \sum_{\tau' \mapsto \tau} S_{\mathrm{disc},
  \tau'}^{\mathbf{H}}(f') \]
  where the sum is over the set of isomorphism classes of elliptic endoscopic
  data \(\efrak = (\Hbf, \Hcal, s, \xi)\) for \(\mathbf{G}\) which are unramified at
  all finite places of \(\Q\).
  Therefore \(\mathbf{H} \simeq \mathbf{SO}_{2a} \times \mathbf{SO}_{2b}\) with
  \(a+b=2n\) and \(a,b \neq 1\).
  By the previous lemma, if \(a\) and \(b\) are odd then the contribution of
  \(\mathfrak{e}\) is zero: in the above formula \(\tau'\) is not invariant under
  the non-trivial outer automorphism of any of the two factor of \(\Hbf\).
  Note that even without assuming that \(\tau\) is not invariant under the
  non-trivial outer automorphism, at least one of the two factors of \(\tau'\) is
  not invariant.

  By the induction hypothesis, we have for all non-trivial \(\mathfrak{e}\)
  \[ S_{\mathrm{disc}, \tau'}^{\mathbf{H}}(f') = \sum_{\psi' \in \tilde{\Psi}_{\disc}^{\unr, \tilde{\tau}'}(\mathbf{H})} \frac{\epsilon_{\psi'}(s_{\psi'})}{ |\mathcal{S}_{\psi'}| } \Lambda_{\dot{\psi}'_{\tau'} \circ \psi_{\infty}}(f'_{\infty}) \times \prod_p \mathrm{Sat}_{\mathbf{H}_{\Zp}}(f'_p)(\dot{\psi}'_{\tau'}( c_p(\psi')). \]
  %TODO: give more details below?
  By the same argument as in the end of the proof of \cite[Theorem 4.0.1]{Taibi_mult} (which originates from \cite[\S 11]{Kottwitz_STFcusptemp}) and using that \({}^L \xi \circ \dot{\psi}'_{\tau'} = \dpsitau\) for \(\psi = {}^L \xi \circ \psi'\) and \(\tau = \xi(\tau')\), we have
  \begin{align} \label{eq:stab_spec_SOeven}
    \sum_{\text{non-trivial } \mathfrak{e}} \iota(e) \sum_{\tau' \mapsto \tau} S^{\mathbf{H}}_{\disc, \tau'}(f') = & \sum_{\psi \in \tilde{\Psi}_{\disc}^{\unr, \tau}(\mathbf{G})} \frac{1}{|\mathcal{S}_{\psi}|} \sum_{\pi_{\infty} \in \Pi(\mathbf{G}_{\R}, \dpsitau \circ \psi_{\infty})} \mathrm{tr}( \pi_{\infty} (f_{\infty}) ) \\
        & \sum_{s \in \mathcal{S}_{\psi},\, s \neq 1} \epsilon_{\psi}(s s_{\psi}) \langle s s_{\psi}, \pi_{\infty} \rangle \prod_p \mathrm{Sat}_{\Gbf_{\Zp}}(f_p)(\dpsitau( c_p(\psi)) ) \nonumber
  \end{align}
  where the last sum is non-empty only if \(\mathcal{S}_{\psi} \neq 1\), so that it makes sense by induction hypothesis.
  By Lemma \ref{lem:pre_refineSOeven} we have
  \[ I^{\mathbf{G}}_{\mathrm{disc}, \tau}(f) = \sum_{\psi \in
      \tilde{\Psi}_{\disc}^{\unr, \tilde{\tau}}(\Gbf)} \sum_{\substack{\pi
    \in \Pi_{\psi_{\infty}, \tau}(\mathbf{G}_{\R}) \\ \langle \cdot ,
    \pi_{\infty} \rangle |_{\mathcal{S}_{\psi}} = \epsilon_{\psi}}}
    \mathrm{tr}( \pi_{\infty} (f_{\infty}) ) \prod_p
    \mathrm{Sat}_{\mathbf{G}_{\Zp}}(f_p)( c_p'(\Gbf, \psi, \pi_{\infty})). \]
  So we obtain
  \begin{align} \label{eq:Sdiscthetarefined0}
    S^{\Gbf}_{\disc, \tau} (f) = & \sum_{\psi \in \tilde{\Psi}_{\disc}^{\unr, \taut}(\mathbf{G})} \frac{1}{|\mathcal{S}_{\psi}|} \sum_{\pi_{\infty} \in \Pi(\mathbf{G}_{\R}, \dpsitau \circ \psi_{\infty})} \mathrm{tr}( \pi_{\infty} (f_{\infty})) \\
    & \left[ \sum_{s \in \mathcal{S}_{\psi}} \epsilon_{\psi}(s s_{\psi}) \langle s s_{\psi}, \pi_{\infty} \rangle \prod_p \Sat_{\Gbf_{\Zp}}(f_p)(c_p'(\Gbf, \psi, \pi_{\infty}))  \right. \nonumber \\
    & \left. - \sum_{s \in \mathcal{S}_{\psi},\, s \neq 1} \epsilon_{\psi}(s s_{\psi}) \langle s s_{\psi}, \pi_{\infty} \rangle \prod_p \Sat_{\Gbf_{\Zp}}(f_p)(\dpsitau(c_p(\psi))) \right] \nonumber
  \end{align}
  where, for convenience, for \(\pi_{\infty}\) such that \(\langle \cdot,
  \pi_{\infty} \rangle |_{\mathcal{S}_{\psi}} \neq \epsilon_{\psi}\) (which
  imposes \(\mathcal{S}_{\psi} \neq 1\)) we let \(c_p'(\Gbf, \psi, \pi_{\infty}) =
  \dpsitau(c_p(\psi))\).
  Note that for such \(\pi_{\infty}\) the difference between square brackets
  simply equals
  \[ \epsilon_{\psi}(s_{\psi}) \langle s_{\psi}, \pi_{\infty} \rangle \prod_p
  \Sat_{\Gbf_{\Zp}}(f_p)(\dpsitau(c_p(\psi))). \]
  To conclude we have to show that:
  \begin{itemize}
    \item For non-endoscopic \(\psi\), in which case \(\dpsitau : \Mpsi \rightarrow
      \Ghat\) is an isomorphism, \(c_p'(\Gbf, \psi, \pi_{\infty})\) does not
      depend on \(\pi_{\infty} \in \Pi_{\dpsitau \circ
      \psi_{\infty}}(\mathbf{G}_{\R})\), and in this case we are forced to define
      \(c_p(\psi)\) as the the preimage by \(\dpsitau\) of the common value.
    \item For endoscopic \(\psi\), for all \(\pi_{\infty}\) we have \(c_p'(\Gbf,
      \psi, \pi_{\infty}) = \dpsitau(c_p(\psi))\).
      By definition this already holds for \(\pi_{\infty}\) such that \(\langle
      \cdot, \pi_{\infty} \rangle |_{\mathcal{S}_{\psi}} \neq \epsilon_{\psi}\).
  \end{itemize}
  Both cases will follow from stability of the distribution \(S_{\disc,
  \tau}^{\Gbf}\), but in the second case we will need an additional piece
  of information (Lemma \ref{lem:enough_elts_AJ_SOeven} below).
  We know (see Remark \ref{rem:mult_one}) that for \(\psi_1 \neq \psi_2\), there
  exists a prime \(p\), which one could even take outside any given finite set,
  such that \(\tilde{c}_p(\psi_1) \cap \tilde{c}_p(\psi_2) = \emptyset\).
  Thus for any \(\psi \in \tilde{\Psi}^{\unr, \tilde{\tau}}_{\disc}(\Gbf)\), the
  summand in \eqref{eq:Sdiscthetarefined0} corresponding to \(\psi\) is a stable
  linear form in \(f \in \Hcal^{\unr}(\Gbf)\).

  \begin{itemize}
    \item For a non-endoscopic \(\psi \in \tilde{\Psi}^{\unr,
      \tilde{\tau}}_{\disc}(\mathbf{G})\), the linear form on
      \(\Hcal^{\unr}(\Gbf)\)
      \begin{equation*} f = \prod_v f_v \mapsto \sum_{\pi_{\infty} \in
        \Pi^{\AJ}(\mathbf{G}_{\R}, \dpsitau \circ \psi_{\infty})} \Tr
        (\pi_{\infty}( f_{\infty} )) \langle s_{\psi}, \pi_{\infty} \rangle
        \prod_p \mathrm{Sat}_{\mathbf{G}_{\Zp}}(f_p)( c_p'(\Gbf, \psi,
        \pi_{\infty}) )
      \end{equation*}
      is stable.
      Recall the discrete parameter \(\varphi_{\tau} : W_{\R} \rightarrow {}^L
      \Gbf\) having infinitesimal character \(\tau\) and the associated discrete
      series L-packet \(\Pi(\Gbf_\R, \varphi_\tau)\).
      By Proposition \ref{pro:tr_AJ_ps}, pseudo-coefficients \(f_{\infty,
      \sigma}\) of elements \(\sigma\) of this L-packet discriminate between
      elements of \(\Pi^{\AJ}(\mathbf{G}_{\R}, \dpsitau \circ \psi_{\infty})\).
      Since the stable orbital integrals of \(f_{\infty, \sigma}\) at regular semisimple elements of \(\Gbf(\R)\) do not depend on the choice of \(\sigma\) in the L-packet\footnote{As recalled in Section \ref{sec:STF} this follows from \cite[Lemma 5.3]{Shelstad_characters}.}, we obtain that \((c_p'(\Gbf, \psi, \pi_{\infty}))_p\) does not depend on \(\pi_{\infty}\).
    \item Now assume that \(\psi\) is an endoscopic parameter in
      \(\tilde{\Psi}^{\unr, \tilde{\tau}}_{\disc}(\mathbf{G})\).
      We can rewrite the contribution of \(\psi\) in \eqref{eq:Sdiscthetarefined0}
      as
      \begin{multline*}
        \frac{\epsilon_{\psi}(s_{\psi})}{|\Scal_{\psi}|} \sum_{\pi_{\infty} \in
          \Pi^{\AJ}(\mathbf{G}_{\R}, \dpsitau \circ \psi_{\infty})} \langle
          s_{\psi}, \pi_{\infty} \rangle \mathrm{tr} (\pi_{\infty}( f_{\infty}))
          \prod_p \mathrm{Sat}_{\mathbf{G}_{\Zp}}(f_p)(\dpsitau(c_p(\psi))) \\
        + \sum_{\substack{\pi_{\infty} \in \Pi^{\AJ}(\mathbf{G}_{\R}, \dpsitau
        \circ \psi_{\infty}) \\ \langle \cdot , \pi_{\infty} \rangle =
        \epsilon_{\psi}}} \mathrm{tr} (\pi_{\infty}( f_{\infty})) \times 
        \left( \prod_p \mathrm{Sat}_{\mathbf{G}_{\Zp}}(f_p)(c_p'(\Gbf, \psi,
          \pi_{\infty})) \right. \\
          \left. - \prod_p
          \mathrm{Sat}_{\mathbf{G}_{\Zp}}(f_p)(\dpsitau(c_p(\psi))) \right)
      \end{multline*}
      and we need to show that the second term (last two lines) vanishes, i.e.\
      that for any \(\pi_{\infty}\) satisfying \(\langle \cdot, \pi_{\infty}
      \rangle = \epsilon_{\psi}\) we have \(c_p'(\Gbf, \psi, \pi_{\infty}) =
      \dpsitau(c_p(\psi))\).
      The above linear form in \(f\) is stable and so is the first term, therefore
      so is the second term.
      As in the previous case we conclude using pseudo-coefficients of discrete
      series, but now we also need the following lemma.
  \end{itemize}
\end{proof}

\begin{lemm} \label{lem:enough_elts_AJ_SOeven}
  For any \(\psi \in \Psit_{\disc}^{\unr, \taut}(\mathbf{G})\), there exists
  \(\pi_{\infty} \in \Pi(\mathbf{G}_{\R}, \dpsitau \circ \psi_{\infty})\) such
  that \(\langle \cdot, \pi_{\infty} \rangle |_{\mathcal{S}_{\psi}} \neq
  \epsilon_{\psi}\).
\end{lemm}
\begin{proof}
  Recall that \(\varphi_\tau\) denotes a discrete Langlands parameter having
  infinitesimal character \(\tau\), and that if we align it (after conjugation by
  \(\Ghat\)) with \(\dpsitau \circ \psi_\infty\) then we have \(C_{\varphi_\tau}
  \subset C_{\dpsitau \circ \psi_\infty}\) and
  \[ \{ \langle \cdot, \sigma \rangle \,|\, \sigma \in \Pi^\AJ(\Gbf_\R, \dpsitau
      \circ \psi_\infty) \} = \{ \langle \cdot, \sigma \rangle|_{C_{\dpsitau
      \circ \psi_\infty}} \,|\, \sigma \in \Pi(\Gbf_\R, \varphi_\tau \}. \]
  Thus it is enough to show that for any \(s \in \mathcal{S}_{\varphi_\tau}
  \smallsetminus \{1\}\), there are discrete series \(\sigma, \sigma'\) such that
  \(\langle s, \sigma \rangle \neq \langle s, \sigma' \rangle\).
  Given the explicit description on p.\ 315 of \cite{Taibi_dimtrace}, it is
  enough to check the following combinatorial fact: if \(I_+\), \(I_-\) are sets of
  cardinality \(n\), and \(\emptyset \subsetneq S \subsetneq I_+ \sqcup I_-\), then
  there are sets \(A_+ \subset I_+\) and \(A_- \subset I_-\) such that \(|A_+| =
  |A_-|\) and \(|S \cap (A_+ \sqcup A_-)|\) is odd.
  This is easy: in fact we can also impose \(|A_+| = 1\).
\end{proof}

\begin{coro}[to Proposition \ref{pro:refineSOeven}] \label{cor:mult_formula_SO}
  Fix \(n \geq 1\) and denote \(\Gbf = \SObf_{4n}\).
  For \(\tau \in \ICcal(\Gbf)\) we have
  \[ \Acal^2(\Gbf(\Q) \backslash \Gbf(\A) / \Gbf(\Zhat)) \simeq
    \bigoplus_{\psi \in \Psit_\disc^{\unr, \taut}(\Gbf)} \Big(
    \chi_f((\dpsitau(c_p(\psi)))_p) \otimes \bigoplus_{\substack{\pi_\infty \in
    \Pi^\AJ(\Gbf_\R, \dpsitau \circ \psi_\infty) \\ \langle \cdot, \pi_\infty
    \rangle|_{\Scal_\psi} = \epsilon_\psi}} \pi_\infty \Big) \]
  where \(\chi_f(c): \Hcal_f^\unr(\Gbf) \to \C\) is the character corresponding to
  the family \(c\) of Satake parameters in \(\Ghat\).
\end{coro}
\begin{proof}
  This simply follows from going through the proof of Proposition
  \ref{pro:refineSOeven} backwards.
\end{proof}

\subsection{Isogenies in level one and stabilization}

In the Langlands-Kottwitz method applied to degree \(n\) Siegel modular varieties,
the reductive groups over \(\Q\) appearing on the ``spectral side'' are
\(\mathbf{GSp}_{2n}\) and some of its endoscopic groups \(\mathbf{G}(
\mathbf{SO}_{4a} \times \mathbf{Sp}_{2n-4a})\), which are only isogenous to
classical groups as considered in \cite{Arthur}.
We need the analogue of the stable multiplicity formula \cite[Theorem
4.1.2]{Arthur} for these isogenous groups.
Bin Xu \cite{Xu} studied the refined local Langlands correspondence for
\(\mathbf{GSp}_{2n}\) and \(\mathbf{GSO}_{2n}\) using stabilization of trace
formulas (and restriction).
He also obtained a stable multiplicity formula in the tempered non-endoscopic
case \cite[Theorem 1.8]{Xu}.
For the purpose of completely describing the cohomology of local systems on
Siegel modular varieties, this is unfortunately not enough.
For local systems corresponding to singular dominant weights, e.g.\ the trivial
local system, it is also necessary to consider non-tempered Arthur-Langlands
parameters (or substitutes thereof).
As Xu explains after \cite[Theorem 1.8]{Xu}, non-tempered packets will certainly
prove more challenging.
This is the main reason why we impose a ``level one'' condition, i.e.\ consider
the moduli stack \(\Acal_n\) instead of an arbitrary Siegel modular variety.

Instead of \cite{Xu} we will give an ad hoc argument using an elementary lifting
result of Chenevier-Renard \cite[Proposition 4.4]{ChRe} relating discrete
automorphic spectra in level one for isogenous split semisimple groups over
\(\Q\).
This result is conjecturally related (but is somewhat more precise) to
properties of the quotient \(L_{\Z}\) of the hypothetical Langlands group \(L_{\Q}\)
of \(\Q\), namely its connectedness (Minkowski's theorem on unramified extensions
of \(\Q\)) and simple connectedness (see \cite[Appendix B]{ChRe}).
In Proposition \ref{pro:csqChRe} we will deduce from \cite[Proposition
4.4]{ChRe} and the relation between Adams-Johnson packets and discrete series
packets an expansion, involving certain families of lifts of Satake parameters,
for the spectral side of the trace formula for \emph{adjoint} groups, restricted
to pseudo-coefficients of discrete series representations at the real place.
Then we ``stabilize'' this lifting result using arguments similar (but more
intricate) to the arguments of the previous section.
Roughly, Proposition \ref{pro:lift_stab_mult} says that the above lifted
families of Satake parameters only depend on the (formal) Arthur-Langlands
parameter, and are compatible with endoscopy.

First we recall a nice property of restriction of \((\gfrak, K)\)-modules under
restriction via an isogeny.

\begin{prop} \label{pro:res_real_iso}
  Let \(\Gbf \rightarrow \Gbf'\) be a morphism of connected reductive groups over
  \(\R\) having central kernel and such that its image contains \(\Gbf'_\der\).
  Let \(K\) be a maximal compact subgroup of \(\Gbf(\R)\), and let \(K'\) be the
  maximal compact subgroup of \(\Gbf'(\R)\) containing the image of \(K\).
  Denote \(\gfrak = \C \otimes_{\R} \Lie( \Gbf)\) and \(\gfrak' = \C \otimes_{\R}
  \Lie ( \Gbf')\).
  Let \(C = \ker(\Gbf(\R) \to \Gbf'(\R))\).
  \begin{enumerate}
    \item
      The restriction of an irreducible \((\gfrak', K')\)-module to \((\gfrak, K)\)
      is semisimple, has finite length and is multiplicity-free.
    \item
      If \(\pi\) is an irreducible \((\gfrak, K)\)-module then there exists an
      irreducible \((\gfrak', K')\)-module \(\pi'\) such that \(\pi\) occurs in the
      restriction of \(\pi'\) if and only if the restriction of the central
      character of \(\pi\) to \(C\) is trivial.
      Moreover the kernel of \(\Hom(\Gbf'(\R), \C^\times) \to \Hom(\Gbf(\R),
      \C^\times)\) acts transitively on the set of isomorphism classes of such
      \(\pi'\), and if \(\pi\) is essentially discrete series and \(\Gbf'\) is
      semisimple then this set has at most one element.
    \item 
      To simplify the formulation assume that \(\Gbf'\) is quasisplit, and fix a
      Whittaker datum \(\mathfrak{w}\) for \(\Gbf\), normalising Shelstad's
      parametrization of L-packets \cite{Shelstad_temp2} for both \(\Gbf\) and
      \(\Gbf'\).
      Let \(\varphi' : W_{\R} \rightarrow {}^L \Gbf'\) be an essentially discrete
      parameter.
      Let \(\varphi : W_{\R} \rightarrow {}^L \Gbf\) be the essentially discrete
      parameter obtained by composing with \({}^L \Gbf' \rightarrow {}^L \Gbf\).
      We have an injection \(\Scal_{\varphi'} \subset \Scal_{\varphi}\) and for
      any \(\pi' \in \Pi(\Gbf', \varphi')\),
      \[ \pi'|_{\Gbf(\R)} \simeq \bigoplus_{\substack{\pi \in \Pi(\Gbf, \varphi)
           \\
         \langle \cdot, \pi \rangle|_{\mathcal{S}_{\varphi'}} = \langle \cdot,
         \pi' \rangle}} \pi. \]
  \end{enumerate}
\end{prop}
\begin{proof}
  Let \(C_K = C \cap K\).
  Recall that \(K'\) is the group of \(g \in \Gbf'(\R)\) normalising \(K/C_K\) such
  that for any (rational) character \(\chi: \Gbf' \to \GLbf_1\) we have \(\chi(g)^2
  = 1\).
  In particular \(K/C_K\) is a distinguished open subgroup of \(K'\), and it is
  well-known that the quotient \(Q\) of \(K'\) by \((K' \cap Z(\Gbf'(\R)))(K/C_K)\) is
  naturally isomorphic to the finite abelian \(2\)-torsion group
  \(\Gbf'(\R)/(Z(\Gbf'(\R)) \, \Gbf(\R)/C)\).
  Now if \(\pi'\) is an irreducible \((\gfrak', K')\)-module, its restriction to
  \((\gfrak, K)\) is finitely generated since \(Q\) is finite, and so it admits an
  irreducible quotient \(\pi'|_{\gfrak, K} \twoheadrightarrow \pi\).
  Denoting by \(W\) the kernel of this map, the subspace \(\cap_{k \in Q} kW\) is a
  proper submodule of the irreducible \(\pi'\), so \(\pi'|_{\gfrak, K}\) embeds in
  \(\bigoplus_{k \in Q} \pi^k\), which is clearly semisimple of finite length.

  We are left to prove the multiplicity one statement.
  We use the same reduction to the case of discrete series representations as in
  \cite{Langlands}.
  \begin{enumerate}
    \item
      If \(\pi'\) is an essentially discrete series representation then the result
      is implicit in \cite[\S 3]{Langlands} which deduces the classification of
      discrete series representations from the case of connected semisimple Lie
      groups considered by Harish-Chandra.
      Since we will need the details later, let us briefly give the argument.
      Choose a  maximal torus \(\Tbf\) in \(\Gbf\) which such that \(\Tbf(\R)\)
      contains a maximal torus of \(K\), and let \(\Tbf'\) be the unique maximal
      torus of \(\Gbf'\) containing the image of \(\Tbf\).
      In particular they are both stable under the Cartan involutions defined by
      \(K\) and \(K'\).
      Let \(W := N_{\Gbf(\C)}(\Tbf) / \Tbf(\C) = N_{\Gbf'(\C)}(\Tbf') /
      \Tbf'(\C)\), \(W_c := N_{\Gbf(\R)}(\Tbf) / \Tbf(\R) = N_K(\Tbf) / (\Tbf(\R)
      \cap K)\) and \(W_c' := N_{\Gbf'(\R)}(\Tbf') / \Tbf'(\R) = N_{K'}(\Tbf') /
      (\Tbf'(\R) \cap K')\).
      Then \(W_c\) is a normal subgroup of \(W_c'\), which is a subgroup of \(W\).
      If we fix a Borel subgroup \(\Bbf'\) of \(\Gbf'_\C\) containing \(\Tbf'_\C\)
      then essentially discrete series \((\gfrak', K')\)-modules are parametrized
      by infinitesimal character, central character (with a compatibility
      condition; these two parameters correspond to the Langlands parameter) and
      an element of \(W'_c \backslash W\), which can be read on Harish-Chandra's
      formula for the restriction of the trace character to \(\Tbf'(\R)\).
      (The essential points in Langlands' argument to classify essentially
      discrete series are the fact that the obviously injective map \(W_c'/W_c
      \to Q\) is also surjective, which simply follows from the fact that maximal
      tori in \(K\) are all conjugated under \(K^0\), and the well-known fact that
      \(\Gbf(\R)\) is connected if \(\Gbf\) is semisimple and simply connected.)
      From this description it is easy to observe that fixing \(\Bbf'\) as above
      and given an essentially discrete Langlands parameter \(\varphi': W_\R \to
      {}^L \Gbf'\), the restriction of the \((\gfrak', K')\)-module corresponding
      to the coset \(W_c' w' \in W_c' \backslash W\) to \((\gfrak, K)\) is the
      direct sum of the essentially discrete modules corresponding to the
      Langlands parameter \(\varphi: W_\R \to {}^L \Gbf\), obtained by composing
      \(\varphi'\) with \({}^L \Gbf' \to {}^L \Gbf\), and all \(|W_c'/W_c|\) cosets
      \(W_c w \in W_c \backslash W\) mapping to \(W_c' w'\).

      The case of essentially tempered modules follows, using the following
      classification results.
      Let \(\pi'\) be an essentially tempered irreducible \((\gfrak', K')\)-module.
      Up to twisting by a character of \(\Gbf'(\R)\) we may assume that its
      central character is unitary.
      By \cite{Trombi_tempered} (see also \cite[Lemma 4.10]{Langlands} or
      \cite[Proposition 5.2.5]{WallachReal}) there exists a parabolic subgroup
      \(\Pbf'\) of \(\Gbf'\) such that, letting \(\Mbf'\) be the unique Levi factor
      which is stable under the Cartan involution corresponding to \(K'\) and
      denoting \(K'_{\Mbf'} = K' \cap \Mbf'(\R)\), there exists an irreducible
      essentially discrete series \((\mfrak', K'_{\Mbf'})\)-module with unitary
      central character \(\sigma'\) such that \(\pi'\) embeds in the parabolically
      induced \((\gfrak', K')\)-module \(\Ind_{\Pbf'}^{\Gbf'} \sigma'\) (which is
      semisimple by unitarity, and known to have finite length).
      Let \(\Pbf\) (resp.\ \(\Mbf\)) be the preimage of \(\Pbf'\) (resp.\ \(\Mbf'\)) in
      \(\Gbf\) and \(K_{\Mbf} = K \cap \Mbf(\R)\).
      From the previous case we know that the restriction of \(\sigma'\) to
      \((\mfrak', K'_{\Mbf'})\) is isomorphic to the direct sum of irreducible
      non-isomorphic \((\mfrak, K_\Mbf)\)-modules \(\sigma_1, \dots, \sigma_r\).
      Thus the restriction of \(\Ind_{\Pbf'}^{\Gbf'} \sigma'\) to \((\gfrak, K)\) is
      isomorphic to \(\bigoplus_i \Ind_{\Pbf}^{\Gbf} \sigma_i\).
      As observed in \cite{Langlands} it is implicit in Harish-Chandra's work
      that for \(j \neq i\), \(\Ind_{\Pbf}^{\Gbf} \sigma_i\) and \(\Ind_{\Pbf}^{\Gbf}
      \sigma_j\) have no constituent in common.
      More precisely, if they had a consitutent in common then any isomorphism
      between the two constituents would extend to an isomorphism between the
      associated irreducible unitary representations of \(\Gbf(\R)\), which are
      constituents of the \(L^2\) induced unitary representations
      \(\Ind_{\Pbf}^{\Gbf} \widehat{\sigma_i}\) and \(\Ind_{\Pbf}^{\Gbf}
      \widehat{\sigma_j}\), where \(\widehat{\sigma_i}\) is the irreducible
      unitary representation of \(\Mbf(\R)\) obtained by completing \(\sigma_i\) for
      the essentially unique Hermitian inner product making it unitary.
      But the Plancherel formula proved in \cite{HC_harmonic_real3} allows one
      to construct, for any irreducible representation \(\delta\) of \(K\), a
      bi-\(K\)-finite Schwartz function on \(\Gbf(\R)\) acting by \(0\) on
      \(\Ind_{\Pbf}^{\Gbf} \widehat{\sigma_i}\) and as the projector on the
      \(\delta\)-isotypic component on \(\Ind_{\Pbf}^{\Gbf} \widehat{\sigma_j}\)
      (using Corollary 26.1 and Lemma 26.1 loc.\ cit., the argument being
      similar to the one in \S 37 loc.\ cit.).
      Finally, thanks to \cite{KnappComm} we know that each \(\Ind_{\Pbf}^{\Gbf} \sigma_i\) is multiplicity-free, and so the restriction of \(\Ind_{\Pbf'}^{\Gbf'} \sigma'\) to \((\gfrak,K)\) is also multiplicity-free.
      %TODO: point out somewhere that by \cite{Lepowsky_alg}, every irreducible \((\gfrak, K)\)-module is associated to a representation of \(\Gbf(\R)\) on a Hilbert space, and in the case of unitary representations or modules see \cite[Theorem 3.4.11]{WallachReal} and \cite[Theorem 4.4.6.6]{Warner_book1}.
      %TODO: see \cite[\S 3]{Casselman_can_ext} for the fact that matrix coefficients only depend on the \((\gfrak, K)\)-module.
      %TODO: for the relation between parabolic induction for representations and modules, see \cite[Proposition 6.3.5]{Vogan_book}

      Finally the case of arbitrary irreducible \((\gfrak', K')\)-modules follows
      from the essentially tempered case and Langlands' classification (Lemmas
      3.14 and 4.2 in \cite{Langlands}).
      More precisely, \(\pi'\) is the unique irreducible quotient of
      \(\Ind_{\Pbf'}^{\Gbf'} \sigma'\) for an irreducible essentially tempered
      \((\mfrak', K'_{\Mbf'})\)-module \(\sigma'\) whose central character satisfies
      a certain positivity condition with respect to \(\Pbf'\).
      By the previous case we have \(\sigma'|_{\mfrak, K_\Mbf} \simeq
      \bigoplus_{i=1}^r \sigma_i\) with \(\sigma_i \not\simeq \sigma_j\) if \(i \neq
      j\), and so we have an isomorphism \((\Ind_{\Pbf'}^{\Gbf'}
      \sigma')|_{\mfrak, K_\Mbf} \simeq \bigoplus_{i=1}^r \Ind_{\Pbf}^{\Gbf}
      \sigma_i\) (essentially because the natural map \(\Pbf(\R) \backslash
      \Gbf(\R) \to \Pbf'(\R) \backslash \Gbf'(\R)\) is a homeomorphism) and thus
      a surjective morphism \(\bigoplus_{i=1}^r \Ind_{\Pbf}^{\Gbf} \sigma_i
      \twoheadrightarrow \pi'|_{\gfrak, K}\).
      Since the right-hand side is semisimple, this map factors through the
      Langlands quotients of all \(\Ind_{\Pbf}^{\Gbf} \sigma_i\), which are
      irreducible and non-isomorphic \((\gfrak, K)\)-modules.
      (In fact we see that \(\pi'|_{\gfrak, K}\) is the direct sum of these
      Langlands quotients since \(K'_{\Mbf'}\) acts transitively on the set of
      \(\sigma_i\)'s.)

    \item
      The ``only if'' condition is obvious.

      First assume that \(\Gbf \to \Gbf'\) is surjective.
      The the map \(\gfrak \to \gfrak'\) is also surjective, and \(K'/C_K\) is an
      open, and thus finite index, subgroup of \(K'\).
      In particular if the central character of \(\pi\) is trivial on \(C\) then the
      existence of \(\pi'\) is easy: at least one irreducible factor of the
      induction of \(\pi\) to a \((\gfrak', K')\)-module works.
      Now fix such a \((\gfrak', K')\)-module \(\pi'\), and let \(Q''\) be the
      stabilizer of the isomorphism class of \(\pi\) under the action of \(Q\) (by
      conjugation).
      Let \(K''\) be the preimage of \(Q''\) in \(K'\).
      By the previous point we know that the restriction of \(\pi'\) to \((\gfrak',
      K'')\) is multiplicity-free, and this shows that there is a unique factor
      \(\pi''\) of this restriction whose restriction to \((\gfrak, K)\) is
      isomorphic to \(\pi\).
      We also see that \(\pi'\) is isomorphic to the induction of \(\pi''\) to a
      \((\gfrak', K')\)-module.
      If \(\Gbf'\) is semisimple and admits discrete series then \(K' \cap
      Z(\Gbf'(\R))\) is contained in \(K/C_K\) and so \(Q = K'/(K/C_K)\); if moreover
      \(\pi\) is essentially discrete series then we see from the previous point
      that \(Q''\) is trivial and so \(\pi'\) is unique up to isomorphism.
      In general (but still assuming that \(\Gbf \to \Gbf'\) is surjective) if
      \(\pi'_\flat\) is another \((\gfrak', K')\)-module such that \(\pi\) occurs in
      its restriction, then the restriction of \(\pi'_\flat\) to \((\gfrak', K'')\)
      contains a unique irreducible \(\pi''_\flat\) in which \(\pi'\) occurs, and by
      irreducibility of \(\pi\) there exists a unique continuous character \(\chi\)
      of \(K''/(K/C_K)\) such that \(\pi''_\flat\) is isomorphic to \(\pi'' \otimes
      \chi\).
      Let \(\tilde{\chi}\) be a character of \(K'\) extending \(\chi\).
      Such a character exists because \(K'/(K/C_K)\) is a compact abelian group
      and \(K''/(K/C_K)\) is an open subgroup.
      Then \(\pi'_\flat\) is isomorphic to \(\pi' \otimes \tilde{\chi}\).

      Without assuming that \(\Gbf \to \Gbf'\) is surjective, we can apply the
      above to the surjective morphism \(\Gbf \times \Zbf(\Gbf')^0 \to \Gbf'\).
      The details of this reduction are left to the reader.

    \item
      It only remains to reformulate part of the proof of the first point in
      terms of dual groups.
      Let \(\Bbf\) be the Borel subgroup of \(\Gbf_{\C}\) containing \(\Tbf_{\C}\)
      corresponding to the choice of \(\mathfrak{w}\) (\(\Bbf\) has the property
      that all simple roots of \(\Tbf\) in \(\Bbf\) are non-compact, and is
      well-defined up to conjugation by \(K\)).
      Denote \((\Bcal, \Tcal)\) (resp.\ \((\Bcal', \Tcal')\)) the Borel pair which
      is part of the pinning used to form \({}^L \Gbf\) (resp.\ \({}^L \Gbf'\)).
      Using \(\Bbf\) (and its image in \(\Gbf'\)) we obtain compatible
      identifications \(\Tcal \simeq \widehat{\Tbf}\) and \(\Tcal' \simeq
      \widehat{\Tbf'}\).
      Up to conjugating by \(\widehat{\Gbf'}\) we can assume that \(\varphi'\) is in
      diagonal and dominant position with respect to the Borel pair \((\Bcal',
      \Tcal')\) in \(\widehat{\Gbf'}\).
      We have a commutative diagram
      \begin{equation*}
        \begin{tikzcd}
          W_c \backslash W \arrow[r, "{\sim}"] \arrow[d] & \ker(H^1(\R, \Tbf)
            \rightarrow H^1(\R, \Gbf)) \arrow[d] \arrow[r, hook] &
            \Scal_{\varphi}^{\wedge} \arrow[d] \\
          W_c' \backslash W \arrow[r, "{\sim}"] & \ker(H^1(\R, \Tbf')
            \rightarrow H^1(\R, \Gbf')) \arrow[r, hook] &
            \Scal_{\varphi'}^{\wedge}
        \end{tikzcd}
      \end{equation*}
      where the left horizontal bijections are given by \(\mathrm{cl}(g) \mapsto
      (\sigma \mapsto g^{-1} \sigma(g))\) (see \cite[Theorem
      2.1]{Shelstad_characters}) and the right horizontal maps are obtained by
      Tate-Nakayama duality (see \cite{KottEllSing}).
      The two horizontal compositions define the maps \(\Pi(\Gbf, \varphi) \to
      \Scal_{\varphi}^{\wedge}\) and \(\Pi(\Gbf', \varphi') \to
      \Scal_{\varphi'}^{\wedge}\).
  \end{enumerate}
\end{proof}

Note that the third point could easily be generalized to non-quasisplit groups
using \cite[\S 5.4]{Kalri}, and one could certainly prove a similar restriction
formula for arbitrary irreducible \((\gfrak', K')\)-modules by following
Shelstad's and Kaletha's arguments closely.

Also observe that the uniqueness property in the second point is particular to
the discrete series: for example the trivial representation of \(\Gbf(\R)\), which
is an Adams-Johnson representation, can be extended into more than one
character of \(\Gbf'(\R)\).

We will nonetheless apply this indirectly to Adams-Johnson representations.
Recall that standard modules form a basis of the Grothendieck group \(K_0(\gfrak, K)\) of finite length \((\gfrak, K)\)-modules (\cite[Corollary 5.5.3]{WallachReal}).
Using this basis one can project to the subgroup (freely) generated by
irreducible \((\gfrak, K)\)-modules in the discrete series.
This projection that we denote \(\prDS\) is also computed by taking the trace on
all pseudo-coefficients for discrete series.
The proof of Proposition \ref{pro:res_real_iso} immediately implies the
following lemma.

\begin{lemm} \label{lem:prDS_isog}
  Let \(\Gbf \to \Gbf'\) be a morphism between connected reductive groups over
  \(\R\) as in Proposition \ref{pro:res_real_iso} (i.e.\ an isogeny).
  Then restriction of \((\gfrak', K')\)-modules to \((\gfrak, K)\) intertwines the
  maps \(\prDS\) for \(\Gbf\) and \(\Gbf'\).
\end{lemm}

\begin{lemm} \label{lem:parity_orth}
  If \(\Gbf = \mathbf{Sp}_{2n}\) (resp.\ \(\mathbf{SO}_{4n}\)), \(\tau \in \ICcal(\Gbf)\) and \(\tilde{\Psi}^{\unr, \tilde{\tau}}_{\disc}(\Gbf) \neq \emptyset\), then \(\tau\) belongs to \(\ICcal(\Gbf_\mathrm{ad})\), i.e.\ it is C-algebraic for \(\Gbf_\mathrm{ad}\).
  Concretely, using notation as in Definition \ref{def:ICcal}, this means that writing \(\tau\) as the class of \(w_1 > \dots > w_n > 0\) (resp.\ \((w_1 > \dots > w_{2n} > 0)\)), the integer \(\sum_{i=1}^n w_i - n(n+1)/2\) (resp.\ \(\sum_{i=1}^{2n} w_i - n\)) is even.
\end{lemm}
\begin{proof}
  See \cite[Remark 4.1.6]{Taibi_dimtrace} and \cite[Proposition 1.8]{ChRe} (and its proof p.\ 42 loc.\ cit.).
\end{proof}

\begin{prodef} \label{def:Mpsisc}
  Let \(n \geq 1\), \(\Gbf = \mathbf{Sp}_{2n}\) or \(\mathbf{SO}_{4n}\), \(\tau \in \ICcal(\Gbf)\).
  Let \(\psi \in \tilde{\Psi}^{\unr, \tilde{\tau}}_{\disc}(\Gbf)\).
  Let \(\Mpsisc\) be the simply connected cover of \(\Mpsi\) (Definition \ref{def:Mpsi}).

  \begin{enumerate}
    \item
      Let \(\dpsitau: \Mpsi \to \Ghat\) be as in Definition \ref{def:Mpsi}.
      There exists a unique \(\dpsitausc : \Mpsisc \rightarrow
      \Ghat_{\mathrm{sc}}\) lifting \(\dpsitau\).
      Thus, like \(\dpsitau\), the lift \(\dpsitausc\) is well-defined up to
      conjugation by \(\Ghat\).
    \item
      For \(\psi_\infty\) as in Definition \ref{def:Mpsi}, there exists
      \(\psi_{\infty, \mathrm{sc}} : W_{\R} \times \SL_2 \rightarrow \Mpsisc\)
      lifting \(\psi_{\infty}\), unique up to multiplication by an element of
      \(Z^1(W_{\R}, \ker( \Mpsisc \rightarrow \Mpsi) )\).
    \item
      The centralizer \(C_{\dpsitausc \circ \psi_{\infty, \mathrm{sc}}}\) of
      \(\dpsitausc \circ \psi_{\infty, \mathrm{sc}}\) in \(\Ghat_{\mathrm{sc}}\) is
      abelian and stays unchanged if we multiply \(\psi_{\infty, \mathrm{sc}}\) by
      an element of \(Z^1(W_{\R}, \ker( \Mpsisc \rightarrow \Mpsi))\).
      In particular the limits (over the choices described above)
      \(C_{\psi_\infty, \mathrm{sc}} := \lim C_{\dpsitausc \circ \psi_{\infty,
      \mathrm{sc}}}\) and \(C_{\psi, \mathrm{sc}} := \lim C_{\dpsitausc}\) are
      naturally isomorphic to any one of their terms, and \(C_{\psi, \mathrm{sc}}
      \subset C_{\psi_\infty, \mathrm{sc}}\).
    \item
      The centralizer \(C_{\dpsitausc}\) of \(\dpsitausc\) in \(\Ghat_{\mathrm{sc}}\)
      is equal to \(\dpsitausc(Z(\Mpsisc))\).
      For this reason we often simply denote it \(C_{\psi, \mathrm{sc}}\).
      In particular the (obviously injective) map \(C_{\psi, \mathrm{sc}} /
      Z(\Ghat_{\mathrm{sc}}) \to \Scal_{\psi}\) is surjective.
  \end{enumerate}
\end{prodef}
\begin{proof}
  \begin{enumerate}
  \item This follows from \cite[Proposition 2.24 (i)]{BorelTits_compl}.
  \item The discussion around Lemma 4.2.2 in \cite{Taibi_dimtrace} completely describes Adams-Johnson parameters for an arbitrary connected semisimple group \(\Hbf\) over \(\R\) admitting discrete series representations.
    Let \((\Bcal_{\Hhat}, \Tcal_{\Hhat})\) be the Borel pair of \(\Hhat\) used to define \({}^L \Hbf\).
    For any C-algebraic regular infinitesimal character, represented by a strictly dominant (for \(\Bcal_{\Hhat}\)) \(\tau \in \rho^\vee_{\Hhat} + X_*(\Tcal_{\Hhat})\), and any standard parabolic subgroup \(\Qcal = \Lcal \Ncal\) of \(\Hhat\) such that the opposite parabolic is conjugated to \(\Qcal\) by an element of \(\Hhat \rtimes j\) and \(\langle \tau, \alpha \rangle = 1\) for any simple root of \(\Tcal\) in the Levi factor \(\Lcal\), there exist Adams-Johnson parameters having this infinitesimal character and associated parabolic subgroup of \(\Hhat\), and the set of such parameters is a torsor under \(Z(\Lcal)/\{t^2 \,|\, t \in Z(\Lcal)\}\).
    Moreover all Adams-Johnson parameters for \(\Hbf\) are obtained in this way.
    Thus an Adams-Johnson parameter admits a lift along the dual of a central isogeny \(\Hbf \to \Hbf'\) if and only if we have \(\tau-\rho^\vee_{\Hhat}\) belongs to the finite index subgroup \(X_*(\Tcal_{\widehat{\Hbf'}})\) of \(X_*(\Tcal_{\Hhat})\).
    The statement thus follows from Lemma \ref{lem:parity_orth}.
  \item This follows from a general property of Adams-Johnson parameters: their centralizer is contained in the \(2\)-torsion of a maximal torus.
  \item The centralizer \(C_{\dpsitausc}\) in \(\Ghat_{\mathrm{sc}}\) is obviously contained in the preimage of the centralizer \(C_{\dpsitau}\) in \(\Ghat\), so this follows from the equality \(\dpsitau(Z(\Mpsi)) = C_{\dpsitau}\) (see Definition \ref{def:Mpsi}).
  \end{enumerate}
\end{proof}

\begin{prop} \label{pro:csqChRe}
  Let \(n \geq 1\), \(\Gbf = \mathbf{Sp}_{2n}\) or \(\mathbf{SO}_{4n}\), and \(\tau \in \ICcal(\Gbf_\ad)\).
  Let \(\delta' \in \Pi(\Gbf_{\mathrm{ad}, \R}, \varphi'_{\tau})\) and let \(f_{\infty}\) be a pseudo-coefficient for \(\delta'\).
  Let \(\prod_p f_p \in \Hcal^{\unr}_f(\Gbf_{\mathrm{ad}})_\C\), and \(f = f_{\infty} \prod_p f_p\).
  Then we have
  \[ I^{\Gbf_{\mathrm{ad}}}_{\disc, \tau}(f) = \sum_{\substack{\psi \in
        \tilde{\Psi}^{\unr, \tilde{\tau}}_{\disc}(\Gbf) \\ \langle \cdot,
      \delta' \rangle|_{\Scal_{\psi}} = \epsilon_{\psi}}}
      \epsilon_{\psi}(s_{\psi}) (-1)^{q(\Lbf^*_{\psi,\tau})} \prod_p
      \Sat_{\Gbf_{\mathrm{ad}, \Zp}}(f_p)(\cpsc'(\Gbf, \psi, \delta')) \]
  where \(\Lbf^*_{\psi,\tau}\) is the quasi-split twisted Levi subgroup of \(\Gbf\) associated to
  the Adams-Johnson parameter \(\dpsitau \circ \psi_{\infty}\), and \(\cpsc'(\Gbf,
  \psi, \delta')\) is a uniquely determined semisimple conjugacy class in
  \(\Ghat_{\mathrm{sc}}(\C)\) lifting the semisimple conjugacy class
  \(\dpsitau(c_p(\psi))\) in \(\Ghat(\C)\).
  (Note that the condition on \(\langle \cdot, \delta' \rangle|_{\Scal_\psi}\) in
  the sum makes sense thanks to 4.\ in Proposition-Definition \ref{def:Mpsisc}.)
\end{prop}
\begin{proof}
  Fix a maximal compact subgroup \(K_\infty\) of \(\Gbf(\R)\) and denote \(K_\infty'\)
  the corresponding maximal compact subgroup of \(\Gbf_\ad(\R)\).
  The pseudo-coefficient \(f_{\infty}\) selects \((\gfrak, K_{\infty}')\)-modules
  having infinitesimal character \(\tau\).
  Arthur's multiplicity formula for \(\Gbf\) in level one and infinitesimal
  character \(\tau\) (as refined in Corollary \ref{cor:mult_formula_SO} in case
  \(\Gbf = \mathbf{SO}_{4n}\)) reads
  \[ \Acal^2(\Gbf(\Q) \backslash \Gbf(\A) / \Gbf(\Zhat))_\tau \simeq
    \bigoplus_{\psi \in \Psit_{\disc}^{\unr, \taut}(\Gbf)}
    \bigoplus_{\pi_{\infty} \in X(\psi)} \pi_{\infty} \otimes
    \chi_f \left( (\dpsitau(c_p(\psi)))_p \right) \]
  where \(X(\psi)\) is the subset of \(\Pi_{\dpsitau \circ \psi_{\infty}}
  (\Gbf_{\R})\) consisting of those \(\pi_{\infty}\) such that \(\langle \cdot,
  \pi_{\infty} \rangle|_{\mathcal{S}_{\psi}} = \epsilon_{\psi}\).
  Reformulating \cite[Prop.\ 4.4]{ChRe}, and using Remark \ref{rem:mult_one}, we obtain that for any \(\psi \in \Psit_{\disc}^{\unr, \taut}(\Gbf)\) there are uniquely determined objects as follows:
  \begin{enumerate}
    \item a set \(X'(\psi)\) of unitary irreducible admissible \((\gfrak,
      K_{\infty}')\)-modules and a surjective map \(X(\psi) \rightarrow
      X'(\psi)\) such that for any \(\pi_{\infty}' \in X'(\psi)\),
      \[ \pi_{\infty}'|_{(\gfrak, K_\infty)} \simeq
         \bigoplus_{\substack{\pi_{\infty} \in X(\psi) \\ \pi_{\infty} \mapsto
         \pi_{\infty}'}} \pi_{\infty}, \]
    \item for each \(\pi_{\infty}' \in X'(\psi)\), semi-simple conjugacy classes
      \(\cpsc'(\Gbf, \psi, \pi_{\infty}')\) in \(\Ghat_{\sico}\) lifting
      \(\dpsitau(c_p(\psi))\),
  \end{enumerate}
  such that
  \[ \Acal^2(\Gbf_\ad(\Q) \backslash \Gbf_\ad(\A) / \Gbf_\ad(\Zhat))_\tau \simeq
    \bigoplus_{\psi \in \Psit_{\disc}^{\unr, \taut}(\Gbf)}
    \bigoplus_{\pi_{\infty}' \in X'(\psi)} \pi_{\infty}' \otimes \chi_f \left(
    (\cpsc'(\Gbf, \psi, \pi_{\infty}'))_p \right). \]
  We now compare \(\mathrm{pr}_{\mathrm{DS}}\) on \(X(\psi)\) and \(X'(\psi)\).
  Let \(\varphi_\tau: W_\R \to \Ghat\) be a discrete Langlands parameter having
  infinitesimal character \(\tau\) and aligned (see Proposition \ref{pro:tr_AJ_ps}) with \(\dpsitau \circ \psi_\infty\).
  Since \(\tau \in \ICcal(\Gbf_\ad)\) there is a unique lift \(\varphi_\tau': W_\R
  \to \Ghat_{\mathrm{sc}}\).
  Recall from Proposition \ref{pro:tr_AJ_ps} that the Adams-Johnson packet
  attached to \(\dpsitau \circ \psi_{\infty}\) occasions a partition of
  \(\Pi(\Gbf_{\R}, \varphi_{\tau})\):
  \[ \Pi(\Gbf_{\R}, \varphi_{\tau}) = \bigsqcup_{\pi_{\infty} \in
    \Pi^\AJ(\Gbf_{\R}, \dpsitau \circ \psi_{\infty})} Y(\dpsitau \circ
    \psi_{\infty}, \pi_{\infty}) \]
  such that whenever \(\delta \in Y(\dpsitau \circ \psi_{\infty}, \pi_{\infty})\)
  we have \(\langle \cdot, \delta \rangle = \langle \cdot, \pi_{\infty}
  \rangle|_{\mathcal{S}_{\dpsitau \circ \psi_{\infty}}}\), and for any
  \(\pi_{\infty} \in \Pi^\AJ(\Gbf_{\R}, \dpsitau \circ \psi_{\infty})\) we have,
  in the Grothendieck group \(K_0(\gfrak, K_\infty)\),
  \begin{equation} \label{eq: proj AJ on DS}
    \prDS(\pi_{\infty}) = \langle s_{\psi}, \pi_\infty \rangle
    (-1)^{q(\Lbf^*_{\psi,\tau})} \sum_{\delta \in Y(\dpsitau \circ \psi_{\infty},
    \pi_{\infty})} \delta.
  \end{equation}
  In particular the supports of \(\prDS(\pi_{\infty})\) (as \(\pi_{\infty}\) varies
  in \(\Pi^\AJ(\Gbf_{\R}, \dpsitau \circ \psi_{\infty})\)) are disjoint, and the
  sign in \eqref{eq: proj AJ on DS} is the same for all elements of \(X(\psi)\).
  By Proposition \ref{pro:res_real_iso} (uniqueness in 2.), Lemma
  \ref{lem:prDS_isog} and Remark \ref{rem:mult_one} this implies that the
  supports of \(\prDS(\pi_{\infty}')\), as \(\pi_{\infty}'\) varies in \(X'(\psi')\),
  are also disjoint so that we have a partition of a subset of \(\Pi(\Gbf_{\ad,
  \R}, \varphi'_{\tau})\):
  \[ \bigsqcup_{\pi_{\infty}' \in X'(\psi)} Y'(\psi, \pi_{\infty}') \subset
  \Pi(\Gbf_{\ad, \R}, \varphi'_{\tau}), \]
  determined by
  \[ \bigoplus_{\delta' \in Y'(\psi, \pi_{\infty}')} \delta'|_{(\gfrak,
     K_\infty)} \simeq \bigoplus_{\substack{\pi_{\infty} \in X(\psi) \\
     \pi_{\infty} \mapsto \pi_{\infty}'}} \ \bigoplus_{\delta \in Y(\dpsitau
     \circ \psi_{\infty}, \pi_{\infty})} \delta. \]
  Since
  \[ \bigsqcup_{\pi_{\infty} \in X(\psi)} Y(\dpsitau \circ \psi_{\infty},
    \pi_{\infty}) = \left\{ \delta \in \Pi(\Gbf_{\R}, \varphi_{\tau})
    \,\middle|\, \langle \cdot, \delta \rangle|_{\mathcal{S}_{\psi}} =
    \epsilon_{\psi} \right\} \]
  we also have, by 3.\ in Proposition \ref{pro:res_real_iso},
  \[ \bigsqcup_{\pi_{\infty}' \in X'(\psi)} Y(\psi, \pi_{\infty}') = \left\{
    \delta' \in \Pi(\Gbf_{\ad, \R}, \varphi'_{\tau}) \,\middle|\, \langle \cdot,
    \delta' \rangle|_{\mathcal{S}_{\psi}} = \epsilon_{\psi} \right\}. \]
  By \eqref{eq: proj AJ on DS} and the definition of \(X(\psi)\), for any
  \(\pi_{\infty}' \in X'(\psi)\) we have
  \[ \prDS(\pi_{\infty}') = \epsilon_{\psi}(s_{\psi}) (-1)^{q(\Lbf^*_{\psi,\tau})}
    \sum_{\delta' \in Y'(\psi, \pi_{\infty}')} \delta'. \]
  We conclude by letting \(\cpsc'(\Gbf, \psi, \delta') = \cpsc'(\Gbf,
  \psi, \pi_{\infty}')\) for \(\delta' \in Y'(\psi, \pi_{\infty}')\).
\end{proof}

We now show that these lifted Satake parameters \(\cpsc'(\Gbf, \psi, \delta_{\infty})\) do not depend on \(\pi_{\infty}\) and are compatible with endoscopy.

\begin{comment}
\begin{lemm} \label{le:welldefined}
  \begin{enumerate}
    \item
      For \(a \geq 2, b \geq 0\), consider the natural morphism of complex
      semisimple groups \(f : \mathrm{Spin}_{2a} \times \mathrm{Spin}_{2b+1}
      \rightarrow \mathrm{Spin}_{2(a+b)+1}\).
      Let \(\theta\) be a non-inner automorphism of \(\mathrm{Spin}_{2a}\) leaving
      \(\ker ( \mathrm{Spin}_{2a} \rightarrow \mathrm{SO}_{2a})\) invariant (this
      condition is only present for the ``triality'' case \(a=4\)).
      Then for any semisimple \((x, x')\) in \(\mathrm{Spin}_{2a} \times
      \mathrm{Spin}_{2b+1}\), \(f(x,x')\) and \(f(\theta(x), x')\) are conjugated in
      \(\mathrm{Spin}_{2(a+b)+1}\).
    \item
      For \(a,b \geq 2\), consider the natural morphism of complex semisimple
      groups \(f : \mathrm{Spin}_{2a} \times \mathrm{Spin}_{2b} \rightarrow
      \mathrm{Spin}_{2(a+b)}\).
      Let \(\theta_1\) (resp.\ \(\theta_2\), \(\theta\)) be a non-inner automorphism
      of \(\mathrm{Spin}_{2a}\) (resp.\ \(\mathrm{Spin}_{2a}\),
      \(\mathrm{Spin}_{2(a+b)}\)) satisfying the same condition as above.
      Then for any semisimple \((x, x')\) in \(\mathrm{Spin}_{2a} \times
      \mathrm{Spin}_{2b}\), \(f(x,x')\), \(f(\theta_1(x), \theta_2(x'))\),
      \(\theta(f(\theta_1(x), x'))\) and \(\theta(f(x, \theta_2(x')))\) are
      conjugated in \(\mathrm{Spin}_{2(a+b)}\).
  \end{enumerate}
\end{lemm}
\begin{proof}
  This is easily checked on maximal tori.
\end{proof}
\end{comment}

\begin{lemm} \label{lem:unr_endo_gps_isog}
  Let \(\Gbf_\mathrm{cl}\) be a group isomorphic to a product of copies of \(\Spbf_{2n}\)'s and \(\SObf_{4n}\)'s.
  Let \(\Gbf\) be a quotient of \(\Gbf_\mathrm{cl}\) by a (finite) central subgroup.
  The natural map from everywhere unramified elliptic endoscopic data for \(\Gbf\) to everywhere unramified elliptic endoscopic data for \(\Gbf_\mathrm{cl}\) is surjective, and induces a bijection between sets of isomorphism classes.
  Moreover when \(\efrak\) maps to \(\efrak_\mathrm{cl}\) we have \(\iota(\efrak) = \iota(\efrak_\mathrm{cl})\).
\end{lemm}
\begin{proof}
  We have seen in Section \ref{sec:STF} that everywhere unramified elliptic endoscopic data for \(\Gbf\) arise from certain elements \(s \in \Ghat\) such that \(s^2\) maps to \(1 \in \widehat{\Gbf_\mathrm{cl}}\).
  With this description surjectivity and bijectivity are clear.

  Now assume that \(\efrak\) maps to \(\efrak_\mathrm{cl}\).
  Recall that for an endoscopic datum \(\efrak = (\Hcal, \Hbf, s, \xi)\) for \(\Gbf\) we have \(\iota(\efrak) = \tau(\Gbf) \tau(\Hbf)^{-1} |\Out(\efrak)|\).
  Since \(\Gbf\) is split we have \(\tau(\Gbf) = |\pi_0(Z(\Ghat))|\), and similarly for \(\Hbf\), so \(\tau(\Gbf) \tau(\Hbf)^{-1}\) is equal to the corresponding quotient for \(\efrak_\mathrm{cl}\).
  The natural morphism \(\Out(\efrak) \to \Out(\efrak_\mathrm{cl})\) is an isomorphism (this is a general fact which does not use the fact that our groups are split or that the endoscopic data are everywhere unramified).
\end{proof}

We can now state and prove the main result of this section, which is a stable analogue of \cite[Proposition 4.4]{ChRe}.

\begin{prop} \label{pro:lift_stab_mult}
  There is a unique family \((\cpsc(\psi))_{\psi, p}\), for \(p\) a prime number and \(\psi \in \Psit_{\disc, \nonendo}^{\unr, \taut}(\Gbf)\) for some \(\Gbf = \mathbf{SO}_{4n}\) or \(\mathbf{Sp}_{2n}\), \(\tau \in \ICcal(\Gbf)\), such that:
\begin{enumerate}
  \item
    \(\cpsc(\psi)\) is a semisimple conjugacy class in \(\Mpsisc(\C)\) lifting
    \(c_p(\psi)\).
  \item For an endoscopic parameter \(\psi = \oplus_i \psi_i\) define \(\cpsc(\psi) = (\cpsc(\psi_i))_i\), a semisimple conjugacy class in \(\Mpsisc(\C)\).
    Let \(\Gbf\) be a product of groups \(\mathbf{Sp}_{2n}\)'s and
    \(\mathbf{SO}_{4n}\)'s, and \(\Gbf'\) a quotient of  by a central subgroup.
    Then for any \(\tau \in \ICcal(\Gbf')\), the following expansion for the linear form \(S_{\disc}^{\Gbf'}\) on \(\PsCo(\Gbf'_{\R}, \tau) \otimes \Hcal^{\unr}_f(\Gbf')_\C\) holds: if \(f_{\infty}\) is a pseudo-coefficient for some discrete series representations of \(\Gbf'(\R)\) with infinitesimal character \(\tau\), and \(\prod_p f_p \in \Hcal^{\unr}_f(\Gbf')_\C\), we have
    \begin{equation} \label{eq:Sdisclifted}
      S^{\Gbf'}_{\disc}(f_{\infty} \prod_p f_p) = \sum_{\psi \in
      \tilde{\Psi}^{\unr, \tau}_{\disc}(\Gbf)}
      \frac{\epsilon_{\psi}(s_{\psi}) (-1)^{q(\Lbf^*_{\psi,\tau})}}{|\Scal_{\psi}|} \prod_p
      \Sat_{\Gbf'_{\Zp}}(f_p)( \dpsitausc(\cpsc(\psi)) )
    \end{equation}
\end{enumerate}
\end{prop}
\begin{proof}
  The proof is similar to that of Proposition \ref{pro:refineSOeven}, only
  slightly more complicated.
  We show by induction on \(N \geq 0\) that there is a unique family
  \((\cpsc(\psi))_{\psi, p}\) for \(\psi \in \tilde{\Psi}^{\unr,
  \tilde{\tau}}_{\disc, \nonendo}(\Gbf)\) for some \(\Gbf\) of dimension \(\leq N\),
  such that \eqref{eq:Sdisclifted} is satisfied for \(\Gbf\) of a product of
  groups of dimension \(\leq N\).

  The case \(N=0\) is obvious.
  Assume that the induction hypothesis is satisfied for some \(N \geq 0\).
  Let \(\Gbf = \mathbf{SO}_{4n}\) or \(\mathbf{Sp}_{2n}\) be of dimension \(N+1\).
  By Proposition \ref{pro:csqChRe} and a calculation using the stabilization of
  the trace formula (Theorem \ref{thm:STF_qs}) similar to
  \eqref{eq:stab_spec_SOeven} in Proposition \ref{pro:refineSOeven} (i.e.\
  essentially a special case of Kottwitz' stabilization of the spectral side of
  the trace formula; the calculation really is almost identical by Lemma
  \ref{lem:unr_endo_gps_isog} and the last point of Proposition-Definition
  \ref{def:Mpsisc}), we have, for any \(\tau \in \ICcal(\Gbfad)\),
  \(\delta_{\infty}\) a discrete series representation of \(\Gbf_{\mathrm{ad}}(\R)\)
  having infinitesimal character \(\tau\), \(f_{\infty}\) a pseudo-coefficient for
  \(\delta_{\infty}\) and \(\prod_p f_p \in \Hcal^{\unr}_f(\Gbfad)_\C\),
  \begin{align*}
    S^{\Gbfad}_{\disc, \tau} (f) = & \sum_{\psi \in \tilde{\Psi}_{\mathrm{disc}}^{\unr, \tau}(\Gbf)} \frac{(-1)^{q(\Lbf^*_{\psi,\tau})}}{|\Scal_{\psi}|} \\
    & \left[ \sum_{s \in \mathcal{S}_{\psi}} \epsilon_{\psi}(s s_{\psi}) \langle s s_{\psi}, \delta_{\infty} \rangle \prod_p \Sat_{\Gbf_{\ad, \Zp}}(f_p)(\cpsc'(\Gbf, \psi, \delta_{\infty}))  \right. \\
    & \left. - \sum_{s \in \Scal_{\psi},\, s \neq 1} \epsilon_{\psi}(s s_{\psi}) \langle s s_{\psi}, \delta_{\infty} \rangle \prod_p \mathrm{Sat}(f_p)(\dpsitausc( \cpsc(\psi) )) \right]
  \end{align*}
  where \(\cpsc'(\Gbf, \psi, \delta_{\infty})\) is defined in Proposition
  \ref{pro:csqChRe} if \(\langle \cdot, \delta_{\infty} \rangle|_{\Scal_{\psi}} =
  \epsilon_{\psi}\), and for convenience we define \(\cpsc'(\Gbf, \psi,
  \delta_{\infty}) = \dpsitausc(\cpsc(\psi))\) otherwise (using the induction
  hypothesis).
  Note that since we are only considering pseudo-coefficients of discrete series
  at the real place, we know that the only relevant endoscopic groups (those for
  which the transfer of \(f\) does not vanish in \(SI(\Hbf)\)) are everywhere
  unramified and have discrete series at the real place (avoiding the analogue
  of Lemma \ref{lem:SO4n+2}, although it does hold true).
  Note also that these relevant endoscopic groups for \(\Gbfad\) are quotients
  of products of groups of the form \(\SObf_{4m}\) or \(\Spbf_{2m}\) which are not
  adjoint in general, which is why it is necessary to include arbitrary
  quotients by central subgroups in the induction hypothesis.

  We know that \(f \mapsto S_{\disc, \tau}^{\Gbfad}(f)\) is a stable linear form
  on \(\PsCo(\Gbf_{\mathrm{ad}, \R}, \tau) \otimes \Hcal^{\unr}_f(\Gbfad)\), and
  by the same argument as in the proof of Proposition \ref{pro:refineSOeven}
  using Remark \ref{rem:mult_one} we deduce that the contribution of each \(\psi\)
  in the summand above is itself stable.
  If \(\psi\) is not endoscopic, by the same argument as in the proof of
  Proposition \ref{pro:refineSOeven} we see that stability implies that
  \(\cpsc'(\Gbf, \psi, \delta_{\infty})\) does not depend on the choice of
  \(\delta_{\infty}\) in its L-packet, and this defines \(\cpsc(\psi)\).
  If \(\psi\) is endoscopic, we also conclude as in Proposition
  \ref{pro:refineSOeven}, using Lemma \ref{lem:enough_elts_AJ_SOeven} in the
  even orthogonal case and the analogue in the symplectic case, which is obvious
  since \(H^1(\R, \mathbf{Sp}_{2n})\) is trivial.
  We obtain that for \(\psi \in \tilde{\Psi}_{\disc, \mathrm{endo}}^{\unr,
  \tilde{\tau}}(\Gbf)\), for any \(\delta_{\infty}\) such that \(\langle \cdot,
  \delta_{\infty} \rangle|_{\Scal_{\psi}} = \epsilon_{\psi}\), we have
  \(\cpsc'(\Gbf, \psi, \delta_{\infty}) = \dpsitausc(\cpsc(\psi))\), and
  \eqref{eq:Sdisclifted} for \(\Gbfad\).

  To conclude the proof of Proposition \ref{pro:lift_stab_mult} we now have to
  check \eqref{eq:Sdisclifted} for \(\Gbf\) a product of \(\mathbf{Sp}_{2n}\)'s and
  \(\mathbf{SO}_{4n}\)'s each of dimension \(\leq N+1\) and \(\Gbf'\) any quotient of
  \(\Gbf\) by a central subgroup.
  Since we now have \eqref{eq:Sdisclifted} for \(\Gbfad\) and its endoscopic
  groups, using one more time the stabilization of the trace formula for
  \(\Gbfad\) and the analogue of \eqref{eq:stab_spec_SOeven} as above, we obtain
  a weak Arthur multiplicity formula for \(\Gbfad\): for any \(\tau \in
  \ICcal(\Gbfad)\) and any \(f \in \PsCo(\Gbf_{\mathrm{ad}, \R}, \tau) \otimes
  \Hcal^{\unr}_f(\Gbfad)\) we have
  \[ I_{\disc, \tau}^{\Gbfad}(f) = \sum_{\substack{\psi \in \tilde{\Psi}^{\unr, \tilde{\tau}}_{\disc}(\Gbf) \\ \langle \cdot, \delta_{\infty} \rangle|_{\Scal_{\psi}} = \epsilon_\psi}} \epsilon_{\psi}(s_\psi) (-1)^{q(\Lbf^*_{\psi,\tau})} \prod_p \Sat_{\Gbf_{\ad, \Zp}} (f_p) (\dpsitausc(\cpsc(\psi))). \]
  Another application of \cite[Proposition 4.4]{ChRe} applied to \(\Gbf'
  \rightarrow \Gbfad\), along with Lemmas \ref{lem:prDS_isog} and
  \ref{lem:parity_orth}, shows that a similar weak multiplicity formula holds
  with \(\Gbf_\ad\) replaced by \(\Gbf'\).
  Finally, \eqref{eq:Sdisclifted} for \(\Gbf'\) is a consequence of this weak
  multiplicity formula, again by the same argument using the stabilization of
  the trace formula for \(\Gbf'\) and the analogue of \eqref{eq:stab_spec_SOeven}.
\end{proof}

\begin{rema} \label{rem:lift_nonendo}
  The proofs of Proposition \ref{pro:csqChRe} and \ref{pro:lift_stab_mult} show that for \(\Gbf = \SObf_{4n}\) or \(\Spbf_{2n}\), \(\tau \in \IC(\Gbf)\) and \(\psi \in \Psit^{\unr,\taut}_{\disc,\nonendo}(\Gbf)\) the eigenspace
  \[ \Acal^2(\Gbf_\ad(\Q) \backslash \Gbf_\ad(\A) / \Gbf_\ad(\Zhat))_\tau[(\dpsitau(c_p(\psi)))_p] \]
  for the character of the subalgebra \(\Hcal^\unr(\Gbf) \subset \Hcal^\unr(\Gbf_\ad)\) corresponding to the Satake parameters \(\dpsitau(c_p(\psi))\), is also eigen for \(\Hcal^\unr(\Gbf_\ad)\), with character corresponding to the Satake parameters \(\dpsitausc(\cpsc(\psi))\).
\end{rema}

\begin{coro}[Weak multiplicity formula] \label{cor:weak_mult_formula}
  Let \(\widetilde{\Gbf}^*\) be either \(\Spbf_{2n}\) or \(\SObf_{4n}\) and let \((\Gbf, \Xi, c)\) be an inner form of \(\Gbf^* := \widetilde{\Gbf}^*_{\ad}\) (see Section \ref{sec:STF}) which is split at all finite places of \(\Q\).
  Fix a reductive model of \(\Gbf\) over \(\Z\) (see \cite[Proposition 1.1]{Gross_overZ}).
  Let \(\tau \in \ICcal(\Gbf) := \ICcal(\Gbf^*)\), let \(\delta_\infty\) be a discrete series representation of \(\Gbf(\R)\) having infinitesimal character \(\tau\) and \(f_{\infty}\) be a pseudo-coefficient for \(\delta_\infty\).
  Associated to the localization \((\Gbf_{\R}, \Xi_{\R}, c_{\R})\) and \(\delta_\infty\) is a character \(\langle \cdot, \delta_\infty \rangle\) of \(\Scal_{\varphi_\tau}\) (we reviewed the construction in detail in \cite[\S 4.2.1]{Taibi_dimtrace} and \cite[\S 3.2.1]{Taibi_mult}), which we may restrict to \(\Scal_{\dpsitau \circ \psi_\infty}\) after aligning the two parameters (see Proposition \ref{pro:tr_AJ_ps}).
  Then for any \(\prod_p f_p \in \Hcal_f^{\unr}(\Gbf)_\C\) we have
  \begin{equation} \label{eq:weak_mult_formula}
    I^{\Gbf}_{\disc, \tau}(f_{\infty} \prod_p f_p) = \sum_{\substack{\psi \in \tilde{\Psi}^{\unr, \tau}_{\disc}(\widetilde{\Gbf}^*) \\ \langle \cdot, \delta_\infty \rangle |_{\Scal_\psi} = \epsilon_{\psi}}} \epsilon_{\psi}(s_{\psi}) (-1)^{q(\Lbf^*_{\psi,\tau})} \prod_p \Sat_{\Gbf_{\Zp}}(f_p)( \dpsitausc(\cpsc(\psi)) ).
  \end{equation}
\end{coro}
\begin{proof}
  For \(\Gbf = \Gbf^*\) this was proved in the proof of Proposition \ref{pro:lift_stab_mult}.
  For arbitrary \(\Gbf\) the formula follows from the stable expansion \eqref{eq:Sdisclifted} for all relevant elliptic endoscopic groups of \(\Gbf\) and the stabilization of the trace formula for \(\Gbf\) (Theorem \ref{thm:STF_inner}).
\end{proof}

\begin{exam} \label{exam:mult_formula_definite}
  A special case that will be useful in Section \ref{sec:even_spin} is when \(\Gbf\) is the inner form of \(\PGSObf_{8n}\) which is split at all finite places and anisotropic at the real place.
  For example a model of \(\Gbf\) over \(\Z\) is then given by \(\PGSObf(q)\) (defined as in \cite[\S C.3]{Conrad_luminy_sga3}) where \(q\) is the quadratic form on the even unimodular lattice \(E_8^{\oplus n}\).
  For such a group any \(\tau \in \ICcal(\Gbf)\) determines a unique (and finite-dimensional!) discrete series representation \(\pi_{\infty}\) of the compact connected group \(\Gbf(\R)\), so that in this case \eqref{eq:weak_mult_formula} really is a multiplicity formula.
  Let us make explicit the character \(\langle \cdot, \pi_{\infty} \rangle\) of \(\Scal_{\varphi_\tau}\), omitting the details which may be found in \cite[Example 3.2.3]{Taibi_mult}.
  We identify \(\Ghat\) with \(\Spin_{8n}\) and use notation as in Section \ref{sec:not_red_gps}.
  Up to conjugation by \(\Spin_{8n}(\C)\) we may assume \(\varphi_\tau(\C^{\times}) \subset \Tcal_{\Spin_{8n}}(\C)\) and that \(\varphi_\tau\) is dominant.
  Then \(C_{\varphi_\tau}\) is the group of \((z_1, \dots, z_{4n}, s) \in \Tcal_{\Spin_{8n}}(\C)\) satisfying \(z_i \in \{\pm 1\}\) and \(\prod_{i=1}^{4n} z_i = 1\) (i.e.\ \(s \in \{\pm 1\}\) as well).
  The character \(\langle \cdot, \pi_{\infty} \rangle\) of \(C_{\varphi_{\tau}}\) is \((z_1, \dots, z_{4n}, s) \mapsto z_1 z_3 \dots z_{4n-1}\).
\end{exam}

\subsection{Exceptional isomorphisms in small rank}

In small dimension we have exceptional isomorphisms:
\begin{enumerate}
\item for \(\Gbf = \Spbf_2\) we have \(\Ghat_\sico \simeq \Spin_3 \simeq \SL_2\),
\item for \(\Gbf = \SObf_4\) we have \(\Ghat_\sico \simeq \Spin_4 \simeq \SL_2 \times \SL_2\),
\item for \(\Gbf = \Spbf_4\) we have \(\Ghat_\sico \simeq \Spin_5 \simeq \Sp_4\).
\end{enumerate}
In these three cases the families \((\cpsc(\psi))_p\) (defined by Proposition \ref{pro:lift_stab_mult}) are not really new, as the following three propositions show for \(\psi\) of the form \(\pi[1]\), essentially rephrasing \cite[\S 4]{ChRe}.
(The non-tempered cases \(\pi[d]\) with \(d>1\) will be obtained later as special cases of more general Arthur-Langlands functoriality.)

\begin{prop} \label{pro:cpsc_Sym2}
  Let \(w\) be a positive integer.
  There is a bijection denoted by \(\Sym^2\) between
  \begin{itemize}
  \item level one cuspidal automorphic representation for \(\PGLbf_2\) having infinitesimal character \((\pm w/2)\) at the real place,
  \item self-dual level one cuspidal automorphic representation for \(\PGLbf_3\) having infinitesimal character \((w,0,-w)\) at the real place.
  \end{itemize}
  For \(\pi\) in the first set its image \(\Sym^2 \pi\) is characterized by the relation \(c((\Sym^2 \pi)_p) = \Sym^2 c(\pi_p)\) and we further have \(\cpsc((\Sym^2 \pi)[1]) = c(\pi_p)\) via the isomorphism \(\SL_2 \simeq \Spin_3\).
\end{prop}
Note that the first set is naturally in bijection with eigenforms in \(S_{w+1}(\SL_2(\Z))\), in particular it is empty if \(w\) is even, as is the second set by Remark \ref{rem:no_param_unless_triv_res_center}.
\begin{proof}
  The fact that is induces a bijection in level one is \cite[Proposition 4.9]{ChRe}, whose proof also shows the relation \(\cpsc((\Sym^2 \pi)[1]) = c(\pi_p)\) in view of Remark \ref{rem:lift_nonendo}.
  The lift \(\Sym^2\) was first constructed by Gelbart and Jacquet \cite[Theorem 9.3]{GelJac}.
\end{proof}

\begin{prop} \label{pro:cpsc_tensor_PGL2_PGL2}
  Let \(w_1, w_2 \in \frac{1}{2} \Z \smallsetminus \Z\) satisfying \(w_1 > w_2 > 0\).
  There is a bijection denoted by \((\pi_1, \pi_2) \mapsto \pi_1 \otimes \pi_2\) between
  \begin{itemize}
  \item pairs of level one cuspidal automorphic representation for \(\PGLbf_2\) having infinitesimal characters \(\pm w_1\) and \(\pm w_2\) at the real place,
  \item self-dual level one cuspidal automorphic representation for \(\PGLbf_4\) having infinitesimal character \((\pm (w_1+w_2), \pm (w_1-w_2))\) at the real place.
  \end{itemize}
  For \((\pi_1,\pi_2)\) in the first set its image \(\pi_1 \otimes \pi_2\) is characterized by the relation (for all primes \(p\)) \(c((\pi_1 \otimes \pi_2)_p) = c(\pi_{1,p}) \otimes c(\pi_{2,p})\) and we further have \(\cpsc((\pi_1 \otimes \pi_2)[1]) = (c(\pi_{1,p}), c(\pi_{2,p}))\) via the isomorphism \(\Mcal_{(\pi_1 \otimes \pi_2)[1], \sico} \simeq \SL_2 \times \SL_2\) mapping \(\tau_{(\pi_1 \otimes \pi_2)[1]}\) to \(((\pm w_1), (\pm w_2))\).
\end{prop}
\begin{proof}
  This bijection is characterized in \cite[Proposition 4.10]{ChRe} (which also handles the case \(w_1=w_2\) which we ignore in this paper).
  The proof loc.\ cit.\ also shows the relation \(\cpsc((\pi_1 \otimes \pi_2)[1]) = (c(\pi_{1,p}), c(\pi_{2,p}))\) in view of Remark \ref{rem:lift_nonendo}.
  The tensor product lift is a special case of \cite[Theorem M]{Ramakrishnan_SL2}.
\end{proof}

\begin{prop} \label{pro:cpsc_Sp4_SO5}
  For any \(w_1, w_2 \in \frac{1}{2} \Z \smallsetminus \Z\) satisfying \(w_1 > w_2 > 0\) there is a bijection denoted by \(\Lambda^*\) between
  \begin{itemize}
  \item self-dual level one cuspidal automorphic representation for \(\PGLbf_4\) having infinitesimal character \((\pm w_1, \pm w_2)\) at the real place,
  \item self-dual level one cuspidal automorphic representation for \(\PGLbf_5\) having infinitesimal character \((w_1+w_2,w_1-w_2,0,w_2-w_1,-w_1-w_2)\) at the real place
  \end{itemize}
  and we have \(\spin(\cpsc((\Lambda^* \pi)[1])) = c(\pi_p)\) (this determines the semi-simple conjugacy class \(\cpsc((\Lambda^* \pi)[1])\) in \(\Spin_5(\C) \simeq \Sp_4(\C)\) since the spin representation of \(\Spin_5\) is the standard representation of \(\Sp_4\).).
\end{prop}
\begin{proof}
  Again this is proved as \cite[Proposition 4.12]{ChRe} with Remark \ref{rem:lift_nonendo}.
  The lift \(\Lambda^*\) can also be constructed using the exterior square lift constructed by Kim \cite{Kim_exterior_sq_GL4} and \cite[Theorem 1.5.3 (a)]{Arthur}.
\end{proof}

\subsection{Rationality of Satake parameters}

The following rationality property will not surprise any expert.
It implies that all statements seemingly relying on a choice of isomorphism \(\iota: \C \simeq \Qellbar\) only depend on the restriction of \(\iota\) to the algebraic closure of \(\Q\).

\begin{prop} \label{pro:rat_Sat}
  Let \(\Gbf = \Spbf_{2n}\) (resp.\  \(\SObf_{4n}\)) for \(n \geq 1\), and let
  \(\delta = n(n+1)/2\) (resp.\ \(n\)).
  Then for any \(\tau \in \ICcal(\Gbf)\) and any \(\psi \in
  \tilde{\Psi}_{\disc}^{\unr, \tau}(\Gbf)\) there exists a finite extension \(E\)
  of \(\Q\) in \(\C\) such that for any prime \(p\), the semisimple conjugacy class
  \(p^{\delta/2} \dpsitausc(\cpsc(\psi))\) in \(\GSpin_{2n+1}(\C)\) (resp.\
  \(\GSpin_{4n}(\C)\)) is defined over \(E\), i.e.\ its trace in any algebraic
  representation belongs to \(E\).
\end{prop}
\begin{proof}
  First recall that \(p^{\delta/2} \dpsitausc(\cpsc(\psi))\) being defined over
  \(E\) is equivalent to the unramified representation
  \(\pi(\dpsitausc(\cpsc(\psi)))\) of \(\Gbf_{\ad}(\Qp)\) corresponding to
  \(\dpsitausc(\cpsc(\psi))\) being defined over \(E\): see \S 2.2, 5.2 and 5.3 of
  \cite{BuzGee} (using the central extension \(\GSpbf_{2n} \rightarrow
  \PGSpbf_{2n}\) resp.\ \(\GSObf_{4n} \rightarrow \PGSObf_{4n}\)).
  Moreover it is easy to check that this rationality property is compatible with
  endoscopy: if \(\psi = \oplus_i \psi_i\) and each \(\psi_i\) satisfies the
  rationality property, then so does \(\psi\).
  Without loss of generality, we may therefore assume that \(\psi\) is
  non-endoscopic, i.e.\ \(\psi \in \tilde{\Psi}^{\unr, \tilde{\tau}}_{\disc,
  \nonendo}(\Gbf)\).

  It should possible to deduce the rationality property from \cite[Corollary
  2.18]{ShinTemplier_fields}.
  However this deduction does not seem completely obvious to us, so we may as
  well adapt the argument (which goes back to \cite{Waldspurger_proparithGL2}
  and \cite[\S 3.5]{Clozel_AA}).
  Below we will also sketch an alternative argument using ``only'' the trace
  formula.

  Let \(K_{\ad}\) be a maximal compact subgroup of \(\Gbf_{\ad}(\R)\), and let \(V\)
  be the irreducible algebraic representation of \(\Gbf_{\ad}(\R)\) having
  infinitesimal character \(-\tau\).
  Note that it is defined over \(\Q\), i.e.\ it originates from an algebraic
  representation of \(\Gbf_{\ad}\).
  By \cite[Lemma 2.2.]{ArthurL2} summing \eqref{eq:weak_mult_formula} over all
  discrete series representations \(\delta_\infty\) having infinitesimal character
  \(\tau\) gives
  \begin{align} \label{eq:L2_Lef_Gad}
    & (-1)^{q(\Gbf_\R)} \sum_{i=0}^{2 q(\Gbf_\R)} (-1)^i \Tr \left( \prod_p f_p
    \, \middle| \, H^i \left( \gfrak, K_{\ad}, \Acal^2(\Gbf_\ad(\Q) \backslash
    \Gbf_\ad(\A)/ \Gbf_{\ad}(\Zhat)) \otimes V \right) \right) \\
  = & \sum_{\psi \in \tilde{\Psi}^{\unr, \tau}_{\disc}(\widetilde{\Gbf})}
    N_\psi \epsilon_{\psi}(s_{\psi}) (-1)^{q(\Lbf^*_{\psi,\tau})} \prod_p
    \Sat_{\Gbf_{\Zp}}(f_p)( \dpsitausc(\cpsc(\psi)) ) \nonumber
  \end{align}
  where \(N_\psi\) is the number of discrete series representations
  \(\delta_\infty\) of \(\Gbf_{\ad}(\R)\) in the L-packet corresponding to \(\tau\)
  which satisfy \(\langle \cdot, \delta_\infty \rangle|_{\Scal_\psi} =
  \epsilon_\psi\).
  For a non-endoscopic \(\psi\) we have \(\Scal_\psi = 1\) and so \(N_\psi\) is
  positive.
  Using Remark \ref{rem:mult_one} this implies that the unramified
  representation of \(\Gbf_{\ad}(\A_f)\) corresponding to
  \((\dpsitausc(\cpsc(\psi)))_p\) occurs in
  \[ \varinjlim_{K_f} H^i(\gfrak, K_{\ad}, \Acal^2(\Gbf_\ad(\Q) \backslash
  \Gbf_\ad(\A)/K_f) \otimes V) \]
  for some \(i\), where the limit is over compact open subgroups of
  \(\Gbf_{\ad}(\A_f)\)
  \footnote{We could consider only level \(\Gbf_{\ad}(\Zhat)\) but that would
  require dealing with orbifolds \dots}.
  The weighted cohomology groups of \cite{GoKoMPh} for the group \(\Gbf_{\ad}\),
  neat level \(K_f\), coefficient system corresponding to \(V\) and the upper middle
  weight profile are vector spaces over \(\Q\) with a Hecke algebra action, which
  by \cite[Corollary B]{Nair_weighted} give a rational structure to these
  \((\gfrak, K_{\ad})\)-cohomology groups.
  This concludes the proof of the proposition.

  Alternatively to prove the proposition one could use Arthur's \(L^2\) Lefschetz
  trace formula \cite{ArthurL2}, which provides another expansion for
  \eqref{eq:L2_Lef_Gad}.
  Let us sketch the argument.
  First one checks that the geometric side of this trace formula is rational if
  the Hecke operators take only rational values, i.e.\ if we have \(\prod_p f_p
  \in \Hcal_f^{\unr}(\Gbf_{\ad})\).
  Choosing a rational Haar measure on \(\Gbf_{\ad}(\A_f)\), i.e.\ one giving
  rational measure to any compact open subgroup, it is easy to see that that
  orbital integrals of rational Hecke operators at semisimple elements are
  rational: it follows from Harish-Chandra's lemma recalled in \cite[\S
  3.1.2]{Taibi_dimtrace} that these orbital integrals are rational linear
  combinations of values of the Hecke operator under consideration.
  Using the formula for Tamagawa numbers, a comparison of Haar measures
  \cite[Theorem 9.9]{Gross_mot} and rationality of the L-function of the motive
  that Gross associated to \(\Gbf_{\ad}\) (Proposition 9.5 loc.\ cit.), one gets
  that the elliptic terms of the geometric side are rational.
  The rationality of the other (``parabolic'') terms on the geometric side is
  proved similarly, starting from the expression \cite[(3.3.1)]{Taibi_dimtrace}
  and using a formula for \(\Phi_{\Mbf}\) given on p.300 loc.\ cit., noting that
  the potentially irrational factor \(\delta_P^{1/2}\) is compensated by a similar
  factor in the normalized constant term of the Hecke operator.
  To conclude the proof of the proposition, apply Lemma \ref{lem:virt_rep_rat}
  below to the virtual representation of \(\Hcal_f^{\unr}(\Gbf_{\ad})_\C\)
  appearing in \eqref{eq:L2_Lef_Gad}.
\end{proof}

\begin{lemm} \label{lem:virt_rep_rat}
  Let \(A\) be a unital associative algebra over \(\Q\), and let \(\sum_{i \in I}
  \lambda_i [V_i, \rho_i]\) be an element of the Grothendieck group of
  finite-dimensional representations of \(\C \otimes_{\Q} A\), where \(I\) is
  finite, \(\lambda_i \in \Z \smallsetminus \{0\}\) and \(\rho_i: \C
  \otimes_{\Qbar} A \to \End_{\C}(V_i)\) are irreducible finite-dimensional
  pairwise non-isomorphic representations.
  Assume that for any \(a \in A\) we have \(\sum_{i \in I} \lambda_i \Tr \rho_i(a)
  \in \Q\).
  Then for any \(i \in I\) there exists a finite extension \(E_i/\Q\) in \(\C\) and a
  representation \(W_i\) of \(E_i \otimes_{\Q} A\) such that \(V_i\) is isomorphic to
  \(\C \otimes_{E_i} W_i\).
\end{lemm}
\begin{proof}
  By the Artin-Wedderburn theorem the map \(\C \otimes_{\Q} A \to \prod_{i \in I}
  \End_{\C}(V_i)\) is surjective and so \(\cap_{i \in I} \ker \rho_i\) is equal to
  \[ \left\{ a \in \C \otimes_{\Q} A \,\middle|\, \forall b \in \C \otimes_{\Q}
  A,\, \sum_{i \in I} \lambda_i \Tr \rho_i(ab) = 0 \right\}. \]
  This equality and the assumption imply that we have \(\cap_{i \in I} \ker
  \rho_i = \C \otimes_{\Q} J\) where \(J\) is the ideal of \(A\) defined by
  \[ J = \left\{ a \in A \,\middle|\, \forall b \in A,\, \sum_{i \in I}
  \lambda_i \Tr \rho_i(ab) = 0 \right\}. \]
  We may replace \(A\) by \(A/J\), which has finite dimension over \(\Q\).
  The Jacobson radical \(R\) of \(A\) is a nilpotent ideal, so \(\C \otimes_{\Q} R\)
  is a nilpotent ideal of the semi-simple algebra \(\C \otimes_{\Q} A\), and so we
  have \(R = 0\), i.e.\ \(A\) is semi-simple and applying the Artin-Wedderburn
  theorem again, we have \(A \simeq \prod_{k \in K} M_{n_k}(D_k)\) where \(K\) is
  finite and \(D_k\) is a finite-dimensional division algebra over \(\Q\).
  If \(E/\Q\) is a finite extension of \(\Q\) in \(\C\) splitting all the \(D_k\)'s then
  it is clear that any irreducible representation of \(\C \otimes_{\Q} A\) is
  defined over \(E\).
\end{proof}

\subsection{Non-semisimple groups}
\label{sec:non_ss_gps}

In this paper we will be mainly interested in \(\GSpbf_{2n}\) and its endoscopic
groups.
We now justify that we can reduce to their (semisimple) quotients by \(\GLbf_1\).
For a connected reductive group \(\Gbf\) over \(\Q\) which is not semisimple (i.e.\
such that the connected center \(\Zbf\) of \(\Gbf\) is a non-trivial torus), to
formulate the trace formula and its stabilization for \(\Gbf\) it is convenient to
fix a closed subgroup \(Z\) of \(\Zbf(\A)\) such that \(Z \Zbf(\Q)\) is a closed
subgroup of \(\Zbf(\A)\) and \(Z \Zbf(\Q) \backslash \Zbf(\A)\) is compact.
We shall only need two cases: \(Z = \Abf_{\Gbf}(\R)^0\) or \(Z = \Zbf(\A)\).
It is also necessary to fix a unitary character \(\chi\) of \(Z\) which is trivial
on \(Z \cap \Zbf(\Q)\).
Since \(\Gbf(\A) = \Abf_{\Gbf}(\R)^0 \times \Gbf(\A)^1\) where
\[ \Gbf(\A)^1 = \left\{ g \in \Gbf(\A) | \forall \delta : \Gbf \rightarrow
\GLbf_1,\, |\delta(g)| = 1 \right\} \]
we can reduce by twisting to the case where \(\chi|_{\Abf_{\Gbf}(\R)^0} = 1\).
Fix a maximal compact subgroup \(K_\infty\) of \(\Gbf(\R)\).
The volume of \(Z \Gbf(\Q) \backslash \Gbf(\A)\) is finite and we can consider the
space of \(\chi\)-equivariant automorphic forms on \(\Gbf(\Q) \backslash \Gbf(\A)\)
which are square-integrable modulo \(Z\), denoted \(\Acal^2(\Gbf, Z, \chi)\), and
for \(\tau\) a semisimple conjugacy class in \(\widehat{\gfrak}\) we can consider
the eigenspace \(\Acal^2(\Gbf, Z, \chi)_\tau\) for the character of the center of
the enveloping algebra of \(\gfrak = \C \otimes \Lie \Gbf(\R)\) corresponding to
\(\tau\).
When \(Z\) intersects \(\Zbf(\R)\) non-trivially there is an obvious necessary
condition relating \(\chi\) and \(\tau\) for this subspace to be non-zero.
There is a decomposition
\[ \Acal^2(\Gbf, Z, \chi)_\tau \simeq \bigoplus_{\pi \in \Pi_{\disc}(\Gbf, Z,
\chi)} \pi^{\oplus m(\pi)} \]
where \(\Pi_{\disc}(\Gbf, Z, \chi, \tau)\) is a countable set of isomorphism
classes of irreducible admissible unitarizable \((\gfrak, K_\infty) \times
\Gbf(\A_f)\)-modules \(\pi = \pi_\infty \otimes \pi_f\) such that the restriction
to \(Z\) of the central character of \(\pi\) is \(\chi\) and \(\pi_\infty\) has
infinitesimal character \(\tau\), and \(m(\pi) \in \Z_{\geq 1}\).
Only finitely many elements of \(\Pi_{\disc}(\Gbf, Z, \chi, \tau)\) have non-zero
invariants under any given open subgroup of \(\Gbf(\A_f)\).
The contribution of the discrete spectrum
to the discrete (sic) part \(I_{\disc, \tau}^{\Gbf, Z, \chi}(f(g) dg)\) of the
spectral side of the trace formula for \((\Gbf, Z, \chi)\) and fixed infinitesimal
character \(\tau\) is then 
\[ \sum_{\pi \in \Pi_{\disc}(\Gbf, Z, \chi, \tau)} m(\pi) \Tr \pi(f(g) dg) \]
which clearly factors through \(f(g) dg \mapsto f_{Z, \chi}(g) d \bar{g}\) where
\(d \bar{g}\) is the quotient Haar measure \(dg / dz\) on \(\Gbf(\A) / Z\) and \(f_{Z,
\chi}(g) := \int_Z \chi(z) f(z g) dz\) is \(\chi^{-1}\)-equivariant, smooth and
compactly supported modulo \(Z\) and bi-\(K_\infty\)-finite.
Denoting \(\Hcal(\Gbf, Z, \chi)\) the space of \(\chi^{-1}\)-equivariant, smooth, 
compactly supported modulo \(Z\) and bi-\(K_\infty\)-finite distributions on
\(\Gbf(\A)\), the map \(\Hcal(\Gbf) \to \Hcal(\Gbf, Z, \chi)\), \(f(g) dg \mapsto
f_{Z, \chi} (\bar{g}) d \bar{g}\) is easily seen to be surjective.
Denote \(I(\Gbf, Z, \chi)\) and \(SI(\Gbf, Z, \chi)\) the quotients of \(Hcal(\Gbf,
Z, \chi)\) obtained by considering orbital integrals and stable orbital integrals
at regular semisimple elements.
We do not recall the other contributions to the discrete part of the trace
formula, since they do not play any role in the present article (as recalled in
the ``proof'' of Theorem \ref{thm:STF_qs} these other contributions vanish if
\(\tau\) is regular).

One can easily compare the linear forms \(I_{\disc, \tau}^{\Gbf, Z, \chi}\) for
varying \(Z\).
For the two cases in which we are interested (\(Z = \Abf_{\Gbf}(\R)^0\) or
\(\Zbf(\A)\)) we have
\begin{equation} \label{eq:compare_TF_central}
  I_{\disc, \tau}^{\Gbf, \Abf_{\Gbf}(\R)^0, \chi}(f(g) dg) = \sum_{\chi'}
  I_{\disc, \tau}^{\Gbf, \Zbf(\A), \chi'}(f(g) dg)
\end{equation}
where the sum is over all characters \(\chi'\) extending \(\chi\).
For a given level (i.e.\ compact open subgroup \(K_f\) of \(\Gbf(\A_f)\) under which
\(f\) is bi-invariant) only finitely many \(\chi'\) may have a non-zero contribution
on the right-hand side of \eqref{eq:compare_TF_central}.
The special case which is relevant for this paper is for \(\Gbf\) a reductive
group over \(\Z\) (this implies \(\Zbf = \Abf_{\Gbf}\)) and level \(\Gbf(\Zhat)\),
then there is at most one \(\chi'\) such that the corresponding term in
\eqref{eq:compare_TF_central} may be non-zero, since \(\Abf_{\Gbf}(\A) =
\Abf_{\Gbf}(\R)^0 \Abf_{\Gbf}(\Q) \Abf_{\Gbf}(\Zhat)\).
In particular if \(\chi = 1\) then only \(\chi' = 1\) may have a non-vanishing
contribution.
Moreover for any \(R \in \{\Q, \Q_v, \Zp, \A \}\) the morphism \(\Gbf(R)
\rightarrow (\Gbf / \Abf_{\Gbf})(R)\) is surjective and it is easy to check that
\[ I_{\disc, \tau}^{\Gbf, \Abf_{\Gbf}(\A), 1}(f) = I_{\disc, \tau}^{\Gbf /
\Abf_{\Gbf}}(f_{\Abf_{\Gbf}(\A), 1}). \]
Note that this reduction is compatible with the Satake isomorphism: the
integration map \(i : \Hcal^{\unr}(\Gbf_{\Zp}) \rightarrow
\Hcal^{\unr}((\Gbf / \Abf_{\Gbf})_{\Zp})\) is surjective and satisfies
\(\Sat_{(\Gbf / \Abf_{\Gbf})_{\Zp}}(i(f))(c) = \Sat_{\Gbf_{\Zp}}(f)(\iota(c))\)
where \(\iota\) is the natural map \({}^L (\Gbf / \Abf_{\Gbf}) \rightarrow {}^L
\Gbf\) (or rather the induced injection from \(\widehat{\Gbf /
\Abf_{\Gbf}}\)-conjugacy classes to \(\Ghat\)-conjugacy classes).

A similar reduction to semisimple groups holds for the stabilization of the
trace formula, replacing \(I_{\disc, \tau}\) by \(S_{\disc, \tau}\) above.
The proof is an obvious induction using the fact that \(\Abf_{\Hbf} =
\Abf_{\Gbf}\) if \(\Hbf\) is an elliptic endoscopic group of \(\Gbf\).

\section{Intersection cohomology of \(\Acal_n^*\)}

In this section we first recall definitions from \cite{MorelSiegel1}: Siegel
modular varieties \(\Acal_{n,K}\), automorphic \(\ell\)-adic étale sheaves
\(\Fcal^K(V)\), Hecke correspondences, arithmetic minimal compactifications
\(\Acal_{n,K} \hookrightarrow \Acal_{n,K}^*\), intermediate extensions \(\IC^K(V)\)
of \(\Fcal^K(V)\) (over \(\Fp\) for \(p \neq \ell\), the case considered in
\cite{MorelSiegel1}, or \(\Q\)) and the canonical extension of Hecke
correspondences.
We check that the specialization isomorphism between cohomology groups over \(\Fp\) and \(\Q\) provided by \cite{Stroh_Kottwitz} (see also \cite{LanStroh_nearby1} and \cite{LanStroh_nearby2}) is compatible with Hecke operators.
We will use general facts about cohomological correspondences gathered in Appendix \ref{app:corr}.

We then recall the main result from \cite{MorelSiegel1}, that is the computation
of the ``intersection cohomology groups'' \(H^{\bullet}((\Acal_{n,K})_{\Fpbar},
\IC^K(V))\), for \(K\) an open compact subgroup of \(\GSpbf_{2n}(\A_f)\) of the form
\(K^p \times \GSpbf_{2n}(\Zp)\), as a Hecke and Galois module, in terms of the
stabilization of the trace formula for certain elliptic endoscopic groups of
\(\mathbf{GSp}_{2n}\).
Combined with results of the previous section, in the case \(K =
\GSpbf_{2n}(\Zhat)\) we obtain an endoscopic formula for these modules in terms
of Arthur's substitute parameters for \(\mathbf{Sp}_{2n}\) and the lifted Satake
parameters of Proposition \ref{pro:lift_stab_mult} (an unconditional
reformulation of Kottwitz' conjecture in \cite{Kottwitz_AA}).

We conclude this section with Corollary \ref{cor:sigma_IH_crys}, essentially saying that the Galois modules considered above have crystalline semisimplification, which will be used in Galois representation-theoretic arguments in the next section.
To this end we prove Corollary \ref{coro:IH_vs_Hc_Q} relating intersection and compactly supported cohomology.
This relation will be further simplified (in both directions) in Sections \ref{sec:IH_from_Hc} and \ref{sec:Hc_from_IH}.

\subsection{Siegel modular varieties}
\label{sec:def_An}

First we recall definitions from \cite{MorelSiegel1}.
Fix an integer \(n \geq 0\) and a free \(\Z\)-module \(\Lambda\) of rank \(2n\) endowed
with a non-degenerate (over \(\Z\)) alternate bilinear form \(\langle \cdot, \cdot
\rangle\).
Let \(\Gbf\) be the associated general symplectic group (see Section
\ref{sec:not_red_gps}), a reductive group over \(\Z\).

For \(M \geq 3\) and integer we consider the functor \(\Acal_{n,M}\) from the
category of schemes over \(\Z[1/M]\) to the category of sets defined as follows:
\(\Acal_{n,M}(S)\) is the set of isomorphism classes of quadruples \((A, \lambda,
\eta, c)\) where \(A\) is an abelian scheme of constant dimension \(n\), \(\lambda\) is
a principal polarization of \(A\), \(\eta : \underline{\Lambda /M \Lambda}_S \simeq
A[M] \) is an isomorphism of finite étale commutative group schemes and \(c :
\underline{\Z / M \Z}_S \simeq \mu_{M,S}\) are such that the following diagram
commutes:
\begin{equation} \label{eq:diag_levelM}
  \begin{tikzcd}[column sep=6em]
    \underline{\Lambda / M \Lambda \times \Lambda / M \Lambda}_S \arrow[d,
      "{\eta \times \eta}"] \arrow[r, "{\langle \cdot, \cdot \rangle}"] &
      \underline{\Z / M \Z}_S \arrow[d, "{c}"] \\
    A[M] \times A[M] \arrow[r] & \mu_{M,S}
  \end{tikzcd}
\end{equation}
where the bottom horizontal map is the Weil pairing induced by \(\lambda\).
If \(n>0\) then \(c\) is determined by \(\eta\) and is thus redundant.
This functor is representable by a smooth quasi-projective scheme over \(\Z[1/M]\)
(see \cite{Mumford_GIT}, \cite[Corollary 1.4.1.12]{Lan_cpctif_PEL}) which we
still denote \(\Acal_{n,M}\).
There is a natural free action of \(\Gbf(\Z/M\Z)\) on the right of \(\Acal_{n,M}\),
by precomposition of \(\eta\).
For \(M' \geq M\) divisible by \(M\) the forgetful functor \(\Acal_{n,M'} \rightarrow
\Acal_{n,M} \times_{\Z[1/M]} \Z[1/M']\) is finite étale, Galois with Galois group
naturally identified to \(K(M) / K(M')\) where \(K(M) := \ker \Gbf(\Zhat)
\rightarrow \Gbf(\Z / M \Z)\).

More generally, if \(K\) is an open compact subgroup of \(\Gbf(\A_f)\), there is a
finite set of primes \(S\) such that \(K = K_S \prod_{p \not\in S} \Gbf(\Zp)\), one
can consider an analogous moduli problem (as a category fibered in groupoids
over the category of schemes over \(\Z[1/M]\), where \(M = \prod_{p \in S} p\)),
which is a smooth Deligne-Mumford stack over \(\Z[1/M]\) (see \cite[Theorem
1.4.1.11]{Lan_cpctif_PEL}).
If \(K\) is neat (in the sense of \cite[\S 0]{Pink_dissertation}), as is the case
if \(K = K(M)\) with \(M \geq 3\), then \(\Acal_{n, K}\) (here we suppress \(S\) from
the notation) is a quasi-projective scheme over \(\Z[1/M]\) (see \cite[Corollary
7.2.3.10]{Lan_cpctif_PEL}, as well as \cite[\S 5]{KottwitzPoints}).
The action of \(\Gbf(\Z/M\Z)\) on \(\Acal_{n,M}\) and the forgetful functor
considered above also generalize, as follows.
Suppose that \(K', K\) are open compact subgroups of \(\Gbf(\A_f)\), choose \(S\) as
above suitable for both \(K\) and \(K'\) and suppose that \(g \in \prod_{p \in S}
\Gbf(\Qp) \times \prod_{p \not\in S} \Gbf(\Zp)\) is such that \(K' \subset g K
g^{-1}\).
Then there is a finite étale morphism \(T_{K', g, K} : \Acal_{n, K'} \rightarrow
\Acal_{n,K}\) of constant degree \(|gKg^{-1}/K'|\), which only depends on \(g\) via
\(K'gK\) (see \cite[\S 1.1]{MorelSiegel1}, \cite[\S 6]{KottwitzPoints}; note that
the construction of \(T_{K',g,K}\) is more natural once the moduli problem is
reformulated in a more flexible manner using quasi-isogenies: see \cite[\S 1.3.8
and \S 1.4.3]{Lan_cpctif_PEL}).
If \(K'\) is normal in \(K\) then \(T_{K', 1, K}\) is a Galois cover which realizes
\(\Acal_{n,K}\) as the quotient of \(\Acal_{n,K'}\) by \(K/K'\).
We have \(T_{K', g, K} \circ T_{K'', h, K'} = T_{K'', hg, K}\) when this makes
sense.

In the present paper we are particularly interested in the ``level one'' Siegel
modular variety \(\Acal_n := \Acal_{n, \Gbf(\Zhat)}\), a smooth Deligne-Mumford
stack over \(\Z\).
There is no need to rely on the general theory of stacks however, since
everything can be formulated using \(\Gbf(\Z/M\Z)\)-equivariant objects on
\(\Acal_{n, M}\), varying \(M\) in a finite set so that \(\Spec \Z[1/M]\) cover \(\Spec
\Z\).

For the purpose of introducing the cocharacter \(\mu\) below, let us recall the
usual description of the orbifold \(\Acal_{n,K}(\C)\) as a double quotient.
For simplicity we only recall the case of a principal level structure \(K(M)\).
It is well-known that a principally polarized abelian variety \((A, \lambda)\)
over \(\C\) is canonically determined by
\begin{itemize}
  \item \(\Lie A\), a vector space over \(\C\) of dimension \(n\),
  \item \(\Gamma := H_1(A(\C), \Z)\), a lattice in \(\Lie A\) such that \(A(\C)
    \simeq \Lie A / \Gamma\), and
  \item \(B : \Gamma \times \Gamma \rightarrow \Z(1) := 2 \sqrt{-1} \pi \Z = \ker
    (\exp)\) an alternating perfect pairing (coming from the polarization
    \(\lambda\)) such that the pairing \((v,w) \mapsto \sqrt{-1} B(v, \sqrt{-1} w)\)
    on the real vector space \(\Gamma_\R = \Lie A\) is symmetric positive definite
    (this condition does not depend on the choice of \(\sqrt{-1}\)).
\end{itemize}
There exists \(\eta_0 : \Lambda \simeq \Gamma\) and \(c_0 : \Z \simeq \Z(1)\) making
the following diagram commute.
\begin{equation}
  \begin{tikzcd}[column sep=6em]
    \Lambda \times \Lambda \arrow[d, "{\eta_0 \times \eta_0}"] \arrow[r,
      "{\langle \cdot, \cdot \rangle}"] & \Z \arrow[d, "{c_0}"] \\
    \Gamma \times \Gamma \arrow[r, "{B}"] & \Z(1)
  \end{tikzcd}
\end{equation}
Clearly \((\eta_0, c_0)\) is unique up to the action of \(\Gbf(\Z)\) (by
precomposition on \(\eta_0\) and multiplication by the similitude factor on
\(c_0\)), and \(c_0\) is determined by \(\eta_0\) if \(n>0\).
Transporting the complex structure on \(\Lie A = \Gamma_{\R}\) via \(\eta_0\) we get
an \(\R\)-algebra morphism \(\ul{h} : \C \rightarrow \End(\Lambda_{\R})\) such that
for \(z \in \C^\times\) we have \(h(z) := (\ul{h}(z), |z|^2) \in \Gbf(\R)\).
It turns out that \(h\) is induced by a (unique) algebraic morphism \(\Res_{\C/\R}
( \GLbf_{1,\C} ) \rightarrow \Gbf_{\R}\) that we abusively still denote \(h\).
The \(\Gbf(\R)\)-conjugacy class of \(h\) does not depend on \((A, \lambda)\).
Denote \(\Xcal\) the set of \(\Gbf(\R)\)-conjugates of \(h\), a hermitian symmetric
space with two connected components \footnote{One can identify each component
with the usual Siegel upper-half space.}.

We now give a direct (without referring to abelian varieties) construction of an element \(h_0 \in \Xcal\).
Choose \((J, 1) \in \Gbf(\R)\) such that \(J^2 = -1\) and the symmetric bilinear
form \((v,w) \mapsto \langle v, J w \rangle\) on \(\Lambda_{\R}\) is either positive
definite or negative definite.
Choose \(i \in \C\) such that \(i^2 = -1\), thereby giving \(\Lambda_{\R}\) a complex
structure: \(i\) acts by \(J\).
Define \(h_0(z) \in \Gbf(\R)\) to be \((\text{multiplication by } z, |z|^2)\).
The complex vector space \(\Lambda_{\R}\) is equipped with a hermitian form \(H :
(v,w) \longmapsto \langle v, J w \rangle - i \langle v, w \rangle\) which is
either positive definite or negative definite.
Choose a decomposition of the Hermitian space \((\Lambda_{\R}, H)\) as a direct
orthogonal sum of \(n\) lines \((L_j)_{1 \leq j \leq n}\).
The common stabilizer \(\Tbf\) in \(\Gbf_{\R}\) of these lines is a maximal torus of
\(\Gbf_\R\) which is anisotropic modulo center.
Explicitly,
\begin{equation} \label{eq:desc_anis_tor_GSp}
  \mathbf{T}(\R) \simeq \big\{ \ul{t} = ((t_1, \dots, t_n), s) \in (\C^{\times})^n \times \R^\times \ \big|\ |t_1|^2 = \dots = |t_n|^2 = s \big\}
\end{equation}
where \(\ul{t}\) stabilizes \(L_j\) and acts on it by multiplication by \(t_j\).
Then \(h_0\) factors through \(\Tbf\) and in these coordinates \(h_0(z) = ((z, \dots, z), |z|^2)\).
It is easy to check that \(h_0\) belongs to \(\Xcal\), for example using a principally polarized abelian variety which is a product of elliptic curves.
Let \(\Bbf\) be a Borel subgroup of \(\GSpbf_{2n, \C}\) containing \(\Tbf_{\C}\) such that \((\Bbf, \Tbf)\) corresponds to the \emph{generic} discrete series representations of \(\Gbf(\R)\) \footnote{We recalled Harish-Chandra's parametrization of discrete series representations  in \cite[\S 4.2.1]{Taibi_dimtrace}.}.
We now compute \(\mu_{h_0}\) with respect to such a Borel pair, as this will be useful in Section \ref{sec:desc_IH_lift}.
Recall that \(\Bbf\) is unique up to conjugation by the normalizer of \(\Tbf\) in \(\Gbf(\R)\).
By \cite[p.\ 315]{Taibi_dimtrace} one can choose \(\Bbf\) such that the
corresponding simple roots are the following characters of \(\Tbf\), as
characterized by their value on \(\Tbf(\R)\):
\[ \ul{t} \mapsto t_1/\ol{t_2},\ \ \ul{t} \mapsto \ol{t_2}/t_3,\ \ \dots,\ \ \ul{t} \mapsto c^{n-2}(t_{n-1}) / c^{n-1}(t_n),\ \ \ul{t} \mapsto (t_n/\ol{t_n})^{(-1)^{n+1}}, \]
where \(c\) denotes complex conjugation.
There exists \(g \in \Gbf(\C)\) conjugating \((\Tbf_{\C}, \Bbf)\) into
\((\Tbf_{\GSpbf_{2n}, \C}, \Bbf_{\GSpbf_{2n}, \C})\).
Via conjugation by \(g\), the above set of roots corresponds to \(\alpha_1, \dots,
\alpha_n\), in this order.
Therefore, in the parametrization of \(\Tbf_{\GSpbf_{2n}}\) introduced in Section \ref{sec:not_red_gps}, the element \(g \ul{t} g^{-1}\) of \(\Tbf_{\GSpbf_{2n}}(\C)\) is \((\ol{t_1}^{-1}, t_2^{-1}, \dots, c^n(t_n)^{-1}, s)\).
In particular the cocharacter \(\Ad(g) \circ \mu_{h_0}: \GLbf_{1,\C} \to
\Tbf_{\GSpbf_{2n}, \C}\) is characterized by the relations
\[ \langle \mu_{h_0}, \alpha_i \rangle = (-1)^{i-1} \text{ for } 1 \leq i \leq n
  \text{ and } \langle \mu_{h_0}, \nu \rangle = 1. \]
If we see \(\mu_{h_0}\) as a character of \(\Tcal_{\GSpin_{2n+1}}\), with the
parametrization of \(\Tcal_{\GSpin_{2n+1}}\) introduced in Section
\ref{sec:not_red_gps}, it is thus equal to
\begin{equation} \label{eq:comp_muh_gen}
  (z_1, \dots, z_n, s, \lambda) \longmapsto s \lambda \prod_{i \text{ even}} z_i^{-1}.
\end{equation}

Returning to moduli problems, if \(\eta : \Lambda / M \Lambda \simeq A[M] =
M^{-1} \Gamma / \Gamma\) is a level structure for the principally polarized
abelian variety \((A, \lambda)\), there exists \(g \in \Gbf(\Z / M \Z)\) such that
\(\eta = (\eta_0 \mod M) \circ g\).
This gives an identification
\[ \Acal_{n,M}(\C) \simeq \Gbf(\Z) \backslash \left( \Xcal \times \Gbf(\Z / M
  \Z) \right). \]
Using \(\Gbf(\Z / M\Z) \simeq \Gbf(\Zhat) / K(M)\), \(\Gbf(\Z) = \Gbf(\Q) \cap
\Gbf(\Zhat)\) and \(\Gbf(\A_f) = \Gbf(\Q) \Gbf(\Zhat)\) (which follows from the
analogous equality for \(\GLbf_1\) and strong approximation for \(\Gbf_{\der}\)), we
finally have the identification \(\Acal_{n,M}(\C) \simeq \Gbf(\Q) \backslash
(\Xcal \times \Gbf(\A_f) / K(M))\).
This reformulation amounts to considering more generally \(\eta_0 : \Lambda_\Q
\simeq \Gamma_\Q\) and \(c_0: \Q \simeq \Q(1)\), and as above it is better suited
to generalization to arbitrary (possibly non-principal) level structures: we
have identifications \(\Acal_{n,K}(\C) \simeq \Gbf(\Q) \backslash (\Xcal \times
\Gbf(\A_f) / K)\), and it is easy to check that via these identifications the map
induced by \(T_{K',g,K}\) on complex points is simply right multiplication by \(g\).

\subsection{Automorphic local systems}
\label{sec:local_systems}

Let \(S\) be a finite set of prime numbers, \(K = K_S \times \prod_{p \not\in S} \Gbf(\Zp)\) a compact open subgroup of \(\Gbf(\A_f)\) and as in the previous section denote \(M = \prod_{p \in S} p\).
Denote \(\pr_{\ell}: \Gbf(\A_f) \rightarrow \Gbf(\Qell)\) the projection map.
Recall from \cite[\S 2.1]{MorelSiegel1} (a construction going back at least to \cite[\S 3]{Langlands_LK_GL2}) that there is a natural functor \(V \rightsquigarrow \Fcal^K(V)\) from the category of finite-dimensional algebraic representations of \(\Gbf(\Qell)\) to the category of \(\ell\)-adic local systems on \(\Acal_{n, K} \times_{\Z[1/M]} \Z[1/\ell]\).
In fact it is first defined on finitely generated \(\Zell\)-modules endowed with a continuous action of \(\pr_{\ell}(K)\).
The local systems \(\Fcal^K(V)\) are of geometric origin, for example for \(V_{\Std}\) the standard representation of \(\Gbf_{\Qell}\) one can check that \(\Fcal^K(V_{\Std})\) is isomorphic to the relative \(\ell\)-adic Tate module of the universal abelian variety \(\pi : A^\univ \rightarrow \Acal_{n,K}\), and so \(\Fcal^K(V_{\Std}^*) = R^1 \pi_* \ul{\Qell}\).
More generally, for \(V\) an irreducible representation with highest weight \((\lambda_1, \dots, \lambda_n, m)\) (for the parametrization introduced in Section \ref{sec:not_red_gps}), a shift of the local system \(\Fcal^K(V)\) can be cut out inside \(R (\pi^{\times s})_* \ul{\Qell}(m)\) using algebraic correspondences, where \(s = \lambda_1 + \dots + \lambda_n\) and \(\pi^{\times s} : A^{\univ} \times_{\Acal_{n,K}} \dots \times_{\Acal_{n,K}} A^{\univ} \rightarrow \Acal_{n,K}\) (see \cite[p.\ 235]{FaltingsChai}).
A bit more precisely \(\Fcal^K(V)\) is a summand of
\[ \Sym^{\lambda_1 - \lambda_2} \left( \bigwedge^1 R^1 \pi_* \ul{\Qell} \right) \otimes \dots \otimes \Sym^{\lambda_{n-1} - \lambda_n} \left( \bigwedge^{n-1} R^1 \pi_* \ul{\Qell} \right) \otimes \left( \bigwedge^n R^1 \pi_* \ul{\Qell} \right)^{\otimes \lambda_n} (m). \]
In particular (see \cite[Lemma 5.6.6]{Pink_ladic_Shim}) if an algebraic representation \(V\) of \(\Gbf_{\Qell}\) has a central character (e.g.\ if \(V\) is irreducible), say \(z \mapsto z^{-w}\), then for any prime \(p\) not dividing \(M \ell\), the local system \(\Fcal^K(V)\) over \((\Acal_{n,K})_{\Fp}\) is pure of weight \(w\) (equal to \(\sum_i \lambda_i - 2m\) for \(V\) as above).

Over \(\C\), for \(V\) an algebraic representation of \(\Gbf_\Q\) we have a
(topological) local system
\[ \Gbf(\Q) \backslash \left( V \times \Xcal \times \Gbf(\A_f) / K \right)
\rightarrow \Gbf(\Q) \backslash \left( \Xcal \times \Gbf(\A_f) / K \right)
\simeq \Acal_{n,K}(\C). \]
If \(\hat{L}\) is a \(\Zhat\)-lattice in \(\A_f \otimes_{\Q} V\) which is stable under
\(K\), then \(\hat{L}\) is determined by the \(\Z\)-lattice \(L = V \cap \hat{L}\) in
\(V\).
Similarly for \(h \in \Gbf(\A_f)\) consider \(h L := V \cap h \hat{L}\).
There is a local system in finite free \(\Z\)-modules
\[ \Gbf(\Q) \backslash \left( \bigsqcup_{h \in \Gbf(\A_f) / K} h L \times \Xcal
\times \Gbf(\A_f) / K \right) \]
over \(\Acal_{n,K}(\C)\), and extending scalars from \(\Z\) to \(\Z/\ell^N \Z\),
algebraising using \cite[Exposé XI Théorème 4.4]{SGA4-3}, considering the
projective system as \(N\) varies and inverting \(\ell\), one recovers
\(\Fcal^K(V_{\Qell})\).

\subsection{Hecke correspondences}
\label{sec:Hecke_corr}

We recall the definition of Hecke correspondences on the local systems \(\Fcal^K(V)\), which induce a Hecke action on ordinary and compactly supported cohomology.
In the next section these correspondences will be extended to intermediate extensions to the minimal compactification of \(\Acal_{n, K}\) (over a field), and for this purpose we introduce a formalism of Hecke operators.
To simplify the formulation of ulterior statements it is also convenient to recall the relation with smooth representations (Proposition \ref{pro:eq_Hecke_smrep} below).

\begin{defi} \label{def:hecke_cat}
  Let \(\Bcal\) be a preadditive category.
  Let \(G\) be a locally profinite group.
  Let \(C\) be a coinitial\footnote{This means that for any compact open subgroup \(K\) of \(G\) there exists \(K' \in C\) which is contained in \(K\).} set of compact open subgroups of \(G\), stable under conjugation by elements of \(G\) and under finite intersections.
  Let \(\Hecke(G, C, \Bcal)\) be the category of families \((V_K)_{K \in C}\) of objects of \(\Bcal\) endowed with morphisms \([K_2, g, K_1, K'] : V_{K_1} \rightarrow V_{K_2}\) defined for \(K_1, K_2, K' \in C\) and \(g \in G\) such that \(K' \subset g K_1 g^{-1} \cap K_2\), satisfying the following conditions
  \begin{enumerate}
    \item For \(K_1, K_2, K' \in C\) and \(g\) such that \(K' \subset g K_1 g^{-1} \cap K_2\), for any \(h_1 \in K_1\) and \(h_2 \in K_2\), we have \([K_2, h_2 g h_1, K_1, h_2 K' h_2^{-1}] = [K_2, g, K_1, K']\).
    \item For any \(K \in C\) we have \([K, 1, K, K] = \id_{V_K}\).
    \item For \(K_1, K_2, K',K'' \in C\) and \(g \in G\) such that \(K'' \subset K' \subset g K_1 g^{-1} \cap K_2\) we have
      \[ [K_2, g, K_1, K''] = |K'/K''| \times [K_2, g, K_1, K']. \]
    \item For \(K_1, K_2, K_3, K', K'' \in C\) and \(g_1, g_2 \in G\) such that \(K' \subset g_1 K_1 g_1^{-1} \cap K_2\) and \(K'' \subset g_2 K_2 g_2^{-1} \cap K_3\) we have
      \[ [K_3, g_2, K_2, K''] \circ [K_2, g_1, K_1, K'] = \sum_{[h] \in K'' \backslash g_2 K_2 / K'} [K_3, hg_1, K_1, K'' \cap h K' h^{-1}]. \]
      (The right-hand side is well defined thanks to the first axiom.)
  \end{enumerate}
  A morphism from \((V_K)_{K \in C}\) to \((V'_K)_{K \in C}\) is a family of morphisms \(V_K \to V'_K\) (in \(\Bcal\)) intertwining the \([K_2, g, K_1, K']\)'s.
  If \(F\) is a field we simply denote \(\Hecke(G, C, F)\) for \(\Hecke(G, C, \Bcal)\) where \(\Bcal\) is the category of vector spaces over \(F\).
\end{defi}

The third property implies that we could equivalently only specify \([K_2, g, K_1, K']\) when \(K' = g K_1 g^{-1} \cap K_2\), but this would make the last expression less natural.

\begin{prop} \label{pro:eq_Hecke_smrep}
  Assume that \(F\) is a field of characteristic zero.
  Let \(G\) and \(C\) be as in Definition \ref{def:hecke_cat}.
  Let \(\Rep_{\sm}(G, F)\) be the category of smooth representations of \(G\) with coefficients in \(F\).
  The following functors are equivalence of categories between \(\Hecke(G, C, F)\) and \(\Rep_{\sm}(G, F)\) which are inverse of each other (up to isomorphism of functors):
  \begin{enumerate}
  \item To a smooth representation \(V\) of \(G\) over \(F\) associate \((V^K)_K\) and \([K_2, g, K_1, K'] = \sum_{k \in K_2/K'} kg\).
  \item To \(((V_K)_K, ([K_2, g, K_1, K'])_{K_1, K_2, g, K'}) \in \Ob \Hecke(G, C, F)\) associate \(V = \varinjlim_K V_K\) for the transition morphisms \([K', 1, K, K'] : V_K \rightarrow V_{K'}\) when \(K' \subset K\).
    The action of \(G\) is induced by \([gKg^{-1}, g, K, gKg^{-1}] : V_K \rightarrow V_{gKg^{-1}}\).
  \end{enumerate}
\end{prop}
\begin{proof}
  We omit the straightforward verification that the first functor is well-defined.

  Let us check that the second functor is well-defined.
  For \(K_1, K_2 \in C\) and \(g \in G\) satisfying \(K_2 \subset g K_1 g^{-1}\) we may consider \([K_2, g, K_1, K_2]: V_{K_1} \to V_{K_2}\).
  Using the fourth axiom in Definition \ref{def:hecke_cat} we see that this subset of operators is compatible with composition: for \(K_1,K_2,K_3 \in C\) and \(g_1,g_2 \in G\) satisfying \(K_2 \subset g_1 K_1 g_1^{-1}\) and \(K_3 \subset g_2 K_2 g_2^{-1}\) we have
  \[ [K_3, g_2, K_2, K_3] \circ [K_2, g_1, K_1, K_2] = [K_3, g_2 g_1, K_1, K_3]. \]
  It follows that each \(g \in G\) defines an operator on \(V := \varinjlim_K V_K\) and that the resulting map \(G \to \End_F(V)\) is multiplicative, and thanks to the second axiom we have a linear action of \(G\) on \(V\).
  We know from the first axiom that for \(K_1,K_2 \in C\) satisfying \(K_2 \subset K_1\) and \(g \in K_1\) we have \([K_2,g,K_1,K_2] = [K_2,1,K_1,K_2]\), and so this action is smooth.

  If we start from a smooth representation \(V\) of \(G\) then the natural map \(\varinjlim_K V^K \to V\) is an isomorphism of \(F\)-vector spaces, and it is tautologically compatible with the action of \(G\).

  The least formal part of the proof is the remaining direction: starting from an object \(((V_K)_{K \in C}, ([K_2,g,K_1,K'])_{K_2,g,K_1,K'})\) of \(\Hecke(G,F)\), we want to identify (in a natural way) \(\left( \varinjlim_{K'} V_{K'} \right)^K\) with \(V_K\) and check that via these identifications we have, for all \(K_1,K_2,K' \in C\) and \(g \in G\) satisfying \(K' \subset K_2 \cap gK_1g^{-1}\):
  \[ [K_2,g,K_1,K'] = \sum_{k \in K_2/K'} kg. \]
  For \(K\) and \(K'\) in \(C\) satisfying \(K' \subset K\) we compute
  \[ [K,1,K',K'] \circ [K',1,K,K'] = \sum_{[h] \in K' \backslash K' / K'} [K,h,K,K'] = [K,1,K,K'] = |K/K'| \id_{V_K} \]
  using the fourth and third axiom in Definition \ref{def:hecke_cat}, and conclude that \([K',1,K,K']\) is injective.
  If moreover \(K'\) is a normal subgroup of \(K\) then using the fourth axiom again we see that \([K',1,K,K']\) maps \(V_K\) to the subspace of \(K/K'\)-invariants in \(V_{K'}\), and we compute
  \[ [K',1,K,K'] \circ [K,1,K',K'] = \sum_{h \in K/K'} [K',h,K',K'] \]
  from which we deduce that the image of \(V_K\) in \(V_{K'}\) contains the subspace of \(K/K'\)-invariants in \(V_{K'}\).
  For a given \(K \in C\), any \(K' \in C\) satisfying \(K' \subset K\) contains \(K'' \in C\) which is a normal subgroup of \(K\) (namely \(\bigcap_{g \in K/K'} g K' g^{-1}\)), and so we conclude from the above that the natural map \(V_K \to \varinjlim_{K'} V_{K'}\) identifies \(V_K\) with the subspace of \(K\)-invariants.
  Finally for \(K_1,K_2,K' \in C\) and \(g \in G\) satisfying \(K' \subset gK_1g^{-1} \cap K_2\), letting \(K'' = \bigcap_{[h] \in K_2/K'} h K' h^{-1}\) we have thanks to the fourth axiom
  \[ [K'',1,K_2,K''] \circ [K_2,g,K_1,K'] = \sum_{[h] \in K_2/K'} [K'',hg,K_1,K'']. \]
\end{proof}

\begin{rema}
  In order to determine a smooth representation of \(G\) over a field of characteristic zero, it would be enough to only specify embeddings \([K', 1, K, K']\) for \(K' \subset K\) as well as actions \([g K g^{-1}, g, K, gKg^{-1}]\) satisfying natural relations (left to the reader).
  When the spaces \(V_K\) come from pro-étale torsors, as in Proposition \ref{pro:Hecke_corr_sat_formalism} below, this can be easier than verifying all axioms in Definition \ref{def:hecke_cat}.
  We will have to consider situations where it is not so obvious that we have \emph{embeddings} \([K', 1, K, K']\) and that \(V_K\) may be identified with the subspace of \(K\)-invariant vectors, so that the above formalism (or some analogue) is needed.

  The formalism in Definition \ref{def:hecke_cat} could also prove useful when dealing with integral (or positive characteristic) coefficients (not needed in this paper).
\end{rema}

\begin{coro} \label{cor:hecke_contra}
  Let \(G\) and \(C\) be as in Definition \ref{def:hecke_cat}.
  Assume that \(F\) is a field of characteristic zero.
  If \(((V_K)_K, ([K_2, g, K_1, K'])_{K_2, g, K_1, K'})\) is an object of \(\Hecke(G, C, F)\), then we have a contragredient object \(((V_K^*)_K, ([K_2, g, K_1, K']^*)_{K_2, g, K_1, K'})\) of \(\Hecke(G, F)\) defined by \(V_K^* = \Hom_F(V_K, F)\) and \([K_2, g, K_1, K']^*: V_{K_1}^* \to V_{K_2}^*\) equal to \(|K_2/K'| / |K_1/g^{-1} K' g|\) times the transpose of \([K_1, g^{-1}, K_2, g^{-1} K' g]: V_{K_2} \to V_{K_1}\).
\end{coro}
\begin{proof}
  By the previous proposition we have identifications between \(V_K\) and \(V^K\) for a smooth representation \((V, \pi)\) of \(G\) over \(F\).
  Let \((\wt{V}, \wt{\pi})\) be the contragredient representation, and \(((\wt{V}^K)_K, ([K_2, g, K_1, K']^*)_{K_2, g, K_1, K'})\) the object of \(\Hecke(G, F)\) associated by the previous proposition.
  By elementary group theory the restriction morphism \(\wt{V}^K \to \Hom_F(V^K, F)\) is an isomorphism.
  For any compact open subgroups \(K_1, K_2, K'\) of \(G\) and any \(g \in G\) satisfying \(K' \subset K_2 \cap g K_1 g^{-1}\), for any \(v \in V^{K_2}\) and \(\wt{v} \in \wt{V}^{K_1}\) we have
  \begin{align*}
    \langle v, [K_2, g, K_1, K']^* \wt{v} \rangle
    &= \sum_{k \in K_2/K'} \langle v, \wt{\pi}(kg) \wt{v} \rangle \\
    &= |K_2/K'| \langle \pi(g^{-1}) v, \wt{v} \rangle \\
    &= \frac{|K_2/K'|}{|K_1/g^{-1} K' g|}
      \langle [K_1, g^{-1}, K_2, g^{-1} K' g] v, \wt{v} \rangle.
  \end{align*}
\end{proof}

We now recall the definition of Hecke correspondences on \(\Acal_{n, K}\).
The definition below is a minor generalization of \cite[Définition
5.2.1]{MorelSiegel1}.
For \(g \in \Gbf(\Qell)\) simply denote \(V \mapsto g V\) the equivalence of
categories from continuous representations of \(K\) (resp.\ algebraic
representations of \(\Gbf_{\Qell}\)) to continuous representations of \(g K g^{-1}\)
(resp.\ algebraic representations of \(\Gbf_{\Qell}\)) where \(g V\) is \(V\) but with
the action composed with \(\ad (g^{-1})\).
If \(V\) is a representation of \(\Gbf(\Qell)\) (or if \(V\) is a representation of
\(K\) and \(g \in K\)) the action of \(g\) induces an isomorphism of representations
\(i_g: g(V|_K) \simeq V|_{gKg^{-1}}\).

Suppose that \(S\) be a finite set of prime numbers, \(K = K_S \times \prod_{p
\not\in S} \Gbf(\Zp)\) and \(K' = K'_S \times \prod_{p \not\in S} \Gbf(\Zp)\)
compact open subgroups of \(\Gbf(\A_f)\), and \(g \in \prod_{p \in S} \Gbf(\Qp)
\times \prod_{p \not\in S} \Gbf(\Zp)\) such that \(K' \subset g K g^{-1}\).
Staring at the definitions gives us an isomorphism of functors from the category
of continuous representations of \(K\) (be it on finite \(\Z/\ell^N\Z\)-modules,
\(\Zell\)-modules or \(\Qell\)-vector spaces) to the category of suitable local
systems on \(\Acal_{n,K'}\):
\[ T_{K', g, K}^* \circ \Fcal^K \simeq \Fcal^{K'} \circ (g_{\ell} - ). \]
On representations of \(\Gbf(\Qell)\) (and not just an open subgroup which might
not contain \(g_{\ell}\)), composing with the isomorphism \(i_{g_{\ell}}\) recalled
above, we get an isomorphism of functors
\begin{equation} \label{eq:iso_func_pullback}
  T_{K', g, K}^* \circ \Fcal^K \simeq \Fcal^{K'},
\end{equation}
which remains unchanged if \(g\) is multiplied on the right by an element of \(K\).
The isomorphisms \eqref{eq:iso_func_pullback} are compatible with composition
but we refrain from naming them and explicitly writing the formula satisfied
whenever \(K' \subset g K g^{-1}\) and \(K'' \subset g' K' g''^{-1}\).

\begin{defi} \label{def:Hecke_corr}
  Suppose that \(S\) be a finite set of prime numbers, \(K_1 = K_{1,S} \times
  \prod_{p \not\in S} \Gbf(\Zp)\) and \(K_2 = K_{2,S} \times \prod_{p \not\in S}
  \Gbf(\Zp)\) compact open subgroups of \(\Gbf(\A_f)\), \(g \in \prod_{p \in S}
  \Gbf(\Qp) \times \prod_{p \not\in S} \Gbf(\Zp)\) and \(K' = K'_S \times \prod_{p
  \not\in S} \Gbf(\Zp)\) and open compact subgroup of \(\Gbf(\A_f)\) contained in
  \(K_2 \cap g K_1 g^{-1}\).
  Let \(u(K_2, g, K_1, K')\) be the cohomological correspondence (in the sense of
  \cite[Exposé III \S 3.2]{SGA5}) from \(\Fcal^{K_1}(V)\) to \(\Fcal^{K_2}(V)\) with
  support in \((T_{K',g,K_1}, T_{K',1,K_2})\) obtained by composing
  identifications:
  \[ T_{K', g, K_1}^* \Fcal^{K_1}(V) \simeq \Fcal^{K'}(g_{\ell} V) \simeq
    \Fcal^{K'}(V) \simeq T_{K', 1, K_2}^* \Fcal^{K_2}(V) \simeq T_{K', 1, K_2}^!
    \Fcal^{K_2}(V). \]
  At the last step \(T_{K',1,K_2}^!\) is identified to \(T_{K',1,K_2}^*\) because
  \(T_{K',1,K_2}\) is étale \cite[Exposé XVIII Proposition 3.1.8 p.91]{SGA4-3}.
\end{defi}

It is easy to check that for \(h_1 \in K_1\) and \(h_2 \in K_2\), \(u(K_2, h_2 g h_1,
K_1, K')\) is isomorphic to \(u(K_2, g, K_1, h_2 K' h_2^{-1})\).
It is also easy to check that for \(h \in \prod_{p \in S} \Gbf(\Qp) \times
\prod_{p \not\in S} \Gbf(\Zp)\) and \(K' \subset g K_1 g^{-1} \cap h K_2 h^{-1}\)
the correspondence obtained as above but using \(T_{K', h, K_2}\) instead of
\(T_{K',1,K_2}\) is simply isomorphic to \(u(K_2, h^{-1}g, K_1, h^{-1} K' h)\).

Using this one easily checks that the dual correspondence \(\Dbb(u(K_2, g, K_1,
K'))\) from \(\Dbb(\Fcal^{K_2}(V)) \simeq \Fcal^{K_2}(V^*)(n(n+1)/2)[n(n+1)]\) (since
\(\Acal_{n,K}\) is smooth of relative dimension \(n(n+1)/2\) over \(\Z[1/M]\)) to
\(\Fcal^{K_1}(V^*)(n(n+1)/2)[n(n+1)]\) and with support in \((T_{K',1,K_2},
T_{K',g,K_1})\) is isomorphic to \(u(K_1, g^{-1}, K_2, g^{-1} K'
g)(n(n+1)/2)[n(n+1)]\).

Denote \(\pi_i : \Acal_{n,K_i} \rightarrow \Spec \Z[1/M]\).
Since \(T_{K',1,K_2}\) (resp.\ \(T_{K',g,K_1}\)) is proper, \(u(K_2, g, K_1, K')\)
induces
\begin{align} \label{eq:corr_induces_coho}
  u(K_2, g, K_1, K')_* : \pi_{1*} \Fcal^{K_1}(V) & \longrightarrow \pi_{2*}
    \Fcal^{K_2}(V), \\
  \text{resp.\ } u(K_2, g, K_1, K')_! : \pi_{1!} \Fcal^{K_1}(V) & \longrightarrow
    \pi_{2!} \Fcal^{K_2}(V),
\end{align}
see \cite[(1.3.2)]{Fujiwara_Deligne_conj} and \cite[(1.3)]{Pink_localterms},
also recalled in Section \ref{sec:corr_def}.
These two operations are dual to each other, so that
\[ \Dbb(u(K_2, g, K_1, K')_*) = u(K_1, g^{-1}, K_2, g^{-1} K'
  g)_! (n(n+1)/2)[n(n+1)]. \]

\begin{rema} \label{rema:Hecke_central_char}
  If \(V\) is irreducible with central character \(\chi\) then it is easy to check
  that for \(K_1, K_2, g, K'\) as above and \(z\) a central element of \(\Gbf(\Q)\)
  such that \(z_f \in \prod_{p \in S} \Gbf(\Qp) \times \prod_{p \not\in S}
  \Gbf(\Zp)\) we have \(u(K_2, z_f g, K_1, K') = \chi(z_{\ell}) u(K_2, g, K_1,
  K')\).
\end{rema}

\begin{rema} \label{rema:twisting_sim}
  Denote by \(\Qell(\nu)\) the algebraic representation of \(\Gbf_{\Qell}\) on
  \(\Qell\) given by the similitude character \(\nu\).
  There is a canonical isomorphism \(c^K : \Fcal^K(\Qell(\nu)) \simeq
  \ul{\Qell}(1)\) (in principal level \(M \geq 3\) it is given by the morphism \(c\)
  in diagram \eqref{eq:diag_levelM} in level \(M \ell^N\) for varying \(N\)).
  Unwinding the definitions (which is easier using the reformulation of the
  moduli problem using quasi-isogenies) we find that in the setting of
  Definition \ref{def:Hecke_corr} we have a commutative diagram
  \[ \begin{tikzcd}[column sep=6em]
      T_{K',g,K}^* \Fcal^K(\Qell(\nu)) \arrow[d, "{T_{K',g,K}^*(c^K)}"]
        \arrow[r, "{\sim}"] & \Fcal^{K'}(\Qell(\nu)) \arrow[d, "{c^{K'}}"] \\
      T_{K',g,K}^* \ul{\Qell}(1) = \ul{\Qell}(1) \arrow[r, "{|\nu(g)|_f^{-1}}"]
        & \ul{\Qell}(1)
  \end{tikzcd} \]
  where the top horizontal arrow is the composition of the first two
  isomorphisms in Definition \ref{def:Hecke_corr}.

  As a consequence for any algebraic representation \(V\) of \(\Gbf_{\Qell}\) and
  any integer \(m\) we have canonical isomorphisms \(\Fcal^K(V(\nu^m)) \simeq
  \Fcal^K(V)(m)\) and commutative diagrams
  \[ \begin{tikzcd}[column sep=10em]
      T_{K',g,K_1}^* \Fcal^{K_1}(V(\nu^m)) \arrow[d, "{\sim}"] \arrow[r, "{u(K_2,g,K_1,K')}"] & T_{K',1,K_2}^* \Fcal^{K_2}(V(\nu^m)) \arrow[d, "{\sim}"] \\
      T_{K',g,K_1}^* \Fcal^{K_1}(V)(m) \arrow[r, "{|\nu(g)|_f^{-m} u(K_2,g,K_1,K')}"] & T_{K',1,K_2}^* \Fcal^{K_2}(V)(m)
  \end{tikzcd} \]
\end{rema}

Of course for \(i : \Spec F \hookrightarrow \Spec \Z[1/M]\) where \(F\) is a prime
field one can similarly define cohomological correspondences between the local
systems \(\Fcal^K(V)\) pulled back to \((\Acal_{n,K})_F\), and these will be denoted
\(u(K_2,g,K_1,K')_F\).
Note that if \(F\) is finite then duality intertwines \(i^*\) and \(i^!\) but \(i^!
\Fcal^K(V) \simeq i^* \Fcal^K(V)(-1)[-2]\) by absolute purity
\cite{Fujiwara_abs_purity}, since \(\Fcal^K(V)\) is a local system on a smooth
scheme over \(\Z[1/M]\).

The following proposition is well-known, and is included to prepare for the case
of intersection complexes.

\begin{prop} \label{pro:Hecke_corr_sat_formalism}
  Let \(F\) be \(\Q\) (resp.\ \(\Fp\) for some prime number \(p \neq \ell\)).
  Let \(V\) be an algebraic representation of \(\Gbf(\Qell)\).
  Then for any \(0 \leq i \leq n(n+1)/2\) the families of finite-dimensional
  \(\Qell\)-vector spaces \((H^i_c((\Acal_{n,K})_{\ol{F}}, \Fcal^K(V)))_K\) and
  \((H^i((\Acal_{n,K})_{\ol{F}}, \Fcal^K(V)))_K\), where \(K\) varies in the set of
  neat compact open subgroups of \(\Gbf(\A_f)\) (resp.\ neat compact open
  subgroups of \(\Gbf(\A_f)\) of the form \(\Gbf(\Zp) \times K^p\)), equipped with
  the operators
  \[ (u(K_2,g,K_1,K')_F)_! : H^i_c((\Acal_{n,K_1})_{\ol{F}}, \Fcal^{K_1}(V))
  \rightarrow H^i_c((\Acal_{n,K_2})_{\ol{F}}, \Fcal^{K_2}(V)) \]
  \[ (u(K_2,g,K_1,K')_F)_* : H^i((\Acal_{n,K_1})_{\ol{F}}, \Fcal^{K_1}(V))
  \rightarrow H^i((\Acal_{n,K_2})_{\ol{F}}, \Fcal^{K_2}(V)) \]
  satisfy the axioms of Definition \ref{def:hecke_cat}.
\end{prop}
\begin{proof}
  The first two axioms follow directly from the definition.
  Let us explain the dependence on \(K'\) (i.e.\ the third axiom) using
  pushforward and pullback of correspondences, recalled in section
  \ref{sec:corr_push_pull}.
  For \(K'' \subset K'\) we simply have
  \[ \corr (T_{K'', 1, K'})^* u(K_2, g, K_1, K') = u(K_2, g, K_1, K'') \]
  essentially because the isomorphism of functors \(T_{K',g,K}^* \simeq
  T_{K',g,K}^!\) compose (\cite[Exposé XVIII Proposition 3.1.8(iii)]{SGA4-3} and
  \cite[Exposé XVII Théorème 6.2.3 (Var 3)]{SGA4-3}), as do the isomorphism
  of functors \(T_{K', g, K}^* \Fcal^K \simeq \Fcal^{K'}\).
  By Lemma \ref{lemm:fet_push_pull_corr} we have
  \[ \corr (T_{K'', 1, K'})_* u(K_2, g, K_1, K'') = |K'/K''| \times u(K_2, g,
  K_1, K'). \]
  This implies the third axiom in Definition \ref{def:hecke_cat}.

  Let us check the fourth axiom, i.e.\ composition.
  Suppose we have a diagram
  \[ \begin{tikzcd}
      & \Acal_{n,K'} \arrow[dl, "{T_{g_1}}" above] \arrow[dr, "{T_1}"] & &
        \Acal_{n,K''} \arrow[dl, "{T_{g_2}}" above] \arrow[dr, "{T_1}"] \\
      \Acal_{n,K_1} & & \Acal_{n,K_2} & & \Acal_{n,K_3}
  \end{tikzcd} \]
  with \(K' \subset g_1 K_1 g_1^{-1} \cap K_2\) and \(K'' \subset  g_2 K_2 g_2^{-1}
  \cap K_3\).
  Then we have an identification of \(\Acal_{n,K'} \times_{\Acal_{n,K_2}}
  \Acal_{n,K''}\) with \(\bigsqcup_{h \in K'' \backslash g_2K_2 / K'} \Acal_{n,K''
  \cap h K' h^{-1}}\) via the morphisms
  \[ \begin{tikzcd}
      & \Acal_{n, K'' \cap h K' h^{-1}} \arrow[dl, "{T_h}" above] \arrow[dr,
      "{T_1}"] \\
      \Acal_{n,K'} & & \Acal_{n,K''}
  \end{tikzcd} \]
  and using this identification, the equality
  \[ u(K_3, g_2, K_2, K'') \circ u(K_2, g_1, K_1, K') = \sum_{g_3 \in K''
    \backslash g_2 K_2 / K'} u(K_3, g_3 g_1, K_1, K'' \cap g_3 K' g_3^{-1}) \]
  easily follows from the definition.
\end{proof}

In particular we get an admissible representation of \(G = \Gbf(\A_f)\) (resp.\ \(\Gbf(\A_f^p)\)) on
\[ H^{\bullet}_?((\Acal_n)_{\ol{F}}, \Fcal(V)) := \varinjlim_K H^{\bullet}_?((\Acal_{n,K})_{\ol{F}}, \Fcal^K(V)) \]
for \(? \in \{c, \emptyset\}\) and \(F = \Q\) (resp.\ \(\Fp\)), with a commuting continuous action of \(\Gal_F\) (continuous in the sense that for any compact open subgroup \(K\) the action on the space of \(K\)-invariants is continuous).
For any choice of Haar measure \(\mathrm{vol}\) on \(G\) such that any compact open subgroup has rational volume there are canonical \(G \times \Gal_F\)-equivariant pairings
\begin{equation} \label{eq:duality_Hecke}
  H^i_c \left( (\Acal_n)_{\ol{F}}, \Fcal(V) \right) \times H^{n(n+1)/2-i} \left( (\Acal_n)_{\ol{F}}, \Fcal(V^*) \right) \rightarrow \Qell \left( -\frac{n(n+1)}{2} \right)
\end{equation}
obtained by multiplying the usual pairing (``Poincaré duality'') in level \(K\) by \(\mathrm{vol}(K)\).
For each \(i\) the pairing \eqref{eq:duality_Hecke} identifies these two admissible representations of \(G\) to the (\(\Qell(-\frac{n(n+1)}{2})\)-valued) contragredient of each other (see Corollary \ref{cor:hecke_contra}).

\subsection{Minimal compactifications and intermediate extensions}
\label{sec:minimal_cpctif}

If \(K = K_S \prod_{p \not\in S} \Gbf(\Zp)\) is neat (as before \(S\) a finite set
of primes, \(K_S\) a compact open subgroup of \(\prod_{p \in S} \Gbf(\Qp)\) and we
denote \(M = \prod_{p \in S} p\)) Chai and Faltings \cite[Theorem
V.2.5]{FaltingsChai} (see also \cite[Section 7]{Lan_cpctif_PEL}) constructed the
minimal compactification \(\Acal_{n,K}^*\) of \(\Acal_{n,K}\) over \(\Z[1/M]\), using
toroidal compactifications.
More precisely, \(\Acal_{n,K}^*\) is a normal projective scheme over \(\Z[1/M]\)
with an open embedding \(j : \Acal_{n,K} \hookrightarrow \Acal_{n,K}^*\).
In general \(\Acal_{n,K}^*\) is not smooth over (any point of) \(\Z[1/M]\).
At least if \(K\) is a principal level, there is a stratification of
\(\Acal_{n,K}^* \smallsetminus \Acal_{n,K}\) by schemes isomorphic to
\(\Acal_{n',K'}\) for \(n'<n\) (see \cite[Theorem V.2.5]{FaltingsChai}), but we
shall not need this description.

If \(K' = K_S' \prod_{p \not\in S} \Gbf(\Zp)\) (\(K_S'\) a compact open subgroup of
\(\prod_{p \in S} \Gbf(\Qp)\)), and \(g \in \prod_{p \in S} \Gbf(\Qp) \times
\prod_{p \not\in S} \Gbf(\Zp)\) is such that \(K' \subset gKg^{-1}\),
\cite[Proposition 7.2.5.1]{Lan_cpctif_PEL} gives a canonical extension of
\(T_{K', g, K}\) as \(\ol{T}_{K',g,K} : \Acal_{n,K'}^* \rightarrow \Acal_{n,K}^*\).
The map \(\ol{T}_{K',g,K}\) is finite (this follows from \cite[Corollary
7.2.5.2]{Lan_cpctif_PEL}), but not necessarily étale.
Since it is canonical it satisfies natural properties similar to \(T_{K',g,K}\):
it only depends on \(g\) via \(K'gK\), and is compatible with composition.

\begin{defi} \label{def:IC}
  For \(K\) as above, \(F = \Q\) or \(\Fp\) for \(p \nmid M\), and \(V\) an algebraic
  representation of \(\Gbf(\Qell)\) define the intersection complex
  \[ \IC^K_\ell(V)_F = j_{!*} ( \Fcal^K(V)_F[n(n+1)/2] ) [-n(n+1)/2] \in D^b_c((\Acal_{n,K})_F, \Qell). \]
\end{defi}
We will also use lighter notation \(\IC^K_\ell(V)\) or \(\IC^K(V)\) when there is no risk of confusion.

If \(F = \Fp\) and \(V\) is irreducible then as recalled in Section \ref{sec:local_systems} \(\Fcal^K(V)_{\Fp}\) is pure of weight determined by the central character of \(V\) and so by \cite[Corollaire 5.4.3]{BBD} \(\IC^K(V)_{\Fp}\) is pure of the same weight.
In this setting Morel identified \(\IC^K(V)_{\Fp}\) with the \emph{weight
truncation} of \(j_* \Fcal^K(V)_{\Fp}\) \cite[Théorème 3.1.4]{MorelSiegel1}.
Using this identification she canonically extended \cite[\S 5]{MorelSiegel1} the
Hecke correspondences \(u(K_2,g,K_1,K')_{\Fp}\) (with \(K_1,K_2,K'\) containing
\(\Gbf(\Zp)\) and \(g_p \in \Gbf(\Zp)\)) of Definition \ref{def:Hecke_corr} (or
rather their base change to \(\Fp\)) to Hecke correspondences from
\(\IC^{K_1}_{\Fp}(V)\) to \(\IC^{K_2}_{\Fp}(V)\), with support in
\(((\ol{T}_{K',g,K_1})_{\Fp}, (\ol{T}_{K',g,K_1})_{\Fp})\).

By Lemma \ref{lemm:uniq_int_ext_corr} (2) in both cases \(F = \Q\) or \(\Fp\) we also have canonical (``geometric'') extensions of Hecke correspondences between the intersection complexes \(\IC^K_F(V)\), that we denote \(u(K_2, g, K_1, K')^{\IC}_F\).
In case \(F = \Fp\) by comparing both characterizations it is clear that this coincides with Morel's canonical (``weight-theoretic'') extension of \(u(K_2, g, K_1, K')\).
Taking cohomology, we have Hecke operators
\begin{equation}
  (u(K_2,g,K_1,K')_F^{\IC})_* : H^{\bullet}((\Acal_{n,K_1}^*)_{\ol{F}},
  \IC^{K_1}_F(V) ) \longrightarrow H^{\bullet}((\Acal_{n,K_2}^*)_{\ol{F}},
  \IC^{K_2}_F(V))
\end{equation}
commuting with the action of \(\Gal(\ol{F}/F)\).

\begin{defi} \label{def:hecke_axioms_derived}
  Let \(F\) be \(\Q\) or \(\Fp\).
  Denote \(G = \GSpbf_{2n}(\A_f)\) if \(F = \Q\) (resp.\ \(G = \GSpbf_{2n}(\A_f^{(p)})\) if \(F = \Fp\)).
  Let \(C\) be a coinitial set of compact open subgroups of \(G\), stable under conjugation and under finite intersections.
  Let \(\PreH(\Acal_{n,?,F}^*, C, \Qell)\) be the category of pairs \(((L_K)_{K \in C}, (v(K_2,g,K_1,K'))_{K_2,g,K_1,K'})\) where \(L_K\) is an object of \(D^b_c(\Acal_{n,K,F}^*, \Qell)\) and \(v(K_2,g,K_1,K'): \ol{T_g}^* L_{K_1} \to \ol{T_1}^! L_{K_2}\) is a cohomological correspondence, defined whenever \(K_1,K_2,K' \in C\) and \(g \in G\) satisfy \(K' \subset gK_1g^{-1} \cap K_2\), subject to the following conditions.
  \begin{enumerate}
  \item For \(K_1, K_2, K' \in C\) and \(g \in G\) such that \(K' \subset g K_1 g^{-1} \cap K_2\), for any \(h_1 \in K_1\) and \(h_2 \in K_2\), the correspondences \(v(K_2, h_2 g h_1, K_1, h_2 K' h_2^{-1})\) and \(v(K_2, g, K_1, K')\) are identified via the isomorphism \(\ol{T_{h_2}}: \Acal_{n,h_2 K' h_2^{-1},F}^* \simeq \Acal_{n,K',F}^*\).
  \item For any \(K \in C\) we have \(v(K,1,K,K) = \id\).
  \item For \(K_1,K_2,K',K'' \in C\) and \(g \in G\) satisfying \(K'' \subset K' \subset gK_1g^{-1} \cap K_2\) we have
    \[ \corr (\ol{T}_{K'',1,K'})_* v(K_2,g,K_1,K'') = |K''/K'| \times v(K_2,g,K_1,K'). \]
  \item For \(K_1,K_2,K_3,K',K'' \in C\) and \(g_1, g_2 \in G\) satisfying \(K' \subset g_1K_1g_1^{-1} \cap K_2\) and \(K'' \subset g_2K_2g_2^{-1} \cap K_3\), denoting
    \[ \ol{f}: \bigsqcup_{[h] \in K'' \backslash g_2 K_2 / K'} \Acal_{n, K'' \cap hK'h^{-1}}^* \longrightarrow \Acal_{n,K'}^* \times_{\Acal_{n,K_2}^*} \Acal_{n,K''}^* \]
    the morphism induced by the pairs of morphisms \((\ol{T_{K'' \cap hK'h^{-1}, h, K'}}, \ol{T_{K'' \cap hK'h^{-1}, 1, K''}})\), we have
    \[ \corr \ol{f}_* \left( (v(K_3, h g_1, K_1, K'' \cap h K' h^{-1}))_{h \in K'' \backslash g_2K_2 / K'} \right) = v(K_3, g_2, K_2, K'') \circ v(K_2, g_1, K_1, K'). \]
  \end{enumerate}
\end{defi}

\begin{lemm} \label{lem:pre_hecke_to_hecke}
  Let \(F\) be \(\Q\) (resp.\ \(\Fp\)), \(\ell\) a prime number (resp.\ a prime number different from \(p\)).
  Denote \(G = \GSpbf_{2n}(\A_f)\) if \(F = \Q\) (resp.\ \(G = \GSpbf_{2n}(\A_f^{(p)})\) if \(F = \Fp\)).
  Let \(C\) be the set of neat\footnote{In case \(F=\Fp\) this means that \(G \times \GSpbf_{2n}(\Zp)\) is neat.} compact open subgroups of \(G\).
  For an object \(((L_K)_{K \in C}, (v(K_2,g,K_1,K'))_{K_2,g,K_1,K'})\) of \(\PreH(\Acal_{n,?,F}^*, C, \Qell)\), denoting \(\pi_K: \Acal_{n,K,F}^* \to \Spec F\), the pair
  \[ ((\pi_{K,*} L_K)_K, (v(K_2,g,K_1,K')_*)_{K_2,g,K_1,K'}) \]
  defines an object of \(\Hecke(G, C, D^b_c(\Spec F, \Qell))\).
  In particular (see Proposition \ref{pro:eq_Hecke_smrep}) for any \(i \in \Z\) we have a smooth admissible action of \(G\) on
  \[ \varinjlim_{K \in C} H^i(\Acal_{n,K,\ol{F}}^*, \IC^K(V)) \]
  with a commuting continuous action of \(\Gal(\ol{F}/F)\).
\end{lemm}
\begin{proof}
  The first two axioms in Definition \ref{def:hecke_cat} clearly follow from the corresponding axioms in Definition \ref{def:hecke_axioms_derived}.
  The third and fourth axioms also follow from the corresponding axioms, using the fact (recalled in Section \ref{sec:corr_push_pull}) that pushforward of correspondences along proper morphisms is compatible with cohomological realizations.
\end{proof}

\begin{prop} \label{pro:IH_is_Hecke}
  Let \(F\) be \(\Q\) (resp.\ \(\Fp\)), \(\ell\) a prime number (resp.\ a prime number different from \(p\)).
  Denote \(G = \GSpbf_{2n}(\A_f)\) if \(F = \Q\) (resp.\ \(G = \GSpbf_{2n}(\A_f^{(p)})\) if \(F = \Fp\)).
  Let \(C\) be the set of neat compact open subgroups of \(G\).
  Let \(V\) be an algebraic representation of \(\GSpbf_{2n,\Qell}\).
  The pair
  \[ \left( (\IC_F^K(V))_{K \in C}, (u(K_2,g,K_1,K')_F^{\IC})_{K_2,g,K_1,K'} \right) \]
  defines an object of \(\PreH(\Acal^*_{n,?,F}, C, \Qell)\).
  In particular the family of finite-dimensional \(\Z\)-graded \(\Qell\)-vector spaces \(\left(H^{\bullet}((\Acal_{n,K}^*)_{\ol{F}}, \IC^{K}(V) )\right)_K\) with the Hecke operators \((u(K_2,g,K_1,K')_F^{\IC})_*\) defines a \(\Z\)-graded admissible representation of \(G\) with commuting continuous action of \(\Gal(\ol{F}/F)\).
\end{prop}
\begin{proof}
  The first two axioms in Definition \ref{def:hecke_axioms_derived} follow from the corresponding relations for correspondences \(u(K_2,g,K_1,K')\) on \(\Acal_{n,K,F}\) and injectivity in Lemma \ref{lemm:uniq_int_ext_corr} (2).

  For the third axiom we know thanks to Proposition \ref{pro:corr_push_pull_compat} and Corollary \ref{cor:compat_pull_corr} that the restriction of \(\corr (\ol{T}_{K'',1,K'})_* u(K_2,g,K_1,K'')_F^{\IC}\) to open Shimura varieties is equal to \(\corr (T_{K'',1,K'})_* u(K_2,g,K_1,K'')_F\), which as we saw in the proof of Proposition \ref{pro:Hecke_corr_sat_formalism} is equal to
  \[ |K'/K''| \times u(K_2,g,K_1,K')_F. \]
  By Lemma \ref{lemm:uniq_int_ext_corr} (2) we deduce
  \[ \corr (\ol{T}_{K'',1,K'})_* u(K_2,g,K_1,K'')_F^{\IC} = |K'/K''| \times u(K_2,g,K_1,K')_F^{\IC}. \]

  The fourth axiom is proved similarly: in the setting of this axiom, the restriction of \(\ol{f}\) to open Shimura varieties is an isomorphism
  \[ f: \bigsqcup_{[h] \in K'' \backslash g_2 K_2 / K'} \Acal_{n, K'' \cap hK'h^{-1}} \longrightarrow \Acal_{n,K'} \times_{\Acal_{n,K_2}} \Acal_{n,K''}. \]
  and we saw in the proof of Proposition \ref{pro:Hecke_corr_sat_formalism} that we have
  \[ \corr f_* \left( (u(K_3, h g_1, K_1, K'' \cap h K' h^{-1})_F)_{h \in K'' \backslash g_2K_2 / K'} \right) = u(K_3, g_2, K_2, K'')_F \circ u(K_2, g_1, K_1, K')_F. \]
  By compatibility of pullback of correspondences (in the case at hand, restriction to open Shimura varieties) with pushforward (Proposition \ref{pro:corr_push_pull_compat}) and composition (Proposition \ref{pro:compat_pullback_compo_corr}) we deduce that the correspondences
  \[ \corr \ul{f}_* \left( (u(K_3, h g_1, K_1, K'' \cap h K' h^{-1})_F^{\IC})_{h \in K'' \backslash g_2K_2 / K'} \right) \]
  \[ \text{and} \ \ u(K_3, g_2, K_2, K'')_F^{\IC} \circ u(K_2, g_1, K_1, K')_F^{\IC}. \]
  agree on open Shimura varieties, and thanks to Lemma \ref{lemm:uniq_int_ext_corr} (2) we conclude that they are equal.
\end{proof}

As in the case of ordinary or compactly supported cohomology we have a natural identification of the \(\Qell(-n(n+1)/2)\)-valued contragredient of \(\varinjlim_K H^i( (\Acal_{n,K}^*)_{\ol{F}}, \IC^K_F(V))\) with \(\varinjlim_K H^{n(n+1)/2-i}( (\Acal_{n,K}^*)_{\ol{F}}, \IC^K_F(V^*))\).

\subsection{Hecke and Galois actions over \(\Q\) and \(\Fp\)}

In this section we recall the ``specialization theorem'' for intersection cohomology, which follows from the existence of smooth (toroidal) compactifications (defined over \(\Zp\)) over the minimal one, which implies compatibility between intermediate extensions and nearby cycles (and so vanishing of vanishing cycles).
This is a special case of \cite[Corollaire 4.4]{Stroh_Kottwitz}.
We also check compatibility with Hecke correspondences.

We mostly use the notation introduced in \cite[Exposé XIII]{SGA7-2} for nearby cycles (in our case over the base \(\Spec \Zp\)), see Sections \ref{sec:corr_nearby_cycles} and \ref{sec:corr_nearby_cycles_perverse}.

\begin{prop} \label{pro:spe_gal_hecke_coho}
  For \(V\) an algebraic representation of \(\Gbf(\Qell)\), \(p \neq \ell\) and \(K\) a
  neat compact open subgroup of \(\Gbf(\A)\) such that \(K = K^p \times \Gbf(\Zp)\),
  the \(\Gal_{\Qp}\)-equivariant specialization morphism
  \[ H^{\bullet}_c((\Acal_{n,K})_{\Fpbar}, \Fcal^K(V)) \rightarrow
  H^{\bullet}_c((\Acal_{n,K})_{\Qbar}, \Fcal^K(V)) \]
  is also \(\Hcal(\Gbf(\A_f^p) // K^p)\)-equivariant and is an isomorphism.
  In particular, the representation of \(\GalQ\) on the right-hand side is
  unramified at \(p\).
\end{prop}
\begin{proof}
  This is a consequence of \cite[Corollaire 4.3]{Stroh_Kottwitz} and similar to
  Corollaire 4.6 loc.\ cit., with extra Hecke action.
  Thanks to the existence of toroidal compactifications, smooth extensions of
  the universal principally polarized abelian variety and the formalism of
  plethysms to construct irreducible representations of \(\Gbf\), for \(j\) the open
  immersion of \(\Acal_{n,K}\) into \(\Acal_{n,K}^*\), the base change morphisms
  \(j_{s,!} \Psi_{\eta} \Fcal^K(V)_{\eta} \rightarrow \Psi_{\eta} j_{\eta, !}
  \Fcal^K(V)_{\eta}\) and \(\Psi_{\eta} j_{\eta,*} \Fcal^K(V)_{\eta} \rightarrow
  j_{s,*} \Psi_{\eta} \Fcal^K(V)_{\eta}\) are isomorphisms: see \cite[Proposition
  4.3 and Corollaire 4.3]{Stroh_Kottwitz}.
  Since \(\Fcal^K(V)\) is a local system on \(\Acal_{n, K}\) which is smooth over
  \(\Z_{(p)}\), there is a canonical isomorphism \(\mathrm{sp}^* \Fcal^K(V)_s
  \xrightarrow{\sim} \Psi_{\eta} \Fcal^K(V)_{\eta}\) (\cite[Exposé XIII,
  Reformulation 2.1.5]{SGA7-2}).
  By proper base change (2.1.7.1 and 2.1.8.3 loc.\ cit.) for \(? \in \{*, !\}\) we
  have a canonical isomorphism \(H^{\bullet}((\Acal_{n,K}^*)_{\sbar},
  \Psi_{\eta} j_{\eta,?} \Fcal^K(V)_{\eta}) \simeq
  H^{\bullet}((\Acal_{n,K}^*)_{\etabar}, j_{\eta,?} \Fcal^K(V))\).
  The upshot is that we have canonical \(\Gal(\etabar/\eta)\)-equivariant
  isomorphisms
  \begin{align*}
    H^{\bullet}((\Acal_{n,K})_{\etabar}, \Fcal^K(V)) & \xrightarrow{\sim}
      H^{\bullet}((\Acal_{n,K})_{\sbar}, \Psi_{\eta} \Fcal^K(V))
      \xleftarrow{\sim} H^{\bullet}((\Acal_{n,K})_{\sbar}, \Fcal^K(V)) \\
    H^{\bullet}_c((\Acal_{n,K})_{\sbar}, \Fcal^K(V)) & \xleftarrow{\sim}
       H^{\bullet}_c((\Acal_{n,K})_{\sbar}, \Psi_{\eta} \Fcal^K(V))
       \xleftarrow{\sim} H^{\bullet}_c((\Acal_{n,K})_{\etabar}, \Fcal^K(V))
  \end{align*}
  which are dual to each other up to replacing on one side \(\Fcal^K(V)\) with its
  dual \(\Fcal^K(V^*)(d)[2d]\) where \(d = n(n+1)/2\) (compatibility of nearby
  cycles with duality was proved in \cite[Théorème 4.2]{Illusie_monodromie}).
  Moreover these isomorphisms are compatible with the action of Hecke operators:
  for neat \(K_1 = K_1^p \Gbf(\Zp)\) and \(K_2 = K_2^p \Gbf(\Zp)\), \(g \in
  \Gbf(\A_f^p)\) and \(K' = K'^p \Gbf(\Zp)\) with \(K'^p \subset K_2^p \cap g K_1^p
  g^{-1}\) we have a morphism of correspondences with support in
  \(((T_{K',g,K_1})_s, (T_{K',1,K_2})_s)\)
  \[ \mathrm{sp}^* \left(\Fcal^{K_1}(V)_s, \Fcal^{K_2}(V)_s, u(K_2, g, K_1,
  K')_s \right) \to \Psi_\eta \left(\Fcal^{K_1}(V)_{\eta}, \Fcal^{K_2}(V)_{\Qp},
  u(K_2, g, K_1, K')_{\eta} \right) \]
  showing that the right square of the following diagram is commutative.
  \[ \begin{tikzcd} 
      H^{\bullet}((\Acal_{n,K_1})_{\etabar}, \Fcal^{K_1}(V)) \arrow[r, "{\sim}"]
      \arrow[d] & H^{\bullet}((\Acal_{n,K_1})_{\sbar}, \Psi_{\eta}
      \Fcal^{K_1}(V)) \arrow[d] & \arrow[l, "{\sim}"]
      H^{\bullet}((\Acal_{n,K_1})_{\sbar}, \Fcal^{K_1}(V)) \arrow[d] \\
      H^{\bullet}((\Acal_{n,K_2})_{\etabar}, \Fcal^{K_2}(V)) \arrow[r, "{\sim}"]
      & H^{\bullet}((\Acal_{n,K_2})_{\sbar}, \Psi_{\eta} \Fcal^{K_2}(V)) &
      \arrow[l, "{\sim}"] H^{\bullet}((\Acal_{n,K_2})_{\sbar}, \Fcal^{K_2}(V))
    \end{tikzcd} \]
  By Proposition \ref{pro:compat_coh_corr_nearby} the left square is also
  commutative.
  An almost identical argument can be used to conclude in the case of compactly
  supported cohomology.
\end{proof}

\begin{prop} \label{pro:spe_gal_hecke_IC}
  For \(V\) an algebraic representation of \(\Gbf(\Qell)\), \(p \neq \ell\) and \(K\) a
  neat compact open subgroup of \(\Gbf(\A)\) such that \(K = K^p \times \Gbf(\Zp)\),
  there are canonical isomorphisms
  \[ H^{\bullet}((\Acal_{n,K}^*)_{\Qbar}, \IC^K_{\Q}(V)) \simeq H^{\bullet}((\Acal_{n,K}^*)_{\Fpbar}, \IC^K_{\Fp}(V)) \]
  compatible with the action of \(\Hcal(\Gbf(\A_f^p) // K^p)\), Galois actions and duality.
  In particular, the Galois action on the left-hand side is unramified at \(p\).
\end{prop}
\begin{proof}
  The construction of the isomorphism is a special case of \cite[Corollaire
  4.6]{Stroh_Kottwitz}.
  We need to prove that this isomorphism is Hecke-equivariant and compatible
  with duality on both sides.

  By \cite[Corollaire 4.4]{Stroh_Kottwitz}, for \(j: \Acal_{n,K} \hookrightarrow
  \Acal_{n,K}^*\) there is a canonical isomorphism \(\mathrm{sp}^* j_{s, !*}
  (\Fcal^K(V)_s[n(n+1)/2]) \simeq \Psi_{\eta} j_{\eta, !*}
  (\Fcal^K(V)_{\eta}[n(n+1)/2])\) sitting in a commutative diagram (using the
  isomorphism \(\mathrm{sp}^* \Fcal^K(V)_s \xrightarrow{\sim} \Psi_{\eta}
  \Fcal^K(V)_{\eta}\) which was already used in the previous proof):
  \[ \begin{tikzcd}[column sep=tiny]
      \Psi_{\eta} j_{\eta, !} \Fcal^K(V)_{\eta}[n(n+1)/2] \arrow[r] &
        \Psi_{\eta} j_{\eta, !*} \Fcal^K(V)_{\eta}[n(n+1)/2] \arrow[d, no head,
        "{\sim}"] \arrow[r] & \Psi_{\eta} j_{\eta, *}
        \Fcal^K(V)_{\eta}[n(n+1)/2] \arrow[d, "{\sim}"] \\
      j_{s,!} \Psi_{\eta} \Fcal^K(V)_{\eta}[n(n+1)/2] \arrow[r] \arrow[u,
        "{\sim}"] & j_{s,!*} \Psi_{\eta} \Fcal^K(V)_{\eta}[n(n+1)/2]
        \arrow[r] & j_{s,*} \Psi_{\eta} \Fcal^K(V)_{\eta}[n(n+1)/2]
    \end{tikzcd} \]
  This isomorphism is characterized by the fact that its restriction to
  \((\Acal_{n,K})_s\) is the identity (implicitly using \(j^* \Psi_{\eta}
  \xrightarrow{\sim} \Psi_{\eta} j^*\), \cite[Exposé XIII, 2.1.7.2]{SGA7-2}).
  On the one hand, as recalled above by Lemma \ref{lemm:uniq_int_ext_corr}
  (2) the family of perverse sheaves \(\IC^K_\eta(V) = (j_{\eta, !*}
  (\Fcal^K(V)_{\eta}[n(n+1)/2]))_K\) (for varying \(K\) as in the proposition) is
  equipped with canonical Hecke correspondences with support in
  \(((\overline{T}_{K', g, K_1})_{\eta}, (\ol{T}_{K',1,K_2})_{\eta})\) of Section
  \ref{sec:minimal_cpctif}, and as explained in \ref{sec:corr_nearby} they
  induce correspondences between the \(\Psi_\eta \IC^K_\eta(V)\).
  On the other hand the family of perverse sheaves \((\IC^K_s(V))_K\) is also
  equipped with canonical Hecke correspondences with support in
  \(((\overline{T}_{K', g, K_1})_s, (\ol{T}_{K',1,K_2})_s)\), which induce
  correspondences between the \(\mathrm{sp}^* \IC^K_s(V)\).
  Via the isomorphisms \(\mathrm{sp}^* \IC^K_s(V) \simeq \Psi_\eta \IC^K_\eta(V)\)
  recalled above, these two families of correspondences coincide by Lemma
  \ref{lem:nearby_int_ext_corr} and injectivity in Lemma
  \ref{lemm:uniq_int_ext_corr} (2) and its analogue Lemma
  \ref{lemm:uniq_int_ext_corr2}.

  Compatibility with duality is rather formal: taking the dual of the above
  diagram yields the same diagram for \(\Fcal^K(V^*)(d)[2d]\).
  Details for correspondences are omitted.
\end{proof}

\subsection{Intersection cohomology: Morel's stabilized formula}
\label{sec:IH_Morel}

The following formula was conjectured for arbitrary Shimura varieties by
Kottwitz in \cite{Kottwitz_AA} and proved for Siegel modular varieties by Morel in
\cite[Corollaire 5.3.3]{MorelSiegel2}.

\begin{theo} \label{t:ICMorel}
  Let \(M \geq 3\) be an integer, so that the principal level \(K = K(M)\) is neat.
  Let \(p\) be a prime number which does not divide \(M \ell\).
  Let \(V\) be an irreducible algebraic representation of \(\Gbf = \GSpbf_{2n}\) and
  let \(\chi^{-1}\) be the restriction of central character of \(V\) to
  \(\Abf_{\Gbf}(\R)^0\).
  For any \(f^{\infty} \in \Hcal(\Gbf(\A_f) // K)\), trivial at \(p\), for any large
  enough integer \(j\),
  \begin{equation} \label{eq:ICMorel}
    \mathrm{Tr} \left( \Frob_p^j f^{\infty} \, | \, H^\bullet((\Acal^*_{n,
    M})_{\Fpbar}, \IC^K_{\Fp}(V)) \right) = \sum_{\mathfrak{e} = (\Hbf, \Hcal,
    s, \xi)} \iota(\mathfrak{e}) S_{\disc}^{\Hbf, \Abf_{\Gbf}(\R)^0,
    \chi}(f_{\Hbf}^{(j)}) 
  \end{equation}
  where the sum is over isomorphism classes of elliptic endoscopic data \(\efrak\)
  such that \(\Hbf_{\R}/\Abf_{\Gbf_{\R}}\) has discrete series and \(\Hbf_{\Qp}\) is
  unramified, \(S_{\disc}^{\Hbf, \Abf_{\Gbf}(\R)^0, \chi}\) is the stable linear
  form appearing in the stabilization of Arthur's invariant trace formula for
  \(\Hbf\) (see Section \ref{sec:non_ss_gps}), and the definition of
  \(f_{\Hbf}^{(j)}\) is recalled below.
\end{theo}

\begin{rema} \label{rem:rationality}
  Before explaining the right-hand side we recall that the left-hand side, which
  a priori is an element of \(\Qell\), is rational and independent of \(\ell\).
  These facts are a by-product of the proof: in the case of compactly supported
  cohomology it is visible on the formula proved in \cite{KottwitzPoints}
  \footnote{in fact it is already visible at the first step of the proof, a
  trace formula of Lefschetz-type}, and the case of intersection cohomology is
  visible on the formula proved in \cite[Théorème 1.2.1]{MorelSiegel2}.

  By a well-known argument concerning linear recurrence sequences (see \cite[\S
  2]{Fujiwara_indep_l_IH}) the left-hand side is rational for arbitrary values
  of \(j \in \Z\).
  Combined with purity to separate cohomological degrees (see after Definition
  \ref{def:IC} and \cite[Proposition 6.2.6]{Deligne_Weil2}) this implies that
  the Hecke action on \(H^\bullet((\Acal_{n,M}^*)_{\Fpbar}, \IC^K_{\Fp}(V))\) is
  defined over \(\Qbar\).
  The slightly stronger statement that this action is defined over \(\Q\) can also
  be proved in a manner similar to \cite[\S 3.5]{Clozel_AA}, by considering the
  analogues of \(\IC^K_{\Q}(V)\) over \(\Acal_{n,M}^*(\C)\) (which are defined over
  \(\Q\), as we recalled above) and comparison results summarized in \cite[\S
  6.1.2]{BBD}.

  This rationality property is related to Proposition \ref{pro:rat_Sat}, as will
  be apparent in Theorem \ref{theo:IH_explicit_crude}.
\end{rema}

The distribution \(f_{\Hbf}^{(j)} \in SI(\Hbf, \Abf_{\Hbf}(\R)^0, \chi)\) is
defined as a product \(f_{\Hbf}^{\infty, p} f_{\Hbf, \infty} f_{\Hbf, p}^{(j)}\),
each term being defined using the theory of endoscopy.
Note that we have a canonical identification \(\Abf_{\Gbf} \simeq \Abf_{\Hbf}\)
because \(\efrak\) is elliptic.
An endoscopic datum for \(\PGSpbf_{2n}\) yields one for \(\GSpbf_{2n}\) and this
induces a bijection between equivalence classes of endoscopic data.
Therefore the set of equivalence classes of elliptic endoscopic data for
\(\GSpbf_{2n}\) is also in bijection with the set of split equivalence classes of
elliptic endoscopic data for \(\Spbf_{2n}\).
In particular for every endoscopic datum occurring in \eqref{eq:ICMorel} the
group \(\Hbf\) is split, and we will use the obvious L-embedding \({}^L \xi: {}^L
\Hbf \to {}^L \GSpbf_{2n}\).
At every place \(v\) of \(\Q\) we will use the Whittaker-normalized transfer factors
\(\Delta_{\lambda}'\) as defined in \cite[(5.5.2)]{KottwitzShelstad2012} (see also
p.\ 178 of \cite{Kottwitz_AA}).
Since \(\Gbf\) is of adjoint type there is a unique \(\Gbf(\Q_v)\)-conjugacy class
of Whittaker data so this is unambiguous.

The most familiar term is \(f_{\Hbf}^{\infty, p} \in SI(\Hbf(\A_f^p))\), which is
the transfer (in the sense of \cite[(7.1)]{Kottwitz_AA}) of \(f^{\infty, p}\) seen
as an element of \(I(\Gbf(\A_f^p))\).

At the real place, the stable orbital integrals of \(f_{\Hbf, \infty} \in
SI(\Hbf_{\R}, \Abf_{\Hbf}(\R)^0, \chi)\) are prescribed by \cite[(7.4) on
p.182]{Kottwitz_AA}.
Qualitatively it is known that the distribution \(f_{\Hbf, \infty}\) can be taken
to be a linear combination of pseudo-coefficients of essentially discrete series
(whose central character coincides with \(\chi\) on \(\Abf_{\Hbf}(\R)^0\)) having
Langlands parameter \(\varphi_{\Hbf}\) such that \({}^L \xi \circ \varphi_{\Hbf}\)
is also discrete.
Kottwitz constructs this function rather explicitly (see \cite[p.\
186]{Kottwitz_AA}) using certain transfer factors \(\Delta_{j, B}\) that are
adapted to the study of endoscopy for discrete series representations of
\(\Gbf(\R)\), but this normalization of transfer factors is not very natural in a
global setting.
We will only need the spectral consequence that Kottwitz draws from these
calculations, namely \cite[Lemma 7.1]{Kottwitz_AA}, and the fact that Shelstad
\cite[Theorem 11.5]{She3} pinned down the spectral transfer factors when these
correspond to the geometric transfer factors given by a Whittaker datum.
For \(\varphi_{\Hbf}\) a tempered parameter for \(\Hbf_{\R}\) such that the
composition with \({}^L \Hbf \rightarrow {}^L \Abf_{\Hbf}\) corresponds to a
character \(\Abf_{\Hbf}(\R) \rightarrow \C^{\times}\) whose restriction to
\(\Abf_{\Hbf}(\R)^0\) is \(\chi\), denote by \(\Lambda_{\varphi_{\Hbf}}\) the
associated linear form on \(SI(\Hbf_{\R}, \Abf_{\Hbf}(\R)^0, \chi)\), i.e.\ the
sum of the traces of all elements in the L-packet associated to
\(\varphi_{\Hbf}\).
By \cite[Lemma 5.3]{Shelstad_characters} the numbers
\(\Lambda_{\varphi_{\Hbf}}(f_{\Hbf, \infty})\) determine \(f_{\Hbf, \infty}\).
Let \(\varphi_{\Gbf,V}: W_{\R} \to {}^L \Gbf\) be the discrete Langlands parameter
corresponding to the L-packet \(\Pi_{\varphi_{\Gbf, V}}(\Gbf_{\R})\) consisting of
all essentially discrete series representations having the same infinitesimal
character and central character as \(V^*\).
If \({}^L \xi \circ \varphi_{\Hbf}\) is not conjugated to \(\varphi_{\Gbf,V}\) then \(\Lambda_{\varphi_{\Hbf}}(f_{\Hbf, \infty})\) vanishes.
In particular \(\Lambda_{\varphi_{\Hbf}}(f_{\Hbf, \infty})\) vanishes if \(\varphi_{\Hbf}\) is not discrete, which implies that \(f_{\Hbf,\infty}\) is the image in \(SI(\Hbf_{\R}, \Abf_{\Hbf}(\R)^0, \chi)\) of a linear combination of pseudo-coefficients of essentially discrete series in \(I(\Hbf_{\R}, \Abf_{\Hbf}(\R)^0, \chi)\).
We now assume the existence of \(g \in \Ghat(\C)\) satisfying \(\Ad(g) \circ {}^L \xi \circ \varphi_{\Hbf} = \varphi_{\Gbf,V}\).
Choose an arbitrary element \(\pi\) of the discrete series L-packet \(\Pi_{\varphi_{\Gbf, V}}(\Gbf_{\R})\), and let \((\Bbf, \Tbf)\) be a corresponding Borel pair (see e.g.\ \cite[\S 4.2.1]{Taibi_dimtrace}).
Choose \(h \in \Xcal\) such that \(h\) factors through \(\Tbf\) (such an \(h\) exists), so that \(\mu_h\) belongs to \(X_*(\Tbf)\).
Let \((\Bcal, \Tcal)\) be the standard Borel pair of \(\Ghat\).
Up to conjugation by \(\Ghat\) we can assume that we have \(\varphi_{\Gbf, V}(\C^\times) \subset \Tcal\) and that the holomorphic part of \(\varphi_{\Gbf, V}|_{\C^\times}\) is dominant with respect to \(\Bcal\).
Via the identification of \(\That\) with \(\Tcal\) determined by \(\Bbf\) and \(\Bcal\) we can see \(\mu_h\) as an algebraic character of \(\Tcal\), and let \(\langle \mu_{\pi}, \cdot \rangle\) be its restriction to \(C_{\varphi_{\Gbf, V}}\).
That this is well-defined in terms of \(\pi\), i.e.\ independent of auxiliary choices (among them \(h\)) is a consequence of \cite[Lemma 5.1]{Kottwitz_AA}.
Finally let \(\langle \pi, \cdot \rangle\) be Shelstad's spectral transfer factor, a character of \(C_{\varphi_{\Gbf, V}} / Z(\Ghat)\) \footnote{This spectral transfer factors depends on the choice of a Whittaker datum, but since \(\Gbf/\Abf_{\Gbf}\) is adjoint there is only one \(\Gbf(\R)\)-orbit of Whittaker data.}.
The character \(\langle \mu_{\pi}, \cdot \rangle \langle \pi, \cdot \rangle\) of \(C_{\varphi_{\Gbf, V}}\) does not depend on the choice of \(\pi\) and its restriction to \(Z(\Ghat)\) coincides with that of \(\mu\), where \(\mu\) is any algebraic character of \(\Tcal\) corresponding to \(\mu_h\) for some \(h \in \Xcal\) (the Weyl orbit of \(\mu\) is well-defined).
To compute the character \(\langle \mu_{\pi}, \cdot \rangle \langle \pi, \cdot \rangle\) one can take for example \(\pi = \pi^\mathrm{gen}\), the unique generic element of the L-packet, for which \(\langle \pi, \cdot \rangle\) is trivial.
In this situation we have \cite[Lemma 7.1]{Kottwitz_AA}
\begin{equation} \label{eq:spectral_char_real}
  \Lambda_{\varphi_{\Hbf}}(f_{\Hbf, \infty}) = (-1)^{q(\Gbf_{\R})} \langle \mu_{\pi}, g \xi(s) g^{-1} \rangle \langle \pi, g \xi(s) g^{-1} \rangle.
\end{equation}

At the \(p\)-adic place, fix a hyperspecial maximal compact subgroup of
\(\Hbf(\Qp)\).
Then \(f_{\Hbf,p}^{(j)}\) can be chosen to be the unique element of the
corresponding unramified Hecke algebra such that for any unramified
representation \(\pi\) of \(\Hbf(\Qp)\) with Satake parameter \(c(\pi)\) (a
semi-simple conjugacy class in \(\Hhat\)), we have
\begin{equation} \label{eq:spectral_char_twisted}
  \Tr( \pi(f_{\Hbf,p}^{(j)}) ) = p^{j n(n+1)/4 } \Tr( r_{- \mu}(s \times
  \xi(c(\pi))^j )) 
\end{equation}
where \(r_{-\mu}\) is the irreducible representation of \(\Ghat\) having extremal
weight \(-\mu\).
As in the real case the stable orbital integrals of the distribution
\(f_{\Hbf,p}^{(j)}\) are prescribed (see \cite[(7.3)]{Kottwitz_AA}), and following
Kottwitz the equivalent spectral characterization
\eqref{eq:spectral_char_twisted} is deduced thanks to:
\begin{itemize}
  \item the twisted fundamental lemma (known in the \(p\)-adic case for the whole
    unramified Hecke algebra and without any assumption on the residual
    characteristic: see \cite{TFLsmallchar1}, \cite{TFLsmallchar2}) for base
    change, seeing \(\Hbf\) as a twisted endoscopic group for \(\Rbf_{\Gbf} :=
    \Res_{\Q_{p^j} / \Qp} \Gbf\) and the automorphism ``arithmetic Frobenius''
    (see \cite[Appendix A]{MorelBook}).
    More precisely, the endoscopic datum \((\Hbf, \Hcal, s, \xi)\) and the choice
    of an unramified L-embedding \({}^L \xi : {}^L \Hbf \rightarrow {}^L \Gbf\)
    extending \(\xi\) determine a twisted endoscopic datum \((\Hbf,
    \widetilde{\Hcal}, \tilde{s}, \tilde{\xi})\) for \(\Rbf_{\Gbf}\) endowed with
    the arithmetic Frobenius automorphism of \(\Q_{p^j} / \Qp\) and an unramified
    L-embedding \({}^L \tilde{\xi} : {}^L \Hbf \rightarrow {}^L \Rbf_{\Gbf}\).
  \item formula \cite[Theorem 5.6.2]{KottwitzShelstad2012} relating twisted (for
    base change) and ordinary transfer factors (to be more precise this formula
    relates the factors \(\Delta_0'\), i.e.\ without epsilon factors, but it is
    easy to check that the epsilon factors are simply equal so \(\Delta_0'\) can
    be replaced by \(\Delta_\lambda'\) in this formula),
  \item the simple formula for the Satake isomorphism in the minuscule case
    \cite[Theorem 2.1.3]{Kottwitz_Satake_minu} and the explicit computation of
    \({}^L \tilde{\xi}\) (see p.\ 179 of \cite{Kottwitz_AA} and \S A.2.6 of
    \cite{MorelBook}).
\end{itemize}

\begin{rema}
  Following \cite[Appendix A]{MorelBook} and \cite[\S 5.6]{KottwitzShelstad2012}
  we are using the transfer factors \(\Delta_\lambda'\), and so the morphism of
  Hecke algebras on p.\ 180 of \cite{Kottwitz_AA} has to be defined via the
  classical (as opposed to ``Deligne'') normalization of the Satake
  isomorphisms, i.e.\ mapping \(p\) to the \emph{arithmetic} Frobenius.
  This is necessary for the twisted fundamental lemma to hold.
  %Similarly on p.\ 193 of \cite{Kottwitz_AA}
  %Similarly in Morel's work, normalizations should change

  Note that the definition of \(c(\pi)\) does not involve a choice of
  normalization because \(\Hbf\) is split, only defining an unramified Langlands
  parameter does.
\end{rema}

In \cite{MorelSiegel2} on the right-hand side of \eqref{eq:ICMorel} the linear
forms \(ST^{\Hbf}\), given by explicit geometric expansions and defined by
Kottwitz in unpublished notes, occur instead of Arthur's \(S_{\disc}^{\Hbf}\).
The goal of Kottwitz's notes is the stabilization of the trace formula for the
action of Hecke operators on middle-weighted cohomology of the locally symmetric
space attached to a reductive group \(\Gbf\) over \(\Q\) such that \(\Gbf /
\Abf_{\Gbf}\) admits discrete series at the real place.
This would be independent of (and more direct and explicit than) Arthur's
stabilization of his invariant trace formula (\cite{ArthurSTF1},
\cite{ArthurSTF2}, \cite{ArthurSTF3}), although the trace formulas of
\cite{GoKoMPh} (in the case of an upper middle weight profile) and
\cite{ArthurL2}, which clearly agree on the geometric (``orbital integrals'')
side, are equivalent by \cite{Nair_weighted}.
Kottwitz' notes are unpublished but in \cite{Peng} Zhifeng Peng used an argument
similar to \cite{ArthurL2} to show that Arthur's stable linear form defined in
his stabilization, when applied to distributions that are pseudo-coefficients of
discrete series at the real place, admits the expansion predicted by Kottwitz,
i.e.\ that \(S_{\disc}^{\Hbf} = ST^{\Hbf}\).

\subsection{Description of intersection cohomology using lifted Satake
parameters}
\label{sec:desc_IH_lift}

We now apply Theorem \ref{t:ICMorel} and Proposition \ref{pro:lift_stab_mult}
to precisely describe the \(\Hcal(\GSpbf_{2n}(\A_f)//K)_{\Qell} \times
\GalQ\)-module structure of \(H^\bullet((\Acal^*_{n, K})_{\Qbar}, \IC^K(V))\) in
the particular case of level one (\(K = \GSpbf_{2n}(\Zhat)\)).
We somewhat abusively define
\[ H^\bullet((\Acal^*_n)_{\Qbar}, \IC(V)) := H^\bullet((\Acal^*_{n,M})_{\Qbar},
  \IC^{K(M)}(V))^{K(M)} \]
where \(M \geq 3\) is any integer.
By Propositions \ref{pro:IH_is_Hecke} and \ref{pro:spe_gal_hecke_IC} this
\(\Hcal_f^{\unr}(\GSpbf_{2n})_{\Qell} \times \GalQ\)-module does not depend on the
choice of \(M\), is unramified away from \(\ell\) and for any prime \(p \neq \ell\)
the semi-simplification of its restriction to \(\otimes'_{q \neq \ell, p}
\Hcal^{\unr}(\GSpbf_{2n, \Z_q})_{\Qell} \times \GalQp\) is determined (abstractly
at least) by Theorem \ref{t:ICMorel}.

We will require the following definition.

\begin{defi} \label{def:spin_pm_psi}
  Let \(m \geq 1\) and \(\taut \in \ICcalt(\SObf_{4m})\).
  Let \(\psi \in \Psit^{\unr, \taut}_{\disc, \nonendo}(\SObf_{4m})\).
  Let \(\spin_\psi^+\) (resp.\ \(\spin_\psi^-\)) be the half-spin representation of \(\Mpsisc\) (see Definition \ref{def:Mpsisc}) such that, if the eigenvalues of \(\tau_\psi\) in the standard representation of \(\Mpsi\) are \(\pm x_1, \dots, \pm x_{2m}\) where \(x_1 > \dots > x_{2m} > 0\) are integers, the eigenvalues (counted with multiplicities) of \(\spin_\psi^+(\tau_\psi)\) (resp.\ \(\spin_\psi^-(\tau_\psi)\)) are the \(\frac{1}{2} \sum_{i=1}^{2m} \epsilon_i x_i\) for \((\epsilon_i)_i \in \{\pm 1\}^{2m}\) such that the cardinality of
  \[ \left\{ i \in \{1,\dots,2m\} \, \middle| \, \epsilon_i = + 1 \right\} \]
  is even (resp.\ odd).
\end{defi}

To be more explicit, in this definition we may take \(\Mpsi = \SO_{4m}\) and \(\tau_\psi\) the class of
\[ (x_1, \dots, x_{2m}) \in \Lie(\Tcal_{\SO_{4m}}) \]
using the parametrization \eqref{eq:param_T_SO_even} of \(\Tcal_{\SO_{4m}}\), and then the weights of the maximal torus \(\Tcal_{\Spin_{4m}}\) of \(\Mpsisc\) occurring in the representation \(\spin_\psi^+\) (resp.\ \(\spin_\psi^-\)) are, using the parametrization \eqref{eq:param_T_Spin_even} of \(\Tcal_{\Spin_{4m}}\),
\[ (z_1, \dots, z_{2m}, s) \longmapsto s \prod_{i \in I} z_i^{-1} \]
for all subsets \(I\) of \(\{1,\dots,2m\}\) having even (resp.\ odd) cardinality.

A simplification particular to the level one case is that we have
\(H^\bullet((\Acal^*_n)_{\Qbar}, \IC(V)) = 0\) unless the central character of \(V\)
is a square.
Indeed this follows from Remark \ref{rema:Hecke_central_char} applied to the
central element \(z = -1\).
Thus up to twisting \(V\) by a power of the similitude character (see Remark
\ref{rema:twisting_sim}) we may and will assume that \(V\) is a representation of
\(\PGSpbf_{2n}\), i.e.\ that the highest weight of \(V\) is \((k_1, \dots, k_n,
(\sum_i k_i)/2)\).
Thanks to this assumption we will be able to use trace formulas for semisimple
groups as in the first section (see Section \ref{sec:non_ss_gps}), simplifying
our notation.
Let \(\tau\) be the projection of the dual of the infinitesimal character of \(V\)
to \(\widehat{\mathfrak{g}}_{\mathrm{der}} = \mathfrak{sp}_{2n}\).

To compute the right-hand side of \eqref{eq:ICMorel}, first recall that for
\(\efrak = (\Hbf, \Hcal, s, \xi)\) an elliptic endoscopic datum for \(\Gbf =
\GSpbf_{2n}\) as in Theorem \ref{t:ICMorel}, \(\iota(\efrak)\) is defined as
\(\tau(\Gbf) \tau(\Hbf)^{-1} |\Out(\efrak)|^{-1}\) where \(\tau\) is the Tamagawa
number and \(\Out(\efrak)\) is the outer automorphism group of \(\efrak = (\Hbf,
\Hcal, s, \xi)\).
These constants are easily computed: if \(\Hbf = \Gbf\), \(\iota(\efrak) = 1\), and
otherwise \(\iota(\efrak) = 1/4\) (in the latter case \(\Out(\efrak) = \Z/2\Z\)).
Also recall from Lemma \ref{lem:unr_endo_gps_isog} and the discussion following
Theorem \ref{t:ICMorel} that the natural map \(\efrak \mapsto \ol{\efrak} =
(\ol{\Hbf}, \ol{\Hcal}, \ol{s}, \ol{\xi})\) induces a bijection between sets of
equivalence classes of everywhere unramified endoscopic data for \(\Gbf\) (or
\(\PGSpbf_{2n}\)) and for \(\Spbf_{2n}\).

For \(\tau' \in \ICcal(\ol{\Hbf})\) satisfying \(\ol{\xi}(\tau') = \tau\) and \(\psi'
\in \Psit_{\disc}^{\unr, \tilde{\tau}}(\ol{\Hbf})\) we can consider \(\xi \circ
\dot{\psi}'_{\tau', \mathrm{sc}} : \Mcal_{\psi'} \rightarrow
(\widehat{\Spbf_{2n}})_{\mathrm{sc}}\), which is identified to \(\dpsitausc\) for a
uniquely determined \(\psi \in \Psit_{\disc}^{\unr, \tau}(\Spbf_{2n})\).
Writing \(\psi' = (\psi'_1, \psi'_2)\) for a decomposition \(\ol{\Hbf} \simeq
\SObf_{4a} \times \Spbf_{2b}\), we have \(\psi = \psi'_1 \oplus \psi'_2\) (if
\(\efrak\) is trivial then \(a=0\) and \(\psi'_1\) is understood to be an empty formal
sum).
In other words we have a natural map \((\efrak, \tau', \psi') \mapsto (\psi,
\ol{s})\), which is a bijection between the set of equivalence classes of triples
\((\efrak, \tau', \psi')\) where \(\efrak\) is an elliptic endoscopic datum for
\(\Gbf\), \(\tau' \in \ICcal(\ol{\Hbf})\) maps to \(\tau\) and \(\psi' \in
\Psit_\disc^{\unr, \tilde{\tau'}}(\ol{\Hbf})\) (the equivalence being induced by
the usual notion of equivalence of endoscopic data, which acts on the second
factor \(\tau'\) and acts trivially on the last factor \(\psi'\)) and the set of
pairs \((\psi, s)\) where \(\psi \in \Psit_\disc^{\unr, \tau}(\Spbf_{2n})\) and
\(\ol{s} \in \Scal_{\psi}\).
The fact that this is a bijection is a special case of an observation of Arthur
\cite[Proposition 2.4.1]{Taibi_mult}, which goes back to a general argument due
to Kottwitz \cite[\S 11]{Kottwitz_STFcusptemp}.
If \(\efrak\) is non-trivial then \(|\Out(\efrak)| = 2\) and the non-trivial outer
automorphism does not fix \(\tau'\), as we have \(\tau' = (\tau'_1, \tau'_2)\) where
\(\tau'_1\) belongs to a \(\thetahat\)-orbit having two elements.

Thus we can begin rewriting the right-hand side of \eqref{eq:ICMorel} purely in
terms of parameters for \(\Gbf\) (or rather \(\Spbf_{2n}\)), using the reduction
from \(\Gbf\) to \(\PGSpbf_{2n}\) explained in Section \ref{sec:non_ss_gps} and the
spectral expansion \eqref{eq:Sdisclifted} in Proposition
\ref{pro:lift_stab_mult}, along with the spectral characterizations
\eqref{eq:spectral_char_real} (at the real place),
\eqref{eq:spectral_char_twisted} (at the \(p\)-adic place) and the fundamental
lemma (at all other places): \eqref{eq:ICMorel} is equal to
\begin{multline*}
  \sum_{\psi \in \Psit_{\disc}^{\unr, \tau}(\Spbf_{2n})} \sum_{\ol{s} \in
  \Scal_{\psi}} \iota(\efrak) \sum_{\tau' \mapsto \tau}
  \frac{\epsilon_{\psi'}(s_{\psi'}) (-1)^{q(\Lbf^*_{\psi',\tau'})}}{|\Scal_{\psi'}|}
  (-1)^{q(\Gbf_{\R})} \langle \mu_{\pi_{\infty}^{\gen}}, s \rangle \times \\
  p^{jn(n+1)/4} \Tr \left( r_{-\mu}(s \dpsitausc(\cpsc(\psi))^j) \right)
  \prod_{q \neq p} \Sat_{\Gbf_{\Zq}}(f_q)(\dpsitausc(c_{q, \mathrm{sc}}(\psi)))
\end{multline*}
where \(\efrak = (\Hbf, \Hcal, s, \xi)\) is the elliptic endoscopic datum for
\(\Gbf\) obtained from a lift \(s \in \widehat{\Gbf}(\C)\) of \(\ol{s}\) and \(\psi'
\in \Psi_{\disc}^{\unr, \tilde{\tau'}}(\ol{\Hbf})\) correspond to \((\psi,
\ol{s})\).
Note that the product \(\langle \mu_{\pi_{\infty}^{\gen}}, s \rangle r_{-\mu}(s)\)
does not depend on the choice of \(s\).
In the above expression we also implicitly see \(\dpsitausc(\cpsc(\psi)) \in
\widehat{\PGSpbf_{2n}}(\C)\) as an element of \(\widehat{\Gbf}(\C)\), and
\(\pi_\infty^\mathrm{gen}\) is the generic element of the L-packet as discussed
after \eqref{eq:spectral_char_real}.

We proceed to remove any reference to endoscopic objects in this expression.
First we have \(\epsilon_{\psi'}(s_{\psi'}) = \epsilon_{\psi}(s s_{\psi})\) by \cite[Lemma 4.4.1]{Arthur}.
The quasi-split real connected reductive group \(\Lbf^*_{\psi',\tau'}\) associated to the Adams-Johnson parameter \(\dpsitau' \circ \psi_{\infty}'\) is isomorphic to \(\Lbf^*_{\psi,\tau}\), associated to the parameter \(\dpsitau \circ \psi_{\infty}\).
Finally if \((\ol{\efrak}, \tau', \psi')\) corresponds to \((\psi, \ol{s})\) then the inclusion \(C_{\psi'} \subset C_{\psi}\) induced by \(\xi\) turns out to be an equality because \(\dpsitau\) induces an isomorphism \(Z(\Mpsi) \to C_{\dpsitau}\) and similarly for \(\psi'\) (see third part of Definition \ref{def:Mpsi}).
Thus we have \(|\Scal_\psi|/|\Scal_{\psi'}| = |Z(\widehat{\ol{\Hbf}})|\) which is equal to \(2\) if \(\efrak\) is non-trivial and \(1\) otherwise.
Thus we have
\[ \iota(\efrak) \times \frac{|\{ \tau' | \tau' \mapsto \tau \}|}{
|\Scal_{\psi'}|} = |\Scal_\psi|^{-1}, \]
and the right-hand side of \eqref{eq:ICMorel} is also equal to

\begin{multline*}
  \sum_{\psi \in \tilde{\Psi}_{\disc}^{\unr, \tau}(\Spbf_{2n})} \frac{
  (-1)^{q(\Gbf_{\R}) + q(\Lbf^*_{\psi,\tau})} }{ |\Scal_{\psi}| } \sum_{\ol{s}
  \in \Scal_{\psi}} \epsilon_{\psi}(\ol{s} s_{\psi}) \langle
  \mu_{\pi_{\infty}^{\mathrm{gen}}}, s \rangle \\
  p^{jn(n+1)/4} \Tr \left( r_{-\mu}(s \dpsitausc(\cpsc(\psi))^j ) \right)
  \prod_{q \neq p} \Sat_{\Spbf_{2n,\Z_q}}(f_q)(c_{q, \mathrm{sc}}(\psi)).
\end{multline*}

Recall that in Arthur's formalism \(s_{\psi} \in \Scal_\psi\) is defined in a
global manner, as the image of \(-1 \in \SL_2\) by a morphism \(\Lcal_\psi \times
\SL_2 \rightarrow \widehat{\Spbf_{2n}}\) (denoted \(\tilde{\psi}_{\Spbf_{2n}}\) in
\cite[\S 1.4.4]{Arthur}) corresponding to the formal Arthur-Langlands parameter
\(\psi\).
We shall need a canonical preimage of \(s_{\psi}\) in \(\widehat{\PGSpbf_{2n}}\),
which we cannot define in exactly the same way since we are using a slightly
weaker formalism where \(\Mpsi\) replaces \(\Lcal_\psi \times \SL_2\).
Recall from Definition \ref{def:Mpsisc} that there exists \(\psi_{\infty,
\mathrm{sc}} : W_{\R} \times \SL_2 \rightarrow \Mpsisc(\C)\) lifting
\(\psi_{\infty}\) and that this lift is unique up to \(Z^1(W_\R, \ker( \Mpsisc
\rightarrow \Mpsi ))\).
Therefore \(\tilde{s}_{\psi} := \psi_{\infty, \mathrm{sc}}(1,-1) \in \Mpsisc\) satisfies \(\tilde{s}_{\psi}^2 = 1\), lifts Arthur's \(s_\psi \in C_{\psi}\) (see Lemma \ref{lem:epsilon_psi_s_psi}) and does not depend on the choice of the lift \(\psi_{\infty, \mathrm{sc}}\).
Moreover \(\tilde{s}_{\psi}\) belongs to \(Z(\Mpsisc)\) because \(s_\psi\) belongs to \(Z(\Mpsi)\).
For any \(\psi_{\infty, \mathrm{sc}}\) as above \(\dpsitausc \circ \psi_{\infty,
\mathrm{sc}}\) is an Adams-Johnson parameter,
so we can apply \cite[Lemma 9.1]{Kottwitz_AA} and conclude that
\((-1)^{q(\Lbf^*_{\psi,\tau})} = \langle \mu_{\pi_{\infty}^{\mathrm{gen}}},
\dpsitausc(\tilde{s}_{\psi}) \rangle\).
Using this identity and the change of variable \(\ol{s} \mapsto \ol{s} s_{\psi}\)
in the sum, our expression for the right-hand side of \eqref{eq:ICMorel} becomes

\begin{multline} \label{eq:penult_expr_before_theo_crude}
  (-1)^{q(\Gbf_{\R})} \sum_{\psi \in \tilde{\Psi}_{\disc}^{\unr,
  \tau}(\Spbf_{2n})} |\Scal_{\psi}|^{-1} \sum_{\ol{s} \in \Scal_{\psi}}
  \epsilon_{\psi}(\ol{s}) \langle \mu_{\pi_{\infty}^{\mathrm{gen}}}, s \rangle
  \\
  p^{jn(n+1)/4} \Tr \left( r_{-\mu}(s \dpsitausc(\tilde{s}_{\psi} \cpsc(\psi))^j
  ) \right) \prod_{q \neq p} \Sat_{\Gbf_{\Z_q}}(f_q)(c_{q, \mathrm{sc}}(\psi)).
\end{multline}

We now use the notations introduced in Section \ref{sec:not_red_gps}, in
particular the identification \(\SO_{2n+1} \simeq \widehat{\Spbf_{2n}}\) and the
parametrization \(\Tcal_{\SO_{2n+1}} \simeq \GL_1^n\), and similarly for \(\Gbf =
\GSpbf_{2n}\) and \(\PGSpbf_{2n}\).
Fix \(\psi_{\infty} : W_{\R} \times \SL_2(\C) \to \Mpsi(\C)\) in the conjugacy
class introduced in Definition \ref{def:Mpsi}.
Up to conjugacy \footnote{of \(\psi_{\infty}\) by \(\Mpsi(\C)\) and of \(\dpsitau\) by \(\Ghat(\Qbar)\)} we may assume that \(\dpsitau \circ \varphi_{\psi_\infty} |_{\C^\times}\) takes values in \(\Tcal_{\SO_{2n+1}}(\C)\) and is dominant for \(\Bcal_{\SO_{2n+1}}\), i.e.\ the holomorphic part of \(\dpsitau \circ \varphi_{\psi_\infty} |_{\C^\times}: \C^\times \to \Tcal_{\SO_{2n+1}}(\C)\) is \(z \mapsto (z^{k_1+n}, \dots, z^{k_n+1})\).
Write \(\psi = \oplus_{i=0}^r \psi_i\) as in Definition \ref{def:Mpsi}.
As in Lemma \ref{lem:epsilon_psi_s_psi} for \(1 \leq i \leq r\) let \(s_i \in Z(\Mpsi)\) be the element such that for \(1 \leq i' \leq r\) the \(i'\)-th projection to \(Z(\Mcal_{\psi_i}) = \{ \pm 1 \}\) is non-trivial if and only if \(i'=i\).
Then \((s_1, \dots, s_r)\) is a basis of \(C_{\dpsitau}\) (seen as a vector space over \(\F_2\)), and it determines a partition \(\{1, \dots, n \} = J_0 \sqcup \dots \sqcup J_r\) as follows: for \(1 \leq i \leq r\) and \(1 \leq j \leq n\) we have \(j \in J_i\) if and only if the \(j\)-th component of \(s_i\) (seen as an element of \(\Tcal_{\SO_{2n+1}}(\Qbar) \simeq (\Qbar^\times)^n\)) is \(-1\).
For \(1 \leq i \leq r\) let \(\tilde{s}_i = (x_1, \dots, x_n, 1, 1) \in
\Tcal_{\GSpin_{2n+1}}(\Qbar)\) where
\[ x_j = \begin{cases}
    -1 & \text{ if } j \in J_i \\
    1 & \text{ otherwise.}
\end{cases} \]
In other words \(\tilde{s}_i\) is the image by \(\dpsitausc\) of the element of \(Z(\Mcal_{\psi_i, \mathrm{sc}})\) mapping to the non-trivial element of \(Z(\Mcal_{\psi_i})\) and acting by \(+1\) in \(\spin^+_{\psi_i}\) and by \(-1\) in \(\spin^-_{\psi_i}\) (see Definition \ref{def:spin_pm_psi}).
Clearly \(\tilde{s}_i\) lifts \(s_i\) and we have \(\tilde{s}_i \in C_{\dpsitausc}\).
This gives us a parametrization \((\Z/2\Z)^r \times \GL_1 \simeq
\Cent(\dpsitausc, \widehat{\Gbf})\) mapping \(((\epsilon_i)_{1 \leq i \leq r},
\lambda)\) to \(\nuhat(\lambda) \prod_i \tilde{s}_i^{\epsilon_i}\).

A simple computation of weights shows that we have a decomposition into
irreducible constituents
\begin{equation} \label{eq:decomp_rmu_hspin}
  r_{-\mu} \circ \dpsitausc \simeq \bigoplus_{(u_1, \dots, u_r) \in \{ \pm 1 \}^r} \spin_{\psi_0} \otimes \spin_{\psi_1}^{u_1} \otimes \dots \otimes \spin_{\psi_r}^{u_r}
\end{equation}
and on each factor the group \(\Cent(\dpsitausc, \widehat{\Gbf})\) acts by a
character that we denote \(\alpha_{\psi, (u_1, \dots, u_r)}\).
With the above parametrization of \(\Cent(\dpsitausc, \widehat{\Gbf})\), this
character maps \(((\epsilon_i)_{1 \leq i \leq r}, \lambda)\) to \(\lambda^{-1}
\prod_i u_i^{\epsilon_i}\).
For \(s\) as in \eqref{eq:penult_expr_before_theo_crude} the automorphism \(\langle
\mu_{\pi_{\infty}^{\gen}}, s \rangle r_{-\mu}(s)\) acts on the factor of
\eqref{eq:decomp_rmu_hspin} corresponding to \((u_1, \dots, u_r)\) by the scalar
\(\langle \mu_{\pi_{\infty}^{\gen}}, s \rangle \alpha_{\psi, (u_1, \dots,
u_r)}(s)\).
This scalar clearly does not depend on the choice of \(s\) lifting \(\ol{s}\) and is
the evaluation at \(\ol{s}\) of a character of \(\Scal_\psi\) that we denote
\(\beta_{\psi,(u_1, \dots, u_r)}\).
This last character is easily computed: by \eqref{eq:comp_muh_gen} we have \(\langle \mu_{\pi_{\infty}^{\mathrm{gen}}}, \tilde{s}_i \rangle = (-1)^{N_i}\) where \(N_i\) is the cardinality of \(\left\{ j \in J_i \,\middle|\, j \text{ even} \right\}\).
For \(1 \leq i \leq r\) define
\begin{equation} \label{eq:def_ui_psi}
  u_i(\psi) = \epsilon_{\psi}(s_i) \langle \mu_{\pi_{\infty}^{\gen}},
  \tilde{s}_i \rangle \in \{ \pm 1 \}
\end{equation}
so that \((u_1(\psi), \dots, u_r(\psi))\) is the unique \((u_1, \dots, u_r)\) in \(\{
\pm 1 \}^r\) such that the character \(\epsilon_\psi \beta_{\psi, (u_1, \dots,
u_r)}\) of \(\Scal_\psi\) is trivial.
Using orthogonality relations for characters of \(\Scal_\psi\) we obtain
that \eqref{eq:penult_expr_before_theo_crude} is equal to
\begin{multline}
  (-1)^{q(\Gbf_{\R})} \sum_{\psi \in \tilde{\Psi}_{\disc}^{\unr,
  \tau}(\Spbf_{2n})} p^{jn(n+1)/4} \alpha_{\psi, (u_1(\psi), \dots, u_r(\psi))}
  (\dpsitausc(\tilde{s}_\psi)) \\
  \Tr \left( \spin_{\psi_0} \cpsc(\psi_0)^j \right) \prod_{i=1}^r \Tr \left(
  \spin_{\psi_i}^{u_i(\psi)} \cpsc(\psi_i)^j \right) \prod_{q \neq p}
  \Sat_{\Gbf_{\Z_q}}(f_q)(c_{q, \mathrm{sc}}(\psi)).
\end{multline}

It will also be useful to compute \(\tilde{s}_{\psi} \in Z(\Mpsisc) =
\prod_{i=0}^r Z(\Mcal_{\psi_i, \mathrm{sc}})\) explicitly.
A special case will be used in Proposition \ref{pro:no_amb_symplectic} below.
\begin{enumerate}
  \item For \(\psi_0 = \pi_0[d_0]\) with \(\pi_0 \in O_o(w_1^{(0)}, \dots,
    w_{(n_0-1)/2}^{(0)})\), we need to compute the image of \(-1\) in the lift of
    \(\nu_{d_0}^{\oplus n_0} : \SL_2 \rightarrow \SO_{n_0 d_0}\) to \(\Spin_{n_0
    d_0}\).
    Computing with weights we find that this element of \(Z(\Spin_{n_0 d_0}) = \{
    1, \tilde{s}_0 \}\) is \(\tilde{s}_0^{(d_0^2-1)/8}\).
  \item For \(\psi_i = \pi_i[d_i]\) with \(d_i\) even and \(\pi_i \in S(w_1^{(i)},
    \dots, w_{k_i}^{(i)})\), computing in the same way we find that the image of
    \(-1 \in \SL_2\) is \((-1, \dots, -1, 1) = \tilde{s}_i \in Z(\Spin_{n_i
    d_i}) \subset \Tcal_{\Spin_{n_i d_i}}\).
  \item For \(\psi_i = \pi_i[d_i]\) with \(d_i\) odd and \(\pi_i \in O_e(w_1^{(i)},
    \dots, w_{2 k_i}^{(i)})\), we find that the image of \(-1 \in \SL_2\) in
    \(Z(\Spin_{n_i d_i}) \simeq \{ \pm 1 \}^2\) is trivial.
\end{enumerate}
Thus we have \(\tilde{s}_{\psi} = \tilde{s}_0^{(d_0^2-1)/8} \prod_{i=1}^r
\tilde{s}_i^{d_i-1}\).

We finally obtain that the right-hand side of \eqref{eq:ICMorel} equals
\begin{multline} \label{eq:last_expr_before_theo_crude}
  \sum_{\psi \in \tilde{\Psi}_{\disc}^{\unr, \tau}(\Spbf_{2n})} (-1)^{n(n+1)/2 +
  (d_0^2-1)/8} \prod_{i=1}^r u_i(\psi)^{d_i-1} \\
  \times p^{jn(n+1)/4} \Tr \left( \spin_{\psi_0}(\cpsc(\psi_0))^j) \right)
  \prod_{i=1}^r \Tr \left( \spin_{\psi_i}^{u_i(\psi)}(\cpsc(\psi_i))^j) \right)
  \\
  \times \prod_{q \neq p} \Sat_{\Gbf_{\Z_q}}(f_q)(c_{q, \mathrm{sc}}(\psi)).
\end{multline}

Theorem \ref{theo:IH_explicit_crude} reformulates this slightly
more conceptually.

\begin{theo} \label{theo:IH_explicit_crude}
  Let \(n \geq 1\), \(\ell\) a prime and \(\iota\) be an isomorphism between the algebraic closures of \(\Q\) in \(\C\) and \(\Qellbar\).
  Let \(V\) be an irreducible algebraic representation of \(\PGSpbf_{2n}\), and let \(\tau\) be the infinitesimal character of its dual representation (understood as a representation of \(\PGSpbf_{2n}(\R)\)), a semisimple conjugacy class in \(\mathfrak{so}_{2n+1}\).
  In the Grothendieck group \(K_0 \left( \Rep_{\Qellbar}(\GalQ \times \Hcal_f^{\unr}(\Gbf)_{\Qellbar}) \right)\) we have
  \begin{equation} \label{eq:explicitIH}
    [\Qellbar \otimes_{\Qell} H^\bullet((\Acal_n^*)_{\Qbar}, \IC(V))] = \sum_{\psi \in \tilde{\Psi}^{\unr,\tau}_{\disc}(\Spbf_{2n})} (-1)^{n(n+1)/2 + (d_0^2-1)/8} \prod_{i=1}^r u_i(\psi)^{d_i-1} \sigma_{\psi, \iota}^{\IH} \otimes \iota(\chi_{f, \psi})
  \end{equation}
  where in the sum,
  \begin{itemize}
    \item as above \(\psi = \psi_0 \oplus \dots \oplus \psi_r\), \(\psi_i = \pi_i[d_i]\), and \(u_i(\psi) = \epsilon_{\psi}(s_i) \langle \mu_{\pi_{\infty}^{\gen}}, \tilde{s}_i \rangle \in \{ \pm 1 \}\),
    \item \(\sigma_{\psi, \iota}^{\IH}\) is a continuous semisimple representation \(\GalQ \rightarrow \GL_{2^{n-r}}(\Qellbar)\) characterized by the properties that it is unramified away from \(\ell\) and that for any \(p \neq \ell\), \(\sigma_{\psi, \iota}^{\IH}(\Frob_p)\) is conjugated to
      \begin{equation} \label{eq:explicitIH_cc_Frob}
        \iota \left( p^{n(n+1)/4} \spin_{\psi_0}(\cpsc(\psi_0)) \otimes \spin_{\psi_1}^{u_1(\psi)}(\cpsc(\psi_1)) \otimes \dots \otimes \spin_{\psi_r}^{u_r(\psi)}(\cpsc(\psi_r)) \right).
      \end{equation}
    \item \(\chi_{f, \psi}\) is the character \(\Hcal^{\unr}_f(\Gbf) \rightarrow \Qbar\) determined by \((\dpsitausc(\cpsc(\psi)))_{p \in \Pcal}\) (see Proposition \ref{pro:rat_Sat}), and \(\iota(\chi_{f,\chi})\) abusively denotes the \(\Qellbar\)-linear extension to \(\Hcal^{\unr}_f(\Gbf)_{\Qellbar}\) of its composition with \(\iota\).
  \end{itemize}
\end{theo}
\begin{proof}
  Recall that the equality between \eqref{eq:ICMorel} and
  \eqref{eq:last_expr_before_theo_crude} holds true for \(f^\infty = \prod_q f_q
  \in \Hcal^{\unr}_f(\Gbf)\) with \(f_p = 1\) and \(j\) sufficiently large.
  Recall from Remark \ref{rem:rationality} that the left-hand side, which a
  priori belongs to \(\Qell\), is in fact rational.
  The right-hand side a priori belongs to \(\C\) by definition of endoscopic
  transfer, but thanks to Proposition \ref{pro:rat_Sat} we see that each term
  belongs to \(\Qbar\).
  It follows that for any \(f^\infty = \prod_q f_q \in
  \Hcal^{\unr}_f(\Gbf)_{\Qellbar}\) with \(f_p = 1\) we have, for \(j\) sufficiently
  large,
  \begin{align} \label{eq:IH_explicit_crude}
    & \mathrm{Tr} \left( \Frob_p^j f^{\infty} \, \middle| \, \Qellbar
    \otimes_{\Qell} H^\bullet((\Acal^*_n)_{\Qbar}, \IC(V)) \right) \\
    = & \sum_{\psi \in \tilde{\Psi}_{\disc}^{\unr, \tau}(\Spbf_{2n})}
    (-1)^{n(n+1)/2 + (d_0^2-1)/8} \prod_{i=1}^r u_i(\psi)^{d_i-1} \times
    \iota(\chi_{f, \psi})(f^\infty) \nonumber \\
    & \ \ \ \times \iota \left( p^{jn(n+1)/4} \Tr \left(
    \spin_{\psi_0}(\cpsc(\psi_0))^j) \right) \prod_{i=1}^r \Tr \left(
    \spin_{\psi_i}^{u_i(\psi)}(\cpsc(\psi_i))^j) \right) \right) \nonumber
  \end{align}
  where the expression between the outer parentheses on the last line belongs to
  a finite extension of \(\Q\) in \(\C\) which only depends on \(\psi\).
  By a standard argument using invertibility of a Vandermonde determinant,
  \eqref{eq:IH_explicit_crude} holds true for any \(j \in \Z\).
  
  In particular for \(j=0\) we obtain
  \begin{multline*} 
    \mathrm{Tr} \left( f^{\infty} \, \middle| \, H^\bullet((\Acal^*_n)_{\Qbar},
    \IC(V)) \right) \\
    = \sum_{\psi \in \tilde{\Psi}_{\disc}^{\unr, \tau}(\Spbf_{2n})}
    (-1)^{n(n+1)/2 + (d_0^2-1)/8} 2^{n-r} \prod_{i=1}^r u_i(\psi)^{d_i-1}
    \times \iota(\chi_{f, \psi})(f^\infty).
  \end{multline*} 
  Observe that \(p\) does not occur in this formula except in the assumption that
  we have \(f_p = 1\), which for a given \(f^\infty\) is satisfied for almost all
  \(p\).
  So the formula is satisfied for any \(f^\infty \in \Hcal_f^{\unr}(\Gbf)\),
  determining \([\Qellbar \otimes_{\Qell} H^\bullet(\Acal_n^*, \IC(V))]\) in \(K_0
  \left( \Rep(\Hcal_f^{\unr}(\Gbf)_{\Qellbar}) \right)\).
  By Remark \ref{rem:mult_one} for any \(\psi \in
  \tilde{\Psi}_{\disc}^{\unr,\tau}(\Spbf_{2n})\) and any finite set \(S\) of prime
  numbers there exists \(f^\infty \in \bigotimes'_{q \not\in S}
  \Hcal^{\unr}(\Gbf_{\Z_q})_{\Qellbar}\) such that for any \(\psi' \in
  \tilde{\Psi}_{\disc}^{\unr,\tau}(\Spbf_{2n})\) we have
  \[ \iota(\chi_{f, \psi'})(f^\infty) = \begin{cases}
	1 & \text{ if } \psi' = \psi \\
	0 & \text{ otherwise.}
  \end{cases} \]
  In particular using just \(S = \emptyset\) we obtain the equality in \(K_0 \left(
  \Rep_{\Qellbar}(\GalQ \times \Hcal_f^{\unr}(\Gbf)_{\Qellbar}) \right)\)
  \begin{align*}
    & \ [\Qellbar \otimes_{\Qell} H^\bullet(\Acal_n^*, \IC(V))] \\
    = & \ \sum_{\psi \in \tilde{\Psi}_{\disc}^{\unr, \tau}(\Spbf_{2n})}
    (-1)^{n(n+1)/2 + (d_0^2-1)/8} \prod_{i=1}^r u_i(\psi)^{d_i-1} \times
    \sigma_{\psi,\iota}^{\IH} \otimes \iota(\chi_{f, \psi})
  \end{align*}
  where \(\sigma_{\psi,\iota}^{\IH} \in K_0 \left( \Rep_{\Qellbar}(\Gal_\Q)
  \right)\) is uniquely determined and has virtual dimension \(2^{n-r}\).
  Taking \(S = \{p\}\) where \(p \neq \ell\) and using \eqref{eq:IH_explicit_crude}
  yields
  \[ \Tr \left( \Frob_p^j \middle| \sigma_{\psi,\iota}^{\IH} \right) = \iota
  \left( p^{jn(n+1)/4} \Tr \left( \spin_{\psi_0}(\cpsc(\psi_0)^j) \right)
  \prod_{i=1}^r \Tr \left( \spin_{\psi_i}^{u_i(\psi)}(\cpsc(\psi_i)^j) \right)
  \right) \]
  for all \(j \in \Z\).

  To conclude we must show that \(\sigma_{\psi,\iota}^{\IH}\) is a genuine
  representation.
  This follows from purity (\cite[Proposition 6.2.6]{Deligne_Weil2}) which
  implies that an element of \(\Qellbar^\times\) can occur as an eigenvalue of
  \(\Frob_p\) acting on \(H^k((\Acal_n^*)_{\Fpbar}, \IC(V))\) in at most one degree
  \(k\), and invertibility of a Vandermonde determinant.
\end{proof}

Observing the characterization of \(\sigma_{\psi, \iota}^{\IH}\) (and the general
conjecture \cite{Kottwitz_AA}) we expect that:
\begin{enumerate}
  \item the representation \(\sigma_{\psi, \iota}^{\IH}\) factors as a tensor
    product of Galois representations as suggested by
    \eqref{eq:explicitIH_cc_Frob}, and
  \item each factor is obtained from a Galois representation taking values in a
    \(\mathrm{GSpin}\) group by composing with a spin or half-spin representation.
\end{enumerate}

As explained in the introduction (\S \ref{sec:intro_IH_GSpin_Gal_rep}) we will show the first point in full generality, and the second point in almost all cases.

\begin{coro} \label{coro:seed_ex_sigma}
  Fix a prime number \(\ell\) and an isomorphism \(\iota\) between \(\Qbar\) and the
  algebraic closure of \(\Q\) in \(\Qellbar\).
  \begin{enumerate}
    \item
      Let \(n \geq 1\), \(\tau \in \ICcal(\Spbf_{2n})\) and \(\psi \in \Psi^{\unr,
      \tau}_{\disc, \nonendo}(\Spbf_{2n})\).
      Explicitly, we have \(\psi = \pi[2d+1]\) for some \(\pi \in O_o(w_1, \dots,
      w_k)\) and \(d \geq 0\).
      There exists a continuous semisimple representation \(\sigma_{\psi,
      \iota}^{\spin} : \GalQ \rightarrow \GL_{2^n}(\Qellbar)\) unramified away
      from \(\ell\) such that for any \(p \neq \ell\), \(\sigma_{\psi,
      \iota}^{\spin}(\Frob_p)^{\sesi}\) is conjugated to
      \[ \iota \left( p^{n(n+1)/4} \spin_{\psi} (\cpsc(\psi)) \right). \]
    \item
      Let \(n \geq 1\), \(\tilde{\tau} \in \tilde{\ICcal}(\SObf_{4n})\), \(\psi \in
      \tilde{\Psi}^{\unr, \tilde{\tau}}_{\disc, \nonendo}(\SObf_{4n})\).
      Recall that either \(\psi = \pi[2d+1]\) for some \(\pi \in O_e(w_1, \dots,
      w_{2k})\) and \(d \geq 0\) or \(\psi = \pi[2d]\) for some \(\pi \in S(w_1,
      \dots, w_k)\) and \(d \geq 1\).
      Let \(\epsilon = (-1)^n \epsilon(\frac{1}{2}, \pi) = (-1)^n
      \epsilon(\frac{1}{2}, \pi_{\infty})\).
      In the first case we have \(\epsilon = (-1)^n\), in the second case it
      depends on the infinitesimal character of \(\pi_{\infty}\): we have
      \(\epsilon = (-1)^{n + k/2 + \sum_i w_i}\).
      Then there is a continuous semisimple representation \(\sigma_{\psi,
      \iota}^{\spin, \epsilon} : \GalQ \rightarrow \GL_{2^{2n-1}}(\Qellbar)\)
      which is unramified away from \(\ell\) and such that for any \(p \neq \ell\),
      \(\sigma_{\psi, \iota}^{\spin, \epsilon}(\Frob_p)^{\sesi}\) is conjugated to
      \[ \iota \left( p^{n/2} \spin_{\psi}^{\epsilon} (\cpsc(\psi)) \right). \]
  \end{enumerate}
\end{coro}
\begin{proof}
  Let \(\sigma_{\psi, \iota}^{\spin} = \sigma_{\psi, \iota}^{\IH}\) in the first case and \(\sigma_{\psi, \iota}^{\spin, \epsilon} = \sigma_{1 \oplus \psi, \iota}^{\IH}(n^2)\) (Tate twist) in the second case (note that \(n^2 = (2n(2n+1)/2 - n/2)/2\)).
\end{proof}

\begin{defi} \label{def:GMpsisc}
  Let \(n \geq 1\), \(\taut \in \ICcalt(\SObf_{4n})\) and \(\psi \in \Psit_{\disc, \nonendo}^{\unr, \taut}(\SObf_{4n})\).
  Let \(\GMpsisc\) be quotient of \(\GL_1 \times \Mpsisc\) by the diagonally embedded subgroup \(\mu_2\) (on the second factor, the kernel of \(\Mpsisc \to \Mpsi\)).
  Recall from Definition \ref{def:spin_pm_psi} the two representations \(\spin_\psi^\epsilon\) of \(\Mpsisc\), for \(\epsilon \in \{+,-\}\).
  They can be extended uniquely to \(\GMpsisc\) by letting \(z \in \GL_1\) act by \(z \, \id\).
  We simply denote these representations by \(\spin_\psi^\epsilon\).
  Denote by \(\nu\) the character \(\GMpsisc \to \GL_1\) induced by \(\GL_1 \times \Mpsisc \to \GL_1\), \((z,g) \mapsto z^2\).
\end{defi}

\begin{prodef} \label{def:alpha_m_d}
  Let \(n,m,d \geq 1\) be integers satisfying \(n=md\).
  The representation \(\Std_{\Sp_{2m}} \otimes \Sym^{2d-1} \Std_{\SL_2}\) of \(\Sp_{2m} \times \SL_2\) is irreducible and self-dual of orthogonal type.
  The set of morphisms \(\alpha: \Sp_{2m} \times \SL_2 \to \SO_{4n}\) satisfying \(\Std_{\SO_{4n}} \circ \alpha \simeq \Std_{\Sp_{2m}} \otimes \Sym^{2d-1} \Std_{\SL_2}\) consists of two \(\SO_{4n}(\Qbar)\)-conjugacy classes.
  Exactly one of these two conjugacy classes has a representative \(\alpha_{m,d}\) restricting to (using parametrizations from Section \ref{sec:not_red_gps})
  \begin{align}
    \Tcal_{\Sp_{2m}} \times \Tcal_{\SL_2} & \longrightarrow \Tcal_{\SO_{4n}} \nonumber \\
    ((x_1, \dots, x_m), t) & \longmapsto (x_1 t^{2d-1}, x_1 t^{2d-3}, \dots, x_1 t^{1-2d}, \dots, x_m t^{2d-1}, \dots, x_m t^{1-2d}).  \label{eq:alpha_m_d_on_T}
  \end{align}
  It admits a (unique) lift \(\tilde{\alpha}_{m,d}: \Sp_{2m} \times \SL_2 \to \Spin_{4n}\).
  We have an identification
  \begin{align*}
    (\Sp_{2m} \times \SL_2 \times \GL_1) / \mu_2 & \longrightarrow G(\Sp_{2m} \times \SL_2) := \{ (g_1,g_2) \in \GSp_{2m} \times \GL_2 \,|\, \nu(g_1) = \det g_2 \} \\
    (h_1, h_2, \lambda) & \longmapsto (\lambda h_1, \lambda h_2)
  \end{align*}
  where \(\mu_2\) is diagonally embedded on the left, and the morphism
  \begin{align}
    \Sp_{2m} \times \SL_2 \times \GL_1 & \longrightarrow \GSpin_{4n} \label{eq:Sp2m_SL2_GL1_to_GSpin4n} \\
    (h_1, h_2, \lambda) & \longmapsto \lambda^n \tilde{\alpha}_{m,d}(h_1, h_2) \nonumber
  \end{align}
  induces a morphism \(G(\Sp_{2m} \times \SL_2) \to \GSpin_{4n}\), that we abusively still denote by \(\tilde{\alpha}_{m,d}\).

  Consider \(\taut \in \ICcalt(\SObf_{4n})\) and \(\psi \in \Psit_{\disc}^{\unr, \taut}(\SObf_{4n})\) of the form \(\pi[2d]\) (as in the second case of the second point of Corollary \ref{coro:seed_ex_sigma}).
  The Langlands parameter \(W_\R \rightarrow \GL_{2m}(\C)\) of \(\pi_{\infty}\) is symplectic, and so it factors through \(\varphi_{\infty} : W_\R \rightarrow \Sp_{2m}(\C)\) which is well-defined up to conjugation.
  There exists \(\alpha_\psi: \Sp_{2m} \times \SL_2 \rightarrow \Mpsi\) satisfying \(\Std_{\Mpsi} \circ \alpha_\psi \simeq \Std_{\Sp_{2m}} \otimes \Sym^{2d-1} \Std_{\SL_2}\) and such that \(\psi_\infty\) is conjugated under \(\Mpsi(\C)\) to \(\alpha_\psi \circ (\varphi_\infty, \id_{\SL_2})\), and for these two properties \(\alpha_\psi\) is unique up to conjugation by \(\Mpsi(\Qbar)\).
  Taking \(\Mpsi = \SO_{4n}\) and \(\tau_\psi = (w_1,\dots,w_{2n})\) with \(w_1 > \dots > w_{2n} > 0\) integers, we have \(\alpha_\psi = \alpha_{m,d}\).
  Let \(\tilde{\alpha}_\psi: \Sp_{2m} \times \SL_2 \rightarrow \Mpsisc\) be the unique lift of \(\alpha_\psi\).
  As above it extends to give a morphism \(\tilde{\alpha}_\psi: G(\Sp_{2m} \times \SL_2) \to \GMpsisc\) mapping \(\lambda \in \GL_1\) to \(\lambda^n\).
\end{prodef}
\begin{proof}
  The standard representations of \(\Sp_{2m}\) and \(\SL_2\) are both self-dual of symplectic type and \(2d-1\) is odd so \(\Std_{\Sp_{2m}} \otimes \Sym^{2d-1} \Std_{\SL_2}\) is of orthogonal type.
  A simple weight computation shows that this representation factors through a morphism \(\alpha_{m,d}\) whose restriction to \(\Tcal_{\Sp_{2m}} \times \Tcal_{\SL_2}\) is given by \eqref{eq:alpha_m_d_on_T}.
  The conjugacy class of this morphism \(\Tcal_{\Sp_{2m}} \times \Tcal_{\SL_2} \to \SO_{4n}\) is not fixed by \(\thetahat\), e.g.\ because there exists \((\tau_1, \tau_2) \in \Lie ( \Tcal_{\Sp_{2m}} \times \Tcal_{\SL_2})\) mapping to \((2n, 2n-1, \dots, 1) \in \Lie \Tcal_{\SO_{4n}}\).
  Existence and uniqueness of \(\tilde{\alpha}_{m,d}\) is \cite[Proposition 2.24 (i)]{BorelTits_compl}.
  The restriction of \(\tilde{\alpha}_{m,d}: \Sp_{2m} \times \SL_2 \to \Spin_{4n}\) to \(Z(\Sp_{2m} \times \SL_2) \simeq \mu_2 \times \mu_2\) is easily computed:
  \[ (z_1, z_2) \mapsto (z_1 z_2, \dots, z_1 z_2, z_1^n) \in \Tcal_{\Spin_{4n}}. \]
  It follows that \eqref{eq:Sp2m_SL2_GL1_to_GSpin4n} is trivial on the diagonally embedded \(\mu_2\).

  The case of \(\tilde{\alpha}_\psi\) follows immediately, using the fact that \(\tau_\psi\) is not fixed by \(\thetahat\).
\end{proof}

\begin{prop} \label{pro:no_amb_symplectic}
  For \(n \geq 1\), \(\taut \in \ICcalt(\SObf_{4n})\) and \(\pi[2d] \in \tilde{\Psi}^{\unr, \tilde{\tau}}_{\disc, \nonendo}(\SObf_{4n})\) and for any prime \(p\) we have \(c_p(\pi[2d]) = \alpha_\psi(c_p(\pi), \diag(p^{1/2}, p^{-1/2}))\).
\end{prop}

Note that by \cite[\S 3.5]{Clozel_AA} for a pair \((\pi, d)\) as above there is a
finite extension \(E\) of \(\Q\) in \(\C\) such that for any prime number \(p\) the
semisimple conjugacy class \(c_p(\pi) p^{1/2}\) in \(\GL_{2m}(\C)\) is defined over
\(E\).
In particular the corresponding conjugacy class in \(\GSp_{2m}(\C)\) having
similitude character \(p\) is also defined over \(E\).
Therefore \(\alpha_{\pi[2d]}(c_p(\pi), \diag(p^{1/2}, p^{-1/2}))\) is also defined over \(E\)
\footnote{In fact the map \(\pi_p \mapsto c(\pi_p) \otimes \diag(p^{(2d-1)/2}, \dots, p^{(1-2d)/2}))\) from unramified representations of \(\GL_{2m}(\Qp)\) to semisimple conjugacy classes in \(\GL_{4n}(\C)\) is defined over \(\Q\).}.

\begin{proof}
  We already know that \(c_p(\pi[2d])\) is either equal to \(\alpha_{\pi[2d]}(c_p(\pi), \diag(p^{1/2}, p^{-1/2}))\) or its image by \(\thetahat\) (if it is not \(\thetahat\)-invariant, but the proof below will show that it never is, i.e.\ that no eigenvalue in the standard representation is \(\pm 1\)).

  Consider \(\psi = 1 \oplus \pi[2d] \in \tilde{\Psi}^{\unr, \tilde{\tau}}_{\disc, \nonendo}(\Spbf_{4n})\).
  With notation as in Theorem \ref{theo:IH_explicit_crude} we have \(\psi_0 = 1\), \(\psi_1 = \pi[2d]\) and \(u_1(\psi)\) is the sign \(\epsilon\) made explicit in the second point of Corollary \ref{coro:seed_ex_sigma}.

  As in Definition \ref{def:alpha_m_d} we may assume \(\Mpsi = \SO_{4n}\) and \(\tau_\psi\) equal to the conjugacy class of \((w_1, \dots, w_{2n}) \in \Lie \Tcal_{\SO_{4n}}\) where \(w_1 > \dots > w_{2n} > 0\).
  Using the parametrizations for maximal tori introduced in Section \ref{sec:not_red_gps}, the morphism \(\alpha_{\pi[2d]}\) maps the conjugacy class of \(((x_1, \dots, x_m), t) \in \Tcal_{\Sp_{2m}} \times \Tcal_{\Sp_2}\) to the conjugacy class of
  \[ (x_1 t^{2d-1}, x_1 t^{2d-3}, \dots, x_1 t^{1-2d}, \dots, x_m t^{2d-1}, \dots, x_m t^{1-2d}) \in \Tcal_{\SO_{4n}}. \]
  The two preimages of this element of \(\Tcal_{\SO_{4n}}\) in \(\Tcal_{\Spin_{4n}}\) are
  \[ y_{\pm} = (x_1 t^{2d-1}, x_1 t^{2d-3}, \dots, x_1 t^{1-2d}, \dots, x_m t^{2d-1}, \dots, x_m t^{1-2d}, \pm (x_1 \dots x_m)^d) \]
  Now we fix some prime number \(p\) and take \((x_1, \dots, x_m)\) to be a representative of the conjugacy class \(c_p(\pi)\) in \(\Sp_{2m}(\C)\) and \(t = p^{1/2}\).
  By \cite{Clozel_purity} or \cite{Caraiani} we know that all \(x_i\)'s have absolute value one (this also holds in any embedding of the number field \(\Q(x_1, \dots, x_m)\) in \(\C\)).
  Therefore any eigenvalue of \(\spin^+_{\pi[2d]}(y_+)\) or \(\spin^+_{\pi[2d]}(y_-)\) (resp.\ \(\spin^-_{\pi[2d]}(y_+)\) or \(\spin^-_{\pi[2d]}(y_-)\)) has absolute value \(p^{i/2}\) for some even (resp.\ odd) integer \(i\).
  Consider the sign
  \[ (-1)^{n(n+1)/2+(d_0^2-1)/8} \prod_{i=1}^r u_i(\psi)^{d_i-1} = (-1)^{n(n+1)/2} u_1(\psi) \]
  corresponding to \(\psi\) in \eqref{eq:explicitIH}.
  Choose any prime number \(\ell \neq p\) and an arbitrary isomorphism \(\iota\) as in Theorem \ref{theo:IH_explicit_crude}, and let \(V\) be the irreducible algebraic representation of \(\PGSpbf_{2n}\) corresponding to the infinitesimal character \(\tau\) of \(\psi\).
  By purity of \(H^\bullet((\Acal_n^*)_{\Qbar}, \IC(V))\), if \(u_1(\psi) = +1\) (resp.\ \(u_1(\psi) = -1\)) then any eigenvalue of \(\spin_{\pi[2d]}^{u_1(\psi)}(\cpsc(\pi[2d]))\) has absolute value \(p^{i/2}\) for some even (resp.\ odd) integer \(i\).
  Since \(\cpsc(\pi[2d])\) equals the Weyl orbit of \(y_+\), \(y_-\), \(\thetahat(y_+)\), or \(\thetahat(y_-)\) we can rule out the last two possibilities, and conclude that \(c_p(\pi[2d]) = \alpha_{\pi[2d]}(c_p(\pi), \diag(p^{1/2}, p^{-1/2}))\).
\end{proof}

\begin{rema} \label{rem:pi_2d_lift}
  Of course we expect a stronger relation
  \[ \cpsc(\pi[2d]) = \tilde{\alpha}_{\pi[2d]}(c_p(\pi), \diag(p^{1/2}, p^{-1/2})) \]
  (in the notation of the proof above, \(y_+\) rather than \(y_-\)).
  This could perhaps be proved by realizing a level one representation of \(\PGSObf_{4n}\) corresponding to \(\pi[2d]\) as iterated residues of Eisenstein series.
  We will prove this in Proposition \ref{pro:Gal_GSpin_pi_2d} using Galois-theoretic arguments.
\end{rema}

\subsection{Intersection versus compactly supported cohomology}
\label{sec:IH_vs_Hc_Morel}

The aim of this section is to express, for an irreducible algebraic representation \(V\) of \(\GSpbf_{2n,\Qell}\), the Euler characteristic \(e_{\IH}(\GSpbf_{2n}, \Xcal_n, V)\) (notation as in the introduction, \S \ref{sec:intro_Hc_IH}) of
\[ \varinjlim_K H^\bullet((\Acal_{n,K}^*)_{\Qbar}, \IC^K_\Q(V)) \]
in terms of
\begin{itemize}
\item the analogous Euler characteristics \(e_c(\GSpbf_{2n'}, \Xcal_{n'}, V')\) for compactly supported cohomology on \((\Acal_{n,K})_{\Qbar}\), for \(n' \leq n\) and certain representations \(V'\),
\item the Euler characteristics of the ``compactly supported'' cohomology (with coefficients) of certain arithmetic subgroups of \(\GLbf_{n'}\), for \(n' \leq n\).
\end{itemize}
This relation is a slightly more conceptual reformulation of results of Morel (\cite[Théorème 5.2.2]{MorelSiegel1} or \cite[Proposition 1.5.3]{MorelBook}) that she used to prove Theorem \ref{t:ICMorel}.
It serves two purposes: first to obtain in Corollary \ref{cor:sigma_IH_crys} crystallineness of the representations \(\sigma_{\psi,\iota}^{\IH}\) of Theorem \ref{theo:IH_explicit_crude}, and later to derive an explicit formula for \(e_c(\GSpbf_{2n}, \Xcal_n, V)\) in terms of intersection cohomology (as announced in \S \ref{sec:intro_Hc_IH} of the introduction).

We first recall the adelic version of group cohomology of arithmetic groups in Section \ref{sec:arith_gp_coh}, and spell out the notion of parabolic induction in the setting of Definition \ref{def:hecke_cat} in Section \ref{sec:para_ind}.
We then review boundary strata of the minimal compactifications \(\Acal_{n,K}^*\) and Morel's weight truncation, before following some of Morel's arguments to express the Euler characteristic of intersection cohomology in terms of that of ordinary cohomology (Corollary \ref{cor:IH_vs_H_Q}).
We then dualize to replace ordinary cohomology by compactly supported cohomology (Corollary \ref{coro:IH_vs_Hc_Q}).
Finally we deduce in Corollary \ref{cor:sigma_IH_crys} that each \(\sigma_{\psi,\iota}^{\IH}\) is crystalline from this relation and a theorem of Faltings and Chai.

\subsubsection{Arithmetic group cohomology}
\label{sec:arith_gp_coh}

Let \(\Gbf\) be a connected reductive group over \(\Q\), \(K_{\infty}\) an open
subgroup of a maximal compact subgroup \(K_{\infty}^{\max}\) of \(\Gbf(\R)\) and
denote by \(\Xcal\) the real manifold \(\Gbf(\R) / K_{\infty} \Abf_{\Gbf}(\R)^0\).
Consider a representation of \(\Gbf(\Q)\) on a finite-dimensional vector space \(V\)
over a field \(F\) of characteristic zero.
Exactly like in the case \(\Gbf = \GSpbf_{2n}\) (see the end of Section
\ref{sec:local_systems}), we have local systems \(\Fcal^K(V)\) on the manifolds
\(\Gbf(\Q) \backslash (\Xcal \times \Gbf(\A_f)/K)\) and we get objects \(H^i(\Gbf,
K_{\infty}, V)\) (resp.\ \(H^i_c(\Gbf, K_{\infty}, V)\)) of \(\Hecke(\Gbf(\A_f), F)\)
with
\[ H^i(\Gbf, K_{\infty}, V)^K = H^i(\Gbf(\Q) \backslash (\Xcal \times \Gbf(\A_f)
  / K), \Fcal^K(V)), \]
\[ H^i_c(\Gbf, K_{\infty}, V)^K = H^i_c(\Gbf(\Q) \backslash (\Xcal \times
  \Gbf(\A_f) / K), \Fcal^K(V)) \]
for any neat compact open subgroup \(K\) of \(\Gbf(\A_f)\).
For \(K_{\infty}\) maximal these will be simply denoted by \(H^i(\Gbf, V)\) (resp.\
\(H^i_c(\Gbf, V)\))
\footnote{We hope that this notation will not create any confusion since
cohomology of algebraic groups does not appear in this article.}.
Recall that thanks to the existence of ``nice'' compactifications, such as the
Borel-Serre compactification, the vector spaces \(H^i(\Gbf, K_{\infty}, V)^K\) and
\(H^i_c(\Gbf, K_{\infty}, V)^K\) have finite dimension and vanish for \(i\) outside
an explicit finite interval (see \cite{Raghunatan}, \cite[\S
11]{BorelSerre_corners}).

The Hecke operators between the cohomology groups \(H^i(\Gbf, K_{\infty} V)^K\)
are defined similarly to the algebro-geometric case of Shimura varieties: for \(g
\in \Gbf(\A_f)\) and \(K'\) an open compact subgroup of \(\Gbf(\A_f)\) contained in
\(gKg^{-1}\), multiplication by \(g\) defines a cover \(T_{K', g, K} : \Gbf(\Q)
\backslash (\Xcal \times \Gbf(\A_f) / K') \rightarrow \Gbf(\Q) \backslash (\Xcal
\times \Gbf(\A_f) / K)\) and there is a natural isomorphism \(T_{K', g, K}^*
\Fcal^K(V) \simeq \Fcal^{K'}(V)\), which is analogous to the composition of the
first two isomorphisms in Definition \ref{def:Hecke_corr}.
There is also an integral version that we will not use.

By Poincaré duality we have perfect pairings
\begin{equation} \label{eqn:pairing_Hecke_topo}
  H^i_c(\Gbf(\Q) \backslash (\Xcal \times \Gbf(\A_f) / K),
  \Fcal^K(V)) \times H^{\dim \Xcal - i}(\Gbf(\Q) \backslash (\Xcal \times
  \Gbf(\A_f) / K), \Fcal^K(V^*) \otimes \Ocal ) \rightarrow \Q
\end{equation}
where \(\Ocal\) is the orientation sheaf.
As in the algebro-geometric setting \eqref{eq:duality_Hecke} it is natural to
normalize this pairing using a Haar measure on \(\Gbf(\A_f)\) in order to realize
these two admissible representations of \(\Gbf(\A_f)\) (as \(K\) varies) as the
contragredient of each other.
The connected components of \(\Xcal\) are simply connected, so choosing an
orientation of \(\Xcal\) gives us an isomorphism \(\Ocal \simeq \Fcal^K(\chi)\)
where \(\chi\) is the restriction to \(\Gbf(\Q)\) of the continuous character
\(\Gbf(\R) \rightarrow \{ \pm 1 \}\) whose restriction to \(K_{\infty}^{\max}\) is
the determinant of the adjoint representation \(\Lie \Gbf(\R) / \Lie K_{\infty}\)
(or equivalently of \(\Lie K_{\infty}\), because the adjoint action of \(\Gbf(\R)\)
has trivial determinant).
In general this character \(\chi\) of \(\Gbf(\Q)\) is non-trivial, in particular
non-algebraic.

\begin{exam} \label{exam:orientation_GL}
  For \(\Gbf = \GLbf_{N, \Q}\) we have \(\chi = \sign \det^{N-1}\) and \(\dim \Xcal = N^2-1-N(N-1)/2 = N(N+1)/2-1\).
\end{exam}

\begin{rema} \label{rema:no_can_orientation}
  Using \cite[Proposition 2.2]{RohlfsSchwermer_inters} one can check that, at
  least when the level \(K\) is small enough, the manifold \(\Gbf(\Q) \backslash (
  \Xcal \times \Gbf(\A_f) / K)\) is orientable.
  Unfortunately it is not canonically so, and so as \(K\) varies it does not seem
  to be possible to choose orientations uniformly so that they are compatible
  with all finite étale covers \(T_{K', g, K}\).
  %TODO: \(\GLbf_N\) in level one gives a counter-example?
\end{rema}

Even though \(\chi\) is non-algebraic in general, we can often reduce to the case of an algebraic representation \(V\), by the following remark.
See also Lemma \ref{lem:fd_rep_real_gp_alg}.

\begin{rema} \label{rema:twisting_top_coh}
  Let \(L\) be a one-dimensional representation of \(\Gbf(\Q)\) over \(F\), and denote by \(\chi : \Gbf(\Q) \rightarrow F^{\times}\) the corresponding character.
  Let \(\Gcal \subset \Gbf(\R)\) be the stabilizer of a connected component \(\Xcal^0\) of \(\Xcal\).
  It is a normal subgroup of \(\Gbf(\R)\), in fact the quotient is commutative and \(2\)-torsion, and so \(\Gcal\) does not depend on the choice of a connected component of \(\Xcal\).
  Assume that there exists a locally constant character \(\widetilde{\chi}_f : \Gbf(\A_f) \rightarrow F^{\times}\) such that \(\widetilde{\chi}_f\) and \(\chi\) coincide on \(\Gcal \cap \Gbf(\Q)\).
  Then for a level \(K \subset \Gbf(\A_f)\) such that \(\widetilde{\chi}_f|_K = 1\), choosing a basis \(e\) of \(L\) there is a unique global section \(s(e, \widetilde{\chi}_f)\) of \(\Fcal^K(L)\) which on \((\Gcal \cap \Gbf(\Q)) \backslash (\Xcal^0 \times \Gbf(\A_f) / K)\) is given by \((x,hK) \mapsto \widetilde{\chi}_f(h) e\).
  Cup-product with \(s(e, \widetilde{\chi}_f)\) gives an isomorphism of admissible representations of \(\Gbf(\A_f)\)
  \[ H^i(\Gbf, K_{\infty}, V) \otimes \widetilde{\chi}_f \simeq H^i(\Gbf, K_{\infty}, V \otimes_F L). \]
  Similarly we have an isomorphism \(H^i_c(\Gbf, K_{\infty}, V) \otimes \widetilde{\chi}_f \simeq H^i_c(\Gbf, K_{\infty}, V \otimes_F L)\).

  Note that Remark \ref{rema:twisting_sim} over \(\C\) is a special case of this (via comparison of étale and singular cohomology), and that a simple way to find a pair \((\chi, \widetilde{\chi}_f)\) is to take \(\chi\) to be the
  restriction to \(\Gbf(\Q)\) of \(\widetilde{\chi}_\infty^{-1}\) where \(\widetilde{\chi} = \widetilde{\chi}_{\infty} \widetilde{\chi}_f\) is a character on \(\Gbf(\A)\) which is trivial on \(\Gbf(\Q)\).
\end{rema}

In this section we will only need the case where \(K_{\infty}\) is maximal, and for simplicity we make this assumption from now on.
We simplify the notation by denoting \(H^\bullet(\Gbf,-)\) for \(H^\bullet(\Gbf,K_\infty,-)\), and similarly for \(H^\bullet_c\) (by conjugacy of maximal compact subgroups of \(\Gbf(\R)\) these cohomology groups do not depend on the choice of \(K_\infty\) indeed).
The general case will be resumed in Section \ref{sec:Franke_formula}.

For a given neat compact open subgroup \(K\) of \(\Gbf(\A_f)\) we have a decomposition into (finitely many) connected components
\[ \Gbf(\Q) \backslash (\Xcal \times \Gbf(\A_f) / K) \simeq \bigsqcup_{[h_j] \in \Gbf(\Q) \backslash \Gbf(\A_f) / K} \Gamma_j \backslash \Xcal \]
where \(\Gamma_j = h_j K h_j^{-1} \cap \Gbf(\Q)\).
This gives an isomorphism
\[ H^i(\Gbf, V)^K \simeq \bigoplus_j H^i(\Gamma_j \backslash \Xcal, \Fcal^K(V)). \]
Since \(\Xcal\) is connected and contractible and the action of \(\Gamma_j\) on it is free we have canonical isomorphisms \(H^i(\Gamma_j \backslash \Xcal, \Fcal^K(V)) \simeq H^i(\Gamma_j, V)\) (see \cite[\S 5.3 Cor.\ 3]{Grothendieck_Tohoku}), giving a purely algebraic interpretation of \(H^i(\Gbf, V)\) in terms of group cohomology.
If we change representatives, say \(h_j' = \gamma_j h_j k_j\) with \(\gamma_j \in \Gbf(\Q)\) and \(k_j \in K\), then \(\Gamma_j' := h_j' K (h_j')^{-1} \cap \Gbf(\Q) = \gamma_j \Gamma_j \gamma_j^{-1}\), we have an isomorphism \(V|_{\Gamma_j'} \simeq V|_{\Gamma_j}\) induced by \(\gamma_j^{-1}\) which is compatible with the
isomorphism \(\Ad(\gamma_j) : \Gamma_j \rightarrow \Gamma_j'\), and so together they induce an isomorphism \(H^\bullet(\Gamma_j', V) \simeq H^\bullet(\Gamma_j, V)\) which of course depends on \((h_j, h_j')\) but not on the choice of \((\gamma_j, k_j)\) (see \cite[Ch.\ VII Prop.\ 3]{Serre_corpsloc}).
In particular we have a canonical isomorphism
\begin{equation} \label{eq:iso_coh_arith_gp}
  H^i(\Gbf, V)^K \simeq \bigoplus_{h \in \Gbf(\Q) \backslash \Gbf(\A_f) / K} \colim_{hK \in c} H^i(hKh^{-1} \cap \Gbf(\Q), V).
\end{equation}

Moreover the Hecke operators \([K_2, g, K_1, K']: H^i(\Gbf, V)^{K_1} \to H^i(\Gbf, V)^{K_2}\), at least when \(K_1\), \(K_2\) and \(K'\) are neat, can be rewritten via \eqref{eq:iso_coh_arith_gp} in terms of combinatorics of double quotients \(\Gbf(\Q) \backslash \Gbf(\A_f) / K\) and maps \(T_{K',g,K}\) between them and restriction and corestriction maps in group cohomology.
We do not make this more explicit here.

We denote
\[ e(\Gbf,V) = \sum_{i \geq 0} (-1)^i [H^i(\Gbf,V)],\ \text{resp. } e_c(\Gbf,V) = \sum_{i \geq 0} (-1)^i [H^i_c(\Gbf,V)] \]
in the Grothendieck group of admissible representations of \(\Gbf(\A_f)\) over \(\Q\).

\subsubsection{Parabolic induction}
\label{sec:para_ind}

\begin{defi} \label{def:hecke_ind}
  Let \(G\) be a locally profinite group, \(P\) a closed subgroup, and \(N\) a closed normal subgroup of \(P\) and denote \(M = P/N\).
  For \(p \in P\) denote by \(\ol{p}\) its image in \(M\).
  Assume that \(G/P\) is compact (so that for any compact open subgroup \(K\) of \(G\) the double quotient \(P \backslash G / K\) is finite).
  Let \(F\) be a field.
  For \(K\) an open compact subgroup of \(G\) and \(h \in G\) we define compact open subgroups of \(N\), \(P\) and \(M\) as follows: let \(K_{N,h} = h K h^{-1} \cap N\), \(K_{P,h} = h K h^{-1} \cap P\) and \(K_{M, h} = K_{P,h} / K_{N, h}\), which is isomorphic to the image of \(K_{P, h}\) in \(M\).
  For an object \(\ul{V} = (V_K)_{K \in C(M)}\) of \(\Hecke(M, F)\), let \(\ind_P^G \ul{V}\) be the object of \(\Hecke(G, F)\) with
  \begin{multline*}
    (\ind_P^G \ul{V})_K = \biggl\{ f : G / K \rightarrow \bigsqcup_{hK \in G / K} V_{K_{M,h}} \ \bigg|\  \forall hK \in G/K, f(hK) \in V_{K_{M,h}} \\
      \text{and } \forall hK \in G/K, \forall p \in P, f(p h K) = [K_{M,ph}, \ol{p}, K_{M,h}, K_{M,ph}] f(hK) \biggr\}
  \end{multline*}
  and Hecke operators \([K_2, g, K_1, K'] : (\ind_P^G \ul{V})_{K_1} \rightarrow (\ind_P^G \ul{V})_{K_2}\) defined by
  \begin{align} \label{eq:def_hecke_ind}
    & [K_2, g, K_1, K'](f)( h_2 K_2 ) \\
    = & \sum_{\substack{[h'] \in P \backslash G / K' \\ P h' K_2 = P h_2 K_2}} \left| \frac{(K_2)_{N, h'}}{(K')_{N, h'}} \right| \left[ (K_2)_{M, h_2}, \ol{p_2^{-1} p_1}, (K_1)_{M, h_1}, (K')_{M, p_2^{-1} h'} \right](f(h_1 K_1)) \nonumber
  \end{align}
  where \(h_1\) is any element of \(P h' g K_1\), \(p_1, p_2 \in P\) are such that \(p_1 h_1 \in h' g K_1\) and \(p_2 h_2 \in h' K_2\).
\end{defi}

\begin{rema} \label{rem:def_hecke_ind}
  \begin{enumerate}
  \item Similarly to the case of Definition \ref{def:hecke_cat}, one can replace the category of vector spaces over \(F\) by an arbitrary additive category.
  \item In the setting of the definition, choosing representatives \(h_1, \dots, h_m\) for \(P \backslash G / K\), we have an isomorphism
    \begin{align*}
      (\ind_P^G \ul{V})_K & \longrightarrow \bigoplus_{i=1}^m V_{K_{M, h_{i}}} \\
      f & \longrightarrow (f(h_i K))_{1 \leq i \leq m}
    \end{align*}
    because the Hecke operators \([K_{M,ph}, \ol{p}, K_{M,h}, K_{M,ph}]\) occurring in the definition are isomorphisms which compose in the obvious way and are equal to the identity when \(phK = hK\).
  \item \label{it:rem_def_hecke_ind_simpl}
    The same properties are used to checked that each term in the sum \eqref{eq:def_hecke_ind} does not depend on the choice of \((h', h_1, p_1, p_2)\).
    For simplicity one could take \(h' \in h_2 K_2\), \(h_1 = h' g\) and \(p_1 = p_2 = 1\).
  \item It is not obvious that \(\ind_P^G \ul{V}\) satisfies all axioms of Definition \ref{def:hecke_cat}.
    Of course it is the last axiom that demands more work.
    The proof is relatively straightforward but a bit long, so we leave it to the reader.
    In fact we will only need the case of \(\Q\)-vector spaces and \(\Q\)-linear categories that reduce to this case, in which case the axioms of Definition \ref{def:hecke_cat} for \(\ind_P^G \ul{V}\) follow from Proposition \ref{pro:Hecke_corr_sat_formalism} and Proposition \ref{pro:hecke_ind} below, which of course is the motivation for Definition \ref{def:hecke_ind}.
  \end{enumerate}
\end{rema}

\begin{prop} \label{pro:hecke_ind}
  If \(\ul{V}\) is associated to a smooth representation \(V\) of \(M\) by Proposition \ref{pro:eq_Hecke_smrep} then \(\ind_P^G \ul{V}\) is canonically associated to the (non-normalized) induced representation \(\ind_P^G V\) of \(G\).
\end{prop}
\begin{proof}
  The identification of \((\ind_P^G V)^K\) with \((\ind_P^G \ul{V})_K\) is straightforward, so let us check that Hecke operators match.
  Recall from Proposition \ref{pro:eq_Hecke_smrep} that \([K_2, g, K_1, K']\) is induced by \(\sum_{k \in K_2/K'} kg\).
  Fix \(h_2 \in G\) (and not just \(h_2 K_2 \in G/K_2\)).
  Each \(k \in K_2 / K'\) defines \([h_2 k] \in P \backslash G / K'\) mapping to \([h_2]\) in \(P \backslash G / K_2\), and all such double cosets in \(P \backslash G / K'\) are obtained in this way.
  Moreover for \(h' \in h_2 K_2\) and \(kK' \in K_2 / K'\) we have \(P h' k K' = P h' K'\) if and only if \(k\) belongs to the image of \((K_2)_{P, h_2} / (K')_{P, h'} \xhookrightarrow{\Ad(h')^{-1}} K_2 / K'\).
  Therefore for \(f \in (\ind_P^G V)^{K_1}\) we have
  \begin{align*}
    [K_2, g, K_1, K'](f)(h_2 K_2)
    &= \sum_{k \in K_2/K'} f(h_2 k g K_1) \\
    &= \sum_{[h'] \in (K_2)_{P, h_2} \backslash h_2 K_2 / K'} \ \sum_{k' \in (K_2)_{P, h'} / (K')_{P, h'}} f(k' h' g K_1)
  \end{align*}
  and since \(f\) is left \(P\)-equivariant we can write
  \begin{align*}
    \sum_{k' \in (K_2)_{P, h'} / (K')_{P, h'}} f(k' h' K_1)
    & = \sum_{k' \in (K_2)_{M, h'} / (K')_{M, h'}} \left| \frac{(K_2)_{N, k' h'}}{(K')_{N, k' h'}} \right| f(k' h' g K_1) \\
    & = \sum_{k' \in (K_2)_{M, h'} / (K')_{M, h'}} \left| \frac{(K_2)_{N, h'}}{(K')_{N, h'}} \right| k' \cdot f(h' g K_1) \\
    & = \left| \frac{(K_2)_{N,h'}}{(K')_{N,h'}} \right| \left[ (K_2)_{M,h_2}, 1, (K_1)_{M,h' g}, (K')_{M,h'} \right] f(h' g K_1)
  \end{align*}
  and we recognize the simplification of \eqref{eq:def_hecke_ind} observed in Remark \ref{rem:def_hecke_ind} \ref{it:rem_def_hecke_ind_simpl}.
\end{proof}

\begin{coro} \label{coro:dual_para_ind}
  Assume that we are in the setting of Definition \ref{def:hecke_ind}.
  Assume that \(F\) has characteristic zero and that \(G\) and \(N\) are unimodular.
  Let \(\delta_P : M \rightarrow \Q^{\times}\) be the modulus character, i.e.\ for any \(p \in P\), for any right Haar measure \(\mu\) on \(P\) and for any measurable set \(X \subset P\) we have \(\mu(pX) = \delta_P(\ol{p}) \mu(X)\).
  Fix right Haar measures on \(P\) and \(G\).
  Then the contragredient object (defined in Corollary \ref{cor:hecke_contra}) \((\ind_P^G \ul{V})^*\) is isomorphic to \(\ind_P^G (\ul{V}^* \otimes \delta_P)\).
\end{coro}
\begin{proof}
  Under the assumptions of unimodularity we have a ``quotient measure'' which is a morphism \(\ind_P^G \delta_P \to F\) of smooth representations of \(G\).
  Let \(V\) be the smooth representation of \(M\) corresponding to \(\ul{V}\), and denote by \(\wt{V}\) its contragredient representation, which naturally corresponds to \(\ul{V}^*\).
  The isomorphism in the Corollary is obtained by composing the obvious pairing \(\ind_P^G V \times \ind_P^G (\wt{V} \otimes \delta_P) \to \ind_P^G \delta_P\) with the quotient measure.
\end{proof}

As the notation suggests we will use Proposition \ref{pro:hecke_ind} in the case of parabolic induction, that is for \(G = \Gbf(\A_f)\), \(P = \Pbf(\A_f)\) and \(N = \Nbf(\A_f)\) where \(\Gbf\) is a connected reductive group over \(\Q\), \(\Pbf\) a parabolic subgroup and \(\Nbf\) its unipotent radical.

\subsubsection{Stratifications of minimal compactifications of Siegel modular varieties}

For the rest of Section \ref{sec:IH_vs_Hc_Morel} we denote \(\Gbf = \GSpbf_{2n}\) and \(\Xcal\) denotes the \(\Gbf(\R)\)-homogeneous space introduced in Section \ref{sec:def_An}.
Maximal proper standard parabolic subgroups of \(\Gbf\) are parametrized by integers \(1 \leq m \leq n\), we denote by \(\Pbf_m\) the block upper triangular subgroup of \(\Gbf\) corresponding to the partition \(n = m + 2(n-m) + m\).
Thus standard Levi subgroups of \(\Gbf\) are parametrized by tuples \(1 \leq n_1 < \dots < n_r \leq n\), and \(\Pbf = \Pbf_{n_1} \cap \dots \cap \Pbf_{n_r}\) has standard Levi factor \(\Mbf_\Pbf\) isomorphic to \(\GLbf_{n_1} \times
\dots \times \GLbf_{n_r-n_{r-1}} \times \GSpbf_{2(n-n_r)}\), the inverse isomorphism being
\[ (g_1, \dots, g_r, h) \longmapsto \diag \left( g_1 \nu(h), \dots, g_r \nu(h), h, J_{n_1}^{-1} {}^t g_1^{-1} J_{n_1} \dots, J_{n_r-n_{r-1}}^{-1} {}^t g_r^{-1} J_{n_r-n_{r-1}} \right) \]
where \((J_k)_{i,j} = \delta_{i+j=k+1}\).
In particular we have a decomposition
\begin{equation} \label{eq:decomp_lin_her_Siegel}
  \Mbf_\Pbf = \Mbf_{\Pbf,\lin} \times \Mbf_{\Pbf,\her}
\end{equation}
where \(\Mbf_{\Pbf,\lin} \simeq \GLbf_{n_1} \times \dots \times \GLbf_{n_r-n_{r-1}}\) and \(\Mbf_{\Pbf, \her} = \GSpbf_{2(n-n_r)}\).

We will use the description of boundary strata of minimal compactifications explained in Section \ref{sec:gen_Shim_boundary}, in terms of generalized Shimura varieties introduced in Definitions \ref{def:gen_Shim_datum} and \ref{def:gen_Shim_var}.
Boundary strata of \(\Acal_{n,K}^*\) are parametrized by maximal proper parabolic subgroups of \(\Gbf\).
For such a parabolic subgroup \(\Pbf = \Pbf_m\) (here \(1 \leq m \leq n\)) we have by Proposition \ref{pro:ext_Shim_parabolic} an associated generalized Shimura datum \((\Mbf_\Pbf, \Xcal_\Pbf, h_\Pbf)\) and for a neat level \(K = K^p \times \Gbf(\Zp)\) the stratum of \(\Acal_{n,K,\Fp}^* = \Sh(\Gbf, \Xcal, K)^*_{\Fp}\) corresponding to \(\Pbf\) may be identified with
\[ \underset{gK \in [\Pbf(\A_f^{(p)}) \curvearrowright \Gbf(\A_f^{(p)})/K^p]}{\colim} \Sh(\Mbf_\Pbf, \Xcal_\Pbf, K(\Pbf,gK))_{\Fp} \simeq \bigsqcup_{[g] \in \Pbf(\A_f^{(p)}) \backslash \Gbf(\A_f^{(p)})/K^p} \Sh(\Mbf_\Pbf, \Xcal_\Pbf, K(\Pbf,gK))_{\Fp} \]
and we denote
\[ i_{\Pbf,gK}: \Sh(\Mbf_\Pbf, \Xcal_\Pbf, K(\Pbf,K))_{\Fp} \hookrightarrow \Acal_{n,K,\Fp}^* \]
the locally closed immersion.

More generally for a standard parabolic subgroup (not assumed to be maximal or proper) \(\Pbf = \Pbf_{n_1} \cap \dots \cap \Pbf_{n_r}\) of \(\Gbf\) we may restrict (in the sense of Proposition-Definition \ref{prodef:linres_ext_Shim}) the generalized Shimura datum \((\Mbf_{\Pbf_{n_r}}, \Xcal_{\Pbf_{n_r}}, h_{\Pbf_{n_r}})\) to \(\Mbf_\Pbf\), to obtain a generalized Shimura datum \((\Mbf_\Pbf, \Xcal_\Pbf, h_\Pbf)\).
As explained in Section \ref{sec:gen_Shim_boundary_boundary} the double coset
\[ \Pbf(\A_f^{(p)}) \backslash \Gbf(\A_f^{(p)}) / K^p \simeq \Pbf(\A_f) \backslash \Gbf(\A_f) / K \]
parametrizes sequences of boundary strata corresponding to the parabolic subgroups \(\Pbf_{n_1}\), \(\Pbf_{n_1} \cap \Pbf_{n_2}\), \dots, \(\Pbf\).
In particular for \(g \in \Gbf(\A_f^{(p)}) \times \Gbf(\Zp)\) we have a morphism
\[ T_{\Pbf, gK}: \Sh(\Mbf_\Pbf, \Xcal_\Pbf, K(\Pbf, gK))_{\Fp} \longrightarrow \Sh(\Gbf, \Xcal, K)^*_{\Fp}. \]
If \(r=1\) this map coincides with the locally closed immersion \(i_{\Pbf,gK}\) but in general \(T_{\Pbf,gK}\) is only finite étale over the boundary stratum corresponding to \(\Pbf_{n_r}\) and \([g] \in \Pbf_{n_r}(\A_f^{(p)}) \backslash \Gbf(\A_f^{(p)}) / K^p\).
If \(n_r<n\) (resp.\ \(n_r=n\)) then \(h_\Pbf: \Xcal_\Pbf \to \Hom(\mathbb{S}, \Mbf_{\Pbf,\R})\) is injective (resp.\ two-to-one) and \(h_\Pbf(\Xcal_\Pbf)\) is the \(\Mbf_\Pbf(\R)\)-orbit of
\[ z \mapsto (I_{n_1}, \dots, I_{n_r-n_{r-1}} h_{n-m,0}(z)) \]
where \(h_{n-m,0} : \mathbb{S} \to \GSpbf_{2(n-m)}\) is as described around \eqref{eq:desc_anis_tor_GSp}.
In any case we are in the situation of Section \ref{sec:ext_Shim_direct_prod}: \(\Mbf_{\Pbf,\lin}(\R)\) acts trivially on \(\Xcal_\Pbf\) and \((\Mbf_{\Pbf,\her}, \Xcal_\Pbf, h_\Pbf)\) is a Shimura datum, so each \(\Sh(\Mbf_\Pbf, \Xcal_\Pbf, K)\) is isomorphic to a finite disjoint union of \(\Acal_{n-m, K'}\) for certain neat levels \(K'\).
For \(\ell \neq p\) prime and \(V\) a bounded complex of (finite-dimensional) algebraic representations  of \(\Mbf_{\Pbf,\Qell}\) we also have (see Proposition-Definition \ref{prodef:AFcal}) objects \(\AFcal^K V\) of \(D^+(\Sh(\Mbf_\Pbf, \Xcal_\Pbf, K)_{\Fp}, \Qell)\) and cohomological correspondences between them (Definition \ref{def:corr_au}), yielding objects
\[ H^i(\Sh(\Mbf_\Pbf, \Xcal_\Pbf, ?)_{\Fpbar}, \AFcal^? V) \]
of \(\Rep_{\Qell}^{\adm,\cont}(\Mbf_\Pbf(\A_f^{(p)}) \times \GalQ)\).
As explained in Section \ref{sec:ext_Shim_direct_prod} the corresponding Euler characteristic in \(K_0(\Rep_{\Qell}^{\adm,\cont}(\Mbf_\Pbf(\A_f^{(p)}) \times \GalQ))\) factorizes as follows when \(V\) is concentrated in degree zero and decomposes as \( = V_\lin \otimes V_\her\)
\[ e(\Sh(\Mbf_\Pbf, \Xcal_\Pbf, ?)_{\Fpbar}, \AFcal^? V) = e(\Mbf_{\Pbf,\lin}, V_\lin) \boxtimes e(\Acal_{n-n_r,?,\Fpbar}, \Fcal^? V_\her). \]

\subsubsection{Weight truncation of correspondences}

Morel \cite[\S 3.3]{MorelSiegel1} associated to the stratification of \((\Acal_{n,K}^*)_{\Fp} = \Sh(\Gbf, \Xcal, K)^*_{\Fp}\) recalled in the previous section and a tuple \(\ul{a} \in (\Z \cup \{\pm \infty\})^{n+1}\) a t-structure on the triangulated category of mixed complexes in \(D^b(\Sh(\Gbf, \Xcal, K)^*_{\Fp}, \Qell)\).
We will use the same notation, e.g.\ \(w_{\leq \ul{a}}\) denotes truncation for this t-structure.

\begin{prop} \label{pro:weight_trunc_PreH}
  Let \(\ul{a} \in (\Z \cup \{\pm \infty\})^{n+1}\).
  Let \(((L_K)_{K \in C}, (v(K_2,g,K_1,K'))_{K_2,g,K_1,K'})\) be an object of \(\PreH(\Acal^*_{n,?,\Fp}, C, \Qell)\) (Definition \ref{def:hecke_axioms_derived}) such that each \(L_K\) is bounded and mixed.
  Then
  \begin{align*}
    & ((w_{\leq \ul{a}} L_K)_{K \in C}, (w_{\leq \ul{a}} v(K_2,g,K_1,K'))_{K_2,g,K_1,K'}) \\
    \text{and } & ((w_{> \ul{a}} L_K)_{K \in C}, (w_{> \ul{a}} v(K_2,g,K_1,K'))_{K_2,g,K_1,K'})
  \end{align*}
  are also objects of \(\PreH(\Acal^*_{n,?,\Fp}, C, \Qell)\).
\end{prop}
\begin{proof}
  We only treat the first case (\(w_{\leq \ul{a}}\)) as the other case (\(w_{> \ul{a}}\)) is entirely similar.
  The proof of the first two axioms in Definition \ref{def:hecke_axioms_derived} is straightforward and we omit it.
  For the third and fourth axioms it is enough to check compatibility of weight truncation with composition of correspondences and pushforward of correspondences.

  For composition, we need to prove the following: assume that we have morphisms
  \[ X_1 \xleftarrow{c_1} X' \xrightarrow{c_2} X_2 \xleftarrow{d_1} X'' \xrightarrow{d_2} X_3 \]
  between schemes separated of finite type over \(\Fp\) and endowed with stratifications compatible with these morphisms.
  Let \(L_i \in D^b_m(X_i, \Qell)\) for \(i=1,2,3\), \(u: c_1^* L_1 \to c_2^! L_2\) and \(v: d_1^* L_2 \to d_2^! L_3\).
  Then we have an equality of correspondences between \(w_{\leq \ul{a}} L_1\) and \(w_{\leq \ul{a}} L_3\) supported on \((c_1 \pi', d_2 \pi'')\), where \(\pi': X' \times_{X_2} X'' \to X'\) and \(\pi'': X' \times_{X_2} X'' \to X''\):
  \begin{equation} \label{eq:weight_trunc_compat_comp_corr}
    w_{\leq \ul{a}}(v \circ u) = (w_{\leq \ul{a}} v) \circ (w_{\leq \ul{a}} u).
  \end{equation}
  In the diagram
  \[
    \begin{tikzcd}[column sep=huge]
      (\pi')^* c_1^* L_1 \ar[r, "{(\pi)^* u}"] & (\pi')^* c_2^! L_2 \ar[r] & (\pi'')^! d_1^* L_2 \ar[r, "{(\pi'')^! v}"] & (\pi'')^! d_2^! L_3 \\
      (\pi')^* c_1^* w_{\leq \ul{a}} L_1 \ar[u] \ar[r, "{(\pi')^* w_{\leq \ul{a}} u}"] & (\pi')^* c_2^! w_{\leq \ul{a}} L_2 \ar[u] \ar[r] & (\pi'')^! d_1^* w_{\leq \ul{a}} L_2 \ar[u] \ar[r, "{(\pi'')^! w_{\leq \ul{a}} v}"] & (\pi'')^! w_{\leq \ul{a}} L_3 \ar[u]
    \end{tikzcd}
  \]
  the left (resp.\ right) square is commutative by application of \((\pi')^*\) (resp.\ \((\pi'')^!\)) to the commutative square characterizing \(w_{\leq \ul{a}} u\) (resp.\ \(w_{\leq \ul{a}} v\)), and the middle square is commutative because it is obtained from the morphism of functors \((\pi')^* c_2^! \to (\pi'')^! d_1^*\) applied to the morphism \(w_{\leq \ul{a}} L_2 \to L_2\).
  So the outer square is commutative, implying \eqref{eq:weight_trunc_compat_comp_corr}.

  For pushforward we need to prove the following: suppose we have a diagram
  \[
    \begin{tikzcd}
      & X' \ar[d, "{f}"] & \\
      X_1 & \ar[l, "{c_1}"] X \ar[r, "{c_2}"] & X_2
    \end{tikzcd}
  \]
  of separated schemes of finite type over \(\Fp\) endowed with stratifications compatible with these morphisms, \(L_i \in D^b_m(X_i, \Qell)\) for \(i=1,2\) and \(u: f^* c_1^* L_1 \to f^! c_2^! L_2\).
  Assume that \(f\) is proper.
  Then we have an equality of correspondences from \(w_{\leq \ul{a}} L_1\) to \(w_{\leq \ul{a}} L_2\) supported on \((c_1,c_2)\):
  \begin{equation} \label{eq:weight_trunc_compat_push_corr}
    \corr f_* (w_{\leq \ul{a}} u) = w_{\leq \ul{a}} (\corr f_* u).
  \end{equation}
  In the diagram
  \[
    \begin{tikzcd}[column sep=huge]
      c_1^* L_1 \ar[r, "{\adj}"] & f_* f^* c_1^* L_1 \ar[r, "{f_* u}"] & f_* f^! c_2^! L_2 \ar[r, "{\adj}"] & c_2^! L_2 \\
      c_1^* w_{\leq \ul{a}} L_1 \ar[r, "{\adj}"] \ar[u] & f_* f^* c_1^* w_{\leq \ul{a}} L_1 \ar[r, "{f_* w_{\leq \ul{a}} u}"] \ar[u] & f_* f^! c_2^! w_{\leq \ul{a}} L_2 \ar[r, "{\adj}"] \ar[u] & c_2^! w_{\leq \ul{a}} L_2
    \end{tikzcd}
  \]
  the left and right square are commutative because they arise from a morphism of functors (\(\id \to f_* f^*\) resp.\ \(f_! f^! \to \id\)) applied to a morphism \(w_{\leq \ul{a}} L_i \to L_i\), and the middle square is commutative by application of \(f_*\) to the commutative square characterizing \(w_{\leq \ul{a}} u\).
  So the outer square is commutative, implying \eqref{eq:weight_trunc_compat_push_corr}.
\end{proof}

\begin{coro} \label{coro:weight_trunc_push_PreH}
  Let \(\ul{a}^{(1)}, \dots, \ul{a}^{(r)} \in (\Z \cup \{\pm \infty\})^{n+1}\).
  For each \(1 \leq i \leq r\) let \(\tau_i\) be either \(w_{\leq \ul{a}^{(i)}}\) or \(w_{> \ul{a}^{(i)}}\).
  Then the pair
  \[ (\tau_r \dots \tau_1 j_* \Fcal^K(V))_K, \ (\tau_r \dots \tau_1 \corr \ul{j}_* u(K_2, g, K_1, K')_{\Fp})_{K_2,g,K_1,K'} \]
  is an object of \(\PreH(\Acal^*_{n,?,\Fp}, C, \Qell)\), where \(\corr \ul{j}_*\) denotes the pushforward of correspondences (Definition \ref{def:gen_push_corr}) for the open immersions \(j: \Acal_{n,K,\Fp} \to \Acal_{n,K,\Fp}^*\).
\end{coro}
\begin{proof}
  This will follow from \(r\) applications of Proposition \ref{pro:weight_trunc_PreH} once we prove that
  \[ (j_* \Fcal^K(V))_K, \ (\corr \ul{j}_* u(K_2, g, K_1, K')^{\IC}_{\Fp})_{K_2,g,K_1,K'} \]
  is an object of \(\PreH(\Acal^*_{n,?,\Fp}, C, \Qell)\).
  This is proved exactly as for Proposition \ref{pro:IH_is_Hecke}, using (repeatedly) the first point of Lemma \ref{lemm:uniq_int_ext_corr} instead of the second one.
\end{proof}

\subsubsection{Intersection cohomology from ordinary cohomology}

\begin{defi} \label{def:trunc_alg_rep}
  Let \(\Pbf = \Pbf_{n_1} \cap \dots \cap \Pbf_{n_r}\) be a standard parabolic subgroup of \(\Gbf\).
  Let \(t_{n_1}, \dots, t_{n_r} \in \Z\).
  Let \(V\) be an algebraic (finite-dimensional) representation of \(\Mbf_\Pbf \simeq \GLbf_{n_1} \times \GLbf_{n_2-n_1} \times \dots \times \GLbf_{n_r-n_{r-1}} \times \GSpbf_{2(n-n_r)}\).
  We have a canonical decomposition
  \[ V = \bigoplus_{s_1, \dots, s_r \in \Z} V_{s_1, \dots, s_r} \]
  where \(V_{s_1, \dots, s_r}\) is the eigenspace for \(\Zbf(\GLbf_{n_1}) \times \dots \times \Zbf(\GLbf_{n_r-n_{r-1}})\) such that for each \(1 \leq i \leq r\) and \(\lambda \in \GLbf_1\),
  \[ \left( \lambda I_{n_1}, \dots, \lambda I_{n_i-n_{i-1}}, I_{n_{i+1}-n_i}, \dots, I_{n_r-n_{r-1}} \right) \]
  acts by \(\lambda^{s_i}\) on \(V_{s_1, \dots, s_r}\).
  Denote
  \[ V_{<\ul{t}} = V_{<t_{n_1}, \dots, <t_{n_r}} = \bigoplus_{s_1 < t_{n_1}, \dots, s_r < t_{n_r}} V_{s_1, \dots, s_r} \]
  and similarly for \(V_{\leq \ul{t}}\), \(V_{>\ul{t}}\) and \(V_{\geq \ul{t}}\).
  Being functorial, these truncations extend to complexes of algebraic representations of \(\Mbf_\Pbf\).
\end{defi}

\begin{rema} \label{rem:compare_trunc_alg_rep}
  This definition of truncation differs from the one in \cite[\S 4.2]{MorelSiegel1}.
  More precisely if the center of \(\GSpbf_{2n}\) acts by \(t \mapsto t^m\) on \(V\) then \(V_{s_1,\dots,s_r}\) is the largest subspace of \(V\) on which
  \[ \diag(\lambda^2 I_{n_i}, \lambda I_{2(n-n_i)}, I_{n_i}) \in \Zbf(\Mbf_{\Pbf}) \]
  acts by \(\lambda^{m+s_i}\) for each \(1 \leq i \leq r\).
  So our \(V_{<\ul{t}}\) is Morel's \(V_{<\ul{t'}}\) where \(t'_i = t_i+m\).
  We translated Morel's conditions so that they become invariant under twisting by characters of \(\Gbf\).

  Note that there seems to be a typographical error (certainly related to Remark \ref{rem:cond3_varphi_h_Pink}) in \cite[\S 4.2]{MorelSiegel1} for the case \(n_i = n\) (corresponding to zero-dimensional strata).
\end{rema}

\begin{theo}[Morel] \label{thm:IH_vs_H_Fp}
  Let \(\ell \neq p\) be prime numbers, \(V\) an algebraic representation of \(\Gbf_{\Qell}\) (or a bounded complex of such representations).
  In \(K_0(\Rep_{\Qell}^{\adm,\cont}(\Gbf(\A_f^{(p)}) \times \Gal(\Fpbar/\Fp)))\) we have
  \[ e(\Acal_{n,?,\Fpbar}^*, \IC^?(V)) = \sum_\Pbf (-1)^{r_\Pbf} \ind_{\Pbf(\A_f^{(p)})}^{\Gbf(\A_f^{(p)})} e(\Sh(\Mbf_\Pbf, \Xcal_\Pbf, ?)_{\Fpbar}, \AFcal^? R\Gamma(\Lie \Nbf_\Pbf, V)_{<\ul{t}}) \]
  where the sum ranges over standard parabolic subgroups \(\Pbf = \Pbf_{n_1} \cap \dots \cap \Pbf_{n_r}\) of \(\Gbf\), \(r_\Pbf = r = \dim \Abf_{\Pbf} - \dim \Abf_\Gbf\) and \(t_{n_i} = (n-n_i)(n-n_i+1)/2 - n(n+1)/2\) for \(1 \leq i \leq r\).
\end{theo}
\begin{proof}
  This is a reformulation of \cite[Théorème 5.2.2]{MorelSiegel1} (see also the dual version \cite[Proposition 1.5.3]{MorelBook}), after taking cohomology.
  Since we have adopted a different formulation in order to make the appearance of parabolic induction more obvious, let us briefly explain how this theorem follows from Morel's results.
  We use the set \(C\) of compact open subgroups \(K^p\) of \(\Gbf(\A_f^{(p)})\) such that \(K^p \times \Gbf(\Zp)\) is neat.
  We have an object \(((j_* \Fcal^K V)_K, (\corr \ul{j}_* u(K_2, g, K_1, K')_{\Fp})_{K_2,g,K_1,K'})\) of \(\PreH(\Acal^*_{n,?,\Fp}, C, \Qell)\).
  By decomposing \(V\) we may assume that \(V\) is irreducible, and so we may and do assume that \(V\) is pure of some weight \(m \in \Z\).
  Let \(a = -m + n(n+1)/2\), so that the perverse sheaves \(\IC^K(V)[n(n+1)/2]\) are pure of weight \(a\).
  %we could also take \(a = -m + n(n+1)/2 - 1\)
  Denoting \(\ul{a} = (a, \dots, a) \in \Z^{n+1}\) we have by \cite[Théorème 3.1.4, Lemme 3.3.3, Lemme 5.1.3]{MorelSiegel1} an isomorphism in \(\PreH(\Acal^*_{n,?,\Fp}, C, \Qell)\)
  \begin{align*}
    &\ ((\IC^K(V))_K, (u(K_2, g, K_1, K')^{\IC}_{\Fp})_{K_2,g,K_1,K'}) \\
    \simeq &\ ((w_{\leq \ul{a}} j_* \Fcal^K V)_K, (w_{\leq \ul{a}} \corr \ul{j}_* u(K_2, g, K_1, K')_{\Fp})_{K_2,g,K_1,K'}).
  \end{align*}
  Recall from Section \ref{sec:gen_Shim_boundary_boundary} (originating from \cite[Proposition 1.1.3]{MorelBook}, see also the end of the proof of \cite[Proposition 4.2.3]{MorelSiegel1}) that for a standard parabolic subgroup \(\Pbf = \Pbf_{n_1} \cap \dots \cap \Pbf_{n_r}\) of \(\Gbf\) and \(K=K^p \times \Gbf(\Zp)\) neat, sequences of boundary strata \((S_1,\dots,S_r)\), where \(S_1\) is a boundary stratum of \(\Sh(\Gbf,\Xcal,K)^*_{\Fp}\) and for \(1 \leq j < r\) \(S_{j+1}\) is a boundary stratum of \(S_j^*\) corresponding to the image of \(\Pbf_{n_1} \cap \dots \cap \Pbf_{n_{j+1}}\) in \(\Mbf_{\Pbf_{n_1} \cap \dots \cap \Pbf_{n_j}}\), are in bijection with \(\Pbf(\A_f^{(p)}) \backslash \Gbf(\A_f^{(p)}) / K\), and for \(g \in \Gbf(\A_f^{(p)}) \times \Gbf(\Zp)\) we have a map
  \[ T_{\Pbf,gK}: \Sh(\Mbf_\Pbf, \Xcal_\Pbf, K(\Pbf, gK))_{\Fp} \longrightarrow (\Acal_{n,K}^*)_{\Fp}. \]
  As in \cite{MorelSiegel1} we denote by \(i_m\) the locally closed immersion of
  \[ \underset{gK^p \in [\Pbf_m(\A_f^{(p)}) \curvearrowright \Gbf(\A_f^{(p)}) / K^p]}{\colim} \Sh(\Mbf_{\Pbf_m}, \Xcal_{\Pbf_m}, K(\Pbf_m, gK))_{\Fp} \]
  in \((\Acal_{n,K}^*)_{\Fp}\) obtained by collecting the immersions \(i_{\Pbf_m, gK}\).
  We need a slightly more complicated Grothendieck group than in \cite[\S 5.1]{MorelSiegel1} and we work with the Grothendieck group of \(\PreH(\Acal^*_{n,?,\Fp}, C, \Qell)\) instead, i.e.\ the group generated by isomorphism classes of objects \((L_K)_K, (u(K_2, g, K_1, K'))_{K_2, g, K_1, K'}\) of \(\PreH(\Acal^*_{n,?,\Fp}, C, \Qell)\), with relations
  \begin{multline*}
    [(L_K)_K, (u(K_2, g, K_1, K'))_{K_2, g, K_1, K'}] + [((L''_K)_K, (u''(K_2, g, K_1, K'))_{K_2, g, K_1, K'})] \\
    = [(L'_K)_K, (u'(K_2, g, K_1, K'))_{K_2, g, K_1, K'}]
  \end{multline*}
  whenever there exists a family \((L_K \to L'_K \to L''_K \xrightarrow{+1})_K\) of exact triangles such that every diagram
  \[
    \begin{tikzcd}
      \ol{T_g}^* L_{K_1} \ar[r] \ar[d, "{u(K_2,g,K_1,K')}"] & \ol{T_g}^* L'_{K_1} \ar[r] \ar[d, "{u'(K_2,g,K_1,K')}"] & \ol{T_g}^* L''_{K_1} \ar[r] \ar[d, "{u''(K_2,g,K_1,K')}"] & \ol{T_g}^* L_{K_1}[1] \ar[d, "{u(K_2,g,K_1,K')[1]}"] \\
      \ol{T_1}^* L_{K_2} \ar[r] & \ol{T_1}^* L'_{K_2} \ar[r] & \ol{T_1}^* L''_{K_2} \ar[r] & \ol{T_1}^* L_{K_2}[1]
    \end{tikzcd}
  \]
  commutes.
  The proof of \cite[Proposition 5.1.5]{MorelSiegel1} works just as well in this context thanks to Corollary \ref{coro:weight_trunc_push_PreH}, and so in \(K_0(\PreH(\Acal^*_{n,?,\Fp}, C, \Qell))\) the class of
  \[ ((\IC^K(V))_K, (u(K_2,g,K_1,K')^{\IC}_{\Fp})) \]
  is equal to the sum, over standard parabolic subgroups \(\Pbf = \Pbf_{n_1} \cap \dots \cap \Pbf_{n_r}\) of \(\Gbf\), of \((-1)^r\) times
  \[ ((w_{>\ul{a}^{(n_r)}} \dots w_{>\ul{a}^{(n_1)}} j_* \Fcal^K(V))_K, (w_{>\ul{a}^{(n_r)}} \dots w_{>\ul{a}^{(n_1)}} \corr \ul{j}_* u(K_2,g,K_1,K')_{\Fp})) \]
  where
  \[ a^{(n_i)}_j =
    \begin{cases}
      a & \text{ if } j = n-n_i \\
      +\infty & \text{ otherwise.}
    \end{cases}
  \]
  Recall from Proposition 3.3.4 (ii) loc.\ cit.\ that each truncation functor \(w_{>\ul{a}^{(n_i)}}\) is (canonically) isomorphic to \(i_{n_i *} w_{>a} i_{n_i}^*\), and from Lemme 5.1.4 loc.\ cit.\ that (via this identification) we have
  \begin{align*}
    \ & w_{>\ul{a}^{(n_r)}} \dots w_{>\ul{a}^{(n_1)}} \corr \ul{j}_* u(K_2,g,K_1,K') \\
    =\ & \corr \ul{i_{n_r}}_* w_{>a} \corr \ul{i_{n_r}}^* \dots \ul{i_{n_1}}_* w_{>a} \ul{i_{n_1}}^* \ul{j}_* u(K_2,g,K_1,K')
  \end{align*}
  (pushforward and pullback of correspondences as defined in Section \ref{sec:more_corr_push_pull}).
  Now by Proposition 4.2.3 loc.\ cit.
  \[ i_{n_r *} w_{>a} i_{n_r}^* \dots i_{n_1 *} w_{>a} i_{n_1}^* j_* \Fcal^K V \]
  is identified with\footnote{As usual this direct sum is really a colimit over the groupoid \([\Pbf(\A_f^{(p)}) \curvearrowright \Gbf(\A_f^{(p)})/K^p]\).}
  \[ \bigoplus_{[h] \in \Pbf(\A_f^{(p)}) \backslash \Gbf(\A_f^{(p)}) / K} (T_{\Pbf,gK})_* \AFcal^{K(\Pbf,hK)} R\Gamma(\Lie \Nbf_\Pbf, V)_{<\ul{t}} \]
  where
  \[ t_{n_i} = \frac{(n-n_i)(n-n_i+1)}{2} - a - m \]
  i.e.\ \(\ul{t}\) as in the theorem (here we are using our convention for truncation in Definition \ref{def:trunc_alg_rep}, see Remark \ref{rem:compare_trunc_alg_rep}).
  Translating \cite[Proposition 1.5.3]{MorelBook}, we obtain that via these identifications each correspondence
  \[ \corr \ul{i_{n_r}}_* w_{>a} \corr \ul{i_{n_r}}^* \dots \corr \ul{i_{n_1}}_* w_{>a} \corr \ul{i_{n_1}}^* \corr \ul{j}_* u(K_2,g,K_1,K') \]
  is equal to the matrix \((u^\Pbf_{[h_1], [h_2]})_{[h_1], [h_2]}\) where \(u^\Pbf_{[h_1], [h_2]}\) is the sum over \([h'] \in \Pbf(\A_f^{(p)}) \backslash \Gbf(\A_f^{(p)}) / K'\) satisfying \(h'g \in \Pbf(\A_f^{(p)}) h_1 K_1\) (say, \(h' g \in p_1 h_1 K_1\) where \(p_1 \in \Pbf(\A_f^{(p)})\)) and \(h' \in \Pbf(\A_f^{(p)}) h_2 K_2\) (say, \(h' \in p_2 h_2 K_2\) where \(p_2 \in \Pbf(\A_f^{(p)})\)), of
  \[ |h_2K_2h_2^{-1} \cap \Nbf_\Pbf(\A_f^{(p)}) / h'K'(h')^{-1} \cap \Nbf_\Pbf(\A_f^{(p)})| \]
  times the pushforward along
  \[
    \begin{tikzcd}
      \Sh(\Mbf_\Pbf, \Xcal_\Pbf, K(\Pbf, h_1 K_1)) \ar[d, "{T_{\Pbf, h_1 K_1}}"] & \ar[l, "{T_{\ol{p_1}}}"] \Sh(\Mbf_\Pbf, \Xcal_\Pbf, K(\Pbf, h' K')) \ar[d, "{T_{\Pbf, h' K'}}"] \ar[r, "{T_{\ol{p_2}}}"] & \Sh(\Mbf_\Pbf, \Xcal_\Pbf, K(\Pbf, h_2 K_2)) \ar[d, "{T_{\Pbf, h_2 k_2}}"] \\
      \Acal_{n,K_1}^* & \ar[l, "{\ol{T_g}}"] \Acal_{n,K'}^* \ar[r, "{\ol{T_1}}"] & \Acal_{n,K_2}
    \end{tikzcd}
  \]
  of the correspondence \(au(K(\Pbf, h_2 K_2), \ol{p_2}, \ol{p_1}, K(\Pbf, h_1 K_1), K(\Pbf, h' K'))\) (Definition \ref{def:corr_au}).
  Taking cohomology (see Proposition \ref{pro:compat_gen_push_corr_coh}), we recognize (see Definition \ref{def:hecke_ind})
  \[ \ind_{\Pbf(\A_f^{(p)})}^{\Gbf(\A_f^{(p)})} e(\Sh(\Mbf_\Pbf, \Xcal_\Pbf, ?)_{\Fpbar}, \AFcal^? R\Gamma(\Lie \Nbf_\Pbf, V)_{<\ul{t}}). \]
\end{proof}

In order to formulate the analogous result with \(\Fp\) replaced by \(\Q\) we need a weaker notion of Grothendieck group of admissible representations.

\begin{defi} \label{def:Groth_adm}
  For a characteristic zero field \(F\) and a connected reductive group \(\Gbf\) over \(\Q\), denote by \(K_0^{\Tr}(\Rep_F^{\adm}(\Gbf(\A_f)))\) the quotient of the Grothendieck group of admissible representations of \(\Gbf(\A_f)\) over \(F\) by the following equivalence relation: two virtual admissible representations \(V\) and \(W\) are equivalent if for every compact open subgroup \(K\) of \(\Gbf(\A_f)\), \(V^K = W^K\) in the Grothendieck group of finite-dimensional representations of \(\Hcal(\Gbf(\A_f) // K)\).

  Let \(K_0^{\Tr}(\Rep_{\Qell}^{\adm,\cont}(\Gbf(\A_f) \times \GalQ))\) be the quotient of the Grothendieck group of admissible representations of \(\Gbf(\A_f)\) with commuting continuous action of \(\GalQ\) (meaning that for any compact open \(K \subset \Gbf(\A_f)\), the finite-dimensional representation of \(\GalQ\) on \(K\)-invariant vectors is continuous) by the analogous equivalence relation with \(\Hcal(\Gbf(\A_f) // K)\) replaced by \(\Hcal(\Gbf(\A_f) // K) \times \GalQ\).
\end{defi}

By the Brauer-Nesbitt theorem, the equivalence relations in the definition amount to equality of traces (and in the second case one can restrict to a dense invariant subset of \(\GalQ\)).
The equivalence relation occurring in Definition \ref{def:Groth_adm} is equivalent to the one in \cite[\S I.2]{HarrisTaylor}, although we will not need this fact.
In the setting of Definition \ref{def:Groth_adm}, if \(\Pbf\) is a parabolic subgroup of \(\Gbf\) with unipotent radical \(\Nbf_\Pbf\) and reductive quotient \(\Mbf_\Pbf = \Pbf/\Nbf_\Pbf\) then (non-normalized) parabolic induction defines a morphism
\[ \ind_{\Pbf(\A_f)}^{\Gbf(\A_f)}: K_0^{\Tr}(\Rep_{\Qell}^{\adm,\cont}(\Mbf_\Pbf(\A_f) \times \GalQ)) \to K_0^{\Tr}(\Rep_{\Qell}^{\adm,\cont}(\Gbf(\A_f) \times \GalQ)), \]
essentially because parabolic induction is an exact functor preserving admissibility.

\begin{coro} \label{cor:IH_vs_H_Q}
  Let \(\ell\) be a prime number and \(V\) an algebraic representation of \(\Gbf_{\Qell}\).
  In \(K_0^{\Tr}(\Rep_{\Qell}^{\adm,\cont}(\Gbf(\A_f) \times \GalQ))\) we have
  \[ e(\Acal_{n,?,\Qbar}^*, \IC_\ell(V)) = \sum_\Pbf (-1)^{r_\Pbf} \ind_{\Pbf(\A_f^{(p)})}^{\Gbf(\A_f^{(p)})} e(\Sh(\Mbf_\Pbf, \Xcal_\Pbf, ?)_{\Qbar}, \AFcal^? R\Gamma(\Lie \Nbf_\Pbf, V)_{<\ul{t}}). \]
  where \(\ul{t}\) is as in Theorem \ref{thm:IH_vs_H_Fp}.
\end{coro}
\begin{proof}
  It is enough to show that for \(f \in \Hcal(\Gbf(\A_f))\) and \(\sigma \in \GalQ\) the traces of \(f \times \sigma\) on both sides are equal, but for a given Hecke operator \(f\) we can apply Theorem \ref{thm:IH_vs_H_Fp} for almost all \(p\), so we conclude using Propositions \ref{pro:spe_gal_hecke_coho} and \ref{pro:spe_gal_hecke_IC} and the \v{C}ebotarev density theorem.
\end{proof}

\subsubsection{Intersection cohomology from compactly supported cohomology}

To state the dual version for compactly supported cohomology we first need to recall Kostant's theorem \cite[Theorem 5.14]{Kostant_Liealgcoh} on Lie algebra cohomology.
Recall that we have fixed a maximal split torus \(\Tbf = \Tbf_{\GSpbf_{2n}}\) in \(\GSpbf_{2n} = \Gbf\), and that we consider the order on the root system corresponding to the upper triangular Borel subgroup of \(\Gbf\).
As usual \(\rho\) denotes half the sum of the positive roots of \(\Tbf\) in \(\Gbf\).
Assume that \(V\) is an irreducible algebraic representation of \(\Gbf\) over \(\Qell\), corresponding to the highest weight \(\lambda \in X^*(\Tbf)\).
Then we have
\[ H^i(\Lie \Nbf_\Pbf, V) \simeq \bigoplus_{\substack{w \in W^\Pbf \\ l(w) = i}} V^{\Mbf_\Pbf}_{w \cdot \lambda} \]
where
\begin{itemize}
\item \(W^\Pbf \subset W(\Gbf)\) is the set of Kostant representatives for \(W(\Mbf_\Pbf) \backslash W(\Gbf)\), i.e.\ the set of \(w \in W(\Gbf)\) satisfying \(w^{-1} \alpha > 0\) for all simple roots occurring in \(\Mbf_\Pbf\),
\item \(l: W(\Gbf) \to \Z_{\geq 0}\) is the length function,
\item \(w \cdot \lambda := w(\lambda + \rho) - \rho\), and
\item \(V^{\Mbf_\Pbf}_{\lambda'}\) is the irreducible representation of \(\Mbf_\Pbf\) admitting highest weight \(\lambda'\).
\end{itemize}
Moreover we have a canonical decomposition
\[ R\Gamma(\Lie \Nbf_\Pbf, V) \simeq \bigoplus_{i=0}^{\dim \Nbf_\Pbf} H^i(\Lie \Nbf_\Pbf, V) [-i] \]
because each \(V^{\Mbf_\Pbf}_{w \cdot \lambda}\) occurs with multiplicity one in the Chevalley-Eilenberg complex.
We will denote \(V^{\Mbf_\Pbf}_{\lambda'} = V^{\Mbf_{\Pbf,\lin}}_{\lambda'_\lin} \otimes V^{\Mbf_{\Pbf,\her}}_{\lambda'_\her}\) according to the decomposition \eqref{eq:decomp_lin_her_Siegel}.
For \(\ul{t} = (t_0, \dots, t_n) \in \Z^{n+1}\) denote by \(W^\Pbf_{<\ul{t}}(\lambda)\) the set of \(w \in W^\Pbf\) satisfying \((V^{\Mbf_\Pbf}_{w \cdot \lambda})_{<\ul{t}} = V^{\Mbf_\Pbf}_{w \cdot \lambda}\) (see Definition \ref{def:trunc_alg_rep}), and similarly for \(>\ul{t}\).

\begin{coro} \label{coro:IH_vs_Hc}
  Let \(V\) be an irreducible algebraic representation of \(\Gbf_{\Qell}\), with highest weight \(\lambda\).
  For \(p \neq \ell\) we have an equality in \(K_0(\Rep_{\Qell}^{\adm,\cont}(\Gbf(\A_f^{(p)}) \times \Gal(\Fpbar/\Fp)))\)
  \begin{multline*}
    e(\Acal_{n,?,\Fpbar}^*, \IC^?(V)) \\
    = \sum_\Pbf \sum_{w \in W^\Pbf_{>\ul{t}}(\lambda)} (-1)^{l(w)} \ind_{\Pbf(\A_f^{(p)})}^{\Gbf(\A_f^{(p)})} \left( e_c \left( \Mbf_{\Pbf,\lin}, V^{\Mbf_{\Pbf,\lin}}_{(w \cdot \lambda)_\lin} \right)^{\Mbf_{\Pbf,\lin}(\Zp)} \otimes e_c \left( \Acal_{n-n_r,?,\Fpbar}, \Fcal^? V^{\Mbf_{\Pbf,\her}}_{(w \cdot \lambda)_\her} \right) \right)
  \end{multline*}
  for the same \(\ul{t} \in \Z^{n+1}\) as in Theorem \ref{thm:IH_vs_H_Fp}.
\end{coro}
\begin{proof}
  The proof is very similar to the proof of \cite[Corollary 1.4.6]{MorelBook}.
  We apply Theorem \ref{thm:IH_vs_H_Fp} to \(V^*\) before taking contragredients.
  The dual of the intersection complex \(\IC^K(V^*)\) is identified with \(\IC^K(V)(d)[2d]\) where \(d = \dim \Acal_n = n(n+1)/2\), and via this identification the dual of the correspondences \(u(K_2,g,K_1,K')^{\IC}_{\Fp}\) are equal to \(u(K_1,g^{-1},K_2,g^{-1}K'g)^{\IC}_{\Fp} (d)[2d]\).
  So the contragredient of the left-hand side in Theorem \ref{thm:IH_vs_H_Fp} is
  \[ e(\Acal_{n,?,\Fpbar}^*, \IC(V))(d). \]

  In order to rewrite the contragredient of the right-hand side we first use Corollary \ref{coro:dual_para_ind} to express, for algebraic representations \(W_\lin\) of \(\Mbf_{\Pbf,\lin}\) and \(W_\her\) of \(\Mbf_{\Pbf,\her}\) (over \(\Qell\)), the contragredient of
  \[ \ind_{\Pbf(\A_f^{(p)})}^{\Gbf(\A_f^{(p)})} \left( e(\Mbf_{\Pbf,\lin}, W_\lin)^{\Mbf_{\Pbf,\lin}(\Zp)} \otimes e(\Acal_{n-n_r,?,\Fpbar}, \Fcal^?(W_\her)) \right) \]
  as
  \[ \ind_{\Pbf(\A_f^{(p)})}^{\Gbf(\A_f^{(p)})} \left( e(\Mbf_{\Pbf,\lin}, W_\lin)^{*,\Mbf_{\Pbf,\lin}(\Zp)} \otimes e(\Acal_{n-n_r,?,\Fpbar}, \Fcal^?(W_\her))^* \otimes \delta_{\Pbf(\A_f^{(p)})} \right). \]
  We have \(\delta_{\Pbf(\A_f^{(p)})} = |\delta_\Pbf|_f\) where \(\delta_\Pbf\) is the (algebraic) character by which \(\Mbf_\Pbf\) acts on \(\bigwedge^{\dim \Nbf_\Pbf} \Lie \Nbf_\Pbf\).
  The character \(\delta_\Pbf\) is easily computed (in the decomposition \eqref{eq:decomp_lin_her_Siegel}): on each \(\GLbf_{n_i-n_{i-1}}\) it is equal to \(\det^{2n+1-n_i-n_{i-1}}\) and on \(\GSpbf_{2(n-n_r)}\) it is equal to \(\nu^{d-d_{n_r}}\) where \(d_{n_r} = \dim \Acal_{n-n_r} = (n-n_r)(n-n_r+1)/2\).
  By Example \ref{exam:orientation_GL} and the observation \(2n+1-n_i-n_{i-1} \equiv n_i-n_{i-1}-1 \mod 2\) we obtain
  \[ e(\Mbf_{\Pbf,\lin}, W_\lin)^{*,\Mbf_{\Pbf,\lin}(\Zp)} \otimes \delta_{\Pbf(\A_f)}|_{\Mbf_{\Pbf,\lin}(\A_f)} = (-1)^{2q_\lin(\Pbf)} e_c(\Mbf_{\Pbf,\lin}, W_\lin^* \otimes \delta_{\Pbf,\lin}^{-1})^{\Mbf_{\Pbf,\lin}(\Zp)} \]
  where \(q_\lin(\Pbf)\) is short for \(q(\Mbf_{\Pbf,\lin,\der}(\R))\) and \(\delta_{\Pbf,\lin}\) denotes the restriction of \(\delta_\Pbf\) to \(\Mbf_{\Pbf,\lin}\).
  Using Remark \ref{rema:twisting_sim} we also compute
  \begin{align*}
    e(\Acal_{n-n_r,?,\Fpbar}, \Fcal^? W_\her)^* \otimes \delta_{\Pbf(\A_f^{(p)})}|_{\Mbf_{\Pbf,\her}(\A_f^{(p)})}
    & = e_c(\Acal_{n-n_r,?,\Fpbar}, \Fcal^? W_\her^*)(d_{n_r}) \otimes \delta_{\Pbf(\A_f^{(p)})}|_{\Mbf_{\Pbf,\her}(\A_f^{(p)})} \\
    & = e_c(\Acal_{n-n_r,?,\Fpbar}, \Fcal^? (W_\her^* \otimes \delta_{\Pbf,\her}^{-1}))(d)
  \end{align*}
  where \(\delta_{\Pbf,\her}\) denotes the restriction of \(\delta_\Pbf\) to \(\Mbf_{\Pbf,\her}\).
  We conclude
  \begin{align*}
    & \ind_{\Pbf(\A_f^{(p)})}^{\Gbf(\A_f^{(p)})} \left( e(\Mbf_{\Pbf,\lin}, W_\lin)^{\Mbf_{\Pbf,\lin}(\Zp)} \otimes e(\Acal_{n-n_r,?,\Fpbar}, \Fcal^?(W_\her)) \right)^* \\
    = & (-1)^{2q_\lin(\Pbf)} \ind_{\Pbf(\A_f^{(p)})}^{\Gbf(\A_f^{(p)})} \left( e_c(\Mbf_{\Pbf,\lin}, W_\lin^* \otimes \delta_{\Pbf,\lin}^{-1})^{\Mbf_{\Pbf,\lin}(\Zp)} \otimes e(\Acal_{n-n_r,?,\Fpbar}, \Fcal^?(W_\her^* \otimes \delta_{\Pbf,\her}^{-1})) \right).
  \end{align*}

  By duality for Lie algebra cohomology \cite{Hazewinkel_duality_Lie} we have
  \[ H^i(\Lie \Nbf_\Pbf, V^*)^* \simeq H^{\dim \Nbf_\Pbf - i}(\Lie \Nbf_\Pbf, V) \otimes \bigwedge^{\dim \Nbf_\Pbf} \Lie \Nbf_\Pbf. \]
  We deduce from the above computation of \(\delta_\Pbf\) that
  \[ \diag(\lambda I_{n_i}, I_{2(n-n_i)}, \lambda^{-1} I_{n_i}) \in \Zbf(\Mbf_{\Pbf,\lin}) \]
  acts by \(\lambda^{2(d-d_{n_i})}\) on \(\bigwedge^{\dim \Nbf_\Pbf} \Lie \Nbf_\Pbf\).
  % perfect pairing:
  %\[ H^i(\Lie \Nbf_\Pbf, V^*) \times H^{\dim \Nbf_\Pbf - i}(\Lie \Nbf_\Pbf, V) \longrightarrow \bigwedge^{\dim \Nbf_\Pbf} (\Lie \Nbf_\Pbf)^*. \]
  Denoting \(s_{n_i} = -t_{n_i} - 2(d-d_{n_i})\) we obtain
  \[ \left( H^i(\Lie \Nbf_\Pbf, V^*)_{<\ul{t}} \right)^* \simeq H^{\dim \Nbf_\Pbf - i}(\Lie \Nbf_\Pbf, V)_{>\ul{s}} \otimes \bigwedge^{\dim \Nbf_\Pbf} \Lie \Nbf_\Pbf. \]
  Note that for \(t_{n_i} = d_{n_i} - d\) we have \(s_{n_i} = t_{n_i}\).

  Write a decomposition into irreducible pieces
  \[ H^i(\Lie \Nbf_\Pbf, V^*)_{<\ul{t}} \simeq \bigoplus_j W_{i,j,\lin} \otimes W_{i,j,\her} \]
  so that we have
  \[ H^i(\Lie \Nbf_\Pbf, V)_{>\ul{s}} \simeq \bigoplus_j W'_{i,j,\lin} \otimes W'_{i,j,\her} \]
  with \(W'_{i,j,?} \simeq  W_{\dim \Nbf_\Pbf - i, j,?}^* \otimes \delta_{\Pbf,?}^{-1}\) for \(? \in \{\lin,\her\}\).
  Putting together the above computations we get
  \begin{align*}
    & \ind_{\Pbf(\A_f^{(p)})}^{\Gbf(\A_f^{(p)})} e(\Sh(\Mbf_\Pbf, \Xcal_\Pbf, ?)_{\Fpbar}, \AFcal^? R\Gamma(\Lie \Nbf_\Pbf, V^*)_{<\ul{t}})^* \\
    =& \sum_{i,j} (-1)^i \ind_{\Pbf(\A_f^{(p)})}^{\Gbf(\A_f^{(p)})} \left( e(\Mbf_{\Pbf,\lin}, W_{i,j,\lin})^{*,\Mbf_{\Pbf,\lin}(\Zp)} \otimes e(\Acal_{n-n_r,?,\Fpbar}, \Fcal^? W_{i,j,\her})^* \otimes \delta_{\Pbf(\A_f^{(p)})} \right) \\
    =& \sum_{i,j} (-1)^{i+2q_\lin(\Pbf)} \ind_{\Pbf(\A_f^{(p)})}^{\Gbf(\A_f^{(p)})} \left( e_c(\Mbf_{\Pbf,\lin}, W_{i,j,\lin}^* \otimes \delta_{\Pbf,\lin}^{-1})^{\Mbf_{\Pbf,\lin}(\Zp)} \otimes e_c(\Acal_{n-n_r,?,\Fpbar}, \Fcal^? (W_{i,j,\her}^* \otimes \delta_{\Pbf,\her}^{-1})) \right)(d) \\
    =& \sum_{i,j} (-1)^{i+2q_\lin(\Pbf)+\dim \Nbf_\Pbf} \ind_{\Pbf(\A_f^{(p)})}^{\Gbf(\A_f^{(p)})} \left( e_c(\Mbf_{\Pbf,\lin}, W'_{i,j,\lin}) \otimes e_c(\Acal_{n-n_r,?,\Fpbar}, \Fcal^? (W'_{i,j,\her})) \right)(d)
  \end{align*}
  We conclude by simplifying signs using the equality
  \[ 2q_\lin(\Pbf) + r_\Pbf + 2 d_{n_r} + \dim \Nbf_\Pbf = 2d. \]
  This equality can be checked directly, or deduced from the Iwasawa decomposition which allows us to see \(\Xcal\) as a homogeneous space under \(\Pbf(\R)\) and comparing dimensions.
\end{proof}

\begin{coro} \label{coro:IH_vs_Hc_Q}
  Let \(V = V_\lambda\) be an irreducible algebraic representation of \(\Gbf_{\Qell}\) characterized by its highest weight \(\lambda\).
  We have an equality in \(K_0^{\Tr}(\Rep_{\Qell}^{\adm,\cont}(\Gbf(\A_f) \times \GalQ))\) 
  \begin{multline*}
    e(\Acal_{n,?,\Qbar}^*, \IC_\ell(V)) \\
    = \sum_\Pbf \sum_{w \in W^\Pbf_{>\ul{t}}} (-1)^{l(w)} \ind_{\Pbf(\A_f)}^{\Gbf(\A_f)} \left( e_c \left( \Mbf_{\Pbf,\lin}, V^{\Mbf_{\Pbf,\lin}}_{(w \cdot \lambda)_\lin} \right) \otimes e_c \left( \Acal_{n-n_r,?,\Qbar}, \Fcal^? V^{\Mbf_{\Pbf,\her}}_{(w \cdot \lambda)_\her} \right) \right).
  \end{multline*}
\end{coro}
\begin{proof}
  This follows from taking traces in Corollary \ref{coro:IH_vs_Hc} and using Proposition \ref{pro:spe_gal_hecke_IC}.
\end{proof}

\subsubsection{Crystallineness of intersection cohomology}

Recall the following corollary to a result of Faltings-Chai.
\begin{theo} \label{theo:FC_crys}
  For any \(n \geq 1\), \(M \geq 3\), \(V\) an irreducible algebraic representation of \(\GSpbf_{2n}\) over \(\Qell\) and \(i \geq 0\) the Galois representation \(H^i_c((\Acal_{n,M})_{\Qbar}, \Fcal^{K(M)}_\ell(V))\) is crystalline if \(\ell\) does not divide \(M\).
\end{theo}
\begin{proof}
  The dual statement (for ordinary cohomology) follows from \cite[Theorem 6.2 (ii)]{FaltingsChai} (using also (i) in this theorem and Theorem 5.5 to compare dimensions).
\end{proof}

\begin{coro} \label{cor:sigma_IH_crys}
  Fix an isomorphism \(\iota\) between the algebraic closures of \(\Q\) in \(\C\) and \(\Qellbar\).
  For any \(n \geq 1\), \(\tau \in \ICcal(\Spbf_{2n})\), \(\psi \in \Psit^{\unr, \tau}_{\disc}(\Spbf_{2n})\) the semisimple Galois representation \(\sigma_{\psi, \iota}^{\IH}\) introduced in Theorem \ref{theo:IH_explicit_crude} is crystalline at \(\ell\).
\end{coro}
\begin{proof}
  Let \(V\) be the algebraic representation of \(\GSpbf_{2n}\) corresponding to \(\tau\).
  By Corollary \ref{coro:IH_vs_Hc_Q} and Theorem \ref{theo:FC_crys} (choosing an auxiliary level \(M \geq 3\) not divisible by \(\ell\)) we know that in the Grothendieck group of continuous \(\ell\)-adic representations of \(\GalQ\), \(e(\Acal_n^*, \IC_{\ell}(V))\) is represented by an alternate sum of crystalline representations.
  As was already used in the proof of Theorem \ref{theo:IH_explicit_crude}, by purity any irreducible representation of \(\GalQ\) occurs as a subquotient in at most one degree in \(H^\bullet(\Acal_n^*, \IC(V))\).
\end{proof}

\section{Odd spin Galois representations}
\label{sec:odd_spin}

\subsection{Existence and uniqueness of a lifting in conductor one}

\begin{prop} \label{pro:ex_lift_cond_one}
  Let \(f : H' \rightarrow H\) be a surjective morphism between reductive groups
  over \(\Qellbar\) such that \(\ker f\) is a central torus \(C\) in \(H'\).
  Let \(\rho : \GalQ \rightarrow H(\Qellbar)\) be a continuous Galois
  representation which is unramified away from \(\ell\) and such that
  \(\rho|_{\GalQell}\) is crystalline.
  Let \(\tau : \GL_{1,\Qellbar} \rightarrow H\) be the Hodge-Tate 1-parameter subgroup for \(\rho\) (well-defined up to \(H(\Qellbar)\)-conjugation).
  Assume that there exists \(\tau' : \GL_{1,\Qellbar} \rightarrow H'\) lifting \(\tau\).
  Then there exists a unique continuous Galois representation \(\rho' : \GalQ
  \rightarrow H'(\Qellbar)\) unramified away from \(\ell\), crystalline at \(\ell\)
  with Hodge-Tate cocharacter \(\tau'\), and such that \(\rho = f \circ \rho'\).
\end{prop}
\begin{proof}
  By \cite[Proposition 2.8.2]{PatrikisVarTate} there exists a geometric lift
  \(\rho'_0 : \GalQ \rightarrow H'(\Qellbar)\) of \(\rho\) with Hodge-Tate
  1-parameter subgroup \(\tau'\) (up to \(H'(\Qellbar)\)-conjugation).
  In fact this proposition is stated for isogenies \(\mathrm{GSpin}_{2n+1}
  \rightarrow \mathrm{SO}_{2n+1}\), but the proof obviously applies to the
  general case.

  Now we need to twist \(\rho'_0\) by a finite order character \(\GalQ \rightarrow
  C(\Qellbar)\) to obtain \(\rho'\) of conductor one.
  First we note that for any prime \(p \neq \ell\), the representation
  \(\rho'_0|_{\GalQp}\) is unramified up to a twist by a character \(\chi_p :
  \GalQp \rightarrow C(\Qellbar)\), and we may and do assume that this character has finite order.
  Similarly, \(\rho'_0|_{\GalQell}\) is crystalline up to a twist by finite order
  character \(\chi_{\ell} : \GalQell \rightarrow C(\Qellbar)\).
  This case (\(p = \ell\)) is not quite as straightforward: it follows from the
  existence of a crystalline lift of \(\rho|_{\GalQell}\) having Hodge-Tate
  1-parameter subgroup \(\tau'\) (\cite[Proposition 6.5]{ConradLifting}, which
  relies on \cite[\S 2]{WintenbergerTann}).

  Denote by \(I_p\) the inertia subgroup of \(\GalQp\).
  Since \(\rho'_0\) is unramified at almost all primes, for almost all primes \(p\)
  we have that \(\chi_p\) is unramified, i.e.\ \(\chi_p|_{I_p} = 1\).
  We claim that there exists a unique continuous finite order character \(\chi :
  \GalQ \rightarrow C(\Qellbar)\) such that for all primes \(p\), \(\chi|_{I_p} =
  \chi_p|_{I_p}\).
  This claim follows from the local Kronecker-Weber theorem, which identifies
  the image of \(I_p\) in \(\GalQp^{\mathrm{ab}}\) with \(\Zp^{\times}\): for any \(n
  \geq 1\) and any finite set of primes \(S\), letting \(N = \prod_{p \in S} p^n\),
  the Galois group of the abelian extension \(\mathbb{Q}(\zeta_N) / \mathbb{Q}\)
  decomposes as a product of inertia subgroups \(\prod_{p \in S} (\Z /
  p^n\Z)^{\times}\).

  Uniqueness follows from Minkowski's theorem (or the global Kronecker-Weber
  theorem).
\end{proof}

\begin{prop} \label{pro:un_cond_one}
  Let \(\rho_1 : \GalQ \rightarrow \GL(V_1)\) and \(\rho_2 : \GalQ \rightarrow \GL(V_2)\) be semi-simple continuous representations, where \(V_1, V_2\) are finite-dimensional \(\Qellbar\)-vector spaces having equal dimensions.
  Assume that \(\rho_1\) and \(\rho_2\) have conductor one, i.e.\ are unramified
  away from \(\ell\) and crystalline at \(\ell\).
  Assume also that there exists an integer \(n \geq 1\) such that for all primes
  \(p>n\), there exists an \(n\)-th root of unity \(\zeta\) such that the
  characteristic polynomials of \(\rho_1(\Frob_p)\) and \(\zeta \rho_2(\Frob_p)\)
  are equal.
  Then \(\rho_1\) and \(\rho_2\) are isomorphic.
\end{prop}
\begin{proof}
  Choose an open subgroup \(U \subset \Qellbar^{\times}\) such that \(U
  \smallsetminus \{1\}\) does not contain any \(n\)-th root of unity.
  There exists a finite Galois extension \(K / \Q\) such that for any \(\sigma \in
  \Gal_K\), the eigenvalues of \(\rho_1(\sigma)\) and those of \(\rho_2(\sigma)\)
  belong to \(U\).
  Then for almost all primes \(\mathfrak{p}\) of \(K\), we have that \(\rho_1(\Frob_{\mathfrak{p}})\) and \(\rho_2(\Frob_{\mathfrak{p}})\) have the same characteristic polynomial, and so by the \v{C}ebotarev density theorem the semi-simple representations \(\rho_1|_{\Gal_K}\) and \(\rho_2|_{\Gal_K}\) are isomorphic.
  Choose an isomorphism \(f: V_1 \rightarrow V_2\) intertwining these
  restrictions.
  Then \(f \in \Hom_{\Qellbar}(V_1, V_2)^{\Gal(K / \Q)}\), and this is a
  sub-\(\GalQ\)-representation of \(\Hom_{\Qellbar}(V_1, V_2)\), therefore it is
  also a geometric representation of conductor one.
  Since any finite image crystalline representation of \(\GalQell\) is in fact
  unramified, we obtain using Minkowski's theorem that \(\Gal(K / \Q)\) acts
  trivially on \(\Hom_{\Qellbar}(V_1, V_2)^{\Gal(K / \Q)}\), and we conclude that
  \(f\) intertwines \(\rho_1\) and \(\rho_2\).
\end{proof}

\subsection{Odd spin Galois representations}

In the following theorem we recall what is already known about the existence of Galois representations valued in orthogonal groups in our case.

\begin{theo} \label{thm:existence_rho_SO}
  Fix a prime number \(\ell\) and an isomorphism \(\iota: \C \simeq \Qellbar\).
  Let \(n \geq 1\).
  Let \(\tau \in \ICcal(\Spbf_{2n})\).
  For any \(\psi \in \tilde{\Psi}_{\disc}^{\unr, \tau}(\Spbf_{2n})\), there exists a continuous semisimple representation \(\rho_{\psi, \iota}^{\SO} : \GalQ \rightarrow \SO_{2n+1}(\Qellbar)\) unramified away from \(\ell\) and crystalline at \(\ell\) and such that for any \(p \neq \ell\), \(\rho_{\psi, \iota}^{\SO}(\Frob_p)^{\sesi}\) is conjugated to \(\iota (\dpsitau (c_p(\psi)) )\).
  The Hodge-Tate cocharacter of \(\rho_{\psi,\iota}^{\SO}|_{\GalQell}\), a conjugacy class of morphisms \(\GL_1 \to \SO_{2n+1}\), has differential equal to \(\tau\).

  If \(\rho : \GalQ \rightarrow \SO_{2n+1}(\Qellbar)\) is any continuous semisimple morphism such that for almost all primes \(p\), \(\rho\) is unramified at \(p\) and \(\rho(\Frob_p)^{\sesi}\) is conjugated to \(\iota ( \dpsitau (c_p(\psi)) )\), then \(\rho\) is conjugated to \(\rho_{\psi, \iota}^{\SO}\).
\end{theo}
\begin{proof}
  Write \(\psi = \oplus_i \pi_i[d_i]\) where \(\pi_i\) is a self-dual automorphic representation for \(\GLbf_{n_i}\) of conductor one, satisfying \(2n+1 = \sum_i n_i d_i\).
  The eigenvalues of the infinitesimal character of \(\pi_{i,\infty}\) are distinct and belong to \(\frac{d_i-1}{2} + \Z\), so it follows from \cite[Theorem 4.2]{CheHar} that there exists a continuous semisimple representation \(\rho_{i,\iota}: \GalQ \to \GL_{n_i}(\Qellbar)\) which is unramified away from \(\ell\) and crystalline at \(\ell\) and such that for any prime number \(p \neq \ell\) the characteristic polynomial of \(\rho_{i,\iota}(\Frob_p)\) is equal to that of \(p^{(d_i-1)/2} c_p(\pi_i)\).
  We have \(\rho_{i,\iota}^\vee \simeq \chi_\ell^{d_i-1} \otimes \rho_{i,\iota}\) because \(\pi_i\) is self-dual.
  Define a \(2n+1\)-dimensional linear representation \(\sigma\) of \(\GalQ\) as \(\bigoplus_i \rho_i \otimes (1 \oplus \chi_\ell \oplus \dots \oplus \chi_\ell^{d_i-1})\).
  It is clearly self-dual, and for any index \(i\) such that \(d_i\) is odd we know thanks to \cite[Corollary 1.3]{BC} that the self-dual representation \(\rho_{i,\iota} \otimes \chi_\ell^{(d_i-1)/2}\) is of orthogonal type, i.e.\ it factors through the standard representation of the orthogonal group \(\mathrm{O}_{n_i}(\Qellbar)\).
  It actually factors through \(\SO_{n_i}(\Qellbar)\) because it is unramified away from \(\ell\) and crystalline at \(\ell\), so its determinant is everywhere unramified.
  It follows that \(\sigma\) is also of orthogonal type, i.e.\ it is isomorphic to \(\rho_{\psi,\iota}^{\SO}: \GalQ \to \SO(\Qellbar)\) composed with the standard representation \(\Std_{\SO_{2n+1}}\) of \(\SO_{2n+1}\).
  The representation \(\Std_{\SO_{2n+1}}\) induces an injective map on semi-simple conjugacy classes, so we deduce from local-global compatibility at non-\(\ell\)-adic finite places in \cite[Theorem 4.2]{CheHar} that for all primes \(p \neq \ell\) the conjugacy class of \(\rho_{\psi, \iota}^{\SO}(\Frob_p)^{\sesi}\) is equal to that of \(\iota (\dpsitau (c_p(\psi)) )\).
  Similarly the Hodge-Tate weights of \(\Std_{\SO_{2n+1}} \circ \rho_{\psi,\iota}^{\SO}|_{\GalQell}\) are given in \cite[Theorem 4.2]{CheHar}, and so the Hodge-Tate cocharacter of \(\rho_{\psi,\iota}^{\SO}\) is conjugated under \(\SO_{2n+1}(\Qellbar)\) to \(\tau\).

  If \(\rho': \GalQ \to \SO_{2n+1}(\Qellbar)\) is another continuous semisimple morphism satisfying the same condition at almost all primes \(p \neq \ell\), then \v{C}ebotarev's theorem implies that the traces of \(\rho'\) and \(\rho_{\psi,\iota}^{\SO}\) are equal, and as they are both semisimple it is well-known that this implies that \(\rho'\) is conjugated to \(\rho_{\psi,\iota}^{\SO}\) by some \(g \in \GL_{2n+1}(\Qellbar)\).
  It then follows from \cite[Theorem 2.3]{Griess_conj_G2} (see also \cite[Proposition 2.3]{Larsen_conj1}) that we may take \(g\) in \(\mathrm{O}_{2n+1}(\Qellbar)\).
  We have \(\mathrm{O}_{2n+1} \simeq \mu_2 \times \SO_{2n+1}\) and so we may even take \(g\) in \(\SO_{2n+1}(\Qellbar)\).
\end{proof}

\begin{theo} \label{theo:odd_spin_rep}
  Fix a prime number \(\ell\) and an isomorphism \(\iota: \C \simeq \Qellbar\).
  Let \(n \geq 1\).
  Let \(\tau \in \ICcal(\Spbf_{2n})\).
  For any \(\psi \in \tilde{\Psi}_{\disc, \nonendo}^{\unr, \tau}(\Spbf_{2n})\), there exists a continuous semisimple morphism \(\rho_{\psi, \iota}^{\GSpin} : \GalQ \rightarrow \GSpin_{2n+1}(\Qellbar)\) unramified away from \(\ell\) and crystalline at \(\ell\) and such that for any \(p \neq \ell\), \(\rho_{\psi, \iota}^{\GSpin}(\Frob_p)^{\mathrm{ss}}\) is conjugated to \(\iota ( p^{n(n+1)/4} \dpsitausc (\cpsc(\psi)) )\).

  Moreover the conjugacy class of \(\rho_{\psi,\iota}^{\GSpin}\) admits the following characterizations.
  \begin{enumerate}
  \item If \(\rho : \GalQ \rightarrow \GSpin_{2n+1}(\Qellbar)\) is any continuous semisimple morphism such that for almost all primes \(p\), \(\rho\) is unramified at \(p\) and \(\rho(\Frob_p)^{\sesi}\) is conjugated to \(\iota ( p^{n(n+1)/4} \dpsitausc (\cpsc(\psi)) )\), then \(\rho\) is conjugated to \(\rho_{\psi, \iota}^{\GSpin}\).
  \item If \(\rho: \GalQ \to \GSpin_{2n+1}(\Qellbar)\) is a continuous morphism which lifts \(\rho_{\psi,\iota}^{\SO}\) (up to conjugacy by \(\SO_{2n+1}(\Qellbar)\)), is unramified away from \(\ell\) and crystalline at \(\ell\), and satisfies \(\nu \circ \rho = \chi_\ell^{-n(n+1/2)}\), then \(\rho\) is conjugated to \(\rho_{\psi, \iota}^{\GSpin}\).
  \end{enumerate}
\end{theo}
\begin{proof}
  Write \(\tau = (w_1, \dots, w_n)\) with \(w_1 > \dots > w_n > 0\) integers.
  Since \(\sum_i w_i = n(n+1)/2 \mod 2\), \(\tau' = (w_1, \dots, w_n, (\sum_i w_i)/2, n(n+1)/4)\) defines a conjugacy class of cocharacters \(\GL_1 \to \GSpin_{2n+1}\) above \(\tau\).
  By Proposition \ref{pro:ex_lift_cond_one}, there exists a unique continuous \(\rho_{\psi, \iota}^{\GSpin} : \GalQ \to \GSpin_{2n+1}(\Qellbar)\) lifting \(\rho_{\psi, \iota}^{\SO}\) and which is unramified away from \(\ell\) and crystalline at \(\ell\) with Hodge-Tate cocharacter \(\tau'\).

  We also have the representation \(\sigma_{\psi, \iota}^{\spin}\) from Corollary \ref{coro:seed_ex_sigma}.
  For all \(p \neq \ell\), \(\sigma_{\psi, \iota}^{\spin}(\Frob_p)^{\sesi}\) is conjugated to \(\spin(\rho_{\psi, \iota}^{\GSpin}(\Frob_p))^{\sesi}\) or \(-\spin(\rho_{\psi, \iota}^{\GSpin}(\Frob_p))^{\sesi}\) in \(\GL_{2^n}(\Qellbar)\).
  Applying Proposition \ref{pro:un_cond_one} to \(\sigma_{\psi, \iota}^{\spin}\) and \(\spin \circ \rho_{\psi, \iota}^{\GSpin}\) (which is semisimple by \cite[Corollaire 4.3]{SerreCompRed}), we obtain that these two representations coincide, and thus \(\rho_{\psi, \iota}^{\GSpin}\) satisfies all desired conditions.

  The second characterization of \(\rho_{\psi,\iota}^{\GSpin}\) follows from uniqueness in Proposition \ref{pro:ex_lift_cond_one} and the fact that the surjective morphism \(\GSpin_{2n+1} \to \SO_{2n+1} \times \GL_1\), where the second component is \(\nu\), is a central isogeny (with kernel \(\mu_2\)).

  We are left to prove the first characterization.
  If \(\rho : \GalQ \rightarrow \GSpin_{2n+1}(\Qellbar)\) is continuous semisimple and satisfies the same property as \(\rho_{\psi,\iota}^{\GSpin}\) at almost all primes \(p\), then \(\spin \circ \rho\) is also continuous semisimple and so by \v{C}ebotarev's theorem it is conjugated to \(\spin \circ \rho_{\psi, \iota}^{\GSpin}\).
  In particular \(\spin \circ \rho\) is also unramified away from \(\ell\) and crystalline at \(\ell\).
  Since the representation \(\spin : \GSpin_{2n+1} \rightarrow \GL_{2^n}\) is faithful, this implies that \(\rho\) enjoys the same property.
  By uniqueness in Proposition \ref{pro:ex_lift_cond_one} we have \(\rho = \chi \rho_{\psi, \iota}^{\GSpin}\) for some continuous character \(\chi : \GalQ \rightarrow \Qellbar^{\times}\) which is unramified away from \(\ell\) and crystalline at \(\ell\), and so \(\chi = \chi_\ell^N\) for some integer \(N\).
  By composing \(\rho\) with \(\nu\) and evaluating at \(\Frob_p\) for any \(p \neq \ell\) we see that we have \(N=0\).
\end{proof}

\begin{rema} \label{rem:rho_GSpin_odd_small_rank}
  For \(n \leq 2\) these Galois representations \(\rho_{\psi,\iota}^\GSpin\) are not new.
  Consider the case \(\psi = \pi_0[1]\) for simplicity (see Section \ref{sec:odd_spin_nontemp} for the other cases).
  \begin{enumerate}
  \item If \(n=1\) then \(\pi_0 = \Sym^2 \pi\) for a level one cuspidal automorphic representation \(\pi\) for \(\PGLbf_2\) (see Proposition \ref{pro:cpsc_Sym2}) and via the isomorphism \(\spin: \GSpin_3 \simeq \GL_2\) the morphism \(\rho_{\psi,\iota}^\GSpin\) is (up to conjugation) a Tate twist of the Galois representation \(\rho_{\pi,\iota}: \GalQ \to \GL_2(\Qellbar)\) associated to \(\pi\) by Deligne \cite{Deligne_GalGL2}.
  \item If \(n=2\) then \(\pi_0 = \Lambda^* \pi\) for a self-dual level one cuspidal automorphic representation \(\pi\) for \(\PGLbf_4\) (see Proposition \ref{pro:cpsc_Sp4_SO5}) and via the isomorphism \(\GSpin_5 \simeq \GSp_4\) the morphism \(\rho_{\psi,\iota}^\GSpin\) is (up to conjugation) a Tate twist of the Galois representation \(\rho_{\pi,\iota}^\GSp: \GalQ \to \GSp_4(\Qellbar)\) associated to \(\pi\) by Weissauer \cite{Weissauer_GalGSp4} (now also a special case of Theorem \ref{theo:existence_rho_GSp}).
  \end{enumerate}
\end{rema}

\subsection{Non-tempered case}
\label{sec:odd_spin_nontemp}

In the setting of Theorem \ref{theo:odd_spin_rep}, we consider the non-tempered case, i.e.\ the case where \(\psi = \pi[2d+1]\) with \(d>0\).
We have \(\pi[1] \in \Psi_{\disc,\nonendo}^{\unr,\tau'}(\Spbf_{2m})\) for a unique \(\tau' \in \ICcal(\Spbf_{2m})\) where the integer \(m\) is determined by the relation \(2n+1=(2m+1)(2d+1)\).
Similarly to Definition \ref{def:alpha_m_d} we have a morphism \(\beta: \SO_{2m+1} \times \SL_2 \to \SO_{2n+1}\) such that \(\Std_{\SO_{2n+1}} \circ \beta \simeq \Std_{\SO_{2m+1}} \otimes \Sym^{2d} \Std_{\SL_2}\), and this morphism is unique up to conjugation by \(\SO_{2n+1}(\Qbar)\).
Up to conjugacy we may assume that it restricts to
\begin{align*}
  \Tcal_{\SO_{2m+1}} \times \Tcal_{\SL_2} & \longrightarrow \Tcal_{\SO_{2n+1}} \\
  ((x_1, \dots, x_m), t) & \longmapsto (x_1 t^{2d}, x_1 t^{2d-2}, \dots, x_1 t^{-2d}, \dots, x_m t^{2d}, \dots, x_m t^{-2d}, t^{2d}, t^{2d-2}, \dots, t^2).
\end{align*}
It lifts uniquely to \(\tilde{\beta}: \Spin_{2m+1} \times \SL_2 \to \Spin_{2n+1}\), which restricts to
\begin{align*}
  \Tcal_{\Spin_{2m+1}} \times \Tcal_{\SL_2} & \longrightarrow \Tcal_{\Spin_{2n+1}} \\
  ((x_1, \dots, x_m, s), t) & \longmapsto (\beta((x_1, \dots, x_m), t), s^{2d+1} t^{d(d+1)/2}).
\end{align*}
The group \(\Spin_{2m+1} \times \SL_2\) naturally embeds in the group
\[ G := \left\{ (g_1,g_2) \in \GSpin_{2m+1} \times \GL_2 \,\middle|\, \nu(g_1) = (\det g_2)^{m(m+1)/2} \right\} \]
and we have a short exact sequence
\[ 1 \to \mu_2 \to \Spin_{2m+1} \times \SL_2 \times \GL_1 \to G \to 1 \]
where the first map is
\begin{align*}
  \mu_2 & \longrightarrow Z(\Spin_{2m+1}) \times Z(\SL_2) \times \GL_1 \\
  z & \longmapsto (z^{m(m+1)/2}, z, z)
\end{align*}
and the second map is \((h_1, h_2, \lambda) \mapsto (\lambda^{m(m+1)/2} h_1, \lambda h_2)\).
The morphism
\begin{align*}
  \Spin_{2m+1} \times \SL_2 \times \GL_1 & \longrightarrow \GSpin_{2n+1} \\
  (h_1, h_2, \lambda) & \longmapsto \lambda^{n(n+1)/2} \tilde{\beta}(h_1, h_2)
\end{align*}
induces an extension \(G \to \GSpin_{2n+1}\) of \(\tilde{\beta}\) because \(m(m+1)/2 + d(d+1)/2 + n(n+1)/2\) is even.
We still denote it by \(\tilde{\beta}\).

\begin{prop} \label{pro:rho_odd_GSpin_pi_2dp1}
  As above assume \(\psi = \pi[2d+1]\) with \(d>0\).
  \begin{enumerate}
  \item For any prime \(\ell\) and any \(\iota: \C \simeq \Qellbar\) the \(\GSpin_{2n+1}(\Qellbar)\)-conjugacy classes of \(\rho_{\psi,\iota}^{\GSpin}\) and \(\tilde{\beta} \circ (\rho_{\pi[1],\iota}^{\GSpin}, 1 \oplus \chi_\ell^{-1})\) are equal.
  \item For all primes \(p\) we have \(\cpsc(\psi) = \tilde{\beta}(\cpsc(\pi[1]), \diag(p^{1/2}, p^{-1/2}))\).
  \end{enumerate}
\end{prop}
\begin{proof}
  Note that we already know that for all primes \(p\) we have \(c_p(\psi) = \beta(c_p(\pi[1]), \diag(p^{1/2}, p^{-1/2}))\) by the very definition of \(c_p(\psi)\).
  \begin{enumerate}
  \item The composition of the continuous morphism
    \[ \rho := \tilde{\beta} \circ (\rho_{\pi[1],\iota}^{\GSpin}, 1 \oplus \chi_\ell^{-1}): \GalQ \to \GSpin_{2n+1}(\Qellbar) \]
    with the projection \(\GSpin_{2n+1}(\Qellbar) \to \SO_{2n+1}(\Qellbar)\) is conjugated to \(\rho_{\psi,\iota}^{\SO}\), by uniqueness in Theorem \ref{thm:existence_rho_SO}.
    The composition of \(\rho\) with \(\nu\) is equal to \(\chi_\ell^{-n(n+1)/2}\), again using \v{C}ebotarev's theorem.
    Using the second characterization in Theorem \ref{theo:odd_spin_rep} we conclude that \(\rho\) is conjugated to \(\rho_{\psi,\iota}^{\GSpin}\).
  \item For a prime \(p\) we choose \(\ell \neq p\) and \(\iota\) arbitrarily and apply the first point at \(\Frob_p\).
  \end{enumerate}
\end{proof}

\begin{rema} \label{rem:odd_pi_2dp1_lift}
  As for Proposition \ref{pro:no_amb_symplectic} and Remark \ref{rem:pi_2d_lift}, one could perhaps prove the second point in Proposition \ref{pro:rho_odd_GSpin_pi_2dp1} directly using Eisenstein series.
\end{rema}

\section{Even spin Galois representations}
\label{sec:even_spin}

\subsection{Local-global compatibility for \(\SObf_{4n}\) yields \(\GSpin\)-valued Galois representations}

\begin{theo} \label{thm:existence_rho_O}
  Fix a prime \(\ell\) and a field isomorphism \(\iota : \C \simeq \Qellbar\).
  Let \(n \geq 1\).
  For \(\tilde{\tau} \in \ICcalt(\SObf_{4n})\) and \(\psi \in \tilde{\Psi}_{\disc}^{\unr, \tilde{\tau}}(\SObf_{4n})\) there exists a continuous semisimple morphism \(\rho_{\psi, \iota}^{\mathrm{O}} : \GalQ \rightarrow \SO_{4n}(\Qellbar)\) unramified away from \(\ell\), crystalline at \(\ell\) and such that for any prime number \(p \neq \ell\) the \(\mathrm{O}_{4n}(\Qellbar)\)-conjugacy class of \(\rho_{\psi, \iota}^{\mathrm{O}}(\Frob_p)^{\mathrm{ss}}\) is equal to \(\tilde{c}_p(\psi)\).
  The \(\mathrm{O}_{4n}(\Qellbar)\)-conjugacy class of the Hodge-Tate cocharacter of \(\rho_{\psi, \iota}^{\mathrm{O}}\) (recall that this is a \(\SO_{4n}(\Qellbar)\)-conjugacy class of morphisms \(\GL_{1,\Qellbar} \to \SO_{4n,\Qellbar}\)) is equal to \(\tilde{\tau}\).
  Any another continuous semisimple morphism \(\GalQ \to \SO_{4n}(\Qellbar)\) satisfying the same condition at almost all primes \(p \neq \ell\) is conjugated by \(\mathrm{O}_{4n}(\Qellbar)\) to \(\rho_{\psi,\iota}^\mathrm{O}\).
\end{theo}
\begin{proof}
  Identical to that of Theorem \ref{thm:existence_rho_SO}, except that for even dimensions conjugacy under the orthogonal group does not imply conjugacy under the special orthogonal group.
\end{proof}

We will refine this theorem in most cases in Theorem \ref{thm:loc_glob_param_SO_even} below, which will be proved in the next sections.

\begin{coro} \label{cor:Sen_not_thetahat_inv}
  Let \(n \geq 1\), \(\tilde{\tau} \in \ICcalt(\SObf_{4n})\) and \(\psi \in \tilde{\Psi}_{\disc}^{\unr, \tilde{\tau}}(\SObf_{4n})\).
  The set of primes \(p\) such that \(\tilde{c}_p(\psi)\) consists of two conjugacy classes in \(\SO_{4n}(\C)\) (swapped by \(\thetahat\)) has (natural) density \(1\), in particular it is infinite.
\end{coro}
\begin{proof}
  Choose any prime \(\ell\) and \(\iota: \C \simeq \Qellbar\).
  By assumption on \(\taut\) and Theorem \ref{thm:existence_rho_O} the Hodge-Tate cocharacter of \(\rho_{\psi,\iota}^{\mathrm{O}}\) is not fixed by the outer automorphism group.
  The corollary thus follows from Sen's theorem \cite{Sen} and \v{C}ebotarev's density theorem.
\end{proof}

\begin{coro} \label{cor:two_SO_cc_rho_O}
  Let \(n \geq 1\), \(\taut \in \ICcalt(\SObf_{4n})\) and \(\psi \in \Psit_{\disc}^{\unr, \taut}(\SObf_{4n})\).
  For any prime \(\ell\) and any \(\iota: \C \simeq \Qellbar\) the \(\mathrm{O}_{4n}(\Qellbar)\)-conjugacy class of \(\rho_{\psi,\iota}^\mathrm{O}\) in Theorem \ref{thm:existence_rho_O} consists of two conjugacy classes under \(\SO_{4n}(\Qellbar)\), distinguished by the conjugacy class of the semisimple part of the image of \(\Frob_p\) for any \(p \neq \ell\) such that \(\tilde{c}_p(\psi)\) consists of two conjugacy classes in \(\SO_{4n}(\C)\).
\end{coro}

\begin{defi} \label{def:bad_tau}
  Let \(n \geq 1\) and \(\tilde{\tau} \in \ICcalt(\SObf_{4n})\).
  There are (uniquely determined) integers \(w_1 > \dots > w_{2n} > 0\) such that \(\tilde{\tau}\) is represented by \((w_1, \dots, w_{2n})\) using the parametrization \eqref{eq:param_T_SO_even}.
  We say that \(\tilde{\tau}\) is \emph{bad} if \(n\) is odd and for all \(1 \leq i \leq n\) we have \(w_{2i-1} = w_{2i} + 1\).
\end{defi}

\begin{theo} \label{thm:loc_glob_param_SO_even}
  Let \(n \geq 1\), \(\taut \in \ICcalt(\SObf_{4n})\) and \(\psi \in \tilde{\Psi}_{\disc, \nonendo}^{\unr, \taut}(\SObf_{4n})\).
  Assume either \(n=1\), \(n\) even, \(\psi = \pi[2d]\) or that \(\taut\) is not bad (Definition \ref{def:bad_tau}).
  Then there exists a continuous semisimple morphism \(\rho_{\psi,\iota}^{\SO}: \GalQ \to \Mpsi(\Qellbar)\) which is unramified away from \(\ell\) and crystalline at \(\ell\) and such that for any prime number \(p \neq \ell\) we have \(\rho_{\psi,\iota}^{\SO}(\Frob_p)^{\sesi} \in \iota(c_p(\psi))\).
  Moreover \(\rho_{\psi,\iota}^{\SO}\) is unique in the following sense: any continuous semisimple \(\rho : \GalQ \rightarrow \Mpsi(\Qellbar)\) such that for almost all primes \(p\), \(\rho\) is unramified at \(p\) with \(\rho(\Frob_p)^{\sesi} \in \iota( c_p(\psi) )\), is conjugated to \(\rho_{\psi, \iota}^{\SO}\) by an element of \(\Mpsi(\Qellbar)\).
\end{theo}
We will prove existence in Proposition \ref{pro:loc_glob_SO4} (\(n=1\)), Proposition \ref{pro:Gal_GSpin_pi_2d} (\(\psi=\pi[2d]\)), Corollary \ref{cor:loc_glob_param_SO8n} (\(n\) even) and Proposition \ref{pro:loc_glob_SO8nm4} (\(n>1\) odd and \(\taut\) not bad).
\begin{proof}[Proof of uniqueness]
  Uniqueness follows from uniqueness of \(\rho_{\psi,\iota}^{\mathrm{O}}\) up to conjugation by \(\mathrm{O}_{4n}(\Qellbar)\) in Theorem \ref{thm:existence_rho_O} and Corollary \ref{cor:two_SO_cc_rho_O}.
\end{proof}

\begin{theo} \label{theo:even_GSpin_Gal_rep}
  Let \(n \geq 1\), \(\tilde{\tau} \in \ICcalt(\mathbf{SO}_{4n})\) and \(\psi \in \tilde{\Psi}_{\disc, \nonendo}^{\unr, \tilde{\tau}}(\mathbf{SO}_{4n})\).
  If Theorem \ref{thm:loc_glob_param_SO_even} holds for \(\psi\) then there exists a continuous semisimple morphism \(\rho_{\psi, \iota}^{\GSpin} : \GalQ \rightarrow \GMpsisc(\Qellbar)\) unramified away from \(\ell\), crystalline at \(\ell\) and such that for any \(p \neq \ell\) we have \(\rho_{\psi, \iota}^{\GSpin}(\Frob_p)^{\sesi} \in \iota(p^{n/2} \cpsc(\psi))\).
  \begin{enumerate}
  \item If \(\rho: \GalQ \to \GMpsisc(\Qellbar)\) is a continuous semisimple morphism also satisfying this property at almost all primes \(p\) then \(\rho\) is conjugated to \(\rho_{\psi, \iota}^{\GSpin}\).
  \item If \(\rho: \GalQ \to \GMpsisc(\Qellbar)\) is a continuous lift of \(\rho_{\psi,\iota}^{\SO}\) (up to conjugacy) which is unramified away from \(\ell\), crystalline at \(\ell\) and satisfying \(\nu \circ \rho = \chi_\ell^{-n}\) then \(\rho\) is conjugated to \(\rho_{\psi, \iota}^{\GSpin}\).
  \end{enumerate}
\end{theo}
\begin{proof}
  The proof is very similar to that of Theorem \ref{theo:odd_spin_rep}.
  Note that Theorem \ref{thm:loc_glob_param_SO_even} does \emph{not} state that the Hodge-Tate cocharacter \(\tau_{\mathrm{HT}}\) of \(\rho_{\psi,\iota}^{\SO}|_{\GalQell}\) is equal to \(\tau_\psi\) (as would be expected), so we only know from Theorem \ref{thm:existence_rho_O} that \(\tau_\mathrm{HT}\) is conjugated to \(\tau_\psi\) under the outer automorphism group of \((\Mpsi)_{\Qellbar}\).
  There are two natural \(\SO_{4n}(\Qbar)\)-conjugacy classes of identifications \(\Mpsi \simeq \SO_{4n}\), and we fix one (arbitrarily) for this proof.
  This also gives an identification \(\GMpsisc \simeq \GSpin_{4n}\) (see Section \ref{sec:not_red_gps}).
  Under the parametrization \eqref{eq:param_T_SO_even} there is a unique representative \((w_1, \dots, w_{2n})\) of \(\tau_\mathrm{HT}\) where \(w_1 > \dots > w_{2n-1} > |w_{2n}| > 0\) are integers.
  Thus \(\tau_\psi\) is the conjugacy class of \((w_1, \dots, w_{2n-1}, \pm w_{2n})\).
  Using the parametrization \eqref{eq:param_TGSpin_even} of the maximal torus \(\Tcal_{\GSpin_{4n}}\) of \(\GSpin_{4n}\), define a lift \(\tau' = (w_1, \dots, w_{2n}, \frac{1}{2} \sum_i w_i, n/2)\) of \(\tau_\mathrm{HT}\).
  This defines a conjugacy class of cocharacters \(\GL_{1,\Qellbar} \to (\GMpsisc)_{\Qellbar}\) because \(n-\sum_i w_i\) is even (see Lemma \ref{lem:parity_orth}).
  By Proposition \ref{pro:ex_lift_cond_one} there exists a unique geometric lift \(\rho_{\psi, \iota}^{\GSpin}: \GalQ \to \GMpsisc(\Qellbar)\) of \(\rho_{\psi,\iota}^{\SO}\) which is unramified away from \(\ell\), crystalline at \(\ell\) and having Hodge-Tate cocharacter \(\tau'\).
  The composition \(\nu \circ \rho_{\psi,\iota}^{\GSpin}\) is a continuous character unramified away from \(\ell\) and crystalline at \(\ell\) with Hodge-Tate weight \(n\) and so it is equal to \(\chi_\ell^{-n}\).
  Thus for any \(p \neq \ell\) we have \(\rho_{\psi,\iota}^{\GSpin}(\Frob_p)^{\sesi} \in \pm \iota(p^{n/2} \cpsc(\psi))\).
  Applying Proposition \ref{pro:un_cond_one} to \(\spin_\psi^\epsilon \circ \rho_{\psi,\iota}^{\GSpin}\) and \(\sigma_{\psi,\iota}^{\spin,\epsilon}\) for \(\epsilon\) as in the second part of Corollary \ref{coro:seed_ex_sigma} allows us to conclude \(\rho_{\psi,\iota}^{\GSpin}(\Frob_p)^{\sesi} \in \iota(p^{n/2} \cpsc(\psi))\) for all \(p \neq \ell\).

  The first characterization of \(\rho_{\psi,\iota}^{\GSpin}\) is proved essentially as in Theorem \ref{theo:odd_spin_rep}.
  Let \(\rho: \GalQ \to \GMpsisc(\Qellbar)\) be another continuous morphism satisfying the same condition at almost all primes \(p\).
  Composing \(\rho\) with the faithful representation \(\spin_\psi^\epsilon \oplus \Std\) of \(\GMpsisc\) we obtain a representation isomorphic to \(\sigma_{\psi,\iota}^{\spin,\epsilon} \oplus \Std \circ \rho_{\psi,\iota}^{\mathrm{O}}\) and we deduce that \(\rho\) is also unramified away from \(\ell\) and crystalline at \(\ell\).
  Composing with the projection \(\GMpsisc \to \Mpsi\) we deduce from Theorem \ref{thm:loc_glob_param_SO_even} (uniqueness part) that up to conjugation by \(\Mpsi(\Qellbar)\), \(\rho\) is also a lift of \(\rho_{\psi,\iota}^{\SO}\).
  We deduce as in the proof of Theorem \ref{theo:odd_spin_rep} that up to conjugation by \(\GMpsisc(\Qellbar)\) we have \(\rho = \chi_\ell^N \rho_{\psi,\iota}\) for some integer \(N\), and then that we have \(N=0\) by composing with \(\nu: \GMpsisc \to \GL_1\) (see Definition \ref{def:GMpsisc}) and considering the image of a Frobenius element at any prime \(p \neq \ell\).

  The second characterization of \(\rho_{\psi,\iota}^{\GSpin}\) follows from uniqueness in Proposition \ref{pro:ex_lift_cond_one}.
\end{proof}

\subsection{Local-global compatibility for \(\SObf_4\)}
\label{sec:SO4}

We first handle the case \(n=1\) (i.e.\ the case of \(\SObf_4\)) directly.
The adjoint group \(\PGSObf_4\) of \(\SObf_4\) is isomorphic to \(\PGLbf_2 \times \PGLbf_2\), say with the first simple root \(\alpha_1\) (see labelling after \eqref{eq:param_T_SO_even}) corresponding to the second factor \(\PGLbf_2\).
The dual morphism \(\SL_2 \times \SL_2 \to \SO_4\) restricts to
\begin{align*}
  \Tcal_{\SL_2} \times \Tcal_{\SL_2} & \longrightarrow \Tcal_{\SO_4} \\
  (t_1, t_2) & \longmapsto (t_1 t_2, t_1/t_2)
\end{align*}
and we recognize the morphism \(\alpha_{1,1}\) of Definition \ref{def:alpha_m_d}.
Note that \(\tilde{\alpha}_{1,1}: G(\SL_2 \times \SL_2) \to \GSpin_4\) is an isomorphism.

Recall that for an integer \(k \geq 1\), discrete level one automorphic representations \(\pi\) for \(\PGLbf_2\) such that \(\pi_\infty\) has infinitesimal character \(k-1/2 \in \Lie \Tcal_{\SL_2}\) are either trivial (if \(k=1\)) or are in bijection with eigenforms in \(S_{2k}(\SLbf_2(\Z))\), the latter ones being cuspidal.
In both cases we know the existence of continuous semi-simple (irreducible in the latter case) Galois representations, for any prime \(\ell\) and any \(\iota: \C \simeq \Qellbar\), \(\rho_{\pi,\iota}: \GalQ \to \GL_2(\Qellbar)\) unramified away from \(\ell\), crystalline at \(\ell\) and such that for any \(p \neq \ell\) we have \(\rho_{\pi,\iota}(\Frob_p)^{\sesi} \in \iota(p^{1/2} c_p(\pi))\), unique up to conjugation.
In the first case we simply have \(\rho_{1,\iota} = 1 \oplus \chi_\ell^{-1}\), in the second case \(\rho_{\pi,\iota}\) is up to a twist the Galois representation associated by Deligne to \(\pi\) and \(\iota\).
Note that we have \(\det \rho_{\pi,\iota} = \chi_\ell^{-1}\).

\begin{prop} \label{pro:loc_glob_SO4}
  Let \(\taut \in \ICcalt(\SObf_4)\), say represented by \(\tau = (w_1, w_2)\) with \(w_1 > w_2 > 0\) integers.
  Then \(\psi \in \Psit^{\unr,\taut}_{\disc}(\SObf_4)\) is either
  \begin{itemize}
  \item \((\pi_1 \otimes \pi_2)[1]\) (see Proposition \ref{pro:cpsc_tensor_PGL2_PGL2} for the notation) where \(\pi_1\) and \(\pi_2\) level one cuspidal automorphic representation for \(\PGLbf_2\) having infinitesimal characters \(\pm (w_1+w_2)/2\) and  \(\pm (w_1 - w_2)/2\) at the real place,
  \item \(\pi[2]\) for some level one cuspidal automorphic representation for \(\PGLbf_2\) having infinitesimal character \(\pm (w_1 - \tfrac{1}{2})\) at the real place (this case can only occur if \(w_1 = w_2+1\)).
  \end{itemize}
  In both cases Theorem \ref{thm:loc_glob_param_SO_even} holds.
  Up to conjugation there is a unique isomorphism \(\Mpsisc \simeq \SL_2 \times \SL_2\) mapping \(\tau_\psi\) to \((\pm (w_1+w_2)/2, \pm (w_1-w_2)/2)\).
  It extends uniquely to an isomorphism \(f: \GMpsisc \simeq G(\SL_2 \times \SL_2)\) restricting to the identity on \(\GL_1\).
  Denoting \((\pi_1,\pi_2) = (\pi, 1)\) in the second case (here \(1\) denotes the trivial automorphic representation for \(\PGLbf_2\)) we have in both cases \(f \circ \rho_{\psi,\iota}^{\GSpin} = (\rho_{\pi_1,\iota}, \rho_{\pi_2, \iota})\) (up to conjugation) in Theorem 6.1.6.
\end{prop}
\begin{proof}
  This follows directly from Proposition \ref{pro:cpsc_tensor_PGL2_PGL2}.
\end{proof}

\begin{coro} \label{cor:spin_pm_SO4}
  Let \(\psi=\pi[d]\) be a parameter for \(\SObf_4\) as in Proposition \ref{pro:loc_glob_SO4}, associated to a pair \((\pi_1, \pi_2)\) of discrete automorphic representations for \(\PGLbf_2\) (so \(\pi_2\) may be the trivial representation).
  We have
  \[ \spin_\psi^+ \circ \rho_{\psi,\iota}^\GSpin \simeq \rho_{\pi_1, \iota} \]
  and
  \[ \spin_\psi^- \circ \rho_{\psi,\iota}^\GSpin \simeq
    \begin{cases}
      \rho_{\pi_2, \iota} & \text{ if } \pi_2 \text{ is cuspidal}, \\
      1 \oplus \chi_\ell^{-1} & \text{ if } \pi_2 \text{ is trivial}.
    \end{cases}
  \]
\end{coro}
\begin{proof}
  Compare images of Frobenius elements.
\end{proof}

\subsection{Local-global compatibility for parameters \(\pi[2d]\)}
\label{sec:even_spin_pi_2d}

Let \(n \geq 1\), \(\tilde{\tau} \in \ICcal(\mathbf{SO}_{4n})\) and \(\psi \in \tilde{\Psi}_{\disc, \nonendo}^{\unr, \tilde{\tau}}(\mathbf{SO}_{4n})\).
Assume that we are in the situation of Proposition \ref{pro:no_amb_symplectic}, i.e.\ \(\psi = \pi[2d]\) where \(\pi\) is a self-dual automorphic cuspidal representation for \(\GLbf_{2m}\) (necessarily of symplectic type and with \(n=md\)).
In this case it turns out that we can prove Theorem \ref{thm:loc_glob_param_SO_even} using Galois-theoretic arguments.

\begin{theo} \label{theo:existence_rho_GSp}
  Let \(m \geq 1\).
  Let \(\pi\) be a self-dual cuspidal automorphic representation for \(\GLbf_{2m}\) such that the infinitesimal character of \(\pi_\infty\) is regular with eigenvalues in \(\tfrac{1}{2} \Z\).
  Let \(\iota: \C \simeq \Qellbar\) be a field isomorphism.
  There exists a continuous semisimple morphism \(\rho_{\pi,\iota}^{\GSp}: \GalQ \to \GSp_{2m}(\Qellbar)\), which is unramified away from \(\ell\), crystalline at \(\ell\) and such that for all \(p \neq \ell\) we have \(\rho_{\pi,\iota}^{\GSp}(\Frob_p)^{\sesi} \in p^{1/2} c_p(\pi)\) (here \(c_p(\pi)\) is considered as a conjugacy class in \(\Sp_{2m}(\C)\)), unique up to conjugation.
\end{theo}
\begin{proof}
  This follows from \cite[Theorem 4.2]{CheHar} and \cite[Corollary 1.3]{BC} as in the proof of Theorem \ref{thm:existence_rho_O}, uniqueness follows from \cite[Lemma 6.1]{GanTakeda_LLC_GSp4}.
\end{proof}

\begin{prop} \label{pro:Gal_GSpin_pi_2d}
  Recall the morphism \(\tilde{\alpha}_\psi: G(\Sp_{2m} \times \SL_2) \to \GMpsisc\) from Definition \ref{def:alpha_m_d}.
  \begin{enumerate}
  \item For all primes \(p\) we have \(\cpsc(\psi) = \tilde{\alpha}_\psi(c_p(\pi), \diag(p^{1/2}, p^{-1/2}))\).
  \item Let \(\iota: \C \simeq \Qellbar\) be a field isomorphism.
    Theorem \ref{thm:loc_glob_param_SO_even} holds for \(\psi = \pi[2d]\) and \(\rho_{\psi,\iota}^{\GSpin}\) (defined in Theorem \ref{theo:even_GSpin_Gal_rep}) is conjugated to \(\tilde{\alpha}_{\psi} \circ (\rho_{\pi,\iota}^{\GSp}, 1 \oplus \chi_\ell^{-1})\).
  \end{enumerate}
\end{prop}
\begin{proof}
  Thanks to Proposition \ref{pro:no_amb_symplectic} we already know \(c_p(\psi) = \alpha_\psi(c_p(\pi), \diag(p^{1/2}, p^{-1/2}))\) for all primes \(p\).
  It follows that the composition \(\rho_{\psi,\iota}^{\SO}\) of \(\tilde{\alpha}_\psi \circ (\rho_{\pi,\iota}^{\GSp}, 1 \oplus \chi_\ell^{-1}): \GalQ \to \GMpsisc(\Qellbar)\) with the projection \(\GMpsisc(\Qellbar) \to \Mpsi(\Qellbar)\) satisfies the conditions in Theorem \ref{thm:loc_glob_param_SO_even}.
  The morphism \(\tilde{\alpha}_\psi \circ (\rho_{\pi,\iota}^{\GSp}, 1 \oplus \chi_\ell^{-1})\) is a continuous lift of \(\rho_{\psi,\iota}^{\SO}\) which is unramified away from \(\ell\), crystalline at \(\ell\) and whose composition with \(\nu\) is \(\chi_\ell^{-n}\), so the second characterization of \(\rho_{\psi,\iota}^{\GSpin}\) in Theorem \ref{theo:even_GSpin_Gal_rep} implies that they are conjugated.
\end{proof}

\subsection{The other half spin Galois representation for \(\SObf_{8n}\)}
\label{sec:other_half_spin}

Let \(n \geq 1\) be an integer.
Let \((L,q)\) be an even unimodular lattice of rank \(8n\) and \(\Gbf\) the corresponding special orthogonal group \(\SObf(L,q)\).
In particular \(\Gbf\) is a semi-simple connected reductive group over \(\Z\), \(\Gbf(\R)\) is compact and connected and for any prime \(p\) the group \(\Gbf_{\Zp}\) is split.
The choice of \((L,q)\) will be irrelevant in the applications below.
We may realize \(\Gbf_{\Q}\) as a pure inner form of \(\SObf_{8n,\Q}\) by choosing a basis of \(\Qbar \otimes_{\Z} L\) splitting \(q\), giving identifications \(\Ghat \simeq \SO_{8n}\) and \(\Ghat_{\sico} \simeq \Spin_{8n}\) (see Section \ref{sec:not_red_gps}).
Changing the realization of \(\Gbf_{\Q}\) as a pure inner form of \(\SObf_{8n,\Q}\) yields the same identifications up to \(\SO_{8n}(\Qbar) \rtimes \{1,\thetahat\}\).

For each prime \(\ell\) we choose a Borel pair \((\Bbf_\ell, \Tbf_\ell)\) of the split group \(\Gbf_{\Qell,\mathrm{ad}}\).
Let \(\Delta \subset X^*(\Tbf_\ell)\) be the set of simple roots of \(\Tbf_\ell\) in \(\Gbf_{\Qell,\mathrm{ad}}\) with respect to \(\Bbf_\ell\).
The maximal compact subgroup \(T_0\) of \(\Tbf_\ell(\Qell)\) is isomorphic to \((\Zell^\times)^\Delta \simeq (\Zell^\times)^{4n}\).
The weight space \(\Wscr\) is the rigid analytic space over \(\Qell\) parametrizing locally analytic (equivalently, continuous) characters of \(T_0\).
It is isomorphic to the product of an open polydisc of dimension \(4n\) and a rigid space finite over \(\Qell\).
We will be particularly interested in the subset of \(\Wscr(\Qell)\) consisting of algebraic and dominant characters, i.e.\ characters \(T_0 \to \Qell^\times\) induced by elements of \(X^*(\Tbf_\ell)\) which are dominant with respect to \(\Bbf_\ell\).

\begin{defi} \label{def:simple_ladic_family}
  Let \(\ell\) be a prime, \(\iota : \C \simeq \Qellbar\) a field isomorphism and \(E\) a finite subextension of \(\Qellbar/\Qell\).
  A simple \(\ell\)-adic family of level one automorphic representations for \(\Gbf_\mathrm{ad}\) is a smooth rigid analytic curve\footnote{Since we will only be interested in local properties around a point \(x_0\) of \(\Cscr(E)\), one could take for \(\Cscr\) the open unit disk.} \(\Cscr\) over \(E\) endowed with
  \begin{enumerate}
  \item a morphism \(w: \Cscr \to \Wscr\) such that \(\Cscr\) is finite over an open affinoid of \(\Wscr\),
  \item a point \(x_0 \in \Cscr(E)\) and a subset \(Z\) of the set of points of \(\Cscr\), accumulating at \(x_0\) (i.e.\ any neighborhood of \(x_0\) contains a point of \(Z\) different from \(x_0\)),
  \item a morphism of \(\Ocal_E\)-algebras \(\Xi : \Hcal^{(\infty,\ell)} \rightarrow \Ocal(\Cscr)\) where \(\Hcal^{(\infty, \ell)} = \bigotimes'_{p \neq \ell} \Hcal(\Gbf_{\mathrm{ad},\Zp}, \Ocal_E)\) is the unramified Hecke algebra away from \(\ell\) over \(\Ocal_E\),
  \end{enumerate}
  such that for any finite subextension \(E'\) of \(\Qellbar/E\) and any \(x \in \Cscr(E')\) taking values in \(\{x_0\} \cup Z\), the character \(w(x): T_0 \to \Ocal_{E'}^\times\) is algebraic dominant and there exists a level one automorphic representation \(\pi \simeq \pi_\infty \otimes \pi_\ell \otimes \pi^{(\infty,\ell)}\) such that \(\pi_\infty\) is isomorphic to the restriction to \(\Gbf_\mathrm{ad}(\R)\) of the irreducible algebraic representation of \(\Gbf_{\mathrm{ad},\C}\) corresponding to \(w(x)\), and the action of \(\Hcal^{(\infty, \ell)}\) on the one-dimensional \(\Qellbar\)-vector space \(\iota^* \left(\pi^{(\infty, \ell)}\right)^{\prod_{p \neq \ell} \Gbf_\mathrm{ad}(\Zp)}\) is by the character \(\Xi_x\).
  We will denote it by \((\Cscr, x_0, Z, w, \Xi)\), or simply \((\Cscr, x_0, Z)\).
  If \(\pi\) corresponds to the point \(x_0 \in \Cscr(E)\) we will also say that the simply \(\ell\)-adic family \((\Cscr, x_0, Z, w, \Xi)\) interpolates \(\pi\) (for the given field isomorphism \(\iota\)).

  For a point \(x \in \Cscr(\Qellbar)\) and \(p \neq \ell\) a prime number we denote by \(c_p(\Xi_x)\) the semi-simple conjugacy class in \(\Ghat_{\sico}(\Qellbar)\) corresponding to the restriction of the specialization \(\Xi_x: \Hcal^{(\infty,\ell)} \to \Qellbar\) to the unramified Hecke algebra \(\Hcal(\Gbf_{\mathrm{ad},\Zp}, \Qellbar)\).

  Let \(x \in \Cscr(E')\) taking values in \(Z\) be such that \(w(x)\) is not invariant under \(\theta\).
  Let \((\Bcal, \Tcal)\) be a Borel pair in \(\Ghat_{\sico}\), in particular \(\Tcal\) is naturally isomorphic to \(\widehat{\Tbf_\ell}\).
  Similarly to Definition \ref{def:spin_pm_psi}, we distinguish the two half-spin representations of \(\Ghat_{\sico}\) as follows.
  The (algebraic dominant) weight \(w(x)\) may be seen as a \(\Ghat_{\mathrm{ad}}\)-orbit of cocharacters taking values in \(\Ghat_{\sico}\), and there is a unique tuple of integers \(k_1 \geq \dots \geq k_{4n} > 0\) such that in the standard representation of \(\Ghat_{\sico}\), \(w(x)\) becomes a direct sum of the \(8n\) characters \(\pm k_i\) of \(\GL_1\) (as usual we identify characters of \(\GL_1\) with \(\Z\)).
  The integer \(\sum_{i=1}^{4n} k_i\) is even.
  We denote by \(\spin^+_x\) (resp.\ \(\spin_x^-\)) the half-spin representation of \(\Ghat_{\sico}\) such that the representation \(\spin^+ \circ w(x)\) (resp.\ \(\spin^- \circ w(x)\)) of \(\GL_1\) is the direct sum, over all \((\epsilon_i)_i \in \{\pm 1\}^{4n}\) such that the cardinality of \(\{i \in \{1,\dots,4n\} \,|\, \epsilon_i = +1\}\) is even (resp.\ odd), of the character \(\frac{1}{2} \sum_i \epsilon_i k_i\).
\end{defi}

\begin{prop} \label{pro:interpol_Gal_rep}
  Let \(\ell\) be a prime, \(\iota: \C \simeq \Qellbar\) a field isomorphism, \(E\) a finite subextension of \(\Qellbar/\Qell\) and \((\Cscr, x_0, Z, w, \Xi)\) a simple \(\ell\)-adic family of level one automorphic representations for \(\Gbf_\mathrm{ad}\), interpolating an automorphic representation \(\pi\).
  Let \(r\) be an algebraic representation of \(\Ghat_{\sico,\Qellbar}\) of dimension \(d\).
  Assume that for any \(x \in \Cscr(\Qellbar)\) taking values in \(Z\) there exists a continuous semisimple representation \(\rho_x: \GalQ \to \GL_d(\Qellbar)\) which is unramified away from \(\ell\) and such that for any prime number \(p \neq \ell\) the semi-simplification \(\rho_x(\Frob_p)^{\sesi}\) belongs to the image of \(c_p(\Xi_x)\) under \(r\).
  Then there exists a continuous semisimple representation \(\rho_{x_0}: \GalQ \to \GL_d(\Qellbar)\) satisfying the same properties at \(x_0\).
\end{prop}
\begin{proof}
  This follows from \cite[Proposition 7.1.1]{TheseG} as explained in \S 7.4 loc.\ cit.
\end{proof}

\begin{defi}
  Let \(H\) be a split semisimple algebraic group over \(\Qell\), \(\Gamma\) a profinite topological group, and \(\rho : \Gamma \rightarrow H(\Qellbar)\) a continuous morphism.
  Then \(\Lie \left( \rho(\Gamma) \right)\) is a finite-dimensional sub-\(\Qell\)-vector space of \(\Qellbar \otimes_{\Qell} \Lie(H)\).
  We say that \(\rho\) has \emph{maximal infinitesimal image} if \(\Lie \left( \rho(\Gamma) \right)\) spans the \(\Qellbar\)-vector space \(\Qellbar \otimes_{\Qell} \Lie (H)\).
\end{defi}

\begin{theo} \label{theo:TaiCC}
  Let \(\pi\) be a level one automorphic representation for \(\Gbf_\mathrm{ad}\), whose infinitesimal character is not invariant under \(\theta^*\).
  Let \(E / \Qell\) be a finite subextension of \(\Qellbar / \Qell\) containing the image by \(\iota\) of a number field \(\subset \C\) over which \(\pi_f\) is defined.
  Let \(r\) be one of the two half-spin representations of \(\Ghat_{\sico}\).

  There exist
  \begin{itemize}
  \item a simple \(\ell\)-adic family \((\Cscr, x_0, Z)\) interpolating \(\pi\),
  \item for each \(x_1 \in Z\), a simple \(\ell\)-adic family \((\Cscr_{x_1}, x_1', Z_{x_1})\) interpolating the representation \(\pi_1\) corresponding to \(x_1\),
  \item for each \(x_1 \in Z\) and each \(x_2 \in Z_{x_1}\), a simple \(\ell\)-adic family \((\Cscr_{x_1, x_2}, x_2', Z_{x_1, x_2})\) interpolating the representation \(\pi_2\) corresponding to \(x_2\),
  \end{itemize}
  such that
  \begin{itemize}
  \item for any \(x \in Z\) (resp.\ \(x_1 \in Z\) and \(x \in Z_{x_1}\), resp.\ \(x_1 \in Z\), \(x_2 \in Z_{x_1}\) and \(x \in Z_{x_2}\)), \(w(x)\) is not invariant under \(\theta^*\) and we have \(\spin_x^+ \simeq r\) (recall that \(\spin_x^+\) was introduced in Definition \ref{def:simple_ladic_family}),
  \item for any \(x_1 \in Z\), \(x_2 \in Z_{x_1}\) and \(x_3 \in Z_{x_1, x_2}\), the Galois representation \(\rho_{\pi_3, \iota}^{\mathrm{O}}\) attached to the level one automorphic representation \(\pi_3\) for \(\Gbf_\mathrm{ad}\) corresponding to \(x_3\) has maximal infinitesimal image.
  \end{itemize}
\end{theo}
\begin{proof}
  The proof is almost identical to that of \cite[Corollary 4.0.2]{TaiCC}, which itself is a simpler version of the proof of Theorem 3.2.2 and Corollary 3.2.3 loc.\ cit.\ so we simply highlight the differences.

  In \S 4 loc.\ cit.\ we worked over a totally real number field of even degree.
  This assumption was made only to guarantee the existence of a non-degenerate quadratic form in dimension \(4\) which is anisotropic at all real places and split at all finite places, but everything works the same over \(\Q\) in dimension divisible by \(8\).

  In loc.\ cit.\ we worked with the group \(\Gbf\), working with \(\Gbf_\mathrm{ad}\) only gives us more Hecke operators.

  In loc.\ cit.\ we worked with eigenvarieties over the whole weight space \(\Wscr\), now we restrict to (i.e.\ take the fiber product over \(\Wscr\) with) certain curves in \(\Wscr\) to simplify our families.
  More precisely, using the parametrization \eqref{eq:param_T_SO_even}, for \(Z\) we restrict to weights of the form
  \[ \ul{t} = (t_1, \dots, t_{4n}) \mapsto w(x_0)(\ul{t}) \prod_{i=1}^{4n} \chi(t_i)^{4n+1-i} \]
  for some (uniquely determined) character \(\chi\) of \(\Zell^\times\), and similarly for the other families.
  Note that this does define a weight for \(\Gbf_\mathrm{ad}\), despite being described as a weight for \(\Gbf\), because the integer \(\sum_{i=1}^{4n} 4n+1-i = 2n(4n+1)\) is even.
  In the proof of Corollary 4.0.2 loc.\ cit.\ the arguments only require classical weights \((k_1, \dots, k_{4n})\) such that \(k_1 - k_2\), \dots, \(k_{4n-1} - k_{4n}\) and \(k_{4n}\) are bigger than some constant, and it is clear that weights of the form above satisfying this condition abound.
  %TODO: better explanation
\end{proof}

\begin{coro}
  Let \(\pi = \otimes'_v \pi_v\) be a level one automorphic representation for \(\Gbf_\mathrm{ad}\), and let \(r : \Ghat \rightarrow \GL_{2^{4n-1}}(\C)\) be any of the two half-spin representations.
  Then for any prime number \(\ell\), and any \(\iota : \C \simeq \Qellbar\), there exists a unique continuous semisimple representation \(\sigma_{\pi, \iota}^r : \GalQ \rightarrow \GL_{2^{4n-1}}(\Qellbar)\) which is unramified away from \(\ell\) and such that for any prime \(p \neq \ell\) we have \(\sigma_{\pi, \iota}^r(\Frob_p)^{\sesi} \in \iota( p^n r(c_p(\pi_p)) )\).
\end{coro}
\begin{proof}
  Note that in Theorem \ref{theo:TaiCC}, for any \(x_1 \in Z, x_2 \in Z_{x_1}\) and \(x_3 \in Z_{x_1, x_2}\), the parameter in \(\tilde{\Psi}^{\unr, \tilde{\tau}}_{\disc}(\mathbf{SO}_{8n})\) attached to the automorphic representation \(\pi^{(3)}\) corresponding to \(x_3\) (here \(\tau = w(x_3) + \rho\)) is a single automorphic cuspidal self-dual level one representation for \(\mathbf{GL}_{8n}\), since the associated Galois representation \(\sigma_{\pi^{(3)}, \iota}^{\Std}\) is irreducible.

  For any \(x_1 \in Z\), \(x_2 \in Z_{x_1}\), and \(x_3 \in Z_{x_1, x_2}\) we have \(\spin_{x_3}^+ \simeq r\), so that the existence of \(\sigma_{x_3, \iota}^r\) as in the Corollary follows from the second point of Corollary \ref{coro:seed_ex_sigma} (we are in the first case with \(d=0\) and so \(\epsilon=(-1)^{2n}=1\)).
  Thanks to Proposition \ref{pro:interpol_Gal_rep} these Galois representations interpolate and we obtain in turn the existence of the representations \(\sigma_{x_2, \iota}^r\), \(\sigma_{x_1, \iota}^r\) and finally \(\sigma_{\pi, \iota}^r\).
\end{proof}

\begin{coro} \label{coro:ex_both_spin}
  Let \(\taut \in \ICcalt(\mathbf{SO}_{8n})\) and \(\psi \in \tilde{\Psi}^{\unr, \tilde{\tau}}_{\disc, \nonendo}(\mathbf{SO}_{8n})\).
  For any prime \(\ell\), any \(\iota : \C \simeq \Qellbar\), and any \(\epsilon \in \{ +,-\}\), there is a unique continuous semisimple \(\sigma_{\psi, \iota}^{\spin, \epsilon} : \GalQ \rightarrow \GL_{2^{4n-1}}(\Qellbar)\) which is crystalline at \(\ell\) and unramified away from \(\ell\) and such that for any prime number \(p \neq \ell\) we have \(\sigma_{\psi, \iota}^{\spin, \epsilon}(\Frob_p)^{\sesi} \in \iota (p^n \spin_{\psi}^{\epsilon}(\cpsc(\psi)))\).
\end{coro}
\begin{proof}
  Since \(\Scal_{\psi} = 1\), so if we choose \(\tau \in \ICcal(\SObf_{8n})\) mapping to \(\taut\) the multiplicity formula for \(\Gbf_\mathrm{ad}\) (see\footnote{Using a realization of \(\Gbf_\mathrm{ad}\) as an inner form of \(\PGSObf_{8n}\) as explained at the beginning of Section \ref{sec:other_half_spin}.} Example \ref{exam:mult_formula_definite}) provides us with a unique automorphic level one representation \(\pi\) for \(\Gbf_\mathrm{ad}\) such that \(\pi_\infty\) has infinitesimal character \(\tau\) and for any prime \(p\) we have \(c(\pi_p) = \dpsitausc(c_p(\psi))\).
  The previous corollary provides us with two continuous semisimple representations \(\sigma_{\psi,\iota}^{\spin,\epsilon}\) unramified away from \(\ell\) and with Frobenius elements at all \(p \neq \ell\) having characteristic polynomial as stated.
  We are left to check that both these representations are crystalline at \(\ell\).
  Thanks to Theorem \ref{theo:IH_explicit_crude} we already know that one is, and their tensor product can be written as a sum of Schur functors applied to \(\sigma_{\psi, \iota}^{\Std}\), which we already know to be crystalline at \(\ell\) (see Theorem \ref{thm:existence_rho_O}).
\end{proof}

\subsection{Local-global compatibility for \(\mathbf{SO}_{8n}\)}

\(n \geq 1\) is still fixed, and \(\Gbf = \mathbf{SO}_{8n}\) (definite).

Can now lift the whole spin.
when big image, get local-global compatibility.
Deform this.

\begin{lemm} \label{lemm:dist_hspin_tr}
  Let \(m \geq 2\) be an integer.
  \begin{enumerate}
  \item Let \(c\) be a semisimple conjugacy class in \(\SO_{2m}(\C)\).
    Recall from \cite[\S 19]{FultonHarris_rep} that the representation \(\bigwedge^m \Std\) of \(\SO_{2m}\) decomposes as the sum of two inequivalent irreducible representations permuted by \(\theta\), say \(\delta_1\) and \(\delta_2\).
    If \(c\) does not admit \(\pm 1\) as an eigenvalue in the standard representation then the traces of \(\delta_1(c)\) and \(\delta_2(c)\) are distinct.
  \item Let \(r_1, r_2\) be the two half-spin representations of \(\Spin_{2m}\).
  Let \(c\) be a semisimple conjugacy class in \(\Spin_{2m}(\C)\).
  Assume that \(c\) does not admit \(\pm 1\) as an eigenvalue in the standard representation.
  Then we have
  \[ \Tr \left( r_1(c) \right)^2 \neq \Tr \left( r_2(c) \right)^2. \]
  \end{enumerate}
\end{lemm}
\begin{proof}
  \begin{enumerate}
  \item Also recall from \cite[\S 19]{FultonHarris_rep} that for \(1 \leq j \leq m-1\) the representation \(\bigwedge^j \Std\) of \(\SO_{2m}\) is irreducible and invariant under \(\theta\).
    Moreover the analysis loc.\ cit.\ shows that any irreducible representation of \(\SO_{2m}\) occurs in a tensor product of representations \(\bigwedge^j \Std\) for \(1 \leq j \leq m-1\), \(\delta_1\) and \(\delta_2\).
    This implies (see \cite[\S 6]{Steinberg}) that the \(\C\)-algebra of conjugation-invariant functions on \(\SO_{2m,\C}\) is generated by the traces in the representations \(\bigwedge^j \Std\) for \(1 \leq j \leq m-1\) and the representations \(\delta_1\) and \(\delta_2\).
    Also recall from Corollary 6.6 loc.\ cit.\ that the traces in these representations determine semi-simple elements in \(\SO_{2m}(\C)\) up to conjugacy.
    The assumption on the eigenvalues of \(c\) is equivalent to \(\theta(c) \neq c\), and the first point follows.
  \item The representations \(r_i^{\otimes 2}\) factor through \(\Spin_{2m}(\C) \rightarrow \SO_{2m}(\C)\).
    A simple calculation with dominant weights shows that, up to swapping \(\delta_1\) and \(\delta_2\), each tensor square \(r_i^{\otimes 2}\) decomposes as \(\delta_i \oplus \bigoplus_{1 \leq j \leq m-1} \left( \bigwedge^j \Std \right)^{\oplus e(i,j)}\) where \(e(i,j) \in \Z_{\geq 0}\).
    So this point follows from the previous one.
  \end{enumerate}
\end{proof}

\begin{lemm} \label{lem:Gal_rep_SO_simple_family}
  Let \(\iota: \C \simeq \Qellbar\) be a field isomorphism.
  Let \(\pi\) be a level one automorphic representation for \(\Gbf_\mathrm{ad}\), and \((\Cscr, x_0, Z, w, \Xi)\) a simple \(\ell\)-adic family interpolating \(\pi\) (for \(\iota\)).
  Assume that for any finite subextension \(E'\) of \(\Qellbar/E\) and any \(x \in \Cscr(E')\) taking values in \(Z\) the following conditions are satisfied:
  \begin{itemize}
  \item the weight \(w(x)\) is not invariant under \(\theta\), and
  \item there exists a continuous semisimple morphism \(\rho_x : \GalQ \rightarrow \Ghat(\Qellbar)\) which is unramified away from \(\ell\) and such that for any \(p \neq \ell\) we have \(\rho_x(\Frob_p)^{\sesi} \in \operatorname{pr} \left( c_p(\Xi_x) \right)\), where \(\operatorname{pr}: \Ghat_{\sico} \to \Ghat\) is the natural projection.
  \end{itemize}
  Then there exists a continuous semisimple morphism \(\rho_{x_0} : \GalQ \rightarrow \Ghat(\Qellbar)\) which is unramified away from \(\ell\) and such that for any \(p \neq \ell\) we have
  \[ \rho_{x_0}(\Frob_p)^{\sesi} \in \operatorname{pr} \left( c_p(\Xi_{x_0}) \right) = \operatorname{pr} \left( \iota(c(\pi_p)) \right). \]
\end{lemm}
\begin{proof}
  For any \(x \in Z\) we have \(\sigma_{\pi(x), \iota}^{\Std} \simeq \Std \circ \rho_x\) where the automorphic representation \(\pi(x)\) for \(\Gbf_\mathrm{ad}\) corresponds to \(x\) (using \(\iota\)) as in Definition \ref{def:simple_ladic_family}.
  By \cite[Proposition 7.1.1]{TheseG} there is a unique continuous pseudocharacter \(T : \GalQ \rightarrow \Ocal(\Cscr)\) of dimension \(8n\) interpolating the traces of these representations.
  By normalizing and restricting to a connected component we can assume that \(\Cscr\) is smooth and connected.
  By \cite[Theorem 1.2]{Tay2} there exists a finite extension \(K\) of \(\op{Frac} \Ocal(\Cscr)\) and a representation of \(\GalQ\) over \(K\) having trace \(T\).
  Up to replacing \(K\) by a finite extension we can assume that this representation decomposes as a direct sum of absolutely irreducible representations.
  Using (the proof of) \cite[Lemma 7.8.11]{BCbook} we obtain that, up to replacing \(\Cscr\) by a finite curve over \(\Cscr\), we can assume that we have a decomposition \(T = \sum_i T_i\) where each \(T_i\) is a generically irreducible pseudocharacter and that there are finite projective \(\Ocal_{\Cscr}\)-modules \(L_i\) and continuous representations \(\tilde{\rho}_i : \GalQ \rightarrow \GL_{\Ocal(\Cscr)}(L_i)\) such that \(\Tr (\tilde{\rho}_i) = T_i\).
  Normalizing again, we can still assume that \(\Cscr\) is smooth and connected.
  Using Burnside's theorem we obtain that, up to restricting to a neighbourhood of \(x_0\), for any \(i\) the specialization of \(\tilde{\rho}_i\) at any point of \(\Cscr\) different from \(x_0\) is absolutely irreducible.
  Since generically the eigenvalues of the Sen polynomial of \(\bigoplus_i \tilde{\rho}_i|_{\GalQell}\) are distinct thanks to the hypothesis that \(w(x)\) is not invariant under \(\theta\) for \(x \in Z\), the pseudocharacters \(T_i\)'s are (generically) distinct.
  Moreover \(T\) is self-dual, so for any index \(i\), either \(T_i\) is self-dual or its dual is isomorphic to \(T_j\) for a uniquely determined index \(j \neq i\).
  In the latter case we denote \(i^\vee = j\).
  By \cite[Corollary 1.3]{BC} any self-dual \(T_i\) is generically of orthogonal type, i.e.\ the (unique up to scalar) non-degenerate \(\tilde{\rho}_i(\GalQ)\)-invariant bilinear form \(\langle \cdot, \cdot \rangle_i\) on \(\Frac \left( \Ocal(\Cscr) \right) \otimes L_i\) is symmetric.
  Up to multiplying it by a non-zero element of \(\Ocal(\Cscr)\) and restricting \(\Cscr\), we can assume that \(\langle \cdot, \cdot \rangle_i\) takes values in \(\Ocal(\Cscr)\) on \(L_i \times L_i\) and is residually non-degenerate away from \(x_0\).
  For an index \(i\) such that \(T_i\) is not self-dual, we can replace \(L_{i^{\vee}}\) by \(L_i^{\vee}\), and we endow \(L_i \oplus L_{i^{\vee}}\) with the non-degenerate symmetric bilinear form \(\langle (v_1, v_1^{\vee}), (v_2, v_2^{\vee}) \rangle = v_1^{\vee}(v_2) + v_2^{\vee}(v_1)\).
  
  Let \(L := \bigoplus_i L_i\), endowed with \(\tilde{\rho} = \bigoplus_i \tilde{\rho}_i : \GalQ \rightarrow \GL_{\Ocal(\Cscr)}(L)\) and \(\langle \cdot, \cdot \rangle\).
  Let \(A = \Ocal_{\Cscr, x_0}\), a discrete valuation ring, and let \(\varpi\) be a uniformizer.
  Let \(n \geq 1\) be minimal and such that \(\varpi^n L^{\vee} \subset L\).
  If \(n>1\), let \(L' = L + \varpi^{n-1} L^{\vee}\), which properly contains \(L\), is such that \(\langle \cdot, \cdot \rangle\) is still \(A\)-valued on \(L' \times L'\), and is stable under \(\tilde{\rho}(\GalQ)\).
  Replacing \(L\) by \(L'\) and iterating this procedure, we can assume that we have \(n=1\), i.e.\
  \[ \varpi L \subset \varpi L^{\vee} \subset L. \]
  Let \(\pi = \sqrt{\varpi}\) and \(A' = A[\pi]\).
  Replacing \(A\) by \(A'\) (thereby replacing \(\Cscr\) by a finite flat cover) and \(L\) by \(A'L + \pi A' L^{\vee}\), we can finally assume that the \(\tilde{\rho}(\GalQ)\)-stable lattice \(L\) is self-dual, i.e.\ \(L^{\vee} = L\).
  Up to restricting \(\Cscr\) we can assume that \(L\) is free over \(\Ocal_{\Cscr}\), and after replacing \(\Cscr\) by a finite étale cover we can assume that the form \(\langle \cdot, \cdot \rangle\) on \(L\) is split, i.e.\ that there exists a basis of \(L\) such that the Gram matrix of \(\langle \cdot, \cdot \rangle\) is equal to \(J_{2n}\).
  Choosing such a basis identifies \(\tilde{\rho}\) with a morphism \(\GalQ \to \SObf_{8n}(\Ocal(\Cscr))\).
  We can specialize \(\tilde{\rho}\) at \(x_0\) to obtain a continuous morphism \(\tilde{\rho}_{x_0}: \GalQ \to \SObf_{8n}(\Qellbar)\), and we have \(\Std \circ \tilde{\rho}_{x_0} \simeq \sigma_{\pi, \iota}^{\Std}\).
  In this proof it is convenient to define the split connected reductive group \(\Ghat\) over \(\Q\), instead of \(\Qbar\) as in Section \ref{sec:not_red_gps}.
  We have (see Section \ref{sec:not_red_gps} and the beginning of Section \ref{sec:other_half_spin}) two natural \(\Ghat_\mathrm{ad}(\Q)\)-orbits of isomorphisms \(\alpha: \SObf_{8n,\Q} \simeq \Ghat\), swapped by \(\thetahat\), and we need to show that choosing \(\alpha\) suitably yields \(\rho_{x_0} := \alpha \circ \tilde{\rho}_{x_0}\) satisfying the condition of the lemma.
  We will also denote \(\rho = \alpha \circ \tilde{\rho} : \GalQ \to \Ghat(\Ocal(\Cscr))\), again leaving \(\alpha\) implicit.

  We know that for any prime number \(p \neq \ell\) and any \(x \in Z\) the characteristic polynomial of \(\tilde{\rho}_x(\Frob_p)\) (in the standard representation of \(\SObf_{8n}\)) is given by the image via the standard representation of \(\Ghat\) of the Satake parameter corresponding to the restriction of \(\Xi_x\) to \(\Hcal(\Gbf_{\Zp})\).
  Any isomorphism \(\alpha: \SObf_{8n,\Q} \simeq \Ghat\) in one the two \(\Ghat_\mathrm{ad}(\Q)\)-orbits recalled above induces an isomorphism of \(\Q\)-algebras
  \[ \Ocal(\Ghat)^{\Ghat, \thetahat} \xrightarrow[\sim]{\alpha^*} \Ocal(\SObf_{8n,\Q})^{\Obf_{8n,\Q}} \]
  for the conjugation actions of \(\Ghat\) and \(\Obf_{8n,\Q}\), and this isomorphism does not depend on the choice of \(\alpha\).
  Moreover the Satake isomorphism \(\Sat_{\Gbf_{\Zp}}: \Hcal(\Gbf_{\Zp}) \simeq \Ocal(\Ghat)^{\Ghat}\), which is defined over \(\Q\) because the sum of the positive roots for \(\SObf_{8n}\) is divisible by \(2\), induces an isomorphism
  \[ \Hcal(\Gbf_{\Zp})^{\theta^*} \simeq \Ocal(\Ghat)^{\Ghat, \thetahat}. \]
  It is classical that the morphism of \(\Q\)-algebras
  \[ \Ocal(\GLbf_{8n,\Q})^{\GLbf_{8n,\Q}} \longrightarrow \Ocal(\SObf_{8n,\Q})^{\Obf_{8n,\Q}} \]
  induced by the standard representation of \(\Obf_{8n}\) is surjective.
  It follows that the above compatibility at \(x \in Z\) and \(p \neq \ell\) may be reformulated as saying that for any \(f \in \Hcal(\Gbf_{\Zp})^{\theta^*}\) we have
  \[ \Sat_{\Gbf_{\Zp}}(f) (\rho_x(\Frob_p)) = \Xi_x(f). \]
  By Zariski density of \(Z\) in the reduced curve \(\Cscr\) we have, for any \(p \neq \ell\) and any \(f \in \Hcal(\Gbf_{\Zp})^{\theta^*}\), the equality in \(\Ocal(\Cscr)\):
  \begin{equation} \label{eq:compat_Xi_Sat_Oeven}
    \Sat_{\Gbf_{\Zp}}(f) (\rho(\Frob_p)) = \Xi(f).
  \end{equation}
  Specializing at \(x_0\) we see in particular that, for any choice of \(\alpha\) as above, we have \(\rho_{x_0}(\Frob_p)^{\sesi} \in \iota(c(\pi_p))\) or \(\rho_{x_0}(\Frob_p)^{\sesi} \in \iota(\thetahat(c(\pi_p)))\).
  %We are left to show that for a suitable \(\alpha\) the equality TODO holds for any \(p \neq \ell\) and any \(f \in \Hcal(\Gbf_{\Zp})\).

  If for all \(p \neq \ell\) the semi-simple conjugacy class \(c(\pi_p)\) is invariant under \(\thetahat\) then clearly any choice of \(\alpha\) works.
  So we may assume that there exists a prime number \(p \neq \ell\) such that \(c(\pi_p)\) is not invariant under \(\thetahat\), and we fix such a prime number \(p\).
  This determines a choice of isomorphism \(\alpha\) (up to composing with \(\Ghat_\mathrm{ad}(\Q)\)) mapping the conjugacy class of \(\tilde{\rho}_{x_0}(\Frob_p)^{\sesi}\) to \(c(\pi_p)\), and we use this isomorphism to form \(\rho := \alpha \circ \tilde{\rho}\).
  Recall from \cite[\S 19]{FultonHarris_rep} that the algebraic representation \(\bigwedge^{4n} \Q^{8n}\) of \(\SObf_{8n,\Q}\) (exterior power of the standard representation) decomposes as \(r_1 \oplus r_2\), where \(r_1\) and \(r_2\) are absolutely irreducible.
  For \(q\) a prime distinct from \(\ell\) let \(T_{q,1}, T_{q,2} \in \Hcal(\Gbf_{\Z_q})\) be the Hecke operators corresponding, via the Satake isomorphism, to the traces of \(r_1 \alpha^{-1}\) and \(r_2 \alpha^{-1}\).
  As a special case of \eqref{eq:compat_Xi_Sat_Oeven} we have
  \[ \sum_{i=1}^2 \Xi(T_{q,i}) = \sum_{i=1}^2 \Tr \left( r_i \tilde{\rho} (\Frob_q) \right) \text{ and } \prod_{i=1}^2 \Xi(T_{q,i}) = \prod_{i=1}^2 \Tr \left( r_i \tilde{\rho} (\Frob_q) \right) \]
  because \(T_{q,1}+T_{q,2}\) and \(T_{q,1} T_{q,2}\) are both invariant under \(\thetahat\).
  Up to swapping \(r_1\) and \(r_2\), at \(p\) we have \(\Xi(T_{p,i})_{x_0} = \Tr \left( r_i \tilde{\rho}_{x_0} (\Frob_p) \right)\), and \(\Tr \left( r_1 \tilde{\rho}_{x_0} (\Frob_p) \right) \neq \Tr \left( r_2 \tilde{\rho}_{x_0} (\Frob_p) \right)\) by the first point in Lemma \ref{lemm:dist_hspin_tr}.
  Since \(\Cscr\) is integral this implies the equalities \(\Xi(T_{p,i}) = \Tr \left( r_i \tilde{\rho} (\Frob_p) \right)\) in \(\Ocal(\Cscr)\).
  As recalled in the proof of Lemma \ref{lemm:dist_hspin_tr}, the \(\Q\)-algebra \(\Hcal(\Gbf_{\Zp})\) is generated by \(\Hcal(\Gbf_{\Zp})^{\theta^*}\) and \(T_{p,1}\), and so the equality \eqref{eq:compat_Xi_Sat_Oeven} holds for all \(f \in \Hcal(\Gbf_{\Zp})\).
  By continuity, up to removing finitely many points from \(Z\) we can assume that for all \(x \in \Cscr(E')\) taking values in \(Z\) we have \(\Tr \left( r_1 \tilde{\rho}_x (\Frob_p) \right) \neq \Tr \left( r_2 \tilde{\rho}_x (\Frob_p) \right)\).
  This implies that \(\rho_x\) is, up to conjugation by \(\Ghat(\Qellbar)\), the unique morphism \(\GalQ \rightarrow \Ghat(\Qellbar)\) satisfying both \(\Std \circ \rho_x \simeq \sigma_{\pi(x), \iota}^{\Std}\) and \(\Tr \left( r_1 \rho_x (\Frob_p) \right) = \Xi(T_{p,1})_x\).
  By assumption this then holds with \(p\) replace by any prime number distinct from \(\ell\).
  It follows that for any \(x \in Z\), any prime \(q \neq \ell\) and any \(f \in \Hcal(\Gbf_{\Z_q})\) we have \(\Sat_{\Gbf_{\Z_q}}(f)(\rho_x(\Frob_q)) = \Xi_x(f)\).
  By Zariski density of \(Z\) in the reduced curve \(\Cscr\) we deduce the equality \(\Sat_{\Gbf_{\Z_q}}(f)(\rho(\Frob_q)) = \Xi(f)\) in \(\Ocal(\Cscr)\).
  Specializing at \(x_0\) yields the result.
\end{proof}

\begin{prop} \label{pro:loc_glob_SO8n}
  Let \(\pi\) be a level one automorphic representation for \(\Gbf_\mathrm{ad}\), \(\ell\) a prime number and \(\iota : \C \simeq \Qellbar\).
  Then there exists a continuous semisimple morphism \(\rho_{\pi, \iota}^{\SO} : \GalQ \rightarrow \Ghat(\Qellbar)\) which is unramified away from \(\ell\) and such that for any prime number \(p \neq \ell\) we have \(\rho_{\pi, \iota}^{\SO}(\Frob_p)^{\sesi} \in \operatorname{pr} \left( \iota( c_p(\pi_p) ) \right)\), where \(\operatorname{pr}: \Ghat_{\sico} \to \Ghat\) is the natural projection.
  If the infinitesimal character of \(\pi_{\infty}\) is not invariant under \(\thetahat\) then \(\rho_{\pi,\iota}^{\SO}\) is crystalline at \(\ell\).
\end{prop}
\begin{proof}
  First we assume that the infinitesimal character \(\tau\) of \(\pi_{\infty}\) is not invariant under \(\thetahat\), and that the representation \(\rho_{\pi, \iota}^{\mathrm{O}}\) (see Theorem \ref{thm:existence_rho_O}) has maximal infinitesimal image.
  Thanks to the multiplicity formula (see Example \ref{exam:mult_formula_definite}) we know that \(\pi\) corresponds to a parameter \(\psi \in \Psit_{\disc}^{\unr, \taut}(\Gbf)\): for all primes \(p\) we have \(c(\pi_p) = \dpsitausc(\cpsc(\psi))\).
  The composition \(\sigma_{\pi, \iota}^{\Std}\) of \(\rho_{\pi, \iota}^{\mathrm{O}}\) with the standard representation of \(\mathrm{O}_{4n}(\Qellbar)\) is irreducible and so \(\psi\) is non-endoscopic.
  Fix a prime \(p \neq \ell\) such that \(c_p(\psi)\) is not fixed by \(\thetahat\) (see Corollary \ref{cor:Sen_not_thetahat_inv}).
  By the first point of Lemma \ref{lemm:dist_hspin_tr} we have \(\thetahat(\operatorname{pr}(c(\pi_p))) \neq \operatorname{pr}(c(\pi_p))\).
  There exists a unique morphism \(\rho : \GalQ \rightarrow \Ghat(\Qellbar)\) such that \(\Std \circ \rho \simeq \sigma_{\pi, \iota}^{\Std}\) and \(\rho(\Frob_p)^{\sesi} \in \operatorname{pr}(\iota(c(\pi_p)))\), uniquely determined up to conjugation by \(\Ghat(\Qellbar)\).
  We are left to show that for any \(q \not\in \{p, \ell\}\) we also have \(\rho(\Frob_q)^{\sesi} \in \operatorname{pr} (\iota(c_q(\pi_q)))\), i.e.\ that \(\rho\) does not depend on the choice of \(p\) as above.
  By Proposition \ref{pro:ex_lift_cond_one}, there exists a unique geometric lift \(\tilde{\rho} : \GalQ \rightarrow \Ghat_{\mathrm{sc}}(\Qellbar)\) of \(\rho\) of conductor one.
  We can see \(1 \oplus \rho_{\pi, \iota}^{\mathrm{O}}\) as a morphism \(\GalQ \rightarrow \SO_{8n+1}(\Qellbar)\), well-defined up to conjugation.
  The composition of \(\tilde{\rho}\) with the natural embedding \(\Ghat_{\sico}(\Qellbar) \hookrightarrow \Spin_{8n+1}(\Qellbar)\) is the unique geometric lift of \(1 \oplus \rho_{\pi, \iota}^{\mathrm{O}}\) of conductor one.
  Denote by \(r^+, r^-\) the two half-spin representations of \(\Ghat_{\sico}\), distinguished using the infinitesimal character \(\tau\) of \(\pi_{\infty}\) by the condition \(r^{\epsilon} \circ \dpsitausc \simeq \spin_{\psi}^{\epsilon}\).
  The composition of \(\Ghat_{\sico} \hookrightarrow \Spin_{8n+1}\) with the spin representation of \(\Spin_{8n+1}\) is isomorphic to \(r^+ \oplus r^-\).
  Applying Proposition \ref{pro:un_cond_one} as in the proof of Theorem \ref{theo:odd_spin_rep}, we obtain
  \[ (r^+ \circ \tilde{\rho}) \oplus ( r^- \circ \tilde{\rho} ) \simeq \sigma_{\psi, \iota}^{\spin, +} \oplus \sigma_{\psi, \iota}^{\spin, -}. \]
  The representations on the right hand side are those defined in Corollary \ref{coro:ex_both_spin}.
  Since by assumption \(\rho\) has maximal infinitesimal image, both \(r^+ \circ \tilde{\rho}\) and \(r^- \circ \tilde{\rho}\) are irreducible, and so \(r^+ \circ \tilde{\rho}\) is isomorphic to \(\sigma_{\psi, \iota}^{\spin, +}\) or \(\sigma_{\psi, \iota}^{\spin, -}\).
  We will show that the latter is impossible.
  By construction we have \(\tilde{\rho}(\Frob_p)^{\sesi} \in Z(\Ghat_{\mathrm{sc}}) \iota(\dpsitausc(\cpsc(\psi)))\), and so
  \[ \Tr \left( (r^+ \circ \tilde{\rho})(\Frob_p) \right)^2 = \iota \left( \Tr \left( \spin_{\psi}^+ (\cpsc(\psi)) \right) \right)^2 = \Tr \left( \sigma_{\psi, \iota}^{\spin, +}(\Frob_p) \right)^2. \]
  By the second point in Lemma \ref{lemm:dist_hspin_tr} applied to \(\cpsc(\psi)\), this is not equal to
  \[ \iota \left( \Tr \left( \spin_{\psi}^- (\cpsc(\psi)) \right) \right)^2 = \Tr \left( \sigma_{\psi, \iota}^{\spin, -}(\Frob_p) \right)^2 \]
  which implies \(r^+ \circ \tilde{\rho} \not\simeq \sigma_{\psi, \iota}^{\spin, -}\), and so \(r^+ \circ \tilde{\rho} \simeq \sigma_{\psi, \iota}^{\spin, +}\).
  Now consider a prime \(q \not\in \{\ell, p\}\).
  The equality 
  \[ \Tr \left( (r^+ \circ \tilde{\rho})(\Frob_q) \right)^2 = \Tr \left( \sigma_{\psi, \iota}^{\spin, +}(\Frob_q) \right)^2 = \iota \left( \Tr \left( \spin_{\psi}^+ (c_{q, \mathrm{sc}}(\psi)) \right) \right)^2 \]
  implies \(\rho(\Frob_q)^{\sesi} \in \iota( \dpsitau( c_q(\psi))) = \iota (c_q(\pi_q))\), again using Lemma \ref{lemm:dist_hspin_tr}.
  This concludes the proof of the existence of \(\rho_{\pi, \iota}^{\SO}\) under the assumptions that the infinitesimal character of \(\pi_{\infty}\) is not invariant under \(\theta\) and \(\rho_{\pi, \iota}^{\mathrm{O}}\) has maximal infinitesimal image.

  Existence in the general case now follows from Lemma \ref{lem:Gal_rep_SO_simple_family} and Theorem \ref{theo:TaiCC}.

  Now assume that the infinitesimal character of \(\pi_{\infty}\) is not invariant under \(\theta\).
  As above the multiplicity formula (Example \ref{exam:mult_formula_definite}) tells us that \(\pi\) corresponds to a parameter \(\psi \in \tilde{\Psi}_{\disc}^{\unr, \tilde{\tau}}(\Gbf)\), and Theorem \ref{thm:existence_rho_O} implies that \(\rho_{\pi,\iota}^{\SO}\) is crystalline.
\end{proof}

\begin{coro} \label{cor:loc_glob_param_SO8n}
  Let \(n \geq 1\), \(\tilde{\tau} \in \ICcalt(\mathbf{SO}_{8n})\) and \(\psi \in \tilde{\Psi}_{\disc, \nonendo}^{\unr, \tilde{\tau}}(\mathbf{SO}_{8n})\).
  Theorem \ref{thm:loc_glob_param_SO_even} holds for \(\psi\) (for any \(\iota: \C \simeq \Qellbar\)).
\end{coro}
\begin{proof}
  Recall from Definition \ref{def:Mpsi} that we may take \(\Mpsi = \Ghat\).
  Thanks to the multiplicity formula (Example \ref{exam:mult_formula_definite}) there exists a (unique) level one automorphic representation \(\pi\) for \(\Gbf_\mathrm{ad}\) corresponding to \(\psi\), i.e.\ such that \(\pi_\infty\) has infinitesimal character \(\tau_\psi\) and for any prime number \(p\) we have \(c(\pi_p) = \cpsc(\psi)\).
  The corollary thus follows from Proposition \ref{pro:loc_glob_SO8n}.
\end{proof}

\subsection{\(\SObf_{8n-4}\): using endoscopy}
\label{sec:using_endo}

We would like to extend Corollary \ref{cor:loc_glob_param_SO8n} to the case of \(\mathbf{SO}_{8n-4}\), for any \(n \geq 1\).
Unfortunately we do not know how to \(\ell\)-adically deform an arbitrary element of \(\tilde{\Psi}_{\disc, \nonendo}^{\unr, \tilde{\tau}}(\mathbf{SO}_{8n-4})\) as in Theorem \ref{theo:TaiCC}.
We circumvent this problem by using endoscopic parameters for \(\mathbf{SO}_{8n}\).
This is only possible in ``almost all'' cases.

For \(f \in S_{2k}(\SLbf_2(\Z))\) an eigenform, denote by \(a_p(f)\) the corresponding (real) eigenvalue for the Hecke operator \(T_p\), and let \(a_p'(f) = a_p(f) / p^{k-1/2}\).
In particular the corresponding level one cuspidal automorphic representation \(\pi\) for \(\PGLbf_2\) (as in Section \ref{sec:SO4}) is such that the semisimple conjugacy class \(c(\pi_p)\) in \(\SL_2(\C)\) has characteristic polynomial \(X^2 - a_p'(f) X + 1\).
The Ramanujan conjecture (proved by Deligne) implies \(|a_p'(f)| \leq 2\).

\begin{lemm} \label{lem:8nm4_to_8n}
  Let \(n \geq 1\), \(\tilde{\tau}_1 \in \tilde{\ICcal}(\mathbf{SO}_{8n-4})\), and \(\psi_1 \in \tilde{\Psi}_{\disc, \nonendo}^{\unr, \tilde{\tau}_1}(\mathbf{SO}_{8n-4})\).
  We consider a definite inner form \(\Gbf\) of \(\SObf_{8n}\) as in Section \ref{sec:other_half_spin}.
  Assume that \(\psi_1 = \pi_1[d_1]\) with \(d_1\) odd, and that \(\tilde{\tau}_1\) is not bad (Definition \ref{def:bad_tau}).
  Let \(S\) be a finite set of prime numbers.
  Then there exists \(\tilde{\tau}_2 \in \tilde{\ICcal}(\SObf_4)\) and \(\psi_2 \in \tilde{\Psi}_{\disc, \nonendo}^{\unr, \tilde{\tau}_2}(\SObf_4)\) of shape \(\pi_2[1]\) satisfying
  \begin{enumerate}
  \item for all primes in \(S\) the conjugacy class \(c_p(\psi_2)\) is not fixed by \(\thetahat\),
  \item we have \(\tilde{\tau}_1 \oplus \tilde{\tau}_2 \in \tilde{\ICcal}(\SObf_{8n})\) and there exists an automorphic level one representation \(\pi\) for \(\Gbf_\mathrm{ad}\) with corresponding parameter \(\psi_1 \oplus \psi_2\) (for the multiplicity formula in Example \ref{exam:mult_formula_definite}).
  \end{enumerate}
\end{lemm}
\begin{proof}
  We may assume that \(\tilde{\tau}_1\) is the class of \((w_1, \dots, w_{4n-2})\) where \(w_1 > \dots > w_{4n-2}\) are integers.
  Since \(\tilde{\tau}_1\) is not bad, there exists \(1 \leq i \leq 2n-1\) satisfying \(w_{2i-1} > w_{2i}+1\).
  Let \(w'_2 = w_{2i}+1\).
  We are looking for level one eigenforms \(f_1, f_2\) of respective weights \(2(k+w'_2), 2k\) such that the associated (by Proposition \ref{pro:loc_glob_SO4}) \(\psi_2 \in \tilde{\Psi}_{\disc, \nonendo}^{\unr, \tilde{\tau}_2}(\SObf_4)\), for \(\tilde{\tau}_2\) the class of \((2k+w_2'-1, w_2')\), satisfies the first condition in the lemma for all \(p \in S\).
  In terms of the Hecke eigenvalues, this condition can be restated as \(a_p'(f_1)^2 \neq a_p'(f_2)^2\).
  Let \(N = |S|\).
  Using \cite[Corollaire 1, p.\ 80]{SerreEquidist} we see that there exists a family \((I_j)_{1 \leq j \leq N+1}\) of closed intervals included in \([0,4]\) and covering \([0,4]\), and \(k_0 \geq 8\) such that for all \(k \geq k_0\), \(p \in S\) and \(1 \leq j \leq N+1\), the proportion of eigenforms \(f \in S_{2k}(\SLbf_2(\Z))\) satisfying \(a_p'(f)^2 \in I_j\) is less than \(1/N\).
  Fix \(k \geq k_0\) satisfying \(2k+w_2'-1>w_1\), and an eigenform \(f_2 \in S_{2k}(\SLbf_2(\Z))\).
  For any \(p \in S\), there exists \(j\) such that \(a_p'(f_2)^2 \in I_j\) and so the proportion of eigenforms \(f \in S_{2(k+w_2')}(\SLbf_2(\Z))\) satisfying \(a_p'(f)^2 = a_p'(f_2)^2\) is less than \(1/N\).
  Therefore the proportion of \(f\) satisfying \(a_p'(f)^2 \neq a_p'(f_2)^2\) for all \(p \in S\) is positive, and we let \(f_1\) be such an eigenform and \(\psi_2\) be the parameter corresponding to \((f_1, f_2)\) by Proposition \ref{pro:loc_glob_SO4}.

  We are left to check that the level one representation of \(\Gbf_\mathrm{ad}(\A)\) associated to \(\psi = \psi_1 \oplus \psi_2\) and any choice of \(\tau \in \ICcal(\Gbf)\) mapping to \(\taut := \tilde{\tau}_1 \oplus \tilde{\tau}_2\) is automorphic, i.e.\ that we have \(\langle \cdot, \pi_{\infty} \rangle|_{\Scal_{\psi}} = \epsilon_{\psi}\), by Example \ref{exam:mult_formula_definite}.
  Let \(v_1 > \dots > v_m > -v_m > \dots > -v_1\) be the eigenvalues of the infinitesimal character of \((\pi_1)_{\infty}\), so that
  \[ \left( w_1, \dots, w_{4n-2} \right) = \left( v_1 + \frac{d_1-1}{2}, \dots, v_1 + \frac{1-d_1}{2}, \dots, v_m + \frac{d_1-1}{2}, \dots, v_m + \frac{1-d_1}{2} \right). \]
  We have \(w_{2i} = v_j + (d_1-1)/2\) for some \(1 < j \leq m\), and we have \(2i = d_1(j-1)+1\) with \(d_1\) odd so \(j\) is even.
  So the indices corresponding to \(2k+w_2'-1\) and \(w_2'\) in \(\taut = (w_1'' > \dots > w_{4n}'')\) are \(1\) and \(2+d_1(j-1)\), which are both odd.
  It follows from the recipe in Example \ref{exam:mult_formula_definite} that the character \(\langle \cdot, \pi_{\infty} \rangle\) of the group (with two elements) \(\Scal_\psi\) is trivial.
  By Lemma \ref{lem:epsilon_psi_s_psi}, more precisely the second formula in the last point, we have \(\epsilon_{\psi}(s) = 1\) because \(d_1\) and \(1\) are odd.
\end{proof}

\begin{prop} \label{pro:loc_glob_SO8nm4}
  Let \(n \geq 2\), \(\taut \in \ICcalt(\mathbf{SO}_{8n-4})\) and \(\psi \in \tilde{\Psi}_{\disc, \nonendo}^{\unr, \tilde{\tau}}(\mathbf{SO}_{8n})\).
  Assume that \(\taut\) is not bad.
  Then Theorem \ref{thm:loc_glob_param_SO_even} holds for \(\psi\) (for any \(\iota: \C \simeq \Qellbar\)).
\end{prop}
\begin{proof}
  Fix a prime \(p \neq \ell\) such that \(c_p(\psi)\) is not fixed by \(\thetahat\) (see Corollary \ref{cor:Sen_not_thetahat_inv}).
  By Corollary \ref{cor:two_SO_cc_rho_O} there is up to conjugation by \(\Mpsi(\Qellbar)\) a unique morphism \(\rho: \GalQ \to \Mpsi(\Qellbar)\) such that \(\alpha \circ \rho\) is conjugated by \(\mathrm{O}_{8n-4}(\Qellbar)\) to the morphism \(\rho_{\psi,\iota}^\mathrm{O}\) of Theorem \ref{thm:existence_rho_O} and \(\rho(\Frob_p)^{\sesi}\) is conjugated to \(\iota(c_p(\psi))\).
  We are left to show that for any \(q \neq \ell\) the conjugacy class of \(\rho(\Frob_q)^{\sesi}\) is also equal to \(c_q(\psi)\).
  Fix \(q \neq \ell\).
  Apply Lemma \ref{lem:8nm4_to_8n} to \(\psi\) and \(S=\{p,q\}\) to obtain \(\tau' \in \ICcal(\SObf_{8n})\), a parameter \(\psi' = \psi \oplus \psi_2 \in \Psit_{\disc}^{\unr, \tilde{\tau'}}(\SO_{8n})\) and a level automorphic representation \(\pi\) for \(\Gbf_\mathrm{ad}\) such that \(\pi_\infty\) has infinitesimal character \(\tau'\) and for any prime \(q\) we have \(c(\pi_q) = \dot{\psi}'_{\tau', \mathrm{sc}}(c_{q,\mathrm{sc}}(\psi'))\).
  By Proposition \ref{pro:loc_glob_SO8n} there exists a continuous semisimple morphism \(\rho': \GalQ \to \SO_{8n}(\Qellbar)\) unramified away from \(\ell\) and such that for any prime \(r \neq \ell\) we have \(\rho'(\Frob_r) \in \dot{\psi'}_{\tau'}(c_r(\psi'))\).
  So \(\rho'\) belongs to one of the two \(\SO_{8n}(\Qellbar)\)-conjugacy classes making up the \(\mathrm{O}_{4n}(\Qellbar)\)-conjugacy class of \(\rho_{\psi',\iota}^\mathrm{O}\) of Theorem \ref{thm:existence_rho_O}.
  We claim that \(\rho'\) is conjugated to \(\dot{\psi'}_{\tau'} \circ (\rho, \rho_{\psi_2,\iota}^{\SO})\), where \(\rho_{\psi_2,\iota}^{\SO}\) is defined by Proposition \ref{pro:loc_glob_SO4}.
  By Corollary \ref{cor:two_SO_cc_rho_O} it is enough to check that \(\rho'(\Frob_p)^{\sesi}\) and \(\dot{\psi'}_{\tau'}(\rho(\Frob_p), \rho_{\psi_2,\iota}^{\SO}(\Frob_p))^{\sesi}\) are conjugated in \(\SO_{8n}(\Qellbar)\) and that this conjugacy class is not invariant under \(\thetahat\).
  The first property follows from the definition of \(\rho\), while the second follows from the fact that neither \(c_p(\psi)\) (by choice of \(p\)) nor \(c_p(\psi_2)\) (because \(p \in S\)) is invariant under \(\thetahat\).
  We obtain that \(\iota \dot{\psi'}_{\tau'}(c_q(\psi), c_q(\psi_2))\) and \(\dot{\psi'}_{\tau'}(\rho(\Frob_q)^{\sesi}, \iota(c_q(\psi_2)))\) are conjugated in \(\SO_{8n}(\Qellbar)\), and since \(c_q(\psi_2)\) is not invariant under \(\thetahat\) (because \(q \in S\)) this implies that \(\iota(c_q(\psi_2))\) and \(\rho(\Frob_q)^{\sesi}\) are conjugated in \(\Mpsi(\Qellbar)\).
\end{proof}

\subsection{Non-tempered parameters \(\pi[2d+1]\)}

Similarly to Section \ref{sec:odd_spin_nontemp}, in the case of parameters \(\psi \in \Psit_{\disc,\nonendo}^{\unr,\taut}(\SObf_{4n})\) of the form \(\pi[2d+1]\) with \(d>0\) (other non-tempered parameters were treated in Section \ref{sec:even_spin_pi_2d}) we would like to relate the lifted Satake parameters \((\cpsc(\psi))_p\) of \(\psi\) to those of \(\pi[1] \in \Psi_{\disc,\nonendo}^{\unr,\tau'}(\SObf_{4m})\) where \(n=m(2d+1)\).
Similarly to Definition \ref{def:alpha_m_d} we have a morphism \(\beta: \Mcal_{\pi[1]} \times \SL_2 \to \Mpsi\) such that \(\Std_{\Mpsi} \circ \beta \simeq \Std_{\Mcal_{\pi[1]}} \otimes \Sym^{2d} \Std_{\SL_2}\) and whose differential maps \((\tau_{\pi[1]}, 1/2)\) to \(\tau_\psi\), unique up to conjugation by \(\Mpsi(\Qbar)\).
We may take \(\Mcal_{\pi[1]} = \SO_{4m}\) with \(\tau_{\pi[1]} = (w_1, \dots, w_{2m})\) with \(w_1 > \dots > w_{2m} > 0 \) integers and \(\Mpsi = \SO_{4n}\) with \(\tau_\psi = (\tau_1, \dots, \tau_{2n})\) with \(\tau_1 > \dots > \tau_{2n} > 0\) integers.
Then up to conjugacy we may assume that \(\beta\) restricts to
\begin{align*}
  \Tcal_{\SO_{4m}} \times \Tcal_{\SL_2} & \longrightarrow \Tcal_{\SO_{4n}} \\
  ((x_1, \dots, x_{2m}), t) & \longmapsto (x_1 t^{2d}, x_1 t^{2d-2}, \dots, x_1 t^{-2d}, \dots, x_{2m} t^{2d}, \dots, x_{2m} t^{-2d}).
\end{align*}
It lifts uniquely to \(\tilde{\beta}: \Mcal_{\pi[1], \mathrm{sc}} \times \SL_2 \to \Mpsisc\), which restricts to
\begin{align*}
  \Tcal_{\Spin_{4m}} \times \Tcal_{\SL_2} & \longrightarrow \Tcal_{\Spin_{4n}} \\
  ((x_1, \dots, x_{2m}, s), t) & \longmapsto (\beta((x_1, \dots, x_{2m}), t), s^{2d+1}).
\end{align*}
In particular \(\tilde{\beta}\) is trivial on the center of \(\SL_2\) and so it induces \(\tilde{\beta}: \Mcal_{\pi[1], \mathrm{sc}} \times \PGL_2 \to \Mpsisc\).
The morphism
\begin{align*}
  \GL_1 \times \Mcal_{\pi[1], \mathrm{sc}} \times \PGL_2 & \longrightarrow \GMpsisc \\
  (\lambda, h, x) & \longmapsto \lambda^{2d+1} \tilde{\beta}(h,x)
\end{align*}
induces an extension \(\mathcal{GM}_{\pi[1], \mathrm{sc}} \times \PGL_2 \to \GMpsisc\) of \(\tilde{\beta}\), that we still denote by \(\tilde{\beta}\).

\begin{prop} \label{pro:rho_even_GSpin_pi_2dp1}
  Assume that \(n\) is even or that \(\taut\) is not bad.
  Exactly one of the following holds true.
  \begin{itemize}
  \item
    \begin{enumerate}
    \item For any prime \(\ell\) and any \(\iota: \C \simeq \Qellbar\) the \(\Mpsi(\Qellbar)\)-conjugacy classes of \(\rho_{\psi,\iota}^{\GSpin}\) and \(\tilde{\beta} \circ (\rho_{\pi[1],\iota}^{\GSpin}, 1 \oplus \chi_\ell^{-1})\) are equal.
    \item For all primes \(p\) we have \(\cpsc(\psi) = \tilde{\beta}(\cpsc(\pi[1]), \diag(p^{1/2}, p^{-1/2}))\).
    \end{enumerate}
  \item 
    \begin{enumerate}
    \item For any prime \(\ell\) and any \(\iota: \C \simeq \Qellbar\) the \(\Mpsi(\Qellbar)\)-conjugacy classes of \(\thetahat \circ \rho_{\psi,\iota}^{\GSpin}\) and \(\tilde{\beta} \circ (\rho_{\pi[1],\iota}^{\GSpin}, 1 \oplus \chi_\ell^{-1})\) are equal.
    \item For all primes \(p\) we have \(\thetahat(\cpsc(\psi)) = \tilde{\beta}(\cpsc(\pi[1]), \diag(p^{1/2}, p^{-1/2}))\).
    \end{enumerate}
  \end{itemize}
\end{prop}
\begin{proof}
  The proof is similar to that of Proposition \ref{pro:rho_odd_GSpin_pi_2dp1}.
  We already know that for all primes \(p\) the conjugacy class of \(\beta(c_p(\pi[1]), \diag(p^{1/2}, p^{-1/2}))\) is equal either to \(c_p(\psi)\) or to \(\thetahat(c_p(\psi))\).
  Note that \(m\) is even or \(\taut'\) is not bad.
  By Theorem \ref{thm:existence_rho_O} the representation \(\rho_{\psi,\iota}^{\SO}\) obtained using Theorem \ref{thm:loc_glob_param_SO_even} is conjugated under \(\Mpsi(\Qellbar)\) to either \(\beta \circ (\rho_{\pi[1],\iota}^{\SO}, 1 \oplus \chi_\ell^{-1})\) or \(\thetahat \beta \circ (\rho_{\pi[1],\iota}^{\SO}, 1 \oplus \chi_\ell^{-1})\) (\(\rho_{\pi[1],\iota}^{\SO}\) is also obtained using Theorem \ref{thm:loc_glob_param_SO_even}), the two being exclusive by Corollary \ref{cor:two_SO_cc_rho_O}.
  Distinguishing these two cases, the rest of the proof is identical to that of \ref{pro:rho_odd_GSpin_pi_2dp1}.
\end{proof}

\begin{rema} \label{rem:even_pi_2dp1_lift}
  Similarly to Remarks \ref{rem:pi_2d_lift} and \ref{rem:odd_pi_2dp1_lift}, this could perhaps be proved using Eisenstein series, yielding a stronger result that does not require \(\taut\) to not be bad and without ambiguity under \(\thetahat\).
\end{rema}

\section{Applications}

\subsection{Tensor product decomposition in intersection cohomology}

\begin{defi}
  Let \(\ell\) be a prime number and \(\iota: \C \simeq \Qellbar\).
  Let \(n \geq 1\), \(\taut \in \ICcalt(\SObf_{4n})\) and \(\psi \in \tilde{\Psi}_{\disc, \nonendo}^{\unr, \taut}(\SObf_{4n})\).
  Assume either \(n=1\), \(n\) even, \(\psi = \pi[2d]\) or that \(\taut\) is not bad (Definition \ref{def:bad_tau}).
  For \(\epsilon \in \{+,-\}\) let \(\sigma_{\psi,\iota}^{\spin,\epsilon}\) be the continuous semisimple representation of \(\GalQ\) over \(\Qellbar\) of dimension \(2^{2n-1}\) which is obtained by composing the morphism \(\rho_{\psi, \iota}^{\GSpin}: \GalQ \to \GMpsisc(\Qellbar)\) of Theorem \ref{theo:even_GSpin_Gal_rep} with the representation \(\spin_\psi^\epsilon\) (Definition \ref{def:GMpsisc}).
\end{defi}

\begin{rema} \label{rem:missing_sigma_spin_epsilon}
  This extends the definition in the second part of Corollary \ref{coro:seed_ex_sigma}, which included only one possible value for the sign \(\epsilon\) for a given parameter \(\psi\).
  Thus \(\sigma_{\psi,\iota}^{\spin,\epsilon}\) is defined in all cases except when the following conditions are simultaneously satisfied: \(n>1\) is odd, \(\taut\) is bad, \(\psi=\pi[d]\) with \(d\) odd and \(\epsilon = +1\).
  If these conditions are satisfied we will say that the pair \((\psi, \epsilon)\) is unreachable.
\end{rema}

\begin{theo} \label{thm:sigma_is_tensor}
  In the setting of Theorem \ref{theo:IH_explicit_crude}, for any \(\psi = \psi_0 \oplus \dots \oplus \psi_r\) we have
  \begin{enumerate}
  \item for any \(1 \leq i \leq r\) the pair \((\psi_i, u_i(\psi))\) is not unreachable,
  \item a decomposition
    \[ \sigma_{\psi,\iota}^{\IH} \simeq \sigma_{\psi_0,\iota}^{\spin} \otimes \sigma_{\psi_1,\iota}^{\spin,u_1(\psi)} \otimes \dots \otimes \sigma_{\psi,\iota}^{\spin,u_r(\psi)} \left( n_0 d_0 (n_0 d_0 + 1)/4 + \sum_{i=1}^r n_i d_i / 8 - n(n+1)/4 \right) \]
    where the expression in the Tate twist is an integer.
  \end{enumerate}
\end{theo}
\begin{proof}
  The second point immediately follows from the first point, the \v{C}ebotarev density theorem, \eqref{eq:explicitIH_cc_Frob}, Theorem \ref{theo:odd_spin_rep}, Theorem \ref{theo:even_GSpin_Gal_rep} and the second part of Corollary \ref{coro:seed_ex_sigma}.

  Let us prove the first point.
  For \(1 \leq i \leq r\) it follows from the definition \eqref{eq:def_ui_psi} of \(u_i(\psi)\) that if \(\psi_i \in \Psit_{\disc,\nonendo}^{\unr, \taut_i}(\SObf_{n_i d_i})\) with \(d_i\) odd, \(n_i \equiv 4 \mod 8\) and \(\taut_i\) bad then \(u_i(\psi) = -\epsilon_\psi(s_i)\).
  We can compute \(\epsilon_\psi(s_i)\) using Lemma \ref{lem:epsilon_psi_s_psi}:
  \[ \epsilon_\psi(s_i) = \prod_j \epsilon(1/2, \pi_i \times \pi_j) \]
  where the product is over indices \(1 \leq j \leq r\) such that \(d_j\) is even and \(d_j>d_i\).
  We have \(\pi_i \in O_o(w_1, \dots, w_{n_i/2})\) where the integers \(w_1 > \dots > w_{n_i/2} > 0\) satisfy \(w_{2k-1} = w_{2k} + d_i\) for all \(1 \leq k \leq n_i/4\).
  For \(j\) as above we have \(\pi_j \in S(w_1', \dots, w_{n_j/2}')\) where \(w_k' \in 1/2+\Z\) satisfy \(w_1' > \dots > w_{n_j/2}' > 0\).
  We have (see \cite[\S 3.9]{ChRe})
  \[ \epsilon(1/2, \pi_i \times \pi_j) = \prod_{\substack{1 \leq a \leq n_i/2 \\ 1 \leq b \leq n_j/2}} (-1)^{1+2\max(w_a,w_b')} \]
  and since for any \(1 \leq k \leq n_i/4\) and \(1 \leq b \leq n_j/2\) we have \(w_{2k-1}>w_{2k}\) and either \(w_b'<w_{2k}\) or \(w_b'>w_{2k-1}\) by badness of \(\taut_i\) and disjointness of \(\taut_i\) and \(\taut_j\), we obtain \(\epsilon(1/2, \pi_i \times \pi_j) = 1\) and \(u_i(\psi) = -1\).
  Thus \((\psi_i, u_i(\psi))\) is not unreachable.
\end{proof}

\subsection{Siegel modular forms}
\label{sec:SMF}

Let \(n \geq 1\).
We recall a few facts about the translation between (level one) genus \(n\) vector-valued Siegel modular forms and automorphic representations for \(\PGSpbf_{2n}\) from \cite[\S 5.2]{Taibi_dimtrace}.
For a tuple \(\ul{k} = (k_1 \geq \dots \geq k_n)\) of integers, that we interpret as a highest weight for the complex Lie group \(\GL_n(\C)\), we have a finite-dimensional complex vector space \(S_{\ul{k}}(\Spbf_{2n}(\Z))\) of vector-valued Siegel cusp forms which has an action of the Hecke algebra \(\Hcal_f^{\unr}(\GSpbf_{2n})\).
There are several competing normalizations for this action (see the two normalizations in \cite[Definition 8]{vdG_123}), and we find it convenient to use yet another normalization, the unitary normalization: add a factor \(\eta(\gamma)^{\sum_i k_i/2}\) in \cite[Definition 8]{vdG_123}.
In level one the action of \(\Hcal_f^{\unr}(\GSpbf_{2n})\) then factors through \(\Hcal_f^{\unr}(\PGSpbf_{2n})\), in particular to an eigenform \(f\) is associated a family \((c_p^\mathrm{unit}(f))_p\) of Satake parameters which are semi-simple conjugacy classes in \(\Spin_{2n+1}(\C)\).
The relation with the notation introduced in Section \ref{sec:intro_SMF} is \(c_p^\mathrm{arith}(f) = p^{\sum_i k_i/2 - n(n+1)/4} c_p^\mathrm{unit}(f)\).
Let \(\gfrak = \C \otimes_{\R} \Lie \PGSpbf_{2n}(\R)\) and \(K\) a maximal compact subgroup of \(\PGSpbf_{2n}(\R)\).
By a celebrated theorem of Gelfand, Graev and Piatetski-Shapiro the space of cuspidal automorphic forms for \(\PGSpbf_{2n}\) decomposes discretely, in particular we have a \((\gfrak,K) \times \Hcal_f^{\unr}(\PGSpbf_{2n})\)-equivariant isomorphism
\[ \Acal_\mathrm{cusp}(\PGSpbf_{2n}(\Q) \backslash \PGSpbf_{2n}(\A))^{\PGSpbf_{2n}(\Zhat)} \simeq \bigoplus_{\pi \simeq \pi_\infty \otimes \pi_f} \left( \pi_\infty \otimes \pi_f^{\PGSpbf_{2n}(\Zhat)} \right)^{\oplus m(\pi)} \]
where the sum is over isomorphism classes of irreducible unitary \((\gfrak,K)\)-modules \(\pi_\infty\) (resp.\ admissible representations \(\pi_f\) of \(\PGSpbf_{2n}(\A_f)\)) and \(m(\pi)\) are integers (only countably many of them are non-zero).
Of course we may restrict to irreducible admissible representations \(\pi_f\) of \(\PGSpbf_{2n}(\A_f)\) which are everywhere unramified, i.e.\ such that the space \(\pi_f^{\PGSpbf_{2n}(\Zhat)}\) is non-zero, in which case it has dimension one.
This gives a \(\Hcal_f^{\unr}(\PGSpbf_{2n})\)-equivariant isomorphism
\[ S_{\ul{k}}(\Spbf_{2n}(\Z)) \simeq \bigoplus_{\pi_f} \left( \pi_f^{\PGSpbf_{2n}(\Zhat)} \right)^{m(\pi_\infty(\ul{k}) \otimes \pi_f)} \]
where \(\pi_\infty(\ul{k})\) is a certain explicit irreducible \((\gfrak,K)\)-module associated to \(\ul{k}\) (essentially a generalized Verma module).
It may fail to be unitarizable, in which case by convention we set \(m(\pi_\infty(\ul{k}) \otimes \pi_f)=0\) for all \(\pi_f\).
For \(k_n \geq n+1\) the irreducible \((\gfrak,K)\)-module \(\pi_\infty(\ul{k})\) is unitarizable, in fact it is the holomorphic discrete series \((\gfrak,K)\)-module having infinitesimal character \((k_1-1, \dots, k_n-n)\).
Under this condition \(m(\pi_\infty(\ul{k}) \otimes \pi_f)\) is also the multiplicity of \(\pi_\infty(\ul{k}) \otimes \pi_f\) in the larger space \(\Acal^2(\PGSpbf_{2n}(\Q) \backslash \PGSpbf_{2n}(\A))\) of square-integrable automorphic forms.

\begin{theo} \label{thm:Siegel_formula}
  Let \(n \geq 1\) and \(k_1 \geq \dots \geq k_n \geq n+1\) be integers.
  Denote \(\tau = (k_1-1, \dots, k_n-n) \in \ICcal(\Spbf_{2n})\).
  We have a decomposition in eigenspaces (lines) under \(\Hcal_f^{\unr}(\PGSpbf_{2n})\)
  \[ S_{\ul{k}}(\Spbf_{2n}(\Z)) \simeq \bigoplus_{\psi} \chi_{f,\psi}\]
  where the sum is over \(\psi = \psi_0 \oplus \dots \oplus \psi_r \in \Psi_{\disc}^{\unr,\tau}(\Spbf_{2n})\) satisfying
  \begin{itemize}
  \item the odd-dimensional factor \(\psi_0\) is tempered, i.e.\ of the form \(\pi_0[1]\),
  \item the holomorphic discrete series \((\gfrak, K)\)-module \(\pi_\infty(\ul{k})\) satisfies \(\langle \cdot, \pi_\infty(\ul{k}) \rangle = \epsilon_\psi\),
  \end{itemize}
  and the character \(\chi_{f,\psi}\) of \(\Hcal_f^{\unr}(\PGSpbf_{2n})\) was introduced in Theorem \ref{theo:IH_explicit_crude}.
\end{theo}
\begin{proof}
  This decomposition is already known for the action of the (smaller) Hecke algebra \(\Hcal_f^{\unr}(\Spbf_{2n})\) thanks to Arthur's endoscopic classification \cite{Arthur} and the comparison \cite{AMR} of Arthur packets for \(\Spbf_{2n,\R}\) with the more explicit Adams-Johnson packets when the infinitesimal character is regular algebraic, see \cite[\S 5.2]{Taibi_dimtrace} or \cite[\S 9]{ChRe}.
  More precisely this analysis shows that the eigenspace in \(\Acal^2(\PGSpbf_{2n}(\Q) \backslash \PGSpbf_{2n}(\A))_\tau^{\PGSpbf_{2n}(\Zhat)}\) for the character of \(\Hcal_f^{\unr}(\Spbf_{2n})\) corresponding to \((c_p(\psi))_p\)
  \begin{itemize}
  \item either does not contain the \((\gfrak,K^0)\)-module \(\pi_\infty(\ul{k})\) if any one of the two conditions in Theorem \ref{thm:Siegel_formula} is not satisfied, or
  \item is isomorphic to \(\pi_\infty(\ul{k}) \oplus \bigoplus_{\pi_\infty'} \pi_\infty'\) otherwise, where the last sum ranges over a finite set of irreducible \((\gfrak,K)\)-modules \(\pi_\infty'\) in which the trace of a pseudo-coefficient for \(\pi_\infty(\ul{k})\) vanishes (see Proposition \ref{pro:tr_AJ_ps}).
  \end{itemize}
  Thus Corollary \ref{cor:weak_mult_formula} completely describes the \(\Hcal_f^{\unr}(\PGSpbf_{2n})\)-module
  \[ \Hom_{\gfrak,K}(\pi_\infty(\ul{k}), \Acal^2(\PGSpbf_{2n}(\Q) \backslash \PGSpbf_{2n}(\A))^{\PGSpbf_{2n}(\Zhat)}). \]
\end{proof}

\begin{coro} \label{cor:GSpin_Gal_rep_Siegel}
  Let \(n \geq 1\) and \(k_1 \geq \dots \geq k_n \geq n+1\) be integers.
  Let \(f \in S_{\ul{k}}(\Spbf_{2n}(\Z))\) be an eigenform for the action of \(\Hcal_f^{\unr}(\PGSpbf_{2n})\).
  Let \(\ell\) be a prime number and \(\iota: \C \simeq \Qellbar\).
  Then there exists a continuous semisimple morphism \(\rho_{f,\iota}^{\GSpin}: \GalQ \to \GSpin_{2n+1}(\Qellbar)\) which is crystalline at \(\ell\), unramified away from \(\ell\) and such that for any prime \(p \neq \ell\) the conjugacy class of \(\rho_{f,\iota}^{\GSpin}(\Frob_p)^{\sesi}\) is equal to \(\iota(c_p^\mathrm{arith}(f))\).
  Any other continuous semisimple morphism \(\rho: \GalQ \to \GSpin_{2n+1}(\Qellbar)\) satisfying this property at almost all primes \(p\) is conjugated to \(\rho_{f,\iota}^{\GSpin}\).
\end{coro}
\begin{proof}
  Denoting \(\tau = (k_1-1, \dots, k_n-n) \in \ICcal(\Spbf_{2n})\), the eigenspace \(\C f\) corresponds to a unique
  \[ \psi = \psi_0 \oplus \dots \oplus \psi_r \in \Psi_{\disc}^{\unr,\tau}(\Spbf_{2n}) \]
  satisfying the two conditions in Theorem \ref{thm:Siegel_formula}.
  We claim that for \(1 \leq i \leq r\) the parameter \(\psi_i \in \Psit_{\disc,\nonendo}^{\unr,\taut_i}(\SObf_{n_i d_i})\) satisfies the assumption in Theorem \ref{thm:loc_glob_param_SO_even}.
  More precisely we show that the condition \(\langle \cdot, \pi_\infty(\ul{k}) \rangle|_{\Scal_\psi} = \epsilon_\psi\) implies that \(\taut_i\) is not bad.
  As in Section \ref{sec:desc_IH_lift} up to conjugacy we may assume that \(\dpsitau \circ \varphi_{\psi_\infty} |_{\C^\times}\) takes values in \(\Tcal_{\SO_{2n+1}}(\C)\) and is dominant for \(\Bcal_{\SO_{2n+1}}\), i.e.\ the holomorphic part of \(\dpsitau \circ \varphi_{\psi_\infty} |_{\C^\times}: \C^\times \to \Tcal_{\SO_{2n+1}}(\C)\) is \(z \mapsto (z^{k_1-1}, \dots, z^{k_n-n})\).
  Now the character \(\langle \cdot, \pi_\infty(\ul{k}) \rangle\) is known to be (see \cite[Lemma 9.1]{ChRe} and \cite[\S 5.1]{Taibi_dimtrace})
  \begin{align*}
    \{ (z_1, \dots, z_n, s) \in \Tcal_{\Spin_{2n+1}} \,|\, z_i^2=1,\,s^2=z_1 \dots z_n=1 \} & \longrightarrow \{\pm 1\} \\
    (z_1, \dots, z_n, s) & \longmapsto \prod_{j=1}^{\lfloor n/2 \rfloor} z_{2i}.
  \end{align*}
  Note the similarity\footnote{This is no coincidence since the two are related by the theta correspondence.} with the character occurring in Example \ref{exam:mult_formula_definite} for the case of definite special orthogonal groups.
  The claim follows from essentially the same computations as in the proof of Theorem \ref{thm:sigma_is_tensor}.

  Thus we have \(\GSpin\)-valued Galois representations for all constituents \(\psi_i\) (Theorem \ref{theo:odd_spin_rep} for \(i=0\), Theorem \ref{theo:even_GSpin_Gal_rep} for \(1 \leq i \leq r\)).
  As usual (see Definition \ref{def:alpha_m_d}) we extend \(\dpsitausc: \prod_i \mathcal{M}_{\psi_i,\sico} \to \Spin_{2n+1}\) to obtain \(\dpsitausc: \prod_i \mathcal{GM}_{\psi_i,\sico} \to \GSpin_{2n+1}\) induced by
  \begin{align*}
    \prod_i \GL_1 \times \mathcal{M}_{\psi_i,\sico} & \longrightarrow \GSpin_{2n+1} \\
    (z_i, h_i)_i & \longmapsto \prod_i z_i \times \dpsitausc((h_i)_i).
  \end{align*}
  Define \(\rho_{f,\iota}^{\GSpin}\) as \(\dpsitausc \circ (\rho_{\psi_i,\iota}^{\GSpin})_i\) twisted by \(\chi_\ell^N\) where
  \[ N = \frac{n_0 d_0 (n_0 d_0 + 1)}{4} + \sum_{i=1}^r \frac{n_i d_i}{2} - \frac{\sum_{j=1}^n k_j}{2} + \frac{n(n+1)}{4} \]
  is an integer (as proved in Theorem \ref{thm:sigma_is_tensor}, because \(\sum_j k_j\) is even).
  This clearly satisfies the conditions of the theorem.
  Uniqueness is proved exactly as for Theorem \ref{theo:odd_spin_rep}.
\end{proof}

\section{Explicit formulas for compactly supported Euler characteristics}
\label{sec:formula_e_c}

The first goal of this section is to deduce from Corollary \ref{coro:IH_vs_Hc_Q} a simpler expression (Theorem \ref{theo:IH_vs_Hc_simple}) for \(e(\Acal_{n,?}^*, \IC_{\ell}(V))\) (in \(K_0(\Rep_{\Qell}(\GSpbf_{2n}(\A_f) \times \GalQ))\)) in terms of \(e_c(\Acal_{n',?,\Qbar}, \Fcal_{\ell}(V'))\) (for \(n' \leq n\) and certain algebraic representations \(V'\) of \(\GSpbf_{2n'}\)) and \(L^2\)-cohomology of arithmetic locally symmetric spaces for the groups \(\GLbf_1\) and \(\GLbf_2\), again via parabolic induction.
To prove it we first recall in Section \ref{sec:Franke_formula} Franke's formula expressing, for \(V\) an irreducible algebraic representation of a reductive group \(\Gbf\) over \(\Q\), the Euler characteristic \(e(\Gbf, V)\) (in \(K_0(\Rep_{\Q}(\Gbf(\A_f)))\), notation as in Section \ref{sec:arith_gp_coh}) in terms of \(e_{(2)}(\Lbf, V')\) (see \eqref{eq:intro_e2}) for \(\R\)-cuspidal Levi subgroups \(\Lbf\) of \(\Gbf\).
Since \(\R\)-cuspidal Levi subgroups of \(\GLbf_N\) are products of \(\GLbf_1\)'s and \(\GLbf_2\)'s, Franke's formula plugged in Corollary \ref{coro:IH_vs_Hc_Q} gives a new formula for \(e(\Acal_{n,?,\Qbar}^*, \IC_{\ell}(V))\), now with a double sum.
Simplifying this to obtain a single sum is a combinatorial matter.

In Section \ref{sec:Hc_from_IH} we invert this formula to obtain an explicit formula (Theorem \ref{theo:Hc_from_IH}) expressing compactly supported cohomology of local systems on \(\Acal_{n,?,\Qbar}\) in terms of intersection cohomology (and again \(e_{(2)}(\GLbf_N, V')\) for \(N \in \{1,2\}\)).
This is again a combinatorial matter.
Specializing to level one, together with Theorem \ref{thm:sigma_is_tensor} this gives an explicit (and, we believe, as simple as possible) description of \(e_c(\Acal_{n,\Qbar}, \Fcal(V))\) in terms of formal Arthur-Langlands parameters \(\psi\) and their associated Galois representations.
The case \(n=3\) (forgetting the Hecke action) verifies the main conjecture of Bergström, Faber and van der Geer in \cite{BFG}.

\subsection{Franke's formula}
\label{sec:Franke_formula}

In this section we consider an arbitrary connected reductive group \(\Gbf\) over \(\Q\).
Fix an open subgroup \(K_{\infty}\) of a maximal compact subgroup \(K_{\infty}^{\max}\) of \(\Gbf(\R)\) and denote \(\Xcal = \Gbf(\R) / K_{\infty} \Abf_{\Gbf}(\R)^0\).
Let \(V\) be an irreducible algebraic representation\footnote{More generally we may consider an irreducible finite-dimensional \((\gfrak,K_\infty^{\max})\)-module, see Appendix \ref{sec:irr_fd_gK}.} of \(\Gbf_{\C}\), that we will consider as a representation of \(\Gbf(\R)\) or \(\Gbf(\Q)\).
Franke's spectral sequence (I in \cite[Theorem 19]{Franke}) implies a formula for the Euler characteristic \(e(\Gbf, K_{\infty}, V)\) of the admissible graded representation of \(\Gbf(\A_f)\)
\[ \varinjlim_{K_f} H^\bullet(\Gbf(\Q) \backslash (\Xcal \times \Gbf(\A_f) / K_f), \Fcal^{K_f}(V)), \]
already considered in Section \ref{sec:IH_vs_Hc_Morel}, in terms of \((\lfrak, K_{\infty, \Lbf})\)-cohomology of the \emph{discrete} automorphic spectrum for \(\Lbf\) with respect to certain finite-dimensional representations of \(\Lbf(\R)\), as \(\Lbf\) varies among the Levi subgroups of \(\Gbf\) (over \(\Q\), up to conjugacy) which are \(\R\)-cuspidal.
We found it clearer to reformulate Franke's filtration for the space of automorphic forms (Theorem \ref{thm:Franke_fil}) and to deduce the formula for the Euler characteristic (thanks to Borel's conjecture, also proved by Franke, see Theorem \ref{thm:Franke_Borel_conj} below) following \cite[\S 7.4]{Franke}.
The final result for the Euler characteristic is Corollary \ref{cor:Franke_Euler2}.

In fact we will only require the case where \(\Gbf = \GLbf_N\) and \(K_{\infty}\) is maximal in later sections, but it would also be natural to take \(\Gbf = \GSpbf_{2n}\) and \(K_{\infty}\) connected to compare Franke's formula in this case (or rather its dual for compactly supported cohomology) to Theorem \ref{theo:Hc_from_IH} via Zucker's conjecture (\cite{LoijZucker}, \cite{SaperStern}, \cite{LooijengaRapoport}).
We give the explicit formula for \(\Gbf = \GLbf_N\) in Corollary \ref{coro:Euler_ord_GL}.

\subsubsection{Borel's conjecture}

We recall Borel's conjecture (proved by Franke in \cite[\S 7.4]{Franke}).
Recall from the beginning of Section \ref{sec:IH_vs_Hc_Morel} that any representation \(V\) of \(\Gbf(\Q)\) over \(\C\) gives rise to local systems \(\Fcal^{K_f} V\) on the manifolds \(\Gbf(\Q) \backslash (\Xcal \times \Gbf(\A_f)/K_f)\) for any neat compact open subgroup \(K_f\) of \(\Gbf(\A_f)\).
Taking cohomology and the colimit over all such levels yields admissible representations of \(\Gbf(\A_f)\)
\[ H^i(\Gbf, K_\infty, V) := \varinjlim_{K_f} H^i(\Gbf(\Q) \backslash (\Xcal \times \Gbf(\A_f)/K_f), \Fcal^{K_f} V). \]
Now assume that \(V\) is the restriction of an irreducible finite-dimensional \((\gfrak,K_\infty^{\max})\)-module (by Lemma \ref{lem:fd_rep_real_gp_alg} this is equivalent to an irreducible finite-dimensional representation of \(\Gbf(\R)\)).
Translating the de Rham comparison isomorphism gives
\[ H^i(\Gbf, K_\infty, V) \simeq H^i((\gfrak/\afrak_\Gbf, K_\infty), (C^\infty(\Gbf(\Q) \backslash \Gbf(\A)) \otimes V)^{\Abf_\Gbf(\R)^0}) \]
where
\[ C^\infty(\Gbf(\Q) \backslash \Gbf(\A)) := \varinjlim_{K_f} C^\infty(\Gbf(\Q) \backslash \Gbf(\A)/K_f). \]
Denote by \(\Acal(\Gbf)\) the space of automorphic forms on \(\Gbf(\Q) \backslash \Gbf(\A)\).
The following theorem is \cite[Theorem 18]{Franke}.

\begin{theo}[Franke] \label{thm:Franke_Borel_conj}
  The inclusion \(\Acal(\Gbf) \hookrightarrow C^\infty(\Gbf(\Q) \backslash \Gbf(\A))\) induces an isomorphism of admissible representations of \(\Gbf(\A_f)\)
  \[ H^\bullet((\gfrak/\afrak_\Gbf, K_\infty), (\Acal(\Gbf) \otimes V)^{\Abf_\Gbf(\R)^0}) \simeq H^\bullet((\gfrak/\afrak_\Gbf, K_\infty), (C^\infty(\Gbf(\Q) \backslash \Gbf(\A)) \otimes V)^{\Abf_\Gbf(\R)^0}). \]
\end{theo}

\subsubsection{Euler characteristic of the discrete automorphic spectrum}

Let us first precise our notation for the cohomology of the discrete automorphic spectrum, and make precise the dependence on the finite-dimensional representation of a Lie group that occurs.
As above let \(V\) be an irreducible finite-dimensional \((\gfrak, K_\infty^{\max})\)-module.

Denote by \(\omega_V\) the central character of \(V\) and let \(\xi^{-1}\) its restriction to \(\Abf_{\Gbf}(\R)^0\).
Denote by \(\Acal^2(\Gbf, \xi)\) the space of automorphic forms on \(\Gbf(\Q) \backslash \Gbf(\A)\) transforming under \(\Abf_{\Gbf}(\R)^0\) by \(\xi\) and square-integrable modulo \(\Abf_{\Gbf}(\R)^0\).
Note that this makes sense because \(\Abf_{\Gbf}(\R)^0\) is canonically a direct factor of \(\Gbf(\A)\) (with complement \(\Gbf(\A)^1\), the subgroup of \(g\) such that \(|\chi(g)| = 1\) for all \(\chi \in X^*(\Gbf)^{\GalQ}\)), and that there is a canonical isomorphism \(\Acal^2(\Gbf, \xi) \simeq \Acal^2(\Gbf, 1) \otimes \xi\) as \((\gfrak, K_{\infty}) \times \Gbf(\A_f)\)-modules, where on the right-hand side \(\xi\) is considered as a character of \(\Gbf(\A)\) (trivial on \(\Gbf(\A)^1\)).
These modules are semisimple, with finitely many constituents having given infinitesimal character and level.
Consider for a level \(K_f \subset \Gbf(\A_f)\) the cohomology groups
\[ H^i_{(2)}(\Gbf, K_{\infty}, V)^{K_f} := H^i((\gfrak / \afrak_{\Gbf}, K_{\infty}), \Acal^2(\Gbf, \xi)^{K_f} \otimes V) \]
endowed with the action of \(\Hcal(\Gbf(\A_f) // K_f)\).
By Wigner's lemma \cite[Corollary I.4.2]{BoWa} one can replace \(\Acal^2(\Gbf, \xi)^{K_f}\) by the direct sum of its constituents (for the action of \(\gfrak\)) having infinitesimal character opposite to that of \(V\), so that \(H^i_{(2)}(\Gbf, K_{\infty}, V)^{K_f}\) has finite dimension over \(\C\).
Of course by varying \(K_f\) we obtain an admissible graded object of \(\Hecke(\Gbf(\A_f), \C)\), which justifies the notation.
If \(\Gbf(\R) / \Abf_{\Gbf}(\R)\) admits discrete series representations these cohomology groups can be identified to \(L^2\)-cohomology groups (see \cite{BorelCasselman}), and in the Hermitian case to intersection cohomology groups as in Section \ref{sec:minimal_cpctif} by Zucker's conjecture (see \cite{LoijZucker}, \cite{SaperStern}, \cite{LooijengaRapoport}).
For \(\chi \in \C \otimes X^*(\Gbf)^{\GalQ}\) there is a canonical isomorphism \(H^i_{(2)}(\Gbf, K_{\infty}, V \otimes |\chi|) \simeq H^i_{(2)}(\Gbf, K_{\infty}, V) \otimes |\chi|_f^{-1}\), where for \(\chi = s \otimes \chi_0\) we denote \(|\chi|_f( (g_p)_p ) = \prod_p |\chi(g_p)|_p^s\) (compare with Remarks \ref{rema:twisting_sim} and \ref{rema:twisting_top_coh}).
\begin{defi} \label{def:e_2_G_Kinfty_V}
  In the setting above we denote \(e_{(2)}(\Gbf, K_{\infty}, V)\) the Euler characteristic of \(H_{(2)}^\bullet(\Gbf, K_{\infty}, V)\) in the Grothendieck group of admissible representations of \(\Gbf(\A_f)\).
\end{defi}
We know that this Euler characteristic vanishes if \(\Gbf\) is not \(\R\)-cuspidal (see \cite[bottom of p.\ 266]{Franke}).
Recall that \(\Gbf\) is said to be \(\R\)-cuspidal if \((\Gbf / \Abf_{\Gbf})(\R)\) has discrete series, or equivalently if there exists a Langlands parameter \(\varphi: W_{\R} \rightarrow {}^L \Gbf\) such that \(\Cent( \varphi, \Lhat) / Z(\Lhat)^{\Gal_{\Q}}\) is finite.
Note that this notion depends on the \(\Q\)-structure of \(\Gbf\), not just on \(\Gbf_{\R}\).
For example \(\GLbf_{N, \Q}\) is \(\R\)-cuspidal if and only if \(N \leq 2\).

\begin{prodef} \label{prodef:varphi_V}
  Let \(\Hbf\) be a connected reductive group over \(\R\) and let \(K\) be a maximal compact subgroup of \(\Hbf(\R)\).
  Let \(V\) be an irreducible finite-dimensional \((\hfrak,K)\)-module and denote by \(\tau_V\) its infinitesimal character and by \(\omega_V\) its central character.
  Assume that \(\Hbf(\R)\) admits essentially discrete series representations (equivalently, that modulo its center it admits an anisotropic maximal torus).
  There exists a unique (up to conjugation by \(\Hhat(\C)\)) essentially discrete Langlands parameter \(\varphi_V: W_\R \to {}^L \Hbf(\C)\) such that the L-packet of \(\varphi_V\) is the set of essentially discrete \((\hfrak,K)\)-modules with infinitesimal character \(-\tau_V\) and central character \(\omega_V^{-1}\).
\end{prodef}
\begin{proof}
  Up to a reduction using a z-extension this follows from the definition of discrete L-packets in \cite{Langlands} and the description of irreducible finite-dimensional \((\hfrak,K)\)-modules in Appendix \ref{sec:irr_fd_gK}.
\end{proof}

Now assume that \(\Gbf\) is \(\R\)-cuspidal.
Denote by \(\tau_V\) the infinitesimal character of \(V\).
In general \(V\) is not determined by \(\tau_V\) and so \(H_{(2)}^\bullet(\Gbf, K_{\infty}, V)\) really depends on \(V\).
However, it turns out (Corollary \ref{coro:Euler2_dep_phi} below) that the Euler characteristic \(e_{(2)}(\Gbf, K_{\infty}, V)\) only depends on the discrete Langlands parameter \(\varphi_V : W_{\R} \rightarrow {}^L \Gbf(\C)\).
This is remarkable since it is not true that the individual cohomology groups only depend on \(\varphi_V\), nor is it true that Euler characteristics for ordinary or compactly supported cohomology only depend on \(\varphi_V\).
Let us now prove this independence statement for \(e_{(2)}(\Gbf, K_{\infty}, V)\).
The question is clearly local at the real place, i.e.\ it is enough to show that \(e((\gfrak / \afrak, K_{\infty}), \pi_{\infty} \otimes V) = e((\gfrak / \afrak, K_{\infty}), \pi_{\infty} \otimes V')\) for any finite-dimensional irreducible continuous representations \(V\) and \(V'\) of \(\Gbf(\R)\) having equal infinitesimal characters and central characters and any finite length admissible \((\gfrak, K_{\infty}^{\max})\)-module \(\pi_{\infty}\) with opposite central character.
Of course both Euler characteristics only depends on the image of \(\pi_{\infty}\)
in the Grothendieck group, and so they may be computed in the basis of standard
modules.
The argument on p.\ 214 of \cite{CloDel} also applies to show that for
\(\pi_{\infty}\) a standard module corresponding to a proper Levi subgroup of
\(\Gbf\) we have \(e((\gfrak / \afrak_{\Gbf}, K_{\infty}^0), \pi_{\infty} \otimes
V) = 0\) as a representation of the finite \(2\)-torsion group
\(K_{\infty}^{\max}/K_{\infty}^0 = \Gbf(\R) / \Gbf(\R)^0\).
Wigner's lemma implies that cohomology (and thus Euler characteristic) vanishes
for tempered representations which are not part of the discrete series (modulo
\(\Abf_{\Gbf}\)) L-packet corresponding to \(\varphi_V\).
We are thus reduced to the case of discrete series \(\pi_{\infty}\), for which
cohomology is completely computed in \cite[\S II.5]{BoWa}.
We recall (a slight variation of) this result in the following theorem.

\begin{theo} \label{thm:coh_DS}
  Let \(\Gbf\) be a connected reductive group over \(\R\), \(\Zbf_{\Gbf}\) its center and \(\Abf_{\Gbf}\) its maximal split central torus.
  Fix a maximal compact subgroup \(K^{\max}\) of \(\Gbf(\R)\), and let \(q(\Gbf) = \frac{1}{2} \dim \Gbf(\R) / K^{\max} \Abf_{\Gbf}(\R)^0\).
  Denote by \(K_{\Zbf_{\Gbf}(\R)}\) the maximal compact subgroup of the center \(\Zbf_{\Gbf}(\R)\) of \(\Gbf(\R)\).
  Let \(V\) be a finite-dimensional \((\gfrak, K^{\max})\)-module and let \(\omega_V\) be its central character, a character of \(\Zbf_{\Gbf}(\R)\).
  Let \(K\) be an open subgroup of \(K^{\max}\), and \(K' = K K_{\Zbf_{\Gbf}(\R)}\), also an open subgroup of \(K^{\max}\).
  Let \(\pi\) be a discrete series \((\gfrak, K^{\max})\)-module with central character \(\omega_{\pi}\) which coincides with \(\omega_V^{-1}\) on \(\Abf_{\Gbf}(\R)^0\).
  We have
  \[ \dim H^i( (\gfrak/ \afrak_{\Gbf}, K), \pi \otimes V) = \begin{cases}
    |K^{\max}/K'| & \text{ if } i = q(\Gbf) \text{ and } \omega_{\pi} \omega_V |_{\Zbf_{\Gbf}(\R) \cap K} = 1 \text{ and } \tau_{\pi} = -\tau_V \\
    0 & \text{ otherwise.}
  \end{cases}\]
\end{theo}
\begin{proof}
  This follows from \cite[Theorem II.5.3]{BoWa} and the proof of Proposition
  II.5.7 loc.\ cit., using the precise description of discrete series (from the
  case where \(\Gbf\) is semisimple and simply connected studied by
  Harish-Chandra) in \cite[pp.\ 134-135]{Langlands}.
  Let us give the details.
  First note that \(V\) is irreducible as a \(\gfrak\)-module.
  If \(\Gbf_{\der}\) is simply connected this was already observed in the proof of
  Lemma \ref{lem:fd_rep_real_gp_alg}, and the general case follows by taking a
  z-extension of \(\Gbf\).
  Let \(K_0 = K^0 K_{\Zbf_{\Gbf}(\R)}\).
  The restriction of \(\pi\) to \(K\) has \(r := |K^{\max}/K'|\) irreducible
  constituents \(\pi'_1, \dots, \pi'_r\) which are distinct discrete series
  \((\gfrak, K)\)-modules.
  We have
  \begin{align*}
    \Hom_K \left( \bigwedge^i (\gfrak / (\afrak_{\Gbf} \oplus \kfrak)), \pi
    \otimes V \right) & \simeq \bigoplus_{j=1}^r \Hom_K \left( \bigwedge^i
    (\gfrak / (\afrak_{\Gbf} \oplus \kfrak)), \pi'_j \otimes V \right) \\
    & \simeq \bigoplus_{j=1}^r \Hom_{K \cap K_0} \left( \bigwedge^i (\gfrak /
    (\afrak_{\Gbf} \oplus \kfrak)), \pi'_{0,j} \otimes V \right)
  \end{align*}
  where \(\pi'_{0,j}\) is any of the constituents of the restriction of \(\pi'_j\)
  to \(K \cap K_0\), so that \(\pi'_j \simeq \Ind_{K \cap K_0}^{K}( \pi'_{0,j})\).
  Let \(K_{\sico}\) be the preimage of \(K^{\max}\) in the connected semisimple Lie
  group \(\Gbf_{\sico}(\R)\) (note that \(K_{\sico}\) is connected, so it is also
  the preimage of \(K^0\)).
  Then \((K \cap K_0) / (K \cap K_{\Zbf_{\Gbf}(\R)}) = K^0 / (K^0 \cap
  K_{\Zbf_{\Gbf}(\R)})\) is naturally a quotient of \(K_{\sico}\).
  Taking cohomology, we obtain
  \begin{align*}
    H^i \left( (\gfrak / \afrak_{\Gbf}, K), \pi \otimes V \right) & \simeq
    \bigoplus_{j=1}^r H^i \left( (\gfrak / \afrak_{\Gbf}, K \cap K_0),
    \pi'_{0,j} \otimes V \right) \\
    & = \bigoplus_{j=1}^r H^i \left( (\gfrak / \zfrak_{\Gbf}, (K \cap K_0)/(K
    \cap \Zbf_{\Gbf}(\R))), (\pi'_{0,j} \otimes V)^{K \cap \Zbf_{\Gbf}(\R)}
    \right) \\
    & = \bigoplus_{j=1}^r H^i \left( (\gfrak_{\sico}, K_{\sico}), (\pi'_{0,j}
    \otimes V)^{K \cap \Zbf_{\Gbf}(\R)} \right)
  \end{align*}
  where \cite[Corollary I.6.6]{BoWa} is used for the second line.
  The last expression is computed by \cite[Theorem II.5.3]{BoWa}.

\end{proof}

We can finally deduce the independence result.

\begin{coro} \label{coro:Euler2_dep_phi}
  Let \(\Gbf\) be a connected reductive group over \(\Q\).
  Assume that \(\Gbf\) is \(\R\)-cuspidal.
  Let \(V\) and \(V'\) be finite-dimensional continuous irreducible representations
  of \(\Gbf(\R)\) such that \(\varphi_V\) and \(\varphi_{V'}\) are conjugated by
  \(\Ghat\).
  Let \(K_{\infty}\) be an open subgroup of a maximal compact subgroup of
  \(\Gbf(\R)\).
  Then \(e_{(2)}(\Gbf, K_{\infty}, V) = e_{(2)}(\Gbf, K_{\infty}, V')\) in
  \(K_0(\Rep_{\C}(\Gbf(\A_f)))\).
\end{coro}

Using a z-extension and following the proof of Lemma
\ref{lem:fd_rep_real_gp_alg} it is easy to see that any discrete Langlands
parameter \(\varphi : W_{\R} \rightarrow {}^L \Gbf\) can be written as \(\varphi_V\)
for some \(V\) as above.
Together with Corollary \ref{coro:Euler2_dep_phi} this justifies the following
definition/notation: \(e_{(2)}(\Gbf, K_{\infty}, \varphi) := e_{(2)}(\Gbf,
K_{\infty}, V)\) for any \(V\) such that \(\varphi_V \sim \varphi\).

\begin{exam} \label{exam:Euler_L2_coh_GL}
  For \(\Gbf = \GLbf_1\) this provides no simplification.
  For \(a \in \Z\) we denote
  \begin{equation} \label{eq:e_GL1}
    e(\GLbf_1, a) := e_{(2)}(\GLbf_1, \{ \pm 1 \}, x \mapsto x^a) = \sum_{\substack{\chi : \Q^{\times} \backslash \A^{\times} \rightarrow \C^{\times} \\ \chi|_{\R^{\times}} = x \mapsto x^{-a}}} \chi_f.
  \end{equation}
  Note that all other characters of \(\R^{\times}\) are obtained by twisting by \(x \mapsto |x|^s\) for some \(s \in \C\).
  Note also that \(e(\GLbf_1, a)^{\GLbf_1(\Zhat)} = 0\) if \(a\) is odd.

  For \(\Gbf = \GLbf_2\) (over \(\Q\)) the center of \(\Gbf(\R)\) is included in its identity component and so \(\varphi_V \sim \varphi_{V'}\) if and only if \(\tau_V = \tau_{V'}\).
  In the sequel we will only need to consider \(V\) with real infinitesimal character, which can then be written \((a+1/2,b-1/2)\) with \(a, b \in \R\) and \(a-b \in \Z_{\geq 0}\) (i.e.\ \(V \simeq \Sym^{a-b} \Std \otimes |\det|^b\) or \(\Sym^{a-b} \Std \otimes \det |\det|^{b-1}\)), and we will simply denote \(e_{(2)}(\GLbf_2, a, b)\) for \(e_{(2)}(\GLbf_2, K_{\infty}, V)\) with \(K_{\infty}\) maximal.
  This can be easily described in terms of modular cusp forms.
  First note that for \((a,b)\) as above and \(s \in \C\) we have
  \[ e_{(2)}(\GLbf_2, a+s, b+s) = e_{(2)}(\GLbf_2, a, b) \otimes |\det|_f^{-s}. \]
  Denote \(S_k = \varinjlim_{\Gamma} S_k(\Gamma)\) where \(\Gamma\) ranges over congruence subgroups of \(\SLbf_2(\Z)\), with the usual action of \(\GLbf_2(\A_f)\).
  Then for \(k \in \Z_{\geq 2}\) we have
  \begin{equation} \label{eq:e2_GL2}
    e_{(2)}(\GLbf_2, k-2, 0) =
    \begin{cases}
      -S_k & \text{if } k>2, \\
      -S_2 + e(\GLbf_1, 0) \circ \det & \text{if } k=2.
    \end{cases}
  \end{equation}
  In particular \(e_{(2)}(\GLbf_2, a, b)^{\GLbf_2(\Zhat)} = 0\) if \(a-b\) is odd.
\end{exam}

\subsubsection{Franke's filtration of the space of automorphic forms}

We recall in Theorem \ref{thm:Franke_fil} below Franke's filtration of the space of automorphic forms, and his description of the associated graduated pieces.
This description uses the Langlands positivity condition on characters of parabolic subgroups (see \cite[\S 5.4.1]{WallachReal} in the real case, \cite[p.233 l.-4]{Franke} for the case at hand), indeed Franke's filtration may be interpreted as the global analogue of Langlands' classification of irreducible \((\gfrak,K)\)-modules as Langlands quotients of standard modules \cite{Langlands}.
This positivity condition is usually expressed using relative root systems.
We favor formulating this condition using Langlands dual groups (and thus absolute root systems), which require some preparation.
Although we give proofs for this reformulation in the general case, we will ultimately only use the case of \(\GLbf_{n,\Q}\) for which the following lemmas are essentially trivial.
The reader only interested in this case may safely skip the proofs and focus on the examples.

In this section \(\Gbf\) is a connected reductive group over \(\Q\).
Following \cite[\S 1]{Kottwitz_STFcusptemp} we consider the Langlands dual group of \(\Gbf\) as an extension \(\Ghat \rightarrow {}^L \Gbf \rightarrow W_{\Q}\) of \(W_{\Q}\) by \(\Ghat\) together with a \(\Ghat(\Qbar)\)-conjugacy class of splittings, i.e.\ quadruples \((\Bcal,\Tcal,(X_\alpha)_{\alpha \in \Delta(\Tcal,\Bcal)},s)\) where \((\Bcal,\Tcal,(X_\alpha)_\alpha)\) is a pinning of \(\Ghat\) and \(s: W_{\Q} \to {}^L \Gbf\) is a section such that each \(s(\sigma)\) stabilizes this pinning.
We will call such quadruples \emph{distinguished splittings}.
Recall from \cite[\S I.3]{BorelCorvallis} the notion of parabolic and Levi subgroups of \({}^L \Gbf\), and the fact that there is a natural injective map
\begin{equation} \label{eq:map_para_conj_to_para_Lgp}
  \{ \text{parabolic subgroups of } \Gbf \} / \Gbf(\Q)\mathrm{-conj} \hookrightarrow \{ \text{parabolic subgroups of } {}^L \Gbf \} / \Ghat\mathrm{-conj}.
\end{equation}
This injection is a bijection if \(\Gbf\) is quasi-split.
If \(\Pbf\) is a parabolic subgroup of \(\Gbf\) with reductive quotient \(\Lbf\), choosing distinguished splittings \(\zeta_\Gbf\) and \(\zeta_\Lbf\) for \({}^L \Gbf\) and \({}^L \Lbf\) yields an embedding
\[ \iota[\Pbf,\zeta_\Gbf,\zeta_\Lbf]: {}^L \Lbf \to {}^L \Gbf \]
with image the standard (for \(\zeta_\Gbf\)) Levi subgroup of \({}^L \Gbf\), and denoting by \(\Bcal\) the Borel subgroup of \(\Ghat\) occurring in \(\zeta_\Gbf\) we have that \(\iota[\Pbf,\zeta_\Gbf,\zeta_\Lbf]({}^L \Lbf) \Bcal\) is the standard parabolic subgroup of \({}^L \Gbf\) corresponding to \(\Pbf\) in \eqref{eq:map_para_conj_to_para_Lgp}.
For \(g \in \Ghat(\Qbar)\) and \(l \in \Lhat(\Qbar)\) we have
\[ \iota[\Pbf, \Ad(g) \zeta_\Gbf, \Ad(l) \zeta_\Lbf] = \Ad(g) \circ \iota[\Pbf, \zeta_\Gbf, \zeta_\Lbf] \circ \Ad(l)^{-1}. \]
We can slightly change the point of view: if \(\Pbf\) and \(\Pcal\) are parabolic subgroups of \(\Gbf\) and \({}^L \Gbf\) corresponding to each other in \eqref{eq:map_para_conj_to_para_Lgp}, and if we choose a Levi factor \(\Lcal\) of \(\Pcal\), then we have an embedding
\[ \iota[\Pbf,\Pcal,\Lcal]: {}^L \Lbf \to {}^L \Gbf \]
well-defined up to composing with \(\Ad(l)\) for some \(l \in \Lhat(\Qbar)\) and satisfying \(\iota[\Pbf,\Pcal,\Lcal]({}^L \Lbf) = \Lcal\) (choose any distinguished splitting \(\zeta_\Lbf\) for \({}^L \Lbf\) and a distinguished splitting \(\zeta_\Gbf\) for \({}^L \Gbf\) for which \((\Pcal,\Lcal)\) is standard).
It is clear that the \(\Ghat(\Qbar)\)-conjugacy class of \(\iota[\Pbf,\Pcal,\Lcal]\) does not depend on the choice of \((\Pcal,\Lcal)\).

\begin{lemm} \label{lem:Levi_para_dual}
  Let \(\Lbf\) be a Levi subgroup of \(\Gbf\).
  \begin{enumerate}
  \item Let \(\Pbf\) and \(\Pbf'\) be parabolic subgroups of \(\Gbf\) admitting \(\Lbf\) as a Levi factor.
    Then \(\iota[\Pbf,\Pcal,\Lcal]\) and \(\iota[\Pbf',\Pcal',\Lcal']\) are conjugated under \(\Ghat(\Qbar)\).
    In particular we have a well-defined conjugacy class \(\Ecal(\Lbf, \Gbf)\) of embeddings \({}^L \Lbf \to {}^L \Gbf\).
  \item For \(\iota_{\Lbf} \in \Ecal(\Lbf, \Gbf)\) we have a bijection \(\Pbf \mapsto {}^L \Pbf\), compatible with \eqref{eq:map_para_conj_to_para_Lgp}, between the set of parabolic subgroups of \(\Gbf\) having \(\Lbf\) as a Levi factor and the set of parabolic subgroups of \({}^L \Gbf\) having \(\Lcal := \iota_{\Lbf}({}^L \Lbf)\) as a Levi factor, compatible with \eqref{eq:map_para_conj_to_para_Lgp}.
    The bijection is determined by the property that \(\iota_\Lbf\) belongs to the \(\Lhat(\Qbar)\)-orbit of \(\iota[\Pbf, {}^L \Pbf, \Lcal]\).
  \end{enumerate}
\end{lemm}
\begin{proof}
  The first statement is \cite[Lemma 2.5]{Langlands}.
  In order to prove the second statement we briefly recall the construction.
  Fix a Borel pair \((\Bbf_\Lbf, \Tbf)\) in \(\Lbf_{\Qbar}\).
  This defines a Borel subgroup \(\Bbf\) (resp.\ \(\Bbf'\)) of \(\Pbf_{\Qbar}\) (resp.\ \(\Pbf'_{\Qbar}\)) by the relation \(\Bbf \cap \Lbf_{\Qbar} = \Bbf_\Lbf\) (resp.\ \(\Bbf' \cap \Lbf_{\Qbar} = \Bbf_\Lbf\)).
  Let \(\zeta_\Gbf = (\Bcal,\Tcal,(X_\alpha)_\alpha,s)\) be a distinguished splitting for \(\Gbf\).
  The pairs \((\Bbf,\Tbf)\) and \((\Bcal,\Tcal)\) yield identifications \(X^*(\Tcal) \simeq X_*(\Tbf)\) and \(W(\Tbf,\Gbf_{\Qbar}) \simeq W(\Tcal,\Ghat)\).
  There is a unique \(x \in W(\Tbf,\Gbf_{\Qbar})\) satisfying \(\Ad(x) \Bbf = \Bbf'\), and denoting by \(n: W(\Tcal,\Ghat) \to N(\Tcal,\Ghat)\) the set-theoretic section induced by \((\Bcal,\Tcal,(X_\alpha)_\alpha)\) (see \cite[\S 9.3.3]{Springer}) we have
  \[ \Ad(n(x)) \circ \iota[\Pbf',\zeta_\Gbf,\zeta_\Lbf] = \iota[\Pbf,\zeta_\Gbf,\zeta_\Lbf]. \]

  We now turn to the second point.
  Fix \(\iota_\Lbf \in \Ecal(\Lbf, \Gbf)\).
  It is straightforward to define a map \(\Pbf \mapsto {}^L \Pbf\) such that \(\iota_\Lbf\) belongs to the \(\Lhat(\Qbar)\)-orbit of \(\iota[\Pbf,{}^L \Pbf, \Lcal]\): starting from any parabolic subgroup \(\Pcal\) of \({}^L \Gbf\) corresponding to \(\Pbf\) via \eqref{eq:map_para_conj_to_para_Lgp} we can choose a Levi factor \(\Lcal_0\) of \(\Pbf\) and \(\iota[\Pbf,\Pcal,\Lcal_0]\) in the corresponding \(\Lhat(\Qbar)\)-orbit, and there exists \(g \in \Ghat(\Qbar)\) satisfying \(\Ad(g) \circ \iota[\Pbf,\Pcal,\Lcal_0] = \iota_\Lbf\) and we define \({}^L \Pbf = \Ad(g) \Pcal\), which does not depend on the choice of \(\Pcal\), \(\Lcal_0\) and \(g\) (\(g\) is unique up to left multiplication by \(Z(\Lcal) \cap \Lcal^0\)).
  We now define the inverse map.
  We may assume \(\iota_\Lbf = \iota[\Pbf,\zeta_\Gbf,\zeta_\Lbf]\) and as above we fix a Borel pair \((\Bbf_\Lbf, \Tbf)\) in \(\Lbf_{\Qbar}\) and denote by \(\Bbf\) the corresponding Borel subgroup of \(\Pbf_{\Qbar}\).
  Let \(\Pcal'\) be a parabolic subgroup of \({}^L \Gbf\) admitting \(\Lcal := \iota_\Lbf({}^L \Lbf)\) as a Levi factor.
  There is a unique \(g \in W(\Tcal,\Ghat)\) such that \(\Ad(g^{-1}) \Pcal'\) is a standard (for \(\zeta_\Gbf\)) parabolic subgroup of \({}^L \Gbf\) and \(\Ad(g^{-1})(\Bcal \cap \Lcal^0) \subset \Bcal\).
  Let \(x \in W(\Tbf,\Gbf_{\Qbar})\) corresponding to \(g\).
  We now check that \(\Ad(x) \Bbf\) contains \(\Bbf_\Lbf\) and that the simple roots for \((\Bbf_\Lbf, \Tbf)\) are also simple for \((\Ad(x) \Bbf, \Tbf)\) (this means that \(\Lbf_{\Qbar}\) is standard for \((\Ad(x) \Bbf, \Tbf)\) and so the corresponding \((\Ad(x) \Bbf, \Tbf)\)-standard parabolic subgroup \(\Pbf'\) of \(\Gbf_{\Qbar}\) admits \(\Lbf_{\Qbar}\) as a Levi factor).
  Because \(\Ad(g^{-1}) \Lcal^0\) is \((\Bcal,\Tcal)\)-standard and \(\Ad(g^{-1})(\Bcal \cap \Lcal^0) \subset \Bcal\) we know
  \[ \Ad(g^{-1}) \left( \Delta^\vee(\Tcal, \Bcal \cap \Lcal^0) \right) \subset \Delta^\vee(\Tcal, \Bcal \cap \Ad(g^{-1}) \Lcal^0) \]
  (where \(\Delta^\vee\) denotes the set of simple coroots).
  Thus \(\Ad(x)^*: \alpha \mapsto \alpha \circ \Ad(x)\) maps \(\Delta(\Tbf,\Bbf_\Lbf)\) to \(\Delta(\Tbf,\Bbf)\).
  Since \(\Ad(x)^*\) induces a bijection \(\Delta(\Tbf, \Ad(x) \Bbf) \simeq \Delta(\Tbf, \Bbf)\) we deduce \(\Delta(\Tbf, \Bbf_\Lbf) \subset \Delta(\Tbf, \Ad(x) \Bbf)\).
  We thus have a parabolic subgroup \(\Pbf'\) of \(\Gbf_{\Qbar}\) admitting \(\Lbf_{\Qbar}\) as a Levi factor and \(\Ad(x) \Bbf\) as a Borel subgroup, and we now check that it is defined over \(\Q\).
  For this we use the fact that \(W(\Tbf,\Gbf_{\Qbar}) \simeq W(\Tcal, \Ghat)\) is \(\GalQ\)-equivariant in the following sense.
  For \(\sigma \in \GalQ\) write \(\sigma(\Bbf,\Tbf) = \Ad(g_\sigma)^{-1}(\Bbf,\Tbf)\) where \(g_\sigma \in \Lbf(\Qbar)\) because \(\Pbf\) is defined over \(\Q\).
  Then \(g_\sigma \sigma(x) g_\sigma^{-1} \in W(\Tbf,\Gbf_{\Qbar})\) corresponds to \(s(\sigma) g s(\sigma)^{-1} \in W(\Tcal, \Ghat)\).
  Now \(\sigma(\Pbf')\) contains
  \[ \sigma(\Ad(x) \Bbf) = \Ad(\sigma(x) g_\sigma^{-1}) \Bbf = \Ad(g_\sigma^{-1} g_\sigma \sigma(x) g_\sigma^{-1}) \Bbf. \]
  We have \(s(\sigma) \in \Lcal\) (because \(\Lcal = \iota_\Lbf({}^L \Lbf)\) is formed using \(\zeta_\Gbf\)) and \(s(\sigma) \in \Ad(g^{-1}) \Lcal\) (because \(\Ad(g^{-1} \Lcal)\) is standard for \(\zeta_\Gbf\)) and so we have
  \[ s(\sigma) g s(\sigma)^{-1} g^{-1} \in \Lcal^0 \]
  which translates to
  \[ g_\sigma \sigma(x) g_\sigma^{-1} x^{-1} \in N(\Tbf, \Lbf(\Qbar)) \]
  and so \(\sigma(\Pbf')\) contains \(\Ad(nx) \Bbf\) for some \(n \in \Lbf(\Qbar)\), which implies \(\sigma(\Pbf') = \Pbf'\) and so \(\Pbf'\) is defined over \(\Q\).
  Comparing with the proof of the first point we see that we have just defined the inverse of \(\Pbf' \mapsto {}^L \Pbf'\).
\end{proof}

\begin{exam} \label{ex:Levi_para_dual_GL}
  Let us describe in simpler terms the content of Lemma \ref{lem:Levi_para_dual} in the case \(\Gbf=\GLbf_n\).
  A Levi subgroup \(\Lbf\) of \(\Gbf\) corresponds to a family \((V_i)_{i \in I}\) of non-zero subspaces of \(\Q^n\) satisfying \(\Q^n = \bigoplus_{i \in I} V_i\), with \(\Lbf\) equal to the intersection of the stabilizers of the \(V_i\)'s.
  A parabolic subgroup \(\Pbf\) of \(\Gbf\) corresponds to a total order on the index set \(I\), with \(\Pbf\) equal to the intersection of the stabilizers of the \(\bigoplus_{i \leq j} V_i\) (as \(j \in I\) varies).
  In other words using the identification of \(\Abf_{\Lbf}\) with \(\GLbf_1^I\), where for a commutative \(\Q\)-algebra \(R\) we let \(\ul{\lambda} = (\lambda_i)_i \in (R^\times)^I\) act by \(\lambda_i\) on \(R \otimes_{\Q} V_i\), the set of roots of \(\Abf_{\Lbf}\) in the unipotent radical of \(\Pbf\) is \((\ul{\lambda} \mapsto \lambda_i / \lambda_j)_{i<j}\).

  There is a natural identification of \({}^L \Lbf\) with \(\prod_{i \in I} \GL_{n_i,\Qbar} \times W_{\Q}\) where \(n_i = \dim_\Q V_i\), and for any distinguished splitting \((\Bcal,\Tcal,(X_\alpha)_\alpha,s)\) the section \(s: W_{\Q} \to {}^L \Lbf\) is the obvious one.
  This also applies to \(\Gbf\) instead of \(\Lbf\), and the parabolic subgroups of \({}^L \Gbf\) are simply the ones of the form \(\Pcal^0 \times W_{\Q}\) where \(\Pcal^0\) is a parabolic subgroup of \(\GL_{n,\Qbar}\).
  With these identifications the orbit \(\Ecal(\Lbf, \Gbf)\) is the obvious one: if \(I = \{i_1, \dots, i_k\}\) it contains
  \begin{align*}
    {}^L \Lbf & \longrightarrow {}^L \Gbf \\
    ((g_i)_i, w) & \longmapsto (\diag(g_{i_1}, \dots, g_{i_k}), w).
  \end{align*}
  In particular for \(\iota_\Lbf \in \Ecal(\Lbf, \Gbf)\) the set of parabolic subgroups of \({}^L \Gbf\) admitting \(\Lcal := \iota_\Lbf({}^L \Lbf)\) as a Levi factor is also parametrized by total orders on \(I\), and with these parametrizations the bijection in the second part of Lemma \ref{lem:Levi_para_dual} is simply the identity map (on total orders on \(I\)).
\end{exam}

Levi subgroups of \({}^L \Gbf\) of the form \(\iota_{\Lbf}({}^L \Lbf)\) (for some Levi subgroup \(\Lbf\) of \(\Gbf\) and some \(\iota_{\Lbf} \in \Ecal(\Lbf, \Gbf)\)) are called \emph{relevant}.

\begin{lemm} \label{lem:Levi_para_dual_overLevi}
  Let \(\Lbf\) be a Levi subgroup of \(\Gbf\).
  Let \(\iota_\Lbf \in \Ecal(\Lbf, \Gbf)\) and denote \(\Lcal = \iota_\Lbf({}^L \Lbf)\).
  \begin{enumerate}
  \item For a Levi subgroup \(\Mbf\) of \(\Gbf\) containing \(\Lbf\) and \(\iota_{\Lbf,\Mbf} \in \Ecal(\Lbf, \Mbf)\) there exists \(\iota_\Mbf \in \Ecal(\Mbf, \Gbf)\) satisfying \(\iota_\Mbf \circ \iota_{\Lbf,\Mbf} = \iota_\Lbf\), and \(\iota_\Mbf\) is unique up to composing with \(\Ad(g)\) for some \(g \in Z(\Lcal) \cap \Ghat(\Qbar)\).
    The class \(\{ \iota_\Mbf \circ \Ad(m) \,|\, m \in \Mhat(\Qbar) \}\) does not depend on the choice of \(\iota_{\Lbf,\Mbf}\), in particular \(\iota_\Mbf({}^L \Mbf)\) does not depend on this choice.
  \item For a Levi subgroup \(\Mcal\) of \({}^L \Gbf\) containing \(\Lcal\) there exists a unique Levi subgroup \(\Mbf\) of \(\Gbf\) containing \(\Lbf\) satisfying \(\iota_\Mbf({}^L \Mbf) = \Mcal\).
  \end{enumerate}
\end{lemm}
\begin{proof}
  \begin{enumerate}
  \item Choose a parabolic subgroup \(\Qbf\) of \(\Gbf\) admitting Levi factor \(\Mbf\) and a parabolic subgroup \(\Pbf_\Mbf\) of \(\Mbf\) admitting Levi factor \(\Lbf\).
    Then there is a unique parabolic subgroup \(\Pbf\) of \(\Qbf\) satisfying \(\Pbf \cap \Mbf = \Pbf_\Mbf\), and it admits Levi factor \(\Lbf\).
    Fix a distinguished splitting \(\zeta_\Lbf\) for \({}^L \Lbf\).
    There exists a distinguished splitting \(\zeta_\Mbf\) for \({}^L \Mbf\) (resp.\ \(\zeta_\Gbf\) for \({}^L \Gbf\)) satisfying \(\iota_{\Lbf,\Mbf} = \iota[\Pbf \cap \Mbf, \zeta_\Mbf, \zeta_\Lbf]\) (resp.\ \(\iota_\Lbf = \iota[\Pbf, \zeta_\Gbf, \zeta_\Lbf]\)).
    Checking that the composition
    \[ {}^L \Lbf \xrightarrow{\iota[\Pbf \cap \Mbf, \zeta_\Mbf, \zeta_\Lbf]} {}^L \Mbf \xrightarrow{\iota[\Qbf, \zeta_\Gbf, \zeta_\Mbf]} {}^L \Gbf \]
    is equal to \(\iota[\Pbf,\zeta_\Gbf,\zeta_\Lbf]\) is a formality.
    This proves the existence of \(\iota_\Mbf\).
    The uniqueness statements are easily checked and we leave the details to the reader.
  \item Let \(\Mcal \supset \Lcal\) be a Levi subgroup of \({}^L \Gbf\).
    Recall from \cite[Lemma 3.5]{BorelCorvallis} that \(\Lcal\) may be recovered as the centralizer of the torus \(Z(\Lcal)^0\) in \({}^L \Gbf\), and similarly for \(\Mcal\).
    Choose \(x_\Mcal \in \R \otimes X_*(Z(\Mcal))\) generic, i.e.\ such that the eigenvalues of \(x_\Mcal\) on \(\gfrakhat / \Lie \Mcal^0\) (which are real) are all non-zero.
    Choose \(x_\Lcal \in \R \otimes X_*(Z(\Lcal))\) similarly.
    Let \(\Ncal\) be the unipotent subgroup of \(\Ghat\) such that \(\Lie \Ncal\) is the direct sum of the eigenspaces corresponding to positive eigenvalues for the adjoint action of \(x_\Mcal\) on \(\gfrakhat\).
    Then \(\Mcal \Ncal\) is a parabolic subgroup of \({}^L \Gbf\) admitting Levi factor \(\Mcal\).
    Let \(\Ucal\) be the unipotent subgroup of \(\Ghat\) such that \(\Lie \Ucal\) is the direct sum of the eigenspaces corresponding to positive eigenvalues for the adjoint action of \(x_\Mcal + \epsilon x_\Lcal\) on \(\gfrakhat\), for small enough \(\epsilon>0\).
    Then \(\Lcal \Ucal\) is a parabolic subgroup of \({}^L \Gbf\) admitting Levi factor \(\Lcal\), and \(\Lcal \Ucal\) is contained in \(\Mcal \Ncal\).
    By the second part of Lemma \ref{lem:Levi_para_dual} \(\Lcal \Ucal\) corresponds to a parabolic subgroup \(\Pbf\) of \(\Gbf\) admitting Levi factor \(\Lbf\), more precisely we have \(\iota_\Lbf = \iota[\Pbf, \zeta_\Gbf, \zeta_\Lbf]\) for some splitting \(\zeta_\Gbf = (\Bcal, \Tcal, (X_\alpha)_\alpha, s)\) for \({}^L \Gbf\) satisfying \(\Bcal \subset \Lcal^0 \Ucal\).
    It is now easy to check that \(\Mcal^0 \Ncal\) corresponds to a parabolic subgroup \(\Qbf\) of \(\Gbf\) containing \(\Pbf\) (and defined over \(\Q\)), and that its Levi factor \(\Mbf\) containing \(\Lbf\) maps to \(\Mcal\).

    For uniqueness we observe that if we choose a maximal torus \(\Tbf\) in \(\Lbf_{\Qbar}\) and a maximal torus \(\Tcal\) in \(\Lcal^0\) then we have an identification \(\Tcal \simeq \That\) which is well-defined up to \(W(\Tbf,\Lbf_{\Qbar})\) and fixing such an identification the set of coroots for \(\Tbf \subset \Mbf_{\Qbar}\) corresponds to the set of roots for \(\Tcal \subset \Mcal^0\).
  \end{enumerate}
\end{proof}

\begin{exam}
  Again for \(\Gbf=\GLbf_n\) this lemma is almost a tautology.
  Reusing the notation of Example \ref{ex:Levi_para_dual_GL} if \(\Lbf\) corresponds to \((V_i)_{i \in I}\) (satisfying \(\Q^n = \bigoplus_{i \in I} V_i\)) then Levi subgroups of \(\Gbf\) containing \(\Lbf\) are parametrized by partitions of \(I\): if \(P\) is such a partition then \(\Mbf \supset \Lbf\) corresponds to \((\oplus_{i \in S} V_i)_{S \in P}\).
\end{exam}

\begin{defi} \label{def:tau_A}
  Let \(\tau\) be a semi-simple conjugacy class in \(\gfrakhat_\C = \Lie \Ghat_{\C}\).
  The Harish-Chandra isomorphism (see e.g.\ \cite[\S 3.2]{WallachReal}) \(Z(U(\gfrak)) \simeq \Ocal(\gfrakhat_\C)^{\Ghat_\C}\) allows us to see \(\tau\) as a morphism of \(\C\)-algebras \(Z(U(\gfrak)) \to \C\).
  We have a natural map from \(\afrak_\Gbf = \C \otimes X_*(\Abf_\Gbf) = \C \otimes_\R \Lie \Abf_\Gbf(\R)\) to \(Z(U(\gfrak))\), so \(\tau\) induces an element \(\tau_A\) of \(\afrak_\Gbf^* = \Hom_\C(\afrak_\Gbf, \C)\).
  Let \(\Cbf_\Gbf\) be the largest quotient of \(\Gbf\) which is a split torus, i.e.\ \(X^*(\Cbf_\Gbf) = \Hom(\Gbf, \GLbf_1) = X^*(\Gbf)^{\GalQ}\).
  The restriction map \(\C \otimes X^*(\Cbf_\Gbf) \to \C \otimes X^*(\Abf_\Gbf)\) is an isomorphism so we also see \(\tau_A\) as an element of \(\C \otimes X^*(\Cbf_\Gbf)\).

  This can also be interpreted dually: we have a natural map \(\Ghat \to \widehat{\Abf_\Gbf}\) and it is easy to check on the definition of the Harish-Chandra isomorphism that
  \[ \tau_A \in \afrak_\Gbf^* \simeq \C \otimes X^*(\Abf_\Gbf) \simeq \C \otimes X_*(\widehat{\Abf_\Gbf}) \simeq \Lie \widehat{\Abf_\Gbf}_\C \]
  is the image of \(\tau\) under the differential of this map.
  We also have an identification \(\widehat{\Cbf_\Gbf} \simeq Z(\Ghat)^{\GalQ,0}\) and we may also see \(\tau_A\) as an element of \(\Lie Z(\Ghat)^{\GalQ}_\C\).
  We thus has a canonical decomposition \(\tau = \tau_A + \tau_0\) where \(\tau_A \in \Lie Z(\Ghat)^{\GalQ}_\C\) and \(\tau_0\) is a semi-simple \(\Ghat(\C)\)-conjugacy class in to the kernel of \(\Lie \Ghat_\C \to \Lie \widehat{\Abf_\Gbf}_\C\).
\end{defi}

\begin{defi} \label{def:Pi_L_tauL}
  Let \(\Lbf\) be a Levi subgroup of \(\Gbf\) and \(\tau_{\Lbf}\) a semi-simple conjugacy class in \(\lfrakhat_\C = \Lie \Lhat_{\C}\).
  Choose \(\iota_\Lbf \in \Ecal(\Lbf, \Gbf)\) and denote \(\Lcal = \iota_\Lbf({}^L \Lbf)\).
  Let \(\tau_{\Lcal,A} \in \Lie Z(\Lcal_\C)^0\) be the image of \(\tau_{\Lbf,A}\) (Definition \ref{def:tau_A}) under the differential of \(\iota_{\Lbf,\C}\).
  Let \(\Mcal\) be the centralizer of \(\Re \tau_{\Lcal,A} \in \R \otimes X_*(Z(\Lcal)^0)\) in \({}^L \Gbf\), a Levi subgroup of \({}^L \Gbf\) (essentially by \cite[Lemma 3.5]{BorelCorvallis}).
  Let \(\Ncal\) be the unipotent subgroup of \(\Ghat\) such that \(\Lie \Ncal_{\C}\) is the direct sum of the eigenspaces for positive eigenvalues for the adjoint action of \(\Re \tau_{\Lcal, A}\).
  Then \(\Mcal \Ncal\) is a parabolic subgroup of \({}^L \Gbf\).
  Let \(\Pi(\Lcal, \tau_{\Lcal})\) be the set of parabolic subgroups \(\Pcal\) of \({}^L \Gbf\) admitting \(\Lcal\) as a Levi factor and which are contained in \(\Mcal \Ncal\).
  Let \(\Pi(\Lbf, \tau_{\Lbf})\) be the set of parabolic subgroups of \(\Gbf\) admitting \(\Lbf\) as a Levi factor which corresponds by the second part of Lemma \ref{lem:Levi_para_dual} to \(\Pi(\Lcal, \tau_{\Lcal})\).
  As the notation suggests \(\Pi(\Lbf, \tau_{\Lbf})\) does not depend on the choice of \(\iota_\Lbf\).
\end{defi}

\begin{rema} \label{rem:direct_def_Pi_L_tauL}
  In the setting of Definition \ref{def:Pi_L_tauL} we have a bijection \(\Pcal \mapsto \Pcal \cap \Mcal\) between \(\Pi(\Lcal, \tau_{\Lcal})\) and the set of parabolic subgroups of \(\Mcal\) admitting \(\Lcal\) as a Levi factor.
  We have a Levi subgroup \(\Mbf \supset \Lbf\) of \(\Gbf\) corresponding to \(\Mcal\) by the second part of Lemma \ref{lem:Levi_para_dual_overLevi}, and a parabolic subgroup \(\Qbf = \Mbf \Nbf\) corresponding to \(\Mcal \Ncal\) by the second part of Lemma \ref{lem:Levi_para_dual}.
  Then \(\Pi(\Lbf, \tau_\Lbf)\) is the set of parabolic subgroups \(\Pbf\) of \(\Qbf\) admitting \(\Lbf\) as a Levi factor.

  We may also give an equivalent definition of \(\Pi(\Lbf, \tau_\Lbf)\) which does not use dual groups: if we choose a maximal torus \(\Tbf\) of \(\Lbf_{\Qbar}\) then a parabolic subgroup \(\Pbf\) of \(\Gbf\) admitting \(\Lbf\) as a Levi factor belongs to \(\Pi(\Lbf, \tau_\Lbf)\) if and only if for any root \(\alpha\) of \(\Tbf\) in \(\Pbf_{\Qbar}\) we have
  \[ \Re \langle \alpha^\vee, \res^{\Cbf_\Lbf}_\Tbf \tau_{\Lbf,A} \rangle \geq 0. \]
\end{rema}

\begin{exam}
  Let us work out these definitions in the case where \(\Gbf = \GLbf_n\).
  Consider as in Example \ref{ex:Levi_para_dual_GL} a Levi subgroup \(\Lbf\) corresponding to \((V_i)_{i \in I}\).
  The semi-simple conjugacy class \(\tau_\Lbf\) in \(\C \otimes_{\Qbar} \lfrakhat \simeq \prod_{i \in I} \mathfrak{gl}_{n_i}(\C)\) is given by the family \(([x_{i,1}, \dots, x_{i,n_i}])_{i \in I}\) of multisets of eigenvalues.
  The factor \(\tau_{\Lbf,A}\) equals
  \[ \left( \left[ \frac{\sum_{a=1}^{n_i} x_{i,a}}{n_i}, \dots, \frac{\sum_{a=1}^{n_i} x_{i,a}}{n_i} \right] \right)_{i \in I}. \]
  The parabolic subgroups in \(\Pi(\Lbf, \tau_\Lbf)\) are the ones corresponding to the total orders on \(I\) satisfying \(i<j\) whenever \(\Re (\sum_a x_{i,a})/n_i > \Re (\sum_b x_{j,b})/n_j\).
\end{exam}

In Lemma \ref{lem:cmp_Langpos_dual} below we compare this definition with a notion that is more commonly used to express the Langlands classification.
First we need to recall a few facts about relative root systems.
Let \(\Pbf_0\) be a minimal parabolic subgroup of \(\Gbf\), and let \(\Abf_0\) be a maximal split torus in \(\Pbf_0\).
Let \(R(\Abf_0,\Gbf)\) be the set of roots of \(\Abf_0\) in \(\Gbf\).
By \cite[Exposé XXVI Théorème 7.4]{SGA3-III}\footnote{See also \cite[Corollaire 5.8]{BorelTits_gpesred}, although the proof seems to be incomplete in the non-reduced case.} there is a unique root datum (possibly non-reduced)
\[ (X^*(\Abf_0), R(\Abf_0,\Gbf), X_*(\Abf_0), R^\vee(\Abf_0,\Gbf)) \]
such that the associated Weyl group, seen as a group of automorphisms of \(\Abf_0\), is equal to the image of the normalizer of \(\Abf_0\) in \(\Gbf(\Q)\).
The parabolic subgroup \(\Pbf_0\) yields an order on the underlying root system and we denote by \(\Delta(\Abf_0, \Pbf_0) \subset R(\Abf_0, \Pbf_0)\) the set of simple roots.
Now consider a parabolic subgroup \(\Pbf\) of \(\Gbf\) which contains \(\Pbf_0\).
There is a unique Levi factor \(\Lbf\) of \(\Pbf\) which contains \(\Abf_0\), and \(\Abf_\Lbf\) is the centralizer of \(\Lbf\) in \(\Abf_0\).
As in Definition \ref{def:tau_A} let \(\Cbf_\Lbf\) be the largest quotient of \(\Lbf\) which is a split torus.
For \(\alpha \in \Delta(\Abf_0, \Pbf_0)\) not occurring in \(\Lbf\) (i.e.\ occurring in \(R_u(\Pbf)\)) with corresponding coroot \(\alpha^\vee \in R^\vee(\Abf_0, \Gbf) \subset X_*(\Abf_0)\) we denote by \(\alpha^\vee_{\Pbf}\) its image in \(X_*(\Cbf_\Lbf)\).

\begin{lemm} \label{lem:cmp_Langpos_dual}
  Let \(\Pbf\) be a parabolic subgroup of \(\Gbf\).
  Let \(\Lbf\) be a Levi factor of \(\Pbf\).
  Let \(\tau_\Lbf\) be a semi-simple conjugacy class in \(\C \otimes_{\Qbar} \lfrakhat\).
  By Definition \ref{def:tau_A} it yields \(\tau_{\Lbf,A} \in \C \otimes X^*(\Cbf_\Lbf)\).

  Let \(\Pbf_0\) be a minimal parabolic subgroup of \(\Pbf\).
  Let \(\Abf_0\) be a maximal split torus in \(\Pbf_0 \cap \Lbf\).
  We have \(\Pbf \in \Pi(\Lbf, \tau_\Lbf)\) if and only if for any simple root \(\alpha \in \Delta(\Abf_0, \Pbf_0)\) occurring in \(R_u(\Pbf)\) we have \(\langle \alpha^\vee_\Pbf, \Re \tau_{\Lbf,A} \rangle \geq 0\).
\end{lemm}
\begin{proof}
  See \cite[Lemma 3.8]{Taibi_notesIHES22}.
\end{proof}

Fix a maximal compact subgroup \(K_\infty^{\max}\) of \(\Gbf(\R)\).
Recall that \(\Gbf(\A)^1\) denotes the subgroup of \(g \in \Gbf(\A)\) satisfying \(|\chi(g)| = 1\) for all \(\chi \in X^*(\Gbf)^{\GalQ}\).
We have an isomorphism \(m_\Gbf: \Gbf(\A)/\Gbf(\A)^1 \simeq \Lie \Abf_\Gbf(\R)\) characterized by the relation
\[ \exp \langle \chi, m_\Gbf(g) \rangle = |\chi(g)| \]
for all \(\chi \in X^*(\Gbf)^{\GalQ}\).
We have \(\Gbf(\A) = \Gbf(\A)^1 \times \Abf_\Gbf(\R)^0\) and the restriction of \(\exp_{\Abf_\Gbf(\R)} \circ m_\Gbf\) to \(\Abf_\Gbf(\R)^0\) is the identity.
We denote by \(\Acal(\Gbf)\) the space of automorphic forms for \(\Gbf\), which are functions \(\Gbf(\Q) \backslash \Gbf(\A) \to \C\).
We have (see \cite[\S I.3.2]{MoeglinWaldspurger_bookspec}) an isomorphism of \((\gfrak,K_\infty^{\max},\Gbf(\A_f))\)-modules
\begin{align} \label{eq:decomp_Acal}
  \bigoplus_\nu \Acal(\Gbf)^{\Abf_\Gbf(\R)^0} \otimes \C(\nu) \otimes \Sym \afrak_\Gbf^* & \longrightarrow \Acal(\Gbf) \\
  f \otimes a \otimes P & \longmapsto a\nu \times (P \circ m_\Gbf) \times f \nonumber
\end{align}
where the sum ranges over \(\nu \in \C \otimes X^*(\Gbf)^{\GalQ}\) (seen as a character of \(\Gbf(\A)/\Gbf(\A)^1\) via \(s \otimes \chi \mapsto (g \mapsto |\chi(g)|^s)\)) and \(\Sym \afrak_\Gbf^*\) is the space of complex polynomial functions on \(\Lie \Abf_\Gbf(\R)\).
For a semi-simple conjugacy class \(\tau\) in \(\C \otimes_{\Qbar} \gfrakhat\), corresponding via Harish-Chandra's isomorphism to a maximal ideal \(\mfrak_\tau\) of \(Z(U(\C \otimes_\R \Lie \Gbf(\R)))\), we denote by \(\Acal(\Gbf)_\tau\) the subspace of automorphic forms which are killed by some power of \(\mfrak_\tau\).
Thus we have \(\Acal(\Gbf) = \bigoplus_\tau \Acal(\Gbf)_\tau\).
Note that \(\Acal(\Gbf)_\tau\) is contained in a single factor in \eqref{eq:decomp_Acal}, corresponding to the character \(\nu = \tau_A\) (Definition \ref{def:tau_A}).
Let \(\Acal^2(\Gbf)\) be the space of automorphic forms in \(\Acal(\Gbf)\) whose restriction to \(\Gbf(\Q) \backslash \Gbf(\A)^1\) is square-integrable.
Restricting \eqref{eq:decomp_Acal} we obtain an isomorphism of \((\gfrak,K_\infty^{\max},\Gbf(\A_f))\)-modules
\begin{equation} \label{eq:decomp_Acal2}
  \bigoplus_\nu \Acal^2(\Gbf)^{\Abf_\Gbf(\R)^0} \otimes \C(\nu) \otimes \Sym \afrak_\Gbf^* \simeq \Acal^2(\Gbf).
\end{equation}
Note that \(\Acal^2(\Gbf)^{\Abf_\Gbf(\R)^0}\) is the space of square-integrable automorphic forms on \(\Abf_\Gbf(\R)^0 \Gbf(\Q) \backslash \Gbf(\A)\), in particular it is semi-simple.
We also combine these two notations: \(\Acal^2(\Gbf)_\tau\) is the subspace of \(\Acal^2(\Gbf)\) consisting of forms killed by some power of \(\mfrak_\tau\).

For a parabolic subgroup \(\Pbf\) of \(\Gbf\) with Levi factor \(\Lbf\) we let \(K_{\Pbf,\infty}^{\max} := \Pbf(\R) \cap K_\infty^{\max}\) and we denote by \(\ind_\Pbf^\Gbf\) the parabolic induction functor from \((\pfrak, K_{\Pbf,\infty}^{\max}, \Pbf(\A_f))\)-modules to \((\gfrak, K_\infty^{\max}, \Gbf(\A_f))\)-modules, obtained by composing the smooth induction functor \(\ind_{\Pbf(\A_f)}^{\Gbf(\A_f)}\) and the induction functor \(\ind_{(\pfrak, K_{\Pbf,\infty}^{\max})}^{\gfrak, K_\infty^{\max}}\) (recalled \cite[p.208]{Franke}).
Let \(2 \rho_\Pbf \in X^*(\Cbf_\Lbf)\) be the determinant of the adjoint action of \(\Lbf\) on \(R_u(\Pbf)\).
We denote \(\Ind_\Pbf^\Gbf(-) := \ind_\Pbf^\Gbf(- \otimes |\rho_\Pbf|)\) for normalized induction.

\begin{theo}[Franke] \label{thm:Franke_fil}
  Let \(\tau\) be a \emph{regular} semi-simple conjugacy class in \(\C \otimes_{\Qbar} \gfrakhat\).
  There exists a finite, separated and exhaustive filtration of \(\Acal(\Gbf)_\tau\) by sub-\((\gfrak, K_\infty^{\max}, \Gbf(\A_f))\)-modules with associated graded pieces
  \[ \Ind_{\Pbf}^{\Gbf} \Acal^2(\Lbf)_{\tau_\Lbf} \]
  where
  \begin{itemize}
  \item \((\Lbf, \tau_\Lbf)\) ranges over \(\Gbf(\Q)\)-conjugacy classes of pairs consisting of a Levi subgroup \(\Lbf\) of \(\Gbf\) and a semi-simple \(\Lhat(\C)\)-conjugacy class \(\tau_\Lbf\) in \(\C \otimes_{\Qbar} \lfrakhat\) mapping to \(\tau\) via the differential of \(\iota_{\Lbf,\C}\) (for any \(\iota_\Lbf \in \Ecal(\Lbf, \Gbf)\)),
  \item \(\Pbf\) is an arbitrary element of \(\Pi(\Lbf, \tau_\Lbf)\) (we recall in Remark \ref{rem:indep_P_Pi_L_tauL} why the choice of \(\Pbf\) is irrelevant).
  \end{itemize}
\end{theo}
Note that we do not specify the order in which these occur as graded pieces, as we will not require this knowledge.
\begin{proof}
  This follows from \cite[Theorem 14]{Franke}.
  Since our notation differs substantially we explain in more detail.
  Franke restricts to automorphic forms invariant under \(\Abf_\Gbf(\R)^0\), using the decomposition \eqref{eq:decomp_Acal} it is easy to reduce to this case.
  Our maximal ideal \(\mfrak_\tau\) of \(Z(U(\C \otimes_\R \Lie \Gbf(\R)))\) corresponds to Franke's \(\mathcal{J}\).
  Franke's result is more refined than Theorem \ref{thm:Franke_fil} in a number of ways.
  \begin{itemize}
  \item In \cite[Theorem 14]{Franke} a certain parameter \(\tau\) in a cone occurs, corresponding to a certain growth condition imposed on automorphic forms (Franke's spaces \(\mathfrak{Fin}_\mathcal{J} S_{\rho_{-\tau}+\log}(\dots)\)), taking it to be very far in this cone (or taking the union over all such parameters \(\tau\)) simply yields all automorphic forms.
  \item We have imposed that \(\tau\) be regular in order to simplify the statement: in general the terms occurring on the left-hand side in \cite[Theorem 14]{Franke} are colimits over certain groupoids, but as observed in \cite[Theorem 19 I]{Franke} for regular infinitesimal characters these groupoids are simply equivalent to sets.
  \item Franke fixes a conjugacy class of Levi subgroups of \(\Gbf\) (in the terminology of his paper, a class of associate parabolic subgroups, denoted \(\{\mathcal{P}\}\) by Franke, see \cite[p.201]{Franke}) and restricts to forms whose cuspidal support corresponds to this class.
    We simply sum over these classes.
  \end{itemize}

  The groupoids \(\mathcal{M}^{k,T,i}_{\mathcal{J},\{\mathcal{P}\},\tau}\) appearing in \cite[Theorem 14]{Franke} are defined at the bottom of p.233, and we now translate between this definition and our formulation.
  As explained above Franke's \(\tau\) is irrelevant for us (we take \(\tau=\infty\)).
  Franke's function \(T\) and integer \(i\) index the filtration, we do not need to make this precise.
  His parabolic subgroup \(\mathcal{R}\) is our \(\Pbf\).
  His continuous character \(\Lambda\) (of what we denote \(\Abf_\Lbf(\A) / \Abf_\Lbf(\Q) \Abf_\Gbf(\R)^0\)) is given by two independent pieces: the continuous character \(\lambda_t\) (of what we denote \(\Abf_\Lbf(\R)^0 / \Abf_\Gbf(\R)^0\)) and the unitary character \(\widetilde{\Lambda}\).
  We do not refine automorphic forms in \(\Acal(\Lbf)^{\Abf_\Lbf(\R)^0}\) by central character \(\widetilde{\Lambda}\), so one has to group the terms corresponding to various \(\widetilde{\Lambda}\) in \cite[Theorem 14]{Franke} to obtain the statement in Theorem \ref{thm:Franke_fil}.
  Franke's \(\lambda_t\) and \(\chi\) together correspond to our \(\tau_\Lbf\), in fact his \(\lambda_t\) corresponds to our component \(\tau_{\Lbf,A}\) of \(\tau_\Lbf\) (Definition \ref{def:tau_A}).
  His condition ``\(\lambda_t \in \mathrm{supp}_{u_t} \mathcal{J}\)'' on p.234 (see also p.230 for the definition of \(\mathrm{supp}\)) is equivalent to our condition that \(\tau_\Lbf\) maps to \(\tau\).
  Finally Franke's condition \(\Re \lambda_t \in \ol{\check{\afrak}_{\Rcal}^+}\) at the bottom of p.233 is equivalent to our condition \(\Pbf \in \Pi(\Lbf, \tau_\Lbf)\) (Definition \ref{def:Pi_L_tauL}) by Lemma \ref{lem:cmp_Langpos_dual}.
  %Recall the definition (p.186): \(\Rcal\) contains a fixed minimal parabolic subgroup \(\Pcal_0\), for a simple root \(\alpha \in \check{\afrak}_0\) occuring in the unipotent radical of \(\Rcal\) (denote by \(\Delta_\Rcal\) the set of such simple roots), let \(\check{\alpha}_\Rcal\) be the projection to \(\afrak_{\Rcal}\) of the corresponding coroot \(\check{\alpha} \in \afrak_0\), then \(\ol{\check{\afrak}_{\Rcal}^+}\) is the set of \(x \in \check{\afrak}_{\Rcal}\) such that \(\langle x, \check{\alpha}_\Rcal \rangle \geq 0\) for any \(\alpha \in \Delta_\Rcal\).

  Finally the \((\gfrak, K_\infty^{\max}, \Gbf(\A_f))\)-modules \(M(t)\) appearing in \cite[Theorem 14]{Franke} are defined on p.234 as \(W(u_t) \otimes D_t\).
  The space \(W(u_t)\) is defined on p.218 as (in our notation) \(\Ind_\Pbf^\Gbf V(u_t)\) for a certain subspace \(V(u_t)\) of \(\Acal^2(\Lbf)^{\Abf_\Lbf(\R)^0}\).
  The space \(D_t\) is precisely \(\C(\nu) \otimes \Sym \afrak_\Lbf^*\) for \(\nu = \tau_{\Lbf,A}\) in \eqref{eq:decomp_Acal2} (for \(\Lbf\) instead of \(\Gbf\)).
\end{proof}

\begin{rema} \label{rem:indep_P_Pi_L_tauL}
  In the setting of Theorem \ref{thm:Franke_fil}, for a pair \((\Lbf, \tau_\Lbf)\) and \(\Pbf, \Pbf' \in \Pi(\Lbf, \tau_\Lbf)\) we have an isomorphism
  \[ \Ind_\Pbf^\Gbf \Acal^2(\Lbf)_{\tau_\Lbf} \simeq \Ind_{\Pbf'}^\Gbf \Acal^2(\Lbf)_{\tau_\Lbf} \]
  given by the standard intertwining operator (defined by meromorphic continuation) which is holomorphic \cite[Lemma 2 p.234]{Franke}, which is why the choice of \(\Pbf \in \Pi(\Lbf, \tau_\Lbf)\) is irrelevant.
  This is perhaps clearer using the parabolic subgroup \(\Qbf = \Mbf \Nbf\) associated to \(\tau_\Lbf\) introduced in Remark \ref{rem:direct_def_Pi_L_tauL}: denoting \(\Pbf_\Mbf = \Pbf \cap \Mbf\) we have \(\Ind_\Pbf^\Gbf \simeq \Ind_\Qbf^\Gbf \circ \Ind_{\Pbf_\Mbf}^\Mbf\) and similarly for \(\Pbf'\), and we have \(\tau_{\Lbf,A} \in \R \otimes X^*(\Cbf_\Mbf) \oplus i\R \otimes X^*(\Cbf_\Lbf)\).
\end{rema}

\begin{rema} \label{rem:L_tauL_fix_L_first}
  In Theorem  \ref{thm:Franke_fil} instead of considering classes of pairs \((\Lbf, \tau_\Lbf)\) we could just as well choose a Levi \(\Lbf\) in each \(\Gbf(\Q)\)-conjugacy class and for each such Levi consider orbits of \(\Lhat(\C)\)-conjugacy classes \(\tau_\Lbf\) under the normalizer of \(\Lbf\) in \(\Gbf(\Q)\).
\end{rema}

\subsubsection{Euler characteristic of \((\gfrak,K_\infty)\)-cohomology of automorphic forms}

Following \cite[\S 7.4]{Franke} we deduce from Theorem \ref{thm:Franke_fil} a formula for the Euler-characteristic of \((\gfrak,K_\infty)\)-cohomology.
Let \(K_\infty\) be an open subgroup of \(K_\infty^{\max}\).
For a parabolic subgroup \(\Pbf\) of \(\Gbf\) we denote \(K_{\Pbf,\infty} = \Pbf(\R) \cap K_\infty\), and similarly for \(K_\infty^{\max}\).
Let \(V\) be an irreducible finite-dimensional \((\gfrak,K_\infty^{\max})\)-module.
Let \(\tau_V\) be the infinitesimal character of \(V\).
By Wigner's lemma\footnote{More precisely we use a straightforward generalization to generalized eigenspaces (instead of eigenspaces) for the action of \(Z(U(\C \otimes_\R \Lie \Gbf(\R)))\).} \cite[Corollary I.4.2]{BoWa} the natural map
\[ H^\bullet((\gfrak,K_\infty), \Acal(\Gbf)_{-\tau_V} \otimes V) \longrightarrow H^\bullet((\gfrak,K_\infty), \Acal(\Gbf) \otimes V) \]
is an isomorphism.
The following lemma provides a further simplification of these cohomology spaces but we defer using it until later.

\begin{lemm} \label{lem:HS_afrak_G}
  We have an isomorphism
  \[ H^\bullet((\gfrak,K_\infty), \Acal(\Gbf)_{-\tau_V} \otimes V) \simeq H^\bullet((\gfrak/\afrak_\Gbf,K_\infty), (\Acal(\Gbf)_{-\tau_V} \otimes V)^{\Abf_\Gbf(\R)^0}). \]
\end{lemm}
\begin{proof}
  This follows from the Hochschild-Serre spectral sequence (see \cite[Theorem I.6.5]{BoWa} and \cite[Appendix A]{Franke}) associated to the ideal \(\afrak_\Gbf\) of \(\gfrak\) and the fact (explained at the bottom of p.256 of \cite{Franke}) that \(H^i(\afrak_\Gbf, \Sym \afrak_\Gbf^*)\) vanishes for \(i>0\).
  Note that the Chevalley-Eilenberg complex computing these cohomology groups may be identified with the (algebraic) de Rham complex of the complex affine space associated to \(\afrak_\Gbf\), so this vanishing may be interpreted as the algebraic Poincaré lemma.
\end{proof}

Now Theorem \ref{thm:Franke_fil} for \(\tau = -\tau_V\) implies
\[ e((\gfrak,K_\infty), \Acal(\Gbf) \otimes V) = \sum_{[\Lbf, \tau_\Lbf]} e \left( (\gfrak,K_\infty), \Ind_\Pbf^\Gbf \Acal^2(\Lbf)_{\tau_\Lbf} \otimes V \right). \]
Separating parabolic induction at the real and finite places we write
\begin{align*}
  & H^\bullet((\gfrak,K_\infty), \Ind_\Pbf^{\Gbf} \Acal^2(\Lbf)_{\tau_\Lbf} \otimes V) \\
  =& \Ind_{\Pbf(\A_f)}^{\Gbf(\A_f)} H^\bullet((\gfrak,K_\infty), \ind_{(\pfrak, K_{\Pbf,\infty}^{\max})}^{(\gfrak,K_\infty^{\max})}( \Acal^2(\Lbf)_{\tau_\Lbf} \otimes V \otimes |\rho_\Pbf|_\infty))
\end{align*}
where \(V\) on the second line is really \(\res_{(\pfrak, K_{\Pbf,\infty}^{\max})}^{(\gfrak,K_\infty^{\max})} V\).
Recall that for a \((\pfrak,K_{\Pbf,\infty}^{\max})\)-module \(M\) the parabolically induced \((\gfrak,K_\infty^{\max})\)-module \(\ind_{(\pfrak,K_{\Pbf,\infty}^{\max})}^{(\gfrak,K_\infty^{\max})} M\) is defined as the largest sub-\(U(\gfrak)\)-module of \(\Hom_{U(\pfrak)}(U(\gfrak), M)\) (where \(U(\gfrak)\) is seen as a left \(U(\pfrak)\)-module and a right \(U(\gfrak)\)-module by multiplication) on which the action of \(U(\kfrak_\infty)\) is locally finite and integrates into an action of \(K_\infty^0\) (it then extends to an action of \(K_\infty^{\max}\) because \(\pi_0(K_{\Pbf,\infty}^{\max}) \to \pi_0(K_\infty^{\max})\) is surjective).
We see that there is an isomorphism of functors
\[ \res_{(\gfrak, K_\infty^{\max})}^{(\gfrak,K_\infty)} \circ \ind_{(\pfrak, K_{\Pbf,\infty}^{\max})}^{(\gfrak,K_\infty^{\max})} \simeq \ind_{(\pfrak, K_{\Pbf,\infty})}^{(\gfrak,K_\infty)} \circ \res_{(\pfrak, K_{\Pbf,\infty})}^{(\pfrak,K_{\Pbf,\infty}^{\max})} \]
and we deduce
\[ H^\bullet((\gfrak,K_\infty), \ind_{(\pfrak, K_{\Pbf,\infty}^{\max})}^{(\gfrak,K_\infty^{\max})}( \Acal^2(\Lbf)_{\tau_\Lbf} \otimes V \otimes |\rho_\Pbf|_\infty)) \simeq H^\bullet((\gfrak,K_\infty), \ind_{(\pfrak, K_{\Pbf,\infty})}^{(\gfrak,K_\infty)}( \Acal^2(\Lbf)_{\tau_\Lbf} \otimes V \otimes |\rho_\Pbf|_\infty)). \]
Staring at the definition also makes evident that \(\ind_{(\pfrak, K_{\Pbf,\infty})}^{(\gfrak,K_\infty)}\) is right adjoint to \(\res_{(\pfrak, K_{\Pbf,\infty})}^{(\gfrak,K_\infty)}\), which implies
\begin{equation} \label{eq:gKinfty_coh_adj}
  H^\bullet((\gfrak,K_\infty), \ind_{(\pfrak, K_{\Pbf,\infty})}^{(\gfrak,K_\infty)}( \Acal^2(\Lbf)_{\tau_\Lbf} \otimes V \otimes |\rho_\Pbf|_\infty)) \simeq H^\bullet((\pfrak,K_{\Pbf,\infty}), \Acal^2(\Lbf)_{\tau_\Lbf} \otimes V \otimes |\rho_\Pbf|_\infty).
\end{equation}
Denoting \(\Nbf := R_u(\Pbf)\), we use the Hochschild-Serre spectral sequence for \(\nfrak \subset \pfrak\), but since we have a semi-direct product \(\pfrak \simeq \lfrak \rtimes \nfrak\) we can be more precise.
The cohomology groups on the right-hand side of \eqref{eq:gKinfty_coh_adj} may be computed with Chevalley-Eilenberg complexes \(C^\bullet_\mathrm{CE}((\pfrak, K_{\Pbf,\infty}), -)\) and a simple computation identifies \eqref{eq:gKinfty_coh_adj} with the cohomology of the total complex associated to the double complex
\[ C^\bullet_\mathrm{CE}((\lfrak, K_{\Pbf,\infty}), \Acal^2(\Lbf)_{\tau_{\Lbf}} \otimes |\rho_\Pbf|_\infty \otimes C^\bullet_\mathrm{CE}(\nfrak, V)). \]
Here we have used the identification of \(\Lbf_\R\) with the Levi factor \(\Pbf_\R \cap \theta(\Pbf_\R)\) of \(\Pbf_\R\), where \(\theta\) is the Cartan involution of \(\Gbf_\R\) satisfying \(K_\infty^{\max} = \Gbf(\R)^\theta\).
Now assume temporarily that \(V\) is (the restriction of) an irreducible algebraic representation of \(\Gbf_\C\).
The proof of Kostant's theorem (already recalled before Corollary \ref{coro:IH_vs_Hc}) computing \(H^\bullet(\nfrak, V)\) shows that \(C^\bullet_\mathrm{CE}(\nfrak, V)\) is quasi-isomorphic to a direct sum of complexes concentrated in one degree, so we simply have
\[ H^i((\gfrak,K_\infty), \ind_{(\pfrak, K_{\Pbf,\infty})}^{(\gfrak,K_\infty)}( \Acal^2(\Lbf)_{\tau_\Lbf} \otimes V \otimes |\rho_\Pbf|_\infty)) \simeq \bigoplus_{a+b=i} H^a((\lfrak, K_{\Pbf,\infty}), \Acal^2(\Lbf)_{\tau_\Lbf} \otimes |\rho_\Pbf|_\infty \otimes H^b(\nfrak, V)). \]
It turns out that there is only one degree \(b\) for which \(H^b(\nfrak, V)\) can be non-zero, to see this we briefly recall Kostant's theorem.
Choose a Borel subgroup \(\Bbf_\Lbf\) of \(\Lbf_{\Qbar}\) and denote by \(\Bbf\) the Borel subgroup of \(\Pbf_{\Qbar}\) satisfying \(\Bbf \cap \Lbf_{\Qbar} = \Bbf_\Lbf\).
Let \(\Tbf\) be a maximal torus of \(\Bbf_\Lbf\).
The irreducible algebraic representation \(V\) of \(\Gbf_\C\) is parametrized by its dominant (for \(\Bbf\)) weight \(\lambda \in X^*(\Tbf)\).
As usual denote
\[ \rho = \rho_\Bbf = \frac{1}{2} \sum_{\alpha \in R(\Tbf,\Bbf)} \alpha \in \frac{1}{2} X^*(\Tbf). \]
The infinitesimal character \(\tau_V\) of \(V\) is represented by \(\lambda + \rho\).
A Kostant representative (for \(\Pbf\), \(\Lbf\) and \((\Bbf_\Lbf,\Tbf)\)) is \(w \in W(\Tbf, \Gbf_{\Qbar})\) such that \(w(\lambda+\rho) \in \frac{1}{2} X^*(\Tbf)\) is dominant for \(\Bbf_\Lbf\).
The set of Kostant representatives is in bijection with \(W(\Tbf, \Lbf_{\Qbar}) \backslash W(\Tbf, \Gbf_{\Qbar})\).
Kostant's theorem says
\[ H^b(\nfrak, V) \simeq \bigoplus_{l(w) = b} V^\Lbf_{w(\lambda+\rho)-\rho} \]
where the sum ranges over Kostant representatives of length \(b\) and \(V^\Lbf_{\lambda'}\) is the irreducible algebraic representation of \(\Lbf_\C\) with highest weight \(\lambda'\).
\begin{defi} \label{def:def_W_V_P_tauL}
  In the situation above there is a unique Kostant representative \(w\) satisfying \(-w^{-1}(\tau_\Lbf) = \tau_V\), that we denote \(w_{\Pbf,\tau_\Lbf}\).
  Define \(W_{V,\Pbf,\tau_\Lbf}\) as the piece of the semi-simple \((\lfrak, K_{\Pbf,\infty})\)-module \(H^{l(w_{\Pbf,\tau_\Lbf})}(\nfrak,V) \otimes |\rho_\Pbf|_\infty\) having infinitesimal character \(-\tau_\Lbf\).
\end{defi}
More generally it can be useful to start from an irreducible finite-dimensional \((\gfrak,K_\infty)\) module \(V\) of the form \(V^\Gbf_\lambda \otimes \chi\) for some \(\chi \in \C \otimes X^*(\Gbf)^{\Gal(\C/\R)}\) (seen as a character of \(\Gbf(\R)\)).
With the same definition we thus have
\[ W_{V,\Pbf,\tau_\Lbf} \simeq V^\Lbf_{w_{\Pbf,\tau_\Lbf}(\lambda+\rho)-\rho} \otimes |\rho_\Pbf|_\infty \otimes \chi. \]
(The general case of an arbitrary irreducible finite-dimensional \((\gfrak,K_\infty)\)-module \(V\) can be reduced to this case using a z-extension and Lemma \ref{lem:fd_rep_real_gp_alg}.)
Using Wigner's lemma and Lemma \ref{lem:HS_afrak_G} (for \(\Lbf\) instead of \(\Gbf\)) we conclude
\begin{align*}
  &\ H^i((\gfrak,K_\infty), \ind_{(\pfrak, K_{\Pbf,\infty})}^{(\gfrak,K_\infty)}( \Acal^2(\Lbf)_{\tau_\Lbf} \otimes V \otimes |\rho_\Pbf|_\infty)) \\
  \simeq &\ H^{i-l(w_{\Pbf,\tau_\Lbf})}((\lfrak, K_{\Pbf,\infty}), \Acal^2(\Lbf)_{\tau_\Lbf} \otimes W_{V,\Pbf,\tau_\Lbf}) \\
  \simeq &\ H^{i-l(w_{\Pbf,\tau_\Lbf})}((\lfrak / \afrak_\Lbf, K_{\Pbf,\infty}), (\Acal^2(\Lbf)_{\tau_\Lbf} \otimes W_{V,\Pbf,\tau_\Lbf})^{\Abf_\Lbf(\R)^0}).
\end{align*}
Recall that the subspace of \(\Acal^2(\Lbf)\) on which \(\Abf_\Lbf(\R)^0\) acts by the inverse of the central character of \(W_{V,\Pbf,\Lbf}\) is semi-simple.
Returning to Euler characteristics we deduce
\begin{equation} \label{eq:final_simpl_e_piece_fil}
  e \left( (\gfrak,K_\infty), \Ind_\Pbf^\Gbf \Acal^2(\Lbf)_{\tau_\Lbf} \otimes V \right) = \epsilon(w_{\Pbf,\tau_\Lbf}) \Ind_{\Pbf(\A_f)}^{\Gbf(\A_f)} e \left( (\lfrak / \afrak_\Lbf, K_{\Pbf,\infty}), (\Acal^2(\Lbf)_{\tau_\Lbf} \otimes W_{V,\Pbf,\tau_\Lbf})^{\Abf_\Lbf(\R)^0} \right)
\end{equation}
where \(\epsilon(w) = (-1)^{l(w)}\).
As explained on \cite[p.266]{Franke} the Euler characteristic vanishes if \(\Lbf\) is not \(\R\)-cuspidal, i.e.\ if \((\Lbf/\Abf_\Lbf)_\R\) does not admit an anisotropic maximal torus.
Using the notation introduced in Definition \ref{def:e_2_G_Kinfty_V} the right-hand side of \eqref{eq:final_simpl_e_piece_fil} is
\[ \Ind_{\Pbf(\A_f)}^{\Gbf(\A_f)} e_{(2)}(\Lbf, K_{\Pbf,\infty}, W_{V,\Pbf,\tau_\Lbf}). \]

We state in the following corollary what we have deduced from Franke's Theorem \ref{thm:Franke_fil} following \cite[\S 7.7]{Franke}.

\begin{coro}[Franke] \label{cor:Franke_Euler1}
  Let \(\Gbf\) be a connected reductive group over \(\Q\).
  Let \(K_\infty\) be an open subgroup of a maximal compact subgroup of \(\Gbf(\R)\).
  Let \(V\) be an irreducible finite-dimensional \((\gfrak,K_\infty^{\max})\)-module, and denote by \(\tau_V\) its infinitesimal character.
  The Euler characteristic \(e((\gfrak,K_\infty), \Acal(\Gbf) \otimes V)\) is equal to
  \[ \sum_{[\Lbf,\tau_\Lbf]} \epsilon(w_{\Pbf,\tau_\Lbf})  \Ind_{\Pbf(\A_f)}^{\Gbf(\A_f)} e_{(2)}(\Lbf, K_{\Pbf,\infty}, W_{V,\Pbf,\tau_\Lbf}) \]
  where
  \begin{itemize}
  \item the sum ranges over \(\Gbf(\Q)\)-conjugacy classes of pairs \((\Lbf,\tau_\Lbf)\) where \(\Lbf\) is an \emph{\(\R\)-cuspidal} Levi subgroup of \(\Gbf\) and \(\tau_\Lbf\) is a semi-simple conjugacy class in \(\lfrakhat_\C\) mapping to \(-\tau_V\),
  \item \(\Pbf\) is an arbitrary element of \(\Pi(\Lbf, \tau_\Lbf)\),
  \item \(w_{\Pbf,\tau_\Lbf}\) and \(W_{V,\Pbf,\tau_\Lbf}\) are defined in Definition \ref{def:def_W_V_P_tauL},
  \item \(e_{(2)}(\Lbf, K_{\Pbf,\infty}, -)\) is defined in Definition \ref{def:e_2_G_Kinfty_V}.
  \end{itemize}
\end{coro}

It follows from Remark \ref{rem:indep_P_Pi_L_tauL} and the computation above that choosing another \(\Pbf \in \Pi(\Lbf, \tau_\Lbf)\) does not change the term in the sum, but this is less obvious now.
In the next section we make this more evident.

\subsubsection{On coefficient systems for Levi subgroups}

\begin{exam} \label{ex:coeff_syst_GL3}
  Consider the case where \(\Gbf = \GLbf_3\), \(V\) is the irreducible algebraic representation of \(\Gbf_\C\) with infinitesimal character \(\tau_V = (k,0,-k)\) for some integer \(k>0\), \(\Lbf\) is the block diagonal Levi subgroup \(\GLbf_1 \times \GLbf_2\) of \(\Gbf\) and \(\tau_\Lbf = (0, {k,-k})\).
  We have \(\tau_{\Lbf,A}=0\), so any parabolic subgroup of \(\Gbf\) with Levi factor \(\Lbf\) belongs to \(\Pi(\Lbf, \tau_\Lbf)\).
  We will use Borel subgroups of \(\Gbf\) containing the diagonal maximal torus \(\Tbf\).
  We choose the upper triangular Borel subgroup \(\Bbf_\Lbf\) of \(\Lbf\).
  For \(\Pbf \in \Pi(\Lbf, \tau_\Lbf)\) the upper block triangular parabolic subgroup the corresponding Borel subgroup \(\Bbf\) of \(\Gbf\) is the upper triangular subgroup (i.e.\ \(\rho_\Bbf = (1,0,-1)\)) and we have \(\lambda = (k-1,0,-k+1)\) and
  \[ w_{\Pbf,\Lbf}(\lambda_\Bbf+\rho_\Bbf)-\rho_\Bbf = (0,k,-k) - (1,0,-1) = (-1,k,-k+1) \]
  and so \(W_{V,\Pbf,\tau_\Lbf} = V^\Lbf_{(-1,(k,-k+1))} \otimes |\rho_\Pbf|\).
  For the lower block triangular parabolic subgroup \(\Pbf' \in \Pi(\Lbf, \tau_\Lbf)\) we have \(\rho_{\Bbf'}=(-1,1,0)\) and thus \(\lambda_{\Bbf'}=(-k+1,k-1,0)\) and we compute
  \[ w_{\Pbf,\Lbf}(\lambda_{\Bbf'}+\rho_{\Bbf'})-\rho_{\Bbf'} = (0,k,-k) - (-1,1,0) = (-1,k-1,-k)\]
  and so
  \[ W_{V,\Pbf',\tau_\Lbf} = V^\Lbf_{(-1,(k-1,-k))} \otimes |\rho_{\Pbf'}| \simeq W_{V,\Pbf,\tau_\Lbf} \otimes (1, \sign \det). \]
  In particular the representations \(W_{V,\Pbf,\tau_\Lbf}\) and \(W_{V,\Pbf',\tau_\Lbf}\) are \emph{not} isomorphic.
  Lemma \ref{lem:Lparam_Levi_welldef} below implies \(l(w_{\Pbf,\tau_\Lbf}) = l(w_{\Pbf',\tau_\Lbf})\).
  We can check directly that for any open subgroup \(K_{\Lbf,\infty}\) of a maximal compact subgroup of \(\Lbf(\R)\) we have
  \[ H^\bullet(\Lbf, K_{\Lbf,\infty}, W_{V,\Pbf,\tau_\Lbf}) \simeq H^\bullet(\Lbf, K_{\Lbf,\infty}, W_{V,\Pbf',\tau_\Lbf}) \]
  because \(\tau_\Lbf\) is very regular (on the factor \(\GLbf_2\) of \(\Lbf\)) and so any \((\lfrak,K_{\Lbf,\infty}^{\max})\)-module occurring in \(\Acal^2(\Lbf,\tau_{\Lbf,A})_{\tau_\Lbf}\) is part of the discrete series, so Theorem \ref{thm:coh_DS} shows that its \((\lfrak,K_{\Lbf,\infty})\)-cohomologies relative to \(W_{V,\Pbf,\tau_\Lbf}\) and \(W_{V,\Pbf',\tau_\Lbf}\) are equal.
\end{exam}

\begin{rema}
  Example \ref{ex:coeff_syst_GL3} seems to contradict the first part of \cite[\S 7.7 Lemma 1]{Franke}.
  The issue seems to be the statement on the first line of p.272: \((\gfrak,K_\infty^{\max})\)-modules with infinitesimal character equal to that of the trivial representation are parametrized not just by \(\Gbf(\R)\)-orbits of pairs \((H,\Delta^+)\) (where \(H = \Tbf(\R)\) with (\(\Tbf\) a maximal torus of \(\Gbf_\R\) and \(\Delta^+\) an order on \(R(\Tbf_\C,\Gbf_\C)\)), but by triples consisting in addition of a character of \(\pi_0(H)\) (see \cite[Theorem 2.2.4]{Vogan_book}).
\end{rema}

\begin{lemm} \label{lem:Lparam_Levi_welldef}
  Assume we are in the setting of Corollary \ref{cor:Franke_Euler1}: \(V\) is an irreducible finite-dimensional \((\gfrak,K_\infty^{\max})\)-module, \(\Lbf\) is an \(\R\)-cuspidal Levi subgroup of \(\Gbf\) and \(\tau_\Lbf\) a semisimple conjugacy class in \(\lfrakhat_\C\) mapping to \(-\tau_V\).
  Then neither the length \(l(w_{\Pbf,\tau_\Lbf})\) nor the essentially discrete Langlands parameter \(\varphi_{V,\Pbf,\tau_\Lbf}: W_\R \to {}^L \Lbf(\C)\) corresponding to \(W_{V,\Pbf,\tau_\Lbf}\) (i.e.\ with corresponding L-packet the set of discrete series representations of \(\Lbf(\R)\) having infinitesimal character \(\tau_\Lbf\) and central character equal to the inverse of that of \(W_{V,\Pbf,\tau_\Lbf}\)) depend on the choice of \(\Pbf \in \Pi(\Lbf, \tau_\Lbf)\).
\end{lemm}
\begin{proof}
  We first consider the length \(l(w_{\Pbf,\tau_\Lbf})\).
  Choose a distinguished splitting \(\zeta_\Gbf = (\Bcal,\Tcal,(X_\alpha)_\alpha,s)\) for \({}^L \Gbf\) and a distinguished splitting \(\zeta_\Lbf\) for \({}^L \Lbf\).
  This gives us an embedding \(\iota_\Lbf := \iota[\Pbf, \zeta_\Gbf, \zeta_\Lbf] \in \Ecal(\Lbf, \Gbf)\) and we denote \(\Lcal := \iota_\Lbf({}^L \Lbf)\) and \(\Pcal := \Lcal \Bcal\).
  Let \(\Ucal\) be the unipotent radical of \(\Pcal\).
  As in Definition \ref{def:Pi_L_tauL} we denote by \(\tau_\Lcal\) the image of \(\tau_\Lbf\) by the differential of \(\iota_\Lbf\), and recall the parabolic subgroup \(\Qcal = \Mcal \Ncal \supset \Pcal\) of \({}^L \Gbf\).
  We identify \(\tau_V\) with an element of \(\C \otimes X_*(\Tcal)\) which is strictly dominant for \(\Bcal\).
  There is a unique representative \(\tau_{\Lcal} \in \C \otimes X_*(\Tcal)\) (in its \(\Lcal^0(\C)\)-conjugacy class) which is stricly anti-dominant for \(\Bcal \cap \Lcal\).
  Using \((\Bbf,\Tbf)\) and \((\Bcal,\Tcal)\) to identify based root data we obtain an identification of the Weyl group of \(\Tbf\) in \(\Gbf_{\C}\) with the Weyl group of \(\Tcal\) in \(\Ghat\), and our Kostant representative \(w = w_{\Pbf,\tau_\Lbf}\) is determined by \(w(\tau_V) = -\tau_{\Lcal}\).
  The length \(l(w)\) is equal to the number of roots \(\alpha \in R(\Tcal, \Bcal)\) satisfying \(\langle w(\tau_V), \alpha \rangle < 0\).
  Recall from Definition \ref{def:tau_A} the decomposition \(\tau_\Lcal = \tau_{\Lcal,A} + \tau_{\Lcal,0}\).
  For any \(\alpha \in R(\Tcal, \Ghat)\) the (a priori complex) pairings \(\langle \tau_{\Lcal,A}, \alpha \rangle\) and \(\langle \tau_{\Lcal,0}, \alpha \rangle\) are both real.
  Let \(w_{0, \Lcal}\) be the longest element of the Weyl group \(W(\Tcal, \Lcal^0)\) (for \(\Bcal \cap \Lcal^0\)).
  Let \(j \in W_{\R}\) be any element of \(W_{\R} \smallsetminus W_{\C}\) satisfying \(j^2 = -1\).
  Then conjugation by \(x := w_{0, \Lcal} s(j)\) is an involution of \(\Tcal\) which leaves \(R(\Tcal, \Ucal)\) invariant.
  It maps \(\tau_{\Lcal, 0}\) to \(-\tau_{\Lcal, 0}\) because \(\Lbf\) is \(\R\)-cuspidal.
  So for \(\alpha \in R(\Tcal, \Ghat)\) we have
  \begin{align*}
    \langle \tau_{\Lcal}, \alpha \rangle &= \langle \tau_{\Lcal,A}, \alpha \rangle + \langle \tau_{\Lcal,0}, \alpha \rangle \\
    \langle \tau_{\Lcal}, x(\alpha) \rangle &= \langle \tau_{\Lcal,A}, \alpha \rangle - \langle \tau_{\Lcal,0}, \alpha \rangle
  \end{align*}
  By regularity of \(\tau_\Lcal\) we have \(|\langle \tau_{\Lcal,A}, \alpha \rangle| \neq |\langle \tau_{\Lcal,0}, \alpha \rangle|\).
  For \(\alpha \in R(\Tcal,\Ucal)\), which satisfies \(\langle \tau_{\Lcal,A}, \alpha \rangle \geq 0\), we distinguish two cases:
  \begin{itemize}
  \item if \(\langle \tau_{\Lcal,A}, \alpha \rangle < |\langle \tau_{\Lcal,0}, \alpha \rangle|\) then in particular \(\langle \tau_{\Lcal,0}, \alpha \rangle \neq 0\) and so \(x(\alpha) \neq \alpha\).
    Up to swapping \(\alpha\) and \(x(\alpha)\) we have \(\langle \tau_{\Lcal}, \alpha \rangle > 0\) and \(\langle \tau_{\Lcal}, x(\alpha) \rangle < 0\).
  \item if \(\langle \tau_{\Lcal,A}, \alpha \rangle > |\langle \tau_{\Lcal,0}, \alpha \rangle|\) then \(\langle \tau_{\Lcal}, \alpha \rangle > 0\).
  \end{itemize}
  We conclude
  \begin{equation} \label{eq:l_w_indep_P}
    l(w) = \frac{\dim \Ghat - \dim \Lcal^0}{2} - \frac{1}{4} \left| \left\{ \alpha \in R(\Tcal,\Ghat) \smallsetminus R(\Tcal,\Lcal^0) \,\middle|\, |\langle \tau_{\Lcal,A}, \alpha \rangle| < |\langle \tau_{\Lcal,0}, \alpha \rangle| \right\} \right|
  \end{equation}
  which clearly does not depend on the choice of \(\Pcal\).

  We now consider the Langlands parameter \(\varphi_{V, \Pbf, \tau_{\Lbf}}\).
  Using a z-extension of \(\Gbf_\R\) we can reduce to the case where \(V = V_\lambda\) is (the restriction of) an irreducible algebraic representation of \(\Gbf_\C\).
  The infinitesimal character of \(\varphi_{V,\Pbf,\tau_\Lbf}\) is clearly \(\tau_{\Lbf}\), so we are left to check that the restriction of\footnote{Unfortunately we have to mix additive and multiplicative notation here, we hope no confusion will arise.} \((w(\lambda + \rho) - \rho) \otimes |\rho_{\Pbf}|\) to \(\Zbf_{\Lbf}^0(\R)\) does not depend on this choice.
  The restriction to \(\Zbf_{\Lbf}(\R)^0\) is imposed by the infinitesimal character, and the finite \(2\)-torsion group of torsion elements in \(\Abf_{\Lbf}(\R) \simeq (\R^{\times})^{\dim \Abf_{\Lbf}}\) surjects onto \(\pi_0(\Zbf_{\Lbf}^0(\R))\), so it is enough to check that the image of \(w(\lambda + \rho) - \rho\) in \(X^*(\Abf_{\Lbf}) / 2X^*(\Abf_{\Lbf})\) does not depend on the choice of \(\Pbf\).
  We have in \(X^*(\Abf_\Lbf)\)
  \[ (w(\lambda + \rho) - \rho)|_{\Abf_{\Lbf}} = (-\tau_{\Lbf} - \rho_{\Pbf}) |_{\Abf_{\Lbf}} \]
  so we are left to show that the image of \(\rho_{\Pbf}\) in \(\frac{1}{2} X^*(\Abf_{\Lbf}) / 2X^*(\Abf_{\Lbf})\) does not depend on the choice of \(\Pbf \in \Pi(\Lbf, \tau_\Lbf)\).
  This is similar to the previous proof, but for a change we do not argue on the dual side.
  Let \(\Tbf_{\mathrm{an}}\) be a maximal torus of \(\Lbf_\R\) such that \(\Tbf_{\mathrm{an}} / \Abf_{\Lbf}\) is anisotropic.
  Denoting \(\{1, \sigma \} = \Gal(\C/\R)\) the action of \(\sigma\) on \(R(\Tbf_{\mathrm{an}, \C}, \Gbf_\C)\) (which corresponds to the action of \(x\) on the dual side considered above) preserves \(R(\Tbf_{\mathrm{an}, \C}, \Ubf_\C)\) where \(\Ubf\) is the unipotent radical of \(\Pbf\), in particular for \(\alpha \in R(\Tbf_{\mathrm{an}, \C}, \Gbf_\C) \smallsetminus R(\Tbf_{\mathrm{an}, \C}, \Lbf_\C)\) we have \(\sigma(\alpha) \neq -\alpha\).
  For \(\alpha \in R(\Tbf_{\mathrm{an}, \C}, \Ubf_\C)\) we distinguish two cases.
  \begin{itemize}
  \item If \(\sigma(\alpha) \neq \alpha\) then since \(\sigma(\alpha)|_{\Abf_\Lbf} = \alpha|_{\Abf_\Lbf}\)the contribution of \(\alpha\) and \(\sigma(\alpha)\) to the restriction of
    \[ \rho_\Pbf = \frac{1}{2} \sum_{\beta \in R(\Tbf_{\mathrm{an}, \C}, \Ubf_\C)} \beta \]
    is simply \(\alpha|_{\Abf_\Lbf}\), which is equal \(\mod 2X^*(\Abf_\Lbf)\) to its opposite.
  \item If \(\sigma(\alpha) = \alpha\) then we have \(\langle \alpha^\vee, \tau_{\Lbf,0} \rangle = 0\) and thus \(\langle \alpha^\vee, \tau_{\Lbf} \rangle = \langle \alpha^\vee, \tau_{\Lbf,A} \rangle > 0\).
  \end{itemize}
  So if we partition the \(\{\pm \id \} \times \Gal(\C/\R)\)-orbits in \(R(\Tbf_{\mathrm{an}, \C}, \Gbf_\C) \smallsetminus R(\Tbf_{\mathrm{an},\C}, \Lbf_\C)\) as \(\Ocal_1 \sqcup \Ocal_2\) where \(\Ocal_1\) is the set of orbits \([\alpha]\) satisfying \(\sigma(\alpha) \neq \alpha\) (i.e.\ orbits with \(4\) elements, so \(\Ocal_2\) is the set of orbits with two elements \(\{ \pm \alpha\}\)) then we have
  \begin{equation} \label{eq:res_rhoP_mod2_indep_P}
    \rho_\Pbf|_{\Abf_\Lbf} = \sum_{[\alpha] \in \Ocal_1} \alpha|_{\Abf_\Lbf} + \frac{1}{2} \sum_{\substack{[\alpha] \in \Ocal_2 \\ \langle \alpha^\vee, \tau_\Lbf \rangle > 0}} \alpha|_{\Abf_\Lbf} \mod 2X^*(\Abf_\Lbf)
  \end{equation}
  which clearly does not depend on the choice of \(\Pbf\).
\end{proof}

The lemma allows us to unambiguously define a sign \(\epsilon_{\tau_\Lbf} := \epsilon(w_{\Pbf,\tau_\Lbf})\) and an essentially discrete Langlands parameter (up to conjugation by \(\Lhat(\C)\)) \(\varphi_{V,\Lbf,\tau_\Lbf} := \varphi_{V,\Pbf,\tau_\Lbf}\), where \(\Pbf\) is any element of \(\Pi(\Lbf, \tau_\Lbf)\).
Note that the proof of Lemma \ref{lem:Lparam_Levi_welldef} gives us a relatively simple way of computing these two objects.

We require two small remarks before reformulating Corollary \ref{cor:Franke_Euler1}.
\begin{itemize}
\item For \(K_{\infty}\) an open subgroup of a maximal compact subgroup \(K_{\infty}^{\max}\) of \(\Gbf(\R)\), \(\Pbf\) a parabolic subgroup of \(\Gbf\) and \(\Lbf\) a Levi factor of \(\Pbf\), denote by \(K_{\Lbf, \infty}\) the image of \(K_{\infty} \cap \Pbf(\R)\) in \(\Lbf(\R)\) (realized as a quotient of \(\Pbf(\R)\)).
  This subgroup of \(\Lbf(\R)\) actually depends on the choice of \(\Pbf\), but its \(\Lbf(\R)\)-conjugacy class does not: there exists \(g \in \Gbf(\R) = \Pbf(\R) K_{\infty}^{\max}\) such that \(g K_{\infty}^{\max} g^{-1}\) contains a maximal compact subgroup of \(\Lbf(\R)\), and we deduce that the \(\Lbf(\R)\)-conjugacy class of \(K_{\infty, \Lbf}\) only depends on the \(\Gbf(\R)\)-conjugacy class of \(K_{\infty}\) and on \(\Lbf\) (not on \(\Pbf\)).
\item For a Levi subgroup \(\Lbf\) of \(\Gbf\) and \(M \in K_0^{\Tr}(\Rep_\C^{\adm}(\Lbf(\A_f)))\) (see Definition \ref{def:Groth_adm}) the element \(\Ind_{\Pbf(\A_f)}^{\Gbf(\A_f)} M\) of \(K_0^{\Tr}(\Rep_\C^{\adm}(\Gbf(\A_f)))\) does not depend on the choice of a parabolic subgroup \(\Pbf\) of \(\Gbf\) admitting \(\Lbf\) as a Levi factor.
  We denote this element of \(K_0^{\Tr}(\Rep_\C^{\adm}(\Gbf(\A_f)))\) by \(\Ind_{\Lbf(\A_f)}^{\Gbf(\A_f)} M\).
\end{itemize}

\begin{coro} \label{cor:Franke_Euler2}
  Let \(\Gbf\) be a connected reductive group over \(\Q\).
  Let \(K_{\infty}\) be an open subgroup of a maximal compact subgroup \(K_\infty^{\max}\) of \(\Gbf(\R)\), \(V\) an irreducible finite-dimensional \((\gfrak,K_\infty)\)-module.
  Then we have the equality in \(K_0^{\Tr}(\Rep^{\adm}_{\C}(\Gbf(\A_f))\)
  \begin{equation} \label{eq:Franke_Euler2}
    e(\Gbf, K_{\infty}, V) = \sum_{[\Lbf, \tau_\Lbf]} \epsilon_{\tau_{\Lbf}} \Ind_{\Lbf(\A_f)}^{\Gbf(\A_f)} e_{(2)}(\Lbf, K_{\Lbf, \infty}, \varphi_{V, \Lbf, \tau_{\Lbf}})
  \end{equation}
  where the sum is over \(\Gbf(\Q)\)-conjugacy classes of pairs \((\Lbf, \tau_{\Lbf})\) with \(\Lbf\) an \(\R\)-cuspidal Levi subgroup of \(\Gbf\) and \(\tau_{\Lbf}\) a semisimple conjugacy class in \(\lfrakhat_\C\) mapping to \(-\tau_V\).
\end{coro}

\begin{rema} \label{rem:Franke_Euler_uniform_lambda}
  The proof of Lemma \ref{lem:Lparam_Levi_welldef} actually shows a bit more than what we stated, and shows that Formula \eqref{eq:Franke_Euler2} is ``uniform in \(V\)'' in a sense to be made precise below.
  Let \((\Bbf,\Tbf)\) be a Borel pair in \(\Gbf_\C\).
  For an irreducible algebraic representation \(V\) parametrized by a dominant (for \(\Bbf\)) weight \(\lambda \in X^*(\Tbf)\) we let \(\tau_V = \lambda + \rho_\Bbf\).
  Fix a representative \(\Lbf\) of a \(\Gbf(\Q)\)-conjugacy class of Levi subgroups of \(\Gbf\).
  The set of \(\Lhat(\C)\)-conjugacy classes \(\tau_\Lbf\) mapping to \(\tau_V\) is parametrized by
  \[ \Nbf(\Lbf,\Gbf)(\C) \backslash \{ g \in \Gbf(\C) \,|\, \Ad(g) \Tbf \subset \Lbf_\C \} \]
  where \(\Nbf(\Lbf,\Gbf)\) denotes the normalizer of \(\Lbf\) in \(\Gbf\).
  We fix a class in this quotient, and we even fix a class in
  \begin{equation} \label{eq:left_cosets_L_transp}
    \Lbf(\C) \backslash \{ g \in \Gbf(\C) \,|\, \Ad(g) \Tbf \subset \Lbf_\C \}
  \end{equation}
  mapping to this class.
  Now fix a Borel pair \((\Bbf_\Lbf,\Tbf_\Lbf)\) in \(\Lbf_\C\).
  In our chosen class in \eqref{eq:left_cosets_L_transp} there is a unique left coset \(\Tbf_\Lbf(\C) g\) such that \(\Ad(g) \Tbf = \Tbf_\Lbf\) and such that
  \[ \Ad(g^{-1})^*: X^*(\Tbf) \simeq X^*(\Tbf_\Lbf) \]
  maps the Weyl chamber
  \[ C := \{ x \in \R \otimes X^*(\Tbf) \,|\, \forall \alpha \in R(\Tbf,\Bbf),\, \langle \alpha^\vee, x \rangle > 0\} \]
  to the Weyl chamber \(C_\Lbf\) for \((\Bbf_\Lbf, \Tbf_\Lbf)\).
  Denote \(\tau_\Lbf = -\Ad(g^{-1})^* \tau_V\), which we will see alternatively as a function of \(\tau_V \in C\) or a function of the dominant weight \(\lambda\).
  For \(\Pbf \in \Pi(\Lbf, \tau_\Lbf)\) we then have \eqref{eq:l_w_indep_P}
  \[ l(w_{\Pbf,\tau_\Lbf}) = \frac{\dim \Gbf - \dim \Lbf}{2} - \frac{1}{4} \left| \left\{ \alpha \in R(\Tbf_\Lbf,\Gbf_\C) \smallsetminus R(\Tbf_\Lbf,\Lbf_\C) \,\middle|\, |\langle \alpha^\vee, \tau_{\Lbf,A} \rangle| < |\langle \alpha^\vee \tau_{\Lbf,0} \rangle| \right\} \right| \]
  and this set of roots does not depend on the choice of a dominant weight \(\lambda\): in fact it makes sense for any \(\tau_V \in C\) and the continuous function
  \begin{align*}
    C & \longrightarrow \R \\
    \tau_V & \longmapsto |\langle \alpha^\vee, \tau_{\Lbf,A} \rangle| - |\langle \alpha^\vee \tau_{\Lbf,0} \rangle|
  \end{align*}
  does not vanish because we have \(\langle \alpha^\vee, \tau_\Lbf \rangle \neq 0\) for all \(\alpha \in R(\Tbf_\Lbf, \Gbf_\C)\).
  Thus \(l(w_{\Pbf,\tau_\Lbf})\), which does not depend on the choice of \(\Pbf \in \Pi(\Lbf, \tau_\Lbf)\) does not depend on \(\lambda\) either.
  We now consider the central characters of the local systems \(W_{V,\Pbf,\tau_\Lbf}\).
  For \(\Pbf \in \Pi(\Lbf, \tau_\Lbf)\) we have
  \[ W_{V,\Pbf,\tau_\Lbf} \simeq V^{\Lbf}_{-\tau_\Lbf - \rho_{\Bbf_\Lbf} - \rho_\Pbf} \otimes |\rho_\Pbf|_\infty. \]
  Define an involution \(\sigma\) of \(X^*(\Tbf_\Lbf)\) by requiring that it acts by \(-\id\) on \(X^*(\Tbf_\Lbf/\Zbf_\Lbf^0)\) and as the complex conjugation on the quotient \(X^*(\Zbf_\Lbf^0)\).
  (Note that if \(\Tbf_\Lbf = \Tbf_{\mathrm{an},\C}\) where \(\Tbf_\mathrm{an}\) is a maximal torus of \(\Lbf_\R\) which is anisotropic modulo center then this action is the natural one.)
  As in the proof of Lemma \ref{lem:Lparam_Levi_welldef} we have a decomposition \(\Ocal_1 \sqcup \Ocal_2\) of the set of \(\{\pm \id\} \times \{1, \sigma \}\)-orbits of \(R(\Tbf_\Lbf,\Gbf_\C) \smallsetminus R(\Tbf_\Lbf,\Lbf_\C)\).
  By \eqref{eq:res_rhoP_mod2_indep_P} we have
  \[ \rho_\Pbf|_{\Abf_\Lbf} = \sum_{[\alpha] \in \Ocal_1} \alpha|_{\Abf_\Lbf} + \frac{1}{2} \sum_{\substack{[\alpha] \in \Ocal_2 \\ \langle \alpha^\vee, \tau_\Lbf \rangle > 0}} \alpha|_{\Abf_\Lbf} \mod 2X^*(\Abf_\Lbf). \]
  and the index sets, which we argued do not depend on the choice of \(\Pbf\) in the proof of Lemma \ref{lem:Lparam_Levi_welldef}, do not depend on the choice of a dominant weight \(\lambda\) either.
  Let \(\delta \in \frac{1}{2} X^*(\Abf_\Lbf)\) be any representative of this class (e.g.\ obtained by choosing a representative in each orbit in \(\Ocal_1\)).
  Then \(varphi_{V,\tau_\Lbf}: W_\R \to {}^L \Lbf(\C)\) is the discrete Langlands parameter with infinitesimal character \(\tau_\Lbf\) such that composing with \({}^L \Lbf \to {}^L \Abf_\Lbf\) yields the parameter of the character \((\tau_\Lbf|_{\Abf_\Lbf} + \delta) \otimes |\delta|^{-1}\) of \(\Abf_\Lbf(\R)\) (note that \(\tau_\Lbf|_{\Abf_\Lbf} + \delta\) belongs to \(X^*(\Abf_\Lbf)\)).
\end{rema}

\subsubsection{Example: \(\GLbf_n\)}

Let us make the formula in Corollary \ref{cor:Franke_Euler2} more explicit for \(\Gbf = \GLbf_{n,\Q}\).
We actually deduce it from the earlier Corollary \ref{cor:Franke_Euler1}, but the formula exemplifies the irrelevance of \(\Pbf\) proved for Corollary \ref{cor:Franke_Euler2}, as well as the ``uniformity in \(V\)'' explained in Remark \ref{rem:Franke_Euler_uniform_lambda}.

\begin{coro} \label{coro:Euler_ord_GL}
  Let \(n \geq 1\).
  For \(a,b \in \Z_{\geq 0}\) such that \(a+2b=n\) denote by \(\Lbf_{a,b} \simeq \GLbf_1^a \times \GLbf_2^b\) the corresponding standard Levi subgroup of \(\GLbf_n\), and let \(\mathfrak{S}(a,b)\) be the subset of \(\mathfrak{S}_n\) consisting of \(\sigma\) such that
  \begin{enumerate}
    \item \(\sigma^{-1}(1) < \dots < \sigma^{-1}(a)\),
    \item \(\sigma^{-1}(a+1) < \sigma^{-1}(a+2)\), \dots, \(\sigma^{-1}(a+2b-1) < \sigma^{-1}(a+2b)\),
    \item \(\sigma^{-1}(a+1) < \sigma^{-1}(a+3) < \dots < \sigma^{-1}(a+2b-1)\).
  \end{enumerate}
  Consider a dominant weight \(\lambda = (\lambda_1 \geq \dots \geq \lambda_n)\) for \(\GLbf_n\) and let
  \[ \tau = (\tau_1 > \dots > \tau_n) := \lambda + \rho \]
  so that \(\tau_i = \lambda_i + \frac{n+1}{2} - i\).
  For \(\sigma \in \mathfrak{S}_n\) denote \(\sigma(\tau)_i = \tau_{\sigma^{-1}(i)}\) and \((\sigma \cdot \lambda)_i = \lambda_{\sigma^{-1}(i)} - \sigma^{-1}(i) + i\), i.e.\ \(\sigma \cdot \lambda = \sigma(\tau) - \rho\).
  Using notation introduced in Example \ref{exam:Euler_L2_coh_GL} we have
  \begin{multline*}
    e(\mathbf{GL}_n, V_{\lambda}) = \sum_{a+2b=n} (-1)^{a(a-1)/2} \sum_{\sigma \in \mathfrak{S}(a,b)} \epsilon(\sigma) \\
    \Ind_{\Lbf_{a,b}(\A_f)}^{\GLbf_n(\A_f)} \Big( \bigotimes_{i=1}^a e(\GLbf_1, (\sigma \cdot \lambda)_i + a + 1) |\cdot|_f^{(\sigma \cdot \lambda)_i + a + 1 - \sigma(\tau)_i} \\
      \otimes \bigotimes_{i=1}^b e_{(2)}(\GLbf_2, \sigma(\tau)_{a+2i-1}-1/2, \sigma(\tau)_{a+2i}+1/2) \Big).
  \end{multline*}
\end{coro}

\begin{proof}
  Conjugacy classes of Levi subgroups of \(\GLbf_n\) are parametrized by unordered partitions of \(n\), \(n = n_1 + \dots + n_r\) (\(n_i \in \Z_{\geq 1}\)) corresponding to Levi subgroups isomorphic to \(\GLbf_{n_1} \times \dots \times \GLbf_{n_r}\) (see Example \ref{ex:Levi_para_dual_GL}).
  Such a Levi subgroup is \(\R\)-cuspidal if and only if \(n_i \in \{1,2\}\) for all \(i\).
  So in \eqref{eq:Franke_Euler2} \(\GLbf_n(\Q)\)-conjugacy classes of pairs \((\Lbf, \tau_{\Lbf})\) are parametrized by partitions \(n = a + 2b\) together with an unordered partition of \(\tau\) (equivalently, of \(\{1, \dots, n\}\)) into \(a\) singletons and \(b\) pairs.
  Note that \(\mathfrak{S}(a,b)\) parametrizes unordered partitions of \(\{1, \dots, n\}\) into \(a\) singletons and \(b\) pairs.
  Fix a pair \((\Lbf, \tau_{\Lbf})\) as in \eqref{eq:Franke_Euler2}, and let \(\sigma \in \mathfrak{S}(a,b)\) be the unique element such that \(\Lbf\) can be conjugated to \(\Lbf_{a,b}\), identifying \(\tau_{\Lbf}\) to
  \[ \left( -\tau_{\sigma^{-1}(1)}, \dots, -\tau_{\sigma^{-1}(a)}, \{-\tau_{\sigma^{-1}(a+1)}, -\tau_{\sigma^{-1}(a+2)}\}, \dots, \{ -\tau_{\sigma^{-1}(n-1)}, -\tau_{\sigma^{-1}(n)}\} \right) . \]
  It remains to compute \(\epsilon_{\tau_\Lbf}\) and \(\varphi_{V,\Lbf,\tau_{\Lbf}}\) in \eqref{eq:Franke_Euler2}.
  By Example \ref{exam:Euler_L2_coh_GL} on \(\GLbf_2\) factors of \(\Lbf\) the
  parameter \(\varphi_{V,\Lbf,\tau_{\Lbf}}\) is determined by \(\tau_{\Lbf}\), so we
  only need to consider \(\GLbf_1\) factors.

  Choose a parabolic subgroup \(\Pbf\) as explained before Lemma
  \ref{lem:Lparam_Levi_welldef}.
  Up to conjugation by \(\GLbf_n(\Q)\), we may assume that \(\Pbf\) and \(\Lbf\) are
  standard (with \(\Lbf \neq \Lbf_{a,b}\) in general).
  There is a unique enumeration
  \[ \{ x_1, \dots, x_r \} = \{ \{\sigma^{-1}(1)\}, \dots, \{\sigma^{-1}(a)\},
  \{ \sigma^{-1}(a+1), \sigma^{-1}(a+2) \}, \dots, \{ \sigma^{-1}(n-1),
  \sigma^{-1}(n) \} \} \]
  such that, denoting \(c_k = |x_k|\) and \(\tau_{x_k} = \{ \tau_i | i \in x_k \}\),
  we have that \(\Lbf\) is the block diagonal \(\GLbf_{c_1} \times \dots \times
  \GLbf_{c_r}\) and \(\tau_{\Lbf} = (-\tau_{x_1}, \dots, -\tau_{x_r})\).
  By definition of \(\Pbf\), denoting
  \[ s(\tau_{x_k}) = \begin{cases}
    \tau_i & \text{ if } x_k = \{i\}, \\
    (\tau_i + \tau_j)/2 & \text{ if } x_k = \{i, j\} \text{ with } i \neq j, \\
  \end{cases} \]
  we have \(-s(\tau_{x_1}) \geq \dots \geq -s(\tau_{x_r})\).
  Let \(w \in \mathfrak{S}_n\) be the Kostant representative for the pair \((\Pbf,
  \Lbf)\) mapping \(\tau\) to \(-\tau_{\Lbf}\).
  Then \(w = w_2 w_1 \sigma\) where \(w_1\) stabilizes \(\{1, \dots, a \}\) and
  satisfies \(w_1(a+2i) = w_1(a+2i-1)+1\) for \(1 \leq i \leq b\) (in particular
  \(w_1\) permutes the pairs \((a+1,a+2)\), \dots, \((n-1,n)\)), and \(w_2\) moves some
  pairs between the singletons, i.e.\ \(w_2(1) < \dots < w_2(a)\), \(w_2(a+2i) =
  w_2(a+2i-1)+1\) for \(1 \leq i \leq b\) and \(w_2(a+2i) < w_2(a+2i+1)\) for \(1 \leq
  i < b\).
  By definition of \(\Pbf\) we have \(w_1(i) = a+1-i\) for \(1 \leq i \leq a\), and so
  the first \(a\) entries of \(w_1 \sigma \cdot \lambda\) are
  \[ (\lambda_{\sigma^{-1}(a+1-i)} - \sigma^{-1}(a+1-i) + i)_{1 \leq i \leq a} =
    \left( (\sigma \cdot \lambda)_{a+1-i} - a - 1 + 2i \right)_{1 \leq i \leq a}
  \]
  Now the set of indices corresponding to the \(\GLbf_1\) factors of \(\Lbf\) is
  \(w_2(\{1, \dots, a\})\) and we have \(\rho_{w_2(i)} - \rho_i \in 2\Z\) for any \(1
  \leq i \leq a\), so that \((w \cdot \lambda)_{w_2(i)} - (w_1 \sigma \cdot
  \lambda)_i \in 2 \Z\) for any \(1 \leq i \leq a\).
  Finally it is easy to see that \(\epsilon(w_1) = (-1)^{a(a-1)/2}\) and
  \(\epsilon(w_2) = 1\).
\end{proof}

\begin{exam}
  For \(n \leq 7\) the dimension of all cohomology groups \(H^i(\GLbf_n(\Z),\Q)\) are known \cite{LeeSzczarba} \cite{ElbazVincentGanglSoule} (except for \(n=4\) for which only \(H^i(\SLbf_4(\Z),\Q)\) is explicitly computed in \cite{LeeSzczarba}).
  Let us check that our formula specialized at \(\lambda=0\) agrees with these computations.
  The Euler characteristic \(e(\GLbf_n,V_0)^{\GLbf_n(\Zhat)}\) belongs to the Grothendieck group of finite-dimensional representations of the unramified Hecke algebra \(\Hcal_f^{\unr}(\GLbf_n)\), and we compute the multiplicity \(e_n\) of the ``trivial'' character (the one corresponding to the trivial representation of \(\GLbf_n(\A_f)\)) in this Euler characteristic.
  Thanks to \cite{JacSha_Euler2} this corresponds to retaining only the power of \(|\cdot|_f\) in the formula \eqref{eq:e_GL1} for \(e(\GLbf_1,-)\) and discarding \(S_k\) from the formula \eqref{eq:e2_GL2} for \(e_{(2)}(\GLbf_2,-)\).
  For \(n<12\) we have \(S_k=0\) for all \(k \leq n\) and we deduce \(e(\GLbf_n(\Z),\C) = e_n\).
  We see that the only terms (corresponding to \(\sigma \in \Sfrak(a,b)\)) contributing in Corollary \ref{coro:Euler_ord_GL} are the ones satisfying \(\sigma(\tau)_{a+2i-1} = \sigma(\tau)_{a+2i}+1\), i.e.
  \begin{equation} \label{eq:cond_sigma_e_GLn_triv}
    \sigma^{-1}(a+2i) = \sigma^{-1}(a+2i-1)+1 \text{ for all } 1 \leq i \leq b.
  \end{equation}
  This implies \(\epsilon(\sigma) = +1\) and \(i - \sigma^{-1}(i)\) even for all \(1 \leq i \leq a\).
  We have \((\sigma \cdot \lambda)_i+a+1 = i - \sigma^{-1}(i) + a + 1 = n+1 \mod 2\).

  First consider the case where \(n\) is even.
  The contribution of \(\sigma\) vanishes if \(a>0\), so we have in fact a single contribution corresponding to \(a=0\) and \(\sigma = \id\).
  We conclude \(e_n = 1\).
  For \(n\) even we may also consider the multiplicity \(\tilde{e}_n\) of the trivial character of \(\Hcal_f^{\unr}(\GLbf_n)\) to \(e(\GLbf_n,\tilde{\C})^{\GLbf_n(\Zhat)}\) where \(\tilde{\C}\) is \(\C\) endowed with the character \(\sign (\det)\) of \(\GLbf_n(\R)\).
  By Shapiro's lemma we have \(e_n + \tilde{e}_n = e(\SLbf_n(\Z),\C)\).
  Taking e.g.\ \(\lambda = (1, \dots, 1)\) allows us to compute \(\tilde{e}_n\) using Corollary \ref{coro:Euler_ord_GL}.
  We still have the condition \eqref{eq:cond_sigma_e_GLn_triv} but now \((\sigma \cdot \lambda)_i+a+1 = 0 \mod 2\) for any \(1 \leq i \leq a\) and so we find
  \[ \tilde{e}_n = \sum_{a+2b=n} (-1)^{a/2} \sum_\sigma 1 \]
  where the sum is over all \(\sigma \in \Sfrak(a,b)\) satisfying \eqref{eq:cond_sigma_e_GLn_triv}.
  The number of such permutations is easily computed, and equal to \(\binom{n-b}{b}\).
  Define for an integer \(m \geq 0\) (with the usual convention \(\binom{x}{y}=0\) for \(y<0\))
  \[ f(m) = \sum_{b \geq 0} (-1)^b \binom{m-b}{b}. \]
  We have \(f(0)=f(1)=1\) and for \(m>0\)
  \[ f(m+1) = \sum_{b \geq 0} (-1)^b \left( \binom{m-b}{b} + \binom{m-b}{b-1} \right) = f(m)-f(m-1) \]
  and we deduce
  \[ f(m) =
    \begin{cases}
      1 & \text{ if } m = 0 \mod 6 \\
      1 & \text{ if } m = 1 \mod 6 \\
      0 & \text{ if } m = 2 \mod 6 \\
      -1 & \text{ if } m = 3 \mod 6 \\
      -1 & \text{ if } m = 4 \mod 6 \\
      0 & \text{ if } m = 5 \mod 6
    \end{cases}
  \]
  and for even \(n>0\)
  \[ \tilde{e}_n = (-1)^{n/2} f(n) =
    \begin{cases}
      1 & \text{ if } m = 0 \mod 12 \\
      0 & \text{ if } m = 2 \mod 12 \\
      -1 & \text{ if } m = 4 \mod 12 \\
      -1 & \text{ if } m = 6 \mod 12 \\
      0 & \text{ if } m = 8 \mod 12 \\
      1 & \text{ if } m = 10 \mod 12.
    \end{cases}
  \]
  This is consistent with
  \[ H^i(\GLbf_2(\Z),\Q) = H^i(\SLbf_2(\Z),\Q) \simeq
    \begin{cases}
      \Q & \text{if } i=0, \\
      0 & \text{if } i>0,
    \end{cases}
  \]
  with the computation in \cite{LeeSzczarba}
  \[ H^i(\SLbf_4(\Z),\Q) \simeq
    \begin{cases}
      \Q & \text{if } i \in \{0,3\} \\
      0 & \text{otherwise}.
    \end{cases}
  \]
  and with the computation of \(H^\bullet(\GLbf_6(\Z),\Q)\) and \(H^\bullet(\SLbf_6(\Z),\Q)\) in \cite[\S 7.3]{ElbazVincentGanglSoule}.
  We also deduce that \(\GLbf_4(\Z)/\SLbf_4(\Z)\) acts by \(-\id\) on the line \(H^3(\SLbf_4(\Z),\Q)\).
  Presumably the method in \cite{LeeSzczarba} could be used to compute \(H^\bullet(\GLbf_4(\Z),\Q)\) as well.

  Now consider the case where \(n\) is odd.
  For any pair \((a,b)\) satisfying \(a+2b=n\) and any \(\sigma \in \Sfrak(a,b)\) satisfying \eqref{eq:cond_sigma_e_GLn_triv} we have \((\sigma \cdot \lambda)_i+a+1\) even for all \(1 \leq i \leq a\), so we find
  \[ e_n = \sum_{a+2b=n} (-1)^{a(a-1)/2} \sum_\sigma 1 = (-1)^{(n-1)/2} f(n) =
    \begin{cases}
      1 & \text{ if } m = 1 \mod 12 \\
      1 & \text{ if } m = 3 \mod 12 \\
      0 & \text{ if } m = 5 \mod 12 \\
      -1 & \text{ if } m = 7 \mod 12 \\
      -1 & \text{ if } m = 9 \mod 12 \\
      0 & \text{ if } m = 11 \mod 12.
    \end{cases}
  \]
  This is consistent with \cite{Soule_SL3Z}
  \[ H^i(\SLbf_3(\Z),\Q) \simeq
    \begin{cases}
      \Q & \text{if } i=0, \\
      0 & \text{if } i>0,
    \end{cases}
  \]
  and with the computation of \(H^\bullet(\SLbf_n(\Z),\Q)\) for \(n \in \{5,7\}\) in \cite[\S 7.3]{ElbazVincentGanglSoule}.
\end{exam}

\begin{coro} \label{coro:Euler_c_GL}
  For \(n \geq 1\) and \(\lambda = (\lambda_1 \geq \dots \geq \lambda_n)\), using notation as in Corollary \ref{coro:Euler_ord_GL} we have
  \begin{multline*}
    e_c(\GLbf_n, V_{\lambda}) = (-1)^{n(n+1)/2-1} \sum_{a+2b=n} (-1)^{a(a-1)/2} \sum_{\sigma \in \mathfrak{S}(a,b)} \epsilon(\sigma) \\
    \Ind_{\Lbf_{a,b}(\A_f)}^{\GLbf_n(\A_f)} \Big( \left( \bigotimes_{i=1}^a e(\GLbf_1, (\sigma \cdot \lambda)_i)|\cdot|_f^{(\sigma \cdot \lambda)_i - \sigma(\tau)_i} \right) \\
    \otimes \bigotimes_{i=1}^b e_{(2)}(\GLbf_2, \sigma(\tau)_{a+2i-1}-1/2, \sigma(\tau)_{a+2i}+1/2) \Big).
  \end{multline*}
\end{coro}
\begin{proof}
  This follows from the preceding Corollary using Example \ref{exam:orientation_GL}, writing the character \(\op{sign}(\det)^{n-1}\) of \(\GLbf_n(\R)\) as \(\det^{n-1} / |\det|^{n-1}\), using Remark \ref{rema:twisting_top_coh} and the duality between \(\Ind_{\Pbf_{a,b}(\A_f)}^{\GLbf_n(\A_f)} \pi\) and \(\Ind_{\Pbf_{a,b}(\A_f)}^{\GLbf_n(\A_f)} \pi^\vee\) for any admissible representation \(\pi\) of \(\Lbf_{a,b}(\A_f)\).
\end{proof}

\subsection{Intersection in terms of compactly supported cohomology}
\label{sec:IH_from_Hc}

In this section we essentially plug Corollary \ref{coro:Euler_c_GL} into Corollary \ref{coro:IH_vs_Hc_Q} and simplify the resulting expansion, ultimately obtaining Theorem \ref{theo:IH_vs_Hc_simple}.
For calculations we need to choose a basis of \(X^*(\Tbf_{\GSpbf_{2n}})\), our convention consists of writing \(\lambda = (\lambda_1, \dots, \lambda_n, m)\) when \(\lambda\) maps \((t_1, \dots, t_n, s) \in \Tbf_{\GSpbf_{2n}}\) (notation as in Section \ref{sec:not_red_gps}) to
\[ s^m \prod_{i=1}^n t_i^{\lambda_i}. \]
The character \(\lambda\) is dominant if and only if we have \(\lambda_1 \geq \dots \geq \lambda_n \geq 0\).
Composing \(\lambda\) with the cocharacter
\[ x \in \GLbf_1 \mapsto x I_{2n} = (x^{-1}, \dots, x^{-1}, x^2) \in \Tbf_{\GSpbf_{2n}} \]
gives \(2m - \sum_{i=1}^n \lambda_i\).
The Weyl group \(W(\GSpbf_{2n})\) of \(\Tbf_{\GSpbf_{2n}}\) in \(\GSpbf_{2n}\) is identified as usual with \(\{\pm 1\}^n \rtimes \Sfrak_n\):
\begin{itemize}
\item for \(w = \sigma \in \Sfrak_n\) and \(\lambda = (\lambda_1, \dots, \lambda_n, m) \in X^*(\Tbf_{\GSpbf_{2n}})\) we have
  \[ w(\lambda) = (\lambda_{\sigma^{-1}(1)}, \dots, \lambda_{\sigma^{-1}(n)}, m), \]
\item for \(w = (\epsilon_i)_{1 \leq i \leq n} \in \{\pm 1\}^n\) we have
  \[ w(\lambda) = \left( \epsilon_1 \lambda_1, \dots, \epsilon_n \lambda_n, m - \sum_{\substack{1 \leq i \leq n \\ \epsilon_i = -1}} \lambda_i \right). \]
\end{itemize}
We will also identify \(W(\GSpbf_{2n})\) with the group of permutations \(w\) of \(\{\pm 1, \dots, \pm n\}\) satisfying \(w(-i) = -w(i)\) (elements of the subgroup \(\Sfrak_n\) of \(W(\GSpbf_{2n})\) are the ones preserving \(\{1, \dots, n\}\)).
Finally for \(0 \leq h \leq n\) the isomorphism
\begin{align*}
  \GLbf_1^h \times \Tbf_{\GSpbf_{2(n-h)}} & \longrightarrow \Tbf_{\GSpbf_{2n}} \\
  ((t_1, \dots, t_h), (t_{h+1}, \dots, t_n, s)) & \longmapsto (t_1, \dots, t_n, s)
\end{align*}
identifies \(\lambda\) with \(((\lambda_1, \dots, \lambda_h), (\lambda_{h+1}, \dots, \lambda_n, m))\), which will also be denoted by \((\lambda_{\lin}, \lambda_{\her})\).

For the rest of this section we fix a dominant weight \(\lambda\) for \(\GSpbf_{2n}\) and let \(V_{\lambda}\) be an irreducible algebraic representation with highest weight \(\lambda\).
Let
\[ \tau := \lambda + \rho = (\lambda_1 +n, \dots, \lambda_n + 1, m + n(n+1)/4). \]
Our first step consists of grouping terms in Corollary \ref{coro:IH_vs_Hc_Q} corresponding to the same value for \(n_r\).
To this end we introduce some notation.
For a positive integer \(h\) we take the diagonal torus and upper triangular Borel subgroup to define standard parabolic and Levi subgroups of \(\GLbf_h\), and as usual we identify the Weyl group of \(\GLbf_h\) with the symmetric group \(\Sfrak_h\).
For an integer \(r \geq 1\) and integers \(0 < n_1 < \dots < n_{r-1} < h\), denoting \(n_r = h\) and \(\ul{n} = (n_1, \dots, n_r)\), let \(\Lbf_{\ul{n}}\) be the standard Levi subgroup \(\GLbf_{n_1} \times \dots \times \GLbf_{n_r - n_{r-1}}\) of \(\GLbf_h\) and let \(\Qbf_{\ul{n}}\) be the corresponding standard parabolic subgroup of \(\GLbf_h\).
For \(\tau'' \in \R^h\) let \(\Sfrak_h(\ul{n}, \tau'') \subset \Sfrak_h\) be the set of Kostant representatives \(\sigma\) for the standard Levi subgroup \(\Lbf_{\ul{n}}\) of \(\GLbf_h\) which satisfy \(\sum_{i=1}^{n_j} \sigma(\tau'')_i > 0\) for all \(1 \leq j \leq r\).
For \(\lambda' \in \Z^h\) a dominant weight for \(\GLbf_h\) and \(\tau'' \in \R^h\) define \(C(\lambda', \tau'')  \in K_0^{\Tr}(\Rep_{\Q}^{\adm}(\GLbf_h(\A_f)))\) by
\begin{equation} \label{eq:collect_nr_IH_from_Hc}
  C(\lambda', \tau'') = \sum_{\substack{1 \leq r \leq h \\ 0 < n_1 < \dots < n_{r-1} < h}} \sum_{\sigma \in \Sfrak_h(\ul{n}, \tau'')} \epsilon(\sigma) \ind_{\Qbf_{\ul{n}}(\A_f)}^{\GLbf_h(\A_f)} \left( e_c(\Lbf_{\ul{n}}, \sigma \cdot \lambda') \right).
\end{equation}
where the dot action \(\sigma \cdot \lambda'\) is for \(\GLbf_h\).
For \(h=0\) we define \(C(\lambda',\tau'') = 1 \in K_0^{\Tr}(\Rep_{\Q}^{\adm}(\GLbf_0(\A_f))) \simeq K_0(\Q) \simeq \Z\).

\begin{lemm}
  The Euler characteristic
  \[ e(\Acal_{n,?,\Qbar}^*, \IC_{\ell}(V_\lambda)) \in K_0^{\Tr}(\Rep_{\Qell}^{\adm,\cont}(\Gbf(\A_f) \times \GalQ)) \]
  is equal to
  \begin{equation} \label{eq:init_IH_from_Hc}
    \sum_{0 \leq h \leq n} \sum_{w_2 \in W^{\Pbf_h}} \epsilon(w_2) \ind_{\Pbf_h(\A_f)}^{\GSpbf_{2n}(\A_f)} \left( C((w_2 \cdot \lambda)_{\lin}, w_2(\tau)_{\lin}) \otimes e_c(\Acal_{n-h,?,\Qbar}, \Fcal_{\ell}^?(V^{\GSpbf_{2(n-h)}}_{(w_2 \cdot \lambda)_{\her}})) \right).
  \end{equation}
\end{lemm}
\begin{proof}
  Of course we start from the formula in Corollary \ref{coro:IH_vs_Hc_Q}.
  First we make the positivity condition (Definition \ref{def:trunc_alg_rep}) appearing in the indexing set \(W^{\Pbf}_{>\ul{t}}(\lambda)\) more explicit.
  Denoting \(d_m = (n-m)(n-m+1)/2\) this condition reads
  \[ \sum_{i=1}^{n_j} (w \cdot \lambda)_i > d_{n_j} - d_0 \text{ for all } 1 \leq j \leq r. \]
  We have \((w \cdot \lambda)_i = w(\tau)_i - (n+1-i)\), and so these inequalities are equivalent to
  \[ \sum_{i=1}^{n_j} w(\tau)_i > 0 \text{ for all } 1 \leq j \leq r. \]
  
  We rewrite the double sum Corollary \ref{coro:IH_vs_Hc_Q} by first summing over \(0 \leq h \leq n\), then summing over standard parabolic subgroups \(\Pbf = \Pbf_{n_1} \cap \dots \cap \Pbf_{n_r}\) of \(\GSpbf_{2n}\) satisfying \(n_r = h\) (the case where \(\Pbf = \GSpbf_{2n}\) corresponding to \(h=0\)).
  Kostant representatives in \(W^\Pbf\) decompose uniquely as \(\sigma w_2\) where \(w_2 \in W^{\Pbf_h}\) and \(\sigma\) is a Kostant representative for \(\Lbf_{\ul{n}} \subset \GLbf_h\), and we have \((\sigma w_2 \cdot \lambda)_\her = (w_2 \cdot \lambda)_\her\).
  Finally the dot action for \(\Mbf_{\Pbf_h}\) coincides with the restriction of the dot action of \(\GSpbf_{2n}\) because \(\rho-\rho_{\Mbf_{\Pbf_h}}\) is fixed by the Weyl group of \(\Mbf_{\Pbf_h}\).
\end{proof}

We now aim to simplify the expression \eqref{eq:collect_nr_IH_from_Hc} for \(C(\lambda', \tau'')\) using Corollary \ref{coro:Euler_c_GL}.
It will be more convenient to use the more symmetric \emph{normalized} parabolic induction: we start from the equality in \(K_0^{\Tr}(\Rep_{\R}^{\adm}(\GLbf_h(\A_f)))\)
\begin{equation} \label{eq:collect_nr_IH_from_Hc_norm}
  C(\lambda',\tau'') = \sum_{\substack{1 \leq r \leq h \\ 0 < n_1 < \dots < n_{r-1} < h}} \sum_{\sigma \in \Sfrak_h(\ul{n}, \tau'')} \epsilon(\sigma) \Ind_{\Lbf_{\ul{n}}(\A_f)}^{\GLbf_h(\A_f)} \left( e_c(\Lbf_{\ul{n}}, \sigma \cdot \lambda') \otimes |\delta_{\Qbf_{\ul{n}}}|_f^{-1/2} \right)
\end{equation}
where \(\delta_{\Qbf_{\ul{n}}}\) is the determinant of the adjoint action of \(\Lbf_{\ul{n}}\) on the unipotent radical of \(\Qbf_{\ul{n}}\).
Note that for normalized parabolic induction we only indicate the Levi factor of the parabolic subgroup, since the semi-simplification of this induced representation does not depend on the choice of parabolic.
As usual working with normalized induction compromises on algebraicity, although it should be clear that all computations below could be done over \(\Q\).
In any case the Brauer-Nesbitt theorem implies that the natural map
\[ K_0^{\Tr}(\Rep_{\Q}^{\adm}(\GLbf_h(\A_f))) \longrightarrow K_0^{\Tr}(\Rep_{\R}^{\adm}(\GLbf_h(\A_f))) \]
is injective.
We need to introduce more notation.
For \((a,b) \in \Z^2_{\geq 0}\) such that \(h := a+2b\) is positive let \(\mathcal{P}_{a,b}\) be the set of triples \((r, (a_j, b_j)_{1 \leq j \leq r}, \delta)\) where \(r \in \Z_{\geq 1}\), \(a_j, b_j \in \Z_{\geq 0}\) and \(\delta \in \Sfrak_h\) satisfy:
\begin{itemize}
\item \(a = \sum_j a_j\), \(b = \sum_j b_j\) and for any \(1 \leq j \leq r\) we have \(a_j +2b_j >0\),
\item denoting \(n_0 = 0\) and \(n_j-n_{j-1} = a_j + 2b_j\), so that \(0 < n_1 < \dots < n_r = h\), we have for any \(1 \leq j \leq r\):
  \[ \delta^{-1}(n_{j-1} + 1) < \dots < \delta^{-1}(n_{j-1} + a_j) \leq a, \]
  \[ a + 1 \leq \delta^{-1}(n_{j-1} + a_j + 1) < \dots < \delta^{-1}(n_{j-1} + a_j + 2b_j-1), \]
  and for any \(1 \leq i \leq b_j\):
  \[ \delta^{-1}(n_{j-1} + a_j + 2i-1) \in a+1+2\Z_{\geq 0} \ \text{and} \]
  \[ \delta^{-1}(n_{j-1} + a_j +2i) = \delta^{-1}(n_{j-1}+a+2i-1)+1. \]
\end{itemize}
Note that such triples simply parametrize ordered partitions \(\{1, \dots, h\} = I_1 \sqcup \dots \sqcup I_r\) such that for any \(1 \leq i \leq b\), \(a+2i-1\) and \(a+2i\) belong to the same subset \(I_j\), by taking
\begin{equation} \label{eq:interp_Pcal_a_b_part}
  I_j = \left\{ \delta^{-1}(i) \,\middle|\, n_{j-1} < i \leq n_j \right\}.
\end{equation}
For \(\tau'' \in \R^h\) we also define
\[ \mathcal{P}_{a,b}(\tau'') = \left\{ (r, (a_j, b_j)_{1 \leq j \leq r}, \delta) \in \mathcal{P}_{a,b} \ \middle| \ \forall j \in \{1,\dots,r\},\, \sum_{i=1}^{n_j} \tau''_{\delta^{-1}(i)} > 0 \right\}. \]

\begin{lemm} \label{lem:step_C_lambda_tau}
  Let \(h\) be a positive integer.
  Let \(\lambda' = (\lambda'_1 \geq \dots \geq \lambda'_h) \in \Z^h\) be a dominant weight for \(\GLbf_h\) and \(\tau'' \in \R^h\).
  Denote \(\tau' = \lambda' + \rho_{\GLbf_h}\), i.e.\ \(\tau_i' = \lambda_i' + (h+1)/2 - i\) for \(1 \leq i \leq h\).
  Then \(C(\lambda', \tau'')\) is equal to
  \begin{multline} \label{eq:reorder_nr_IH_from_Hc}
    \sum_{\substack{a,b \geq 0 \\ a+2b=h}} \sum_{\eta \in \Sfrak(a,b)} \epsilon(\eta) (-1)^{a+b} \sum_{(r, (a_j, b_j)_{1 \leq j \leq r}, \delta) \in \mathcal{P}_{a,b}(\eta(\tau''))} (-1)^r \epsilon(\delta) \\
    \Ind_{\Lbf_{a, b}(\A_f)}^{\GLbf_h(\A_f)} \Big( \bigotimes_{i=1}^a e(\GLbf_1, (\eta \cdot \lambda')_i + \delta(i) - i) |\cdot|_f^{(\eta \cdot \lambda')_i + \delta(i) - i - \eta(\tau')_i} \\
  \otimes \bigotimes_{i=1}^b e_{(2)}(\GLbf_2, \eta(\tau')_{a + 2i-1} - 1/2, \eta(\tau')_{a + 2i} + 1/2) \Big).
\end{multline}
\end{lemm}
Recall that \(\Sfrak(a,b)\) was defined in Corollary \ref{coro:Euler_ord_GL}.
\begin{proof}
  We start from the expression \eqref{eq:collect_nr_IH_from_Hc_norm} for \(C(\lambda', \tau'')\).
  For \(0 < n_1 < \dots < n_r = h\), a permutation \(\sigma \in \mathfrak{S}_h\) is a Kostant representative for \(\Lbf_{\ul{n}}\) if and only if for every \(1 \leq j \leq r\),
  \[ \sigma^{-1}(n_{j-1} +1) < \dots < \sigma^{-1}(n_j) \]
  where by convention \(n_0=0\).
  We observe that such Kostant representatives correspond bijectively to ordered partitions \(\{1, \dots, h\} = I_1 \sqcup \dots \sqcup I_r\) with \(|I_j| = n_j - n_{j-1}\): set \(I_j = \{ \sigma^{-1}(i) \,|\, n_{j-1} < i \leq n_j \}\).
  %For \(1 \leq j \leq r\) denote by \(\sigma(\tau')^{(j)}\) the element of \((\frac{1}{2} \Z)^{n_j-n_{j-1}}\) given by \(\sigma(\tau')^{(j)}_i = \tau'_{\sigma^{-1}(n_{j-1}+i)}\) for \(1 \leq i \leq n_j-n_{j-1}\), and similarly define \((\sigma \cdot \lambda')^{(j)} \in \Z^{n_j-n_{j-1}}\) by \((\sigma \cdot \lambda')^{(j)}_i = \lambda'_{\sigma^{-1}(n_{j-1}+i)} + n_{j-1} + i - \sigma^{-1}(n_{j-1}+i)\) for \(1 \leq i \leq n_j-n_{j-1}\).

  For such an ordered partition of \(\{1,\dots,h\}\) we will be led to consider families \((a_j, b_j, \gamma_j)_{1 \leq j \leq r}\) where \(a_j, b_j \in \Z_{\geq 0}\) satisfy \(a_j+2b_j = |I_j|\) and \(\gamma_j \in \Sfrak(a_j,b_j)\), which we can think of as a partition of each \(I_j\) into \(a_j\) singletons and \(b_j\) pairs.
  In this situation we define \(\gamma \in \Sfrak_h\) by \(\gamma^{-1}(n_{j-1} + i) = n_{j-1} + \gamma_j^{-1}(i)\) for \(1 \leq j \leq r\) and \(1 \leq i \leq n_j-n_{j-1}\).
  We apply Corollary \ref{coro:Euler_c_GL} for each \(\GLbf_{n_j-n_{j-1}}\), and observing that the restriction of \(\delta_{\Qbf_{\ul{n}}}\) to the diagonal maximal torus of \(\GLbf_h\) is equal (in additive notation) to \(2 \rho_{\GLbf_h} - 2 \rho_{\Lbf_{\ul{n}}}\) we obtain that \(C(\lambda', \tau'')\) equals the sum over \(r \geq 1\), \((n_j)_{1 \leq j \leq r}\), \(\sigma \in \Sfrak_h(\ul{n}, \tau'')\) and \((a_j, b_j, \gamma_j)_{1 \leq j \leq r}\) as above of
  \begin{multline} \label{eq:term_in_expanded_collect_nr_IH_from_Hc}
    \epsilon(\sigma) \left( \prod_j (-1)^{(n_j-n_{j-1})(n_j-n_{j-1}+1)/2-1+a_j(a_j-1)/2} \epsilon(\gamma_j) \right) \times \\
    \Ind_{\prod_j \Lbf_{a_j, b_j}(\A_f)}^{\GLbf_h(\A_f)} \Big( \bigotimes_j \Big( \bigotimes_{i=1}^{a_j} e(\GLbf_1, (\gamma \sigma \cdot \lambda')_{n_{j-1}+i}) |\cdot|_f^{(\gamma \sigma \cdot \lambda')_{n_{j-1}+i} - \gamma \sigma(\tau')_{n_{j-1}+i}} \\
    \otimes \bigotimes_{i=1}^{b_j} e_{(2)}(\GLbf_2, \gamma \sigma(\tau')_{n_{j-1} + a_j + 2i-1} - 1/2, \gamma \sigma(\tau')_{n_{j-1} + a_j + 2i} +1/2) \Big) \Big).
  \end{multline}

  Let us simplify the first line.
  We clearly have \(\prod_{j=1}^r \epsilon(\gamma_j) = \epsilon(\gamma)\), and we claim
  \begin{equation} \label{eq:step_C_lambda_tau_signs}
    \prod_{j=1}^r (-1)^{(n_j-n_{j-1})(n_j-n_{j-1}+1)/2-1+a_j(a_j-1)/2} = (-1)^{a+b+r}.
  \end{equation}
  This follows from the congruence
  \[ \frac{a'(a'-1)}{2} + \frac{(a'+2b')(a'+2b'+1)}{2} \equiv a'+b' \mod 2 \]
  which holds for any pair \((a',b')\) of non-negative integers and is easily proved by induction on \(b'\), applied to each pair \((a_j,b_j)\).

  Consider \(r \geq 1\), \(0 < n_1 < \dots < n_r = h\), \(\sigma \in \Sfrak_h(\ul{n}, \tau'')\) and \((a_j, b_j, \gamma_j)_{1 \leq j \leq r}\) as above, with \(\gamma \in \Sfrak_h\) defined as above.
  We can associate to this datum \(a = \sum_j a_j\), \(b = \sum_j b_j\) and \(\eta \in \Sfrak(a,b)\) characterized by the relations
  \[ \{\eta^{-1}(i) \,|\, 1 \leq i \leq a \} = \{ \sigma^{-1} \gamma^{-1}(n_{j-1}+i) \,|\, 1 \leq j \leq r,\, 1 \leq i \leq a_j \} \]
  \begin{multline*}
    \{ (\eta^{-1}(a+2i-1),\eta^{-1}(a+2i)) \,|\, 1 \leq i \leq b \} = \\
    \{ (\sigma^{-1} \gamma^{-1}(n_{j-1}+a_j+2i-1), \sigma^{-1} \gamma^{-1}(n_{j-1}+a_j+2i-1)) \,|\, 1 \leq j \leq r,\, 1 \leq i \leq b_j \}
  \end{multline*}
  and then \(\delta := \gamma \sigma \eta^{-1}\) is such that \(P := (r,(a_j,b_j)_j,\delta)\) belongs to \(\mathcal{P}_{a,b}(\tau'')\).
  In fact this yields a bijection
  \[ (r, (n_j)_{1 \leq j \leq r}, \sigma, (a_j,b_j,\gamma_j)_{1 \leq j \leq r}) \longmapsto (a,b,\eta,P) \]
  with the set of quadruples satisfying \(a+2b = h\), \(\eta \in \Sfrak(a,b)\) and \(P \in \mathcal{P}_{a,b}(\eta(\tau''))\).
  It is rather tedious to check this formally, but note that we are simply reordering choices here: instead of choosing first an ordered partition \(\{1,\dots,h\} = I_1 \sqcup \dots \sqcup I_r\) and then a partition of each \(I_j\) into singletons and pairs, we can first choose a partition of \(\{1,\dots,h\}\) into singletons and pairs and then choose how to distribute these singletons and pairs in \(r\) packets; the conditions involving \(\tau''\) are equivalent because we have \(\gamma \sigma = \delta \eta\).
  This bijection allows us to reindex the sum of \eqref{eq:term_in_expanded_collect_nr_IH_from_Hc} over \(r \geq 1\), \((n_j)_{1 \leq j \leq r}\), \(\sigma \in \Sfrak_h(\ul{n}, \tau'')\) and \((a_j, b_j, \gamma_j)_{1 \leq j \leq r}\), first summing over \((a,b)\) such that \(a+2b=h\) and \(\eta \in \mathfrak{S}(a,b)\).
  Using \eqref{eq:step_C_lambda_tau_signs} and the equality
  \begin{align*}
    & \left\{ ((\gamma \sigma \cdot \lambda')_{n_{j-1}+i}, \gamma \sigma(\tau'')_{n_{j-1}+i}) \,\middle|\, 1 \leq j \leq r,\, 1 \leq i \leq a_j \right\} \\
    = & \left\{ ((\eta \cdot \lambda')_i +\delta(i)-i, \eta(\tau'')_i) \,\middle|\, 1 \leq i \leq a \right\}
  \end{align*}
  we obtain the formula claimed in the lemma.
\end{proof}

We now aim to simplify the innermost sum in \eqref{eq:reorder_nr_IH_from_Hc}.
Recall from Example \ref{exam:Euler_L2_coh_GL} that for \(k \in \Z\) and \(s \in \C\), the representation \(e(\GLbf_1, k)|\cdot|_f^{k-s}\) of \(\GLbf_1(\A_f)\) only depends on \((s, k \mod 2)\).
For this reason we can formalize our computations using the ring \(\Z[ t_1, \dots, t_a]/(t_i^2-1)_{1 \leq i \leq a}\), which as \(\Z\)-module admits \((\prod_{i=1}^a t_i^{\epsilon_i})_{\epsilon \in \{0,1\}^a}\) as a basis.
Let
\[ f(a,b,\tau'') = \sum_{(r,(a_j,b_j)_j,\delta) \in \mathcal{P}_{a,b}(\tau'')} (-1)^r \epsilon(\delta) \prod_{i=1}^a t_i^{\delta(i)-i} \in \Z[ t_1, \dots, t_a]/(t_i^2-1)_{1 \leq i \leq a} \]
and for \(\epsilon \in \{0,1\}^a\) let \(f(a,b,\tau'')_{\epsilon} \in \Z\) be the coefficient of \(\prod_{i=1}^a t_i^{\epsilon_i}\) in \(f(a,b,\tau'')\).
Then \eqref{eq:reorder_nr_IH_from_Hc} can be rewritten
\begin{multline} \label{eq:reorder_nr_IH_from_Hc2}
  \sum_{a,b,\eta} \epsilon(\eta) (-1)^{a+b} \sum_{\epsilon} f(a,b,\eta(\tau''))_{\epsilon} \times \Ind_{\Lbf_{a, b}(\A_f)}^{\GLbf_h(\A_f)} \Big( \bigotimes_{i=1}^a e(\GLbf_1, (\eta \cdot \lambda')_i + \epsilon_i) |\cdot|_f^{(\eta \cdot \lambda')_i + \epsilon_i - \eta(\tau)_i} \\
  \otimes \bigotimes_{i=1}^b e_{(2)}(\GLbf_2, \eta(\tau')_{a + 2i-1} - 1/2, \eta(\tau')_{a + 2i} + 1/2) \Big).
\end{multline}

\begin{lemm} \label{lem:simplify_f_a_b_tau}
  For \(a,b \in \Z_{\geq 0}\) satisfying \(a+2b>0\) and \(\tau'' \in \R^{a+2b}\) satisfying \(\tau_1'' \geq \dots \geq \tau''_a\) we have
  \[ f(a,b,\tau'') =
    \begin{cases}
      0 & \text{ if } a > 0 \text{ and } \tau''_a \leq 0, \\
      0 & \text{ if } \tau''_{a+2i-1} + \tau''_{a+2i} \leq 0 \text{ for some } 1 \leq i \leq b, \\
      (-1)^{a(a+1)/2 + b} (t_1 \dots t_a)^{a-1} & \text{ otherwise.}
  \end{cases} \]
\end{lemm}
\begin{proof}
  In this proof we use the interpretation of \(\mathcal{P}_{a,b}\) (and its subset \(\mathcal{P}_{a,b}(\tau'')\)) as parametrizing certain partitions (see \eqref{eq:interp_Pcal_a_b_part}).
  \begin{enumerate}
    \item Assume \(b>0\) and that there exists \(1 \leq i_0 \leq b\) for which we have \(\tau''_{a+2i_0-1} + \tau''_{a+2i_0} \leq 0\).
      Let us show that \(f(a,b,\tau'')\) vanishes.
      The set \(\mathcal{P}_{a,b}(\tau'')\) can be partitioned into two subsets \(\mathcal{P}_{a,b}^{(i)}(\tau'')\), \(i \in \{1,2\}\), where \(\mathcal{P}_{a,b}^{(i)}(\tau'')\) is the set of partitions \(P = (I_1, \dots, I_r)\) such that \(I_j = \{a+2i_0-1, a+2i_0\}\) for some \(1 < j \leq r\).
      Note that we cannot have \(I_1 = \{a+2i_0-1, a+2i_0\}\) because \(\tau''_{a+2i_0-1} + \tau''_{a+2i_0} \leq 0\).
      There is a natural bijection \(\mathcal{P}_{a,b}^{(1)}(\tau'') \simeq \mathcal{P}_{a,b}^{(2)}(\tau'')\): if \(P_1 = (I_1, \dots, I_{r_1})\) is such that \(I_j = \{a+2i_0-1 ,a+2i_0\}\) for some \(1 < j \leq r_1\), consider \(P_2 = (I_1, \dots, I_{j-1} \sqcup I_j, I_{j+1}, \dots, I_{r_2})\).
      Note that surjectivity uses the inequality \(\tau''_{a+2i_0-1} + \tau''_{a+2i_0} \leq 0\).
      For \(P_1, P_2\) as above the permutations \(\delta_1, \delta_2 \in \mathfrak{S}_{a+2b}\) associated to \(P_1\) and \(P_2\) differ by \(|I_{j-1} \cap \{ a+2i_0+1, \dots, a+2b \}|\) (an even number) of transpositions and we have \(r_2 = r_1 - 1\).
      By pairwise cancellation we obtain \(f(a,b,\tau'') = 0\) in this case.
    \item Assume now \(a > 0\) and \(\tau''_a \leq 0\).
      We now have a partition \(\mathcal{P}_{a,b}(\tau'') = \mathcal{P}_{a,b}^{[1]}(\tau'') \sqcup \mathcal{P}_{a,b}^{[2]}(\tau'')\) where \(\mathcal{P}_{a,b}^{[1]}(\tau'')\) is the set of partitions \(P = (I_1, \dots, I_r)\) such that \(I_j = \{a\}\) for some \(1 < j \leq r\).
      Again we have a natural bijection \(\mathcal{P}_{a,b}^{[1]}(\tau'') \simeq \mathcal{P}_{a,b}^{[2]}(\tau'')\): if \(P_1 = (I_1, \dots, I_{r_1})\) is such that \(I_j = \{a\}\) for some \(1 < j \leq r_1\), consider \(P_2 = (I_1, \dots, I_{j-1} \sqcup I_j, I_{j+1}, \dots, I_{r_2})\).
      Now \(\delta_1\) and \(\delta_2\) differ by \(|I_{j-1} \cap \{ a+1, \dots, a+2b \}|\) (an even number) of transpositions, \(\delta_1(i) = \delta_2(i)\) for \(1 \leq i < a\), \(\delta_1(a) = n_{j-1}+1\) and \(\delta_2(a) = n_{j-2} + a_{j-1} + 1 = n_{j-1} - 2b_j + 1\) (where \((a_j, b_j)_{1 \leq j \leq r_1}\) and consequently \((n_{j-1})_{1 \leq j \leq r_1}\) are associated to \(P_1\)).
      As before we obtain pairwise cancellation in the sum defining \(f(a,b,\tau'')\).
    \item Finally we assume [\(a=0\) or \(\tau''_a > 0\)] and \(\tau''_{a+2i-1} + \tau''_{a+2i} > 0\) for all \(1 \leq i \leq b\).
      Now we simply have \(\mathcal{P}_{a,b}(\tau'') = \mathcal{P}_{a,b}\).
      If \(b>0\) we have a partition \(\mathcal{P}_{a,b} = \bigsqcup_{i=1}^3 \mathcal{P}_{a,b}^{(i)}\) where
      \begin{itemize}
      \item \(\mathcal{P}_{a,b}^{(1)}\) is the set of partitions \(P = (I_1, \dots, I_r)\) such that \(I_j = \{a+2b-1, a+2b\}\) for some \(1 < j \leq r\).
      \item \(\mathcal{P}_{a,b}^{(2)}\) is the set of partitions \(P = (I_1, \dots, I_r)\) such that \(I_j \neq \{a+2b-1, a+2b\}\) for all \(1 \leq j \leq r\),
      \item \(\mathcal{P}_{a,b}^{(3)}\) is the set of partitions \(P = (I_1, \dots, I_r)\) such that \(I_1 = \{a+2b-1, a+2b\}\).
      \end{itemize}
      We have a bijection \(\mathcal{P}_{a,b}^{(1)} \simeq \mathcal{P}_{a,b}^{(2)}\) defined as in the first step (although it is a bijection for slightly simpler reasons), and we obtain
      \[ f(a,b,\tau'') = \sum_{P \in \mathcal{P}_{a,b}^{(3)}} (-1)^r \epsilon(\delta) \prod_{i=1}^a t_i^{\delta(i)-i} = -f(a, b-1, (\tau''_1, \dots, \tau''_{a+2b-2})) \]
      where the second equality follows from consideration of the bijection \(\mathcal{P}_{a,b}^{(3)} \simeq \mathcal{P}_{a,b-1}\), \(P = (I_1, \dots, I_r) \mapsto P' = (I_2, \dots, I_r)\): if \(\delta \in \mathfrak{S}_{a+2b}\) (resp.\ \(\delta' \in \mathfrak{S}_{a+2b-2}\)) is associated to \(P\) (resp.\ \(P'\)), we have \(\delta'(i) = \delta(i)-2\) for \(1 \leq i \leq a+2b-2\), and so \(\epsilon(\delta) = \epsilon(\delta')\) and \(\delta'|_{\{1, \dots, a\}} \equiv \delta|_{\{1, \dots, a\}} \mod 2\).
      By induction on \(b\) we get
      \[ f(a, b, \tau'') = (-1)^b f(a, 0, (\tau''_1, \dots, \tau''_a)). \]
      Thus we can assume \(b=0\) for the rest of the proof.
      We denote \(\mathcal{P}_a = \mathcal{P}_{a,0}\).
      If \(a=0\) the result is obvious, and we will conclude the proof by induction on \(a\).
      Suppose \(a>0\).
      We now use the decomposition \(\mathcal{P}_a = \bigsqcup_{i=1}^3 \mathcal{P}_a^{[i]}\) where \(\mathcal{P}_a^{[i]}\) are defined as in the second step and \(\mathcal{P}_a^{[3]}\) is the set of partitions \(P = (I_1, \dots, I_r)\) satisfying \(I_1 = \{a\}\).
      \begin{itemize}
      \item \(\mathcal{P}_a^{[1]}\) is the set of partitions \(P = (I_1, \dots, I_r)\) such that \(I_j = \{a\}\) for some \(1 < j \leq r\).
      \item \(\mathcal{P}_a^{[2]}\) is the set of partitions \(P = (I_1, \dots, I_r)\) such that \(I_j \neq \{a\}\) for all \(1 \leq j \leq r\),
      \item \(\mathcal{P}_a^{[3]}\) is the set of partitions \(P = (I_1, \dots, I_r)\) such that \(I_1 = \{a\}\).
      \end{itemize}
      Again the terms for \(\mathcal{P}_a^{[1]}\) and \(\mathcal{P}_a^{[2]}\) cancel and we obtain
      \begin{align*}
        f(a, 0, \tau'')
        &= \sum_{P \in \mathcal{P}_a^{[3]}} (-1)^r \epsilon(\delta) \prod_{i=1}^a t_i^{\delta(i)-i} \\
        &= (-1)^a t_1 \dots t_{a-1} t_a^{a-1} f(a-1, 0, (\tau''_1, \dots, \tau''_{a-1}))
      \end{align*}
      using the bijection \(\mathcal{P}_a^{[3]} \simeq \mathcal{P}_{a-1}\), \(P = (I_1, \dots, I_r) \mapsto P' = (I_2, \dots, I_r)\) and the fact that for \(\delta \in \mathfrak{S}_a\) (resp.\ \(\delta' \in \mathfrak{S}_{a-1}\)) associated to \(P\) (resp.\ \(P'\)) we have \(\delta(a) = 1\) and \(\delta'(i) = \delta(i)-1\) for \(1 \leq i < a\).
      Using the induction hypothesis we conclude
      \[ f(a, 0, \tau'') = (-1)^{a(a+1)/2} (t_1 \dots t_a)^{a-1}. \]
  \end{enumerate}
\end{proof}

Let \(\Sfrak(a,b,\tau'')\) be the set of \(\eta \in \Sfrak(a,b)\) such that
\[ \big( \eta(\tau'')_i \big)_{1 \leq i \leq a}, \, \big( \eta(\tau'')_{a+2i-1} + \eta(\tau'')_{a+2i} \big)_{1 \leq i \leq b} \]
are all positive.
We deduce from \eqref{eq:reorder_nr_IH_from_Hc2} and Lemma \ref{lem:simplify_f_a_b_tau} the following expression for \(C(\lambda', \tau'')\):
\begin{multline*}
  \sum_{\substack{a,b \geq 0 \\ a+2b=h}} \sum_{\eta \in \Sfrak(a,b,\tau'')} \epsilon(\eta) (-1)^{a(a-1)/2} \times \Ind_{\Lbf_{a,b}(\A_f)}^{\GLbf_h(\A_f)} \Big( \bigotimes_{i=1}^a e \left( \GLbf_1, (\eta \cdot \lambda')_i + a-1 \right) |\cdot|_f^{(\eta \cdot \lambda')_i + a-1 - \eta(\tau')_i} \\
  \otimes \bigotimes_{i=1}^b e_{(2)} \left( \GLbf_2, \eta(\tau)_{a+2i-1} - 1/2, \eta(\tau)_{a+2i} + 1/2\right) \Big).
\end{multline*}

This can be slightly simplified further.
For \(\eta \in \Sfrak(a,b,\tau'')\) define \(\eta' \in \Sfrak_h\) by
\[ (\eta')^{-1}(i) =
  \begin{cases}
    \eta^{-1}(a+1-i) & \text{ if } 1 \leq i \leq a \\
    \eta^{-1}(i) & \text{ if } a < i \leq h.
  \end{cases}
\]
Then \(\eta \mapsto \eta'\) defines a bijection \(\Sfrak(a,b,\tau'') \to \Sfrak'(a,b,\tau'')\) where \(\Sfrak'(a,b,\tau'')\) is defined as \(\Sfrak(a,b,\tau'')\) except that the condition \(0 < \eta^{-1}(1) < \dots < \eta^{-1}(a)\) (see Corollary \ref{coro:Euler_ord_GL}) is replaced by the condition \(0 < \eta'^{-1}(a) < \dots < \eta'^{-1}(1)\).
We obtain the following expression for \(C(\lambda', \tau'')\):
\begin{multline*}
  \sum_{\substack{a,b \geq 0 \\ a+2b=h}} \sum_{\eta' \in \Sfrak'(a,b,\tau'')} \epsilon(\eta') \times \Ind_{\Lbf_{a,b}(\A_f)}^{\GLbf_h(\A_f)} \Big( \bigotimes_{i=1}^a e \left( \GLbf_1, (\eta' \cdot \lambda')_i \right) |\cdot|^{(\eta' \cdot \lambda')_i - \eta'(\tau')_i} \\
  \otimes \bigotimes_{i=1}^b e_{(2)} \left( \GLbf_2, \eta'(\tau')_{a+2i-1} - 1/2, \eta'(\tau')_{a+2i} + 1/2 \right) \Big).
\end{multline*}
We convert back to \emph{unnormalized} parabolic induction:
\begin{equation} \label{eq:C_lambda_tau_final}
  C(\lambda', \tau'') = \sum_{\substack{a,b \geq 0 \\ a+2b=h}} \sum_{\eta' \in \Sfrak'(a,b,\tau'')} \epsilon(\eta') \, \ind_{\Qbf_{a,b}(\A_f)}^{\GLbf_h(\A_f)} \Big( e_{(2)}(\Lbf_{a,b}, \eta' \cdot \lambda') \Big)
\end{equation}
where \(\Qbf_{a,b}\) is the standard parabolic subgroup of \(\GLbf_h\) with Levi factor \(\GLbf_1^a \times \GLbf_2^b\).
Note that each individual term on the right-hand side is defined over \(\Q\) (not just \(\R\)).
This formula trivially holds true for \(h=0\) as well.

\begin{theo} \label{theo:IH_vs_Hc_simple}
  For \(n \geq 1\) and \(a,b \geq 0\) such that \(a+2b \leq n\) let \(W'(a,b,n)\) be the set of \(w \in W(\GSpbf_{2n})\) satisfying (as permutations of \(\{\pm 1, \dots, \pm n\}\))
  \[ w^{-1}(1) > \dots > w^{-1}(a) > 0 \]
  \[ 0 < w^{-1}(a+1) < w^{-1}(a+3) < \dots < w^{-1}(a+2b-1) \]
  \[ |w^{-1}(a+2)| > w^{-1}(a+1), \dots, |w^{-1}(a+2b)| > w^{-1}(a+2b-1) \]
  \[ 0 < w^{-1}(a+2b+1) < \dots < w^{-1}(n). \]
  For any dominant weight \(\lambda\) for \(\GSpbf_{2n}\) the Euler characteristic
  \[ e(\Acal_{n,?,\Qbar}^*, \IC_{\ell}(V_\lambda)) \in K_0^{\Tr}(\Rep_{\Qell}^{\adm,\cont}(\Gbf(\A_f) \times \GalQ)) \]
  is equal to
  \begin{multline*}
    \sum_{\substack{a,b \geq 0 \\ a+2b \leq n}} \sum_{w \in W'(a,b,n)} \epsilon(w) \, \ind_{\Pbf_{a,b,n}(\A_f)}^{\GSpbf_{2n}(\A_f)} \Big( e_{(2)}(\GLbf_1^a \times \GLbf_2^b, (w \cdot \lambda)_{\lin}) \\
    \otimes e_c(\Acal_{n-a-2b,?,\Qbar}, \Fcal_{\ell}^?(V^{\GSpbf_{2(n-a-2b)}}_{(w \cdot \lambda)_{\her}})) \Big)
  \end{multline*}
  where \(\Pbf_{a,b,n}\) is the standard parabolic subgroup of \(\GSpbf_{2n}\) with Levi \(\GLbf_1^a \times \GLbf_2^b \times \GSpbf_{2(n-a-2b)}\) and the linear and hermitian parts of \(w \cdot \lambda\) are as defined at the beginning of section \ref{sec:IH_from_Hc}.
\end{theo}
\begin{proof}
  We plug the final expression \eqref{eq:C_lambda_tau_final} for \(C(\lambda', \tau'')\) into \eqref{eq:init_IH_from_Hc}, taking \(\lambda' = (w_2 \cdot \lambda)_{\lin}\) and \(\tau'' = w_2(\tau)_{\lin}\) where \(w_2 \in W^{\Pbf_h}\).
  An element \(w_2 \in W\) is a Kostant representative for \(\Pbf_h \subset \GSpbf_{2n}\) if and only if \(w_2(\tau)_{\lin}\) and \(w_2(\tau)_{\her}\) are both dominant.
  The second condition is equivalent to (seeing \(w_2\) as a permutation of \(\{\pm 1, \dots, \pm n\}\))
  \[ 0 < w_2^{-1}(h+1) < \dots < w_2^{-1}(n). \]
  Translating the first condition, we see that it is equivalent to the existence of (a unique) \(i \in \{0,\dots,h\}\) for which we have
  \[ 0 < w_2^{-1}(1) < \dots < w_2^{-1}(i) \ \ \text{and} \ \ w_2^{-1}(i+1) < \dots < w_2^{-1}(h) < 0. \]
  Then for \(a,b \in \Z_{\geq 0}\) satisfying \(a+2b = h\) and \(\sigma \in \Sfrak_h\), we have \(\sigma \in \Sfrak'(a,b,w_2(\tau)_{\lin})\) if and only if
  \begin{itemize}
  \item \(0 < \sigma^{-1}(1) < \dots < \sigma^{-1}(a) \leq i\),
  \item \(\sigma^{-1}(a+1) < \dots < \sigma^{-1}(a+2b-1)\), and
  \item for any \(1 \leq j \leq b\) we have either
    \begin{itemize}
    \item \(\sigma^{-1}(a+2j-1) < \sigma^{-1}(a+2j) \leq i\) (i.e.\ \(\sigma w_2(\tau)_{a+2j-1}\) and \(\sigma w_2(\tau)_{a+2j}\) are both positive), or
    \item \(\sigma^{-1}(a+2j-1) \leq i < \sigma^{-1}(a+2j)\) and \(-w_2^{-1} \sigma^{-1}(a+2j) > w_2^{-1} \sigma^{-1}(a+2j-1)\) (i.e.\ \(\sigma w_2(\tau)_{a+2j}\) is negative and \(\sigma w_2(\tau)_{a+2j-1} > - \sigma w_2(\tau)_{a+2j}\)).
    \end{itemize}
  \end{itemize}
  We deduce
  \[ \left\{ \sigma w_2 \,\middle|\, w_2 \in W^{\Pbf_h}, \, \sigma \in \Sfrak'(a,b,w_2(\tau)_{\lin}) \right\} = W'(a,b,n) \]
  and the theorem follows.
\end{proof}

\subsection{Compactly supported in terms of intersection cohomology}
\label{sec:Hc_from_IH}

The result of the previous section can be roughly described as saying that the matrix expressing intersection cohomology in terms of compactly supported cohomology is unipotent with coefficients in \(\{-1,0,1\}\).
Somewhat surprisingly, this is also the case for the inverse matrix that we compute in this section.

For \(a, b, n \in \Z_{\geq 0}\) satisfying \(a+2b \leq n\) let \(W(a,b,n)\) be the set of \(w \in W(\GSpbf_{2n})\) satisfying (as permutations of \(\{\pm 1, \dots, \pm n\}\)):
\[ 0 < w^{-1}(1) < \dots < w^{-1}(a) \]
\[ 0 < w^{-1}(a+1) < \dots < w^{-1}(a+2b-1) \]
\[ |w^{-1}(a+2)| > w^{-1}(a+1), \dots, |w^{-1}(a+2b)| > w^{-1}(a+2b-1) \]
\[ 0 < w^{-1}(a+2b+1) < \dots < w^{-1}(n). \]
(The only difference with \(W'(a,b,n)\) is the first line of inequalities.)
To \(w \in W(a,b,n)\) is associated an unordered partition of \(\{1, \dots, n\}\) into \(a\) singletons, \(b\) pairs and a set having \(n-a-2b\) elements.
The fiber of any such partition is parametrized by \(\{ \pm 1 \}^b\) via \(w \mapsto (\sign(w^{-1}(a+2i)))_{1 \leq i \leq b}\).

\begin{theo} \label{theo:Hc_from_IH}
  For an integer \(n \geq 1\), a dominant weight \(\lambda\) for \(\GSpbf_{2n}\) and a prime number \(\ell\) the Euler characteristic
  \[ e_c(\Acal_{n,?,\Qbar}, \Fcal_{\ell}^?(V_\lambda)) \in K_0^{\Tr}(\Rep_{\Qell}^{\adm,\cont}(\GSpbf_{2n}(\A_f) \times \GalQ)) \]
  is equal to
  \begin{multline} \label{equ:Hc_from_IH}
    \sum_{\substack{a,b \geq 0 \\ a+2b \leq n}} \sum_{w \in W(a,b,n)} (-1)^{a+b} \epsilon(w)  \\
    \ind_{\Pbf_{a,b,n}(\A_f)}^{\GSpbf_{2n}(\A_f)} \left( e_{(2)} \left( \GLbf_1^a \times \GLbf_2^b, (w \cdot \lambda)_{\lin} \right) \otimes e(\Acal_{n-a-2b,?,\Qbar}^*, \IC_\ell(V_{(w \cdot \lambda)_{\her}})) \right).
  \end{multline}
\end{theo}
\begin{proof}
  Using Theorem \ref{theo:IH_vs_Hc_simple} we find that the right-hand side of \eqref{equ:Hc_from_IH} equals
  \begin{multline} \label{equ:pf_theo_Hc_from_IH1}
    \sum_{\substack{a_2, b_2 \geq 0 \\ a_2+2b_2 \leq n}} \sum_{w_2 \in W(a_2,b_2,n)} (-1)^{a_2+b_2} \epsilon(w_2) \sum_{\substack{a_1, b_1 \geq 0 \\ a_1+2b_1 \leq n-n_2}} \sum_{w_1 \in W'(a_1,b_1,n-n_2)} \epsilon(w_1) \\
    \ind_{\Pbf_{a_2, b_2, a_1, b_1 ,n}(\A_f)}^{\GSpbf_{2n}(\A_f)} \Big( e_{(2)} \left( \GLbf_1^{a_2} \times \GLbf_2^{b_2}, (w_2 \cdot \lambda)_{\lin} \right) \otimes e_{(2)} \left( \GLbf_1^{a_1} \times \GLbf_2^{b_1}, (w_1 \cdot (w_2 \cdot \lambda)_{\her})_{\lin} \right) \\
    \otimes e_c(\Acal_{n-n_2-n_1,?,\Qbar} \Fcal_{\ell}(V_{(w_1 \cdot (w_2 \cdot \lambda)_\her)_\her})) \Big)
  \end{multline}
  where \(n_i = a_i + 2b_i\) and \(\Pbf_{a_2,b_2,a_1,b_1,n}\) is the standard parabolic subgroup of \(\GSpbf_{2n}\) with Levi factor \(\GLbf_1^{a_2} \times \GLbf_2^{b_2} \times \GLbf_1^{a_1} \times \GLbf_2^{b_1} \times \GSpbf_{2(n-n_2-n_1)}\) (in this order).
  We will reorder the sums, summing over \(a = a_1 + a_2\) and \(b = b_1 + b_2\) first, and we will show that the resulting inner sums for \((a,b) \neq (0,0)\) vanish.

  To this end we first introduce, for \(a_2, a_1, b \geq 0\) satisfying \(a_2 + a_1 + 2b \leq n\) the set \(W''(a_2, a_1, b, n)\) of \(w \in W(\GSpbf_{2n})\) satisfying
  \[ 0 < w^{-1}(1) < \dots < w^{-1}(a_2) \]
  \[ w^{-1}(a_2 + 1) > \dots > w^{-1}(a_2 + a_1) > 0 \]
  \[ 0 < w^{-1}(a+1) < \dots < w^{-1}(a+2b-1) \]
  \[ |w^{-1}(a+2)| > w^{-1}(a+1), \dots, |w^{-1}(a+2b)| >
    w^{-1}(a+2b-1) \]
  \[ 0< w^{-1}(a+2b+1) < \dots < w^{-1}(n) \]
  where \(a = a_1 + a_2\).
  For \(b_2,b_1 \geq 0\) satisfying \(b_1+b_2 = b\) we have a well-defined map
  \[ \xi: W(a_2,b_2,n) \times W'(a_1,b_1,n-n_2) \longrightarrow W''(a_2, a_1, b, n) \]
  characterized by the following conditions.
  For \(w_2 \in W'(a_2, b_2, n)\) and \(w_1 \in W(a_1, b_1, n-n_2)\), denoting \(w = \xi(w_2, w_1) \in W''(a_2, a_1, b, n)\) we have
  \begin{itemize}
  \item \(w^{-1}(i) = w_2^{-1}(i)\) for \(1 \leq i \leq a_2\),
  \item \(w^{-1}(a_2 + i) = w_2^{-1}(n_2 + w_1^{-1}(i))\) for \(1 \leq i \leq a_1\), and
  \item the set
    \[\left\{ (w^{-1}(a+2i-1), w^{-1}(a+2i)) \,\middle|\, 1 \leq i \leq b \right\} \]
    equals
    \begin{equation} \label{equ:pf_theo_Hc_from_IH_parti}
      \begin{aligned}
        & \left\{ (w_2^{-1}(a_2+2i-1), w_2^{-1}(a_2+2i)) \,\middle|\, 1 \leq i \leq b_2 \right\} \\
        \bigsqcup & \left\{ (w_2^{-1}(n_2 + w_1^{-1}(a_1+2i-1)), w_2^{-1}(n_2 + w_1^{-1}(a_1+2i))) \,\middle|\, 1 \leq i \leq b_1 \right\}.
      \end{aligned}
    \end{equation}
  \end{itemize}
  In fact we have \(w = w_0 w_1 w_2\) where \(w_0 \in \Sfrak_n \subset W(\GSpbf_{2n})\) satisfies
  \begin{itemize}
  \item for any \(1 \leq i \leq a_2\) we have \(w_0^{-1}(i) = i\),
  \item for any \(1 \leq i \leq a_1\) we have \(w_0^{-1}(a_2+i) = n_2 + i\),
  \item for any \(1 \leq i \leq b\) we have
    \begin{itemize}
    \item \(w_0^{-1}(a+2i-1)\) is either equal to \(a_2+2j-1\) for some \(j \in \{1,\dots,b_2\}\) or to \(n_2 + a_1 + 2j-1\) for some \(j \in \{1,\dots,b_1\}\) and
    \item \(w_0^{-1}(a+2i) = w_0^{-1}(a+2i-1)+1\),
    \end{itemize}
  \item for any \(a+2b < i \leq n\) we have \(w_0^{-1}(i) = i\).
  \end{itemize}
  We omit the straightforward but tedious verification that the map \(\xi\) is well-defined, satisfies \(\epsilon(\xi(w_2, w_1)) = \epsilon(w_2) \epsilon(w_1)\) for all \(w_2\) and \(w_1\) (this follows from \(\epsilon(w_0) = +1\)), and is surjective with each fibre having \(\binom{b}{b_2}\) elements, corresponding to the possible partitions in \eqref{equ:pf_theo_Hc_from_IH_parti}.
  (These facts are clear when considering elements of \(W(a_2,b_2,n)\), \(W'(a_1, b_1, n-n_2)\) and \(W''(a_2,a_1,b,n)\) as partitions into singletons, pairs and an extra set along with sign changes.)
  For \(1 \leq i \leq a_2\) we have \((w \cdot \lambda)_i = (w_1 w_2 \cdot \lambda)_i\) and for \(1 \leq i \leq a_1\) it is easy to check that we have
  \[ (w \cdot \lambda)_{a_2+i} = (w_1 w_2 \cdot \lambda)_{n_2+i} - 2b_2 \equiv (w_1 w_2 \cdot \lambda)_{n_2+i} \mod 2. \]
  We deduce that \eqref{equ:pf_theo_Hc_from_IH1} is equal to
  \begin{multline*}
    \sum_{\substack{a_2, a_1, b \geq 0 \\ a_2 + a_1 + 2b \leq n}} \sum_{w \in W''(a_2, a_1, b, n)} (-1)^{a_2} \sum_{b_2 = 0}^b (-1)^{b_2} \binom{b}{b_2} \epsilon(w) \\
    \ind_{\Pbf_{a, b, n}(\A_f)}^{\GSpbf_{2n}(\A_f)} \Big( e_{(2)} \left( \GLbf_1^a \times \GLbf_2^b, (w \cdot \lambda)_{\lin} \right) \otimes e_c(\Acal_{n-n_2-n_1,?,\Qbar} \Fcal_{\ell}(V_{(w \cdot \lambda)_\her})) \Big)
  \end{multline*}
  For \(b > 0\) we simply have \(\sum_{b_2=0}^b (-1)^{b_2} \binom{b}{b_2} = (1-1)^b = 0\) so this simplifies as
  \begin{equation} \label{equ:pf_theo_Hc_from_IH2}
    \sum_{0 \leq a \leq n} \sum_{\substack{0 \leq a_2 \leq a \\ w \in W''(a_2, a-a_2, 0, n)}} (-1)^{a_2} \epsilon(w) \ind_{\Pbf_{a, 0, n}(\A_f)}^{\GSpbf_{2n}(\A_f)} \left( e_{(2)} \left( \GLbf_1^a, (w \cdot \lambda)_{\lin} \right) \otimes e_c(\Acal_{n-a,?,\Qbar} \Fcal_{\ell}(V_{(w \cdot \lambda)_\her})) \right).
  \end{equation}
  To conclude it is enough to check that for \(a>0\) the inner sum in \eqref{equ:pf_theo_Hc_from_IH2} vanishes.
  As in the previous section (the proof of Lemma \ref{lem:simplify_f_a_b_tau} in particular) this follows from cancelling pairs of terms.
  Fix \(0 < a \leq n\).
  We define a partition of \(\bigsqcup_{0 \leq a_2 \leq a} \{a_2\} \times W''(a_2, a-a_2, 0, n)\) as \(W_1 \sqcup W_2\) and a bijection \(\beta : W_1 \simeq W_2\).
  Let \(W_1\) be the disjoint union of \(\{0\} \times W''(0, a, 0, n)\) and
  \[ \bigsqcup_{0 < a_2 < a} \{a_2\} \times \left\{ w \in W''(a_2, a-a_2, 0, n) \, \middle| \, w^{-1}(a_2+1) > w^{-1}(a_2) \right\} \]
  and let \(W_2\) be the disjoint union of \(\{a\} \times W''(a, 0, 0, n)\) and
  \[ \bigsqcup_{0 < a_2 < a} \{a_2\} \times \left\{ w \in W''(a_2, a-a_2, 0, n) \, \middle| \, w^{-1}(a_2+1) < w^{-1}(a_2) \right\}. \]
  The bijection \(\beta\) is simply defined by \(\beta(a_2, w) = (a_2+1, w)\).
  Thanks to the sign \((-1)^{a_2}\) we obtain that for any \(0 < a \leq n\) the inner sum in \eqref{equ:pf_theo_Hc_from_IH2} vanishes.
\end{proof}

\section{Special cases}

\subsection{Genus \(n \leq 3\)}
\label{sec:BFG}

We work out Theorems \ref{theo:IH_explicit_crude}, \ref{thm:sigma_is_tensor}, and \ref{theo:IH_vs_Hc_simple} or \ref{theo:Hc_from_IH} for \(n \leq 3\) and deduce \cite[Conjecture 7.1]{BFG} at the level of \(\ell\)-adic Galois representations (for any prime number \(\ell\)).
While comparing two rather large formulas is not terribly exciting, this comparison serves two purposes: it gives us confidence that the formulas in the present article are correct, and because it makes the computation by Bergström, Faber and van der Geer of the traces of certain Hecke operators on certain spaces of Siegel cusp forms in genus \(3\) unconditional (as in Examples 7.5 and 7.6 loc.\ cit., see also \S 9 loc.\ cit.).

We fix a prime number \(\ell\) and \(\iota: \C \simeq \Qellbar\), although the formulas in this section will ultimately not depend on the choice of \(\iota\).
First we note that in \eqref{equ:Hc_from_IH}, forgetting the Hecke action we have
\begin{multline*}
  \ind_{\Pbf_{a,b,n}(\A_f)}^{\GSpbf_{2n}(\A_f)} \left( e_{(2)} \left( \GLbf_1^a \times \GLbf_2^b, (w \cdot \lambda)_{\lin} \right) \otimes e(\Acal_{n-a-2b,?,\Qbar}^*, \IC_\ell(V_{(w \cdot \lambda)_{\her}})) \right)^{\GSpbf_{2n}(\Zhat)} \\
  = \dim e_{(2)} \left( \GLbf_1^a \times \GLbf_2^b, (w \cdot \lambda)_{\lin} \right)^{(\GLbf_1^a \times \GLbf_2)(\Zhat)} \times e(\Acal_{n-a-2b,\Qbar}^*, \IC_\ell(V_{(w \cdot \lambda)_{\her}}))
\end{multline*}
and we have for \(a \in \Z\)
\[ \dim e(\GLbf_1,a)^{\GLbf_1(\Zhat)} = \delta_{a \text{ even}} \]
and for \(a,b \in \R\) satisfying \(a-b \in \Z_{\geq 0}\)
\[ \dim e_{(2)}(\GLbf_2, a, b) = - \dim S_{a-b+2}(\SLbf_2(\Z)) + \delta_{a=b} \]
(see Example \ref{exam:Euler_L2_coh_GL}).
For \(a \in \Z_{\geq 0}\) we denote
\[ s_{a+2} = - \dim e_{(2)}(\GLbf_2, a, 0). \]
This notation is consistent with \cite[\S 2]{BFG}.
Since we work in level one in this section we will keep the level implicit and the notation and simply write \(\Fcal_\ell(V)\) for \(\Fcal^{\GSpbf_{2n}(\Zhat)}_\ell(V)\).
As recalled in Section \ref{sec:local_systems} for \(\lambda = (\lambda_1 \geq \dots \geq \lambda_n)\) the local system denoted by \(\mathbb{V}_{\lambda_1, \dots, \lambda_n}\) loc.\ cit.\ is our \(\Fcal_\ell(V_{\lambda,0})\) where \(V_{\lambda,0}\) is the irreducible representation of \(\GSpbf_{2n,\Q}\) of highest weight \((\lambda_1, \dots, \lambda_n, 0)\) (parametrization of weights as in Section \ref{sec:IH_from_Hc}).
We can reduce to this case by Remark \ref{rema:twisting_sim}.
If \(\sum_i \lambda_i\) is odd then cohomology (ordinary, compactly supported or intersection) vanishes.
If \(\sum_i \lambda_i\) is even then defining \(m \in \Z\) by \(2m = \sum_i \lambda_i\) the representation \(V_{\lambda,m} \simeq V_{\lambda,0} \otimes \nu^m\) has trivial central character (as assumed in Theorem \ref{theo:IH_explicit_crude}) and we have
\[ e(\Acal^*_{n,\Qbar}, \IC_\ell(V_{\lambda,0})) = e(\Acal^*_{n,\Qbar}, \IC_\ell(V_{\lambda,m})) \otimes \chi_\ell^{-m} \]
and similarly for compactly supported cohomology.

For integers \(\lambda_1 \geq \dots \geq \lambda_n \geq 0\) denote
\[ \ul{k}(\lambda) = (\lambda_1+n+1, \dots, \lambda_n+n+1) \in \Z^n, \]
which represents the same highest weight for \(\GLbf_n\) as \(n(\lambda)\) in \cite[Notation 4.3]{BFG}, except that our parametrizations of highest weights differ (we already used our parametrization in Section \ref{sec:SMF}).
Reformulating \cite[\S 5]{BFG} using Corollary \ref{cor:GSpin_Gal_rep_Siegel}, the authors conjectured for any \(n>1\) and integers \(\lambda_1 \geq \dots \geq \lambda_n \geq 0\) the existence of a virtual motive over \(\Q\), \(S[\ul{k}(\lambda)]\) (in their notation, \(S[n(\lambda)]\)) such that for any prime \(\ell\) we have
\begin{equation} \label{eq:S_k_lambda_ell}
  S[\ul{k}(\lambda)]_\ell = \sum_f \spin \circ \rho_{f,\iota}^{\GSpin}
\end{equation}
where the sum is over eigenforms in \(S_{\ul{k}}(\Spbf_{2n}(\Z))\), \(\rho_{f,\iota}^{\GSpin}\) was defined in Corollary \ref{cor:GSpin_Gal_rep_Siegel}.
We simply take \ref{eq:S_k_lambda_ell} as a definition.
We recall the (slightly different in the weight zero case) definition of \(S[-]_\ell\) for \(n=1\) in the next section.

\subsubsection{Genus one}

First consider the case \(n=1\).
Let \(k \geq 0\) be an integer.
Let \(\tau \in \ICcal(\Spbf_2)\) be the orbit of \(k+1 \in \Lie \Tcal_{\SO_3}\).
Any \(\psi \in \tilde{\Psi}_{\disc, \nonendo}^{\unr, \tau}(\Spbf_2)\) is either \([3]\) (only if \(k=0\)) or \((\Sym^2 \pi)[1]\) where \(\pi\) is a level one cuspidal automorphic representation for \(\PGLbf_2\) such that \(\pi_\infty\) has infinitesimal character \(\pm (k+1)/2\) (see Proposition \ref{pro:cpsc_Sym2}).
Such automorphic representations for \(\PGLbf_2\) correspond bijectively to eigenforms in \(S_{k+2}(\SLbf_2(\Z))\), in particular they exist only for \(k \geq 10\) even, \(k \neq 12\).
Deligne \cite{Deligne_GalGL2} proved the existence of a unique (up to conjugation) continuous Galois representation \(\rho_{\pi,\iota}: \GalQ \to \GL_2(\Qellbar)\) unramified away from \(\ell\) such that for all primes \(p \neq \ell\) we have \(\rho_{\pi,\iota}(\Frob_p)^{\sesi} \in \iota(p^{1/2} c(\pi_p))\).
(The more common normalization associates \(\chi_\ell^{-k} \rho_{\pi,\iota}\) to \(\pi\).)
In any case the spin representation \(\spin_\psi: \GMpsisc \to \GL_2\) is an isomorphism and we have (see Remark \ref{rem:rho_GSpin_odd_small_rank})
\[ \sigma_{\psi,\iota}^{\spin} \simeq
  \begin{cases}
    1 + \chi_\ell^{-1} & \text{ if } \psi = [3],\\
    \rho_{\pi,\iota} & \text{ if } \psi = (\Sym^2 \pi)[1].
  \end{cases}
\]
For \(k \geq 0\) an integer, following \cite[\S 2]{BFG} define in \(K_0(\Rep_{\Qellbar}^{\cont}(\GalQ))\)
\[ S[k+2]_\ell =
  \begin{cases}
    0 & \text{ if } k \text{ is odd},\\
    -1 - \chi_\ell^{-1} & \text{ if } k=0,\\
    \sum_\psi \chi_\ell^{-k/2} \rho_{\pi,\iota} & \text{ if } k>0 \text{ even}.
  \end{cases}
\]
For any even integer \(k \geq 0\) by Theorem \ref{theo:IH_explicit_crude} applied to the representation \(V_{k,k/2} \simeq V_{k,0} \otimes \nu^{k/2}\) of \(\PGSpbf_2\) and Remark \ref{rema:twisting_sim} we have (in \(K_0(\Rep_{\Qellbar}^{\cont}(\GalQ))\))
\[ e(\Acal^*_{1,\Qbar}, \IC_\ell(V_{k,0})) = -S[k+2]_\ell \]
which shows in particular that the right-hand side does not depend on the choice of \(\iota\).
Theorem \ref{theo:Hc_from_IH} for \(n=1\) yields
\[ e_c(\Acal_{1,\Qbar}, \Fcal_\ell(V_{k,0})) = e(\Acal^*_{1,\Qbar}, \IC_\ell(V_{k,0})) - \delta_{k \text{ even}} \]
because the only possible pairs \((a,b)\) occurring in the sum \eqref{equ:Hc_from_IH} are \((0,0)\) and \((1,0)\), and \(W(1,0,1) = \{\id\}\).
We thus recover \cite[Theorem 2.3]{BFG}.

\subsubsection{Genus two}
\label{sec:BFG_genus_two}

We now consider the case \(n=2\), which is also already known (see \cite{Petersen}) but is a good sanity check.
Following \cite[Conjecture 6.3]{BFG} we define for integers \(\lambda_1 \geq \lambda_2 \geq 0\) an element \(e_{2,\extr}(\lambda_1,\lambda_2)_\ell\) of \(K_0(\Rep_{\Qellbar}^{\cont}(\GalQ))\)) as
\begin{equation} \label{eq:def_e_2_extr}
  -s_{\lambda_1+\lambda_2+4} \, (S[\lambda_1-\lambda_2+2]_\ell + 1) \chi_\ell^{-\lambda_2-1} + s_{\lambda_1-\lambda_2+2} - S[\lambda_1+3]_\ell + S[\lambda_2+2]_\ell + \delta_{\lambda_1 \text{ even}}
\end{equation}
(We note that this definition for \(\lambda=(0,0)\) will not be used until we consider the genus three case.)
First we compare the contributions of a parameter \(\psi\) to \(S[\ul{k}(\lambda)]_\ell\) and to \(e(\Acal^*_{2,\Qbar}, \IC_\ell(V_{\lambda,0}))\), assuming \(\lambda \neq (0,0)\) and \(\lambda_1+\lambda_2\) even.
Denote \(m = (\sum_i \lambda_i)/2\).
\begin{enumerate}[leftmargin=0.4cm]
\item If \(\psi = (\Lambda^* \pi)[1]\) where \(\pi\) is a self-dual automorphic cuspidal representation for \(\PGLbf_4\) of symplectic type (see Proposition \ref{pro:cpsc_Sp4_SO5} and Remark \ref{rem:rho_GSpin_odd_small_rank}) then it contributes
  \[ - \chi_\ell^{-m} \sigma_{\psi,\iota}^{\spin} = - \chi_\ell^{-m} \spin \circ \rho_{\psi,\iota}^{\GSpin} = \chi_\ell^{-m-1} \std \circ \rho_{\pi,\iota}^{\GSp} \]
  (recall \(\rho_{\pi,\iota}^{\GSp}: \GalQ \to \GSp_4(\Qellbar)\) from Theorem \ref{theo:existence_rho_GSp}) to \(e(\Acal^*_{2,\Qbar}, \IC_\ell(V_{\lambda,0}))\), and the opposite to \(S[\ul{k}(\lambda)]_\ell\).
\item Assume \(\psi = (\pi_1 \otimes \pi_2)[1] \oplus [1]\) where each \(\pi_i\) is a cuspidal automorphic representation for \(\PGLbf_2\), \(\pi_{1,\infty}\) has infinitesimal character \(\pm (\lambda_1+\lambda_2+3)/2\) and \(\pi_{2,\infty}\) has infinitesimal character \(\pm (\lambda_1-\lambda_2+1)\) (see Proposition \ref{pro:cpsc_tensor_PGL2_PGL2} for the lift \(\pi_1 \otimes \pi_2\)).
  We compute (see \eqref{eq:def_ui_psi})
  \[ u_1(\psi) = \langle \mu_{\pi_\infty^\gen}, \tilde{s}_1 \rangle = -1 \]
  and so \(\psi\) contributes
  \[ -\chi_\ell^{-m-1} \rho_{\pi_2,\iota} \]
  to \(e(\Acal^*_{2,\Qbar}, \IC_\ell(V_{\lambda,0}))\), and does not contribute to \(S[\ul{k}(\lambda)]_\ell\) as an application of Arthur's multiplicity formula (see the proof of Corollary \ref{cor:GSpin_Gal_rep_Siegel}).
  Parameters of this shape contribute
  \[ -\delta_{\lambda_1>\lambda_2} \, \chi_\ell^{-\lambda_2-1} s_{\lambda_1+\lambda_2+4} \, S[\lambda_1-\lambda_2+2]_\ell \]
  to \(e(\Acal^*_{2,\Qbar}, \IC_\ell(V_{\lambda,0})) + S[\ul{k}(\lambda)]_\ell\).
\item Assume \(\psi = \pi[2] + 1\) where \(\pi\) is a self-dual automorphic cuspidal representation for \(\GLbf_2\) and \(\pi_\infty\) has infinitesimal character \(\pm (\lambda_1+\lambda_2+3)/2\).
  This imposes \(\lambda_1=\lambda_2\).
  This is rather similar to the previous case except that \(\epsilon_\psi\) is not always trivial and we have \(u_1(\psi) = -\epsilon_\psi(s_1)\).
  \begin{enumerate}[leftmargin=0.4cm]
  \item if \(\epsilon_\psi(s_1)=+1\) then \(\psi\) contributes \(\chi_\ell^{-m-1}+\chi_\ell^{-m-2}\) to \(e(\Acal^*_{2,\Qbar}, \IC_\ell(V_{\lambda,0}))\) and does not contribute to \(S[\ul{k}(\lambda)]_\ell\),
  \item if \(\epsilon_\psi(s_1)=-1\) then \(\psi\) contributes \(-\chi_\ell^{-m-1} \rho_{\pi,\iota}\) to \(e(\Acal^*_{2,\Qbar}, \IC_\ell(V_{\lambda,0}))\) and contributes \(\chi_\ell^{-m-1} (\rho_{\pi,\iota} + 1 + \chi_\ell^{-1})\) to \(S[\ul{k}(\lambda)]_\ell\).
  \end{enumerate}
  Parameters of this shape contribute
  \[ -\delta_{\lambda_1=\lambda_2} \, \chi_\ell^{-\lambda_2-1} s_{\lambda_1+\lambda_2+4} \, S[\lambda_1-\lambda_2+2]_\ell \]
  to \(e(\Acal^*_{2,\Qbar}, \IC_\ell(V_{\lambda,0})) + S[\ul{k}(\lambda)]_\ell\).
\end{enumerate}
We conclude for \(\lambda \neq (0,0)\)
\[ e(\Acal^*_{2,\Qbar}, \IC_\ell(V_{\lambda,0})) + S[\ul{k}(\lambda)]_\ell = - s_{\lambda_1+\lambda_2+4} \chi_\ell^{-\lambda_2-1} S[\lambda_1-\lambda_2+2]_\ell. \]
Now we use Theorem \ref{theo:Hc_from_IH} to express \(e_c(\Acal_{2,\Qbar}, \Fcal_\ell(V_{\lambda,0}))\) as \(e(\Acal^*_{2,\Qbar}, \IC_\ell(V_{\lambda,0}))\) plus the following contributions.
\begin{enumerate}[leftmargin=0.4cm]
\item For \((a,b) = (1,0)\) we have \(W(1,0,2) = \{ \id, (1 2) \}\), contributing respectively
  \[ - \delta_{\lambda_1 \text{ even}} \, e(\Acal^*_{1,\Qbar}, \IC_\ell(V_{\lambda_2,0})) = \delta_{\lambda_1 \text{ even}} \, S[\lambda_2+2]_\ell, \]
  \[ \delta_{\lambda_2 \text{ odd}} \, e(\Acal^*_{1,\Qbar}, \IC_\ell(V_{\lambda_1+1, 0})) = - \delta_{\lambda_2 \text{ odd}} \, S[\lambda_1+3]_\ell. \]
  Note that as we assume \(\lambda_1+\lambda_2\) even the two Kronecker \(\delta\) are superfluous.
\item For \((a,b) = (2,0)\) we have \(W(2,0,2) = \{\id\}\), contributing \(\delta_{\lambda_1 \text{ and } \lambda_2 \text{ even}}\).
\item For \((a,b)=(0,1)\) we have \(W(0,1,2) = \{ \id, (1 \mapsto 1, 2 \mapsto -2) \}\), contributing respectively
  \[ s_{\lambda_1-\lambda_2+2}, \]
  \[ -s_{\lambda_1+\lambda_2+4} \, \chi_\ell^{-\lambda_2-1}. \]
\end{enumerate}
We conclude for \(\lambda \neq (0,0)\)
\begin{align*}
  &e_c(\Acal_{2,\Qbar}, \Fcal_\ell(V_{\lambda,0})) \\
  =\ & -S[\ul{k}(\lambda)]_\ell - s_{\lambda_1+\lambda_2+4} \, \chi_\ell^{-\lambda_2-1} (1+S[\lambda_1-\lambda_2+2]_\ell) + S[\lambda_2+2]_\ell - S[\lambda_1+3]_\ell \\
  & + \delta_{\lambda_1 \text{ even}} + s_{\lambda_1-\lambda_2+2} \\
  =\ & -S[\ul{k}(\lambda)]_\ell + e_{2,\extr}(\lambda_1,\lambda_2)_\ell,
\end{align*}
recovering \cite[Conjecture 6.3]{BFG}, already proved by Petersen \cite{Petersen}.

\subsubsection{Genus three}

\begin{theo} \label{thm:BFG}
  Conjecture 7.1 of \cite{BFG} holds true at the level of \(\ell\)-adic Galois representations, i.e.\ in \(K_0(\Rep_{\Qellbar}^{\cont}(\GalQ))\)) we have for any \(\lambda \neq (0,0,0)\)
  \[ e_c(\Acal_{3,\Qbar}, \Fcal_\ell(V_{\lambda,0})) = S[\ul{k}(\lambda)]_\ell + \sum_{(\eta,a,b,c) \in X(\lambda)} \eta \times \left( e_c(\Acal_{2,\Qbar}, \Fcal_\ell(V_{a,b,0})) + e_{2,\extr}(a,b)_\ell \otimes S[c]_\ell \right) \]
  where
  \[ X(\lambda) = \{ (-1, \lambda_1+1, \lambda_2+1, \lambda_3+2), (1, \lambda_1+1, \lambda_3, \lambda_2+3), (-1, \lambda_2, \lambda_3, \lambda_1+4) \} \]
  and \(e_{2,\extr}(-,-)_\ell\) defined in \eqref{eq:def_e_2_extr}.
\end{theo}

This section is devoted to the proof of Theorem \ref{thm:BFG}.

First we compare \(e_c(\Acal^*_{3,\Qbar}, \IC_\ell(V_{\lambda,0}))\) and \(S[\ul{k}(\lambda)]_\ell\).
Along the way we prove \cite[Conjecture 7.7]{BFG} (thanks to Theorem \ref{thm:Siegel_formula}) when it holds true (see Remark \ref{rem:conj77_BFG}), ``upgrading'' \cite[Proposition 9.5]{ChRe} to Satake parameters in \(\Spin_7(\C)\).
\begin{enumerate}[leftmargin=0.4cm]
\item For \(\psi = \pi[1]\) a single level one self-dual automorphic cuspidal representation for \(\PGLbf_7\) such that the infinitesimal character of \(\pi_\infty\) is \((\pm (\lambda_1+3), \pm (\lambda_2 + 2), \pm (\lambda_1+1), 0)\),  it contributes
  \[ \chi_\ell^{-m} \sigma_{\psi,\iota}^{\spin} = \chi_\ell^{-m} \spin \circ \rho_{\pi,\iota}^{\GSpin} \]
  to both \(e(\Acal^*_{3,\Qbar}, \IC_\ell(V_{\lambda,0}))\) and \(S[\ul{k}(\lambda)]_\ell\).
\item For \(\psi = (\pi_1 \otimes \pi_2)[1] \oplus (\Sym^2 \pi_0)[1]\) where each \(\pi_i\) is a cuspidal automorphic representation for \(\PGLbf_2\), there are three subcases:
  \begin{enumerate}[leftmargin=0.4cm]
  \item infinitesimal characters \(\pi_{1,\infty} \mapsto \pm (\lambda_1+\lambda_2+5)/2\), \(\pi_{2,\infty} \mapsto \pm (\lambda_1-\lambda_2+1)/2\), \(\pi_{3,\infty} \mapsto \pm (\lambda_3+1)/2\).
    We have \(u_1(\psi) = \langle \mu_{\pi_\infty^\gen}, \tilde{s}_1 \rangle = -1\) and so \(\psi\) contributes
    \[ \chi_\ell^{-m-2} \rho_{\pi_2,\iota} \otimes \rho_{\pi_3,\iota} \]
    to \(e(\Acal^*_{3,\Qbar}, \IC_\ell(V_{\lambda,0}))\), and does not contribute to \(S[\ul{k}(\lambda)]_\ell\).
    These parameters contribute
    \[ \delta_{\lambda_1>\lambda_2} \, \delta_{\lambda_3>0} \, s_{\lambda_1+\lambda_2+6} \, \chi_\ell^{-\lambda_2-2} S[\lambda_1-\lambda_2+2]_\ell \otimes S[\lambda_3+2]_\ell \]
    to \(e(\Acal^*_{3,\Qbar}, \IC_\ell(V_{\lambda,0})) - S[\ul{k}(\lambda)]_\ell\).
  \item infinitesimal characters \(\pi_{1,\infty} \mapsto \pm (\lambda_1+\lambda_3+4)/2\), \(\pi_{2,\infty} \mapsto \pm (\lambda_1-\lambda_3+2)/2\), \(\pi_{3,\infty} \mapsto \pm (\lambda_2+2)/2\).
    We have \(u_1(\psi) = \langle \mu_{\pi_\infty^\gen}, \tilde{s}_1 \rangle = +1\) and so \(\psi\) contributes
    \[ \chi_\ell^{-m-2} \rho_{\pi_1,\iota} \otimes \rho_{\pi_3,\iota} \]
    to \(e(\Acal^*_{3,\Qbar}, \IC_\ell(V_{\lambda,0}))\), and contributes
    \[ \chi_\ell^{-m-2} (\rho_{\pi_1,\iota} + \rho_{\pi_2,\iota}) \otimes \rho_{\pi_3,\iota} \]
    to \(S[\ul{k}(\lambda)]_\ell\) (this corresponds to case (i) in \cite[Conjecture 7.7]{BFG}).
    These parameters contribute
    \[ - s_{\lambda_1+\lambda_3+5} \, \chi_\ell^{-\lambda_3-1} S[\lambda_1-\lambda_3+3]_\ell \otimes S[\lambda_2+3]_\ell \]
    to \(e(\Acal^*_{3,\Qbar}, \IC_\ell(V_{\lambda,0})) - S[\ul{k}(\lambda)]_\ell\).
  \item infinitesimal characters \(\pi_{1,\infty} \mapsto \pm (\lambda_2+\lambda_3+3)/2\), \(\pi_{2,\infty} \mapsto \pm (\lambda_2-\lambda_3+1)/2\), \(\pi_{3,\infty} \mapsto \pm (\lambda_1+3)/2\).
    As in the first case we have \(u_1(\psi) = -1\) and so \(\psi\) contributes
    \[ \chi_\ell^{-m-2} \rho_{\pi_2,\iota} \otimes \rho_{\pi_3,\iota} \]
    to \(e(\Acal^*_{3,\Qbar}, \IC_\ell(V_{\lambda,0}))\), and does not contribute to \(S[\ul{k}(\lambda)]_\ell\).
    These parameters contribute
    \[ \delta_{\lambda_2>\lambda_3} \, s_{\lambda_2+\lambda_3+4} \, \chi_\ell^{-\lambda_3-1} S[\lambda_2-\lambda_3+2]_\ell \otimes S[\lambda_1+4]_\ell \]
    to \(e(\Acal^*_{3,\Qbar}, \IC_\ell(V_{\lambda,0})) - S[\ul{k}(\lambda)]_\ell\).
  \end{enumerate}
\item For \(\psi = (\pi_1 \otimes \pi_2)[1] \oplus [3]\) with infinitesimal characters \(\pi_{1,\infty} \mapsto \pm (\lambda_1+\lambda_2+5)/2\) and \(\pi_{2,\infty} \mapsto \pm (\lambda_1-\lambda_2+1)/2\) (and imposing \(\lambda_3=0\)) we have
  \[ u_1(\psi) = \epsilon_\psi(s_1) \langle \mu_{\pi_\infty^\gen}, \tilde{s}_1 \rangle = -\epsilon(\pi_1 \times \pi_2) = -1 \]
  and so \(\psi\) contributes
  \[ -\chi_\ell^{-m-2} \rho_{\pi_2,\iota} \otimes (1 + \chi_\ell^{-1}) \]
  to \(e(\Acal^*_{3,\Qbar}, \IC_\ell(V_{\lambda,0}))\), and does not contribute to \(S[\ul{k}(\lambda)]_\ell\) because of the factor \([3]\) of \(\psi\) (see \cite[Lemma 9.2]{ChRe}).
  Parameters of this shape contribute
  \[ \delta_{\lambda_1>\lambda_2} \, \delta_{\lambda_3=0} \, s_{\lambda_1+\lambda_2+6} \, \chi_\ell^{-\lambda_2-2} S[\lambda_1-\lambda_2+2]_\ell \otimes S[\lambda_3+2]_\ell \]
  to \(e(\Acal^*_{3,\Qbar}, \IC_\ell(V_{\lambda,0})) - S[\ul{k}(\lambda)]_\ell\).
\item For \(\psi = \pi_1[2] \oplus (\Sym^2 \pi_2)[1]\) there are two cases to consider:
  \begin{enumerate}[leftmargin=0.4cm]
  \item for infinitesimal characters \(\pi_{1,\infty} \mapsto \pm (\lambda_1+\lambda_2+5)/2\) and \(\pi_{2,\infty} \mapsto \pm (\lambda_3+1)/2\) (imposing \(\lambda_1=\lambda_2\)) we have
    \[ u_1(\psi) = \epsilon_\psi(s_1) \langle \mu_{\pi_\infty^\gen}, \tilde{s}_1 \rangle = -\epsilon(\pi_1 \times \mathrm{ad}^0 \pi_2) = - (-1)^{1+\max(\lambda_1+\lambda_2+5, 2\lambda_3+2)} (-1)^{\lambda_1+3} = (-1)^{\lambda_1} \]
    so
    \begin{enumerate}[leftmargin=0.4cm]
    \item if \(\lambda_1\) is even then \(\psi\) contributes
      \[ \chi_\ell^{-m-2} \rho_{\pi_1,\iota} \otimes \rho_{\pi_2,\iota} \]
      to \(e(\Acal^*_{3,\Qbar}, \IC_\ell(V_{\lambda,0}))\), and
      \[ \chi_\ell^{-m-2} (\rho_{\pi_1,\iota} + 1 + \chi_\ell^{-1}) \otimes \rho_{\pi_2,\iota} \]
      to \(S[\ul{k}(\lambda)]_\ell\),
    \item if \(\lambda_1\) is odd then \(\psi\) contributes
      \[ -\chi_\ell^{-m-2} (1 + \chi_\ell^{-1}) \otimes \rho_{\pi_2,\iota} \]
      to \(e(\Acal^*_{3,\Qbar}, \IC_\ell(V_{\lambda,0}))\), and does not contribute to \(S[\ul{k}(\lambda)]_\ell\).
    \end{enumerate}
    In particular part (iii) of \cite[Conjecture 7.7]{BFG} holds true if and only if \(\lambda_1\) (denoted by \(a\) loc.\ cit.) is even.
    In any case parameters of this shape thus contribute
    \[ \delta_{\lambda_1=\lambda_2} \, \delta_{\lambda_3>0} \, \chi_\ell^{-\lambda_2-2} S[\lambda_1-\lambda_2+2]_\ell S[\lambda_3+2]_\ell \]
    to \(e(\Acal^*_{3,\Qbar}, \IC_\ell(V_{\lambda,0})) - S[\ul{k}(\lambda)]_\ell\).
  \end{enumerate}
\item for infinitesimal characters \(\pi_{1,\infty} \mapsto \pm (\lambda_2+\lambda_3+3)/2\) and \(\pi_{2,\infty} \mapsto \pm (\lambda_1+3)/2\) (imposing \(\lambda_2=\lambda_3\)) we have \(u_1(\psi) = -\epsilon_\psi(s_1)\) with
  \[ \epsilon_\psi(s_1) = (-1)^{\lambda_2+2} (-1)^{1+\max(2\lambda_2+3,2\lambda_1+6)} = (-1)^{\lambda_2+1} \]
  so
  \begin{enumerate}[leftmargin=0.4cm]
  \item if \(\lambda_2\) is odd then \(\psi\) contributes
    \[ - \chi_\ell^{-m-2} (1+\chi_\ell^{-1}) \otimes \rho_{\pi_2,\iota} \]
    to \(e(\Acal^*_{3,\Qbar}, \IC_\ell(V_{\lambda,0}))\) and does not contribute to \(S[\ul{k}(\lambda)]_\ell\),
  \item if \(\lambda_2\) is even then \(\psi\) contributes
    \[ \chi_\ell^{-m-2} \rho_{\pi_1,\iota} \otimes \rho_{\pi_2,\iota} \]
    to \(e(\Acal^*_{3,\Qbar}, \IC_\ell(V_{\lambda,0}))\) and contributes
    \[ \chi_\ell^{-m-2} (\rho_{\pi_1,\iota}+1+\chi_\ell^{-1}) \otimes \rho_{\pi_2,\iota} \]
    to \(S[\ul{k}(\lambda)]_\ell\).
  \end{enumerate}
  In particular part (ii) of \cite[Conjecture 7.7]{BFG} holds true if and only if \(\lambda_2\) (denoted by \(b\) loc.\ cit.) is even.
  In any case parameters of this shape thus contribute
  \[ \delta_{\lambda_2=\lambda_3} \, s_{\lambda_2+\lambda_3+4} \, \chi_\ell^{-\lambda_3-1} S[\lambda_2-\lambda_3+2]_\ell \otimes S[\lambda_1+4]_\ell \]
  to \(e(\Acal^*_{3,\Qbar}, \IC_\ell(V_{\lambda,0})) - S[\ul{k}(\lambda)]_\ell\).
\item Finally for \(\psi = \pi[2] \oplus [3]\) where \(\pi_{1,\infty}\) has infinitesimal character \(\pm (\lambda_1+\lambda_2+5)/2\) (and imposing \(\lambda_1=\lambda_2\) and \(\lambda_3=0\)) we have
  \[ u_1(\psi) = \epsilon_\psi(s_1) \langle \mu_{\pi_\infty^\gen}, \tilde{s}_1 \rangle = -\epsilon(\pi_1)^2 = -1 \]
  and so \(\psi\) contributes
  \[ \chi_\ell^{-m-2} (1+\chi_\ell^{-1}) \otimes (1+\chi_\ell^{-1}) \]
  to \(e(\Acal^*_{3,\Qbar}, \IC_\ell(V_{\lambda,0}))\) and does not contribute to \(S[\ul{k}(\lambda)]_\ell\).
  Parameters of this shape contribute
  \[ \delta_{\lambda_1=\lambda_2} \, \delta_{\lambda_3=0} \, s_{\lambda_1+\lambda_2+6} \, \chi_\ell^{-\lambda_2-2} S[\lambda_1-\lambda_2+2]_\ell \otimes S[\lambda_3+2]_\ell \]
  to \(e(\Acal^*_{3,\Qbar}, \IC_\ell(V_{\lambda,0})) - S[\ul{k}(\lambda)]_\ell\).
\end{enumerate}
Summing all these contributions we obtain for \(\lambda \neq (0,0,0)\)
\begin{align}
  & e(\Acal^*_{3,\Qbar}, \IC_\ell(V_{\lambda,0})) - S[\ul{k}(\lambda)]_\ell \nonumber \\
  =\ & s_{\lambda_1+\lambda_2+6} \, \chi_\ell^{-\lambda_2-2} S[\lambda_1-\lambda_2+2]_\ell \otimes S[\lambda_3+2]_\ell \nonumber \\
  & - s_{\lambda_1+\lambda_3+5} \, \chi_\ell^{-\lambda_3-1} S[\lambda_1-\lambda_3+3]_\ell \otimes S[\lambda_2+3]_\ell \nonumber \\
  & + s_{\lambda_2+\lambda_3+4} \, \chi_\ell^{-\lambda_3-1} S[\lambda_2-\lambda_3+2]_\ell \otimes S[\lambda_1+4]_\ell \label{eq:eIH_vs_Siegel_genus3}
\end{align}

\begin{rema} \label{rem:conj77_BFG}
  Note that there are parity conditions for the existence of the lifts predicted in cases (ii) and (iii) \cite[Conjecture 7.7]{BFG}, in agreement with \cite[Proposition 9.5]{ChRe}.
  % For example taking \((a,b,c)=(8,7,7)\) in (ii) of \cite[Conjecture 7.7]{BFG} we have eigenforms in \(S_{a+4}\) and \(S_{2b+4}\) (\(\Delta\) and \(E_6 \Delta\)) but no Siegel cusp form of the corresponding weight.
  In \cite[Table 2]{BFG} the authors seem to be aware of these conditions, so it seems that they were simply forgotten in the statement of \cite[Conjecture 7.7]{BFG}.
\end{rema}

Now we use Theorem \ref{theo:IH_vs_Hc_simple} to express \(e_c(\Acal_{3,\Qbar}, \Fcal_\ell(V_{\lambda,0}))\) as \(e(\Acal^*_{3,\Qbar}, \IC_\ell(V_{\lambda,0}))\) plus the following contributions (assuming \(\sum_i \lambda_i\) even)
\begin{enumerate}[leftmargin=0.4cm]
\item for \((a,b)=(1,0)\) we have \(W'(1,0,3) = \{\id, (1 2), (1 2 3)\}\), respectively contributing
  \[ -\delta_{\lambda_1 \text{ even}} \, e_c(\Acal_{2,\Qbar}, \Fcal_\ell(V_{\lambda_2,\lambda_3,0})), \]
  \[ \delta_{\lambda_2 \text{ odd}} \, e_c(\Acal_{2,\Qbar}, \Fcal_\ell(V_{\lambda_1+1,\lambda_3})), \]
  \[ -\delta_{\lambda_3 \text{ even}} \, e_c(\Acal_{2,\Qbar}, \Fcal_\ell(V_{\lambda_1+1,\lambda_2+1,0})). \]
  The Kronecker \(\delta\)'s are superfluous.
\item For \((a,b) = (2,0)\) the set \(W'(2,0,3)\) has three elements:
  \(\{(1 2), (1 2 3), (1 3)\}\), respectively contributing
  \begin{enumerate}[leftmargin=0.4cm]
  \item \(w=\begin{pmatrix} 1 & 2 & 3 \\ 2 & 1 & 3 \end{pmatrix}\) yields
    % \epsilon(w)=-1
    % w \cdot \lambda = (\lambda_2-1, \lambda_1+1, \lambda_3, 0)
    \[ \delta_{\lambda_1 \text{ odd}} \, \delta_{\lambda_2 \text{ odd}} \, e_c(\Acal_{1,\Qbar}, \Fcal_\ell(V_{\lambda_3,0})) = - \delta_{\lambda_1 \text{ odd}} \, S[\lambda_3+2]_\ell - \delta_{\lambda_1 \text{ odd}} \, \delta_{\lambda_2 \text{ odd}}. \]
  \item \(w=\begin{pmatrix} 1 & 2 & 3 \\ 2 & 3 & 1 \end{pmatrix}\) yields
    % \epsilon(w)=+1
    % w \cdot \lambda = (\lambda_3-2, \lambda_1+1, \lambda_2+1, 0)
    \[ -\delta_{\lambda_1 \text{ odd}} \, \delta_{\lambda_3 \text{ even}} \, e_c(\Acal_{1,\Qbar}, \Fcal_\ell(V_{\lambda_2+1,0})) = \delta_{\lambda_1 \text{ odd}} \, S[\lambda_2+3]_\ell + \delta_{\lambda_1 \text{ odd}} \, \delta_{\lambda_2 \text{ odd}}. \]
  \item \(w=\begin{pmatrix} 1 & 2 & 3 \\ 3 & 2 & 1 \end{pmatrix}\) yields
    % \epsilon(w)=-1
    % w \cdot \lambda = (\lambda_3-2, \lambda_2, \lambda_1+2, 0)
    \[ \delta_{\lambda_2 \text{ even}} \, \delta_{\lambda_3 \text{ even}} \, e_c(\Acal_{1,\Qbar}, \Fcal_\ell(V_{\lambda_1+2,0})) = -\delta_{\lambda_2 \text{ even}} \, S[\lambda_1+4]_\ell - \delta_{\lambda_1 \text{ even}} \, \delta_{\lambda_2 \text{ even}} \]
  \end{enumerate}
\item For \((a,b)=(3,0)\) we have \(W'(3,0,3)=\{(1 3)\}\), contributing \(\delta_{\lambda_1 \text{ even}} \, \delta_{\lambda_2 \text{ even}}\).
\item For \((a,b)=(0,1)\) the set \(W'(0,1,3)\) has six elements:
  \begin{enumerate}[leftmargin=0.4cm]
  \item \(w=\id\) yields
    \[ s_{\lambda_1-\lambda_2+2} \, e_c(\Acal_{1,\Qbar}, \Fcal_\ell(V_{\lambda_3,0})) = - s_{\lambda_1-\lambda_2+2} \, S[\lambda_3+2]_\ell - s_{\lambda_1-\lambda_2+2} \]
  \item \(w=\begin{pmatrix} 1 & 2 & 3 \\ 1 & -2 & 3 \end{pmatrix}\) yields
    % \epsilon(w) = -1
    % w \cdot \lambda = (\lambda_1, -\lambda_2-4, \lambda_3, -\lambda_2-2)
    \[ -s_{\lambda_1+\lambda_2+6} \, e_c(\Acal_{1,\Qbar}, \Fcal_\ell(V_{\lambda_3,-\lambda_2-2})) = s_{\lambda_1+\lambda_2+6} \, \chi_\ell^{-\lambda_2-2} S[\lambda_3+2]_\ell + s_{\lambda_1+\lambda_2+6} \, \chi_\ell^{-\lambda_2-2} \]
  \item \(w=\begin{pmatrix} 1 & 2 & 3 \\ 1 & 3 & 2 \end{pmatrix}\) yields
    % \epsilon(w) = -1
    % w \cdot \lambda = (\lambda_1, \lambda_3-1, \lambda_2+1, 0)
    \[ -s_{\lambda_1-\lambda_3+3} \, e_c(\Acal_{1,\Qbar}, \Fcal_\ell(V_{\lambda_2+1,0})) = s_{\lambda_1-\lambda_3+3} \, S[\lambda_2+3]_\ell + s_{\lambda_1-\lambda_3+3} \]
  \item \(w=\begin{pmatrix} 1 & 2 & 3 \\ 1 & 3 & -2 \end{pmatrix}\) yields
    % \epsilon(w) = +1
    % w \cdot \lambda = (\lambda_1, -\lambda_3-3, \lambda_2+1, -\lambda_3-1)
    \[ s_{\lambda_1+\lambda_3+5} \, e_c(\Acal_{1,\Qbar}, \Fcal_\ell(V_{\lambda_2+1,-\lambda_3-1})) = - s_{\lambda_1+\lambda_3+5} \, \chi_\ell^{-\lambda_3-1} S[\lambda_2+3]_\ell - s_{\lambda_1+\lambda_3+5} \, \chi_\ell^{-\lambda_3-1} \]
  \item \(w=\begin{pmatrix} 1 & 2 & 3 \\ 3 & 1 & 2 \end{pmatrix}\) yields
    % \epsilon(w) = +1
    % w \cdot \lambda = (\lambda_2-1, \lambda_3-1, \lambda_1+2, 0)
    \[ s_{\lambda_2-\lambda_3+2} \, e_c(\Acal_{1,\Qbar}, \Fcal_\ell(V_{\lambda_1+2,0})) = - s_{\lambda_2-\lambda_3+2} \, S[\lambda_1+4]_\ell - s_{\lambda_2-\lambda_3+2} \]
  \item \(w=\begin{pmatrix} 1 & 2 & 3 \\ 3 & 1 & -2 \end{pmatrix}\) yields
    % \epsilon(w) = -1
    % w \cdot \lambda = (\lambda_2-1, -\lambda_3-3, \lambda_1+2, -\lambda_3-1)
    \[ -s_{\lambda_2+\lambda_3+4} \, e_c(\Acal_{1,\Qbar}, \Fcal_\ell(V_{\lambda_1+2,-\lambda_3-1})) = s_{\lambda_2+\lambda_3+4} \, \chi_\ell^{-\lambda_3-1} S[\lambda_1+4]_\ell + s_{\lambda_2+\lambda_3+4} \, \chi_\ell^{-\lambda_3-1} \]
  \end{enumerate}
\item For \((a,b)=(1,1)\) the set \(W'(1,1,3)\) has six elements:
  \begin{enumerate}[leftmargin=0.4cm]
  \item \(w=\id\) yields
    \[ \delta_{\lambda_1 \text{ even}} \, s_{\lambda_2-\lambda_3+2} = s_{\lambda_2-\lambda_3+2}. \]
  \item \(w=\begin{pmatrix} 1 & 2 & 3 \\ 1 & 2 & -3 \end{pmatrix}\) yields
    % \epsilon(w) = -1
    % w \cdot \lambda = (\lambda_1, \lambda_2, -\lambda_3-2, -\lambda_3-1)
    \[ -\delta_{\lambda_1 \text{ even}} \, s_{\lambda_2+\lambda_3+4} \chi_\ell^{-\lambda_3-1} = -s_{\lambda_2+\lambda_3+4} \chi_\ell^{-\lambda_3-1}. \]
  \item \(w=\begin{pmatrix} 1 & 2 & 3 \\ 2 & 1 & 3 \end{pmatrix}\) yields
    % \epsilon(w) = -1
    % w \cdot \lambda = (\lambda_2-1, \lambda_1+1, \lambda_3, 0)
    \[ -\delta_{\lambda_2 \text{ odd}} \, s_{\lambda_1-\lambda_3+3} = -s_{\lambda_1-\lambda_3+3}. \]
  \item \(w=\begin{pmatrix} 1 & 2 & 3 \\ 2 & 1 & -3 \end{pmatrix}\) yields
    % \epsilon(w) = +1
    % w \cdot \lambda = (\lambda_2-1, \lambda_1+1, -\lambda_3-2, -\lambda_3-1)
    \[ \delta_{\lambda_2 \text{ odd}} \, s_{\lambda_1+\lambda_3+5} \chi_\ell^{-\lambda_3-1} = s_{\lambda_1+\lambda_3+5} \chi_\ell^{-\lambda_3-1}. \]
  \item \(w=\begin{pmatrix} 1 & 2 & 3 \\ 2 & 3 & 1 \end{pmatrix}\) yields
    % \epsilon(w) = +1
    % w \cdot \lambda = (\lambda_3-2, \lambda_1+1, \lambda_2+1, 0)
    \[ \delta_{\lambda_3 \text{ even}} \, s_{\lambda_1-\lambda_2+2} = s_{\lambda_1-\lambda_2+2}. \]
  \item \(w=\begin{pmatrix} 1 & 2 & 3 \\ 2 & -3 & 1 \end{pmatrix}\) yields
    % \epsilon(w) = -1
    % w \cdot \lambda = (\lambda_3-2, \lambda_1+1, -\lambda_2-3, -\lambda_2-2)
    \[ -\delta_{\lambda_3 \text{ even}} \, s_{\lambda_1+\lambda_2+6} \chi_\ell^{-\lambda_2-2} = - s_{\lambda_1+\lambda_2+6} \chi_\ell^{-\lambda_2-2}. \]
  \end{enumerate}
  Note that these contributions simplify the previous ones (for \((a,b)=(0,1)\)).
\end{enumerate}
Summing all these contributions with \eqref{eq:eIH_vs_Siegel_genus3} and factoring by \(S[\lambda_i+5-i]_\ell\) we obtain that for \(\lambda \neq (0,0,0)\)
\[ e_c(\Acal_{3,\Qbar}, \Fcal_\ell(V_{\lambda,0})) - S[\ul{k}(\lambda)]_\ell \]
is equal to
\begin{multline*}
  - e_c(\Acal_{2,\Qbar}, \Fcal_\ell(V_{\lambda_2,\lambda_3,0})) + e_c(\Acal_{2,\Qbar}, \Fcal_\ell(V_{\lambda_1+1,\lambda_3,0})) - e_c(\Acal_{2,\Qbar}, \Fcal_\ell(V_{\lambda_1+1,\lambda_2+1,0})) \\
  - (e_{2,\extr}(\lambda_2, \lambda_3) + S[\lambda_2+3]_\ell - S[\lambda_3+2]_\ell) \otimes S[\lambda_1+4]_\ell \\
  + (e_{2,\extr}(\lambda_1+1, \lambda_3) + S[\lambda_1+4]_\ell - S[\lambda_3+2]_\ell) \otimes S[\lambda_2+3]_\ell \\
  - (e_{2,\extr}(\lambda_1+1, \lambda_2+1) + S[\lambda_1+4]_\ell - S[\lambda_2+3]_\ell) \otimes S[\lambda_3+2]_\ell
\end{multline*}
and the terms \(\pm S[\lambda_i+5-i]_\ell \otimes S[\lambda_j+5-j]_\ell\) cancel each other out, concluding the proof of Theorem \ref{thm:BFG}.

\subsection{Trivial local systems: \(|\mathcal{A}_n(\Fq)|\) for small \(n\)}

In this section we prove Theorem \ref{thmintro:card_An_Fq_small_n}.
Consider a dominant weight \(\lambda = (\lambda_1 \geq \dots \geq \lambda_n \geq 0)\) for \(\Spbf_{2n}\), and the dominant weight \((\lambda,0) = (\lambda_1, \dots, \lambda_n, 0\) for \(\GSpbf_{2n}\) (recall our choice of parametrization from Section \ref{sec:not_red_gps}).
First we use Theorem \ref{theo:Hc_from_IH} in level one, forgetting the Hecke action.
From Example \ref{exam:Euler_L2_coh_GL} we easily deduce formulas for \(e(\GLbf_1, a)^{\Zhat^\times}\) and \(e_{(2)}(\GLbf_2, a, b)^{\GLbf_2(\Zhat)}\), and we deduce formulas (in the Grothendieck group of continuous finite-dimensional \(\ell\)-adic representations of \(\GalQ\)), for any \(\lambda_1 \geq \dots \lambda_n \geq 0\), expressing \(e_c(\Acal_{n, \Qbar}, \Fcal_\ell(V_{\lambda,0}))\) in terms of
\[ e_\IH(\lambda') := e_\IH(\GSpbf_{2n'}, \Xcal_{n'}, V_{\lambda',0})^{\GSpbf_{2n'}(\Zhat)} \]
where \(n' \leq n\) and \(\lambda'_1 + n' \leq \lambda_1 + n\), more precisely as a linear combination with integral coefficients of \(e_\IH(\lambda') \chi_\ell^{-N}\) (where \(N \geq 0\) is an integer).
Next we apply Theorems \ref{theo:IH_explicit_crude} and \ref{thm:sigma_is_tensor} to decompose
\[ e_\IH(\lambda') = e_\IH(\GSpbf_{2n'}, \Xcal_{n'}, V_{\lambda',m'})^{\GSpbf_{2n'}(\Zhat)} \, \chi_\ell^{-m'}, \]
where \(2m' = \sum_i \lambda'_i\) (so that the central character of \(V_{\lambda',m'}\) is trivial), by parameter \(\psi \in \Psit^{\unr,\tau'}_{\disc}(\Spbf_{2n'})\), where \(\tau' = (\lambda'_1 + n', \dots, \lambda'_{n'}+1)\) is the infinitesimal character for associated to \(\lambda'\).
Each such parameter \(\psi\) decomposes as
\[ \psi_0 \oplus \dots \oplus \psi_r = \pi_0[d_0] \oplus \dots \oplus \pi_r[d_r], \]
and for \(\lambda_1' + n' \leq 12\) the classification theorems \cite[Theorem 9.3.3]{CheLan} and \cite[Theorems 3 and 4]{ChenevierTaibi} tell us that each \(\pi_i\) belongs to an explicit (short) list of possibilities.
For the cases at hand we even have \(\lambda_1' + n' \leq 7\), so one of the following holds:
\begin{itemize}
\item \(\psi = \psi_0 = [2d+1]\) with \(d=n'\), and \(\lambda' = (0, \dots, 0)\),
\item \(\psi = \psi_0 + \psi_1 = [2d+1] + \Delta_{11}[2]\) with \(d=n'-2\) and \(\lambda' = (6-n', 6-n', 0, \dots, 0)\) (this can occur only if \(n' \geq 2\)), where \(\Delta_{11}\) is the unique level one cuspidal automorphic representation for \(\PGLbf_2\) with infinitesimal character \((11/2,-11/2)\) (corresponding to the unique eigenform in \(S_{12}(\SL_2(\Z))\))
\item \(\psi = \psi_0 + \psi_1 = [2d+1] + \Delta_{11}[4]\) with \(d = n' - 4\) and \(\lambda' = (7-n', 7-n', 7-n', 7-n', 0, \dots, 0)\) (this can occur only if \(n' \geq 4\)).
\end{itemize}

We give more details for the case \(n=6\): the cases where \(n<6\) are easier, and we briefly discuss the case \(n=7\) below.
Discarding highest weights \(\lambda'\) which are not of the form \((0,\dots,0)\) or \((6-n',6-n',0,\dots,0)\), we\footnote{In fact we had this computation done by a computer.} obtain
\begin{align*}
  e_c(\Acal_{6,\Qbar}, \Qell) =&\ e_\IH(0,0,0,0,0,0) - e_\IH(0,0,0,0,0) - e_\IH(0,0,0,0) \chi_\ell^{-5} \\
                               &\ + e_\IH(0,0,0) (1 + \chi_\ell^{-5}) - e_\IH(0,0) - e_\IH(0) \chi_\ell^{-5} - e_\IH(1,1,0,0,0) \\
                               &\ + e_\IH(3,3,0) - e_\IH(4,4).
\end{align*}
For \(\lambda' = (0,\dots,0)\) the contribution of the parameter \([2n'+1]\) to \(e_\IH(\lambda')\) is
\[ \prod_{i=1}^{n'} (1 + \chi_\ell^{-i}). \]
For \(n' \in \{2,3,5,6\}\) and \(\lambda' = (6-n',6-n',0,\dots,0)\) we also need to compute the contribution of the parameter \(\psi = [2d+1] \oplus \Delta_{11}[2]\) (here \(d = n'-2\)) to \(e_\IH(\lambda')\).
With notation as in Theorem \ref{theo:IH_explicit_crude} we have
\[ u_1(\psi) = - \epsilon(1/2, \Delta_{11})^{\min(2,2d+1)} = -1. \]
Using Corollary \ref{cor:spin_pm_SO4} we compute (for any \(\iota\))
\[ \sigma_{\psi_1,\iota}^{\spin,u_1(\psi)} = \sigma_{\psi_1,\iota}^{\spin,-} = 1 + \chi_\ell^{-1} \]
and so by Theorem \ref{thm:sigma_is_tensor} the contribution of \(\psi\) to \(e_\IH(\lambda')\) is
\[ \chi_\ell^{n'-6} \sigma_{\psi,\iota}^{\IH} = \chi_\ell^{-5} (1 + \chi_\ell^{-1}) \times \prod_{i=1}^d (1 + \chi_\ell^{-i}). \]
Adding all contributions we obtain \(e_c(\Acal_{6,\Qbar}, \Qell) = P_6(\chi_\ell^{-1})\) (defined in Theorem \ref{thmintro:card_An_Fq_small_n}), which is equivalent (via the Grothendieck-Lefschetz trace formula) to the formula \(|\Acal_6(\Fq)| = P_6(q)\) for all prime powers \(q\).

For \(n=7\) we also have to consider parameters of the form
\[ \psi = \psi_0 \oplus \psi_1 = [2d+1] \oplus \Delta_{11}[4] \]
for which we compute
\[ u_1(\psi) = \epsilon(1/2, \Delta_{11})^{\min(4,2d+1)} = +1. \]
A simple weight computation shows that in Proposition \ref{pro:Gal_GSpin_pi_2d} we have
\[ \spin_{\psi_1}^+ \circ \tilde{\alpha}_{\psi_1}|_{\Sp_2 \times \SL_2} \simeq \left( (\Sym^2 \Std_{\Sp_2}) \otimes 1 \right) \oplus \left( 1 \otimes \Sym^4 \Std_{\SL_2} \right) \]
and we deduce
\[ \sigma_{\psi_1,\iota}^{\spin,u_1(\psi)} = \sigma_{\psi_1,\iota}^{\spin,+} = \Sym^2 \rho_{\Delta_{11}, \iota} + \sum_{i=-1}^3 \chi_\ell^{-i}. \]
The contribution of \(\psi\) to \(e_\IH(\lambda')\) is thus
\[ \chi_\ell^{2n'-14} \sigma_{\psi,\iota}^\IH = \chi_\ell^{-10} \sigma_{\psi_1,\iota}^{\spin,+} \times \prod_{i=1}^d 1 + \chi_\ell^{-i}. \]
We have
\[ \Tr \chi_\ell^{-10} \Sym^2 \rho_{\Delta_{11},\iota}(\Frob_p^m) = \iota \left( p^{11m} \Tr \Sym^2 c_p(\Delta_{11}) \right) \]
where the semi-simple conjugacy class \(c_p(\Delta_{11})\) in \(\SL_2(\C)\) is determined by
\[ \Tr c_p(\Delta_{11}) = p^{-11/2} \tau(p) \]
(recall the Ramanujan \(\tau\) function from Theorem \ref{thmintro:card_An_Fq_small_n}).
The well-known relation \(\tau(p)^2 = \tau(p^2) + p^{11}\) and elementary computations give
\[ \Tr \chi_\ell^{-10} \Sym^2 \rho_{\Delta_{11},\iota}(\Frob_p^m) = a(p^m). \]
We omit the details leading to the formula for \(|\Acal_7(\Fq)|\) in Theorem \ref{thmintro:card_An_Fq_small_n}.

\begin{rema} \label{rem:card_An_Fq_upto12}
  In principle the classification theorems \cite[Theorem 9.3.3]{CheLan} and \cite[Theorems 3 and 4]{ChenevierTaibi} allow us to obtain explicit formulas for \(e_c(\Acal_{n, \Qbar}, \Fcal_\ell(V_{\lambda,0}))\) for \(\lambda_1 + n \leq 12\).
\end{rema}

\begin{rema}
  \begin{enumerate}
  \item We checked that the dimensions of the Euler characteristics \((e_c(\Acal_{n,\Qbar}, \Qell))_{1 \leq n \leq 7}\) (for \(n \leq 6\) this amounts to evaluating the polynomial \(P_n\) at \(q=1\)) coincide with the values that we computed independently using the trace formula (see \cite[Appendix, Proposition 4]{HulekTommasi_survey} and \cite{Taibi_dimtrace}).
  \item The method explained above to compute \(e_c(\Acal_{n,\Qbar}, \Qell)\) for small \(n\) clearly also works to compute \(e_c(\Acal_{n,\Qbar}, \Fcal(V_{\lambda,0}))\) for small \(\lambda_1 + n\).
    For example we checked that we recover \cite[Theorem 8.1]{BFG}.
    %TODO: reference code
  \end{enumerate}
\end{rema}

\newpage

\appendix

\section{Cohomological correspondences}
\label{app:corr}

\subsection{Definitions and induced maps in cohomology}
\label{sec:corr_def}

Let \(\ell\) be a prime number.
Let \(B\) be a scheme over \(\Z[1/\ell]\) assumed\footnote{More generally one could assume that \(B\) is a regular excellent Noetherian scheme, or an excellent Noetherian scheme endowed with a dimension function, see \cite[Exposé XVII]{Gabber_unifloc}.
We will not need this generality.} to be Noetherian, regular and of dimension \(\leq 1\).
Fix a dualising complex on \(B\), thus determining dualising complexes on schemes separated of finite type over \(B\).
Let \(\pi_i : X_i \rightarrow B\) for \(i=1,2\) be separated finite type schemes over \(B\), and \(L_i \in D^b_c(X_i, \Qell)\).
Let \(\pi : X \rightarrow B\) be also separated of finite type, and suppose that \(c_i : X \rightarrow X_i\) for \(i=1,2\) are morphisms over \(B\).
Recall that a correspondence from \(L_1\) to \(L_2\) with support in \((c_1, c_2)\) is a morphism \(u : c_1^* L_1 \rightarrow c_2^! L_2\).

The dual \(\Dbb(u)\) of \(u\) is a correspondence from \(\Dbb(L_2)\) to \(\Dbb(L_1)\)
with support in \((c_2, c_1)\) defined in the obvious way using the
identifications
\begin{align*}
  & \Hom(c_1^* L_1, c_2^! L_2) \simeq H^0(R \Gamma (\RuHom(c_1^* L_1, c_2^! L_2))) \ \ \text{ and} \\
  & \RuHom(c_1^* L_1, c_2^! L_2) \simeq \Dbb( c_1^* L_1 \otimes^L \Dbb(c_2^! L_2)) \simeq \Dbb( c_2^* \Dbb(L_2) \otimes^L c_1^! \Dbb(L_1)).
\end{align*}

If \(c_1\) (resp.\ \(c_2\)) is proper then \(u\) induces a morphism ``in cohomology''
\begin{align}
  u_! : \pi_{1!} L_1 & \rightarrow \pi_{2!} L_2 \label{eq:corr_in_ord_coh} \\
  (\text{resp.\ } u_* : \pi_{1*} L_1 & \rightarrow \pi_{2*} L_2), \label{eq:corr_in_cpct_coh}
\end{align}
see \cite[(1.3.2)]{Fujiwara_Deligne_conj} and \cite[(1.3)]{Pink_localterms} (see also \cite{Varshavsky_LefschetzVerdier}).
For example in the first case \(u_!\) is obtained as the composition
\[ \pi_{1!} L_1 \rightarrow \pi_{1!} c_{1*} c_1^* L_1 \simeq \pi_{2!} c_{2!} c_1^* L_1 \xrightarrow{\pi_{2!} c_{2!} * u} \pi_{2!} c_{2!} c_2^! L_2 \rightarrow \pi_{2!} L_2 \]
where the first map is the unit for the pair of adjunct functors \((c_1^*,
c_{1*})\), the isomorphism follows from the composition rule for \(\cdot_!\) and
the fact that \(c_1\) is proper, and the last map is the counit for \((c_{2!},
c_2^!)\).
There is an obvious notion of morphism of correspondences having the same
support \((L_1, L_2, u) \to (K_1, K_2, v)\), and it is easy to check that if \(c_1\)
(resp.\ \(c_2\)) is proper then the following diagram commutes.
\[
  \begin{tikzcd}
    \pi_{1 !} L_1 \arrow[r, "{u_!}"] \arrow[d] & \pi_{2 !} L_2 \arrow[d] \\
    \pi_{1 !} K_1 \arrow[r, "{v_!}"] & \pi_{2 !} K_2
  \end{tikzcd}
  \text{ resp.\ }
  \begin{tikzcd}
    \pi_{1 *} L_1 \arrow[r, "{u_*}"] \arrow[d] & \pi_{2 *} L_2 \arrow[d] \\
    \pi_{1 *} K_1 \arrow[r, "{v_*}"] & \pi_{2 *} K_2
  \end{tikzcd}
\]

We leave it to the reader to check that the formation of \(u_!\) (resp.\ \(u_*\)) is compatible with duality: \(\Dbb(u_!) : \pi_{2*} \Dbb(L_2) \rightarrow \pi_{1*} \Dbb(L_1)\) coincides with \(\Dbb(u)_*\).

We briefly recall the definition of composition for correspondences (see also \cite[Exposé III \S 5.2]{SGA5}).
Consider a diagram of schemes separated of finite type over \(B\)
\[
  \begin{tikzcd}
    & \ar[dl, "{c_1}"] X' \ar[dr, "{c_2}"] & & \ar[dl, "{d_2}"] X'' \ar[dr, "{d_3}"] & \\
    X_1 & & X_2 & & X_3
  \end{tikzcd}
\]
and correspondences \(u: c_1^* L_1 \to c_2^! L_2\) and \(v: d_2^* L_2 \to d_3^! L_3\).
Denote \(p_1: X' \times_{X_2} X'' \to X'\) and \(p_2: X' \times_{X_2} X'' \to X''\) the two projections.
The composition \(v \circ u\) is the correspondence supported on \((c_1 p_1, d_3 p_2)\) defined as the composition
\[ p_1^* c_1^* L_1 \xrightarrow{p_1^*(u)} p_1^* c_2^! L_2 \to p_2^! p_{2!} p_1^* c_2^! L_2 \xrightarrow{BC^{-1}} p_2^! d_2^* c_{2!} c_2^! L_2 \to p_2^! d_2^* L_2 \xrightarrow{p_2^!(v)} p_2^! d_3^! L_3 \]
where as before unlabelled maps are (co)units and \(BC\) is the base change isomorphism \(d_2^* c_{2!} \simeq p_{2!} p_1^*\) (\cite[Exposé XVII Théorème 5.2.6]{SGA4-3}).
It is formal to check that this notion is compatible with cohomology (\eqref{eq:corr_in_ord_coh} and \eqref{eq:corr_in_cpct_coh}) when this makes sense.

\subsection{Base change}

We now discuss base change.
Let \(f : B' \rightarrow B\) be morphism such that one of the following holds:
\begin{enumerate}
  \item \(f\) is separated of finite type, or
  \item \(f\) is flat with geometrically regular fibres and \(B'\) is Noetherian
    and excellent.
\end{enumerate} 
Denote \(\pi_i' : X_i' \rightarrow B'\) and \(\pi' : X' \rightarrow B'\) the objects
obtained by base change, \(c_i' : X' \rightarrow X\) the natural morphisms
obtained from \(c_i\) and \(g_i : X_i' \rightarrow X_i\), \(g : X' \rightarrow X\)
parallel to \(f\).
There is a notion of base change for correspondences \(f^* : \Hom(c_1^*, c_2^!
L_2) \rightarrow \Hom(c_1'^* g_1^* L_1, c_2'^! g_2^* L_2)\) (resp.\ \(f^! :
\Hom(c_1^*, c_2^! L_2) \rightarrow \Hom(c_1'^* g_1^! L_1, c_2'^!  g_2^! L_2)\)
mapping \(u\) to \(f^* u\) (resp.\ \(f^! u\)) defined as the composition
\[ c_1'^* g_1^* L_1 \simeq g^* c_1^* L_1 \rightarrow c_2'^! c_{2!}' g^* c_1^* L_1 \xrightarrow{BC^{-1}} c_2'^! g_2^* c_{2!} c_1^* L_1 \xrightarrow{u} c_2'^! g_2^* c_{2!} c_2^! L_2 \rightarrow c_2'^! g_2^* L_2, \]
resp.\
\[ c_1'^* g_1^! L_1 \rightarrow c_1'^* g_1^! c_{1*} c_1^* L_1 \xrightarrow{u}
  c_1'^* g_1^! c_{1*} c_2^! L_2 \xrightarrow{\Dbb(BC)} c_1'^* c_{1*}' g^! c_{2!}
  L_2 \rightarrow g^! c_{2!} \simeq c_2'^! g_2^! L_2. \]
Unsurprisingly, \(\Dbb(f^* u) = f^! \Dbb(u)\).

\begin{lemm}
  For \(f : B' \rightarrow B\) and \(u \in \Hom(c_1^* L_1, c_2^! L_2)\) as above, if
  \(c_1\) is proper then \((f^* u)_! = f^* (u_!)\), i.e.\ \((f^* u)_! : \pi_{1!}'
  g_1^* L_1 \rightarrow \pi_{2!}' g_2^* L_2\) is obtained from \(u_!\) by applying
  \(f^*\) and via the base change isomorphisms \(f^* \pi_{i!} \simeq \pi_{i!}'
  g_i^*\).
\end{lemm}
\begin{proof}
  Details are left to the reader, essentially uses compatibility of proper base
  change with composition \cite[Exposé XVII Lemme 5.2.4]{SGA4-3} and the
  unit/counit relations for pairs of adjoint functors.
\end{proof}
Of course if \(c_2\) is proper the dual assertion \((f^! u)_* = f^! (u_*)\) is a
direct consequence.

\subsection{Pushforward and pullback}
\label{sec:corr_push_pull}

We now discuss proper pushforward and étale pullback (for \(X\)), and
compatibility with base change and induced morphisms in cohomology.
We do not consider the most general situation.

Consider a commutative diagram of schemes separated of finite type over \(B\):
\begin{equation} \label{eq:diagram_pushpull_corr}
  \begin{tikzcd}
    & X' \arrow[ddl, bend right=20, "{c_1'}" above] \arrow[ddr, bend left=20,
      "{c_2'}"] \arrow[d, "{f}"] & \\
    & X  \arrow[dl, "{c_1}"] \arrow[dr, "{c_2}"]& \\
  X_1 & & X_2
  \end{tikzcd}
\end{equation}
Under hypotheses on \(f\), correspondences with support in \((c_1, c_2)\) and
\((c_1', c_2')\) can be related in both directions (see \cite[\S
1.4]{Fujiwara_Deligne_conj} for the proper pushforward).
\begin{enumerate}
  \item If \(f\) is proper, we have a pushforward morphism \(\corr f_* :
    \Hom(c_1'^* L_1, c_2'^! L_2) \rightarrow \Hom(c_1^* L_1, c_2^! L_2)\): for \(u
    \in \Hom(c_1'^* L_1, c_2'^! L_2)\), \(\corr f_* u\) is obtained as the
    composition
    \[ c_1^* L_1 \rightarrow f_* f^* c_1^* L_1 = f_* c_1'^* L_1
      \xrightarrow{f_*(u)} f_* c_2'^! L_2 = f_! f^! c_2^! L_2 \rightarrow c_2^!
      L_2 \]
    where the first and last morphisms are given by unit and counit of
    adjunctions.
    As before denote \(\pi_i : X_i \rightarrow B\).
    If \(c_1\) is proper (resp.\ \(c_2\) is proper) then \((\corr f_* u)_* \in
    \Hom(\pi_{1*} L_1, \pi_{2*} L_2)\) equals \(u_*\) (resp.\ \((\corr f_* u)_! \in
    \Hom(\pi_{1!} L_1, \pi_{2!} L_2)\) equals \(u_!\)).
    This follows from the compatibility of adjunctions with compositions.
  \item If \(f\) is étale and \(u \in \Hom(c_1^* L_1, c_2^! L_2)\), define \(\corr
    f^* u\) as the composition
    \[ c_1'^* L_1 \simeq f^* c_1^* L_1 \simeq f^! c_1^* L_1 \xrightarrow{f^!(u)}
    f^! c_2^!  L_2 \simeq c_2'^! L_2. \]
    For \(f\) étale this defines a pullback morphism
    \[ \corr f^* : \Hom(c_1^* L_1, c_2^! L_2) \longrightarrow \Hom(c_1'^* L_1,
    c_2'^! L_2). \]
    When \(f\) is an open immersion \(\corr f^*\) is just the restriction morphism.
\end{enumerate}

\begin{lemm} \label{lemm:fet_push_pull_corr}
  In the situation of \eqref{eq:diagram_pushpull_corr} above, if \(f\) is finite étale of constant degree \(N\) then for any correspondence \(u\) with support in \((c_1, c_2)\) we have \(\corr f_*  (\corr f^* u) = N u\).
\end{lemm}
\begin{proof}
  This follows from the fact that the composition \(\id \rightarrow f_* f^* \simeq f_! f^! \rightarrow \id\) is multiplication by \(N\) (\cite[Exposé XVIII Proposition 3.1.8(iii)]{SGA4-3} and \cite[Exposé XVII Théorème 6.2.3 (Var 4)]{SGA4-3}).
\end{proof}

\subsection{More pushforwards and pullbacks}
\label{sec:more_corr_push_pull}

The pushforward in the previous section admits the following variant.

\begin{defi} \label{def:gen_push_corr}
  Consider a commutative diagram of qcqs schemes
  \begin{equation} \label{diag:gen_push_corr}
    \begin{tikzcd}
      X_1 \ar[d, "{f_1}"] & \ar[l, "{c_1}"] X' \ar[d, "{f}"] \ar[r, "{c_2}"] & X_2 \ar[d, "{f_2}"] \\
      Y_1 & \ar[l, "{d_1}"] Y' \ar[r, "{d_2}"] & Y_2
    \end{tikzcd}
  \end{equation}
  in which all morphisms are separated of finite type and \(d_2\) and \(c_2\) are proper.
  Let \(u: c_1^* L_1 \to c_2^! L_2\) be a correspondence supported on \((c_1,c_2)\).
  Define a correspondence from \(f_{1*} L_1\) to \(f_{2*} L_2\) supported on \((d_1,d_2)\):
  \[ d_1^* f_{1*} L_1 \to f_* c_1^* L_1 \xrightarrow{f_*(u)} f_* c_2^! L_2 \to d_2^! d_{2!} f_* c_2^! L_2 \simeq d_2^! f_{2*} c_{2!} c_2^! L_2 \to d_2^! f_{2*} L_2 \]
  where the first map is obtain from two adjunctions, the third map is an adjunction, the fourth map follows from \(d_{2,!} \simeq d_{2,*}\) and \(c_{2,*} \simeq c_{2,!}\) because \(d_2\) and \(c_2\) are proper and the last map is also an adjunction.
  Denoting \(\ul{f} = (f_1,f,f_2)\), we will denote by \(\corr \ul{f}_* u\) this correspondence.
\end{defi}

It is often more convenient to see the correspondence \(u\) as a morphism \(u: c_{2!} c_1^* L_1 \to L_2\), by adjunction, and then its pushforward \(\corr \ul{f}_* u\) by \(\ul{f} = (f_1,f,f_2)\), seen as a morphism \(d_{2!} d_1^* f_{1*} L_1 \to f_{2*} L_2\), is equal to the composition
\[ d_{2!} d_1^* f_{1*} L_1 \to d_{2!} f_* c_1^* L_1 \simeq f_{2*} c_{2!} c_1^* L_1 \xrightarrow{f_{2*}(u)} f_{2*} L_2. \]
As usual this follows from the adjunction formalism.

\begin{rema}
  Assume that we have a commutative diagram \eqref{diag:gen_push_corr}, and instead of assuming that \(d_2\) and \(c_2\) are proper, assume that \(f\) and \(f_2\) are proper.
  In this situation we may define, for a correspondence \(u: c_1^* L_1 \to c_2^! L_2\), its pushforward along \(\ul{f}\) as above, deriving \(d_{2!} f_* \simeq f_{2*} c_{2!}\) from \(f_* = f_!\) and \(f_{2*} = f_{2!}\).
  If all four morphisms \(c_2\), \(d_2\), \(f\) and \(f_2\) are proper then the two notions of pushforward coincide.
  In particular when \(f_1 = \id\), \(f_2 = \id\) and \(f\) and \(d_2\) are proper the pushforward defined in the present section is equal to the pushforward defined in Section \ref{sec:corr_push_pull}.

  In practice the assumption that \(c_2\) and \(d_2\) are proper is always satisfied (at least in this article), whereas the vertical morphisms are not always proper.
\end{rema}

As in the previous section, pushfoward of correspondences is compatible with cohomology, as the following proposition shows.

\begin{prop}[{Compare \cite[\S 1.4]{Fujiwara_Deligne_conj}}] \label{pro:compat_gen_push_corr_coh}
  In the setting of Definition \ref{def:gen_push_corr}, assume that the diagram is a diagram of schemes over \(B\), and denote \(\pi_i: Y_i \to B\) and \(\pi: Y \to B\).
  Then the morphisms \(u_*: (\pi_1 f_1)_* L_1 \to (\pi_2 f_2)_* L_2\) and \((\corr \ul{f}_* u)_*: \pi_{1*} f_{1*} L_1 \to \pi_{2*} f_{2*} L_2\) are equal.
\end{prop}
\begin{proof}
  Here it is convenient to see \(u\) as a morphism \(c_{2!} c_1^* L_1 \to L_2\) and similarly for \(\corr \ul{f}_* u\).
  Writing the ``base change'' morphism \(d_1^* f_{1*} \to f_* c_1^*\) as the composition
  \[ d_1^* f_{1*} \to d_1^* f_{1*} c_{1*} c_1^* \simeq d_1^* d_{1*} f_* c_1^* \to f_* c_1^* \]
  and plugging this in the definition of \((\corr \ul{f}_* u)_*\), we obtain a long composition where the unit and counits for \((d_1^*, d_{1*})\) both appear and may be eliminated.
  Details are left to the reader.
\end{proof}

We now recall from \cite[\S 5]{MorelSiegel1} a definition of pullback for correspondences.

\begin{defi} \label{def:gen_pull_corr}
  Consider a commutative diagram of qcqs schemes
  \begin{equation} \label{diag:gen_pull_corr}
    \begin{tikzcd}
      X_1 \ar[d, "{f_1}"] & \ar[l, "{c_1}"] X' \ar[d, "{f}"] \ar[r, "{c_2}"] & X_2 \ar[d, "{f_2}"] \\
      Y_1 & \ar[l, "{d_1}"] Y' \ar[r, "{d_2}"] & Y_2
    \end{tikzcd}
  \end{equation}
  in which all morphisms are separated of finite type and assume that the right square is Cartesian up to nilpotents.
  Let \(u: d_1^* L_1 \to d_2^! L_2\) be a correspondence supported on \((d_1,d_2)\).
  Define its pullback \(\corr \ul{f}^* u\) as the composition
  \[ c_1^* f_1^* L_1 \simeq f^* d_1^* L_1 \xrightarrow{f^* u} f^* d_2^! L_2 \to c_2^! f_2^* L_2 \]
  where the morphism of functors \(f^* d_2^! \to c_2^! f_2^*\) is \cite[Exposé XVIII (3.1.14.2)]{SGA4-3}, i.e.\ it is obtained by adjunction from the base change isomorphism \(f_2^* d_{2!} \simeq c_{2!} f^*\).
\end{defi}

As for pushforwards, this notion of pullback coincides with the one defined in the previous section when both make sense, but proving this requires some work.
The general compatibilities in the following lemma are probably folklore.

\begin{lemm} \label{lem:BC_for_et_using_trace}
  Let
  \[
    \begin{tikzcd}
      X' \ar[r, "{f}"] \ar[d, "{c}"] & Y' \ar[d, "{d}"] \\
      X \ar[r, "{g}"] & Y
    \end{tikzcd}
  \]
  be a Cartesian diagram of qcqs schemes in which all morphisms are separated of finite type.
  \begin{enumerate}
  \item Assume that \(d\) is proper (whence also \(c\)).
    Then the isomorphism \(g_! c_* \simeq d_* f_!\) (see \cite[Exposé XVII \S 3.3.2.3]{SGA4-3}) is equal to the composition
    \[ g_! c_* \xrightarrow{\adj} d_* d^* g_! c_* \xrightarrow[\sim]{\BC} d_* f_! c^* c_* \xrightarrow{\adj} d_* f_! \]
  \item Assume that \(g\) is étale, so that we have an adjoint pair \((g_!,g^*)\) (see \cite[Exposé XVII Proposition 6.2.11]{SGA4-3}), and similarly for \(f\).
    Then the composition
    \[ c^* \xrightarrow{\adj} c^* g^* g_! \simeq f^* d^* g_! \xrightarrow[\sim]{\BC} f^* f_! c^* \]
    is equal to the unit \(\id \to f^* f_!\) applied to \(c^*\).
  \item Assume that \(g\) is étale (whence also \(f\)).
    Then the base change isomorphism \(c_! f^* \simeq g^* d_!\) is equal to the composition
    \[ c_! f^* \xrightarrow{\adj} g^* g_! c_! f^* \simeq g^* d_! f_! f^* \xrightarrow{\adj} g^* d_!. \]
  \item Assume that \(g\) is étale (whence also \(f\)).
    Then the morphism of functors \(f^* d^! \to c^! g^*\) defined in \cite[Exposé XVIII (3.1.14.2)]{SGA4-3} is equal to
    \[ f^* d^! \simeq f^! d^! \simeq (d f)^! = (g c)^! \simeq c^! g^! \simeq c^! g^*. \]
  \end{enumerate}
\end{lemm}
\begin{proof}
  \begin{enumerate}
  \item Writing \(g\) as the composition of an open immersion and a proper morphism, we are reduced to proving that our two compositions coincide in the following cases.

    If \(g\) is proper then \(g_! = g_*\) and \(f_! = f_*\) and we are left with an exercise in adjunction.
    Details are left to the reader.

    If \(g\) is an open immersion then the isomorphism \(g_! c_* \simeq d_* f_!\) is defined (see \cite[Exposé XVII (5.1.5.2)]{SGA4-3}) as
    \[ g_! c_* \xrightarrow[\sim]{\adj} g_! c_* f^* f_! \xleftarrow[\sim]{\BC} g_! g^* d_* f_! \xrightarrow[\sim]{\adj} d_* f_! \]
    and it follows readily from the definition (see Lemme 5.1.2 loc.\ cit.\ and its proof) that the base change isomorphism \(f_! c^* \simeq d^* g_!\) is equal to
    \[ f_! c^* \xrightarrow[\sim]{\adj} f_! c^* g^* g_! \simeq f_! f^* d^* g_! \xrightarrow{\adj} d^* g_!. \]
    So we are left to check that the two paths from \(g_! c_*\) to \(d_* f_!\) in the diagram
    \[
      \begin{tikzcd}
        g_! c_* \ar[r, "{\sim}" below, "{\adj}" above] \ar[d, "{\adj}"] & g_! c_* f^* f_! & \ar[l, "{\BC}" above, "{\sim}" below] g_! g^* d_* f_! \ar[r, "{\adj}" above, "{\sim}" below] & d_* f_! \\
        d_* d^* g_! c_* & \ar[l, "{\adj}" above, "{\sim}" below] d_* f_! f^* d^* g_! c_* \ar[r, dash, "{\sim}"] & d_* f_! c^* g^* g_! c_* & \ar[l, "{\adj}" above, "{\sim}" below] d_* f_! c^* c_* \ar[u, "{\adj}"]
      \end{tikzcd}
    \]
    are equal.
    The bottom path is also equal to
    \[ g_! c_* \xrightarrow[\sim]{\adj} g_! c_* f^* f_! \xrightarrow{\adj} d_* d^* g_! c_* f^* f_! \simeq d_* f_! c^* g^* g_! c_* f^* f_! \xleftarrow[\sim]{\adj} d_* f_! c^* c_* f^* f_! \xrightarrow{\adj} d_* f_! f^* f_! \xrightarrow[\sim]{\adj} d_* f_! \]
    so we have to show that the (anti-clockwise) cycle from \(d_* f_!\) to itself in the diagram
    \[
      \begin{tikzcd}
        g_! c_* f^* f_! \ar[d, "{\adj}"] & \ar[l, "{\BC}" above, "{\sim}" below] g_! g^* d_* f_! \ar[r, "{\adj}" above, "{\sim}" below] & d_* f_! & \ar[l, "{\adj}" above, "{\sim}" below] d_* f_! f^* f_! \\
        d_* d^* g_! c_* f^* f_! & \ar[l, "{\adj}" above, "{\sim}" below] d_* f_! f^* d^* g_! c_* f^* f_! \ar[r, dash, "{\sim}"] & d_* f_! c^* g^* g_! c_* f^* f_! & \ar[l, "{\adj}" above, "{\sim}" below] d_* f_! c^* c_* f^* f_! \ar[u, "{\adj}"]
      \end{tikzcd}
    \]
    is the identity.
    We observe that this diagram is obtained by applying a similar diagram to \(f_!\), and so it is enough to show that the composition
    \begin{multline*}
      g_! g^* d_* \xrightarrow[\sim]{\BC} g_! c_* f^* \xrightarrow{\adj} d_* d^* g_! c_* f^* \xleftarrow[\sim]{\adj} d_* f_! f^* d^* g_! c_* f^* \\
      \simeq d_* f_! c^* g^* g_! c_* f^* \xleftarrow[\sim]{\adj} d_* f_! c^* c_* f^* \xrightarrow{\adj} d_* f_! f^* \xrightarrow{\adj} d_*
    \end{multline*}
    is simply equal to the counit \(g_! g^* \to \id\) applied to \(d_*\).
    We may reorder:
    \begin{multline*}
      g_! g^* d_* \xrightarrow{\adj} d_* d^* g_! g^* d_* \xleftarrow[\sim]{\adj} d_* f_! f^* d^* g_! g^* d_* \\
      \simeq d_* f_! c^* g^* g_! g^* d_* \xleftarrow[\sim]{\adj} d_* f_! c^* g^* d_* \xrightarrow[\sim]{\BC} d_* f_! c^* c_* f^* \xrightarrow{\adj} d_* f_! f^* \xrightarrow{\adj} d_*.
    \end{multline*}
    Using the fact that the inverse of \(g^* \xrightarrow{\adj} g^* g_! g^*\) (counit applied to \(g^*\)) is \(g^* g_! g^* \xrightarrow{\adj} g^*\) (\(g^*\) applied to unit), and the equality of compositions
    \[ \left( c^* g^* d_* \xrightarrow{\BC} c^* c_* f^* \xrightarrow{\adj} f^* \right) = \left( c^* g^* d_* \simeq f^* d^* d_* \xrightarrow{\adj} f^* \right) \]
    we obtain
    \begin{multline*}
      g_! g^* d_* \xrightarrow{\adj} d_* d^* g_! g^* d_* \xleftarrow[\sim]{\adj} d_* f_! f^* d^* g_! g^* d_* \\
      \simeq d_* f_! c^* g^* g_! g^* d_* \xrightarrow{\adj} d_* f_! c^* g^* d_* \simeq d_* f_! f^* d^* d_* \xrightarrow{\adj} d_* f_! f^* \xrightarrow{\adj} d_*
    \end{multline*}
    which is simply
    \[ g_! g^* d_* \xrightarrow{\adj} d_* d^* g_! g^* d_* \xleftarrow[\sim]{\adj} d_* f_! f^* d^* g_! g^* d_* \xrightarrow{\adj} d_* f_! f^* d^* d_* \xrightarrow{\adj} d_* f_! f^* \xrightarrow{\adj} d_*. \]
    Moving the counit \(f_! f^* \to \id\) to the left simplifies the composition:
    \[ g_! g^* d_* \xrightarrow{\adj} d_* d^* g_! g^* d_* \xrightarrow{\adj} d_* d^* d_* \xrightarrow{\adj} d_* \]
    and one last reordering and unit/counit relation for \((d^*,d_*)\) yields the result.

  \item
    We need to prove that the composition
    \[ c^* \xrightarrow{\adj} c^* g^* g_! \simeq f^* d^* g_! \xrightarrow[\sim]{\BC} f^* f_! c^* \]
    is equal to the unit \(\id \to f^* f_!\) applied to \(c^*\).
    All functors occurring are derived from exact functors, as are the morphisms between them, so we may check this equality for sheaves, and for this it is enough to prove equality on stalks.
    For \(x'\) a geometric point of \(X'\) and \(\Fcal\) an abelian sheaf on \(X\) we have a commutative diagram
    \[
      \begin{tikzcd}
        (c^* \Fcal)_{x'} \ar[d, dash, "{\sim}"] \ar[r] & (c^* g^* g_! \Fcal)_{x'} \ar[r, "{\sim}"] \ar[d, dash, "{\sim}"] & (f^* d^* g_! \Fcal)_{x'} \ar[d, dash, "{\sim}"] \ar[r, "{\BC}" above, "{\sim}" below] & (f^* f_! c^* \Fcal)_{x'} \ar[d, "{\sim}"] \\
        \Fcal_{c(x')} \ar[r] & \bigoplus\limits_{\substack{x \\ g(x) = g(c(x'))}} \Fcal_x \ar[r, dash, "{\sim}"] & \bigoplus\limits_{\substack{x \\ g(x) = d(f(x'))}} \Fcal_x \ar[r] & \bigoplus\limits_{\substack{t \\ f(t) = f(x')}} \Fcal_{c(t)}
      \end{tikzcd}
    \]
    where \(x\) denotes a geometric point of \(X\) and \(t\) denotes a geometric point of \(X'\), the first map on the bottom line is the obvious inclusion and the third map is induced by the bijection \(t \mapsto c(t)\) from lifts of \(f(x')\) along \(f\) to lifts of \(d(f(x'))\) along \(g\), by the Cartesian property.
    We see that the composition at the bottom coincides with the adjunction \((c^* \Fcal)_{x'} \to (f^* f_! c^* \Fcal)_{x'}\) for \((f_!, f^*)\).

  \item Spelling out the definition of \(d_!\) and \(c_!\) by writing \(d\) as the composition of an open immersion and a proper morphism, we are reduced to proving the statement in two special cases: if \(d\) is an open immersion or if \(d\) is proper.

    If \(d\) is an open immersion then all functors involved are derived from exact functors on sheaves and equality may be checked on stalks.
    (Alternatively we may check equality after applying \(c^*\) because both source and target vanish on the complement of \(X'\) in \(X\), and this is rather formal.)

    So we assume that \(d\) is proper for the rest of the proof, and we need to show that the composition
    \[ c_* f^* \xrightarrow{\adj} g^* g_! c_* f^* \simeq g^* d_! f_! f^* \xrightarrow{\adj} g^* d_* \xrightarrow[\sim]{\BC} c_* f^* \]
    is the identity.
    Thanks to the first point of the lemma this composition may be rewritten as
    \[ c_* f^* \xrightarrow{\adj} g^* g_! c_* f^* \xrightarrow{\adj} g^* d_* d^* g_! c_* f^* \xrightarrow[\sim]{\BC} g^* d_* f_! c^* c_* f^* \xrightarrow{\adj} g^* d_* f_! f^* \xrightarrow{\adj} g^* d_* \xrightarrow[\sim]{\BC} c_* f^*. \]
    Writing the second base change isomorphism as the composition
    \[ g^* d_* \xrightarrow{\adj} c_* c^* g^* d_* \simeq c_* f^* d^* d_* \xrightarrow{\adj} c_* f^*, \]
    reordering and noticing that the unit and counit for \((d^*,d_*)\) cancel each other out, we obtain
    \[ c_* f^* \xrightarrow{\adj} c_* c^* c_* f^* \xrightarrow{\adj} c_* c^* g^* g_! c_* f^* \simeq c_* f^* d^* g_! c_* f^* \xrightarrow[\sim]{\BC} c_* f^* f_! c^* c_* f^* \xrightarrow{\adj} c_* f^* f_! f^* \xrightarrow{\adj} c_* f^*. \]
    The second point of the lemma allows us to simplify the composition of the second, third and fourth morphisms:
    \[ c_* f^* \xrightarrow{\adj} c_* c^* c_* f^* \xrightarrow{\adj} c_* f^* f_! c^* c_* f^* \xrightarrow{\adj} c_* f^* f_! f^* \xrightarrow{\adj} c_* f^*. \]
    This is equal to the identity by unit/counit relations.

  \item The morphism of functors \(f^* d^! \to c^! g^*\) is, by definition, equal to the composition
  \[ f^* d^! \xrightarrow{\adj} c^! c_! f^* d^! \xrightarrow[\sim]{\BC} c^! g^* d_! d^! \xrightarrow{\adj} c^! g^* \]
  and by the previous point in the lemma this is equal to
  \[ f^* d^! \xrightarrow{\adj} c^! c_! f^* d^! \xrightarrow{\adj} c^! g^* g_! c_! f^* d^! \simeq c^! g^* d_! f_! f^* d^! \xrightarrow{\adj} c^! g^* d_! d^! \xrightarrow{\adj} c^! g^*. \]
  The first two units may be combined into the unit for \(((g c)_!, (g c)^!)\), expanded into two units for \((f_!, f^*)\) and \((d_!, d^!)\), compensating the last two counits.
  \end{enumerate}
\end{proof}

\begin{rema}
  It seems likely that in Lemma \ref{lem:BC_for_et_using_trace} one could replace ``étale'' by ``smooth'' (introducing the appropriate shifts and Tate twists), but of course our proof (the second point in particular) does not obviously extend to this case.
  It may be possible to give a more conceptual proof that applies to the smooth case as well, probably using tensor products.
\end{rema}

\begin{coro} \label{cor:compat_pull_corr}
  Assume that we have a commutative diagram
  \begin{equation}
    \begin{tikzcd}
      X_1 \ar[d, "{f_1}"] & \ar[l, "{c_1}"] X' \ar[d, "{f}"] \ar[r, "{c_2}"] & X_2 \ar[d, "{f_2}"] \\
      Y_1 & \ar[l, "{d_1}"] Y' \ar[r, "{d_2}"] & Y_2
    \end{tikzcd}
  \end{equation}
  of qcqs schemes, where all morphisms are separated of finite type and \(f\) and \(f_2\) are étale.
  Let \(u: d_1^* L_1 \to d_2^! L_2\) be a correspondence.
  We obtain a correspondence
  \[ c_1^* f_1^* L_1 \simeq f^* d_1^* L_1 \xrightarrow{f^* u} f^* d_2^! L_2 \simeq c_2^! f_2^* L_2 \]
  where we used \(f^* \simeq f^!\) and \(f_2^! \simeq f_2\), generalizing the notion of pullback defined in Section \ref{sec:corr_push_pull}.

  If the right square in the above diagram is Cartesian up to nilpotent elements, then this pullback coincides with \(\corr \ul{f}^* u\) (Definition \ref{def:gen_pull_corr}).
\end{coro}
\begin{proof}
  This follows directly from the last point in Lemma \ref{lem:BC_for_et_using_trace}.
\end{proof}

The compatibility in Corollary \ref{cor:compat_pull_corr} will be useful to apply the following lemma (compatibility of pullback and pushforward) in a setting where we simply consider pullback along open immersions.

\begin{prop} \label{pro:corr_push_pull_compat}
  Assume that we have commutative diagrams of qcqs schemes
  \[
    \begin{tikzcd}
      X_1 \ar[d, "{f_1}"] & \ar[l, "{c_1}"] X' \ar[d, "{f}"] \ar[r, "{c_2}"] & X_2 \ar[d, "{f_2}"] & X' \ar[r, "{g}"] \ar[d, "{f'}"] & X \ar[d, "{f}"] \\
      Y_1 & \ar[l, "{d_1}"] Y' \ar[r, "{d_2}"] & Y_2 & Y' \ar[r, "{h}"] & Y
    \end{tikzcd}
  \]
  in which all morphisms are separated of finite type, the right square of the left diagram and the right diagram are Cartesian up to nilpotents, and \(h\) is proper (and thus also \(g\)).
  Denote \(\ul{f} = (f_1,f,f_2)\) and \(\ul{f}' = (f_1, f', f_2)\).
  For any correspondence \(u: h^* d_1^* L_1 \to h^! d_2^! L_2\) we have \(\corr \ul{f}^* (\corr h_* u) = \corr g_* (\corr (\ul{f}')^* u)\).
\end{prop}
\begin{proof}
  Unwinding definitions we find that it is enough to prove commutativity for the following diagrams of functors, where the unlabelled morphisms are \cite[Exposé XVIII (3.1.14.2)]{SGA4-3}.
  \[
    \begin{tikzcd}
      f^* \ar[r, "{\adj}"] \ar[d, "{\adj}"] & f^* h_* h^* \ar[d, "{\BC}" right, "{\sim}" left] &
      f^* h_* h^! d_2^! \ar[d, "{\sim}" right, "{\BC}" left] \ar[r, "{\adj}"] & f^* d_2^! \ar[r] & c_2^! f_2^* \ar[d, equal] \\
      g_* g^* f^* \ar[r, dash, "{\sim}"] & g_* (f')^* h^* &
      g_* (f')^* h^! d_2^! \ar[r] & g_* g^! c_2^! f_2^* \ar[r, "{\adj}"] & c_2^! f_2^*
    \end{tikzcd}
  \]
  For the left diagram we expand the base change morphism as
  \[ f^* h_* \xrightarrow{\adj} g_* g^* f^* h_* \simeq g_* (f')^* h^* h_* \xrightarrow{\adj} g_* (f')^* \]
  and note that unit and counit for \((h^*,h_*)\) cancel each other.
  For the right diagram, expanding the definition of \((f')^* h^! d_2^! \to g^! c_2^! f_2^*\) involves the unit \(\id \to g^! g_!\) which cancels with the bottom right counit \(g_* g^! \to \id\), and by compatibility of base change with composition the composition via the bottom path is equal to
  \[ f^* h_* h^! d_2^! \xrightarrow{\adj} c_2^! c_{2!} f^* h_* h^! d_2^! \xrightarrow{\BC} c_2^! f_2^* d_{2!} h_* h^! d_2^! \xrightarrow{\adj} c_2^! f_2^* \]
  i.e.\ the same composition as the top path but with natural transformations applied in a different order.
\end{proof}

\begin{prop} \label{pro:compat_pullback_compo_corr}
  Consider a commutative diagram of qcqs schemes where all morphisms are separated of finite type
  \[
    \begin{tikzcd}
      U_1 \ar[d, "{f_1}"] & \ar[l, "{c_1}"] U' \ar[d, "{f'}"] \ar[r, "{c_2}"] \ar[dr, phantom, "{\square}"] & U_2 \ar[d, "{f_2}"] & \ar[l, "{c_3}"] U'' \ar[d, "{f''}"] \ar[r, "{c_4}"] \ar[dr, phantom, "{\square}"] & U_3 \ar[d, "{f_3}"] \\
      X_1 & \ar[l, "{d_1}"] X' \ar[r, "{d_2}"] & X_2 & \ar[l, "{d_3}"] X'' \ar[r, "{d_4}"] & X_3
    \end{tikzcd}
  \]
  in which the two marked squares are Cartesian.
  Denote \(q_1\), \(q_2\), \(p_1\), \(p_2\) and \(g\) the morphisms
  \[
    \begin{tikzcd}
      U' \times_{U_2} U'' \ar[dd, "{q_1}"] \ar[rr, "{q_2}"] \ar[dr, "{g}"] & & U'' \ar[dd, "{c_3}"] \ar[dr, "{f''}"] & \\
      & X' \times_{X_2} X'' \ar[rr, "{p_2}" near start, crossing over] & & X'' \ar[dd, "{d_3}"] \\
      U' \ar[rr, "{c_2}" near start] \ar[dr, "{f'}" below] & & U_2 \ar[dr, "{f_2}"] & \\
      & X' \ar[<-, uu, "{p_1}" near end, crossing over] \ar[rr, "{d_2}"] & & X_2
    \end{tikzcd}
  \]
  and \(\ul{g} = (f_1, g, f_3)\), \(\ul{f'} = (f_1, f', f_2)\) and \(\ul{f''} = (f_2, f'', f_3)\).
  Then the right square in the commutative diagram
  \[
    \begin{tikzcd}
      U_1 \ar[d, "{f_1}"] & \ar[l, "{c_1 q_1}" above] U' \times_{U_2} U'' \ar[d, "{g}"] \ar[r, "{c_4 q_2}" above] & U_3 \ar[d, "{f_3}"] \\
      X_1 & \ar[l, "{d_1 p_1}" above] X' \times_{X_2} X'' \ar[r, "{d_4 p_2}" above] & X_3
    \end{tikzcd}
  \]
  is Cartesian and for any correspondences \(u: d_1^* L_1 \to d_2^! L_2\) and \(v: d_3^* L_2 \to d_4^! L_3\) we have an equality of correspondences supported on \((c_1 q_1, c_4 q_2)\):
  \[ \corr \ul{g}^* (v \circ u) = (\corr \ul{f''}^* v) \circ (\corr \ul{f'}^* u). \]
\end{prop}
\begin{proof}
  It is easy to deduce that the top square in the cube above is Cartesian from the fact that the left marked square is Cartesian in the first diagram, and by composition we deduce that the right square in the last diagram is also Cartesian.

  We only sketch the rather tedious comparison of correspondences.
  Unraveling definitions, we obtain long compositions for both sides involving units and counits for \((p_{2!}, p_2^!)\) (only for the left-hand side, with both unit and counit appearing), \((d_{2!}, d_2^!)\), \((c_{4!}, c_4^!)\), \((q_{2!}, q_2^!)\) \((d_{4!}, d_4^!)\), \((c_{2!}, c_2^!)\) (only for the right-hand side, with both unit and counit appearing), and base change isomorphisms:
  \begin{itemize}
  \item for the left-hand side, for the squares
    \[
      \begin{tikzcd}
        X' \times_{X_2} X'' \ar[d, "{p_1}"] \ar[r, "{p_2}"] & X'' \ar[d, "{d_3}"] & U' \times_{U_2} U'' \ar[d, "{g}"] \ar[r, "{q_2}"] & U'' \ar[d, "{f''}"] \ar[r, "{c_4}"] & U_3 \ar[d, "{f_3}"] \\
        X' \ar[r, "{d_2}"] & X_2 & X' \times_{X_2} X'' \ar[r, "{p_2}"] & X'' \ar[r, "{d_4}"] & X_3
      \end{tikzcd}
    \]
  \item for the right-hand side, for the three squares in
    \[
      \begin{tikzcd}
        U' \times_{U_2} U'' \ar[d, "{q_1}"] \ar[r, "{q_2}"] & U'' \ar[d, "{c_3}"] & U'' \ar[d, "{f''}"] \ar[r, "{c_4}"] & U_3 \ar[d, "{f_3}"] \\
        U' \ar[d, "{f'}"] \ar[r, "{c_2}"] & U_2 \ar[d, "{f_2}"] & X'' \ar[r, "{d_4}"] & X_3 \\
        X' \ar[r, "{d_2}"] & X_2 & &
      \end{tikzcd}
    \]
  \end{itemize}
  where pullback is for vertical maps and exceptional direct image for horizontal maps, and all base change isomorphisms are directed from top to bottom (i.e.\ left to right).
  One can reorder both compositions to remove the redundant unit/counit pairs \((p_{2!}, p_2^!)\) and \((c_{2!}, c_2^!)\), and use compatibility of base change isomorphisms with composition (both horizontally and vertically) to express on both sides the composition of all base change isomorphisms as the composition of two base change isomorphisms for the same squares.
\end{proof}

\subsection{Compactifications and canonical extensions}

The first part of the following lemma restates \cite[Lemma
1.3.1]{Fujiwara_Deligne_conj} and extends it to intermediate extensions of
perverse sheaves.

\begin{lemm} \label{lemm:uniq_int_ext_corr}
  Suppose that we have a commutative diagram
  \[ \begin{tikzcd}
      & U \arrow[dl, "{c_1}" above] \arrow[dr, "{c_2}"] \arrow[dd, "{j}"] & \\
      U_1 \arrow[dd, "{j_1}"] & & U_2 \arrow[dd, "{j_2}"]& \\
      & X \arrow[dl, "{\bar{c}_1}" above] \arrow[dr, "{\bar{c}_2}"] & \\
      X_1 & & X_2
    \end{tikzcd} \]
  of schemes separated and of finite type over \(B\), where \(j, j_1, j_2\) are open immersions.
  Let \(Z := X \smallsetminus j(U)\) and \(Z_k = X_k \smallsetminus j_k(U_k)\) for \(k=1,2\).
  Assume that for any \(k\) we have \(Z = \bar{c}_k^{-1}(Z_k)\), i.e.\ that both squares in the above diagram are cartesian (the inclusion \(\bar{c}_k^{-1}(Z_k) \subset Z\) holds automatically).
  \begin{enumerate}
  \item For any \(L_1 \in D^b_c(U_1)\) and \(L_2 \in D^b_c(U_2)\), the restriction morphisms
    \begin{align*}
      \Hom(\bar{c}_1^* j_{1 !} L_1, \bar{c}_2^! j_{2 !} L_2) & \longrightarrow \Hom(c_1^* L_1, c_2^! L_2) \\
      \Hom(\bar{c}_1^* j_{1 *} L_1, \bar{c}_2^! j_{2 *} L_2) & \longrightarrow \Hom(c_1^* L_1, c_2^! L_2)
    \end{align*}
    are isomorphisms.
  \item Assume moreover that \(B\) is the spectrum of a field and that the morphisms \(b_k : Z \rightarrow Z_k\) are quasi-finite.
    Then for any perverse sheaves \(L_1, L_2\) on \(U_1, U_2\), the restriction morphism
    \[ \Hom(\bar{c}_1^* j_{1 !*} L_1, \bar{c}_2^! j_{2 !*} L_2) \rightarrow \Hom(c_1^* L_1, c_2^! L_2) \]
    is an isomorphism.
  \end{enumerate}
\end{lemm}
\begin{proof}
  We only prove the second case of intermediate extensions, the other case (proved in \cite[Lemma 1.3.1]{Fujiwara_Deligne_conj} in a slightly different generality) being similar but using \(i_1^* j_{1 !} L_1 = 0\) and \(i_2^! j_{2 *} L_2 = 0\).
  Denote \(i : Z \hookrightarrow X\) and \(i_k : Z_k \hookrightarrow X_k\).
  We have \(\bar{c}_k \circ i = i_k \circ b_k\) where \(b_k : Z \rightarrow Z_k\) is quasi-finite by assumption.
  In particular \(b_1^*\) (resp.\ \(b_2^!\)) is left (resp.\ right) t-exact \cite[Proposition 2.2.5]{BBD}.
  Using the induction formula \cite[Exposé XVIII Corollaire 3.1.12.2]{SGA4-3}\footnote{For an immersion \(i\) the induction formula can be proved directly using the formalism of \cite[Exposé IV \S 14]{SGA4-1}.} we get
  \begin{align*}
    i^! \RuHom ( \bar{c}_1^* j_{1 !*} L_1, \bar{c}_2^! j_{2 !*} L_2) & \simeq \RuHom( i^* \bar{c}_1^* j_{1 !*} L_1, i^! \bar{c}_2^! j_{2 !*} L_2) \\
    & \simeq \RuHom( b_1^* i_1^* j_{1 !*} L_1, b_2^! i_2^! \bar{c}_2^! j_{2 !*} L_2).
  \end{align*}
  Moreover \(i_1^* j_{1 !*} L_1 \in {}^p D^{\leq -1}\) and \(i_2^! j_{2 !*} L_2 \in {}^p D^{\geq 1}\) \cite[Proposition 1.4.14]{BBD}, and by \cite[Proposition 2.1.20]{BBD}\footnote{The proof is given for topological spaces there, but it applies without modification to the categories \(D^b_{S,L}(-)\) for pairs \((S,L)\) as in \S 2.1.10 loc.\ cit.}
  \[ A := \RuHom( b_1^* i_1^* j_{1 !*} L_1, b_2^! i_2^! \bar{c}_2^! j_{2 !*} L_2) \]
  is an object of \(D^{\geq 2}\).
  Applying \(i_* i^! \rightarrow \mathrm{id} \rightarrow j_* j^*\) to \(\RuHom ( \bar{c}_1^* j_{1 !*} L_1, \bar{c}_2^! j_{2 !*} L_2)\) we obtain a distinguished triangle
  \[ i_* A \rightarrow \RuHom ( \bar{c}_1^* j_{1 !*} L_1, \bar{c}_2^! j_{2 !*} L_2) \rightarrow j_* \RuHom( c_1^* L_1, c_2^! L_2) \xrightarrow{+1}. \]
  Taking global sections and cohomology in degree \(0\), we find the desired isomorphism since \(H^0(R \Gamma A) = H^1(R \Gamma A) = 0\).
\end{proof}

In the situation of (a) in this lemma, for a correspondence \(u : c_1^* L_1
\rightarrow c_2^! L_2\), the induced correspondence \(\bar{c}_1^* j_{1!} L_1
\rightarrow \bar{c}_2^! j_{2!} L_2\) can be defined directly as follows:
\[ \bar{c}_1^* j_{1!} L_1 \xrightarrow[\sim]{\mathrm{BC}} j_! c_1^* L_1
\xrightarrow{j_! * u} j_! c_2^! L_2 \xrightarrow{\mathrm{BC}} \bar{c}_2^! j_{2!}
L_2. \]
and dually for \(\bar{c}_1^* j_{1*} L_1 \rightarrow \bar{c}_2^! j_{2*} L_2\).

\subsection{Nearby cycles}
\label{sec:corr_nearby_cycles}

For the rest of this section we consider the case where the base \(B = \{s, \eta\}\) is a Henselian trait, with special point \(s\) and generic point \(\eta\).
We use the notation of \cite[Exposé XIII]{SGA7-2}, except that we continue to
suppress the letter R from our notation (all objects considered are in derived
categories, unless explicitly mentioned otherwise).
For simplicity we fix a ``separable closure'' \(\etabar \to \eta\) and take, for
\(X\) a scheme over \(s\), Construction 1.2.4 loc.\ cit.\ as the definition of \(X
\times_s \eta\).
We will denote by \(F_X: X_{\bar{s}} \to X \times_s \eta\) ``the'' morphism of
toposes for which \(F_X^*\) is ``forgetting the action of \(\Gal(\etabar/\eta)\)''.
We will say that a sheaf in abelian groups \(\Fcal\) of \(X \times_s \eta\) is
constructible if \(F_X^* \Fcal\) is constructible.
The morphism \(F_X\) can be thought as an analogue of the base change morphism of
toposes \(\BC_X: X_{\sbar} \to X\), and the two are related by the specialization
morphism \(\mathrm{sp}: X \times_s \eta \to X\) since we have \(F_X^* \mathrm{sp}^*
= \BC_X^*\).
By \cite[Théorèmes de finitude \S 3]{SGA4etdemi} if \(X\) is a scheme of finite
type over \(B\) and \(\Fcal\) is a constructible sheaf in abelian groups over
\(X_\eta\) then \(\Psi_{\eta} \Fcal\) is constructible.

\subsubsection{Direct image}
\label{sec:app_nearby_direct}

Recall from \cite[Exposé XIII (2.1.6.2)]{SGA7-2} that any morphism \(f: X \to Y\)
of schemes over \(s\) induces a morphism of toposes \(X \times_s \eta \to Y
\times_s \eta\) and thus a derived direct image functor \(f_*: D^+(X \times_s
\eta, \Ocal_E/\mfrak_E^N) \to D^+(Y \times_s \eta, \Ocal_E/\mfrak_E^N)\).
Recall that the formation of \(f_*\) is even a (contravariant) normalized
pseudo-functor (see \cite[Exposé VI \S 8]{SGA1}), i.e.\ we can impose \((\id_X)_*
= \id\) for any \(X\), the obvious isomorphisms of functors \((f g)_* \simeq f_*
g_*\) satisfy a cocycle condition and are obvious when \(f\) or \(g\) is an identity
morphism.
 we have
the ``base change'' morphism of functors \(F_Y^* f_* \to f_{\bar{s}, *} F_X^*\)
(before derivation, this is simply an equality by definition).
If \(f\) is qcqs then this morphism of functors is an isomorphism
\footnote{Deligne claims loc.\ cit.\ that quasi-compact is enough for this to
hold, but we cannot think of an argument that does not require \(f\) to be
quasi-separated.
This can be proved by a familiar argument (analogous to \cite[Exposé VII
Théorème 5.7]{SGA4-2}, for the projective system \((X_{s'} \times_{s'}
\eta')_{\eta'}\), where \(\eta'\) ranges over all the finite étale covers of \(\eta\)
covered by \(\etabar\) and \(s'\) is the corresponding finite étale cover of \(s\)).
}.
Using the analogous property for \(\BC_X^*\) instead of \(F_X^*\) and the equality
\(F_X^* \mathrm{sp}^* = \BC_X^*\) we see that if \(f\) is qcqs then the morphism of
functors \(\mathrm{sp}^* f_* \to f_* \mathrm{sp}^*\) obtained by deriving the
obvious isomorphism of functors between categories of sheaves, is also an
isomorphism.

Recall ((2.1.7.1) loc.\ cit.) that if \(f: X \to Y\) is a morphism of schemes over
\(B\) then we have a morphism of functors \(\Psi_\eta f_{\eta *} \to f_{s *}
\Psi_\eta\) which is an isomorphism if \(f\) is proper.

\begin{lemm} \label{lem:nearby_direct_compo}
The morphisms of functors \(\Psi_\eta f_{\eta *} \to f_{s *}
\Psi_\eta\) are compatible with composition, in the sense
that for qcqs morphisms \(f: X \to Y\) and \(g: Y \to Z\) between schemes over \(B\)
the following diagram of functors is commutative.
\begin{equation} \label{diag:nearby_direct_compo}
  \begin{tikzcd}
    \Psi_\eta g_{\eta *} f_{\eta *} \arrow[r] & g_{s
    *} \Psi_\eta f_{\eta *} \arrow[r] & g_{s *} f_{s *} \Psi_\eta \\
    \Psi_\eta (gf)_{\eta *} \arrow[u, "{\sim}"] \arrow[rr] & & (gf)_{s *}
    \Psi_\eta \arrow[u, "{\sim}"]
  \end{tikzcd}
\end{equation}
\end{lemm}
\begin{proof}
  In this proof we momentarily denote by \(f_{s *}\) etc.\ the \emph{underived}
  functors, and denote by \(R?\) the right derived functor of \(?\).

  Before derivation the commutativity of \eqref{diag:nearby_direct_compo} may be
  checked after forgetting the action of \(\Gal(\etabar/\eta)\).
  This amounts to the compatibility of base change maps (for \(\ol{i}^*\)) with
  composition (see \cite[Exposé XII Proposition 4.4]{SGA4-2}; this can also be
  proved by considering the category fibered and cofibered over the category of
  schemes over \(\ol{B}\), with fiber over a scheme \(X\) the opposite category to
  that of étale sheaves of \(\Ocal_E/\mfrak_E^N\)-modules over \(X\), using the
  characterization of base change maps given by \cite[Exposé XVII Proposition
  2.1.3]{SGA4-2}).

  One can then derive the diagram, and obtain the following commutative diagram
  \[ \begin{tikzcd}[column sep=tiny]
      R\Psi_\eta  Rg_{\eta *}  Rf_{\eta *} & & Rg_{s *} R\Psi_\eta  Rf_{\eta *}
        & & \\
      R(\Psi_\eta g_{\eta *})  Rf_{\eta *} \arrow[u, "{\sim}"] \arrow[r] & 
        R(g_{s *} \Psi_\eta)  Rf_{\eta *} \arrow[ru] & Rg_{s *} R(\Psi_\eta
        f_{\eta *}) \arrow[u, "{\sim}"] \arrow[r] & Rg_{s *} R(f_{s *}
        \Psi_\eta) \arrow[r] & Rg_{s *}  Rf_{s *} R\Psi_\eta \\
      & R(\Psi_\eta g_{\eta *} f_{\eta *}) \arrow[lu] \arrow[r] \arrow[d, equal]
        & R(g_{s *} \Psi_\eta f_{\eta *}) \arrow[r] \arrow[lu] \arrow[u] &
        R(g_{s *} f_{s *} \Psi_\eta) \arrow[r] \arrow[u] \arrow[d, equal] &
        R(g_{s *} f_{s *})  R\Psi_\eta \arrow[u, "{\sim}"] \arrow[d, equal] \\
      R\Psi_\eta R(gf)_{\eta *} & \arrow[l, "{\sim}"] R(\Psi_\eta (gf)_{\eta *})
        \arrow[rr] & & R((gf)_{s *} \Psi_\eta) \arrow[r] & R(gf)_{s *}
        R\Psi_\eta
  \end{tikzcd} \]
  using the general fact that if \(F,G,H\) are composable additive functors
  between abelian categories such that \(F\), \(G\), \(H\), \(FG\), \(GH\) and \(FGH\) are
  everywhere right derivable, the two composite morphisms of functors \(R(FGH)
  \to R(FG) \circ RH \to RF \circ RG \circ RH\) and \(R(FGH) \to R(F) \circ R(GH)
  \to RF \circ RG \circ RH\) coincide.
  The original diagram \eqref{diag:nearby_direct_compo} is easily extracted.
\end{proof}

We see no reason for the forgetful functors \(F_X^*\) between derived categories
to be faithful, which is why the proof begins with the underived case.
All of the above admit parallel statements for the toposes \(X \times_s B\) and
the morphisms of functors \(\Psi f_* \to f_{s *} \Psi\), that we leave to the
reader.

\subsubsection{Inverse image}

Similarly, for \(f: X \to Y\) a morphism of schemes over \(s\) we have \(f^*: D^+(Y
\times_s \eta, \Ocal_E/\mfrak_E^N) \to D^+(X \times_s \eta,
\Ocal_E/\mfrak_E^N)\), which is ``the'' pseudo-functor (as \(f\) varies) left
adjoint to \(f_*\).
More explicitly, it is also obtained by deriving the obvious exact functor on
\(\Ocal_E/\mfrak_E^N\)-modules.
In particular we have isomorphisms of functors \(F_X^* f^* \simeq f_{\sbar}^*
F_Y^*\) and \(f^* \mathrm{sp}^* \simeq \mathrm{sp}^* f^*\).
For \(f: X \to Y\) a morphism of schemes over \(B\) one can define a morphism of
functors \(f_s^* \Psi_\eta \to \Psi_\eta f_\eta^*\), either directly (before
derivation and ignoring the action of \(\Gal(\etabar/\eta)\), this is a base
change map) or by adjunction, as the composition
\[ f_s^* \Psi_\eta \to f_s^* \Psi_\eta f_{\eta *} f_\eta^* \to \Psi_\eta
f_\eta^* f_{\eta *} f_\eta^* \to \Psi_\eta f_\eta^*. \]
Again these morphisms of functors are compatible with composition as \(f\) varies
(a diagram similar to \eqref{diag:nearby_direct_compo} is commutative), and
there are analogous constructions and statements for the toposes \(X \times_s B\)
and \(\Psi\).

\subsubsection{Exceptional direct image}

Recall from \cite[Exposé XIII 2.1.6 c)]{SGA7-2} that, for a separated morphism
\(f: X \to Y\) between schemes of finite type over \(s\), we have \(f_!: D^b_c(X
\times_s \eta, \Qell) \to D^b_c(Y \times_s \eta, \Qell)\) defined by
compactifying \(f\) as in \cite[Exposé XVII]{SGA4-3}.
We have an isomorphism of functors \(F_Y^* f_! \simeq f_{\bar{s}, !} F_X^*\): the
case of open immersions is trivial, the case of proper morphisms was considered
in \ref{sec:app_nearby_direct} and the general case follows.
Similarly we have an isomorphism of functors \(\mathrm{sp}^* f_! \simeq f_!
\mathrm{sp}^*\).
Again the formation of \(f_!\) is a pseudo-functor (by construction: it is
obtained by ``glueing'' the pseudo-functors \(f \mapsto f_*\) for \(f\) proper and
\(f \mapsto f_!\) for \(f\) an open immersion using the analogue of \cite[Exposé
XVII \S 5.1.5]{SGA4-3} for the toposes \(X \times_s \eta\) as ``glueing datum'').
Also by compactification (over \(B\)) and using the proper base change theorem,
for \(f: X \to Y\) a separated morphism between schemes of finite type over \(B\),
we have a morphism of functors \(f_{s !} \Psi_\eta \to \Psi_\eta f_{\eta !}\).

\begin{lemm} \label{lem:nearby_exc_direct_compo}
  For \(f: X \to Y\) and \(g: Y \to Z\) separated morphisms between schemes of
  finite type over \(B\) the following diagram of functors is commutative.
  \begin{equation} \label{diag:nearby_exc_direct_compo}
    \begin{tikzcd}
      \Psi_\eta g_{\eta !} f_{\eta !} \arrow[d, equals, "{\sim}"] &
      g_{s !} \Psi_\eta f_{\eta !} \arrow[l] & g_{s !} f_{s !} \Psi_\eta
      \arrow[d, equals, "{\sim}"] \arrow[l] \\
      \Psi_\eta (gf)_{\eta !} & & (gf)_{s !} \Psi_\eta \arrow[ll] 
    \end{tikzcd}
  \end{equation}
\end{lemm}
\begin{proof}
  \begin{enumerate}
    \item If \(f\) is an open immersion and \(g\) is proper then this is the very
      definition of \((gf)_{s !} \Psi_\eta \to \Psi_\eta (gf)_{\eta !}\).

    \item If \(f\) and \(g\) are proper all horizontal morphisms in the diagram
      \eqref{diag:nearby_exc_direct_compo} are isomorphisms and are the inverses
      of the horizontal morphisms in the diagram
      \eqref{diag:nearby_direct_compo}.

    \item In the case where \(f\) and \(g\) are open immersions we can argue as in
      the proof of Lemma \ref{lem:nearby_direct_compo} (again we momentarily
      consider functors before derivation, and explicitly denote the right
      derived functors by \(R?\)).
      \begin{enumerate}
        \item First we check the commutativity of the analogue of
          \eqref{diag:nearby_direct_compo} before derivation, and for this we
          can forget the action of \(\Gal(\etabar/\eta)\).
          We have a commutative diagram of functors
          \[ \begin{tikzcd}
              \ol{i}^* \ol{g}_! \ol{f}_! \arrow[d, equals] \arrow[r, "{\sim}"] &
              g_{\sbar !} \ol{i}^* \ol{f}_! \arrow[r, "{\sim}"] & g_{\sbar !}
              f_{\sbar !} \ol{i}^* \arrow[d, equals] \\
              \ol{i}^* \ol{gf}_! \arrow[rr, "{\sim}"] & & (gf)_{\sbar !}
              \ol{i}^* \end{tikzcd} \]
          thanks to the commutativity of the analogous diagram with \(!\) replaced
          by \(*\) (already used in the proof of Lemma
          \ref{lem:nearby_direct_compo}) and the characterization of the
          isomorphism of functors \(\ol{i}^* \ol{f}_! \to f_{\sbar !} \ol{i}^*\)
          in \cite[Exposé XVII Lemme 5.1.2]{SGA4-3}.
          We also have a commutative diagram of functors
          \[ \begin{tikzcd}
            \ol{j}_* g_{\etabar !} f_{\etabar !} \arrow[d, equals] & \ol{g}_!
            \ol{j}_* f_{\etabar !} \arrow[l] & \ol{g}_! \ol{f}_! \ol{j}_*
            \arrow[d, equals] \arrow[l] \\
            \ol{j}_* (gf)_{\etabar !} & & \ol{gf}_! \ol{j}_* \arrow[ll] 
          \end{tikzcd} \]
          This can be checked directly on the definition (details left to the
          reader).
          We thus obtain the commutativity of the underived analogue of
          \eqref{diag:nearby_exc_direct_compo}.

        \item As in the proof of Lemma \ref{lem:nearby_direct_compo}, deriving
          gives us a commutative diagram containing
          \eqref{diag:nearby_exc_direct_compo}, using the fact that the morphism
          of functors \(R(f_{s !} \circ \Psi_\eta) \to Rf_{s !} \circ R\Psi_\eta\)
          and the analogues for \(g\) and \(gf\) are isomorphisms (\(f_{s !}\) is
          exact).
      \end{enumerate}

    \item To conclude the proof it remains to consider the case where \(f\) is
      proper and \(g\) is an open immersion.
      Factoring \(gf\) as \(ba\) where \(b: T \to Z\) is proper and \(a: X \to T\) is an
      open immersion, we need to show that the following diagram of functors is
      commutative.
      \begin{equation} \label{diag:nearby_exc_direct_compo_swap}
        \begin{tikzcd}
          g_{s !} \Psi_\eta f_{\eta *} \arrow[r] \arrow[d, "{\sim}"] & \Psi_\eta
          g_{\eta !} f_{\eta *} \arrow[r, "{\sim}"] & \Psi_\eta b_{\eta *}
          a_{\eta !} \arrow[d, "{\sim}"] \\
          g_{s !} f_{s *} \Psi_\eta \arrow[r, "{\sim}"] & b_{s *} a_{s !}
          \Psi_\eta \arrow[r] & b_{s *} \Psi_\eta a_{\eta !}
      \end{tikzcd} \end{equation}
      where the top right and bottom left horizontal isomorphisms are induced by
      the isomorphisms \(g_{\eta !} f_{\eta *} \simeq b_{\eta *} a_{\eta !}\) and
      \(g_{s !} f_{s *} \simeq b_{s *} a_{s !}\) defined in \cite[Exposé XVII \S
      5.1.5]{SGA4-3}.
      Since the morphism \(X \to Y \times_Z T\) is a proper open immersion, we can
      easily reduce the problem to the case where the square
      \[ \begin{tikzcd}
        X \arrow[r, "{a}"] \arrow[d, "{f}"] & T \arrow[d, "{b}"] \\
        Y \arrow[r, "{g}"] & Z
      \end{tikzcd} \]
      is Cartesian.
      As usual to check the commutativity of the underived analogue of
      \eqref{diag:nearby_exc_direct_compo_swap} one can ignore the action of
      \(\Gal(\etabar/\eta)\).
      Details are left to the reader (ingredients are formal properties of
      adjunction, the triviality of base change by open immersions and the
      compatibility of base change maps with composition).
      Deriving and contemplating the resulting commutative diagram, one can
      extract the commutative diagram \eqref{diag:nearby_exc_direct_compo_swap}
      (this is very similar to the second part of the proof of Lemma
      \ref{lem:nearby_direct_compo}).
  \end{enumerate}
\end{proof}

Note that the fourth case with \(f=\id\) is already used to show that the
morphism of functors \(h_{s !} \Psi_\eta \to \Psi_\eta h_{\eta !}\) is
well-defined, i.e.\ does not depend on the choice of a compactification of \(h\).
Again, analogous statements hold for the toposes \(X \times_s B\) and with respect
to the functors \(\Psi\).

\subsubsection{Exceptional inverse image}

According to \cite[p.\ 45]{Illusie_monodromie} following the same method as in
\cite[Exposé XVIII \S 3.1]{SGA4-3} one obtains, for \(f: X \to Y\) a separated
morphism between schemes of finite type over \(s\), ``the'' functor \(f^!: D^+(Y
\times_s \eta, \Ocal_E/\mfrak_E^N) \to D^+(X \times_s \eta, \Ocal_E/\mfrak_E^N)\)
right adjoint to \(f_!\), which maps \(D^+_c(Y \times_s \eta,
\Ocal_E/\mfrak_E^N)\) to \(D^+_c(X \times_s \eta, \Ocal_E/\mfrak_E^N)\).
We have a morphism of functors \(F_X^* f^! \to f_{\sbar}^! F_Y^*\) defined by
adjunction as the composition
\begin{equation} \label{eq:def_forget_upper_shriek}
  F_X^* f^! \to f_{\sbar}^! f_{\sbar !} F_X^* f^! \simeq f_{\sbar}^! F_Y^* f_!
  f^! \to f_{\sbar}^! F_Y^*.
\end{equation}
To check that this is an isomorphism, we will adapt the site-theoretic arguments
of \cite[Exposé XVIII \S 3]{SGA4-3}.
For this we will need a ``nice'' site \(\Ccal_{X,B}\) whose category of sheaves
can be identified with \(X \times_s \eta\).
We model the definition on \cite[Exposé XVIII 3.1.15]{SGA4-3}: consider the
category \(\Ccal_{X,B}\) of pairs \((\eta', U)\) where \(\eta'\) is finite étale over
\(\eta\), with residual scheme \(s' \to s\), and \(U\) is a scheme étale over
\(X_{s'}\).
A morphism \((\eta_1, U_1) \to (\eta_2, U_2)\) is given by a morphism \(\eta_1 \to
\eta_2\) over \(\eta\) and a morphism \(U_1 \to U_2\) compatible with \(X_{s_1} \to
X_{s_2}\).
In particular the morphism \(U_1 \to U_2\) is étale.
One easily checks that \((\eta,X)\) is a final object in this category and that
fiber products exist in this category, thus so do finite projective limits.
Define a family of morphisms \(((\eta_i, U_i) \to (\eta', U))_i\) to be a covering
if the family \((U_i \to U)_i\) is jointly surjective.
It is not difficult to check that this defines a pretopology, and so a site with
underlying category \(\Ccal_{X,B}\).
From \cite[Exposé VII 2.a)]{SGA4-2} one can easily deduce that the site
\(\Ccal_{X,B}\) is subcanonical.

\begin{itemize}
  \item
    If \(\Fcal\) is an object of \(X \times_s \eta\), then
    \[ \Fcal': (\eta', U) \longmapsto \varprojlim_{\alpha: \etabar \to \eta'}
    \Fcal(U \times_{s'} \sbar), \]
    where the projective limit is over the groupoid of morphisms compatible with
    the morphisms to \(\eta\), defines a sheaf on \(\mathcal{C}_{X,B}\).
    Note that if \(\eta'\) is connected (i.e.\ a point) and if we choose
    \(\alpha_0: \etabar \to \eta'\) then the right-hand side is identified with
    \(\Fcal(U \times_{s'} \sbar)^{\Gal(\etabar/\eta')}\), and a general \((\eta',
    U)\) can be written as a finite disjoint union of pairs \((\eta'', U')\) with
    \(\eta''\) a point.
  \item
    Conversely, if \(\Fcal'\) is a sheaf on \(\Ccal_{X,B}\) then we can define an
    object \(\Fcal\) of \(X \times_s \eta\) by
    \[ \Fcal(\ol{U}) = \varinjlim_{(\alpha: \etabar \to \eta', U, \varphi)}
    \Fcal'(\eta', U) \]
    for \(\ol{U}\) a scheme étale over \(X_{\sbar}\) with \(\ol{U}\) affine (or qcqs),
    where the injective limit is over the category opposite to the category
    \(I_{\ol{U}}\) of triples \((\alpha: \etabar \to \eta', U, \varphi)\) where
    \(\eta'\) is finite étale over \(\eta\), \(U\) is a scheme étale over \(X_{s'}\) and
    \(\varphi: \ol{U} \simeq U \times_{s'} \sbar\) (the morphism \(\sbar \to s'\)
    used here being the reduction of \(\alpha\)).
    Note that thanks to the qcqs assumption on \(\ol{U}\) the index category
    \(I_{\ol{U}}^\mathrm{opp}\) is filtered (essentially by \cite[Théorème
    8.8.2]{EGA4-3} and \cite[Proposition 17.7.8]{EGA4-4}), and if \((\alpha_0:
    \etabar \to \eta_0', U_0, \varphi_0)\) is any object then the subcategory of
    \((\alpha: \etabar \to \eta', U_0 \times_{s'_0} s', \varphi_0)\) given by some
    \(\etabar \to \eta' \to \eta_0'\), is cofinal.
    For \(\sigma \in \Gal(\etabar/\eta)\) with reduction \(\sigma_\mathrm{red}\) the
    action \(a_\sigma: \Fcal(\sigma_\mathrm{red}^* \ol{U}) \simeq \Fcal(\ol{U})\)
    is defined by the equivalence of index categories \(I_{\ol{U}} \to
    I_{\sigma_\mathrm{red}^* \ol{U}}\), \((\alpha, U, \varphi) \mapsto (\alpha
    \sigma, U, \varphi \times \sigma_\mathrm{red})\).
\end{itemize}

\begin{prop}
  The first (resp.\ second) construction gives a sheaf \(\Fcal'\) on \(\Ccal_{X,B}\)
  (resp.\ an object of \(X \times_s \eta\)), and the two constructions are inverse
  of each other.
\end{prop}
\begin{proof}
  Let \(\Fcal\) be an object of \(X \times_s \eta\), and let \(((\eta_i, U_i) \to
  (\eta',U))\) be a covering family of morphisms in \(\Ccal_{X,B}\).
  Let \((f_i)_i\) be a family of sections \(f_i \in \Fcal'(\eta_i,U_i)\).
  We may see each \(f_i\) as a family \((f_{i,\alpha_i})_{\alpha_i}\) indexed by all
  morphisms \(\alpha_i: \etabar \to \eta_i\), where \(f_{i,\alpha_i} \in
  \Fcal(U_{i,\alpha_i})\) with \(U_{i,\alpha_i} = U_i \times_{s_i} \sbar\) for the
  morphism \(\sbar \to s_i\) obtained by reduction from \(\alpha_i\).
  The action of \(\Gal(\etabar/\eta)\) on \(\Fcal\) induces isomorphisms, for each
  \(\sigma \in \Gal(\etabar/\eta)\), \(a_{\sigma, \alpha_i}: \Fcal(U_{i, \alpha_i
  \sigma}) \to \Fcal(U_{i, \alpha_i})\) satisfying the cocycle condition \(a_{\tau
  \sigma, \alpha_i} = a_{\tau, \alpha_i} a_{\sigma, \alpha_i \tau}\).
  By definition we have \(a_{\sigma, \alpha_i}(f_{i, \alpha_i \sigma}) = f_{i,
  \alpha_i}\) for all \(\sigma \in \Gal(\etabar/\eta)\) and all \(\alpha_i: \etabar
  \to \eta_i\).
  Now assume that for all indices \(i\) and \(j\) the images of \(f_i\) and \(f_j\) in
  \(\Fcal'(\eta_{i,j}, (U_i \times_U U_j) \times_{(s_i \times_{s'} s_j)}
  s_{i,j})\) coincide, where we have denoted \(\eta_{i,j} = \eta_i \times_\eta
  \eta_j\), with residue scheme \(s_{i,j}\).
  This means that for any \(\alpha: \etabar \to \eta'\) and for any \(\alpha_i\) and
  \(\alpha_j\) above \(\alpha\), the images of \(f_{i,\alpha_i}\) and \(f_{j,\alpha_j}\)
  in \(\Fcal(U_{i,\alpha_i} \times_{U'_\alpha} U_{j,\alpha_j})\) coincide.
  For any \(\alpha: \etabar \to \eta'\) the family \((U_{i,\alpha_i} \to
  U'_\alpha)_{i,\alpha_i}\) is an étale covering, and so the family
  \((f_{i,\alpha_i})_{i,\alpha_i}\) comes from a unique \(f_\alpha \in
  \Fcal(U'_\alpha)\).
  The relation \(a_{\sigma, \alpha}(f_{\alpha \sigma}) = f_\alpha\) for \(\sigma
  \in \Gal(\etabar/\eta)\) is clear.

  Conversely let \(\Fcal'\) be a sheaf on \(\Ccal_{X,B}\).
  The fact that the associated functor \(\Fcal\) is a sheaf on the étale site of
  \(X_{\sbar}\) (restricted to qcqs objects) follows from standard arguments, and
  the fact that the Galois action is continuous is obvious.

  We leave it to the reader to check that the two functors are inverse of each
  other.
\end{proof}

If \((\eta', U)\) is an object of \(\Ccal_{X,B}\) such that \(\eta'\) is a point, with
corresponding trait \(B' \to B\), we have an identification of the sites
\(\Ccal_{X,B}/(\eta',U)\) and \(\Ccal_{U,B'}\) and thus an identification of the
toposes \(\widetilde{\Ccal_{X,B}/(\eta',U)}\) and \(U \times_{s'} \eta'\).
In general \(\eta'\) is a finite disjoint union of points and we define
\(\Ccal_{U,B'}\) and \(U \times_{s'} \eta'\) as products over these points in an
obvious way.

We now describe certain points of the site \(\Ccal_{U,B'}\), similarly to the
usual case \cite[Exposé VIII Proposition 3.9]{SGA4-2}.
Let \(\alpha: \etabar \to \eta'\) be a morphism over \(\eta\) \footnote{one could
more generally consider morphisms from the spectrum of an arbitrary separably
closed (or algebraically closed) field, but this would complicate the
compatibility constraint given by \(\gamma\)}, inducing a morphism
\(\alpha_\mathrm{red}: \sbar \to s'\) over \(s\).
Let \(\beta: \Spec(k) \to U\) be a geometric point, inducing a geometric point of
\(s'\), and let \(\gamma: \Spec(k) \to \sbar\) be a morphism over \(s'\).
We will call such triples \((\alpha, \beta, \gamma)\) geometric points of
\((\eta',U)\).
Consider the category \(\Ccal_{U,B'}^{\alpha, \beta, \gamma}\) with objects
\footnote{and obvious morphisms \dots} \((\eta'', V, \tilde{\alpha},
\tilde{\beta}, \tilde{\gamma})\) where \((\eta'', V)\) is an object of
\(\Ccal_{U,B'}\) (the morphism to \((\eta',U)\) being implied) and \((\tilde{\alpha},
\tilde{\beta}, \tilde{\gamma})\) is a geometric point of \((\eta'',V)\) lifting
\((\alpha,\beta,\gamma)\).
Using a variation of the proof loc.\ cit.\ one easily checks that this category
is cofiltered.
The functor \(\widetilde{\Ccal_{U,B'}} \to \mathrm{Sets}\),
\begin{equation} \label{eq:fiber_func}
  \Fcal' \longmapsto \Fcal'_{\alpha, \beta, \gamma} :=
  \varinjlim_{(\Ccal_{U,B'}^{\alpha, \beta, \gamma})^\mathrm{opp}} \Fcal'
\end{equation}
is a fiber functor.
In fact \(\alpha\) induces a morphism of toposes \(F_{U,\alpha}: U_{\sbar} \to U
\times_{s'} \eta'\) (here \(U_{\sbar} = U \times_{s'} \sbar\) uses
\(\alpha_\mathrm{red}\)), and \eqref{eq:fiber_func} is identified with the fiber
functor of the étale site of \(U_{\sbar}\) obtained by composing \(F_{U,\alpha}^*\)
with the fiber functor corresponding to the geometric point \(\beta \times
\gamma: \Spec(k) \to U_{\sbar}\).
In particular if \((\alpha_i,\beta_i,\gamma_i)_i\) is a family of points such that
the \(\alpha_i\)'s are jointly surjective then it is a conservative family of
points.

If \(f: X \to Y\) is a morphism between schemes over \(s\) then the associated
morphism of toposes \(X \times_s \eta \to Y \times_s \eta\) can also be obtained
from the morphism of sites \(\Ccal_{X,B} \to \Ccal_{Y,B}\) given by the functor
\(\Ccal_{Y,B} \to \Ccal_{X,B}\), \((\eta',V) \mapsto (\eta', V \times_Y X)\) (this
is almost identical to \cite[Exposé XVII \S 1.4]{SGA4-2}).
Let \(j: (\eta',V) \to (\eta,Y)\) be an object of \(\Ccal_{Y,B}\), then we can form
the ``Cartesian diagram''
\footnote{
  In general it only is a Cartesian diagram in a category that we have not
  defined, so this ``diagram'' is only a visual aid.
  One could extend the (at least) 4 functors formalism to this category, but we
  will not need this generality.
  If \(f\) is étale then this is a Cartesian diagram in \(\Ccal_{Y,B}\).
} 
\begin{equation} \begin{tikzcd} \label{eq:fake_Cartesian}
  (\eta',V \times_Y X) \arrow[r, "{k}"] \arrow[d, "{g}"] & (\eta, X) \arrow[d,
  "{f}"] \\
  (\eta',V) \arrow[r, "{j}"] & (\eta, Y)
\end{tikzcd} \end{equation}
We will need various compatibilites for this ``diagram'', generalising the ones
already known when \(\eta'=\eta\).
We obviously have \(j_* g_* = f_* k_*\), and by transposition we have an
isomorphism of functors \(k^* f^* \simeq g^* j^*\), compatibly with composition of
\(f\) or \(j\).
We also have the usual ``tautological base change'' isomorphism \(j^* f_* \simeq
g_* k^*\), also compatible with compositions in two ways.
If \(f\) is an open immersion then the general \cite[Exposé XVII Lemme
5.1.2]{SGA4-3} gives a canonical isomorphism \(j^* f_! \simeq g_! k^*\).
If \(X\) and \(Y\) are separated and of finite type over \(s\) then by
compactifying \(f\) we obtain an isomorphism \(j^* f_! \simeq g_! k^*\), and by the
same arguments as in Lemmas 5.2.3, 5.2.4 and 5.2.5 in \cite[Exposé XVII]{SGA4-3}
this isomorphism does not depend on the compactification and is compatible with
composition (both horizontal and vertical).
Before giving more compatibilities, we recall and compare two constructions,
horizontal or vertical in \eqref{eq:fake_Cartesian}, of trace maps.

\begin{lemm} \label{lem:trace_map_zero}
  If \(f: X \to Y\) is a flat and quasi-finite morphism between separated schemes
  of finite type over \(s\), there is a unique morphism of functors \(\Tr_f: f_!
  f^* \to \id\) on \(\Ocal_E/\mfrak_E^N\)-modules in \(Y \times_s \eta\) which after
  applying \(F_Y^*\) is the usual trace map for \(f_{\sbar}\) \cite[Exposé XVII
  Théorème 6.2.3]{SGA4-3}.
\end{lemm}
\begin{proof}
  For once we are dealing with sheaves and not objects in derived categories,
  and so it is enough to observe that the trace map for \(f_{\sbar}\) is
  \(\Gal(\etabar/\eta)\)-equivariant, which follows from compatibility with base
  change.
\end{proof}
This trace map is also characterized by taking stalks at all geometric points in
\(\Ccal_{X,B}\) as above.

Consider a scheme \(X\) over \(s\) and \(j: (\eta_1, U_1) \to (\eta_2,U_2)\) a
morphism in \(\Ccal_{X,B}\).
For any \(\alpha_1: \etabar \to \eta_1\) over \(\eta\) we have an obvious
isomorphism \(F_{U_1,\alpha_1}^* j^* \simeq j_{\alpha_1}^* F_{U_2,\alpha_2}^*\)
where \(j_{\alpha_1}\) is the morphism from \(U_1 \times_{s_1} \sbar\) (using the
reduction \(\sbar \to s_1\) of \(\alpha_1\)) to \(U_2 \times_{s_2} \sbar\) (using the
reduction \(\sbar \to s_2\) of \(\alpha_2: \etabar \to \eta_2\) obtained by
composing \(\alpha_1\) and \(\eta_1 \to \eta_2\)).
This isomorphism is compatible with the identifications of stalks
\eqref{eq:fiber_func} of pullbacks.
Recall that the restriction functor \(j^*\) (on \(\Ocal_E/\mfrak_E^N\)-modules)
admits a left adjoint that we denote by \(j_!\) (see \cite[Exposé IV Proposition
11.3.1]{SGA4-1}).
The formation of \(j_!\), like \(j^*\), is a pseudo-functor on \(\Ccal_{X,B}\).
This notation does not conflict with the previous one: if \(\eta'\) is finite
étale over \(\eta\) and \(\eta_1 = \eta_2 = \eta'\) then the trace map of Lemma
\ref{lem:trace_map_zero} realizes \(j_!\) (as previously defined by
compactification, replacing \(B\) by the trait \(B'\) corresponding to \(\eta'\)) as
the left adjoint to \(j^*\) (see \cite[Exposé XVII Proposition 6.2.11]{SGA4-3}).
For \(\alpha_1: \etabar \to \eta_1\) over \(\eta\) consider the composition, where
\(\alpha_2: \etabar \to \eta_2\) is \(\alpha_1\) composed with \(\eta_1 \to \eta_2\),
\[ j_{\alpha_1 !} F_{U_1,\alpha_1}^* \to j_{\alpha_1 !} F_{U_1,\alpha_1}^* j^*
j_! \simeq j_{\alpha_1 !} j_{\alpha_1}^* F_{U_2,\alpha_2}^* j_! \to
F_{U_2,\alpha_2}^* j_!. \]
Summing all possible \(\alpha_1\)'s we obtain, for \(\alpha_2: \etabar \to \eta_2\)
over \(\eta\), a morphism of functors
\begin{equation} \label{eq:forget_site_shriek}
  \bigoplus_{\alpha_1} j_{\alpha_1 !} F_{U_1, \alpha_1}^* \longrightarrow
  F_{U_2, \alpha_2}^* j_!
\end{equation}
where the sum is over all morphisms \(\alpha_1: \etabar \to \eta_1\) whose
composition with \(\eta_1 \to \eta_2\) is \(\alpha_2\).
If \((\alpha_2,\beta_2,\gamma_2)\) is a geometric point of \((\eta_2, U_2)\) and
\(\Fcal\) is a sheaf of \(\Ocal_E/\mfrak_E^N\)-modules in \(\widetilde{\Ccal_{U_1,
B_1}}\) then as usual we have an identification
\begin{equation} \label{eq:stalk_site_shriek}
  (j_! \Fcal)_{\alpha_2,\beta_2,\gamma_2} \simeq \bigoplus_{(\alpha_1, \beta_1,
  \gamma_1)} \Fcal_{\alpha_1,\beta_1,\gamma_1}
\end{equation}
where the sum is over all geometric points of \((\eta_1,U_1)\) lifting \((\alpha_2,
\beta_2, \gamma_2)\).
Using this identification one can easily check that the morphism of functors
\eqref{eq:forget_site_shriek} is an isomorphism.
By construction it is compatible with the adjunctions of functors \((j_!, j^*)\)
and \((j_{\alpha_1 !}, j_{\alpha_1}^*)\).
The functor \(j_!\) (on \(\Ocal_E/\mfrak_E^N\)-modules) is exact, and we will still
denote by \(j_!\) the induced functor between derived categories.

We get back to the situation of \eqref{eq:fake_Cartesian}, given by a morphism
\(f: X \to Y\) of schemes over \(s\) and an object \(j: (\eta',V) \to (\eta,Y)\) in
\(\Ccal_{Y,B}\).
Consider the morphism of functors obtained as the following composition
\begin{equation} \label{eq:BC_site_shriek}
  k_! g^* \longrightarrow k_! g^* j^* j_! \simeq k_! k^* f^* j_! \to f^* j_!.
\end{equation}
Applying the forgetful functor \(F_X^*\) and using the identifications above (in
particular \eqref{eq:forget_site_shriek}) we see that this morphism of functors
is an isomorphism.
Under the assumption that \(X\) and \(Y\) be separated and of finite type over \(s\)
we can define a second morphism of functors (now between derived categories)
\(j_! g_! \to f_! k_!\) by abstract nonsense, as the composition
\begin{equation} \label{eq:compose_site_shriek}
  j_! g_! \to j_! g_! k^* k_! \simeq j_! j^* f_! k_! \to f_! k_!.
\end{equation}
Applying the forgetful functor \(F_Y^*\) and using the above identification
\eqref{eq:forget_site_shriek}, we see that this morphism of functors is also an
isomorphism.

For a scheme \(X\) over \(s\) and an integer \(n\) invertible on \(s\) we can pull back
the sheaf \(\mu_n\) via the morphism of toposes \(\mathrm{sp}: X \times_s \eta \to
X\) and then define Tate twists as usual by tensoring with tensor powers of
\(\mathrm{sp}^* \mu_n\).

\begin{lemm} \label{lem:trace_map_one}
  Let \(f: X \to Y\) be a morphism between schemes which are separated and of
  finite type over \(s\).
  Assume that \(f\) is a flat and that all of its fibers are purely of dimension
  \(1\).
  Then there exists, for any abelian group object \(\Fcal\) of \(Y \times_s \eta\)
  killed by an integer invertible on \(B\), a unique trace morphism \(\Tr_f : H^2(
  f_! f^* \Fcal)(1) \to \Fcal\) which after applying \(F_Y^*\) and using the
  identification \(F_Y^* f_! f^* \simeq f_{\sbar !} f_{\sbar}^* F_Y^*\) is the
  usual trace morphism defined in \cite[Exposé XVIII Proposition
  1.1.6]{SGA4-3}.
\end{lemm}
\begin{proof}
  We are considering sheaves rather than objects in the derived category and so
  it is enough to show that the trace morphism for \(f_{\sbar}\) is
  \(\Gal(\etabar/\eta)\)-equivariant, which follows from compatibility of trace
  maps with base change.
\end{proof}

\begin{lemm} \label{lem:compat_trace_maps}
  Let \(f: X \to Y\) be a flat morphism between schemes which are separated and of
  finite type over \(s\).
  Assume that all fibers of \(f\) are purely of dimension one.
  Let \(j: (\eta',V) \to (\eta,Y)\) be an object of \(\Ccal_{Y,B}\), and let \(g\) and
  \(k\) be as in \eqref{eq:fake_Cartesian}.
  Let \(\Fcal\) be an \(\Ocal_E/\mfrak_E^N\)-module in \(V \times_{s'} \eta'\).
  Then the following diagram commutes.
  \[ \begin{tikzcd}[column sep=large]
      H^2(f_! k_! g^* \Fcal) \arrow[r, "{\sim}", "{\eqref{eq:BC_site_shriek}}"']
      & H^2(f_! f^* j_! \Fcal) \arrow[d, "{\Tr_f * j_!}"] \\
      H^2(j_! g_! g^* \Fcal) \arrow[u, "{\sim}",
      "{\eqref{eq:compose_site_shriek}}"'] \arrow[r, "{j_! * \Tr_g}"] & j_!
      \Fcal(-1)
  \end{tikzcd} \]
\end{lemm}
\begin{proof}
  All morphisms in the diagram are compatible with the forgetful functor \(F_Y^*\)
  by definition, and so we are reduced to check that for any \(\alpha: \etabar
  \to \eta'\), the following diagram is commutative.
  \[ \begin{tikzcd}[column sep=large]
    H^2(f_{\sbar !} k_{\alpha !} g_\alpha^* F_{V,\alpha}^* \Fcal) \arrow[r,
    "{\sim}"] & H^2(f_{\sbar !} f_{\sbar}^* j_{\alpha !}
    F_{V,\alpha}^* \Fcal) \arrow[d, "{\Tr_{f_{\sbar}} * j_{\alpha !}}"] \\
    H^2(j_{\alpha !} g_{\alpha !} g_\alpha^* F_{V,\alpha}^* \Fcal) \arrow[u,
    "{\sim}"] \arrow[r, "{j_{\alpha !} * \Tr_{g_\alpha}}"] & j_{\alpha !}
    F_{V,\alpha}^* \Fcal(-1)
  \end{tikzcd} \]
  This can be checked on stalks at any geometric point of \(Y_{\sbar}\), using
  compatibility of trace maps with base change \footnote{the fact that similar
  diagrams commute seems to be used implicitly in the proof of \cite[Exposé
  XVIII Théorème 3.2.5]{SGA4-3}}.
\end{proof}

Let \(f\) be a flat curve as in the previous lemma.
Fix a compactification of \(f\) and a conservative set of points of the topos \(X
\times_s \eta\) (for example one for each geometric point of \(X\)).
As in \cite[Exposé XVIII \S 3.1.15]{SGA4-3} we have a site \(\Gamma_B(f)\)
equivalent to \(\Ccal_{X,B}\) and adapted to \(f^*\).
Namely, the underlying category of \(\Gamma_B(f)\) is that of quadruples \((\eta',
U, V, \varphi)\) where \((\eta',V)\) is an object of \(\Ccal_{Y,B}\), \((\eta',U)\) is
an object of \(\Ccal_{X,B}\) and \(\varphi: U \to V \times_Y X\) is a morphism of
schemes over \(X\) (automatically étale); morphisms are defined in the obvious
manner, and a family of morphisms \((\eta_i, U_i, V_i, \varphi_i) \to (\eta', U,
V \varphi)\) is covering if it induces a covering of \(U\).
For \((\eta',U,V,\varphi)\) an object of \(\Gamma_B(f)\) denote (somewhat abusively)
\[ (\eta', U) \xrightarrow{\varphi} (\eta', V \times_Y X) \xrightarrow{k}
(\eta,X) \]
in \(\Ccal_{X,B}\) and define \(K(\eta',U,V,\varphi)\) as the complex concentrated
in degrees \(0,1,2\) analogous to \cite[Exposé XVIII (3.1.4.7)]{SGA4-3} computing
\(f_! k_! \varphi_! \ul{\Ocal_E/\mfrak_E^N}\).
Denoting \(j\) the morphism \((\eta',V) \to (\eta,Y)\) in \(\Ccal_{Y,B}\), define
\(K''(\eta',U,V,\varphi)\) as \(j_! \ul{\Ocal_E/\mfrak_E^N}(-1)\) concentrated in
degree \(2\).
The composition of the trace map for \(\varphi\), the identification
\eqref{eq:BC_site_shriek} \(k_! g^* \simeq f^* j_!\) and the trace map for \(f\)
gives a morphism of complexes of \(\Ocal_E/\mfrak_E^N\)-modules in \(Y \times_s
\eta\):
\[ K(\eta', U, V, \varphi) \longrightarrow K''(\eta', U, V, \varphi) \]
analogous to \cite[Exposé XVIII (3.2.1.1)]{SGA4-3}.
We obtain a morphism of functors \(K \to K''\) from \(\Gamma_B(f)\) to the category
of complexes of \(\Ocal_E/\mfrak_E^N\)-modules in \(Y \times_s \eta\).
We have seen in Lemma \ref{lem:compat_trace_maps} that this morphism can be
described using the trace map for \(g: V \times_Y X \to V\).
Using this description we obtain the following consequence of the ``fundamental
lemma'' \cite[Exposé XVIII Lemme 1.6.9]{SGA4-3}.

\begin{lemm}
  Let \(f: X \to Y\) be a smooth relative curve between separated schemes of
  finite type over \(s\).
  Let \((\eta', U, V, \varphi)\) be an object of \(\Gamma_B(f)\).
  Let \((\alpha, \beta, \gamma)\) be a geometric point of \(\Ccal_{U,B'}\).
  There exist an object \((\eta', U', V', \varphi')\) of \(\Gamma_B(f)\) over
  \((\eta', U, V, \varphi)\) and a geometric point \((\alpha', \beta', \gamma')\) of
  \((\eta', U', V', \varphi')\) lifting \((\alpha, \beta, \gamma)\) such that
  \begin{enumerate}
    \item in cohomology in degrees \(0\) and \(1\), the morphism
      \[ K(\eta',U',V',\varphi') \longrightarrow K(\eta',U,V,\varphi) \]
      induces zero maps.
    \item in cohomology in degree \(2\), the morphism \(K(\eta', U', V', \varphi')
      \to K''(\eta', U', V', \varphi')\) induces an isomorphism.
  \end{enumerate}
\end{lemm}
\begin{proof}
  This follows from Lemme 1.6.9 loc.\ cit.\ using Lemma
  \ref{lem:compat_trace_maps} twice (``base change'' the trace map for \(f: X \to
  Y\) by \((\eta',V) \to (\eta,Y)\) and the trace map for \(U \to V\) (with \(B\)
  replaced by \(B'\)) by \((\eta',V') \to (\eta',V)\) (this second case is easier
  than Lemma \ref{lem:trace_map_zero})
\end{proof}

Exactly as in \cite[\S 3.2.1]{SGA4-3}, for \(f\) a flat relative curve which is a
morphism between separated schemes of finite type over \(s\) the morphism of
functors \(K \to K''\) induces a morphism of functors \(f^*(1)[2] \to f^!\) (on
derived categories) which, by the same proof as in Lemma 3.2.3 loc.\ cit.,
coincides with the composition
\begin{equation} \label{eq:duality_sm_curve}
  f^*(1)[2] \longrightarrow f^! f_! f^*(1)[2] \xrightarrow{\Tr_f} f^!
\end{equation}
where \(\Tr_f: f_! f^*(1)[2] \to \id\) is deduced from Lemma
\ref{lem:trace_map_one} as in (2.13.2) loc.\ cit.

\begin{prop}
  Let \(f: X \to Y\) be a smooth relative curve between separated schemes of
  finite type over \(s\).
  Then the morphism of functors \(f^*(1)[2] \to f^!\) defined by
  \eqref{eq:duality_Hecke} is an isomorphism.
\end{prop}
\begin{proof}
  Using the previous lemma, the proof is identical to that of Théorème 3.2.5
  loc.\ cit.
\end{proof}

It seems likely that the analogue of Théorème 3.2.5 loc.\ cit.\ (i.e.\ the
generalization of this proposition to smooth morphisms of arbitrary relative
dimension) could be proved using the same strategy, but fortunately we will not
need this statement.

\begin{coro} \label{cor:compat_upper_shriek_curves}
  Under the same assumption, the morphism of functors
  \eqref{eq:def_forget_upper_shriek} is an isomorphism.
\end{coro}
\begin{proof}
  This follows from the proposition and the fact that the following diagram
  commutes.
  \[ \begin{tikzcd}
      F_X^* f^* (1)[2] \arrow[d, "{\sim}"] \arrow[r] & F_X^* f^! f_! f^* (1)[2]
      \arrow[r, "{\Tr_f}"] & F_X^* f^! \arrow[r] & f_{\sbar}^! f_{\sbar !} F_X^*
      f^! \arrow[d, "{\sim}"] \\
      f_{\sbar}^* (1)[2] F_Y^* \arrow[r] & f_{\sbar}^! f_{\sbar !} f_{\sbar}^*
      (1)[2] F_Y^* \arrow[r, "{\Tr_{f_{\sbar}}}"] & f_{\sbar}^! F_Y^* &
      \arrow[l] f_{\sbar}^! F_Y^* f_! f^!
  \end{tikzcd} \]
  This commutativity follows from the compatibility of \(\Tr_f\) and
  \(\Tr_{f_{\sbar}}\) using the adjunction formalism.
\end{proof}

\begin{prop} \label{pro:compat_upper_shriek}
  For \(f: X \to Y\) a morphism between schemes separated and of finite type over
  \(s\) the morphism of functors \eqref{eq:def_forget_upper_shriek} \(F_X^* f^! \to
  f_{\sbar}^! F_Y^*\) is an isomorphism.
\end{prop}
\begin{proof}
  The morphism of functors \eqref{eq:def_forget_upper_shriek} is compatible with
  composition.
  We are thus reduced to proving that this morphism is an isomorphism
  Zariski-locally (or even étale-locally) on \(X\), since for \(g\) an open
  immersion into \(X\) we have an isomorphism \(g^! \simeq g^*\) given by the trace
  map (Lemma \ref{lem:trace_map_zero}) and similarly for \(g_{\sbar}\), and we
  already know the compatibility of \(g^*\) with \(F_?^*\).
  Thus we can assume that \(f\) is a morphism between \emph{affine} schemes of
  finite type over \(s\).
  We can factor \(f\) as \(\pi i\) where \(\pi: \A_Y^d \to Y\) is the typical affine
  space of relative dimension \(d\) and \(i\) is a closed immersion.
  We are reduced to proving the statement with \(f\) replaced by \(i\) or \(\pi\).

  So let us assume first that \(f: X \to Y\) is a closed immersion between affine
  schemes of finite type over \(s\).
  Denote by \(g: U \to Y\) the complementary open immersion.
  The right adjoint functor \(f^!\) also arises as the right derived functor of
  the ``sections with support in \(X\)'' functor \({}^0 f^!\), and similarly for
  \(f_{\sbar}\) and we leave it to the reader to check that the morphism \(F_X^*
  f^! \to f_{\sbar}^! F_Y^*\) which is defined by adjunction coincides with the
  one obtained by deriving the obvious isomorphism of functors \(F_X^* {}^0 f^!
  \simeq {}^0 f_{\sbar}^! F_Y^*\)
  \footnote{Note that, as in the case of the derived direct image functor, there
    is no a priori reason for \(F_Y^*\) to map an injective sheaf of
    \(\Ocal_E/\mfrak_E^N\)-modules to one that is acyclic for \({}^0 f_{\sbar}^!\),
    and a priori so the morphism of functors between derived categories may not
  be an isomorphism.}.
  Thus we have a commutative diagram whose rows are (functors in) distinguished
  triangles.
  \[ \begin{tikzcd}
      F_X^* f_* f^! \arrow[r] \arrow[d] & F_X^* \arrow[d, equal] \arrow[r] &
      F_X^* g_* g^* \arrow[d] \arrow[r, "{+1}"] & {} \\
      f_{\sbar *} f_{\sbar}^! F_X^* \arrow[r] & F_X^* \arrow[r] & g_{\sbar *}
      g_{\sbar}^* F_X^* \arrow[r, "{+1}"] & {}
  \end{tikzcd} \]
  The middle and right vertical morphisms are isomorphisms, so the left one is
  also an isomorphism, and applying \(f_{\sbar}^*\) gives us the lemma in this
  case.

  Finally, assume that \(Y\) is affine and that \(f: X \to Y\) is a typical affine
  space of relative dimension \(d\).
  By composition we are easily reduced to the case where \(d=1\), which is covered
  by Corollary \ref{cor:compat_upper_shriek_curves}.
\end{proof}

\begin{coro}
  For \(f: X \to Y\) a morphism between schemes separated and of finite type over
  \(s\) the composite morphism of functors
  \begin{equation} \label{eq:compat_sp_upper_shriek}
    \mathrm{sp}^* f^! \to f^! f_! \mathrm{sp}^* f^! \simeq f^! \mathrm{sp}^* f_!
    f^! \to f^! \mathrm{sp}^*
  \end{equation}
  is an isomorphism.
\end{coro}
\begin{proof}
  The same arguments used to prove Proposition \ref{pro:compat_upper_shriek}
  show that the composite morphism of functors
  \[ \BC_X^* f^! \to f_{\sbar}^! f_{\sbar !} \BC_X^* f^! \simeq f_{\sbar}^!
    \BC_Y^* f_! f^! \to f_{\sbar}^! \BC_Y^* \]
  is also an isomorphism, and using the fact that the isomorphisms of functors
  \(\BC_?^* \simeq F_?^* \mathrm{sp}^*\), \(\mathrm{sp}^* f_! \simeq f_!
  \mathrm{sp}^*\), \(\BC_Y^* f_! \simeq f_{\sbar !} \BC_X^*\) and \(F_Y^* f_! \simeq
  f_{\sbar !} F_X^*\) are compatible one can deduce formally that
  \eqref{eq:compat_sp_upper_shriek} is also an isomorphism.
\end{proof}

For \(f: X \to Y\) a separated morphism between schemes of finite type over \(B\),
from the morphism of functors \(f_{s !} \Psi_\eta \to \Psi_\eta f_{\eta !}\)
we obtain by adjunction a morphism of functors \(\Psi_\eta f_\eta^! \to f_s^!
\Psi_\eta\).
The analogue of Lemma \ref{lem:nearby_exc_direct_compo} holds true for formal
reasons.

\subsubsection{Correspondences and nearby cycles}
\label{sec:corr_nearby}

Let \(c_i: X \to X_i\) for \(i=1,2\) be a pair of separated morphisms between
schemes of finite type over \(B\), and \(L_i\) an object of \(D^b_c(X_{i, \eta},
\Ocal_E/\mfrak_E^N)\) (resp.\ \(D^b_c(X_i, \Ocal_E/\mfrak_E^N)\)).
Denote \(\pi_i: X_i \to B\).
If \(u\) is a correspondence from \(L_1\) to \(L_2\) with support in \((c_{1, \eta},
c_{2, \eta})\) (resp.\ \((c_1, c_2)\)) then we get a correspondence \(\Psi_{\eta} u\)
(resp.\ \(\Psi u\)) from \(\Psi_{\eta} L_1\) to \(\Psi_{\eta} L_2\) (resp.\ \(\Psi L_1\)
to \(\Psi L_2\)) with support in the pair of morphisms between 2-fibred toposes
\((c_{1,s}, c_{2,s})\) (resp.\ \((c_1, c_2)\)).
In the case of \(\Psi_{\eta}\) it is defined as the composition
\begin{equation} \label{eq:corr_nearby}
  c_{1, s}^* \Psi_\eta L_1 \rightarrow \Psi_{\eta} c_{1, \eta}^* L_1
  \xrightarrow{\Psi_{\eta} * u} \Psi_{\eta} c_{2, \eta}^! L_2 \rightarrow c_{2,
  s}^! \Psi_{\eta} L_2.
\end{equation}
In the case of \(\Psi u\), note that \((\Psi u)_s : c_{1,s}^* i^* L_1 \rightarrow
c_{2,s}^! i^* L_2\) is simply \(i^* u\), defined using \(i^* c_2^! \rightarrow
c_{2,s}^! i^*\), and that we get a morphism of correspondences \(\mathrm{sp}^* i^*
(L_1, L_2, u) \to \Psi_\eta j^* (L_1, L_2, u)\).

\begin{prop} \label{pro:compat_coh_corr_nearby}
  If \(u\) is a correspondence from \(L_1 \in \Ob D^+_c(X_{1,\eta},
  \Ocal_E/\mfrak_E^N)\) to \(L_2 \in \Ob D^+_c(X_{2,\eta}, \Ocal_E/\mfrak_E^N)\)
  and if \(c_1\) (resp.\ \(c_2\)) is proper then the diagram
  \[ \begin{tikzcd}
      \pi_{1, s !} \Psi_{\eta} L_1 \arrow[r, "{(\Psi_{\eta} u)_!}"] \arrow[d] &
        \pi_{2, s !} \Psi_{\eta} L_2 \arrow[d] \\
      \Psi_\eta \pi_{1, \eta !} L_1 \arrow[r, "{u_!}"] & \Psi_\eta \pi_{2, \eta
      !} L_2
  \end{tikzcd}
  \text{ resp.\ }
  \begin{tikzcd}
    \pi_{1, \eta *} L_1 \arrow[r, "{u_*}"] \arrow[d] & \pi_{2, \eta *} L_2
      \arrow[d] \\
    \pi_{1, s *} \Psi_{\eta} L_1 \arrow[r, "{(\Psi_{\eta} u)_*}"] & \pi_{2, s *}
      \Psi_{\eta} L_2
  \end{tikzcd} \]
  in \(D^+_c(\eta, \Ocal_E/\mfrak_E^N)\) is commutative.
\end{prop}
\begin{proof}
  We only prove the commutativity of the first diagram, as the second case is
  very similar.
  For the proof it will be slightly simpler to view \(u\) as a morphism \(c_{2,
  \eta !} c_{1, \eta}^* L_1 \to L_2\) (by adjunction).
  It is a formal exercise in adjunction to check that \(\Psi_\eta u : c_{2, s !}
  c_{1, s}^* \Psi_\eta L_1 \to \Psi_\eta L_2\) is the composite
  \[ c_{2, s !} c_{1, s}^* \Psi_\eta L_1 \to c_{2, s !} \Psi_\eta c_{1, \eta}^*
  L_1 \to \Psi_\eta c_{2, \eta !} c_{1, \eta}^* L_1 \xrightarrow{\Psi_\eta * u}
  \Psi_\eta L_2 \]
  and that \(u_!: \pi_{1, \eta !} L_1 \to \pi_{2, \eta !} L_2\) is the composite
  \[ \pi_{1, \eta !} L_1 \to \pi_{1, \eta !} c_{1, \eta *} c_{1, \eta}^* L_1
  \simeq \pi_{2, \eta !} c_{2, \eta !} c_{1, \eta}^* L_1 \xrightarrow{\pi_{2,
  \eta !} * u} \pi_{2, \eta !} L_2 \]
  and similarly for \((\Psi_\eta u)_!\).
  We want to show that the diagram
  \[ \begin{tikzcd}[column sep=tiny]
    \pi_{1, s!} \Psi_\eta L_1 \arrow[r] \arrow[d] & \pi_{1, s!} c_{1, s *}
    c_{1, s}^* \Psi_\eta L_1 \arrow[r, equal, "{\sim}"] & \pi_{2, s!} c_{2, s !}
    c_{1, s}^* \Psi_\eta L_1 \arrow[r] & \pi_{2, s!} c_{2, s !}
    \Psi_\eta c_{1, \eta}^* L_1 \arrow[r] & \pi_{2, s!} \Psi_\eta c_{2, \eta !}
    c_{1, \eta}^* L_1 \arrow[d] \\
    \Psi_\eta \pi_{1, \eta !} L_1 \arrow[r] & \Psi_\eta \pi_{1, \eta !} c_{1,
    \eta *} c_{1, \eta}^* L_1 \arrow[r, equal, "{\sim}"] & \Psi_\eta \pi_{2,
    \eta !} c_{2, \eta !} c_{1, \eta}^* L_1 \arrow[r] & \Psi_\eta \pi_{2, \eta
    !} L_2 & \arrow[l] \pi_{2, s!} \Psi_\eta L_2
  \end{tikzcd} \]
  commutes.
  Reordering the morphisms of functors used along both paths as we may, we find
  that the diagram commutes thanks to the commutativity of
  \[ \begin{tikzcd}
    \Psi_\eta \arrow[r] \arrow[d] & c_{1, s *} c_{1, s}^* \Psi_\eta \arrow[d] \\
    \Psi_\eta c_{1, \eta *} c_{1, \eta}^* \arrow[r] & c_{1, s*} \Psi_\eta c_{1,
    \eta}^* \\
  \end{tikzcd} \]
  which follows formally from the definition of \(c_{1, s}^* \Psi_\eta \to
  \Psi_\eta c_{1, \eta}^*\) by adjunction, and Lemma
  \ref{lem:nearby_exc_direct_compo} applied with \((g, f) = (\pi_1, c_1)\) and
  \((\pi_2, c_2)\).
\end{proof}

\subsection{Nearby cycles of perverse sheaves}
\label{sec:corr_nearby_cycles_perverse}

We will need the notion of perverse sheaves for the topos \(X \times_s \eta\)
where \(X\) is scheme of finite type over \(s\).
The properties of \S 1.4.3 in \cite{BBD} are satisfied for any open subtopos of
an arbitrary topos, with an arbitrary sheaf of rings.
Note that the open subtoposes of \(X \times_s \eta\) are precisely the \(U \times_s
\eta\) for \(U\) an open subscheme of \(X\) (to check this one is immediately reduced
to finite étale descent of open immersions).
The consequences in \S 1.4.2.1 loc.\ cit.\  are thus satisfied as well, and they
also hold for \(\Ocal_E\)-sheaves and \(E\)-sheaves.
Thus the results of \S 1.4 loc.\ cit.\ can be applied without modification.
The definitions of \S 2.2.9 and 2.2.10 loc.\ cit.\ only need to be slightly
modified: we consider stratifications of \(X\) satisfying the same smoothness
condition over \(\sbar\), and for a stratum \(S\) a finite set of isomorphism
classes of local systems (with coefficients \(\Ocal_E/\mfrak_E\), \(\Ocal_E\) or
\(E\)) over \(S_{\sbar}\).
The argument at the end of \S 2.2.10 loc.\ cit.\ still applies since for \(j: U
\hookrightarrow X\) we have \(F_X^* j_* \simeq j_{\sbar *} F_U^*\).
Assuming as usual that the perversity function and its dual are non-increasing,
the argument using cohomological purity in \S 2.2.11 loc.\ cit.\ also applies to
show that by refining the stratification and enlarging the finite collections of
isomorphism classes of local systems on strata considered, we get compatible
t-structures.
Note that this uses Proposition \ref{pro:compat_upper_shriek}.
This gives us a t-structure on \(D_{\ctf}(X, \Ocal_E/\mfrak_E^N)\) (or \(D^b_c(X,
\Ocal_E)\) or \(D^b_c(X, E)\)) associated to a perversity function.
From now on we will only consider the case of \(E\)-coefficients and the self-dual
perversity function.
To sum up, we have a t-structure on \(D^b_c(X \times_s \eta, E)\) which is
compatible with the usual one on \(D^b_c(X_{\sbar}, E)\): a complex \(K\) in
\(D^c_c(X \times_s \eta, E)\) is in \({}^p D^{\geq 0}\) if and only if \(F_X^* K\) is
in \({}^p D^{\geq 0}\).
In particular the functor \(\mathrm{sp}^*: D^b_c(X, E) \to D^b_c(X \times_s \eta,
E)\) is also compatible with the perverse t-structures.

These compatibilities rely on the compatibility of \(F_?^*\) with the four
operations associated to a morphism between separated schemes of finite type
over \(s\).
We will also need the following compatibility with (derived, as usual)
``internal Hom''.

\begin{lemm} \label{lem:compat_RuHom_forget}
  For a complex \(K\) in \(D_{\ctf}(X \times_s \eta, \Ocal_E/\mfrak_E^N)\) (resp.\
  \(D^b_c(X \times_s \eta, \Ocal_E)\), resp.\ \(D^b_c(X \times_s \eta, E)\)) and a
  complex \(L\) in \(D(X \times_s \eta, \Ocal_E/\mfrak_E^N)\) (resp.\ \(D(X \times_s
  \eta, \Ocal_E)\), resp.\ \(D(X \times_s \eta, E)\)), the morphism \(F_X^*
  \RuHom(K, L) \to \RuHom(F_X^* K, F_X^*L)\) obtained by adjunction from
  \cite[Exposé IV (13.4.2)]{SGA4-2} is an isomorphism.
\end{lemm}
\begin{proof}
  This follows from \cite[Exposé VI Corollaire 8.7.9]{SGA4-2}.
\end{proof}

\begin{lemm} \label{lemm:uniq_int_ext_corr2}
  The conclusion of Lemma \ref{lemm:uniq_int_ext_corr} 2 holds true if we consider schemes over the base \(s\) and perverse sheaves \(L_k\) in \(D^b_c(U_k \times_s \eta, E)\).
\end{lemm}
\begin{proof}
  The ingredients in the proof of Lemma \ref{lemm:uniq_int_ext_corr} are still
  valid over toposes \(? \times_s \eta\):
  \begin{enumerate}
    \item The induction formula
      \[ \RuHom(i^* K, i^!L) \simeq i^! \RuHom(K, L) \]
      holds for any inclusion \(i: F \hookrightarrow X\) of a closed subtopos \(F\)
      (complementary to an open subtopos \(U\)) of a topos \(X\).
    \item The inclusion (\cite[Proposition 2.1.20]{BBD} in the usual setting)
      \[ \RuHom({}^p D_c^{\geq a}, {}^p D_c^{\leq b}) \subset D^{\leq b-a} \]
      holds true thanks to the previous lemma
      \footnote{this could certainly also be proved for an arbitrary topos
      endowed with a stratification.}.
  \end{enumerate}
\end{proof}

\begin{lemm} \label{lem:nearby_int_ext_corr}
  Suppose that we have a commutative diagram of schemes separated and of finite
  type over \(B\)
  \[ \begin{tikzcd}
    & U \arrow[dl, "{c_1}" above] \arrow[dr, "{c_2}"] \arrow[dd, "{j}"] & \\
    U_1 \arrow[dd, "{j_1}"] & & U_2 \arrow[dd, "{j_2}"]& \\
    & X \arrow[dl, "{\bar{c}_1}" above] \arrow[dr, "{\bar{c}_2}"] & \\
    X_1 & & X_2
  \end{tikzcd} \]
  where \(j, j_1, j_2\) are open immersions and both squares are cartesian.
  For \(k=1,2\) let \(L_k\) be an object of \(D^b_c(X_{k,\eta}, E)\).
  Then the diagram of groups of correspondences
  \[ \begin{tikzcd}
    \Hom(\bar{c}_{1, \eta}^* L_1, \bar{c}_{2, \eta}^! L_2) \arrow[r] \arrow[d] &
    \Hom(c_{1, \eta}^* j_{1,\eta}^* L_1, c_{2, \eta}^! j_{2,\eta}^* L_2)
    \arrow[d] \\
    \Hom(\bar{c}_{1,s}^* \Psi_\eta L_1, \bar{c}_{2, s}^!  \Psi_\eta L_2)
    \arrow[r] & \Hom(c_{1,s}^* \Psi_\eta j_{1,\eta}^* L_1, c_{2, s}^!  \Psi_\eta
    j_{2,\eta}^* L_2)
    \end{tikzcd} \]
  is commutative.
\end{lemm}
\begin{proof}
  This follows from Lemmas \ref{lem:nearby_direct_compo} and
  \ref{lem:nearby_exc_direct_compo} and their adjoints.
  Details are left to the reader.
\end{proof}

\section{Irreducible finite-dimensional \((\gfrak, K)\)-modules}
\label{sec:irr_fd_gK}

In this appendix we denote \(\Gamma = \Gal(\C/\R)\).
Let \(\Gbf\) be a (connected) reductive group over \(\R\).
Choose a maximal compact subgroup \(K\) of \(\Gbf(\R)\).
Let \(\gfrak\) be the complexification of \(\Lie \Gbf(\R)\), i.e.\ the Lie algebra
of \(\Gbf_\C\) in the algebraic sense.
Denote by \(X^*(\Gbf)^\Gamma\) the group of morphisms \(\Gbf \to \GLbf_{1,\R}\)
(defined over \(\R\)).
We have an injective morphism
\begin{align*}
  \C \otimes_\Z X^*(\Gbf)^\Gamma & \longrightarrow \Hom_\cont(\Gbf(\R),
  \C^\times) \\
  s \otimes \chi & \longmapsto (g \mapsto |\chi(g)|^s)
\end{align*}
and we will implicitly identify \(\C \otimes_\Z X^*(\Gbf)^\Gamma\) with a subgroup
of \(\Hom_\cont(\Gbf(\R), \C^\times)\).
Our next goal in Lemma \ref{lem:fd_rep_real_gp_alg} is to show that if
\(\Gbf_\der\) is simply connected then irreducible finite-dimensional \((\gfrak,
K)\)-modules are obtained by twisting an algebraic representation of \(\Gbf_\C\) by
a character in \(\C \otimes_\Z X^*(\Gbf)^\Gamma\), but first we consider the case
of tori.

\begin{lemm} \label{lem:twist_alg_char}
  Let \(\Tbf\) be a torus over \(\R\).
  \begin{enumerate}
    \item
      Any continuous character \(\chi: \Tbf(\R) \to \C^\times\) can be written as
      the product of an element of \(\C \otimes_\Z X^*(\Tbf)\) and the restriction
      of an algebraic character of \(\Tbf(\C)\).
    \item
      A character in \(\C \otimes_\Z X^*(\Tbf)^\Gamma\) is the restriction of an
      algebraic character \(\Tbf(\C) \to \C^\times\) if and only if it belongs to
      \((1+\sigma) X^*(\Tbf)\).
  \end{enumerate}
\end{lemm}
\begin{proof}
  \begin{enumerate}
    \item
      Assume first that \(\Tbf = \GLbf_{1,\R}\).
      Up to multiplying \(\chi\) by the restriction of the identity character we
      can assume that \(\chi(-1) = -1\), so that \(\chi\) belongs to \(\C \otimes_\Z
      X^*(\Tbf)^\Gamma\).

      Assume next that \(\Tbf\) is anisotropic of dimension one, i.e.\ \(\Tbf(\C)
      \simeq \C^\times\) and for \(z \in \Tbf(\C)\) we have \(\sigma(z) =
      \ol{z}^{-1}\).
      The character \(\chi\) of the circle \(\Tbf(\R)\) can be written \(z \mapsto
      z^a\) for some integer \(a\), and is clearly the restriction of an algebraic
      character of \(\Tbf(\C)\).

      Assume that \(\Tbf = \Res_{\C/\R} \GL_{1,\C}\), so that \(\Tbf(\C) \simeq
      \C^\times \times \C^\times\) with Galois action given by \(\sigma(z_1, z_2)
      = (\ol{z_2}, \ol{z_1})\).
      In particular \(\Tbf(\R) \simeq \{ (z,\ol{z}) \,|\, z \in \C^\times \}\).
      The character \(\chi\) of \(\Tbf(\R)\) can be written \((z,\ol{z}) \mapsto
      (z/|z|)^a |z|^b\) for \(a \in \Z\) and \(b \in \C\).
      We have \(\chi(z,\ol{z}) = z^a |z|^{b-a}\).
      The character \(\Tbf(\R) \to \R_{>0}\), \((z,\ol{z}) \mapsto |z|^{b-a}\)
      belongs to \(\C \otimes_\Z X^*(\Tbf)^\Gamma\) and the character \(T(\R) \to
      \C^\times\), \((z,\ol{z}) \mapsto z^a\) is the restriction of the algebraic
      character \(\Tbf(\C) \to \C^\times\), \((z_1,z_2) \mapsto z_1^a\).

      Finally if \(\Tbf\) is an arbitrary torus then \(\Tbf\) decomposes as a
      product of indecomposable tori which are isomorphic to one of the three
      tori considered above.

    \item
      As for the previous point this can be checked easily for each isomorphism
      class of indecomposable tori.
      Details are left to the reader.
  \end{enumerate}
\end{proof}

\begin{lemm} \label{lem:fd_rep_real_gp_alg}
  Assume that the derived subgroup \(\Gbf_\der\) of \(\Gbf\) is simply connected.
  Any irreducible finite-dimensional \((\gfrak, K)\)-module is, up to twisting by
  a character in \(\C \otimes_\Z X^*(\Gbf)^\Gamma\), the restriction of an
  irreducible algebraic representation of \(\Gbf(\C)\).
\end{lemm}
\begin{proof}
  Let \(V\) be an irreducible finite-dimensional \((\gfrak, K)\)-module.
  More precisely, we have morphisms \(\pi: \gfrak \to \End(V)\) and \(\rho: K \to
  \GL(V)\) satisfying the usual conditions.
  The semisimple Lie algebra \([\gfrak, \gfrak]\) acts via \(\pi\) on \(V\), and
  comparing the classifications of representations of \([\gfrak, \gfrak]\) and of
  \(\Gbf_\der\) (or using the fact that \(\Gbf_\der(\C)\) is simply connected and
  that any holomorphic representation of \(\Gbf_\der(\C)\) is algebraic) we see
  that it integrates to an algebraic representation of \(\Gbf_\der(\C)\), that we
  still denote \(\pi\).

  Let \(\Abf_\Gbf\) be the largest split central torus in \(\Gbf\), and denote \(A =
  \Abf_\Gbf(\R)^0\).
  Similarly via \(\pi\) and Schur's lemma we get a scalar action of the vector
  group \(\Gbf\) on \(V\), again denoted by \(\pi\).

  Next we check that we can glue the actions of \(A\), \(\Gbf_\der(\C)\) and \(K\) to
  get an action of the subgroup \(A \cdot \Gbf_\der(\C) \cdot K\) of \(\Gbf(\C)\).
  First we observe that \(A\) is central in \(\Gbf(\C)\) and \(\Gbf_\der(\C)\) is
  normal in \(\Gbf(\C)\).
  Any element of \(A \cdot \Gbf_\der(\C) \cdot K\) can be written \(a g k\) with \(a
  \in A\), \(g \in \Gbf_\der(\C)\) and \(k \in K\), and the triple \((a,g,k)\) is
  unique up to the equivalence relation \((a,g,k) \sim (a,gx,x^{-1}k)\) for all \(x
  \in \Gbf_\der(\C) \cap K\).
  Now \(\Gbf_\der(\C) \cap K\) is a maximal compact subgroup in \(\Gbf_\der(\R)\)
  which is connected because \(\Gbf_\der\) is simply connected \cite[Corollaire
  4.7]{BorelTits_compl}, and so \(\Gbf_\der(\C) \cap K\) is also connected, which
  implies that \(\pi\) (defined by integration on \(\Gbf_\der(\C)\)) and \(\rho\)
  agree on \(\Gbf_\der(\C) \cap K\).
  It follows that the map
  \begin{align*}
    A \cdot \Gbf_\der(\C) \cdot K & \longrightarrow \GL(V) \\
    a g k & \longmapsto \pi(a) \pi(g) \rho(k)
  \end{align*}
  is well defined.
  Using the relation \(\pi(\Ad(k)(X)) = \Ad(\rho(k))(\pi(X))\) for \(k \in K\) and
  \(X \in \gfrak\) it is easy to check that we get an action of \(A \cdot
  \Gbf_\der(\C) \cdot K\) on \(V\), both extending \(\rho\) and integrating
  \(\pi|_{[\gfrak,\gfrak] \oplus \Lie A}\).

  Letting \(\Dbf = \Gbf/\Gbf_\der\) we have a short exact sequence
  \[ 1 \to \Gbf_\der(\C) \to \Gbf(\C) \xrightarrow{p} \Dbf(\C) \to 1. \]
  The inclusion \(p(A \cdot K) \subset p(\Gbf(\R))\) is an equality: \(A \cdot K\)
  contains the center of \(\Gbf(\R)\) and so \(p(A \cdot K)\) is an open subgroup of
  \(\Dbf(\R)\), and \(K\) meets every connected component of \(\Gbf(\R)\).
  Therefore we have \(A \cdot \Gbf_\der(\C) \cdot K = \Gbf_\der(\C) \cdot
  \Gbf(\R)\).
  If \(\Tbf\) is a maximal torus of \(\Gbf\) then \(p(\Tbf(\R))\) is an open subgroup
  of \(\Dbf(\R)\).
  If moreover \(\Tbf\) contains a maximal split torus of \(\Gbf\) (i.e.\ if \(\Tbf\)
  is a maximal torus in the centralizer of a maximal split torus of \(\Gbf\)) then
  \(\Tbf(\R)\) meets every connected component of \(\Gbf(\R)\) \cite[Théorème
  14.4]{BorelTits_gpesred} and so we have \(p(\Tbf(\R)) = p(\Gbf(\R))\), i.e.\
  \(\Gbf_\der(\C) \cdot \Gbf(\R) = \Gbf_\der(\C) \Tbf(\R)\).
  Let \(\Bbf\) be a Borel subgroup of \(\Gbf_\C\) containing \(\Tbf_\C\), let \(\ufrak
  \subset \gfrak\) be the (algebraic) Lie algebra of its unipotent radical and
  let \(\ufrak^-\) be the Lie algebra of the unipotent radical of the opposite
  Borel subgroup of \(\Gbf_\C\) with respect to \(\Tbf_\C\), so that we have \(\gfrak
  = \ufrak^- \oplus \tfrak \oplus \ufrak\) where \(\tfrak = \Lie \Tbf_\C = \C
  \otimes_\R \Lie \Tbf(\R)\).
  Let \(V^\ufrak = \{ v \in V \,|\, \forall X \in \gfrak,\, \pi(X) v = 0 \}\) be
  the subspace of maximal vectors in \(V\), so that \(V\) is generated as a
  representation of \(\ufrak^-\) by \(V^\ufrak\), i.e.\ \(\pi(U(\ufrak^-)) V^\ufrak =
  V\) where \(U(-)\) denotes the universal enveloping algebra.
  Denoting \(\Tbf_\der = \Gbf_\der \cap \Tbf\), the subgroup \(\Tbf_\der(\C) \cdot
  \Tbf(\R)\) of \(\Tbf(\C)\) preserves \(V^\ufrak\) and the action of \(\Tbf_\der(\C)\)
  on \(V^\ufrak\) is algebraic.
  There exists a line \(L \subset V^\ufrak\) preserved by the commutative group
  \(\Tbf(\R)\), on which it acts by a continuous character \(\chi\).
  If \(X_1, \dots, X_n\) are eigenvectors for the roots \(\alpha_1, \dots,
  \alpha_n\) of \(\Tbf_\C\) acting on \(\ufrak^-\) then for \(t \in \Tbf(\R)\) and \(v
  \in L\) the action of \(t\) on \(\pi(X_1) \dots \pi(X_n) v\) is multiplication by
  \(\chi(t) \alpha_1(t) \dots \alpha_n(t)\).
  The subspace \(\pi(U(\ufrak^-)) L\) is an irreducible
  sub-\([\gfrak,\gfrak]\)-representation of \(V\), i.e.\ an irreducible
  sub-\(\Gbf_\der(\C)\)-representation.
  Therefore \(\pi(U(\ufrak^-)) L\) is an irreducible sub-\(\Gbf_\der(\C) \cdot
  \Tbf(\R)\)-representation of \(V\), and so it is equal to \(V\) (this also implies
  the equality \(V^\ufrak = L\)).
  By the first point of Lemma \ref{lem:twist_alg_char} there exists a character
  \(\nu \in \C \otimes_\Z X^*(\Tbf)^\Gamma\) of \(\Tbf(\R)\) such that \(\chi \nu\) is
  the restriction of an algebraic character of \(\Tbf(\C)\).
  Since the restriction of \(\chi\) to \(\Tbf_\der(\R)\) is already algebraic, by
  the second point of Lemma \ref{lem:twist_alg_char} there exists \(\lambda \in
  X^*(\Tbf_\der)\) such that the restriction of \(\nu\) to \(\Tbf_\der\) equals
  \((1+\sigma) \lambda\).
  We have a short exact sequence
  \[ 0 \to X^*(\Dbf) \to X^*(\Tbf) \to X^*(\Tbf_\der) \to 0. \]
  Let \(\tilde{\lambda} \in X^*(\Tbf)\) be any lift of \(\lambda\).
  Up to dividing \(\nu\) by \((1+\sigma) \tilde{\lambda}\), which preserves the
  algebraicity of \(\chi \nu\), we can assume that \(\nu\) belongs to
  \[ \C \otimes_Z X^*(\Dbf)^\Gamma = \C \otimes_\Z X^*(\Gbf)^\Gamma. \]
  So up to twisting \(V\) by an element of \(\C \otimes_\Z X^*(\Gbf)^\Gamma\) we can
  assume that \(\chi\) is the restriction of an algebraic character of \(\Tbf(\C)\).
  It is now clear that \(V\) extends uniquely to an (irreducible) algebraic
  representation of \(\Gbf(\C)\).
\end{proof}

\newpage
\section{Hecke formalism for boundary strata of minimal compactifications of Shimura varieties}
\label{sec:gen_Shim}

A technical nuisance when working with minimal compactifications of Shimura varieties is the fact that the boundary strata are not quite Shimura varieties, but only quotients by certain finite groups of Shimura varieties.
Even in cases where they are (disjoint union of) Shimura varieties, keeping track of Hecke operators quickly becomes a notational burden: the analogues of automorphic local systems on boundary strata, which appear thanks to Pink's theorem \cite{Pink_ladic_Shim}, are defined using group cohomology of arithmetic groups, so that we are led to consider a mix of group cohomology and étale cohomology.
We find convenient to isolate concerns as follows:
\begin{itemize}
\item introduce a slight generalization of the notion of Shimura datum, associate ``generalized Shimura varieties'' to them, define the analogue of automorphic local systems on them and check that their cohomology groups are naturally endowed with Hecke operators (commuting with the action of the Galois group in the \(\ell\)-adic setting),
\item show that the boundary components of such generalized Shimura varieties are themselves generalized Shimura varieties.
\end{itemize}

Most of our efforts will be spent proving that our generalized automorphic local systems induce a Hecke formalism in cohomology.

\subsection{Generalized Shimura varieties}

\begin{defi} \label{def:gen_Shim_datum}
  A generalized Shimura datum is a triple \((\Gbf, \Xcal, h)\) where \(\Gbf\) is a connected reductive group over \(\Q\), \(\Xcal\) is a homogeneous space under \(\Gbf(\R)\) and \(h: \Xcal \to \Hom(\mathbb{S}, \Gbf_{\R})\) is a finite-to-one \(\Gbf(\R)\)-equivariant map satisfying the following weakening of \cite[\S 2.1.1]{Deligne_ShimuraCorvallis}, for any (equivalently, one) \(x \in \Xcal\):
  \begin{enumerate}
  \item the adjoint action of \(h(x)\) on \(\Lie \Gbf_{\R}\) is of type \((-1,1)\), \((0,0)\), \((1,-1)\).
  \item for any simple (over \(\Q\)) factor \(\Hbf\) of \(\Gbf_{\ad}\), either \(h(x)\) acts trivially on \(\Lie \Hbf_{\R}\) or \(\Hbf_{\R}\) is isotropic and conjugation by \(h(x)(i)\) is a Cartan involution of \(\Hbf_{\R}\).
  \end{enumerate}
\end{defi}
The first condition implies that for any \(x \in \Xcal\) the image of \(h(x)\) by \(\Hom(\mathbb{S}, \Gbf_{\R}) \to \Hom(\mathbb{S}, \Gbf_{\ad, \R})\) is trivial on the split one-dimensional subtorus in \(\mathbb{S}\).
Let \(\Gbf_\her\) be the smallest algebraic (over \(\Q\)) subgroup of \(\Gbf\) satisfying \(h(\Xcal) \subset \Hom(\mathbb{S}, \Gbf_{\her, \R})\).
It is a normal subgroup of \(\Gbf\) and we have a canonical decomposition \(\Gbf_\ad = \Gbf_{\ad,\lin} \times \Gbf_{\her,\ad}\).
Letting \(\Xcal_1\) be a \(\Gbf_{\her}(\R)\)-orbit in \(\Xcal\), we see that \((\Gbf_\her, \Xcal_1)\) is a pure Shimura datum in the sense of \cite[\S 2.1]{Pink_dissertation}.

Following \cite[\S 3.1]{Pink_ladic_Shim} we also make the following assumption.
\begin{assu} \label{ass:gen_Shim_easy_center}
  The connected center of \(\Gbf / \Abf_{\Gbf}\) stays anisotropic after base changing along \(\Q \hookrightarrow \R\).
\end{assu}

\begin{rema} \label{rem:isog_Gher_gen_Shim}
  One could replace \(\Gbf_{\her}\) by any connected group between \(\Gbf_{\her}\) and \(\Gbf_{\her} Z(\Gbf)^0\).
  It is very formal to check that the subsequent definitions (e.g.\ Definition \ref{def:gen_Shim_var} and Proposition-Definition \ref{prodef:AFcal}) do not depend on this choice.
\end{rema}

\begin{exam} \label{ex:direct_prod_gen_Shim}
  Starting from a Shimura datum \((\Gbf_{\her}, \Xcal_1, h)\) (in the sense of \cite[\S 2.1]{Pink_dissertation}) and a connected reductive group \(\Gbf_{\lin}\) over \(\Q\) we can form a generalized Shimura datum with \(\Gbf = \Gbf_{\lin} \times_{\Spec \Q} \Gbf_{\her}\) and \(\Xcal = \Xcal_1\).
\end{exam}
This is in fact the only case needed in the main part of the paper.

Let \(\Xcal' = \Gbf(\R) / K_\infty \Abf_{\Gbf}(\R)^0\), a real manifold isomorphic to \(\R^N\) for some integer \(N\).
In the following discussion it could be replaced by any ``large enough'' contractible space with an action of \(\Gbf(\Q)\), for example any space \(E \Gbf(\Q)\) with a proper free action of \(\Gbf(\Q)\), giving a model \(\Gbf(\Q) \backslash E \Gbf(\Q)\) of \(B \Gbf(\Q)\).
Consider the projections, for \(K\) a compact open subgroup of \(\Gbf(\A_f)\)
\begin{equation} \label{eq:proj_to_gen_Shim}
  \Gbf(\Q) \backslash (\Xcal \times \Xcal' \times \Gbf(\A_f) / K) \longrightarrow \Gbf(\Q) \backslash (\Xcal \times \Gbf(\A_f) / K).
\end{equation}
We will see below that the right-hand side is a naturally a quasi-projective complex variety when \(K\) is neat (Deligne-Mumford stack in general) even though in general the action of \(\Gbf(\Q)\) is not proper.
The action of \(\Gbf(\Q)\) on \(\Xcal \times \Xcal' \times \Gbf(\A_f) / K\) is proper because it is already propre on \(\Xcal' \times \Gbf(\A_f) / K\), and if \(K\) is neat this action is moreover free.
The source of \eqref{eq:proj_to_gen_Shim} is not a complex manifold (it may have odd real dimension) in general however, in particular it cannot be algebraized.

The spaces \(\Gbf(\Q) \backslash (\Xcal \times \Gbf(\A_f)/K) \) may be described using Shimura varieties.
Let \(S\) be the stabilizer of \(\Xcal_1\) in \(\Gbf(\R)\), an open subgroup of \(\Gbf(\R)\).
Denote \(S_{\Q} = S \cap \Gbf(\Q)\).
Similarly let \(C\) be the centralizer of \(\Xcal_1\), i.e.\ the subgroup of elements of \(S\) fixing \(\Xcal_1\) pointwise, and \(C_{\Q} = C \cap \Gbf(\Q)\).
We have an isomorphism
\[ S_{\Q} \backslash (\Xcal_1 \times \Gbf(\A_f)/K) \simeq \Gbf(\Q) \backslash (\Xcal \times \Gbf(\A_f)/K) \]
because \(\Gbf(\Q)\) meets every connected component of \(\Gbf(\R)\).
The left-hand side decomposes as
\begin{equation} \label{eq:decomp_gen_Shim}
  \bigsqcup_{[g] \in S_\Q \Gbf_{\her}(\A_f) \backslash \Gbf(\A_f) / K} S_{\Q} \backslash (\Xcal_1 \times S_{\Q} \Gbf_{\her}(\A_f) g K / K).
\end{equation}
Denoting \(P(gK) := gKg^{-1} \cap \Gbf_{\her}(\A_f)\) and
\[ \Sh(\Gbf_{\her}, \Xcal_1, P(gK))(\C) := \Gbf_{\her}(\Q) \backslash (\Xcal_1 \times \Gbf_{\her}(\A_f) / P(gK)) \]
we have surjective maps (right multiplication by \(g\))
\begin{equation} \label{eq:surj_to_piece_generalized_Shim}
  \Sh(\Gbf_{\her}, \Xcal_1, P(gK))(\C) \longrightarrow S_{\Q} \backslash (\Xcal_1 \times S_{\Q} \Gbf_{\her}(\A_f) g K / K).
\end{equation}
In order to describe them as torsors under finite groups we introduce a slight generalization of the action of \(\Gbf_{\her}(\A_f)\) on the tower of Shimura varieties for the Shimura datum \((\Gbf_{\her}, \Xcal_1)\).
Note that \(C\) is an open subgroup of the centralizer of \(\Gbf_{\her,\R}\) in \(\Gbf(\R)\), and that \(S_{\Q}/C_{\Q}\) is a group of automorphisms of this Shimura datum.

\begin{prodef} \label{def:action_gen_Shim_tower}
  The group \(S_{\Q} \Gbf_{\her}(\A_f)\) has a right action on the tower \((\Sh(\Gbf_{\her}, \Xcal_1, K_{\her})(\C))_{K_\her}\) of Shimura varieties, where \(K_{\her}\) ranges over neat compact open subgroup of \(\Gbf_{\her}(\A_f)\).
  For \(t = s g_{\her} \in S_{\Q} \Gbf_{\her}(\A_f)\) and \(K'_{\her} \subset t K_{\her} t^{-1}\) we have a map
  \begin{align*}
    T_t: \Sh(\Gbf_{\her}, \Xcal_1, K'_{\her})(\C) & \longrightarrow \Sh(\Gbf_{\her}, \Xcal_1, K_{\her})(\C) \\
    [x, gK'_{\her}] & \longmapsto [s^{-1} \cdot x, s^{-1} g t K_{\her}]
  \end{align*}
  In other words \(T_t\) is the composition \(T_{g_{\her}} \circ T_s\) with \(T_s\) the isomorphism induced by the automorphism of the Shimura datum \((\Gbf_{\her}, \Xcal_1)\) induced by \(s^{-1}\).
  This map does not depend on the choice of decomposition \(t = s g_{\her}\).
\end{prodef}
\begin{proof}
  First we note that \(s^{-1} g t = (s^{-1} g s) g_{\her}\) belongs to \(\Gbf_{\her}(\A_f)\).
  For \(k \in K'_{\her}\) we have \(g k t K_{\her} = g t (t^{-1} k t) K_{\her} = g t K_{\her}\).
  For \(\gamma \in S_{\Q} \cap \Gbf_{\her}(\A_f) = \Gbf_{\her}(\Q)\), replacing \((s,g_{\her})\) by \((s \gamma^{-1}, \gamma g_{\her})\) we find \([\gamma s^{-1} \cdot x, \gamma s^{-1} g t K_{\her}] = [s^{-1} \cdot x, s^{-1} g t K_{\her}]\).
\end{proof}
Note that the distinguished subgroup \(C_{\Q}\) of \(S_{\Q} \Gbf_{\her}(\A_f)\) acts trivially on the tower.
For \(K'_{\her} \subset tK_{\her}t^{-1}\) and \(K''_{\her} \subset t' K'_{\her} (t')^{-1}\) we have \(T_t \circ T_{t'} = T_{t't}\) as morphisms \(\Sh(\Gbf_{\her}, \Xcal_1, K''_{\her}) \to \Sh(\Gbf_{\her}, \Xcal_1, K_{\her})\): this follows formally from the functoriality of the maps \(T_{g_1}\) with respect to Shimura datum isomorphisms.
In particular we have \(T_{tk} = T_t\) for any \(k \in K_{\her}\).

Returning to \(K\) a compact open subgroup of \(\Gbf(\A_f)\) and \(gK \in \Gbf(\A_f)/K\) denote also \(Q(gK) = g K g^{-1} \cap S_{\Q} \Gbf_{\her}(\A_f)\) and \(C(gK) = g K g^{-1} \cap C_{\Q}\).

\begin{lemm} \label{lem:CgK_PgK_QgK}
  Let \(\Gbf_{\lin,\der} := \Cent(\Gbf_{\her}, \Gbf)^0_{\der}\) (derived subgroup of neutral component of centralizer).
  Assume that \(K\) is neat.
  Then \(C(gK)\) is a torsion-free arithmetic subgroup of \(\Gbf_{\lin,\der}(\Q)\) and the multiplication map identifies \(C(gK) \times P(gK)\) (with the product topology) with an open subgroup of \(Q(gK)\) having finite index.
\end{lemm}
\begin{proof}
  We show first that \(C(gK)\) is an arithmetic subgroup of \(\Gbf_{\lin,\der}(\Q)\).
  For \(x_0 \in h(\Xcal)\) the \((\Gbf/\Gbf_{\lin,\der})(\R)\)-orbit of \(x_0\) in \(\Hom(\mathbb{S}, (\Gbf/\Gbf_{\lin,\der})_{\R})\) is a Shimura datum still satisfying Assumption \ref{ass:gen_Shim_easy_center} and so the image of \(C(gK)\) in \((\Gbf/\Gbf_{\lin,\der})(\Q)\) is trivial, i.e.\ \(C(gK)\) is contained in the arithmetic subgroup \(\Gbf_{\lin,\der}(\Q) \cap gKg^{-1}\) of \(\Gbf_{\lin,\der}(\Q)\).
  It is torsion-free because \(\Gbf_{\lin,\der}(\A_f) \cap gKg^{-1}\) is neat.
  Let \(\Gbf_{\lin,\sico}\) be the simply connected cover of \(\Gbf_{\lin,\der}\), then \(\Gbf_{\lin,\sico}(\R)\) is connected and so its image in \(\Gbf(\R)\) is contained in \(C\).
  It follows that \(C(gK)\) contains the image of the intersection of \(\Gbf_{\lin,\sico}(\Q)\) with the preimage of \(gKg^{-1}\), so by \cite[Corollary 6.11]{BorelHC_arith_subgps} it contains an arithmetic subgroup of \(\Gbf_{\lin,\der}(\Q)\).

  The intersection \(C_{\Q} \cap \Gbf_{\her}(\A_f)\) is contained in the center of \(\Gbf_{\her}(\Q)\) and so by neatness of \(K\) and Assumption \ref{ass:gen_Shim_easy_center} we have \(C(gK) \cap P(gK) = \{1\}\).
  Note that \(S_{\Q} \Gbf_{\her}(\A_f)\) is the preimage of \(S_{\Q} / \Gbf_{\her}(\Q) \subset (\Gbf/\Gbf_{\her})(\Q)\) under the projection \(\Gbf(\A_f) \to (\Gbf/\Gbf_{\her})(\A_f)\), and so \(S_{\Q} \Gbf_{\her}(\A_f)\) is a closed subgroup of \(\Gbf(\A_f)\) and \(\Gbf_{\her}(\A_f)\) is an open subgroup of \(S_{\Q} \Gbf_{\her}(\A_f)\).
  It follows that \(P(gK)\) is an open subgroup of \(Q(gK)\), and so the multiplication map \(C(gK) \times P(gK) \to Q(gK)\), identifies \(C(gK) \times P(gK)\) with an open subgroup.
  We are left to show that this subgroup has finite index.
  We have a map
  \begin{align*}
    Q(gK) / C(gK) P(gK) & \longrightarrow \Gbf_{\her}(\Q) \backslash \Gbf_{\her}(\A_f) / P(gK) \\
    [\gamma g_{\her}] & \longmapsto [g_{\her}]
  \end{align*}
  where \(\gamma \in S_{\Q}\) and \(g_{\her} \in \Gbf_{\her}(\A_f)\), which is well-defined because \(C(gK)\) commutes with \(\Gbf_{\her}(\A_f)\).
  The target of this map is finite \cite[Théorème 5]{Godement_BBki_domfond} and we are left to show that its fibers are finite.
  If \(\gamma_1 g_{\her,1}\) and \(\gamma_2 g_{\her,2}\) map to the same class then up to right multiplication by an element of \(P(gK)\) we may assume \(g_{\her,2} \in \Gbf_{\her}(\Q) g_{\her,1}\) for some \(\gamma_{\her} \in \Gbf_{\her}(\Q)\), and since \(\Gbf_{\her}(\Q)\) is contained in \(S_{\Q}\) we may even assume \(g_{\her,1} = g_{\her,2}\).
  Thus it is enough to show that for any \(g_{\her} \in \Gbf_{\her}(\A_f)\) the quotient
  \[ \{ \gamma \in S_{\Q} \,|\, \gamma g_{\her} \in gKg^{-1} \} / C(gK) (\Gbf_{\her}(\Q) \cap g_{\her} gK g^{-1} g_{\her}^{-1}) \]
  is finite.
  This quotient is either empty or in bijection with
  \[ (S_{\Q} \cap (g_\her g K g^{-1} g_{\her}^{-1})) / C(gK) (\Gbf_{\her}(\Q) \cap g_{\her} gK g^{-1} g_{\her}^{-1}), \]
  so it is enough to show that
  \[ (\Gbf(\Q) \cap (g_\her g K g^{-1} g_{\her}^{-1})) / C(gK) (\Gbf_{\her}(\Q) \cap g_{\her} gK g^{-1} g_{\her}^{-1}) \]
  is finite.
  There exists a central torus \(\Tbf\) of \(\Gbf\) such that the multiplication morphism \(\Tbf \times \Gbf_{\lin,\der} \times \Gbf_{\her} \to \Gbf\) is a central isogeny.
  Assumption \ref{ass:gen_Shim_easy_center} implies \(\Tbf(\Q) \cap K = \{1\}\), so another application of \cite[Corollary 6.11]{BorelHC_arith_subgps} tells us that \(C(gK) (\Gbf_{\her}(\Q) \cap g_{\her} gK g^{-1} g_{\her}^{-1})\) is an arithmetic subgroup of \(\Gbf(\Q)\).
\end{proof}

Let \(R(gK) := Q(gK)/C(gK)\), which by Lemma \ref{lem:CgK_PgK_QgK} is a profinite topological group.
The natural map \(P(gK) \to R(gK)\) identifies \(P(gK)\) with an open subgroup of \(R(gK)\).
The finite quotient \(R(gK) / P(gK)\) acts on \(\Sh(\Gbf_{\her}, \Xcal_1, P(gK))(\C)\) via the maps \(T_t\) defined above.

\begin{prop} \label{pro:gen_Shim_var_no_fixed}
  For a neat compact open subgroup \(K\) of \(\Gbf(\A_f)\) and \(g \in \Gbf(\A_f)\) the left action of \((S_{\Q} \cap gKg^{-1}) / C(gK)\) on \(\Xcal_1\) is proper and free.
  Similarly the right action of \(R(gK) / P(gK)\) on \(\Sh(\Gbf_{\her}, \Xcal_1, P(gK))(\C)\) is free and exhibits \eqref{eq:surj_to_piece_generalized_Shim} as a right \(R(gK) / P(gK)\)-torsor.
\end{prop}
\begin{proof}
  See \cite[Proposition 3.7.5]{Pink_ladic_Shim}.
\end{proof}
We also see that over the piece corresponding to \([g]\) in the decomposition \eqref{eq:decomp_gen_Shim}, the map \eqref{eq:proj_to_gen_Shim} is a fibration with fiber \(C(gK) \backslash \Xcal'\), in fact base changing \eqref{eq:proj_to_gen_Shim} along the finite étale map \(\Sh(\Gbf_{\her}, \Xcal_1, P(gK)) \to \Gbf(\Q) \backslash (\Xcal \times \Gbf(\A_f)/K)\) given by \eqref{eq:surj_to_piece_generalized_Shim} and \eqref{eq:decomp_gen_Shim} yields a trivial fibration.

Let \(F\) be the reflex field of \((\Gbf_{\her}, \Xcal_1)\).
We would like to have a tower of smooth quasi-projective schemes \(\Sh(\Gbf, \Xcal, K)\) over \(F\), where \(K\) ranges over neat compact open subgroups of \(\Gbf(\A_f)\), and a right action of \(\Gbf(\A_f)\) on this tower, along with identifications (isomorphisms of complex manifolds) \(\Sh(\Gbf, \Xcal, K)(\C) \simeq \Gbf(\Q) \backslash (\Xcal \times \Gbf(\A_f)/K)\) compatible with transition maps and the action of \(\Gbf(\A_f)\).
It should be possible to redo the whole theory with this more general definition of Shimura datum, but one can simply reduce to ``classical'' Shimura varieties as follows.
For \(K_{\her}\) a neat compact open subgroup of \(\Gbf_{\her}(\A_f)\) denote by \(\Sh(\Gbf_{\her}, \Xcal_1, K_{\her})\) the canonical model, a smooth quasi-projective scheme over \(F\).
The action of \(S_{\Q} \Gbf_{\her}(\A_f)\) on the tower \((\Sh(\Gbf_{\her}, \Xcal_1, K_{\her}))_{K_{\her}}\) introduced in Proposition-Definition \ref{def:action_gen_Shim_tower} descends to canonical models: by \cite[Proposition 11.10]{Pink_dissertation} for any \(s \in S_\Q\) its inverse induces isomorphisms
\[ T_s: \Sh(\Gbf_{\her}, \Xcal_1, K_{\her}) \xrightarrow{\sim} \Sh(\Gbf_{\her}, \Xcal_1, s^{-1}K_{\her}s) \]
and the arguments of Proposition-Definition \ref{def:action_gen_Shim_tower} adapt to show that we have well-defined finite étale maps \(T_t\) between canonical models.
In particular for \(K\) a neat compact open subgroup of \(\Gbf(\A_f)\) and \(gK \in \Gbf(\A_f)/K\) we have a right action of \(R(gK)/P(gK)\) on \(\Sh(\Gbf_{\her}, \Xcal_1, P(gK))\) and this action is free because it is so on complex points.
We obtain \cite[\href{https://stacks.math.columbia.edu/tag/07S7}{Tag 07S7}]{stacks-project} smooth quasi-projective quotients
\[ U(gK) := \Sh(\Gbf_{\her}, \Xcal_1, P(gK)) / (R(gK)/P(gK)) \]
and morphisms \(T_t\) between them, for \(t \in S_{\Q} \Gbf_{\her}(\A_f)\).
We have a functor
\begin{align*}
  F(\Gbf, \Xcal, h, \Gbf_{\her}, \Xcal_1, K): [S_\Q \Gbf_{\her}(\A_f) \curvearrowright \Gbf(\A_f)/K] & \longrightarrow \mathrm{Sch} \\
  gK & \longmapsto U(gK) \\
  \left( gK \xrightarrow{t} tgK \right) & \longmapsto \left( U(gK) \xrightarrow{T_{t^{-1}}} U(tg K) \right).
\end{align*}

\begin{defi} \label{def:gen_Shim_var}
  Let \((\Gbf, \Xcal, h)\) be a generalized Shimura datum satisfying assumption \ref{ass:gen_Shim_easy_center}.
  For \(K\) a neat compact open subgroup of \(\Gbf(\A_f)\) define the (quasi-projective and smooth) scheme \(\Sh(\Gbf, \Xcal, K)\) over the reflex field \(F\) as the colimit of the functor \(F(\Gbf, \Xcal, h, \Gbf_{\her}, \Xcal_1, K)\).
\end{defi}

For \(gK \in \Gbf(\A_f)/K\) and \(t \in \Aut_{[S_\Q \Gbf_{\her}(\A_f) \curvearrowright \Gbf(\A_f)/K]}(gK) = Q(gK)\) we have \(T_{t^{-1}} = \id\) and so choosing representatives yields an isomorphism
\begin{equation} \label{eq:gen_Shim_var_explicit}
  \Sh(\Gbf, \Xcal, K) \simeq \bigsqcup_{[g] \in S_{\Q} \Gbf_{\her}(\A_f) \backslash \Gbf(\A_f) / K} U(gK).
\end{equation}
We denote by
\[ \iota_g: U(gK) \hookrightarrow \Sh(\Gbf, \Xcal, K) \]
the canonical embedding.
It is very formal to check that choosing a different \(\Gbf_{\her}(\R)\)-orbit \(\Xcal_1\) in \(\Xcal\) yields a canonically isomorphic colimit, which is why \(\Xcal_1\) is absent from the notation \(\Sh(\Gbf, \Xcal, K)\).

\begin{exam} \label{ex:direct_prod_gen_Shim_var}
  Consider the case of a direct product (Example \ref{ex:direct_prod_gen_Shim}) where \(\Gbf = \Gbf_{\lin} \times \Gbf_{\her}\) and \((\Gbf_{\her}, \Xcal, h)\) is a Shimura datum.
  We have \(S = \Gbf(\R)\), \(S_{\Q} \Gbf_{\her}(\A_f) = \Gbf_{\lin}(\Q) \times \Gbf_{\her}(\A_f)\) and \(\Gbf_{\lin}(\R) \times Z(\Gbf_{\her}(\R))^0\) is an open subgroup of finite index in \(C\).
  For any neat compact open subgroup \(K\) of \(\Gbf(\A_f)\) and any \(gK \in \Gbf(\A_f)/K\) we have \(C(gK) = \Gbf_{\lin}(\Q) \cap gKg^{-1}\) (thanks to Assumption \ref{ass:gen_Shim_easy_center}) and so \(R(gK)\) is identified with a neat compact open subgroup of \(\Gbf_{\her}(\A_f)\), and \(U(gK)\) is identified with \(\Sh(\Gbf_{\her}, \Xcal, R(gK))\).
  In the case where \(K = K_{\lin} \times K_{\her}\) the situation is even simpler: we have \(R(gK) \simeq P(gK) = g_{\her} K_{\her} g_{\her}^{-1}\) and \(S_{\Q} \Gbf_{\her}(\A_f) \backslash \Gbf(\A_f) / K \simeq \Gbf_{\lin}(\Q) \backslash \Gbf_{\lin}(\A_f) / K_{\lin}\).
\end{exam}

\begin{rema} \label{rem:models_gen_Shim}
  When we have integral models for the Shimura varieties \(\Sh(\Gbf_{\her}, \Xcal_1, K_{\her})\), say for \(K_{\her}\) hyperspecial at \(p\) for some prime \(p\), we would like to construct integral models of the generalized Shimura varieties \(\Sh(\Gbf, \Xcal, K)\) as well.
  For this one needs to check the analogue of the second part of Proposition \ref{pro:gen_Shim_var_no_fixed} for the special fibers (in characteristic \(p\)) of \(\Sh(\Gbf_{\her}, \Xcal_1, P(gK))\).
  In the direct product case (Example \ref{ex:direct_prod_gen_Shim_var}) this reduces to the fact that \(\Sh(\Gbf_{\her}, \Xcal_1, K'_{\her}) \to \Sh(\Gbf_{\her}, \Xcal_1, K_{\her})\) is a right \(K_{\her}/K'_{\her}\)-torsor whenever \(K'_{\her}\) is a distinguished open subgroup of the neat level \(K_{\her}\).

  In the rest of this appendix we discuss the case of generalized Shimura varieties over the reflex field.
  Under the above hypotheses (hyperspecial level, direct product case) everything generalizes to integral models, but we leave the details implicit.
\end{rema}

For \(h \in \Gbf(\A_f)\) and neat compact open subgroups \(K\) and \(K'\) of \(\Gbf(\A_f)\) satisfying \(K' \subset h K h^{-1}\) we have a well-defined morphism \(T_h\) making the following diagrams commutative, whenever \(g,g' \in \Gbf(\A_f)\) and \(t \in S_\Q \Gbf_{\her}(\A_f)\) satisfy \(g'h \in tgK\):
\begin{equation} \label{eq:def_T_h_Pink}
  \begin{tikzcd}
    U(g'K') \ar[r, hook, "{\iota_{g'}}"] \ar[d, "{T_t}"] & \Sh(\Gbf, \Xcal, K') \ar[d, "{T_h}"] \\
    U(gK) \ar[r, hook, "{\iota_g}"] & \Sh(\Gbf, \Xcal, K)
  \end{tikzcd}
\end{equation}
For any \(k \in K\) we have \(T_{hk} = T_h\).
For any \(h' \in \Gbf(\A_f)\) we have \(T_h T_{h'} = T_{h' h}\).

\begin{lemm} \label{lem:Pink_Shim_as_quot}
  Let \(K'\) be a distinguished open subgroup of a neat compact open subgroup \(K\) of \(\Gbf(\A_f)\).
  We have an action of the finite group \(K/K'\) on \(\Sh(\Gbf, \Xcal, K')\) via the maps \((T_h)_{h \in K/K'}\) defined above (see \eqref{eq:def_T_h_Pink}).
  The finite étale map
  \[ T_1: \Sh(\Gbf, \Xcal, K') \longrightarrow \Sh(\Gbf, \Xcal, K). \]
  is surjective and \(K/K'\) acts transitively on its geometric fibers.
  The stabilizer of a geometric point of the component of \(\Sh(\Gbf, \Xcal, K')\) corresponding to \([g] \in S_{\Q} \Gbf_{\her}(\A_f) \backslash \Gbf(\A_f) / K'\) is the image of \(K \cap g^{-1} C_{\Q} g = g^{-1} C(gK) g\) in \(K/K'\).
\end{lemm}
\begin{proof}
  It is clear that \(K/K'\) acts transitively on the fibers of the surjective map
  \[ S_{\Q} \Gbf_{\her}(\A_f) \backslash \Gbf(\A_f) / K' \longrightarrow S_{\Q} \Gbf_{\her}(\A_f) \backslash \Gbf(\A_f) / K. \]
  The stabilizer in \(K/K'\) of \([g] \in S_{\Q} \Gbf_{\her}(\A_f) \backslash \Gbf(\A_f) / K'\) is the image of \(K \cap g^{-1} S_{\Q} \Gbf_{\her}(\A_f) g = g^{-1} Q(gK) g\) in \(K/K'\).
  So we may fix \(g \in \Gbf(\A_f)\) and we are left to check that the finite étale map
  \begin{equation} \label{eq:}
    U(gK') \xrightarrow{T_1} U(gK)
  \end{equation}
  is surjective, that \(Q(gK)/Q(gK')\) acts transitively on its geometric fibers and that the stabilizer of any geometric point of \(U(gK')\) is the image of \(C(gK)\).
  This is clear because the natural map
  \[ \Sh(\Gbf_{\her}, \Xcal_1, P(gK')) \longrightarrow U(gK) \]
  is a right \(R(gK) / P(gK')\)-torsor thanks to Proposition \ref{pro:gen_Shim_var_no_fixed}.
\end{proof}

\subsection{Hecke formalism}
\label{sec:gen_Shim_Hecke}

Let \(\ell\) be a prime and \(E\) a finite extension of \(\Qell\).
When working with integral models we also assume that \(\ell\) is invertible over the base.
Let \(K\) be a neat compact open subgroup of \(\Gbf(\A_f)\) and \(g \in \Gbf(\A_f)\).
For any open subgroup \(H\) of \(P(gK)\) which is distinguished in \(R(gK)\) we have a finite étale map
\[ T_1: \Sh(\Gbf_{\her}, \Xcal_1, H) \longrightarrow U(gK) \]
and an action of \(R(gK)/H\) on the source which identifies \(U(gK)\) with the quotient.
These maps and actions are compatible with change of \(H\).
This gives us, as a special case of \cite[(1.10)]{Pink_ladic_Shim}, a morphism from the étale topos of \(U(gK)\) to the topos \(\Sets_{R(gK)}\) of (discrete) sets with continuous action of \(R(gK)\).
More precisely,
\begin{itemize}
\item the pullback functor \(\Fcal^{R(gK)}\) maps a set \(A\) with continuous action of \(R(gK)\) to the colimit over \(H\) (as above) of the sheaf on \(U(gK)\) associated to the \(R(gK)/H\)-equivariant constant sheaf on \(\Sh(\Gbf_{\her}, \Xcal_1, H)\) corresponding to \(A^H\),
\item the pushforward functor maps an étale sheaf \(\Gcal\) on \(U(gK)\) to
  \[ \colim_H \Gcal \left( \Sh(\Gbf_{\her}, \Xcal_1, H) \xrightarrow{T_1} U(gK) \right). \]
\end{itemize}
Ekedahl \cite{Ekedahl_adic} defined an \(\Ocal_E\)-linear triangulated category \(D^+(R(gK), \Ocal_E)\) associated to the topos \(\Sets_{R(gK)}\), the constant ring object \(\Ocal_E\) in this topos and the maximal ideal \(\mfrak_E\) of \(\Ocal_E\).
We only recall that this is a certain quotient of a certain subcategory of the derived category of \((\Ocal_E)_\bullet\)-modules, where \((\Ocal_E)_\bullet\) is the ring object \((\Ocal_E/\mfrak_E^i)_{i \geq 0}\) in the topos \(\Sets_{R(gK)}^\N\).
As in \cite{Taibi_vanEst} we denote by \(D^+(R(gK), E)\) the triangulated \(E\)-linear category obtained from \(D^+(R(gK), \Ocal_E)\) by inverting \(\ell\).
Similarly we consider the triangulated \(E\)-linear category \(D^+(Q(gK), E)\), and we have a triangulated functor \(R\Gamma(C(gK), -): D^+(Q(gK),E) \to D^+(R(gK), E)\).

Now let \(V\) be a bounded object in the derived category of the abelian category of continuous representations of \(\Gbf(\Qell)\) over finite-dimensional \(E\)-vector spaces.
For any neat compact open subgroup \(K\) of \(\Gbf(\A_f)\) and any \(t \in S_\Q \Gbf_{\her}(\A_f)\) we have a composite isomorphism \(s(Q(K),t)\)
\begin{align} \label{eq:T_t_pullback_Pink}
  T_t^* \Fcal^{R(K)} R\Gamma(C(K), V)
  &\simeq \Fcal^{R(tK)} \ol{\Ad(t^{-1})}^* R\Gamma(C(K), V) \nonumber \\
  &\simeq \Fcal^{R(tK)} R\Gamma(C(tK), \Ad(t^{-1})^* V) \nonumber \\
  &\simeq \Fcal^{R(tK)} R\Gamma(C(tK), V).
\end{align}
Here
\begin{itemize}
\item We abusively still denote by \(V\) its image in \(D^+(Q(K),E)\) (see \cite[Corollary 2.5]{Taibi_vanEst}),
\item the first isomorphism follows from the definition of \(\Fcal\) and the fact that the maps \(T_t\) induce a morphism between towers of Shimura varieties above \(T_t: U(tK) \to U(K)\), intertwining the actions of \(R(tK)\) and \(R(K)\) via the isomorphism \(\ol{\Ad(t^{-1})}: R(tK) \simeq R(K)\) induced by \(\Ad(t^{-1}): Q(tK) \simeq Q(K)\),
\item the second isomorphism is completely formal using \(C(tK) = tC(K)t^{-1}\),
\item the third isomorphism is the action of \(t\) on \(V\).
\end{itemize}

\begin{prodef} \label{prodef:AFcal}
  Let \(K\) be a neat compact open subgroup of \(\Gbf(\A_f)\).
  There exists an object \(\AFcal^K(V)\) of \(D^+(\Sh(\Gbf, \Xcal, K), E)\) together with isomorphisms, for all \(gK \in \Gbf(\A_f)/K\),
  \[ \iota_g^* \AFcal^K(V) \simeq \Fcal^{R(gK)} R\Gamma(C(gK), V) \]
  such that for any \(t \in S_{\Q} \Gbf_{\her}(\A_f)\) the following diagram is commutative.
  \[
    \begin{tikzcd}
      T_t^* \iota_g^* \AFcal^K (V) \ar[r, dash, "{\sim}"] \ar[d, equal, "{\sim}"] & T_t^* \Fcal^{R(gK)} R\Gamma(C(gK),V) \ar[d, "{s(Q(gK),t)}" right, "{\sim}" left] \\
      \iota_{tg}^* \AFcal^K (V) \ar[r, dash, "{\sim}"] & \Fcal^{R(tgK)} R\Gamma(C(tgK),V)
    \end{tikzcd}
  \]
  This characterizes \(\AFcal^K(V)\).
\end{prodef}
\begin{proof}
  Uniqueness follows from \eqref{eq:gen_Shim_var_explicit}.
  It follows from Lemma \ref{lem:compo_s_Pink} below (we are in the case where the inclusions between levels are equalities) that the formation of \(s(Q(K),t)\) is compatible with composition, i.e.\ for \(t_1,t_2 \in S_{\Q} \Gbf_{\her}(\A_f)\) the composition
  \begin{align*}
    & T_{t_2}^* T_{t_1}^* \Fcal^{R(K)} R\Gamma(C(K),V) \\
    \xrightarrow{T_{t_2}^* s(Q(K),t_1)} & T_{t_2}^* \Fcal^{R(t_1K)} R\Gamma(C(t_1K), V) \\
    \xrightarrow{s(Q(t_1K),t_2)} & \Fcal^{R(t_2t_1K)} R\Gamma(C(t_2 t_1 K), V)
  \end{align*}
  is equal to \(s(Q(K), t_2 t_1)\) (implicitly using \(T_{t_2}^* T_{t_1}^* \simeq T_{t_2 t_1}^*\)).
  We can use this to show the analogue of \eqref{eq:iso_func_pullback}: for any \(q \in Q(K)\) we have \(T_{tq}=T_t\), \(R(tqK)=R(tK)\), \(C(tqK)=C(tK)\) and \(s(Q(K),tq)=s(Q(K),t)\) (the cocycle relation reduces it to the case where \(t=1\), which is formal).
  This invariance property and the cocycle relation imply existence.
\end{proof}

Our next goal is to prove a finiteness property of \(\AFcal^K(V)\) in Proposition \ref{pro:RGamma_CgK_finite} below.
Denote \(\pi_{\lin}: \Gbf \to \Gbf/\Gbf_{\her}\) and \(D(gK) = \pi_{\lin}(gKg^{-1}) \cap (\Gbf/\Gbf_{\her})(\Q)\), a neat arithmetic subgroup of \((\Gbf/\Gbf_{\her})(\Q)\) (the morphism \(\Gbf(\A_f) \to (\Gbf/\Gbf_{\her})(\A_f)\) induced by \(\pi_{\lin}\) is open because \(\Gbf_{\her}\) is connected).

\begin{lemm} \label{lem:QgK_DgK_RgK}
  We keep the assumption that \(K\) is neat.
  The natural morphism \(Q(gK) \to D(gK) \times R(gK)\) is injective and identifies \(Q(gK)\) with an open finite index subgroup.
\end{lemm}
\begin{proof}
  The kernel of this morphism is \(C(gK) \cap \Gbf_{\her}(\Q)\), which is trivial thanks to Assumption \ref{ass:gen_Shim_easy_center}.
  By Lemma \ref{lem:CgK_PgK_QgK} it is enough to show that \(C(gK)\) has finite index in \(D(gK)\).
  Thanks to this lemma we know that \(C(gK)\) is an arithmetic subgroup of \(\Gbf_{\lin,\der}(\Q)\), and the natural map \(\Gbf_{\lin,\der} \times \Zbf(\Gbf/\Gbf_{\her})^0 \to \Gbf/\Gbf_{\her}\) is a central isogeny.
  Thanks to Assumption \ref{ass:gen_Shim_easy_center} we have \(\Zbf(\Gbf/\Gbf_{\her})^0(\Q) \cap K = \{1\}\) and yet another application of \cite[Corollary 6.11]{BorelHC_arith_subgps} shows that the image of \(C(gK)\) in \((\Gbf/\Gbf_{\her})(\Q)\) is an arithmetic subgroup.
\end{proof}

\begin{prop} \label{pro:RGamma_CgK_finite}
  Assume that \(K\) is neat and that \(V\) is a finite complex of continuous finite-dimensional representations of \(K\).
  Then the object \(R\Gamma(C(gK),V)\) of \(D^+(R(gK),E)\) is isomorphic to the image\footnote{See \cite[Corollary 2.5]{Taibi_vanEst}} of a finite complex of finite free \(\Ocal_E\)-modules with continuous action of \(R(gK)\).
\end{prop}
\begin{proof}
  There exists a finite complex \(\Lambda\) of finite free \(\Ocal_E\)-modules with continuous action of \(K\) and an isomorphism \(V \simeq E \otimes_{\Ocal_E} \Lambda\), and by definition \(R\Gamma(C(gK),V)\) is represented by \(R\Gamma(C(gK),-)\) applied to the image of \(\Lambda\) in \(D^+(S_{Q(gK)}^{\N}, (\Ocal_E)_\bullet)\).
  The functor \(\Gamma(C(gK),-)\) from \((\Ocal_E)_\bullet\)-modules in \(\Sets_{Q(gK)}^\N\) to \((\Ocal_E)_\bullet\)-modules in \(\Sets_{R(gK)}^\N\) is isomorphic to \(\Gamma(D(gK), \ind_{Q(gK)}^{D(gK) \times R(gK)}(-))\).
  It follows from Lemma \ref{lem:QgK_DgK_RgK} that the induction functor \(\ind_{Q(gK)}^{D(gK) \times R(gK)}(-)\) is right adjoint to the restriction functor, which is obviously exact, and so this induction functor preserves injective objects.
  We deduce an isomorphism of functors from \(D^+(Q(gK),E)\) to \(D^+(R(gK),E)\)
  \[ R\Gamma(C(gK),-) \simeq R\Gamma(D(gK), \ind_{Q(gK)}^{D(gK) \times R(gK)}(-)). \]
  By \cite[11.1.c]{BorelSerre_corners} there exists a finite resolution \(P^\bullet\) of \(\Z\), considered as a \(\Z[D(gK)]\)-module with trivial action of \(D(gK)\), by finite free \(\Z[D(gK)]\)-modules.
  By the same argument as in \cite[Lemma 4.1]{Taibi_vanEst} we have an isomorphism of functors from \(D^+(D(gK) \times R(gK),E)\) to \(D^+(R(gK),E))\) between \(R\Gamma(D(gK),-)\) and \(\Tot^\bullet \left( \Hom_{\Z[D(gK)]}(P^\bullet, -) \right)\), where \(\Tot^\bullet\) denotes taking the total complex of a double complex.
  Finally the complex (considered as an object of Ekedahl's category \(D^+(R(gK), \Ocal_E)\))
  \[ \Tot^\bullet \Hom_{\Z[D(gK)]}(P^\bullet, \ind_{Q(gK)}^{D(gK) \times R(gK)} \Lambda) \]
  is clearly the image of a finite complex of finite free \(\Ocal_E\)-modules with continuous action of \(R(gK)\).
\end{proof}

We shall need a slight generalization of \eqref{eq:T_t_pullback_Pink}: if \(K\) and \(K'\) are neat compact open subgroups of \(\Gbf(\A_f)\) and \(t \in S_{\Q} \Gbf_{\her}(\A_f)\) satisfies \(K' \subset tKt^{-1}\) we have a natural isomorphism (defined as a composition like \eqref{eq:T_t_pullback_Pink})
\begin{equation} \label{eq:T_t_pullback_Pink_gal}
  T_t^* \Fcal^{R(K)} R\Gamma(C(K), V) \xrightarrow{s(Q(K), t, Q(K'))} \Fcal^{R(K')} R\Gamma(C(tK), V).
\end{equation}
On the right-hand side the object \(R\Gamma(C(tK),V)\) of \(D^+(R(tK),E)\) is implicitly restricted along \(R(K') \hookrightarrow R(tK)\) to obtain an object of \(D^+(R(K'),E)\).

\begin{lemm} \label{lem:compo_s_Pink}
  Assume \(K' \subset t_1 K t_1^{-1}\) and \(K'' \subset t_2 K' t_2^{-1}\).
  The following diagrams are commutative
  \[
    \begin{tikzcd}[column sep=3em]
      T_{t_2}^* T_{t_1}^* \Fcal^{R(K)} R\Gamma(C(K),V) \ar[d, "{\sim}"] \ar[r, "{\sim}" below, "{T_{t_2}^* s(Q(K),t_1,Q(K'))}" above] & T_{t_2}^* \Fcal^{R(K')} R\Gamma(C(t_1K),V) \ar[r, "{\res}"] & T_{t_2}^* \Fcal^{R(K')} R\Gamma(C(K'),V) \ar[d, "{\sim}" left, "{s(Q(K'),t_2,Q(K''))}"] \\
      T_{t_2 t_1}^* \Fcal^{R(K)} R\Gamma(C(K),V) \ar[d, "{s(Q(K),t_2t_1,Q(K''))}"] & & \Fcal^{R(K'') R\Gamma(C(t_2K'),V)} \ar[d, "{\res}"] \\
      \Fcal^{R(K'')} R\Gamma(C(t_2t_1K),V) \ar[rr, "{\res}"] & & \Fcal^{R(K'')} R\Gamma(C(K''), V)
    \end{tikzcd}
  \]
  \[
    \begin{tikzcd}[column sep=3em]
      T_{t_2}^* T_{t_1}^* \Fcal^{R(K)} R\Gamma(C(K),V) \ar[d, "{\sim}"] \ar[r, "{\sim}" below, "{T_{t_2}^* s(Q(K),t_1,Q(K'))}" above] & T_{t_2}^* \Fcal^{R(K')} R\Gamma(C(t_1K),V) & \ar[l, "{\cores}"] T_{t_2}^* \Fcal^{R(K')} R\Gamma(C(K'),V) \ar[d, "{\sim}" left, "{s(Q(K'),t_2,Q(K''))}"] \\
      T_{t_2 t_1}^* \Fcal^{R(K)} R\Gamma(C(K),V) \ar[d, "{s(Q(K),t_2t_1,Q(K''))}"] & & \Fcal^{R(K'') R\Gamma(C(t_2K'),V)} \\
      \Fcal^{R(K'')} R\Gamma(C(t_2t_1K),V) & & \ar[ll, "{\res}"] \Fcal^{R(K'')} R\Gamma(C(K''), V) \ar[u, "{\cores}"]
    \end{tikzcd}
  \]
\end{lemm}
\begin{proof}
  To prove the lemma it is useful to introduce yet another slight generalization of \eqref{eq:T_t_pullback_Pink}: if \(K\), \(\tilde{K}\) and \(K'\) are neat compact open subgroups of \(\Gbf(\A_f)\) and \(t \in S_{\Q} \Gbf_{\her}(\A_f)\) satisfy \(K' \subset tKt^{-1}\) and \(K \subset \tilde{K}\) we have a natural isomorphism
  \[ T_t^* \Fcal^{R(K)} \res^{R(\tilde{K})}_{R(K)} R\Gamma(C(\tilde{K}), V) \xrightarrow{s(Q(\tilde{K}), Q(K), t, Q(K'))} \Fcal^{R(K')} \res^{R(t\tilde{K})}_{R(K')} R\Gamma(C(t\tilde{K}), V) \]
  defined as the composition (being explicit with restriction functors)
  \begin{align*}
    T_t^* \Fcal^{R(K)} \res^{R(\tilde{K})}_{R(K)} R\Gamma(C(\tilde{K}), V)
    &\simeq \Fcal^{R(K')} \res^{R(tK)}_{R(K')} \Ad(t^{-1})^* \res^{R(\tilde{K})}_{R(K)} R\Gamma(C(K), V) \\
    &\simeq \Fcal^{R(tK)} \res^{R(tK)}_{R(K')} \res^{R(t\tilde{K})}_{R(tK)} \Ad(t^{-1})^* R\Gamma(C(\tilde{K}), V) \\
    &\simeq \Fcal^{R(tK)} \res^{R(t\tilde{K})}_{R(K')} R\Gamma(C(t\tilde{K}), \Ad(t^{-1})^* V) \\
    &\simeq \Fcal^{R(tgK)} \res^{R(t\tilde{K})}_{R(K')} R\Gamma(C(t\tilde{K}), V)
  \end{align*}
  where (again) the first isomorphism comes from the isomorphism of functors \(T_t^* \Fcal^{R(K)} \simeq \res^{R(\tilde{tK})}_{R(K')} \Ad(t^{-1})^*\), the middle isomorphisms are very formal and the last isomorphism is the action of \(t\) on \(V\).
  First we claim that these isomorphisms are compatible with composition, in other words they satisfy a cocyle relation: for neat compact open subgroups \(\tilde{K}\), \(K\), \(K'\) and \(K''\) of \(\Gbf(\A_f)\) and \(t_1,t_2 \in S_{\Q} \Gbf_{\her}(\A_f)\) satisfying \(K \subset \tilde{K}\), \(K' \subset t_1 K t_1^{-1}\) and \(K'' \subset t_2 K' t_2^{-1}\) the composition
  \begin{align*}
    & T_{t_2}^* T_{t_1}^* \Fcal^{R(K)} \res^{R(\tilde{K})}_{R(K)} R\Gamma(C(\tilde{K}),V) \\
    \xrightarrow{T_{t_2}^* s(Q(\tilde{K}),Q(K),t_1,Q(K'))} & T_{t_2}^* \Fcal^{R(K')} \res^{R(t_1\tilde{K})}_{R(K')} R\Gamma(C(t_1 \tilde{K}), V) \\
    \xrightarrow{s(Q(t_1\tilde{K}),Q(K'),t_2,Q(K''))} & \Fcal^{R(K'')} \res^{R(t_2 t_1 \tilde{K})}_{R(K'')} R\Gamma(C(t_2 t_1 \tilde{K}), V)
  \end{align*}
  is equal to \(s(Q(\tilde{K}),Q(K), t_2 t_1, Q(K''))\) (implicitly using \(T_{t_2}^* T_{t_1}^* \simeq T_{t_2 t_1}^*\)).
  We do not give details for the rather tedious proof, which simply goes through each step in the definition of \(s(Q(\tilde{K}),Q(K),t,Q(K'))\).

  To prove commutativity of the first diagram we rewrite the composition of the two arrows in the top right corner as the composition (here we leave restrictions implicit)
  \[ T_{t_2}^* \Fcal^{R(K')} R\Gamma(C(t_1K),V) \xrightarrow{s(Q(t_1K),Q(K'),t_2,Q(K''))} \Fcal^{R(K'')} R\Gamma(C(t_2t_1K),V) \xrightarrow{\res} \Fcal^{R(K'')} R\Gamma(C(K''),V) \]
  which again follows from general functoriality arguments for which we do not give details.
  This reduces commutativity of the first diagram to the cocycle relation for \(s(-,-,-,-)\) above.

  We proceed similarly for the second diagram: the top right corner sits in a commutative diagram
  \[
    \begin{tikzcd}
      T_{t_2}^* \Fcal^{R(K')} R\Gamma(C(t_1K),V) \ar[d, "{s(Q(t_1K),Q(K'),t_2,Q(K''))}"] & \ar[l, "{\cores}"] T_{t_2}^* \Fcal^{R(K')} R\Gamma(C(K'),V) \ar[d, "{\sim}" left, "{s(Q(K'),t_2,Q(K''))}"] \\
      \Fcal^{R(K'')} R\Gamma(C(t_2t_1K),V) & \ar[l, "{\cores}"] \Fcal^{R(K'') R\Gamma(C(t_2K'),V)}
    \end{tikzcd}
  \]
  (again, we omit the details) and this reduces commutativity of the second diagram to the cocycle relation for \(s(-,-,-,-)\).
\end{proof}

\begin{prodef} \label{prodef:ru_cu}
  Let \(h \in \Gbf(\A_f)\) and let \(K\) and \(K'\) be neat compact open subgroups of \(\Gbf(\A_f)\) satisfying \(K' \subset hKh^{-1}\).
  There is a unique morphism
  \[ T_h^* \AFcal^K(V) \xrightarrow{ru(h,K,K')} \AFcal^{K'}(V) \]
  such that for any \(g,g' \in \Gbf(\A_f)\) and \(t \in S_\Q \Gbf_{\her}(\A_f)\) satisfying \(g' h \in t g K\) we have a commutative diagram
  \begin{equation} \label{eq:diag_def_ru_Pink}
    \begin{tikzcd}
      T_t^* \iota_g^* \AFcal^K(V) \ar[r, "{\sim}"] \ar[d, dash, "{\sim}"] & \iota_{g'}^* T_h^* \AFcal^K(V) \ar[r, "{\iota_{g'}^* \left( ru(h,K,K') \right)}" above] & \iota_{g'}^* \AFcal^{K'}(V) \ar[d, dash, "{\sim}"] \\
      T_t^* \Fcal^{R(gK)} R\Gamma(C(gK), V) \ar[r, "{\sim}" above, "{s(Q(gK),t)}" below] & \Fcal^{R(g'K')} R\Gamma(C(tgK), V) \ar[r, "{\res}"] & \Fcal^{R(g'K')} R\Gamma(C(g'K'), V)
    \end{tikzcd}
  \end{equation}
  where the top left horizontal isomorphism comes from \eqref{eq:def_T_h_Pink}.

  We also have a unique morphism
  \[ \AFcal^{K'}(V) \xrightarrow{cu(K,h,K')} T_h^* \AFcal^K(V) \]
  characterized by commutative diagrams
  \begin{equation} \label{eq:diag_def_cu_Pink}
    \begin{tikzcd}
      T_t^* \iota_g^* \AFcal^K(V) \ar[r, "{\sim}"] \ar[d, "{\sim}"] & \iota_{g'}^* T_h^* \AFcal^K(V) & \ar[l, "{\iota_{g'}^* \left( cu(K,h,K') \right)}" above] \iota_{g'}^* \AFcal^{K'}(V) \ar[d, "{\sim}"] \\
      T_t^* \Fcal^{R(gK)} R\Gamma(C(gK), V) \ar[r, "{\sim}"] & \Fcal^{R(g'K')} R\Gamma(C(tgK), V) & \ar[l, "{\cores}" above] \Fcal^{R(g'K')} R\Gamma(C(g'K'), V)
    \end{tikzcd}
  \end{equation}
\end{prodef}
\begin{proof}
  Uniqueness again follows from \eqref{eq:gen_Shim_var_explicit}, existence follows from Lemma \ref{lem:compo_s_Pink} (in the case where one of the two inclusions between levels is an equality).
\end{proof}

\begin{prop} \label{pro:compo_ru_cu}
  For \(K' \subset h_1Kh_1^{-1}\) and \(K'' \subset h_2K'h_2^{-1}\) the following diagram is commutative
  \[
    \begin{tikzcd}[column sep=5em]
      T_{h_2}^* T_{h_1}^* \AFcal^K(V) \ar[r, "{T_{h_2}^* ru(h_1,K,K')}"] \ar[d, dash, "{\sim}"] & T_{h_2}^* \AFcal^{K'}(V) \ar[d, "{ru(h_2,K',K'')}"] \\
      T_{h_2h_1}^* \AFcal^K(V) \ar[r, "{ru(h_2h_1,K,K'')}"] & \AFcal^{K''}(V)
    \end{tikzcd}
  \]
  and similarly for the maps \(cu\) with reversed directions of arrows:
  \[ cu(K,h_2h_1,K'') = T_{h_2}^* cu(K,h_1,K') \circ cu(K', h_2, K''). \]
\end{prop}
\begin{proof}
  Of course this is proved by restricting to an arbitrary \(U(g''K'') \xhookrightarrow{\iota_{g''}} \Sh(\Gbf,\Xcal,K'')\).
  We do not even need Lemma \ref{lem:compo_s_Pink} now, as we may choose \(g'=g''h_2\) and \(g=g'h_1\) to form \(U(g'K') \xhookrightarrow{\iota_{g'}} \Sh(\Gbf,\Xcal,K')\) and \(U(gK) \xhookrightarrow{\iota_g} \Sh(\Gbf,\Xcal,K)\).
\end{proof}

We deduce formally \(ru(hk,K,K')=ru(h,K,K')\) and \(cu(K,kh,K')=cu(K,h,K')\) for \(k \in K\).

\begin{defi} \label{def:corr_au}
  Let \(h_1,h_2 \in \Gbf(\A_f)\), \(K_1\), \(K_2\) and \(K'\) neat compact open subgroups of \(\Gbf(\A_f)\) such that \(K' \subset h_i K_i h_i^{-1}\) for any \(i \in \{1,2\}\).
  We define a correspondence
  \[ au(K_2, h_2, h_1, K_1, K'): T_{h_1}^* \AFcal^{K_1} (V) \longrightarrow T_{h_2}^* \AFcal^{K_2}(V) \]
  supported on
  \[ \Sh(\Gbf, \Xcal, K_1) \xleftarrow{T_{h_1}} \Sh(\Gbf, \Xcal, K') \xrightarrow{T_{h_2}} \Sh(\Gbf, \Xcal, K_2) \]
  as the composition
  \[ T_{h_1}^* \AFcal^{K_1}(V) \xrightarrow{ru(h_1,K_1,K')} \AFcal^{K'}(V) \xrightarrow{cu(K_2,h_2,K')} T_{h_2}^* \AFcal^{K_2}(V). \]
\end{defi}

\begin{prop} \label{pro:easy_prop_corr_au}
  In the setting of Definition \ref{def:corr_au} we have
  \[ au(K_2, h_2 k_2, h_1 k_1, K_1, K') = au(K_2, h_2, h_1, K_1, K') \]
  for any \(k_i \in K_i\) and
  \[ au(K_2, x h_2, x h_1, K_1, xK'x^{-1}) = \corr (T_x)_* au(K_2, h_2, h_1, K_1, K') \]
  for any \(x \in \Gbf(\A_f)\) (note that \(T_x: \Sh(\Gbf,\Xcal,xK'x^{-1}) \to \Sh(\Gbf,\Xcal,K')\) is an isomorphism).
\end{prop}
\begin{proof}
  This follows from Proposition \ref{pro:compo_ru_cu} and the simple facts:
  \begin{itemize}
  \item \(ru(k_1,K_1,K_1) = \id\) for \(k_1 \in K_1\),
  \item \(cu(K_2, k_2, K_2) = \id\) for \(k_2 \in K_2\),
  \item \(ru(x,K',xK'x^{-1})\) and \(cu(K',x,xK'x^{-1})\) are isomorphisms which are inverse of each other.
  \end{itemize}
\end{proof}

The second part of Proposition \ref{pro:easy_prop_corr_au} explains why we may take \(h_1=1\) without loss of generality.

We have the following analogue of Proposition \ref{pro:Hecke_corr_sat_formalism}.

\begin{prop}
  The correspondences \(au(K_2, 1, h, K_1, K')\) induce a Hecke formalism in cohomology, i.e.\ the maps
  \[ (au(K_2, 1, h, K_1, K'))_? \]
  where \(? \in \{!,*\}\), satisfy the axioms in Definition \ref{def:hecke_cat}.
\end{prop}

This section \ref{sec:gen_Shim_Hecke} is devoted to the proof of this proposition.
The first axiom in Definition \ref{def:hecke_cat} follows from Proposition \ref{pro:easy_prop_corr_au}.
The second axiom is trivial.

The proof of the third axiom proceeds as in the proof of Proposition \ref{pro:Hecke_corr_sat_formalism}, with a slight complication.
Assuming \(K'' \subset K' \subset h_iK_ih_i^{-1}\), we compute using Proposition \ref{pro:compo_ru_cu}
\begin{align*}
  au(K_2,h_2,h_1,K_1,K'')
  &= cu(K_2,h_2,K'') \circ ru(h_1,K_1,K'') \\
  &= T_1^* cu(K_2,h_2,K') \circ cu(K',1,K'') \circ ru(1,K',K'') \circ T_1^* ru(h_1,K_1,K')
\end{align*}
where \(T_1: \Sh(\Gbf,\Xcal,K'') \to \Sh(\Gbf,\Xcal,K')\).
The composition
\[ cu(K',1,K'') \circ ru(1,K',K'') : T_1^* \AFcal^{K'} (V) \longrightarrow T_1^* \AFcal^{K'} (V) \]
is equal to multiplication by the locally constant function \(\delta(K',K''): \Sh(\Gbf,\Xcal,K'') \to \Z_{>0}\) such that \(\delta(K',K'') \circ \iota_{gK''} = |C(gK')/C(gK'')|\) for any \(gK'' \in \Gbf(\A_f)/K''\), because the composition of morphisms of functors
\[ \res^{R(gK')}_{R(gK'')} R\Gamma(C(gK'),-) \xrightarrow{\res} R\Gamma(C(gK''), \res^{Q(gK')}_{Q(gK'')} -) \xrightarrow{\cores} \res^{R(gK')}_{R(gK'')} R\Gamma(C(gK'),-) \]
is equal to multiplication by \(|C(gK')/C(gK'')|\) (this follows from \cite[Lemma 3.1]{Taibi_vanEst}).
Thus we have
\[ au(K_2,h_2,h_1,K_1,K'') = \delta(K',K'') \corr (T_1)^* au(K_2,h_2,h_1,K_1,K'). \]
The map \(T_1: U(gK'') \to U(gK')\) is finite étale of constant degree
\[ \frac{|P(gK')/P(gK'')| \times |R(gK')/P(gK')|}{|R(gK'')/P(gK'')|} = \frac{|Q(gK')/Q(gK'')|}{|C(gK')/C(gK'')|}. \]
For \([g'] \in S_\Q \Gbf_{\her}(\A_f) \backslash \Gbf(\A_f) / K'\) we have
\[ \sum_{[g] \in S_\Q \Gbf_{\her}(\A_f) \backslash \Gbf(\A_f) / K'' \mapsto [g']} |Q(gK')/Q(gK'')| = K'/K'' \]
and together with a variant of Lemma \ref{lemm:fet_push_pull_corr} this allows us to conclude
\[ \corr (T_1)_* au(K_2,h_2,h_1,K_1,K'') = |K'/K''| au(K_2,h_2,h_1,K_1,K'). \]

We are left to check the last axiom in Definition \ref{def:hecke_cat}.
We want to compute the composition \(au(K_3,h_3,h_2',K_2,K'') \circ au(K_2,h_2,h_1,K_1,K')\), which in the following diagram is supported on \((T_{h_1} \circ \pr_1, T_{h_3} \circ \pr_2)\).
\[
  \begin{tikzcd}[column sep=1em]
    & & \Sh(\Gbf,\Xcal,K') \underset{\Sh(\Gbf, \Xcal, K_2)}{\times} \Sh(\Gbf, \Xcal, K'') \ar[dl, "{\pr_1}"] \ar[dr, "{\pr_2}"] & & \\
    & \Sh(\Gbf, \Xcal, K') \ar[dl, "{T_{h_1}}"] \ar[dr, "{T_{h_2}}"] & & \Sh(\Gbf, \Xcal, K'') \ar[dl, "{T_{h_2'}}"] \ar[dr, "{T_{h_3}}"] & \\
    \Sh(\Gbf, \Xcal, K_1) & & \Sh(\Gbf, \Xcal, K_2) & & \Sh(\Gbf, \Xcal, K_3)
  \end{tikzcd}
\]
As in the classical case we consider the colimit of the functor
\begin{align*}
  F_{K_2,K',h_2,K'',h_2'}: [K'' \curvearrowright h_2' K_2 h_2^{-1} / K'] & \longrightarrow \mathrm{Sch} \\
  hK' & \longmapsto \Sh(\Gbf, \Xcal, K'' \cap hK'h^{-1}) \\
  \left( hK' \xrightarrow{x \in K''} xhK' \right) & \longmapsto \left( \Sh(\Gbf, \Xcal, K'' \cap hK'h^{-1}) \xrightarrow{T_{x^{-1}}} \Sh(\Gbf, \Xcal, K'' \cap xhK'h^{-1}) \right)
\end{align*}
which again may simply be expressed as a disjoint union over \(K'' \backslash h_2'K_2h_2^{-1} / K'\) because this functor maps automorphisms to identity morphisms.
For an object \(hK'\) of \([K'' \curvearrowright h_2' K_2 h_2^{-1} / K']\) we have a map
\[ f_{hK'} := T_h \times T_1: \Sh(\Gbf, \Xcal, K'' \cap hK'h^{-1}) \longrightarrow \Sh(\Gbf,\Xcal,K') \underset{\Sh(\Gbf, \Xcal, K_2)}{\times} \Sh(\Gbf, \Xcal, K'') \]
and these maps induce a well-defined map
\begin{equation} \label{eq:f_from_colim_to_fibprod_Sh}
  f: \colim F_{K_2,K',h_2,K'',h_2'} \longrightarrow \Sh(\Gbf,\Xcal,K') \underset{\Sh(\Gbf, \Xcal, K_2)}{\times} \Sh(\Gbf, \Xcal, K'').
\end{equation}
Thanks to Proposition \ref{pro:easy_prop_corr_au} there is a unique correspondence
\[ v(K_3, h_3, K'', h_2', K_2, h_2, K', h_1, K_1) \]
supported on \((T_{h_1} \circ \pr_1 \circ f, T_{h_3} \circ \pr_2 \circ f)\) from \(\AFcal^{K_1}(V)\) to \(\AFcal^{K_3}(V)\), which is equal to \(au(K_3, h_3, hh_1, K_1, K'' \cap hK'h^{-1})\) on \(\Sh(\Gbf, \Xcal, K'' \cap hK'h^{-1})\) for any \(hK'\).
The last axiom of Definition \ref{def:hecke_cat} follows from the next lemma (again using Section \ref{sec:corr_push_pull}).

\begin{lemm} \label{lem:compo_corr_au}
  The composition of correspondences
  \[ au(K_3,h_3,h_2',K_2,K'') \circ au(K_2,h_2,h_1,K_1,K') \]
  is equal to
  \[ (\corr f_*)v(K_3, h_3, K'', h_2', K_2, h_2, K', h_1, K_1). \]
\end{lemm}

We devote the rest of Section \ref{sec:gen_Shim_Hecke} to the proof of this lemma.
We will proceed in three steps.
This first step is a rather formal reduction to a similar but simpler statement involving fewer objects.
In the second step we exhibit representatives (explicit morphisms between explicit complexes) for our two morphisms (originally in the localization at \(\ell\) of Ekedahl \(\ell\)-adic categories).
In the third step we check that these representatives are equal (as morphisms of complexes, not just in localized categories) by checking that they agree on all stalks.

\subsubsection{Proof of Lemma \ref{lem:compo_corr_au}: step 1}

The composition of correspondences decomposes as
\begin{align*}
  &\ \pr_1^* T_{h_1}^* \AFcal^{K_1}(V) \\
  \xrightarrow{\mathmakebox[3cm]{\pr_1^*(ru(h_1,K_1,K'))}} &\ \pr_1^* \AFcal^{K'}(V) \\
  \xrightarrow{\mathmakebox[3cm]{\pr_1^*(cu(K_2,h_2,K'))}} &\ \pr_1^* T_{h_2}^* \AFcal^{K_2}(V) \\
  \xrightarrow{\mathmakebox[3cm]{\sim}} &\ \pr_2^* T_{h_2'}^* \AFcal^{K_2}(V) \\
  \xrightarrow{\mathmakebox[3cm]{\pr_2^*(ru(h_2',K_2,K''))}} &\ \pr_2^* \AFcal^{K''}(V) \\
  \xrightarrow{\mathmakebox[3cm]{\pr_2^*(cu(K_3,h_3,K''))}} &\ \pr_2^* T_{h_3}^* \AFcal^{K_3}(V).
\end{align*}

The pushforward \((\corr f_*)v(K_3, h_3, K'', h_2', K_2, h_2, K', h_1, K_1)\) is the sum over \([h] \in K'' \backslash h_2' K_2 h_2^{-1} / K'\) of
\begin{multline*}
  \pr_1^* T_{h_1}^* \AFcal^{K_1}(V) \xrightarrow{\adj} (f_{hK'})_* f_{hK'}^* \pr_1^* T_{h_1}^* \AFcal^{K_1}(V) \xrightarrow{\sim} (f_{hK'})_* T_{hh_1}^* \AFcal^{K_1}(V) \\
  \xrightarrow{(f_{hK'})_*(ru(hh_1,K_1, K'' \cap hK'h^{-1}))} (f_{hK'})_* \AFcal^{K'' \cap hK'h^{-1}}(V) \\
  \xrightarrow{(f_{hK'})_*(cu(K_3,h_3, K'' \cap hK'h^{-1}))} (f_{hK'})_* T_{h_3}^* \AFcal^{K_3}(V) \\
  \xrightarrow{\sim} (f_{hK'})_* f_{hK'}^* \pr_2^* T_{h_3}^* \AFcal^{K_3}(V) \xrightarrow{\Tr} \pr_2^* T_{h_3}^* \AFcal^{K_3}(V)
\end{multline*}
Using the identities (Proposition \ref{pro:compo_ru_cu})
\[ ru(hh_1,K_1, K'' \cap hK'h^{-1}) = ru(h,K', K'' \cap hK'h^{-1}) \circ T_h^*(ru(h_1,K_1,K')) \]
\[ cu(K_3,h_3, K'' \cap hK'h^{-1}) = T_{h_3}^*(cu(K_3,h_3,K'')) \circ cu(K'',1, K'' \cap hK'h^{-1}) \]
we decompose \((\corr f_*)v(K_3, h_3, K'', h_2', K_2, h_2, K', h_1, K_1)\) as the sum over \([h]\) of
\begin{multline*}
  \pr_1^* T_{h_1}^* \AFcal^{K_1}(V) \xrightarrow{\pr_1^*(ru(h_1,K_1,K'))} \pr_1^* \AFcal^{K'}(V) \xrightarrow{\adj} (f_{hK'})_* f_{hK'}^* \pr_1^* \AFcal^{K'}(V) \\
  \xrightarrow{\sim} (f_{hK'})_* T_h^* \AFcal^{K'}(V) \xrightarrow{(f_{hK'})_*(ru(h,K', K'' \cap hK'h^{-1}))} (f_{hK'})_* \AFcal^{K'' \cap hK'h^{-1}}(V) \\
  \xrightarrow{(f_{hK'})_*(cu(K'',1, K'' \cap hK'h^{-1}))} (f_{hK'})_* T_1^* \AFcal^{K''}(V) \\
  \xrightarrow{\sim} (f_{hK'})_* f_{hK'}^* \pr_2^* \AFcal^{K''}(V) \xrightarrow{\Tr}  \pr_2^* \AFcal^{K''}(V) \xrightarrow{\pr_2^*(cu(K_3,h_3, K''))} \pr_2^* T_{h_3}^* \AFcal^{K_3}(V).
\end{multline*}
Now the first and last maps are the same so we are left to check that the correspondence
\begin{equation} \label{eq:comp_au_corr1}
  \pr_1^* \AFcal^{K'}(V) \xrightarrow{\pr_1^*(cu(K_2,h_2,K'))} \pr_1^* T_{h_2}^* \AFcal^{K_2}(V) \simeq \pr_2^* T_{h_2'}^* \AFcal^{K_2}(V) \xrightarrow{\pr_2^*(ru(h_2',K_2,K''))} \pr_2^* \AFcal^{K''}(V)
\end{equation}
supported on \((\pr_1,\pr_2)\) is equal to the pushforward along \(f\) of the correspondence
\begin{equation} \label{eq:comp_au_corr2}
  w(K'',K_2,K') := \underset{[h] \in K'' \backslash h_2' K_2 h_2 / K'}{\boxplus} au(K'',1,h,K',K'' \cap hK'h^{-1}).
\end{equation}
We may also reduce to the case where \(h_2=h_2'=1\) by replacing \(K'\) and \(K''\) by conjugate subgroups of \(\Gbf(\A_f)\).

\subsubsection{Proof of Lemma \ref{lem:compo_corr_au}: step 2}

It suffices to prove that for any \(g_2 \in \Gbf(\A_f)\) our two morphisms (\eqref{eq:comp_au_corr1} and the pushforward of \eqref{eq:comp_au_corr2}, for \(h_2=h_2'=1\)) coincide over the preimage of the component \(\iota_{g_2}(U(g_2 K_2))\) of \(\Sh(\Gbf, \Xcal, K_2)\).
We fix \(g_2\) for the rest of this section.
Recall that \(V\) is (represented by) a bounded complex of algebraic representations of \(\Gbf_E\), which we consider as continuous representations of \(\Gbf(\Qell)\) over finite-dimensional vector spaces over \(E\).
We fix a subcomplex \(\Lambda\) consisting of \(\Ocal_E\)-lattices stable under the action of \(g_2 K_2 g_2^{-1}\).
To such a complex we may associate a complex of \((\Ocal_E)_\bullet\)-modules in the topos \(S_{Q(g_2K_2)}^{\N}\), and we fix a quasi-isomorphic bounded below complex \(I^\bullet\) of injective \((\Ocal_E)_\bullet\)-modules in \(S_{Q(g_2K_2)}^{\N}\).

We will need an integral version of the morphisms \(s(-,-,-,-)\) defined in the proof of Lemma \ref{lem:compo_s_Pink}.
Assume that \(L \subset \tilde{L}\) are open subgroups of \(g_2K_2g_2^{-1}\), \(t \in Q(g_2K_2)\) and \(L'\) an open subgroup of \(tLt^{-1}\).
Now \(t\) induces an isomorphism \(\Ad(t^{-1})^* \res_{Q(\tilde{L})} \Lambda \simeq \res_{Q(t\tilde{L})} \Lambda\) and similarly for \(I^\bullet\), and so we have an isomorphism of complexes of \((\Ocal_E)_\bullet\)-modules defined as in the proof of Lemma \ref{lem:compo_s_Pink}:
\begin{equation} \label{eq:T_t_pullback_Pink_gal_integral}
  T_t^* \Fcal^{R(L)} (I^\bullet)^{C(\tilde{L})} \xrightarrow{s_{g_2K_2}(Q(\tilde{L}),Q(L),t,Q(L'))} \Fcal^{R(L')} (I^\bullet)^{C(t\tilde{L})}
\end{equation}
Note that \((I^\bullet)^{C(\tilde{L})}\) is naturally a complex of \((\Ocal_E)_\bullet\)-modules in \(S_{R(\tilde{L})}^{\N}\) and is considered in \(S_{R(L)}^{\N}\) by restriction, and similarly for \((I^\bullet)^{C(t\tilde{L})}\).
We have a cocycle relation (same proof as for Lemma \ref{lem:compo_s_Pink}): for \(t' \in Q(g_2K_2)\) and \(L'' \subset t'L'(t')^{-1}\) we have
\[ s_{g_2K_2}(Q(t\tilde{L}),Q(L'),t',Q(L'')) \circ T_{t'}^* s_{g_2K_2}(Q(\tilde{L}),Q(L),t,Q(L')) = s_{g_2K_2}(Q(\tilde{L}),Q(L),t't,Q(L'')). \]
By the same argument as in Proposition-Definition \ref{prodef:AFcal} replacing \(t\) by another element of \(tQ(L)\) leaves \(s_{g_2K_2}(Q(\tilde{L}),Q(L),t,Q(L'))\) unchanged.
Now if \(L\) is an open subgroup of \(K_2\), the full subgroupoid of \([S_\Q \Gbf_{\her}(\A_f) \curvearrowright \Gbf(\A_f)/L]\) (used in Definition \ref{def:gen_Shim_var}) whose objects are \(gL\) for \(g \in g_2K_2\) is equivalent to \([Q(g_2K_2) \curvearrowright g_2K_2/L]\), and the preimage of \(\iota_{g_2}(U(g_2 K_2))\) (in other words, the fiber product of \(\iota_{g_2}: U(g_2K_2) \hookrightarrow \Sh(\Gbf,\Xcal,K_2)\) and \(T_1: \Sh(\Gbf,\Xcal,L) \to \Sh(\Gbf,\Xcal,K_2)\)) in \(\Sh(\Gbf, \Xcal, L)\) may be identified with
\[ \underset{gL \in [Q(g_2K_2) \curvearrowright g_2K_2/L]}{\colim} U(g L). \]
Similarly to Proposition-Definition \ref{prodef:AFcal} we have a well-defined complex \(\AFcal^L_{g_2 K_2} I^\bullet\) of \((\Ocal_E)_\bullet\)-modules on this colimit\footnote{More accurately, in \(S^\N\) where \(S\) is the étale topos of this colimit.} together with isomorphisms
\[ \iota_g^* \AFcal^L_{g_2 K_2} I^\bullet \simeq \Fcal^{R(gL)} (I^\bullet)^{C(gL)}. \]
For any \(h \in K_2\) and \(L'\) an open subgroup of \(hLh^{-1}\) we have well-defined morphisms of complexes of \((\Ocal_E)_\bullet\)-modules on \(\Sh(\Gbf,\Xcal,L') \times_{\Sh(\Gbf,\Xcal,K_2)} U(g_2K_2)\)
\[ T_h^* \AFcal^L_{g_2 K_2} (I^\bullet) \xrightarrow{ru_{g_2 K_2}(h, L, L')} \AFcal^{L'}_{g_2 K_2} (I^\bullet) \]
and
\[ \AFcal^{L'}_{g_2 K_2} (I^\bullet) \xrightarrow{cu_{g_2 K_2}(L, h, L')}  T_h^* \AFcal^L_{g_2 K_2} (I^\bullet) \]
fitting in commutative diagrams (now of complexes of \((\Ocal_E)_\bullet\)-modules) similar to \eqref{eq:diag_def_ru_Pink} and \eqref{eq:diag_def_cu_Pink} with \(I^\bullet\) replacing \(V\), \((-)^H\) replacing \(R\Gamma(H, -)\) and the obvious inclusion (resp.\ the norm map) replacing \(\res\) (resp.\ \(\cores\)).
The proof is the same as for Proposition-Definition \ref{prodef:ru_cu}, and the same proof as that of Proposition \ref{pro:compo_ru_cu} shows that \(ru_{g_2K_2}\) and \(cu_{g_2K_2}\) satisfy cocycle relations.

The integral analogue of the morphism \eqref{eq:comp_au_corr1} is the morphism \(A_{g_2K_2}\) of complexes of \((\Ocal_E)_\bullet\)-modules defined as the composition
\begin{equation} \label{diag:def_A_g2_K2}
  \begin{tikzcd}[column sep=huge]
    \pr_1^* \AFcal^{K'}_{g_2 K_2}(I^\bullet) \ar[r, "{A_{g_2 K_2}}"] \ar[d, "{\pr_1^*(cu_{g_2K_2}(K_2,1,K'))}" left] & \pr_2^* \AFcal^{K''}_{g_2 K_2}(I^\bullet) \\
    \pr_1^* T_1^* \AFcal^{K_2}_{g_2 K_2}(I^\bullet) \ar[r, dash, "{\sim}"] & \pr_2^* T_1^* \AFcal^{K_2}_{g_2 K_2}(I^\bullet) \ar[u, "{\pr_2^*(ru_{g_2K_2}(1,K_2,K''))}" right]
  \end{tikzcd}
\end{equation}
More precisely for any \(g',g'' \in g_2 K_2\), denoting
\[ U(g'K') \underset{U(g_2K_2)}{\times} U(g''K'') \xhookrightarrow{\iota_{g',g''} := \iota_{g'} \times \iota_{g''}} \Sh(\Gbf, \Xcal, K') \underset{\Sh(\Gbf, \Xcal, K_2)}{\times} \Sh(\Gbf, \Xcal, K'') \]
we have a commutative diagram of complexes of \((\Ocal_E)_\bullet\)-modules
\begin{equation} \label{diag:char_A_g2_K2}
  \begin{tikzcd}[column sep=8em]
    \iota_{g',g''}^* \pr_1^* \AFcal^{K'}_{g_2 K_2}(I^\bullet) \ar[r, "{\iota_{g',g''}^*(A_{g_2 K_2})}"] \ar[d, dash, "{\sim}"] & \iota_{g',g''}^* \pr_2^* \AFcal^{K''}_{g_2 K_2}(I^\bullet) \ar[d, dash, "{\sim}"] \\
    \pr_1^* \Fcal^{R(g'K')} (I^\bullet)^{C(g'K')} \ar[d, "{\pr_1^* \Fcal^{R(g'K')}(N)}"] \ar[r, "{A_{g_2K_2,g',g''}}"] & \pr_2^* \Fcal^{R(g''K'')} (I^\bullet)^{C(g''K'')} \\
    \pr_1^* \Fcal^{R(g'K')} (I^\bullet)^{C(g_2K_2)} & \pr_2^* \Fcal^{R(g''K'')} (I^\bullet)^{C(g_2K_2)} \ar[u, hook] \\
    \pr_1^* T_1^* \Fcal^{R(g_2K_2)} (I^\bullet)^{C(g_2K_2)} \ar[r, dash, "{\sim}"] \ar[u, "{\sim}" left, "{\pr_1^* \left( s_{g_2K_2}(Q(g_2K_2),Q(g_2K_2),1,Q(g'K')) \right)}" right] & \pr_2^* T_1^* \Fcal^{R(g_2K_2)} (I^\bullet)^{C(g_2K_2)} \ar[u, "{\sim}" left, "{\pr_2^* \left( s_{g_2K_2}(Q(g_2K_2),Q(g_2K_2), 1, Q(g''K'')) \right)}" right]
  \end{tikzcd}
\end{equation}
where \(N\) denotes the norm map for \(C(g'K') \subset C(g_2K_2)\) and \(A_{g_2K_2,g',g''}\) is defined here by commutativity (this notation will be convenient later).
By design the morphism \(A_{g_2K_2}\) is compatible with \eqref{eq:comp_au_corr1}: the following diagram in
\[ D^+ \left( \Sh(\Gbf,\Xcal,K') \underset{\Sh(\Gbf, \Xcal, K_2)}{\times} \Sh(\Gbf, \Xcal, K''), E \right) \]
is commutative, denoting by \(i_{g_2 K_2}\) the base change along the open immersion (with closed image) \(\iota_{g_2}: U(g_2 K_2) \hookrightarrow \Sh(\Gbf, \Xcal, K_2)\) of \(\Sh(\Gbf,\Xcal,K') \underset{\Sh(\Gbf, \Xcal, K_2)}{\times} \Sh(\Gbf, \Xcal, K'') \to \Sh(\Gbf, \Xcal, K_2)\).
\[
  \begin{tikzcd}[column sep=huge]
    i_{g_2 K_2}^* \pr_1^* \AFcal^{K'}(V) \ar[d, dash, "{\sim}"] \ar[r, "{i_{g_2 K_2}^* \eqref{eq:comp_au_corr1}}"] & i_{g_2 K_2}^* \pr_2^* \AFcal^{K''}(V) \ar[d, dash, "{\sim}"] \\
    \pr_1^* \AFcal^{K'}_{g_2 K_2}(I^\bullet) \ar[r, "{A_{g_2 K_2}}"] & \pr_2^* \AFcal^{K''}_{g_2 K_2}(I^\bullet)
  \end{tikzcd}
\]

We proceed similarly for our second morphism \((\corr f_*) w(K'', K_2, K')\) (defined in \eqref{eq:comp_au_corr2}).
We have an analogue \(F_{K_2,K',K''}'\) of the functor \(F_{K_2,K',1,K'',1}\) used to define \(f\) (see \ref{eq:f_from_colim_to_fibprod_Sh}), mapping \(hK' \in [K'' \curvearrowright K_2/K']\) to
\[ \underset{gK(h) \in [Q(g_2K_2) \curvearrowright g_2K_2/K(h)]}{\colim} U(gK(h)) \]
where we have abbreviated \(K(h) := K'' \cap hK'h^{-1}\), and a morphism
\[ \colim F_{K_2,K',K''}' \xrightarrow{f_{g_2K_2}} \underset{g'K' \in [Q(g_2K_2) \curvearrowright g_2K_2/K']}{\colim} U(g'K') \ \underset{U(g_2K_2)}{\times} \ \underset{g''K'' \in [Q(g_2K_2) \curvearrowright g_2K_2/K'']}{\colim} U(g''K'') \]
induced by the maps \(T_1 \times T_1: U(gK(h)) \to U(ghK') \times_{U(g_2K_2)} U(gK'')\).
Note that \(f_{g_2K_2}\) is just the restriction of \(f\) to the preimage of \(\iota_{g_2}(U(g_2K_2))\).
We have a morphism of complexes of \((\Ocal_E)_\bullet\)-modules
\[ f_{g_2K_2}^* \pr_1^* \AFcal^{K'}_{g_2 K_2} I^\bullet \simeq \underset{[h] \in K'' \backslash K_2 / K'}{\boxplus} T_h^* \AFcal^{K'}_{g_2 K_2} I^\bullet \xrightarrow{(ru_{g_2K_2}(h,K',K(h)))_{[h]}} \underset{[h] \in K'' \backslash K_2 / K'}{\boxplus} \AFcal_{g_2 K_2}^{K(h)} I^\bullet \]
which is well-defined thanks to the integral analogue of Proposition \ref{pro:compo_ru_cu}.
Similarly we have a well-defined morphism
\[ \underset{[h] \in K'' \backslash K_2 / K'}{\boxplus} \AFcal_{g_2 K_2}^{K(h)} I^\bullet  \xrightarrow{(cu_{g_2K_2}(K'',1,K(h)))_{[h]}} \underset{[h] \in K'' \backslash K_2 / K'}{\boxplus} T_1^* \AFcal^{K''}_{g_2 K_2} I^\bullet \simeq f^* \pr_2^* \AFcal^{K''}_{g_2 K_2} I^\bullet. \]
Define a morphism \(B_{g_2K_2}\) of \((\Ocal_E)_\bullet\)-modules as the composition
\begin{equation} \label{diag:def_B_g2_K2}
  \begin{tikzcd}
    \pr_1^* \AFcal^{K'}_{g_2 K_2} I^\bullet \ar[rr, "{B_{g_2 K_2}}"] \ar[d, "{\adj}"] & & \pr_2^* \AFcal^{K''}_{g_2 K_2} I^\bullet \\
    (f_{g_2K_2})_* f_{g_2K_2}^* \pr_1^* \AFcal^{K'}_{g_2 K_2} I^\bullet \ar[dr, "{(f_{g_2K_2})_*((ru_{g_2 K_2}(h, K', K(h)))_{[hK']})}" left, near end] & & (f_{g_2K_2})_* f_{g_2K_2}^* \pr_2^* \AFcal^{K''}_{g_2 K_2} I^\bullet \ar[u, "{\Tr}"] \\
    & (f_{g_2K_2})_* \left( \underset{[hK'] \in K'' \backslash K_2 / K'}{\boxplus} \AFcal^{K(h)}_{g_2 K_2} I^\bullet \right) \ar[ur, "{(f_{g_2K_2})_*((cu_{g_2 K_2}(K'', 1, K(h)))_{[hK']})}" right, near start] &
  \end{tikzcd}
\end{equation}
Again by design we have a commutative diagram in
\[ D^+ \left( \Sh(\Gbf,\Xcal,K') \underset{\Sh(\Gbf, \Xcal, K_2)}{\times} \Sh(\Gbf, \Xcal, K''), E \right): \]
\[
  \begin{tikzcd}[column sep=10em]
    i_{g_2 K_2}^* \pr_1^* \AFcal^{K'}(V) \ar[d, dash, "{\sim}"] \ar[r, "{i_{g_2 K_2}^* \left( (\corr f_*) w(K'', K_2, K') \right)}"] & i_{g_2 K_2}^* \pr_2^* \AFcal^{K''}(V) \ar[d, dash, "{\sim}"] \\
    \pr_1^* \AFcal^{K'}_{g_2 K_2} I^\bullet \ar[r, "{B_{g_2 K_2}}"] & \pr_2^* \AFcal^{K''}_{g_2 K_2} I^\bullet
  \end{tikzcd}
\]
In order to compute \(B_{g_2 K_2}\) more explicitly we need some preparation.
First we observe that \(\colim F_{K_2,K',1,K'',1}\) (the source of \(f\)), being a colimit (disjoint union) of colimits (also disjoint unions), may be expressed as a single colimit, over the groupoid
\[ [K'' \times S_\Q \Gbf_{\her}(\A_f) \curvearrowright K_2 / K' \times \Gbf(\A_f)] \]
for the (left) action \((x,s) \cdot (hK', g) = (xhK', sgx^{-1})\), of the functor
\begin{align*}
  (hK',g) & \longmapsto U(g K(h)) \\
  \left( (hK', g) \xrightarrow{(x,s)} (xhK', sgx^{-1}) \right) & \longmapsto \left( U(g K(h)) \xrightarrow{T_{s^{-1}}} U(sgx^{-1} K(xh)) \right).
\end{align*}
We denote
\[ \iota_{hK', g}: U(g K(h)) \hookrightarrow \colim F_{K_2,K',1,K'',1}. \]
Restricting over our fixed component \(\iota_{g_2}(U(g_2 K_2))\) of \(\Sh(\Gbf, \Xcal, K_2)\), we restrict to the full subgroupoid of objects mapping to \([g_2] \in S_\Q \Gbf_{\her}(\A_f) \backslash \Gbf(\A_f) / K_2\), which is easily checked to be equivalent to
\[ [ K'' \times Q(g_2K_2) \curvearrowright K_2 / K' \times g_2K_2 ] .\]
For \(g', g'' \in g_2 K_2\) we consider the full subgroupoid \(\Gcal_{g',g''}\) of \((hK', g)\) satisfying \(gh \in Q(g_2K_2) g' K'\) and \(g \in Q(g_2 K_2) g'' K''\), say \(gh \in t' g' K'\) and \(g \in t'' g'' K''\).
We have a Cartesian diagram
\begin{equation} \label{diag:cart_prod_Sh_U1}
  \begin{tikzcd}
    \underset{(hK',g) \in \Gcal_{g',g''}}{\colim} U(g K(h)) \ar[d, "{f_{g',g''}}"] \ar[r, hook] & \colim F_{K_2,K',1,K'',1} \ar[d, "{f}"] \\
    U(g'K') \underset{U(g_2K_2)}{\times} U(g''K'') \ar[r, hook, "{\iota_{g', g''}}"] & \Sh(\Gbf, \Xcal, K') \underset{\Sh(\Gbf, \Xcal, K_2)}{\times} \Sh(\Gbf, \Xcal, K'')
  \end{tikzcd}
\end{equation}
where the top open immersion is induced by the maps \(\iota_{hK',g}\) and \(f_{g',g''}\) is induced by the maps
\[ U(g' K') \xleftarrow{T_{t'}} U(g K(h)) \xrightarrow{T_{t''}} U(g''K''). \]
Now we have a commutative diagram
\begin{equation} \label{diag:char_B_g2_K2}
  \begin{tikzcd}
    \iota_{g',g''}^* \pr_1^* \AFcal^{K'}_{g_2 K_2}(I^\bullet) \ar[r, "{\iota_{g',g''}^*(B_{g_2 K_2})}"] \ar[d, dash, "{\sim}"] & \iota_{g',g''}^* \pr_2^* \AFcal^{K''}_{g_2 K_2}(I^\bullet) \ar[d, dash, "{\sim}"] \\
    \pr_1^* \Fcal^{R(g'K')} (I^\bullet)^{C(g'K')} \ar[d, "{\adj}"] \ar[r, "{B_{g_2K_2,g',g''}}"] & \pr_2^* \Fcal^{R(g''K'')} (I^\bullet)^{C(g''K'')} \\
    (f_{g',g''})_* f_{g',g''}^* \pr_1^* \Fcal^{R(g'K')} (I^\bullet)^{C(g'K')} \ar[d, dash, "{\sim}"] & (f_{g',g''})_* f_{g',g''}^* \pr_2^* \Fcal^{R(g''K'')} (I^\bullet)^{C(g''K'')} \ar[u, "{\Tr}"] \ar[d, dash, "{\sim}"] \\
    (f_{g',g''})_* \left( \underset{[hK',g]}{\boxplus} T_{t'}^* \Fcal^{R(g'K')} (I^\bullet)^{C(g'K')} \right) \ar[d, "{\left( s(Q(g'K'),Q(g'K'),t',Q(gK(h))) \right)_{hK',g}}" right, "{\sim}" left] & (f_{g',g''})_* \left( \underset{[hK',g]}{\boxplus} T_{t''}^* \Fcal^{R(g''K'')} (I^\bullet)^{C(g''K'')} \right) \ar[d, "{\left( s(Q(g''K''),Q(g''K''),t'',Q(gK(h))) \right)_{hK',g}}" right, "{\sim}" left] \\
    (f_{g',g''})_* \left( \underset{[hK',g]}{\boxplus} \Fcal^{R(gK(h))} (I^\bullet)^{C(ghK'h^{-1})} \right) \ar[d, hook] & (f_{g',g''})_* \left( \underset{[hK',g]}{\boxplus} \Fcal^{R(gK(h))} (I^\bullet)^{C(gK'')} \right) \\
    (f_{g',g''})_* \left( \underset{[hK',g]}{\boxplus} \Fcal^{R(gK(h))} (I^\bullet)^{C(gK(h))} \right) \ar[ur, "{\left( N_{C(gK'')/C(gK(h))} \right)_{[hK',g]}}" below] &
  \end{tikzcd}
\end{equation}
where as above \(B_{g_2K_2,g',g''}\) is defined for later use by requiring that the top square be commutative.

\subsubsection{Proof of Lemma \ref{lem:compo_corr_au}: step 3}

Now we will check that the morphisms of complexes \(A_{g_2 K_2}\) and \(B_{g_2 K_2}\) defined in \eqref{diag:def_A_g2_K2} and \eqref{diag:def_B_g2_K2} are equal by comparing them on stalks, using their more concrete descriptions \eqref{diag:char_A_g2_K2} and \eqref{diag:char_B_g2_K2}.
We may fix a geometric point \(p\) of \(U(g_2K_2)\), and only consider lifts of the point \(\iota_{g_2} \circ p\) to the various spaces above \(\Sh(\Gbf, \Xcal, K_2)\).
We fix a compatible family \(\ul{p} = (p_L)_L\) where the family runs over open subgroups \(L\) of \(R(g_2K_2)\) and \(p_L\) is a lift of \(p\) via \(T_1: \Sh(\Gbf_{\her}, \Xcal_1, L) \to U(g_2 K_2)\).
Note that the set of such compatible families is a torsor under \(R(g_2K_2)\).
For any such family \(\ul{p}\) we have an isomorphism of functors \(j_{\ul{p}}: p^* \Fcal^{R(g_2K_2)} \simeq \id\), and more generally \(p_L^* \Fcal^L \simeq \id\) for any open subgroup \(L\) of \(R(g_2K_2)\), abusively also denoted by \(j_{\ul{p}}\).

Looking at the diagrams \eqref{diag:char_A_g2_K2} and \eqref{diag:char_B_g2_K2} we see that it would be useful to compute the morphisms \eqref{eq:T_t_pullback_Pink_gal_integral} on stalks.
Unwinding the definitions we find that for open subgroups \(L \subset \tilde{L}\) of \(g_2K_2g_2^{-1}\), \(t \in Q(g_2K_2)\) and \(L'\) an open subgroup of \(tLt^{-1}\) we have a commutative diagram
\begin{equation} \label{eq:s_on_stalks}
  \begin{tikzcd}[column sep=14em]
    p_{R(L')}^* T_t^* \Fcal^{R(L)}(I^\bullet)^{C(\tilde{L})} \ar[r, "{\sim}" below, "{p_{R(L')}^* \left( s_{g_2K_2}(Q(\tilde{L}),Q(L), t, Q(L')) \right)}" above] \ar[d, "{j_{T_t(\ul{p})}}" right, "{\sim}" left] & p_{R(L')}^* \Fcal^{R(L')}(I^\bullet)^{C(t\tilde{L})} \ar[d, "{j_{\ul{p}}}" right, "{\sim}" left] \\
    (I^\bullet)^{C(\tilde{L})} \ar[r, "{t}" above, "{\sim}" below] & (I^\bullet)^{C(t\tilde{L})}
  \end{tikzcd}
\end{equation}
Note that, unlike the map \(s(\dots)\), the family \(T_t(\ul{p})\) of geometric points depends on the choice of \(t\) in \(t Q(L)\), thus so does \(j_{T_t(\ul{p})}\).
The fact that we are using two different families of points to compute stalks may seem problematic now, but ultimately we will be able to take \(t \in C(g_2K_2)\), implying \(T_t(\ul{p}) = \ul{p}\).

In order to describe the various lifts of \(\iota_{g_2} p\) we also fix \(g_0 \in g_2 K_2\), so that for any open subgroup \(L\) of \(K_2\) we have a geometric point \(\iota_{g_0} \circ p_{R(g_0L)}\) of \(\Sh(\Gbf, \Xcal, L)\) above \(\iota_{g_2} p\).
Denote \(\Gamma = g_0^{-1} C(g_0K_2) g_0 = K_2 \cap g_0^{-1} C_{\Q} g_0\).
Let \(S = \Gamma \backslash K_2\) considered as a set with right action of \(K_2\).
Denote \(s_0 = \Gamma \in S\).
Thanks to Lemma \ref{lem:Pink_Shim_as_quot} we have an identification of the base change to \(\iota_{g_2} p\) of the diagram
\[
  \begin{tikzcd}[row sep=large]
    & \underset{hK' \in [K'' \curvearrowright K_2 / K']}{\colim} \, \Sh(\Gbf, \Xcal, K(h)) \ar[dl, "{(T_h)_{hK'}}" above left] \ar[dr, "{(T_1)_{hK'}}" above right] \ar[d, "{f}"] & \\
    \Sh(\Gbf, \Xcal, K') \ar[dr, "{T_1}" below left] & \Sh(\Gbf, \Xcal, K') \underset{\Sh(\Gbf, \Xcal, K_2)}{\times} \Sh(\Gbf, \Xcal, K'') \ar[l] \ar[r] & \Sh(\Gbf, \Xcal, K'') \ar[dl, "{T_1}" below right] \\
    & \Sh(\Gbf, \Xcal, K_2)
  \end{tikzcd}
\]
with
\begin{equation} \label{diag:desc_pts_prod_Sh_Pink}
  \begin{tikzcd}[row sep=large]
    & \underset{hK' \in [K'' \curvearrowright K_2 / K']}{\colim} S/K(h)  \ar[dl, "{(T_h)_{hK'}}" above left] \ar[dr, "{(T_1)_{hK'}}" above right] \ar[d, "{f_{\ul{p}}}"] & \\
    S / K' \ar[dr, "{T_1}" below left] & S/K' \times S/K'' \ar[l] \ar[r] & S / K'' \ar[dl, "{T_1}" below right] \\
    & *  
  \end{tikzcd}
\end{equation}
More precisely for \(K \in \{K',K'',K(h)\}\) and \(k \in K_2\) the element \(s_0 \cdot kK\) of \(S / K\) corresponds to the geometric point \(T_k \iota_{g_0} p_{R(g_0L)}\) of \(\Sh(\Gbf, \Xcal, K)\) where \(L\) is any distinguished open subgroup of \(K_2\) contained in \(K\) and \(T_k: \Sh(\Gbf, \Xcal, L) \to \Sh(\Gbf, \Xcal, K)\).
Just like the colimit at the top of the previous diagram, the colimit of \(S/K(h)\) at the top of this diagram may be written as a disjoint union over \([h] \in K'' \backslash K_2 / K'\) if one chooses representatives.
Let us compute the morphism \(A_{g_2 K_2}\) defined in \eqref{diag:def_A_g2_K2} on stalks at the geometric point corresponding to \((s_0 \cdot x'K', s_0 \cdot x''K'')\) for some \(x',x'' \in K_2\) (i.e.\ an arbitrary point above \(p\)).
Denote \(g' = g_0x'\) and \(g'' = g_0x''\).
For \(L\) a distinguished open subgroup of \(K_2\) contained in \(K' \cap K''\) we have a commutative diagram
\[
  \begin{tikzcd}
    U(g_0L) \ar[r, hook, "{\iota_{g_0}}"] \ar[d, "{T_1}"] & \Sh(\Gbf, \Xcal, L) \ar[d, "{T_{x'}}"] \\
    U(g'K') \ar[r, hook, "{\iota_{g'}}"] & \Sh(\Gbf, \Xcal, K')
  \end{tikzcd}
\]
and similarly for \(K'',x'',g''\).
We consider the stalk of \eqref{diag:char_A_g2_K2} at the geometric point \(p_{g',g''} := p_{R(g'K')} \times p_{R(g''K'')}\) and using \eqref{eq:s_on_stalks} we obtain a commutative diagram
\begin{equation} \label{diag:A_g2_K2_stalks}
  \begin{tikzcd}[column sep=large]
    p_{g',g''}^* \pr_1^* \Fcal^{R(g'K')}(I^\bullet)^{C(g'K')} \ar[rr, "{p_{g',g''}^* (A_{g_2K_2,g',g''})}"] \ar[d, dash, "{\sim}"]  & & p_{g',g''}^* \pr_2^* \Fcal^{R(g''K'')}(I^\bullet)^{C(g''K'')} \ar[d, dash, "{\sim}"] \\
    p_{R(g'K')}^* \Fcal^{R(g'K')}(I^\bullet)^{C(g'K')} \ar[d, "{\sim}" left, "{j_{\ul{p}}}" right] & & p_{R(g''K'')}^* \Fcal^{R(g''K'')}(I^\bullet)^{C(g''K'')} \ar[d, "{\sim}" left, "{j_{\ul{p}}}" right] \\
    (I^\bullet)^{C(g'K')} \ar[r, "{N_{C(g_2K_2)/C(g'K')}}"] & (I^\bullet)^{C(g_2K_2)} \ar[r, hook] & (I^\bullet)^{C(g''K'')}
  \end{tikzcd}
\end{equation}

The case of \(B_{g_2K_2,g',g''}\) (diagrams \eqref{diag:def_B_g2_K2} and \eqref{diag:char_B_g2_K2}) is more complicated.
We have an isomorphism of functors
\begin{equation} \label{eq:stalk_f_direct}
  p_{g',g''}^* (f_{g',g''})_* \simeq \bigoplus_{\tilde{p}} \tilde{p}^*
\end{equation}
where the sum is over geometric points \(\tilde{p}\) of \(\colim_{(hK',g) \in \Gcal_{g',g''}} U(gK(h))\) satisfying \(f_{g',g''} \circ \tilde{p} = p_{g',g''}\).
We compute on stalks the first two maps in the definition of \(B_{g_2K_2,g',g''}\) in \eqref{diag:char_B_g2_K2} using the following commutative diagram
\begin{equation} \label{diag:partial_B_g2_K2_stalks}
  \begin{tikzcd}[column sep=tiny]
    p_{g',g''}^* \pr_1^* \Fcal^{R(g'K')} (I^\bullet)^{C(g'K')} \ar[r, dash, "{\sim}"] \ar[d, "{\adj}"] & p_{R(g'K')}^* \Fcal^{R(g'K')} \ar[r, "{\sim}" above, "{j_{\ul{p}}}" below] &  (I^\bullet)^{C(g'K')} \ar[d, hook, "{\diag}"] \\
    p_{g',g''}^* (f_{g',g''})_* f_{g',g''}^* \pr_1^* \Fcal^{R(g'K')} (I^\bullet)^{C(g'K')} \ar[r, dash, "{\sim}"] \ar[d, dash, "{\sim}"] & \underset{\tilde{p}}{\bigoplus} p_{g',g''}^* \pr_1^* \Fcal^{R(g'K')} (I^\bullet)^{C(g'K')} \ar[dr, dash, "{\sim}"] & \underset{\tilde{p}}{\bigoplus} (I^\bullet)^{C(g'K')} \\
    p_{g',g''}^* (f_{g',g''})_* \left( \underset{[hK',g]}{\boxplus} T_{t'}^* \Fcal^{R(g'K')} (I^\bullet)^{C(g'K')} \right) \ar[r, dash, "{\sim}"] & \underset{\tilde{p}}{\bigoplus} \tilde{p}^* T_{t'(\tilde{p})}^* \Fcal^{R(g'K')} (I^\bullet)^{C(g'K')} \ar[r, dash, "{\sim}"] & \underset{\tilde{p}}{\bigoplus} p_{R(g'K')}^* \Fcal^{R(g'K')} (I^\bullet)^{C(g'K')} \ar[u, "{\sim}" left, "{(j_{\ul{p}})_{\tilde{p}}}" right]
  \end{tikzcd}
\end{equation}
where we choose for each \(\tilde{p}\) a pair \((h(\tilde{p})K', g(\tilde{p}))\) such that \(\tilde{p}\) maps to the component of \(\colim_{(hK',g) \in \Gcal_{g',g''}} U(gK(h))\) corresponding to \((h(\tilde{p})K',g(\tilde{p}))\), the coset \(t'(\tilde{p}) Q(g'K') \in Q(g_2K_2)/Q(g'K')\) is determined as usual by \(g(\tilde{p}) h(\tilde{p}) \in t'(\tilde{p}) g' K'\), and we abusively still denote by \(\tilde{p}\) this geometric point of \(U(g(\tilde{p}) K(h(\tilde{p})))\) (i.e.\ the inclusion of the latter in the colimit is kept implicit).
Note that we have \(T_{t'(\tilde{p})} \circ \tilde{p} = p_{R(g'K')}\) by definition.
Commutativity of the top right part is the usual computation of the unit \(\id \to (f_{g',g''})_* f_{g',g''}^*\) on stalks.
Commutativity of the bottom left part is formal using the definition of \(f_{g',g''}\).

We now face two obstacles.
First, in order to compute on stalks the next map \(\left( s(Q(g'K'),Q(g'K'),t,Q(gK(h))) \right)_{hK',g}\) occurring in the definition of \(B_{g_2K_2,g',g''}\) in \eqref{diag:char_B_g2_K2} we want to use \eqref{eq:s_on_stalks}.
This requires expressing each \(\tilde{p}\) as \(T_{q'} \circ p_{R(q' g(\tilde{p}))}\) for some \(q' \in Q(g_2K_2)\) and using \(j_{T_{q'}(\ul{p})}\) and \(j_{T_{t'} T_{q'}(\ul{p})}\).
We also need to do the same for the right side of the diagram \eqref{diag:char_B_g2_K2}.
This would merely complicate computations.
More importantly, the index set of lifts \(\tilde{p}\) of the geometric point \(p_{g',g''}\) is defined only implicitly.
We need a more explicit description to compare with \eqref{diag:A_g2_K2_stalks}.
Since the diagram \eqref{diag:cart_prod_Sh_U1} is Cartesian the set of lifts \(\tilde{p}\) of \(p_{g',g''}\) is in bijection with the set of lifts of \(\iota_{g',g''} p_{g',g''}\) along \(f\), which may be described explicitly using the diagram \eqref{diag:desc_pts_prod_Sh_Pink}.
This is the purpose of the following lemma, which fortunately will also make our first obstacle disappear.

\begin{lemm}
  \begin{enumerate}
  \item Consider the (left) action of  \(K'' \times K' \times \Gamma\) on \(K_2 \times K_2\) defined by
    \[ (k'',k',\gamma) \cdot (h,z) = (k'' h (k')^{-1}, \gamma z (k'')^{-1}). \]
    We have a well-defined bijection
    \begin{align*}
      i: (K'' \times K' \times \Gamma) \backslash (K_2 \times K_2) & \longrightarrow \bigsqcup_{[h] \in K'' \backslash K_2 / K'} S/K(h) \\
      [h,z] & \longmapsto [hK', s_0 \cdot zK(h)]
    \end{align*}
    and we have \(f_{\ul{p}} i([h,z]) = (s_0 \cdot zh K', s_0 \cdot z K'')\).
  \item For \(x',x'' \in K_2\) we have a bijection
    \begin{align} \label{eq:geom_fiber_prod_Pink_Shim}
      C(g_0x'' K'') \backslash C(g_0K_2) / C(g_0 x' K') & \longrightarrow (f_{\ul{p}} i)^{-1}(s_0 \cdot x'K', s_0 \cdot x''K'') \\
      [\alpha] & \longmapsto [(g_0x'')^{-1} \alpha g_0 x', x'']. \nonumber
    \end{align}
  \end{enumerate}
\end{lemm}
\begin{proof}
  \begin{enumerate}
  \item We leave this elementary verification to the reader.
  \item An element of the fiber is a class \([h,z]\) satisfying \(g_0 z \in C(g_0K_2) g_0 x'' K''\) and \(g_0 z h \in C(g_0K_2) g_0 x' K'\).
    Up to translating by an element of \(K'' \times K' \times \Gamma\) we may assume both \(g_0 z = g_0 x''\) and \(g_0 z h = \alpha g_0 x'\) for some \(\alpha \in C(g_0K_2)\).
    So we have a surjective map \(C(g_0K_2) \to (f_{\ul{p}} i)^{-1}(s_0 \cdot x'K', s_0 \cdot x''K'')\), and an uneventful computation shows that it induces a bijective map \eqref{eq:geom_fiber_prod_Pink_Shim}.
  \end{enumerate}
\end{proof}

The second point in the lemma allows us the parametrize (with \eqref{diag:desc_pts_prod_Sh_Pink}) the direct sum on the right-hand side of \eqref{eq:stalk_f_direct} by \([\alpha] \in C(g'' K'') \backslash C(g_0 K_2) / C(g'K')\).
More precisely for \(\alpha \in C(g_2K_2)\) the corresponding geometric point \(\tilde{p}\) of \(\colim_{(hK',g) \in \Gcal_{g',g''}} U(gK(h))\) is \(p_{R(g(\alpha)K(h(\alpha)))}\) (via the canonical clopen immersion of \(U(g(\alpha) K(h(\alpha)))\) in the colimit) where \(g(\alpha) = g''\) and \(h(\alpha) = (g_0 x'')^{-1} \alpha g_0 x'\).
Indeed for \(L\) small enough we have \(\iota_{g(\alpha)} \tilde{p} = T_{x''} \, \iota_{g_0} \, p_{R(g_0 L)}\) and a commutative diagram
\[
  \begin{tikzcd}
    U(g_0L) \ar[r, hook, "{\iota_{g_0}}"] \ar[d, "{T_1}"] & \Sh(\Gbf, \Xcal, L) \ar[d, "{T_{x''}}"] \\
    U(g(\alpha) K(h(\alpha))) \ar[r, hook, "{\iota_{g(\alpha)}}"] & \Sh(\Gbf, \Xcal, K(h(\alpha))).
  \end{tikzcd}
\]
So in the diagram \eqref{diag:partial_B_g2_K2_stalks} we may index by a set of such representatives \(\alpha\) instead of \(\{\tilde{p}\}\), and replace \((h(\tilde{p}), g(\tilde{p}), t'(\tilde{p}))\) by \((h(\alpha), g'', \alpha)\).
We note that the restriction of \(f_{g',g''}\) to the component \(U(g'' K(h))\) is
\[ T_{\alpha} \times T_1: U(g'' K(h)) \longrightarrow U(g'K') \underset{U(g_2K_2)}{\times} U(g''K'') \]
and that we simply have \(T_{\alpha} = T_1\) because \(\alpha \in C_\Q\).
Continuing the diagram \eqref{diag:partial_B_g2_K2_stalks} to compute the left side of \eqref{diag:char_B_g2_K2} on stalks using \eqref{eq:s_on_stalks}, and proceeding similarly for the right side of \eqref{diag:char_B_g2_K2}, we finally obtain a commutative diagram
\begin{equation} \label{diag:B_g2_K2_stalks}
  \begin{tikzcd}[column sep=large]
    p_{g',g''}^* \pr_1^* \Fcal^{R(g'K')} (I^\bullet)^{C(g'K')} \ar[rr, "{p_{g',g''}^*(B_{g_2K_2,g',g''})}"] \ar[d, dash, "{\sim}"] & & p_{g',g''}^* \pr_2^* \Fcal^{R(g''K'')} (I^\bullet)^{C(g''K'')} \ar[d, dash, "{\sim}"] \\
    p_{R(g'K')}^* \Fcal^{R(g'K')} (I^\bullet)^{C(g'K')} \ar[d, "{\sim}" left, "{j_{\ul{p}}}" right] & & p_{R(g''K'')}^* \Fcal^{R(g''K'')} (I^\bullet)^{C(g''K'')} \ar[d, "{\sim}" left, "{j_{\ul{p}}}" right] \\
    (I^\bullet)^{C(g'K')} \ar[d, hook, "{\diag}"] & & (I^\bullet)^{C(g''K'')} \\
    \underset{[\alpha]}{\bigoplus} (I^\bullet)^{C(g'K')} \ar[r, "{(\alpha)_{\alpha}}"] & \underset{\alpha}{\colim} (I^\bullet)^{C(g(\alpha) K(h(\alpha)))} \ar[r, "{(N_{\alpha})_{\alpha}}"] & \underset{[\alpha]}{\bigoplus} (I^\bullet)^{C(g''K'')} \ar[u, "{\Sigma}"]
  \end{tikzcd}
\end{equation}
where \(N_{\alpha}\) is the norm map for \(C(g(\alpha) K(h(\alpha))) \subset C(g''K'')\) and the colimit is over the groupoid \([ C(g''K'') \curvearrowright C(g_2K_2)/C(g'K') ]\).

To conclude we are left to compare the diagrams \eqref{diag:A_g2_K2_stalks} and \eqref{diag:B_g2_K2_stalks}.
The following lemma implies that the bottom paths in these diagrams are in fact equal, by decomposing \(C(g_2K_2) / C(g' K')\) into \(C(g'' K'')\)-orbits.

\begin{lemm}
  For \(\alpha \in C(g_2K_2)\), denoting \(h = (g'')^{-1} \alpha g'\) we have a bijection
  \begin{align*}
    C(g''K'') / C(g''K(h)) & \longrightarrow C(g''K'') \backslash C(g''K'') \alpha C(g'K') / C(g'K') \\
    [\delta] & \longmapsto \delta \alpha
  \end{align*}
\end{lemm}
\begin{proof}
  This follows from the equality \(C(g''K(h)) = C(g''K'') \cap \alpha C(g'K') \alpha^{-1}\).
\end{proof}

\subsection{Direct product case}
\label{sec:ext_Shim_direct_prod}

Assume \(\Gbf = \Gbf_\lin \times \Gbf_\her\) where \(\Gbf_\lin(\R)\) acts trivially on \(\Xcal\).
In particular we have \(S_\Q \Gbf_\her(\A_f) = \Gbf_\lin(\Q) \times \Gbf_\her(\A_f)\).
For neat \(K\) we have \(C(gK) = \Gbf_\lin(\Q) \cap gKg^{-1}\) and \(Q(gK) = (\Gbf_\lin(\Q) \times \Gbf_\her(\A_f) ) \cap gKg^{-1}\).
If \(K\) is neat then \(C(gK)\) is a neat arithmetic subgroup of \(\Gbf_\lin(\Q)\), and if \(K\) factors as \(K_\lin \times K_\her\) then we have \(Q(gK) = C(gK) \times P(gK)\) and we can compute
\[ R\Gamma(C(gK),-): D^+(Q(gK), E) \longrightarrow D^+(P(gK), E) \]
as \(\Tot^\bullet \left( \Hom_{\Z[C(gK)]}(F^\bullet,-) \right)\) using a resolution \(F^\bullet\) of the \(\Z[C(gK)]\)-module \(\Z\) consisting of finite free \(\Z[C(gK)]\)-modules (see Proposition \ref{pro:RGamma_CgK_finite}).
For such \(K\) we also have a bijection
\[ \Gbf_\lin(\Q) \backslash \Gbf_\lin(\A_f) / K_\lin \xrightarrow{\sim} S_\Q \Gbf_\her(\A_f) \backslash \Gbf(\A_f) / K. \]
If \(K' = K'_\lin \times K'_\her\) is an open (also factorizable) subgroup of \(K\) then we may restrict the resolution \(F^\bullet\) above from \(C(gK)\) to \(C(gK')\) to compute \(R\Gamma(C(gK'),-)\) as well.
By double complex arguments the restriction and corestriction maps may be computed using \(F^\bullet\): we have commutative diagrams of functors
\[
  \begin{tikzcd}
    R\Gamma(C(gK), -) \ar[r, "{\res}"] \ar[d, dash, "{\sim}"] & R\Gamma(C(gK'), -) \ar[d, dash, "{\sim}"] \\
    \Tot^\bullet(\Hom_{\Z[C(gK)]}(F^\bullet), -) \ar[r] & \Tot^\bullet(\Hom_{\Z[C(gK')]}(F^\bullet), -)
  \end{tikzcd}
\]
where the bottom morphism of functors is given by obvious inclusions, and
\[
  \begin{tikzcd}
    R\Gamma(C(gK'), -) \ar[r, "{\cores}"] \ar[d, dash, "{\sim}"] & R\Gamma(C(gK), -) \ar[d, dash, "{\sim}"] \\
    \Tot^\bullet(\Hom_{\Z[C(gK')]}(F^\bullet), -) \ar[r] & \Tot^\bullet(\Hom_{\Z[C(gK)]}(F^\bullet), -)
  \end{tikzcd}
\]
where the bottom morphism of functors is induced by the norm maps
\[ \Hom_{\Z}(F^i, M)^{C(gK')} \xrightarrow{N_{C(gK)/C(gK')}} \Hom_{\Z}(F^i, M)^{C(gK)}. \]

Let us temporarily denote by \(S\) any topos.
We have an ``external tensor product'' bifunctor
\begin{align*}
  \boxtimes: \Perf(\Ocal_E) \times D^+(S, \Ocal_E) & \longrightarrow D^+(S, \Ocal_E) \\
  (M^\bullet, N^\bullet) & \longmapsto \Tot^\bullet( M^\bullet \otimes_{\Ocal_E} N^\bullet)
\end{align*}
where \(\Tot^n(M^\bullet \otimes_{\Ocal_E} N^\bullet) = \bigoplus_{a+b=n} M^a \otimes_{\Ocal_E} N^b\), with differentials \((-1)^b d_M^a \otimes \id_{N^b} + \id_{M^a} \otimes d_N^b\), inducing a bifunctor
\[ \boxtimes: D^b(E) \times D^+(S, E) \longrightarrow D^+(S, E), \]
using the identification \(\Perf(\Ocal_E)[\ell^{-1}] \simeq D^b(E)\).
Now \(D^b(E)\) is abelian, equivalent to its full subcategory of complexes with vanishing differentials.

We come back to the setting above: \(K = K_\lin \times K_\her\) is a neat (factorizable) compact open subgroup of \(\Gbf(\A_f)\).
Let \(\Lambda_\lin\) be an \(\Ocal_E[C(gK)]\)-modules which is finite free as \(\Ocal_E\)-module.
For \(\Lambda_\her \in D^+(P(gK), \Ocal_E)\) we define \(\Lambda_\lin \otimes_{\Ocal_E} \Lambda_\her \in D^+(Q(gK), \Ocal_E)\) in the obvious way.
The computation of group cohomology for \(C(gK)\) recalled above yields
\[ R\Gamma(C(gK), \Lambda_\lin \otimes_{\Ocal_E} \Lambda_\her) \simeq R\Gamma(C(gK), \Lambda_\lin) \boxtimes \Lambda_\her \]
compatible with restriction and corestriction (commutative diagrams left to the reader).
Similarly we can replace \(\Lambda_{\lin}\) by an \(E[C(gK)]\)-module having finite dimension over \(E\) and admitting \(C(gK)\)-stable \(\Ocal_E\)-lattice and \(\Lambda_{\her}\) by any object of \(D^+(P(gK),E)\).
It follows that for \(K_1\), \(K_2\) and \(K'\) neat factorizable compact open subgroups of \(\Gbf(\A_f)\) satisfying \(K' \subset K_1 \cap K_2\) we can compute the correspondence \(au(K_2,1,1,K_1,K')\) (Definition \ref{def:corr_au}) as a kind of tensor product using the following commutative diagram, where \(g \in \Gbf(\A_f)\) is arbitrary and with morphisms \(s(\dots)\) defined in \eqref{eq:T_t_pullback_Pink_gal}
\[
  \begin{tikzcd}
    T_1^* \Fcal^{R(gK_1)} R\Gamma(C(gK_1), V_\lin \otimes V_\her) \ar[r, dash, "{\sim}"] \ar[d, "{\sim}" right, "{s(Q(gK_1),1,Q(gK'))}" left] & R\Gamma(C(gK_1), V_\lin) \boxtimes T_1^* \Fcal^{P(gK_1)} V_\her \ar[dddd, "{(\cores \circ \res) \boxtimes u(K_2,1,K_1,K')}"] \\
    \Fcal^{R(gK')} R\Gamma(C(gK_1), V_\lin \otimes V_\her) \ar[d, "{\Fcal^{R(gK')}(\res)}"] & \\
    \Fcal^{R(gK')} R\Gamma(C(gK'), V_\lin \otimes V_\her) \ar[d, "{\Fcal^{R(gK')}(\cores)}"] & \\
    \Fcal^{R(gK')} R\Gamma(C(gK_2), V_\lin \otimes V_\her) & \\
    T_1^* \Fcal^{R(gK_1)} R\Gamma(C(gK_1), V_\lin \otimes V_\her) \ar[r, dash, "{\sim}"] \ar[u, "{\sim}" right, "{s(Q(gK_2),1,Q(gK'))}" left] & R\Gamma(C(gK_2), V_\lin) \boxtimes T_1^* \Fcal^{P(gK_2)} V_\her
  \end{tikzcd}
\]
where the correspondence \(u(K_2,1,K_1,K')\) is as in Definition \ref{def:Hecke_corr}.
It follows that for \(? \in \{*,!\}\) we have a commutative diagram in \(D^+(B, E)\):
\[
  \begin{tikzcd}
    \pi_? \AFcal^{K_1} (V_\lin \otimes V_\her) \ar[r, dash, "{\sim}"] \ar[d, "{au(K_2,1,K_1,K')_?}"] & R\Gamma(\Gbf_\lin, V_\lin, K_1) \boxtimes \pi_? \Fcal^{K_{1,\her}} V_\her \ar[d, "{\left[ K_{2,\lin},1,K_{1,\lin},K'_{\lin} \right] \boxtimes u(K_{2,\her}, 1, K_{1,\her}, K'_\her)}"] \\
    \pi_? \AFcal^{K_2} (V_\lin \otimes V_\her) \ar[r, dash, "{\sim}"] & R\Gamma(\Gbf_\lin, V_\lin, K_2) \boxtimes \pi_? \Fcal^{K_{2,\her}} V_\her.
  \end{tikzcd}
\]
The reduction of the general case (correspondences \(au(K_2,g,K_1,K')\) with \(g = g_\lin g_\her \in \Gbf(\A_f)\) arbitrary) is quite formal and we omit the details.

To conclude with explicit formulas in the direct product case it is useful to be more explicit with the \(\ell\)-adic formalism that we use.
As in \cite[\S 6]{Ekedahl_adic} we denote by \(c\) be the category of constructible sheaves of \(\Ocal_E/\mfrak_E\)-vector spaces in \(\wt{B_\mathrm{et}}\), and we denote by \(D^b_c(B, \Ocal_E)\) the associated (\S 3 loc.\ cit.) triangulated category, and by \(D^b_c(B, E)\) the triangulated category obtained by inverting \(\ell\) in the \(\Ocal_E\)-modules of morphisms.
The external tensor product maps \(\Perf(\Ocal_E) \times D^b_c(B, \Ocal_E)\) to \(D^b_c(B, \Ocal_E)\), as one checks using ``stupid truncations'' of \(M \in \Perf(\Ocal_E)\) to reduce to the case where \(M\) is concentrated in one degree.
Recall from \cite[Theorem 3.6.v]{Ekedahl_adic} that \(D^b_c(B,E)\) admits a \(t\)-structure whose heart may be identified with the \(E\)-linear abelian category \(\Acal(c,E)\) of \(c-\Ocal_E\)-modules in \(\wt{B_\mathrm{et}}^\N\) with \(\ell\) inverted in morphisms\footnote{In the case where \(B=\Spec k\), if we choose a separable closure of \(k\) then we get an identification of this category with the category of finite-dimensional continuous representations of the absolute Galois group of \(k\) over \(E\).}.
We denote by \(H^i_t\) the corresponding homological functors \(D^b_c(B,E) \to \Acal(c,E)\)\footnote{More concretely the functors \(H^i_t\) may be identified with the restriction of the usual functor \(H^i\) on \(D^-(\wt{B_\mathrm{et}}^\N, (\Ocal_E)_\bullet)\) composed with the quotient functor from \((\Ocal_E)_\bullet\)-modules in \(\wt{B_\mathrm{et}}^\N\) to the quotient by the Serre subcategory of essentially zero systems.
It is easy to check that the restriction of this quotient functor to \(c-\Ocal_E\) sheaves (in the sense of \cite[before Theorem 3.6]{Ekedahl_adic}) is fully faithful.}.
In the Grothendieck group of \(\Hecke(\Gbf(\A_f), \Acal(c,E))\) we have an equality
\[ e_?(\Gbf, \Xcal, B) = e(\Gbf_\lin) \boxtimes e_?(\Gbf_\her, \Xcal_1, B) \]
where
\[ e_?(\Gbf, \Xcal, B) = \sum_{i \geq 0} (-1)^i \left[ \left( H^i_t(\pi_? \AFcal^K(V_\lin \otimes V_\her)) \right)_K, \left( au(K_2, g, K_1, K')_? \right)_{K_2,g,K_1,K'} \right] \]
\[ e_?(\Gbf_\her, \Xcal_1, B) = \sum_{i \geq 0} (-1)^i \left[ \left( H^i_t(\pi_? \Fcal^{K_\her}(V_\her)) \right)_{K_\her}, \left( u(K_{2,\her}, g_\her, K_{1,\her}, K'_\her)_? \right)_{K_{2,\her}, g_\her, K_{1,\her}, K'_\her} \right] \]
\[ e(\Gbf_\lin) = \sum_{i \geq 0} (-1)^i \left[ (H^i(\Gbf_\lin, V_\lin, K_\lin))_{K_\lin}, \left( \left[ K_{2,\lin}, g_\lin, K_{1,\lin}, K'_\lin \right] \right)_{K_{2,\lin}, g_\lin, K_{1,\lin}, K'_\lin} \right], \]
the latter in the Grothendieck group of
\[ \Hecke(\Gbf_\lin(\A_f), \text{finite-dimensional } E \text{-vector spaces}) \]
i.e.\ admissible representations of \(\Gbf_\lin(\A_f)\) over \(E\).

\subsection{Minimal compactifications}

We will see how generalized Shimura varieties occur naturally as boundary strata in minimal compactifications of (usual, or even generalized) Shimura varieties.

\subsubsection{Real connected case}

First we recall from \cite[\S III]{AMRT} the description of boundary components.
Let \(\Gbf\) be a connected reductive group over \(\R\), assumed to be adjoint and simple, and let \(D\) be a \(\Gbf(\R)^0\)-orbit in \(\Hom(\mathbb{S}, \Gbf)\).
We assume as usual that \(\Gbf\) is isotropic, and that for any (equivalently, one) \(h \in D\) we have
\begin{itemize}
\item the restriction of \(h\) via \(\GLbf_{1,\R} \hookrightarrow \mathbb{S}\) is trivial,
\item letting \(\mathbb{S}_\C\) act on \(\Lie \Gbf_\C\) via \(h\) and the adjoint representation, only the trivial character and the characters \(z \mapsto (z/\ol{z})^{\pm 1}\) occur in the diagonalization of this action,
\item the involution \(\Ad h(i)\) of \(\Gbf\) is a Cartan involution.
\end{itemize}
The map \(h \in D \mapsto h(i)\) is injective and so we may see \(D\) as parametrizing Cartan involutions of \(\Gbf\) (equivalently, maximal compact subgroups of \(\Gbf(\R)\) or \(\Gbf(\R)^0\)).
For \(h \in D\) denote by \(\Kbf_h\) the centralizer of \(h\) in \(\Gbf\).
In particular \(K_h = \Kbf_h(\R)\) is a maximal compact subgroup of \(\Gbf(\R)^0\).
For any \(h \in D\) we have an associated morphism \(\mu_h: \GLbf_{1,\C} \to \Gbf_\C\), determining a parabolic subgroup of \(\Gbf_\C\).
This yields (see e.g.\ \cite[Theorem III.2.1]{AMRT}) a natural \(\Gbf(\R)^0\)-equivariant open embedding of \(D\) in the complex points of a Grassmannian \(\check{D}\), mapping \(h \in D\) to the parabolic subgroup \(\Qbf\) of \(\Gbf_\C\) whose Lie algebra is the sum of the nonpositive eigenspaces for the adjoint action of \(\mu_h\).
Following \cite[\S III.3]{AMRT} the closure \(\ol{D}\) of \(D\) in \(\check{D}\) may be described as the disjoint union of \(D\) and of lower-dimensional simple hermitian symmetric domains.
Choose \(h_0 \in D\), choose a maximal torus of \(\Kbf_{h_0}\), and choose (following Harish-Chandra) a maximal set of strongly orthogonal non-compact roots (see \cite[\S III.2.3]{AMRT}).
This yields a morphism \(\SLbf_{2,\R}^r \to \Gbf\) with finite kernel, compatible with Cartan involutions, mapping the diagonal maximal split torus of \(\SLbf_{2,\R}^r\) to a maximal split torus of \(\Gbf\) which is \(\Ad h_0(i)\)-stable.
Taking the \(\Gbf(\R)^0\)-orbit, we obtain a homogeneous space \(\Sigma\) under \(\Gbf(\R)^0\) and a \(\Gbf(\R)^0\)-equivariant map \(\Sigma \to D\) such that the fiber of \(h\) is a \(K_h\)-orbit of morphisms \(f: \SLbf_{2,\R}^r \to \Gbf\) as above (intertwining the standard Cartan involution on \(\SLbf_{2,\R}^r\) with \(\Ad h(i)\)).
Following \cite[\S 4.3]{Pink_dissertation} let \(H_r\) be the subgroup of \(\mathbb{S} \times \GLbf_{2,\R}^r\) defined by
\begin{equation} \label{eq:def_Hr}
  H_r(T) = \left\{ (z,g_1,\dots,g_r) \in \mathbb{S}(T) \times \GLbf_2(T)^r \,\middle|\, z \ol{z} = \det g_1 = \dots = \det g_r \right\}
\end{equation}
for any \(\R\)-scheme \(T\).
Denote
\begin{align*}
  h_\std: \mathbb{S} & \longrightarrow \GLbf_{2,\R} \\
  z = a+ib & \longmapsto \begin{pmatrix} a & b \\ -b & a \end{pmatrix}.
\end{align*}
It follows from \cite[Theorem III.2.4]{AMRT} that for any \(f \in \Sigma\) above \(h \in D\) we have a unique morphism \(\varphi_f: H_r \longrightarrow \Gbf\) satisfying
\begin{itemize}
\item the restriction of \(\varphi_f\) to \(\SLbf_{2,\R}^r \subset H_r\) is \(f\),
\item for any \(z \in \mathbb{S}(\R)\) we have \(\varphi_f(z, h_{\std}(z), \dots, h_{\std}(z)) = h(z)\).
\end{itemize}
More precisely \cite[Theorem III.2.4]{AMRT} gives a morphism \(\mathbb{S}^1 \times \SLbf_{2,\R}^r \to \Gbf\) satisfying these two conditions, where \(\mathbb{S}^1\) is the one-dimensional anisotropic subtorus of \(\mathbb{S}\), and it extends uniquely to \(H_r\) by requiring it to be trivial on the diagonally embedded \(\GLbf_{1,\R} \hookrightarrow H_r\).
The map \(f \mapsto \varphi_f\) is clearly \(\Gbf(\R)^0\)-equivariant.
Recall that a parabolic subgroup of \(\Gbf\) is called \emph{admissible} if it is equal to \(\Gbf\) or maximal among proper parabolic subgroups.
Let \(S\) be a subset of \(\{1, \dots, r\}\) of cardinality \(b\).
For a set \(X\) denote
\[ r_S: X^2 \to X^r,\ \ \ r_S(x,y)_i =
  \begin{cases}
    x & \text{ if } i \in S \\
    y & \text{ if } i \not\in S.
  \end{cases}
\]
For \(S \subset \{1,\dots,r\}\) having cardinality \(b\) we have an admissible parabolic subgroup \(\Pbf_{f,S}\) of \(\Gbf\) associated to the cocharacter
\begin{align}
  \GLbf_{1,\R} & \longrightarrow \Gbf \label{eq:cochar_f_S_to_P} \\
  t & \longmapsto \varphi_f \left(1, r_S(\diag(t,t^{-1}), I_2) \right), \nonumber
\end{align}
i.e.\ the Lie algebra of \(\Pbf_{f,S}\) is the sum of the nonnegative eigenspaces for the adjoint action of this cocharacter.
We have \(\Pbf_{f,S} = \Gbf\) if and only if \(b=0\).
\begin{rema}
  For \(h \in D\), for any admissible parabolic subgroup of \(\Gbf\) there exists \(f \in \Sigma\) above \(h\) and \(0 \leq b \leq r\) such that we have \(\Pbf = \Pbf_{f,S}\) for \(S=\{1,\dots,b\}\): this follows from the fact that the relative root system of \(\Gbf\) is of type \(C_r\) or \(BC_r\) \cite[\S I.1.2]{BailyBorel} and the fact that \(K_h\) acts transitively on each \(\Gbf(\R)\)-orbit of parabolic subgroups of \(\Gbf\).
\end{rema}
It turns out that, by explicit computation in bounded symmetric domains, \(\ol{D}\) decomposes as the disjoint union of complex submanifolds \(F_\Pbf\), one for each admissible parabolic subgroup \(\Pbf\) of \(\Gbf\), with \(F_{\Gbf} = D\), \(\mathrm{Stab}(F_\Pbf, \Gbf(\R)^0) = \Pbf(\R) \cap \Gbf(\R)^0\) and \(g \cdot F_\Pbf = F_{\Ad(g)(\Pbf)}\) for any admissible parabolic \(\Pbf\) and \(g \in \Gbf(\R)^0\).
We briefly recall how to associate \(F_\Pbf\) to an admissible parabolic subgroup \(\Pbf\).
First choose \(h \in D\), and choose \(f \in \Sigma\) above \(h\) and \(S \subset \{1, \dots, r\}\) such that we have \(\Pbf = \Pbf_{f,S}\).
There is a unique Levi factor \(\Mbf_{\Pbf,h}\) of \(\Pbf\) which is stable under \(\Ad h(i)\), namely \(\Pbf \cap \Ad(h(i))(\Pbf)\), also equal to the centralizer of \eqref{eq:cochar_f_S_to_P}.
In \cite[\S III.3.1 and Theorem III.3.10]{AMRT} a semi-simple ideal \(\lfrak_{f,S}\) of \(\Lie \Mbf_{\Pbf,h}\) is defined.
We have \(\lfrak_{f,S} = \Lie \Lbf_{f,S}\) for a unique connected semi-simple subgroup \(\Lbf_{f,S}\) of \(\Mbf_{\Pbf,h}\), and the adjoint quotient \(\Lbf_{f,S,\ad}\) of \(\Lbf_{f,S}\) is a simple factor of \(\Mbf_{\Pbf,h,\ad}\).
Denote \(S^c := \{1,\dots,r\} \smallsetminus S\).
The group \(\Lbf_{f,S}\) contains the image of \(\SLbf_{2,\R}^{S^c}\)  via \(\varphi_{f,S}\), and commutes with the image of \(\SLbf_{2,\R}^S\).
Thus the subgroup \(\Hbf_{f,S}\) of \(\Gbf\) generated by \(\Lbf_{f,S}\) and the image of \(\varphi_f\) is reductive, with derived subgroup a quotient of \(\SLbf_{2,\R}^S \times \Lbf_{f,S}\) by a finite central subgroup.
We have a factorization of \(\varphi_f\) through \(\varphi_{f,S} : H_r \to \Hbf_{f,S}\), in particular \(h\) factors through \(h_{f,S} : \mathbb{S} \to \Hbf_{f,S}\).
These constructions are obviously \(\Gbf(\R)^0\)-equivariant, in particular \(\Hbf_{f,S}(\R)^0\)-equivariant.
The \(\Hbf_{f,S,\ad}(\R)^0\)-orbit of the image of \(h_{f,S}\) in
\[ \Hom(\mathbb{S}, \Hbf_{f,S,\ad}) = \Hom(\mathbb{S}, \PGLbf_{2,\R}^S \times \Lbf_{f,S,\ad}) \]
clearly satisfies the same conditions as \(D\), and is naturally identified with \(\Hcal^S \times D_{f,S}\) where \(\Hcal\) is the \(\PGLbf_2(\R)^0\)-orbit of \(h_\std\) (also known as the upper half-plane) and \(D_{f,S}\) is the \(\Lbf_{f,S,\ad}(\R)^0\)-orbit of the image of \(h_{f,S}\) in \(\Hom(\mathbb{S}, \Lbf_{f,S,\ad})\).
The holomorphic embedding \(\Hcal^S \times D_{f,S} \hookrightarrow D\) extends uniquely to a holomorphic embedding \(\iota_{f,S}: \mathbb{P}^1(\C)^S \times \check{D_{f,S}} \hookrightarrow \check{D}\), essentially because \(\mu_h\) factors through \(\Hbf_{f,S,\C}\).
The boundary component associated to \((f,S)\) is defined as \(F_{f,S} := \iota_{f,S} \left( \{\infty\}^S \times D_{f,S} \right)\).
With this description it seems to depend on the choice of \(f\) and \(S\) satisfying \(\Pbf_{f,S} = \Pbf\), but in fact it does not so we denote it by \(F_\Pbf\) (see \cite[Theorem 3.7 and Proposition 3.9]{AMRT}, which also shows that \(\Pbf \mapsto F_\Pbf\) is a bijection between admissible parabolic subgroups and boundary components and that \(\Pbf(\R) \cap \Gbf(\R)^0\) is the stabilizer of \(F_\Pbf\) in \(\Gbf(\R)^0\)).
We have \(h = \iota_{f,S}(i,\dots,i,h')\) for some \(h' \in D_{f,S}\) so we have an associated point \(\iota_{f,S}(\infty,\dots,\infty,h') \in F_\Pbf\).
It will be useful to recall a more intrinsic description of this point.
The maximal split central torus \(\Abf_{\Mbf_{\Pbf,h}}\) in \(\Mbf_{\Pbf,h}\) has dimension one, and there is a unique isomorphism \(u_{\Pbf,h}: \GLbf_{1,\R} \simeq \Abf_{\Mbf_{\Pbf,h}}\) such that in the adjoint action of \(\GLbf_{1,\R}\) on \(\Lie \Nbf_{\Pbf}\), only positive characters occur.
For any \(h' \in D_{f,S}\) we have
\[ \iota_{f,S}(\infty, \dots, \infty, h') = \lim_{t \to + \infty} u_{\Pbf,h}(t) \cdot \iota_{f,S}(i,\dots,i,h'). \]
Denoting
\begin{align*}
  \pi_\Pbf: D & \longrightarrow \ol{D} \\
  h & \longmapsto \lim_{t \to + \infty} u_{\Pbf,h}(t) \cdot h
\end{align*}
we see that for any \(g \in \Gbf(\R)^0\) and \(h \in D\) we have \(\pi_{\Ad(g) \Pbf}(g \cdot h) = g \cdot \pi_\Pbf(h)\), in particular \(\pi_\Pbf\) is \(\Pbf(\R) \cap \Gbf(\R)^0\)-equivariant.
By \cite[Theorem III.3.10 (2)]{AMRT} the unipotent radical \(\Nbf_\Pbf(\R)\) of \(\Pbf(\R)\) acts trivially on \(F_\Pbf\) (so does the connected centralizer of \(\Lbf_{f,S}\) in \(\Mbf_{\Pbf,h}(\R)\)).
We also deduce \(\pi_\Pbf(D) = F_\Pbf\).
The map \(\pi_\Pbf\) defined above is the \emph{geodesic projection} denoted by \(\pi_F\) in \cite[\S III.3.4]{AMRT}, see p.140 loc.\ cit.
In particular it is holomorphic and its image \(F_\Pbf\) is a submanifold of \(\check{D}\), isomorphic to \(D_{f,S}\) for any pair \((f,S)\) satisfying \(\Pbf_{f,S} = \Pbf\).
By \(\Pbf(\R)^0\)-equivariance (or from the description in \cite[\S III.4.1]{AMRT}) it is also clear that the image of the subgroup \(\Lbf_{f,S}\) of \(\Mbf_{\Pbf,h}\) in the reductive quotient \(\Mbf_\Pbf\) of \(\Pbf\) does not depend on the choice of \((f,S)\) mapping to \(\Pbf\), and we denote it by \(\Mbf_{\Pbf,\her,\der}\).
Also denote by \(\Mbf_{\Pbf,\her,\ad}\) its adjoint quotient, so that \(F_\Pbf\) is a hermitian symmetric domain with automorphism group \(\Mbf_{\Pbf,\her,\ad}(\R)^0\) (now independently of a choice of \((f,S)\)).
Pink observed (see \cite[Proposition 4.6]{Pink_dissertation} and \cite[\S 3.6]{Pink_ladic_Shim}) that the geodesic projection \(\pi_\Pbf\) may be interpreted à la Deligne (i.e.\ with morphisms from \(\mathbb{S}\)) and more intrinsically as follows: there is a unique identification of \(F_\Pbf\) with a \(\Mbf_\Pbf(\R)^0\)-orbit in \(\Hom(\mathbb{S}, \Mbf_\Pbf)\) such that for any \(h \in D\) and \(z \in \mathbb{S}(\R)\) we have
\begin{equation} \label{eq:pi_P_simple_real}
  \pi_\Pbf(h)(z) = \varphi_{f,S}(z,r_S(\diag(z\ol{z},1), h_\std(z)))
\end{equation}
for any \(f \in \Sigma\) and \(S \subset \{1,\dots,r\}\) mapping to \((h,\Pbf)\).
Indeed it is clear that projecting the right-hand side of \eqref{eq:pi_P_simple_real} to \(\Mbf_{\Pbf,\her,\ad}\) recovers \(\iota_{f,S}(\infty, \dots, \infty, h')\) where \(h = \iota_{f,S}(i, \dots, i, h')\).
Moreover
\begin{align*}
  \varphi_{h,\Pbf}: H_1 & \longrightarrow \Gbf \\
  (z,g) & \longmapsto \varphi_{f,S}(z,r_S(g,h_\std(z)))
\end{align*}
is the unique morphism satisfying
\begin{enumerate}
\item for all \(z \in \mathbb{S}\) we have \(\varphi_{h,\Pbf}(z, h_\std(z)) = h(z)\),
\item \(h_\infty: z \mapsto \varphi_{h,\Pbf}(z, \diag(z\ol{z},1))\) takes values in \(\Pbf\),
\item \label{it:cond3_varphi_h}
  the restriction of \(h_\infty\) to \(\GLbf_{1,\R} \hookrightarrow \mathbb{S}\) has adjoint action \(t \mapsto t^2\) on the center \(\ufrak\) of \(\nfrak_\Pbf\) and by \(t \mapsto t\) on \(\nfrak_\Pbf/\ufrak\) (see \cite[\S III.4]{AMRT}).
\end{enumerate}
With this characterization we thus have
\begin{equation} \label{eq:pi_P_simple_real_can}
  \pi_\Pbf(h) = \varphi_{h,\Pbf}(z, \diag(z\ol{z},1)).
\end{equation}

\begin{rema} \label{rem:cond3_varphi_h_Pink}
  Condition \eqref{it:cond3_varphi_h} is stated differently in \cite[Proposition 4.6]{Pink_dissertation} and \cite[\S 3.6]{Pink_ladic_Shim}, but it does not seem to be the correct condition for \(0\)-dimensional boundary components (in this case \(\ufrak = \nfrak_\Pbf\)), e.g.\ for \(\Gbf=\PGLbf_{2,\R}\) and \(\Pbf\) a Borel subgroup.
\end{rema}

\subsubsection{Boundary strata of minimal compactifications}
\label{sec:gen_Shim_boundary}

We now turn to the global setting and consider a connected reductive group \(\Gbf\) over \(\Q\) and a generalized Shimura datum \(h: \Xcal \to \Hom(\mathbb{S}, \Gbf_\R)\).
Again we fix a \(\Gbf(\R)\)-equivariant family of strongly orthogonal non-compact roots, i.e.\ choose any \(x_0 \in \Xcal\), a maximal torus of the centralizer \(K_{h(x_0)}\) of \(h(x_0)\) in \(\Gbf_\R\) and a maximal set of strongly orthogonal non-compact roots, and then take the \(\Gbf(\R)\)-orbit.
We obtain a \(\Gbf(\R)\)-equivariant surjective map \(\Scal \to \Xcal\), and at each point of \(\Scal\) we have a morphism \(\SLbf_{2,\R}^r \to \Gbf_\R\) with finite kernel.
As in the previous section for any \(s \in \Scal\) above \(x \in \Xcal\) there is a unique morphism \(\varphi_s: H_r \longrightarrow \Gbf_\R\) (recall the group \(H_r\) from \eqref{eq:def_Hr}) satisfying
\begin{itemize}
\item the restriction to \(\SLbf_{2,\R}^r \subset H_r\) is as above,
\item for any \(z \in \mathbb{S}(\R)\) we have \(\varphi_s(z, h_{\std}(z), \dots, h_{\std}(z)) = h(x)(z)\).
\end{itemize}
The map \(s \mapsto \varphi_s\) is \(\Gbf(\R)\)-equivariant.

\begin{defi}
  A parabolic subgroup \(\Pbf\) of \(\Gbf\) is \emph{admissible} (with respect to the generalized Shimura datum \((\Xcal,h)\)) if its image in each simple factor (over \(\Q\)) \(\Hbf\) of \(\Gbf_{\ad}\) is either equal to \(\Hbf\) or is maximal among proper parabolic subgroups of \(\Hbf\), the latter being allowed only if there exists \(x \in \Xcal\) such that \(h(x)\) acts non-trivially on \(\Hbf_\R\) (by conjugation).
\end{defi}

We have a natural map from \(\Xcal\) to the complex points of a Grassmannian \(\mathrm{Gr}_\mu\):
\[ x \in \Xcal \mapsto h(x) \mapsto (\mu_{h(x)}: \GLbf_{1,\C} \to \Gbf_{\C}) \mapsto \Qbf \]
where \(\Qbf\) is the parabolic subgroup of \(\Gbf_\C\) whose Lie algebra is the sum of the nonpositive eigenspaces for the adjoint action of \(\GLbf_{1,\C}\) via \(\mu_{h(x)}\).
This embeds each connected component of \(\Xcal\) as an open subset of \(\mathrm{Gr}_\mu(\C)\).
Define \(\ol{\Xcal}\) as the disjoint union, over all connected components \(\Xcal'\) of \(\Xcal\), of the closure of \(\Xcal'\) in \(\mathrm{Gr}_\mu(\C)\).
Each such \(\Xcal'\) is a product of domains \(D\)'s as in the real adjoint case considered in the previous section.
More precisely \(\Gbf_{\ad,\R}\) decomposes as \(\prod_{i \in I} \Gbf_i\) where each \(\Gbf_i\) is a simple adjoint group over \(\R\), and there is a subset \(I_\her\) of \(I\) such that the image of \(h(x)\) (for one or any \(x \in \Xcal\)) in \(\Hom(\mathbb{S}, \Gbf_i)\) is trivial if and only if \(i \not\in I_\her\), and the map
\[ \Xcal \longrightarrow \prod_{i \in I_\her} \Hom(\mathbb{S}, \Gbf_i) \]
identifies each connected component of \(\Xcal\) with \(\prod_{i \in I_\her} D_i\) where \(D_i\) is a simple hermitian symmetric domain as in the previous case.
For \(\Pbf\) an admissible parabolic subgroup of \(\Gbf_\R\) define \(\Xcal_\Pbf\) as the disjoint union over \(\Xcal'\) of the boundary component \(\Xcal'_\Pbf\) of \(\Xcal'\) corresponding to \(\Pbf\).
From the real simple case considered before we know that \(\Xcal_\Pbf\) is stable under \(\Pbf(\R)\) and that \(\Nbf_\Pbf(\R)\) acts trivially on \(\Xcal_\Pbf\), so \(\Xcal_\Pbf\) has a natural action of \(\Mbf_\Pbf(\R)\).
We also have a natural \(\Mbf_\Pbf(\R)\)-equivariant embedding of \(\Xcal_\Pbf\) in \(\pi_0(\Xcal) \times \Hom(\mathbb{S}, \Mbf_{\Pbf,\R})\) characterized by Pink's interpretation \eqref{eq:pi_P_simple_real_can} of the geodesic projection: for any \(x \in \Xcal\) there is a unique \(\varphi_{h(x),\Pbf}: H_1 \to \Gbf_{\R}\) satisfying the same three conditions, allowing us to define
\[ \pi_\Pbf(x) := ([x], \varphi_{h(x),\Pbf}(z, \diag(z \ol{z}, 1))) \in \pi_0(\Xcal) \times \Hom(\mathbb{S}, \Mbf_{\Pbf,\R}), \]
and we have an identification of \(\Xcal_\Pbf\) with \(\pi_\Pbf(\Xcal)\).
The second projection then gives a \(\Mbf_\Pbf(\R)\)-equivariant map \(h_\Pbf: \Xcal_\Pbf \to \Hom(\mathbb{S}, \Mbf_{\Pbf,\R})\).

\begin{prop} \label{pro:ext_Shim_parabolic}
  Let \((\Gbf, \Xcal, h)\) be a generalized Shimura datum.
  Let \(\Pbf\) be an admissible parabolic subgroup of \(\Gbf\).
  Then \((\Mbf_\Pbf, \Xcal_\Pbf, h_\Pbf)\) is also a generalized Shimura datum.
\end{prop}
Compare \cite[Lemma 4.8 and Corollary 4.10]{Pink_dissertation}.
\begin{proof}
  We can argue for each factor of \(\Gbf_\ad\) separately.
  The factors centralized by all \(h(x)\) for \(x \in \Xcal\) trivially yield factors of \(\Mbf_{\Pbf,\ad}\) centralized by all \(h_\Pbf(x)\) for \(x \in \Xcal_\Pbf\).
  We are reduced to showing that for a simple adjoint group \(\Gbf\) over \(\Q\), a hermitian symmetric domain \(D = \Gbf(\R)^0/K\) and a boundary component \(F_\Pbf\) corresponding to an admissible parabolic subgroup of \(\Gbf_\R\) defined over \(\Q\), for each simple factor \(\Hbf\) of \(\Mbf_{\Pbf,\ad}\) we have either
  \begin{itemize}
  \item for all \(x \in D\) conjugation by \(\pi_\Pbf(x)\) acts trivially on \(\Hbf_\R\), or
  \item for all \(x \in D\) conjugation by \(\pi_\Pbf(x)\) on \(\Lie \Hbf_\R\) is of type \((0,0), (1,-1), (-1,1)\) and the involution \(\Ad \pi_\Pbf(x)(i)\) of \(\Hbf_\R\) is a Cartan involution.
  \end{itemize}
  This follows from \cite[\S III.3.5]{AMRT}: \(\Mbf_{\Pbf,\ad}\) has a decomposition \(\Mbf_{\Pbf,\ad,\lin} \times \Mbf_{\Pbf,\ad,\her}\) where \(\Mbf_{\Pbf,\ad,\lin,\R}\) (corresponding to \(\Gcal_\ell(F)\) loc.\ cit.) is centralized by all \(\pi_\Pbf(x)\) and \(F\) is isomorphic to the quotient of \(\Mbf_{\Pbf,\ad,\her}(\R)^0\) by a maximal compact subgroup (up to compact factors \(\Mbf_{\Pbf,\ad,\her,\R}\) is the image in \(\Mbf_{\Pbf,\ad,\R}\) of the product over all simple factors of \(\Gbf_\R\) of the group denoted by \(\Lbf_{f,S}\) in the real simple case).
\end{proof}

Let \(\Xcal^*\) be the union of the components \(\Xcal_\Pbf\), for \(\Pbf\) ranging over all admissible subgroups of \(\Gbf\).
We have an action of \(\Gbf(\Q)\) on \(\Xcal^*\), with \(g \in \Gbf(\Q)\) mapping \(\Xcal_{\Pbf}\) to \(\Xcal_{\Ad(g) \Pbf}\), in particular the stabilizer of \(\Xcal_\Pbf\) is \(\Pbf(\Q)\).
We do not recall the definition of the Satake topology on \(\Xcal^*\), see \cite[\S III.6]{AMRT}.
For a neat compact open subgroup \(K\) of \(\Gbf(\A_f)\) the minimal (or Satake-Baily-Borel) compactification \(\Sh(\Gbf,\Xcal,K)^*(\C)\) of \(\Sh(\Gbf,\Xcal,K)(\C)\) is defined as
\[ \Gbf(\Q) \backslash \Xcal^* \times \Gbf(\A_f)/K \]
and as the notation suggests may naturally be identified with the complex points of a projective variety over \(\C\).
The decomposition (over admissible parabolic subgroups) \(\Xcal^* = \bigsqcup_\Pbf \Xcal_\Pbf\) corresponds to a stratification (by locally closed analytic subsets)
\begin{align*}
  \Sh(\Gbf,\Xcal,K)^*(\C)
  &= \bigsqcup_{[\Pbf]} \Pbf(\Q) \backslash \left( \Xcal_\Pbf \times \Gbf(\A_f)/K \right) \\
  &= \bigsqcup_{[\Pbf]} \underset{gK \in [\Pbf(\A_f) \curvearrowright \Gbf(\A_f)/K]}{\colim} \Sh(\Mbf_\Pbf, \Xcal_\Pbf, \pi(gKg^{-1} \cap \Pbf(\A_f)))
\end{align*}
where the disjoint union ranges over \(\Gbf(\Q)\)-conjugacy classes of admissible parabolic subgroups of \(\Gbf\) (equivalently, admissible parabolic subgroups containing a fixed minimal parabolic subgroup of \(\Gbf\)) and \(\pi: \Pbf \to \Mbf_\Pbf\) is the canonical projection.
The second equality follows from the fact that \(\Nbf_\Pbf(\Q)\) acts trivially on \(\Xcal_\Pbf\) and is dense in \(\Nbf_\Pbf(\A_f)\), so we have well-defined homeomorphisms
\begin{align*}
  \Mbf_\Pbf(\Q) \backslash \left( \Xcal_P \times \Mbf(\A_f) / \pi(gKg^{-1} \cap \Pbf(\A_f)) \right) & \longrightarrow \Pbf(\Q) \backslash \left( \Xcal_\Pbf \times \Pbf(\A_f) g K / K \right) \\
  [x, m] & \longmapsto [x, \tilde{m}g]
\end{align*}
where \(\tilde{m} \in \Pbf(\A_f)\) is any lift of \(m \in \Mbf_\Pbf(\A_f)\).
The functor implied in the colimit maps an isomorphism \(gK \xrightarrow{p \in \Pbf(\A_f)} pgK\) to
\[ \Sh(\Mbf_\Pbf, \Xcal_\Pbf, \pi(gKg^{-1} \cap \Pbf(\A_f))) \xrightarrow{T_{\pi(p)^{-1}}} \Sh(\Mbf_\Pbf, \Xcal_\Pbf, \pi(pgK(pg)^{-1} \cap \Pbf(\A_f))), \]
which is the identity morphism when \(pgK = gK\), so each colimit may be (non-canonically) identified with a disjoint union over \(\Pbf(\A_f) \backslash \Gbf(\A_f) / K\).
We denote
\[ i_{\Pbf,gK}: \Sh(\Mbf_\Pbf, \Xcal_\Pbf, \pi(gKg^{-1} \cap \Pbf(\A_f))) \hookrightarrow \Sh(\Gbf, \Xcal, K)^* \]
the locally closed immersion.
The stratification described in \cite[\S 6.3]{Pink_dissertation} is simply obtained by spelling out the generalized Shimura varieties as colimits.
By \S 12.3 loc.\ cit.\ there is a unique model of \(\Sh(\Gbf, \Xcal, K)^*\) over the reflex field \(E\) extending the canonical model of \(\Sh(\Gbf, \Xcal, K)\), each stratum is also defined over \(E\) and is identified with the canonical model of \(\Sh(\Mbf_\Pbf, \Xcal_\Pbf, \pi(gKg^{-1} \cap \Pbf(\A_f)))\).

\subsubsection{Iterated boundary strata}
\label{sec:gen_Shim_boundary_boundary}

We recall how boundary strata of boundary strata map to boundary strata.
We will need the following notion.

\begin{prodef} \label{prodef:linres_ext_Shim}
  Let \((\Gbf, \Xcal, h)\) be a generalized Shimura datum, with corresponding decomposition \(\Gbf_\ad = \Gbf_{\ad,\lin} \times \Gbf_{\ad,\her}\).
  Let \(\Pbf\) be a parabolic subgroup of \(\Gbf\).
  Assume that the image of \(\Pbf\) in \(\Gbf_{\ad,\her}\) is \(\Gbf_{\ad,\her}\).
  We obtain a generalized Shimura datum \((\Mbf_\Pbf, \Xcal, h)\) by first restricting the action of \(\Gbf(\R)\) on \(\Xcal\) to \(\Pbf(\R)\) and observing that \(\Nbf_\Pbf(\R)\) acts trivially and that for any \(x \in \Xcal\) the morphism \(h(x): \mathbb{S} \to \Gbf_\R\) factors through \(\Pbf_\R\).

  For \(K\) a neat compact open subgroup of \(\Gbf(\A_f)\) we have a unique morphism of schemes over \(E\)
  \[ \Sh(\Mbf_\Pbf, \Xcal, \pi(K \cap \Pbf(\A_f))) \longrightarrow \Sh(\Gbf, \Xcal, K) \]
  given on complex points by
  \begin{align*}
    & \Mbf_\Pbf(\Q) \backslash \left( \Xcal \times \Mbf_\Pbf(\A_f) / \pi(K \cap \Pbf(\A_f)) \right) \\
    \xleftarrow{\sim} & \Pbf(\Q) \backslash \left( \Xcal \times \Pbf(\A_f) / K \cap \Pbf(\A_f) \right) \\
    \to & \Gbf(\Q) \backslash \left( \Xcal \times \Gbf(\A_f) / K \right)
  \end{align*}
  This morphism is finite étale.
\end{prodef}
\begin{proof}
  The fact that \(\Mbf_{\Pbf}(\R)\) acts transitively on \(\Xcal\) easily follows from the Iwasawa decomposition of \(\Gbf(\R)\), so \((\Mbf_\Pbf, \Xcal, h)\) is a generalized Shimura datum.
  The remaining statements are clear once we go back to the definition of generalized Shimura varieties as colimits of Shimura varieties (Definition \ref{def:gen_Shim_var}).
\end{proof}

Let \((\Gbf, \Xcal, h)\) be a generalized Shimura datum and \(\Pbf\) an admissible parabolic subgroup of \(\Gbf\).
The embedding of \(\Xcal_\Pbf\) in \(\Xcal^*\) extends to \((\Xcal_\Pbf)^* \hookrightarrow \Xcal^*\) (or even \(\ol{\Xcal_\Pbf} \hookrightarrow \ol{\Xcal}\) which is a homeomorphism onto a closed subset of \(\ol{\Xcal}\)).
We obtain for any neat compact open subgroup \(K\) of \(\Gbf(\A_f)\) and any \(gK \in \Gbf(\A_f)/K\) a map
\[ \Sh(\Mbf_\Pbf, \Xcal_\Pbf, \pi(gKg^{-1} \cap \Pbf(\A_f)))^*(\C) \xrightarrow{\ol{i_{\Pbf,gK}}} \Sh(\Gbf, \Xcal, K)^*(\C) \]
extending (the complexification of) the embedding \(i_{\Pbf,gK}\).
There is a unique morphism of schemes over \(E\) inducing \(\ol{i_{\Pbf,gK}}\), that we still denote by \(\ol{i_{\Pbf,gK}}\).
This morphism is finite because it is a morphism between proper schemes over \(E\) whose fibers at closed points are finite.

Let us describe the restriction of \(\ol{i_{\Pbf,gK}}\) to a boundary stratum of the source.
For an admissible parabolic subgroup \(\Qbf\) of \(\Mbf_\Pbf\) there is a unique admissible parabolic subgroup \(\Pbf'\) of \(\Gbf\) such that the preimage of \(\Qbf\) in \(\Pbf\) is \(\Pbf \cap \Pbf'\).
Let \(\Qbf'\) be the image of \(\Pbf \cap \Pbf'\) in \(\Mbf_{\Pbf'}\).
The generalized Shimura datum \((\Mbf_\Qbf, (\Xcal_\Pbf)_\Qbf, (h_\Pbf)_\Qbf)\) may be identified with the restriction (in the sense of Proposition-Definition \ref{prodef:linres_ext_Shim}) of the generalized Shimura datum \((\Mbf_{\Pbf'}, \Xcal_{\Pbf'}, h_{\Pbf'})\) to \(\Mbf_{\Qbf'}\).
More precisely the subsets \((\Xcal_\Pbf)_\Qbf\) and \(\Xcal_{\Pbf'}\) of \(\pi_0(\Xcal) \times \Hom(\mathbb{S}, \Mbf_{\Pbf \cap \Pbf', \R})\) are equal: this may be checked using \eqref{eq:pi_P_simple_real}.
Denote \(K(\Pbf,gK) = \pi_\Pbf(gKg^{-1} \cap \Pbf(\A_f))\).
For \(g' K(\Pbf,gK) \in \Mbf_\Pbf(\A_f) / K(\Pbf,gK)\) the composition
\begin{align*}
  & \Sh(\Mbf_\Qbf, (\Xcal_\Pbf)_\Qbf, \pi_\Qbf(g' K(\Pbf,gK) (g')^{-1} \cap \Qbf(\A_f))) \\
  \xrightarrow{i_{\Qbf,g'K(\Pbf,gK)}} \ & \Sh(\Mbf_\Pbf, \Xcal_\Pbf, K(\Pbf,gK))^* \\
  \xrightarrow{\ol{i_{\Pbf,gK}}} & \Sh(\Gbf, \Xcal, K)^*
\end{align*}
is equal to the composition
\begin{align*}
  & \Sh(\Mbf_{\Qbf'}, \Xcal_{\Pbf'}, \pi_{\Pbf \cap \Pbf'}(g' g K (g' g)^{-1} \cap (\Pbf \cap \Pbf')(\A_f))) \\
  \to \ & \Sh(\Mbf_{\Pbf'}, \Xcal_{\Pbf'}, \pi_{\Pbf'}(g' g K (g' g)^{-1} \cap \Pbf'(\A_f))) \\
  \xrightarrow{i_{\Pbf',g'gK}} \ & \Sh(\Gbf, \Xcal, K)^*
\end{align*}
where the first map is the finite étale map introduced in Proposition-Definition \ref{prodef:linres_ext_Shim}.

We have the following analogue of \cite[Proposition 1.1.3]{MorelBook}.
First choose an order \(\{\Pbf_1, \dots, \Pbf_n\}\) on the set of standard maximal proper parabolic subgroups of \(\Gbf\) mapping onto \(\Gbf_{\ad,\lin}\) satisfying \(\Ubf_i \subsetneq \Ubf_{i+1}\) for \(1 \leq i < n\), where \(\Ubf_i\) is the center of the unipotent radical of \(\Pbf_i\) (see \cite[Theorem III.4.8 (i)]{AMRT}).
We have a bijection \(I \mapsto \Pbf_I := \bigcap_{i \in I} \Pbf_i\) between subsets of \(\{1, \dots, n\}\) and standard parabolic subgroups of \(\Gbf\) mapping onto \(\Gbf_{\ad,\lin}\) (these are in bijection with standard parabolic subgroups of \(\Gbf_{\ad,\her}\)).
For \(I = \{ i_1 < \dots < i_r \}\) a non-empty subset of \(\{1,\dots,n\}\) we have a bijection between \(\Pbf_I(\A_f) \backslash \Gbf(\A_f) / K\) and the set of tuples \((S_1, \dots, S_r)\) where \(S_1\) is a boundary stratum of \(\Sh(\Gbf, \Xcal, K)^*\) corresponding to \(\Pbf_{i_1}\) and for \(1 \leq j < r\) \(S_{j+1}\) is a boundary stratum of \(S_j^*\) corresponding to the image \(\Qbf_{j+1}\) of \(\Pbf_{i_1} \cap \dots \cap \Pbf_{i_{j+1}}\) in \(\Mbf_{\Pbf_{i_1} \cap \dots \cap \Pbf_{i-1}}\).
Under this bijection the \(\Pbf_I(\A_f)\)-orbit of \(gK \in \Gbf(\A_f)/K\) corresponds to
\[ S_1 := \Sh(\Mbf_{\Pbf_{i_1}}, \Xcal_1, K(\Pbf_{i_1}, gK)) \xrightarrow{i_{\Pbf_{i_1},gK}} \Sh(\Gbf, \Xcal, K)^* \]
\[ S_{j+1} := \Sh(\Mbf_{\Qbf_{j+1}}, \Xcal_{j+1}, K_{j+1}) \xrightarrow{i_{\Qbf_{j+1},K_j}} S_j^* \]
where \(\Xcal_1 = \Xcal_{\Pbf_{i_1}}\), \(\Xcal_{j+1} = (\Xcal_j)_{\Qbf_j}\), \(K_1 = K(\Pbf_{i_1}, gK)\) and \(K_{j+1} = K(\Qbf_{j+1}, K_j)\).

\newpage

\printbibliography

\end{document}